\documentclass[]{amsart}

\usepackage[letterpaper, hmarginratio=1:1]{geometry}
\usepackage{stmaryrd}
\usepackage{graphicx} 
\usepackage[utf8]{inputenc}
\usepackage{amsmath,amsthm,amsfonts,amssymb,bm,bbm,mathrsfs,comment}
\usepackage{tikz}
\usepackage{ytableau}

\usepackage{mathtools}

\usepackage{caption}
\usepackage{subcaption}
\usepackage{enumerate}
\usepackage{hyperref}
\usepackage{cleveref}

\usepackage{soul}

\numberwithin{equation}{section}

\theoremstyle{plain}
\newtheorem{thm}{Theorem}[section]
\newtheorem{theorem}{Theorem}[section]
\newtheorem*{thm*}{Theorem}
\newtheorem{prop}[thm]{Proposition}
\newtheorem*{prop*}{Proposition}

\newtheorem*{conj*}{Conjecture}
\newtheorem{cor}[thm]{Corollary}
\newtheorem{lem}[thm]{Lemma}
\newtheorem*{lem*}{Lemma}

\newtheorem{innerthm}{Theorem}
\newenvironment{theorema}[1]
{\renewcommand\theinnerthm{#1}\innerthm}
{\endinnerthm}

\theoremstyle{definition}

\newtheorem{definition}[thm]{Definition}
\newtheorem{remark}[thm]{Remark}
\newtheorem{ass}[thm]{Assumption}

\newcommand{\varrhobetabf}{\boldsymbol{\varrho}_{\beta}}

\renewcommand{\i}{\mathrm{i}}
\newcommand{\T}{\mathbf{T}}
\newcommand{\G}{\mathbf{G}}
\newcommand{\M}{\mathbf{M}}
\renewcommand{\L}{\mathbf{L}}

\renewcommand{\S}{\mathbf{S}}

\newcommand{\ve}{v_{\mathsf{eff}}}

\DeclareMathOperator{\Id}{Id}
\DeclareMathOperator{\Real}{Re}
\DeclareMathOperator{\Imaginary}{Im} 
\DeclareMathOperator{\dr}{dr}
\DeclareMathOperator{\eig}{eig}
\DeclareMathOperator{\supp}{supp}
\DeclareMathOperator{\Tr}{Tr}

\DeclareMathOperator{\var}{Var}
\DeclareMathOperator{\Var}{Var}

\DeclareMathOperator{\Cov}{\mathrm{Cov}}
\DeclareMathOperator{\sgn}{sgn}
\DeclarePairedDelimiter{\floor}{\lfloor}{\rfloor}
\DeclareMathOperator{\dist}{dist}
\DeclareMathOperator{\op}{op}

\title{Fluctuations for the Toda lattice}
\author{Amol Aggarwal and Matthew Nicoletti}

\begin{document}

\begin{abstract}

   In this paper we consider the Toda lattice $(\bm{p}(t);\bm{q}(t))$ at thermal equilibrium, meaning that its variables $(p_j)$ and $(e^{q_j - q_{j+1}})$ are independent Gaussian and Gamma random variables, respectively. We show under diffusive scaling that the space-time fluctuations for the model's currents converge to an explicit Gaussian limit. As consequences, we deduce, (i) the scaling limit for the trajectory of a single particle $q_0$ is a Brownian motion; (ii) space-time two-point correlation functions for the model decay inversely with time, with explicit scaling distributions predicted by Doyon \cite{ECIS} and Spohn \cite{Spo20}. Our starting point is the notion that the Toda lattice can be thought of as a dense collection of many ``quasi-particles'' that interact through scattering. The core of our work is to establish that the full joint scaling limit of the fluctuations for these quasi-particles is given by a Gaussian process, called a dressed L\'{e}vy--Chentsov field.

\end{abstract}

\maketitle

\tableofcontents

\section{Introduction}

\label{Introduction}

\subsection{Overview} 

A central question in mathematical physics is to describe the long-time behavior of many-body Hamiltonian systems. An important aspect of this is to determine universal scaling limits for their space-time fluctuations at equilibrium. Much has now been established in this direction if the system's density of particles is vanishingly low; see the review of Bodineau--Gallagher--Saint-Raymond--Simonella \cite{DDG}. However, little is known if this density is strictly positive, even in one spatial dimension.

The broad belief in this context is that the limiting fluctuations for these systems should separate into universality classes, depending on whether the underlying dynamics are chaotic or integrable. In the chaotic case, physics works of van Beijeren \cite{vBe12} and Spohn \cite{Spo14} posit Kardar--Parisi--Zhang (KPZ) \cite{DSGI} type fluctuations; after a long time $T \gg 1$, certain space-time two-point correlation functions should decay as $T^{-2/3}$ along particular characteristic directions, and the system's fluctuations should be of order $T^{1/3}$ and scale to non-Gaussian limits. For various stochastic interacting particle systems, this (and substantially more) is now established; see the surveys of Corwin \cite{Cor16} and Quastel \cite{Qua11}. However, proving such phenomena for any deterministic Hamiltonian system remains open.

In the integrable case, physics works of De Nardis--Bernard--Doyon \cite{BDD18,BDD19} and Gopalakrishnan--Huse--Khemani--Vasseur \cite{GHKV18} instead predict that diffusive fluctuations typically arise; after a long time $T \gg 1$, space-time two-point correlation functions should decay as $T^{-1}$ along all directions, and the system's fluctuations should be of order $T^{1/2}$ and scale to a non-trivially correlated Gaussian process. The only model for which such phenomena has previously been proven (as stochastic dynamics no longer serve as suitable proxies) is the hard rods model; see the earlier works of Boldrighini--Dobrushin--Suhov \cite{OHRH,OHRHLS}, Spohn \cite{HETCF}, Presutti--Wick \cite{MSFMS}, and the more recent developments by Ferrari--Franceschini--Grevino--Spohn \cite{HRHF} and Ferrari--Olla \cite{FO25}.

While the hard rods model provided a useful initial test case, it has simplifying features (such as conservation of each individual particle's momentum, and a constant scattering shift) that raise the question of whether its fluctuation behavior fully captures that of Hamiltonian integrable systems more broadly. This warrants understanding the above predictions for a more representative such system with richer interactions. The purpose of this paper is to do so for the Toda lattice. \\

The \emph{Toda lattice} \cite{Tod67} is a Hamiltonian system $(\bm{p}(t); \bm{q}(t))$, where $\bm{p}(t) = (p_i(t))$ and $\bm{q}(t) = (q_i(t))$ are indexed by $i \in \mathbb{Z}$, evolving under the equations
\begin{flalign*}
	\partial_t q_i (t) = p_i (t), \quad \text{and} \quad \partial_t p_i (t) = e^{q_{i-1} (t) - q_i (t)} - e^{q_i (t) - q_{i+1} (t)}.
\end{flalign*}

\noindent As these dynamics are invariant under a simultaneous shift of all $(q_i)$, we normalize $q_0(0) = 0$. This model may be thought of as a system of particles moving on the real line, with locations $(q_i)$ and momenta $(p_i)$. Since the works of Flaschka \cite{Fla74} and Manakov \cite{Man74}, which exhibited its full set of conserved quantities, and that of Moser \cite{Mos75}, which determined its scattering shift, the Toda lattice has become recognized as an archetypal example of an integrable system.

The long-time behavior of the Toda lattice has been extensively studied when its initial data either decays or approximates an almost everywhere smooth profile (so that its solutions have a well-defined background). In these cases, detailed results have been attained (often by analyzing associated Riemann--Hilbert problems) by Venakides--Deift--Oba \cite{DOV91}, Deift--Kamvissis--Kriecherbauer--Zhou \cite{DKKZ96}, Deift--McLaughlin \cite{CLL}, Bloch--Golse--Paul--Uribe \cite{DAO}, and Kr\"{u}ger--Teschl \cite{LADI}.

In this work, we study the Toda lattice under equilibrium initial conditions, which are necessarily rough and fully supported on $\mathbb{Z}$. Specifically, we consider the model under perhaps its most natural invariant measure,\footnote{This is a random initial condition $(\bm{p}(0); \bm{q}(0))$, so that the law of $(\bm{p}(t); \bm{q}(t))$ equals that of $(\bm{p}(0); \bm{q}(0))$ for all $t$.} sometimes called \emph{thermal equilibrium}, which is when the $(p_j)$ and $(e^{(q_j - q_{j+1})/2})$ are independent Gaussian and chi random variables, respectively. Fixing two parameters $\beta, \theta > 0$, this means that we sample each $p_j$ and $e^{(q_j - q_{j+1})/2}$ independently, under the densities $C_{\beta} e^{-\beta x^2/2}$ and $C_{\beta,\theta} x^{2\theta-1} e^{-\beta x^2/2} \cdot \mathbbm{1}_{x>0}$, respectively (where $C_{\beta}$ and $C_{\beta,\theta}$ are normalization constants).

Now we state several sample results. To lighten the notation in this introductory section, these results may be written informally (often without giving full definitions or equations) or in a special case. However, we will provide specific references to Section \ref{Results} for the precise statements in full generality. Throughout, we will assume that $\theta$ is sufficiently small with respect to $\beta$.\footnote{This corresponds to assuming that density of particles in the system is not too large (but is still strictly positive). Its main source is our use of certain estimates from \cite{Agg25}, where this assumption was imposed to verify that a specific matrix was (quantitatively) invertible; see \cite[Remark 6.7]{Agg25}. We expect this invertibility to hold for all $\theta$, so it is likely that this condition is artificial. However, we will make no effort to remove it here.}

Our first result identifies the scaling limit for the trajectory of a single particle, say $q_0$, in the Toda lattice, as a Brownian motion.

\begin{theorema}{A}[Corollary \ref{cor:q0bm}]
	
	\label{q0t00}
	
	The process $T^{-1/2} \cdot q_0 (T\tau)$ converges to a Brownian motion in $\tau$ with an explicit variance (given by \eqref{eqn:q0var} below), as $T$ tends to $\infty$.
\end{theorema}

Next, it is natural to ask for the asymptotics of two-point correlation functions after a long time $T \gg 1$. Indeed, these were initially sought after in the physics literature soon after the original discovery of the Toda lattice \cite{ESLMD,DCL,CSM}, but a conjectural answer did not emerge at the time. More recently, Doyon \cite{ECIS} and then Spohn \cite{Spo20} predicted that these correlations decay as $T^{-1}$ with an explicit scaling function, which was later tested numerically to high accuracy by Mazzuca--Grava--Kriecherbauer--McLaughlin--Mendl--Spohn \cite{ESCL}. Corollary \ref{cor:twopoint_intro} below verifies this prediction in generality (for arbitrary local charges). Here we only informally state a special case of this result, addressing the two-point correlation function $\mathbb{E} [p_0 (0) p_{\lfloor xT\rfloor} (T)]$ between particle momenta.

\begin{theorema}{B}[Corollary \ref{cor:twopoint_intro}, special case]
	\label{expectationp}
	
	The two-point correlation function $T \cdot \mathbb{E} [p_0 (0) p_{\lfloor xT \rfloor} (T)]$ converges, as a distribution in $x$, to an explicit scaling function (given by the $(m,n)=(1,1)$ case of \eqref{eqn:two_pt_asymp_intro} below), as $T$ tends to $\infty$.
	
\end{theorema}

Theorems \ref{q0t00} and \ref{expectationp} will be consequences of a broader space-time fluctuation scaling limit result for currents of the Toda lattice. Theorem \ref{thm:current_fluct_intro} below provides the general form of this result (for currents of arbitrary local charges); here we only informally state a special case of it, addressing the particle current. Below, given real numbers $t \ge 0$ and $x, y \in \mathbb{R}$, let $J_t (x, y)$ denote the flux of particles across the interval $[x,y]$ after time $t$, namely, $J_t (x,y) = \# \{ j : q_j (0) < x, q_j (t) \ge y \} - \# \{ j : q_j (0) \ge x, q_j (t) < y \}$.

\begin{theorema}{C}[Theorem \ref{thm:current_fluct_intro}, special case]
	\label{limitj}
	The process $T^{-1/2} \cdot (J_{T\tau} (T \mathfrak{q}, T\mathfrak{q}') - \mathbb{E}[J_{T\tau} (T \mathfrak{q}, T\mathfrak{q}')])$ converges to an explicit Gaussian limit in $(\tau, \mathfrak{q}, \mathfrak{q}')$ (given by \eqref{w1} below), as $T$ tends to $\infty$.  
\end{theorema}

We will deduce Theorem \ref{limitj} (and thus Theorems \ref{q0t00} and \ref{expectationp}) from a ``quasi-particle'' fluctuation result, given by Theorem \ref{qjt0} below. Before setting up the latter, let us comment on how the above three theorems relate to the literature on kinetic theory and stochastic dynamics. 

First, for the hard spheres model under the low-density, Boltzmann--Grad limit, the direct counterparts (witnessing Gaussian limits under similar scaling) of Theorems \ref{q0t00}, \ref{expectationp}, and \ref{limitj} were established by Bodineau--Gallagher--Saint-Raymond--Simonella in \cite{MLDS}, \cite{LCGE}, and \cite{LDEF}, respectively. For the weakly nonlinear Schr\"{o}dinger equation, the counterpart of Theorem \ref{expectationp} was proven by Lukkarinen--Spohn in \cite{WNER} (where the question of understanding two-point correlation functions for strongly nonlinear integrable systems, such as the Toda lattice, was incidentally also mentioned).

At positive density and strong interaction, the above results had only been shown before for the hard rods model; the analogs of Theorems \ref{expectationp} and \ref{limitj} were proven in \cite{HETCF} and \cite{MSFMS,OHRHLS,FO25}, respectively. Partial results have also been obtained for the box-ball model, an integrable cellular automaton. For this system, Croydon--Kato--Sasada--Tsujimoto \cite{DBSR} showed tightness of a particle's trajectory under diffusive scaling (in the direction of Theorem \ref{q0t00}), and identified the particle current fluctuations at a specific space-time point (a marginal for the counterpart of Theorem \ref{limitj} at this point).

For various stochastic interacting particle systems, such as the totally asymmetric simple exclusion process (TASEP), much is now known. After the initial current fluctuation results for the TASEP by Johansson \cite{Joh00}, Ferrari, Pr\"{a}hofer, and Spohn \cite{ESFG,SLSTC} showed that its two-point function (as in Theorem \ref{expectationp}) decay as $T^{-2/3}$ along a characteristic line, and limit to a universal scaling function. Matetski--Quastel--Remenik \cite{TF} later showed that its full space-time current fluctuations (as in Theorem \ref{limitj}) are of order $T^{1/3}$ and rescale to a non-Gaussian process called the KPZ fixed point.

Theorems \ref{q0t00}, \ref{expectationp}, and \ref{limitj} give the fluctuation scaling limit for the Toda lattice and show it differs (both in scaling exponents and limiting distributions) from that of stochastic dynamics such as the TASEP, in agreement with the predictions described earlier. This is the first such result for a dense, strongly interacting, deterministic Hamiltonian system beyond the hard rods model. 

\subsection{Quasi-particles} 

\label{ParticleQ}

The starting point for our analysis is the notion that the Toda lattice can be viewed as a dense collection of objects called \emph{quasi-particles}, each possessing a time-independent \emph{spectral parameter} $\lambda_j$ and a time-dependent \emph{location} $Q_j (t)$. To define them, it is convenient\footnote{It is also possible to set up this formalism directly for the Toda lattice on $\mathbb{Z}$, without changing any of the asymptotic results to be stated below. This is implemented in Appendix \ref{ParticlesInfinite} (see Definition \ref{qlambdanlimit0}).} to truncate the index $j \in \mathbb{Z}$, of the Toda lattice variables $(p_j; q_j)$, to a long interval $[N_1, N_2]$.

As the quasi-particle spectral parameters are time-independent, they are given by the conserved quantities of the Toda lattice, which are the eigenvalues $(\lambda_1, \lambda_2, \ldots , \lambda_N)$ of its \emph{Lax matrix}. This is the tridiagonal, symmetric matrix $\mathbf{L}(t) = [L_{ij}(t)]$, whose diagonal and off-diagonal entries are the $(p_i(t))$ and $(e^{(q_i(t) - q_{i+1}(t))/2})$, respectively. When $N = N_2 - N_1 + 1$ is large and $(\bm{p} (t); \bm{q} (t))$ is random, $\mathbf{L}(t)$ becomes a high-dimensional random matrix, whose eigenvalue density prescribes the distribution of quasi-particle spectral parameters in the Toda lattice under thermal equilibrium. After initial work of Opper \cite{ASCEC}, Spohn \cite{Spo20GGE} predicted formulas for this limiting density (under a broader class of invariant measures), by relating random Lax matrices to high temperature $\beta$ ensembles. These predictions were proven in works of Grava, Guionnet, Kozlowski, Little, Mazzuca, and Memin \cite{Maz22,GM22,MM24,LDPC}, with \cite{GM22,LDPC} also providing large deviations principles for these densities. 

A definition for the quasi-particle positions $(Q_j(t))$ was provided in \cite{Agg25a}, in terms of the eigenvectors of $\mathbf{L}(t)$.  To explain it (simplifying slightly; see Definition \ref{def:quasi_intro} below for its full description), let $\mathbf{u}_j (t) = (u_j (i;t))_{i \in [N_1, N_2]}$ denote a unit eigenvector of $\mathbf{L}(t)$ with eigenvalue $\lambda_j$. Results on random tridiagonal matrices due to Kunz--Souillard \cite{SSO} and Schenker \cite{LRBM}   imply that, if $\mathbf{L}(t)$ is under thermal equilibrium, then $\mathbf{u}_j (t)$ is \emph{exponentially localized}. This means that it admits some ``center'' $\varphi_t (j) \in [N_1, N_2]$ such that $|u_j (i; t)| \le C e^{-c|i-\varphi_t (j)|}$ likely holds for any $i \in [N_1, N_2]$. We view this center $\varphi_t (j)$ as the index of the Toda particle associated with the $j$-th quasi-particle. So, we define the $j$-th quasi-particle's location on $\mathbb{R}$ to be this particle's position $Q_j (t) = q_{\varphi_t (j)} (t)$.

In \cite{Agg25a}, it was proven that the above quasi-particle data $(\lambda_j, Q_j (t))$ has two properties that could potentially be used together to study asymptotics for the Toda lattice. First, natural  quantities describing the original Toda lattice, such as currents, are closely approximable by simple expressions of the quasi-particles (for example, the number of Toda particles in an interval is equal to the number of quasi-particles in it). Second, these quasi-particles satisfy the approximate equation
\begin{flalign}
	\label{qktqk0} 
	\begin{aligned} 
		Q_k (t) = Q_k (0) + \lambda_k & t + 2 \displaystyle\sum_{j=1}^N \log |\lambda_k - \lambda_j| \cdot (\mathbbm{1}_{Q_j (0) < Q_k (0)} - \mathbbm{1}_{Q_j(t) < Q_k (t)}) + O( (\log N)^C). 
	\end{aligned} 
\end{flalign}

In the physics literature, the concept of viewing many-body integrable systems as dense collections of quasi-particles was posited (though without providing a microscopic definition of their locations $(Q_j)$) by Bertini--Collura--De Nardis--Fagotti \cite{TEP} and Castro-Alvaredo--Doyon--Yoshimura \cite{EHSOE}, and variants of \eqref{qktqk0} were widely predicted by Doyon--H\"{u}bner--Yoshimura \cite{NCIS,DGHG,GDCF}. One may informally interpret \eqref{qktqk0} as follows. The $k$-th quasi-particle $Q_k$ moves with velocity $\lambda_k$ until it meets another one $Q_j$. At that moment, $Q_k$ instantaneously moves forward or backward by $2 \log |\lambda_k - \lambda_j|$, depending on whether it met $Q_j$ from the right or left, respectively (this could make it pass another quasi-particle, producing a ``cascade'' of such interactions). Then, $Q_k$ proceeds at velocity $\lambda_k$ until meeting another quasi-particle, when the procedure repeats. Here, $2 \log |\lambda_k - \lambda_j|$ is the Toda lattice's scattering shift, describing the phase displacement of just two quasi-particles passing through each other. As observed by Zakharov \cite{KES} and El--Kamchatnov--Pavlov--Zykov \cite{TLE,KESG}, this description parallels the collision structure between solitons in classical integrable systems. Indeed, solitons also proceed independently until they get close, after which they undergo an intricate interaction, and then emerge as independent solitons. Their shapes (spectral parameters) are the same as before the interaction, but their locations are translated (relative to their free evolutions) by a scattering shift. Thus, we refer to \eqref{qktqk0} as the \emph{asymptotic scattering relation}.\footnote{It has also been termed the ``flea-gas algorithm'' \cite{SGGH} or ``semiclassical Bethe equations'' \cite{DGHG}, since it was predicted by taking a semi-classical limit of the Bethe equations underlying quantum integrable systems.}

Given the definition of the $(Q_j)$ and that they asymptotically determine the Toda lattice currents, the next question is to understand the limiting behavior of these quasi-particles. In \cite{Agg25}, it was proven by analyzing \eqref{qktqk0} that a linear trajectory governs their law of large numbers, namely, 
\begin{flalign}
	\label{qjt} 
	Q_j (T) = Q_j (0) + T \cdot \ve (\lambda_j) + O( T^{1/2} (\log T)^C), 
\end{flalign}

\noindent likely holds, where the \emph{effective velocity} $\ve: \mathbb{R} \rightarrow \mathbb{R}$ is an explicit function (see Definition \ref{def:v_eff_intro} below) predicted by Doyon \cite{GHCS,Doy20} and Spohn \cite{Spo23}. However, \eqref{qjt} is insufficient to obtain fine distributional results (Theorems \ref{q0t00}, \ref{expectationp}, and \ref{limitj}) for the original Toda lattice. To do the latter, we must further pinpoint the limiting fluctuations for these quasi-particles.  \\

The joint scaling limit for the fluctuations of all quasi-particles in time is given by the following theorem (stated informally here), which we view as the core of this paper; Theorems \ref{q0t00}, \ref{expectationp}, and \ref{limitj} will follow as consequences of it. 

\begin{theorema}{D}[Theorem \ref{thm:quasi_fluct_intro} and Corollary \ref{q200}]
	
	\label{qjt0} 
	
	The fluctuations $T^{-1/2} \cdot (Q_j (T \tau) - Q_j (0) - T \tau \cdot \ve (\lambda_j))$ of all quasi-particles converge jointly to an explicit Gaussian process called a dressed L\'{e}vy--Chentsov field (given by Definition \ref{def:Z_intro} below), as $T$ tends to $\infty$. 
	
\end{theorema}

Fluctuation results as in Theorem \ref{qjt0} are also known versions of the Navier--Stokes corrections (to Euler scale ballistic transport) \cite[Chapter 15]{Spo23}. We have not seen the full dressed L\'{e}vy--Chentsov field describing these corrections, in the precise form written in Definition \ref{def:Z_intro}, appear previously. However, asymptotics for closely related correlation functions were predicted earlier in the physics literature in \cite{BDD18,BDD19,GHKV18}, including a formula for the limiting diffusivity of a tracer quasi-particle. We recover the latter from the dressed L\'{e}vy--Chentsov field as Proposition \ref{prop:tracer_fluct_intro} below. 

The term, ``dressed L\'{e}vy--Chentsov field,'' is used since it is a variant of a L\'{e}vy--Chentsov field originating in \cite{PSM,MSP}, composed with a nontrivial ``dressing operator'' (see Definition \ref{def:dress_intro} and Remark \ref{rmk:LCrmk} below). Forms of Theorem \ref{qjt0} for the hard rods model were proven in \cite{HRHF,FO25}. There, it is the ordinary L\'{e}vy--Chentsov field that arises as a scaling limit, as the associated dressing operator is multiplication by a scalar (essentially since the scattering shift for the hard rods model is constant). In \cite[Section 3.2]{FO25}, it was further explained how a L\'{e}vy--Chentsov type convergence (for quasi-particles initially separated by distances of order $T^{1/2}$) can be recast as a fluctuation result for the hard rods model on diffusive time scales. It is likely possible to show an analog of that statement for the Toda lattice, as a corollary of Theorem \ref{qjt0}, but we will not pursue this here.

Certain marginals of the counterpart of Theorem \ref{qjt0} were also established by Olla--Sasada--Suda \cite{SLSBS} for the box-ball model. They showed that the fluctuations of one box-ball soliton converge to a Brownian motion, and those of two solitons of the same length (and of distance at most $T^{1/2}$ from each other) converge to shifts of the same Brownian motion. These can be seen to be consistent with dressed L\'{e}vy--Chentsov field fluctuations. Indeed, as a consequence of Theorem \ref{qjt0}, we deduce the corresponding results for the Toda lattice as Proposition \ref{prop:tracer_fluct_intro} and Corollary \ref{zqi2} below (see also Corollary \ref{zqi200}).

\subsection{Proof discussion} 

\label{ProofFluctuationZ}

In this section we comment on the proof of Theorem \ref{qjt0}. Since  a detailed exposition will be given in Section \ref{ProofQ} below, we will be brief here. Below, we let $\alpha \coloneqq \mathbb{E}[q_{i+1}] - \mathbb{E}[q_i]$ denote the average spacing between consecutively indexed particles under thermal equilibrium (which is an explicit function of $\beta$ and $\theta$; see \eqref{eq:alpha_def}). Then, we likely have $q_i (0) = \alpha i + o(i)$. 

Recall the approximate ``center'' $\varphi_t (j)$ of the eigenvector $\mathbf{u}_j (t)$ of $\mathbf{L}(t)$ with eigenvalue $\lambda_j$, which was used to define the quasi-particle $Q_j (t) = q_{\varphi_t (j)} (t)$. For convenience, we take $\varphi_0 : [1,N] \rightarrow [N_1, N_2]$ to be a bijection (one can verify \cite{Agg25a} that this can always be arranged). We denote the quasi-particle fluctuations by
\begin{flalign}
\label{ziqt}
Z_{i}^{\mathcal Q}(t) \coloneqq  T^{-1/2} \cdot \big(Q_{\varphi_0^{-1}(i)}(t) - q_i(0)- t \ve(\Lambda_i) \big), \qquad \text{where} \qquad \Lambda_i = \lambda_{\varphi_0^{-1} (i)}. 
\end{flalign}

\noindent Here, we have indexed the quasi-particles $(Q_j)$ by $j=\varphi_0^{-1} (i)$, instead of $j=i$, since this indexing is consistent with their initial positions; indeed, these initial positions are given by $Q_{\varphi_0^{-1} (i)} = q_i (0) = \alpha i + o(i)$. By \eqref{qjt} (at $j = \varphi_0^{-1} (i)$), we have with high probability that
\begin{flalign}
\label{ziqt0}
|Z_i^{\mathcal{Q}}(t)| \le (\log T)^C.
\end{flalign}

Now, the proof of Theorem \ref{qjt0} is based on an analysis of the asymptotic scattering relation \eqref{qktqk0}, but it will first be convenient to smoothen the indicator functions there. So, fix a parameter $T^{1/2} \ll \mathfrak{M} \ll T$, and let $\chi(x)$ denote a smooth function approximating $\mathbbm{1}_{x>0}$ on scale $\mathfrak{M}$, meaning that $\chi(x) = \mathbbm{1}_{x > 0}$ for $|x| \ge \mathfrak{M}$. It was shown in \cite{Agg25} (see Lemma \ref{lem:asymptotic-scattering}) that we can replace the indicator functions in \eqref{qktqk0} with $\chi$, while incurring an error of at most $\mathfrak{M}^{1/2+o(1)} \ll T^{1/2}$, which is negligible on the fluctuation scale. This, with a change of indices in \eqref{qktqk0} from $(j,k)$ to $(\varphi_0^{-1} (i), \varphi_0^{-1} (k))$, gives
\begin{flalign}
\label{ziqt00} 
\begin{aligned} 
Q_{\varphi_0^{-1}(k)} (t) & \approx Q_{\varphi_0^{-1} (k)} (0) + \Lambda_k  t \\ 
& \quad + 2 \displaystyle\sum_{i \ne k} \log |\Lambda_k - \Lambda_i| \cdot \big( \chi (Q_{\varphi_0^{-1} (k)} (0) - Q_{\varphi_0^{-1} (i)} (0)) - \chi (Q_{\varphi_0^{-1} (k)} (t) - Q_{\varphi_0^{-1} (i)} (t)) \big).
\end{aligned} 
\end{flalign}

\noindent Since $\chi$ is smooth and $Q_{\varphi_0^{-1} (j)} (t) = q_j (0) + t \ve (\Lambda_j) + T^{1/2} Z_j^{\mathcal{Q}} (t)$ by \eqref{ziqt}, the approximation \eqref{ziqt00} can be Taylor expanded to first order around $Q_{\varphi_0^{-1} (j)} (t) \approx q_j (0) + t \ve (\Lambda_j)$, using \eqref{ziqt0}. After some computation (using $q_j (0) = \alpha j + o(j)$), one obtains from this expansion that (see Lemma \ref{eqn:Zeq_1})
\begin{flalign}
\label{qktqk01}
\begin{aligned}
Z_{k}^{\mathcal Q} (& t) + 2 \sum_{i \ne k} ( Z_{k}^{\mathcal Q} (t)- Z_{i}^{\mathcal Q} (t)) \cdot \log |\Lambda_{k}-\Lambda_i| \cdot \chi' \big(\alpha( k- i) +  t  (\ve(\Lambda_{k})  - \ve(\Lambda_{i}) ) \big) \approx  \Xi(\Lambda_{k}, k, t),
\end{aligned}
\end{flalign}

\noindent where throughout we view $t = T\tau$ as of order $T$, and we have denoted 
\begin{multline}
\label{xilambda}
    \begin{aligned} 
 \Xi & (\Lambda, k, t) \\
 & =  \displaystyle\frac{2}{T^{1/2}} \sum_{i = N_1}^{N_2} \log |\Lambda- \Lambda_i| \cdot
\big(  
    \chi(q_{k}(0) - q_i(0)) -
    \chi \big(q_{k}(0) - q_i(0) + t  \left(\ve(\Lambda) - \ve(\Lambda_i) \right) \big)
   \big) - \mathbb{E} [\cdots ],
    \end{aligned} 
 \end{multline}

 \noindent where the quantity within the expectation $\mathbb{E}[\cdots]$ is the first term on the right side of \eqref{xilambda}. 
 
 Observe that \eqref{qktqk01} is now linear in the $(Z_j^{\mathcal{Q}})$ variables (though with coefficients that are correlated both with the $(Z_j^{\mathcal{Q}})$ and each other). Thus, we must first pinpoint the limit of the right side $\Xi$ of \eqref{qktqk01}, and then solve \eqref{qktqk01} for the $(Z_j^{\mathcal{Q}})$. This amounts to showing the following two statements.

\begin{enumerate} 

\item \emph{Scaling limit for $\Xi$} (Theorem \ref{thm:couplethm_outline}): $\Xi$ converges, as a process in $(\Lambda,k,t)$, to a Gaussian limit $\mathfrak{X}$ (given by an ordinary L\'{e}vy--Chentsov field). 
\item \emph{Approximating the $(Z_j^{\mathcal{Q}})$ through $\Xi$} (Theorem \ref{prop:dressing_Z_approx_outline}): The $(Z_j^{\mathcal{Q}})$ approximate an explicit linear operator (given by the dressing operator) applied to $\Xi$.

\end{enumerate} 

Together, the above two statements indicate that the $(Z_j^{\mathcal{Q}})$ converge to the composition of $\mathfrak{X}$ with a certain operator, which will be the dressed L\'{e}vy--Chentsov field of Theorem \ref{qjt0}. Their proofs occupy most of this work; we remark on them here briefly, leaving further discussion to Section \ref{ProofQ}. 

For the former, observe that $\Xi$ is not a linear statistic of $\mathbf{L}(0)$, as it is a function of both its eigenvalues $(\Lambda_i)$ and off-diagonal entries, through the $(q_i)$. One aspect of accessing its scaling limit will involve realizing an independence structure within the sum in \eqref{qktqk01}; this uses an approximate independence between the pairs $(\Lambda_i, q_{i+1}(0)-q_i(0))$ and $(\Lambda_j, q_{j+1}(0)-q_j(0))$ if $|i-j|$ is large, which was also used in \cite{Agg25a,Agg25} and can be traced to the localization of eigenvectors of $\mathbf{L}(0)$. See the end of Section \ref{ProofXLambda} for further commentary, and also for other aspects in identifying the limit of $\Xi$. 

For the second, observe in \eqref{qktqk01} that the integer $k$ is involved in the coefficients there both through the spatial coordinate $\alpha k$ and the spectral variable $\Lambda_k$. Since it will happen that the latter is discontinuous in the former, we will generalize \eqref{qktqk01} by replacing every appearance of $\Lambda_k$ there by an arbitrary parameter $\Lambda$. We will then guess a solution $Z$ to this new equation, dependent on two independent parameters (a spectral and a spatial one), obtained by applying an explicit linear operator (the dressing operator) to $\Xi$ (see \eqref{eqn:Z_outline}). Verifying that $Z$ is indeed a solution (Proposition \ref{prop:sumZlambdaQ}) requires a concentration bound (Proposition \ref{prop:sumZconc}) to account for the correlation between $Z$ and the coefficients in this equation. See Section \ref{ProofZ00} for further commentary. \\ 

Before proceeding, let us make two remarks. First, it is possible that, on a high level, the above framework may apply to other integrable systems that have microscopically defined quasi-particles $(Q_j)$ satisfying a form of the asymptotic scattering relation \eqref{qktqk0}.\footnote{These are known for the box-ball system (due to works of Takahashi--Satsuma \cite{SCA}, Ferrari--Nguyen--Rolla--Wang \cite{SDBS}, and Croydon--Sasada \cite{CS21}) and believed for various other integrable models.} In particular, we expect it should be correct that, (i) such a relation can be smoothened with negligible error into a form of \eqref{ziqt00}, and the latter can be Taylor expanded into a linear equation analogous to \eqref{qktqk01}; (ii) the right side of this equation converges to a L\'{e}vy--Chentsov field; and (iii) the solution to the linear equation from (i) converges to a dressed L\'{e}vy--Chentsov field. If so, proving these statements might still require model-specific arguments (such as in realizing the independence structure required to show (ii)).

Second, works of Saitoh \cite{TCLE} and Matetski--Quastel--Remenik \cite{PGL} predict under a certain limit that the Toda lattice converges to the Korteweg--de Vries (KdV) equation; in this limit, thermal equilibrium (when suitably scaled) can be seen to map to white noise. Such a convergence has not yet been justified, partly since even the well-definedness of the KdV equation under white noise initial data is quite subtle, proven by Kappeler--Topalov \cite{GW} and Quastel--Valk\'{o} \cite{PN} on the torus and by Killip--Murphy--Visan \cite{INL} on the line (the latter \cite{INL} also used ideas from Anderson localization). It would be interesting to verify this convergence and see whether the above framework for the Toda lattice could be adapted to access the long-time asymptotics for the KdV equation.

\subsection*{Acknowledgements}

The authors thank Alexei Borodin, Percy Deift, Benjamin Doyon, Tamara Grava, Ken McLaughlin, Stefano Olla, and Herbert Spohn for valuable conversations. The authors are also grateful to Makiko Sasada for suggesting to add Appendix \ref{ParticlesInfinite}. The work of Amol Aggarwal was partially supported by a Clay Research Fellowship and a Packard Fellowship for Science and Engineering. The work of Matthew Nicoletti was partially supported by an NSF postdoctoral fellowship, grant No. DMS 2402237. This work was started at the ``Program on Classical, Quantum, and Probabilistic Integrable Systems,'' held in the spring of 2025 at the Center of Mathematical Sciences and Applications in Harvard University.

\section{Results} 

\label{Results} 

\subsection{Toda lattice}
\label{subsec:model_and_objs}

In this section we recall the Toda lattice on an interval, and its thermal equilibrium. Throughout, we fix integers $N_1 \le N_2$ and set $N = N_2 - N_1 + 1$ (which will prescribe the interval's endpoints and length, respectively). We may have $N_1 = \infty$ or $N_2 = -\infty$; in either case, $N = \infty$. We also let $\llbracket a, b \rrbracket = [a,b] \cap \mathbb{Z}$ for any real numbers $a \le b$. 

The state space of the \emph{Toda lattice} on the interval $\llbracket N_1, N_2 \rrbracket$ is given by a pair $( \mathbf{p} (t); \mathbf{q}(t) ) \in \mathbb{R}^N \times \mathbb{R}^N$, where $\mathbf{p}(t) = ( p_{N_1} (t), p_{N_1+1} (t), \ldots , p_{N_2} (t) )$ and $\mathbf{q}(t) = ( q_{N_1} (t), q_{N_1+1}(t), \ldots , q_{N_2} (t) )$ are $N$-tuples, both indexed by a real number $t \ge 0$ called the time. Given any \emph{initial data} $( \mathbf{p}(0); \mathbf{q}(0) ) \in \mathbb{R}^N \times \mathbb{R}^N$, the joint evolution of $( \mathbf{p}(t); \mathbf{q}(t) )$ for $t \ge 0$ is prescribed by the system of ordinary differential equations 
\begin{flalign}
	\label{eqn:eom1}
	\partial_t q_j (t) = p_j (t), \quad \text{and} \quad \partial_t p_j (t) = e^{q_{j-1} (t) - q_j (t)} - e^{q_j (t) - q_{j+1} (t)},
\end{flalign}

\noindent for all $(j, t) \in \llbracket N_1, N_2 \rrbracket \times \mathbb{R}_{\ge 0}$; here, we set $q_{N_1-1}(t) = -\infty$ and $q_{N_2+1}(t) = \infty$ for all $t \ge 0$. One might interpret this as the dynamics for $N$ points (indexed by $\llbracket N_1, N_2 \rrbracket$) moving on the real line, whose locations and momenta at time $t \ge 0$ are given by the $(q_i(t))$ and $(p_i(t))$, respectively. 	

The system of differential equations \eqref{eqn:eom1} is equivalent to the Hamiltonian dynamics generated by the Hamiltonian $\mathfrak{H} : \mathbb{R}^N \times \mathbb{R}^N \rightarrow \mathbb{R}$ that is defined, for any $\mathbf{p} = (p_{N_1}, p_{N_1+1}, \ldots , p_{N_2}) \in \mathbb{R}^N$ and $\mathbf{q} = (q_{N_1}, q_{N_1+1}, \ldots , q_{N_2}) \in \mathbb{R}^N$, by setting
\begin{flalign*}
	\mathfrak{H} (\mathbf{p}; \mathbf{q}) = \displaystyle\frac{1}{2} \displaystyle\sum_{j=N_1}^{N_1} p_j^2 + \displaystyle\sum_{j=N_1}^{N_2} e^{q_j - q_{j+1}}.
\end{flalign*} 

\noindent When $N$ is finite, the existence and uniqueness of solutions to \eqref{eqn:eom1} for all time $t \ge 0$, under arbitrary initial data $(\mathbf{p}; \mathbf{q}) \in \mathbb{R}^N \times \mathbb{R}^N$, is thus a consequence of the Picard--Lindel\"{o}f theorem (see, for example, the proof of \cite[Theorem 12.6]{OINL}). 

It will often be useful to reparameterize the variables of the Toda lattice, following \cite{Fla74}. To that end, for any $(i, t) \in \llbracket N_1, N_2 \rrbracket \times \mathbb{R}_{\ge 0}$, define 
\begin{flalign}
	\label{abr} 
	r_i (t) = q_{i+1} (t) - q_i (t); \qquad a_i (t) = e^{-r_i(t)/2}; \qquad b_i (t) = p_i (t).
\end{flalign}

\noindent Denoting $\mathbf{a}(t) = ( a_{N_1} (t), a_{N_1+1} (t), \ldots , a_{N_2} (t) ) \in \mathbb{R}_{\ge 0}^N$ and $\mathbf{b}(t) = ( b_{N_1} (t), b_{N_1+1} (t), \ldots , b_{N_2} (t) ) \in \mathbb{R}^N$, the $( \mathbf{a} (t); \mathbf{b} (t) )$ are sometimes called \emph{Flaschka variables}; they satisfy $r_{N_2}(t) = q_{N_2+1}(t) - q_{N_2} (t) = \infty$ and $a_{N_2} (t) = e^{-r_{N_2}(t) / 2} = 0$. Then, \eqref{eqn:eom1} is equivalent to the system  
\begin{flalign}
	\label{eqn:Flaschka_ev}
	\partial_t a_j (t) = \displaystyle\frac{a_j (t)}{2} \cdot \big(b_j (t) - b_{j+1} (t) \big), \qquad \text{and} \qquad \partial_t b_j (t) = a_{j-1} (t)^2 - a_j (t)^2,
\end{flalign} 

\noindent for each $(j, t) \in \llbracket N_1, N_2 \rrbracket \times \mathbb{R}_{\ge 0}$. Since the dynamics \eqref{eqn:eom1} are well-posed when the interval $\llbracket N_1, N_2 \rrbracket$ is finite, so are the dynamics \eqref{eqn:Flaschka_ev}. It was shown as \cite[Proposition 4.7]{Agg25a} that the dynamics \eqref{eqn:Flaschka_ev} are also well-posed when $\llbracket N_1, N_2 \rrbracket = \mathbb{Z}$ is infinite (by taking a limit of these dynamics on $\llbracket N_1, N_2 \rrbracket$, as $(N_1, N_2)$ tends to $(-\infty, \infty)$), under the assumption that 
\begin{flalign}
	\label{ajbjlimit}
	\displaystyle\limsup_{j \in \mathbb{Z}} \displaystyle\frac{|a_j(0)| + |b_j(0)|}{(|j|+1)^{\mathfrak{p}}} < \infty, \qquad \text{for some $\mathfrak{p} \in (0, 1)$}. 
\end{flalign}

\noindent Thus, if \eqref{ajbjlimit} holds, then we may consider the Toda lattice \eqref{eqn:Flaschka_ev} on the line, when $\llbracket N_1, N_2 \rrbracket = \mathbb{Z}$, as will be our setting when stating our current fluctuation results in Section \ref{subsec:main_res} below. 

It will at times be necessary to define the original Toda state space variables $( \mathbf{p}(t); \mathbf{q}(t) )$ from the Flaschka variables $( \mathbf{a}(t); \mathbf{b}(t) )$; it suffices to do this at $t = 0$, as $( \mathbf{p} (t); \mathbf{q}(t) )$ is determined from $( \mathbf{p}(0); \mathbf{q}(0) )$ and $(\mathbf{a}(s); \mathbf{b}(s) )_{s \ge 0}$, by \eqref{eqn:eom1} and \eqref{abr}. We explain how to do this when $0 \in \llbracket N_1, N_2 \rrbracket$ (as will be the case for us). By \eqref{abr}, the Flaschka variables $\mathbf{a}(0)$ specify the differences between consecutive entries in $\mathbf{q}(0)$, so the former determines the latter up to an overall shift. We fix this shift by setting $q_0 (0) = 0$. Then, $( \mathbf{p} (0); \mathbf{q} (0) )$ is called the \emph{Toda state space initial data} associated with $(\mathbf{a} (0); \mathbf{b}(0) )$. The evolution $( \mathbf{p}(t); \mathbf{q}(t) )$ of this initial data under \eqref{eqn:eom1} is called the \emph{Toda state space dynamics} associated with $( \mathbf{a}(t); \mathbf{b}(t) )$; observe that we may have $q_0 (t) \ne 0$ if $t \ne 0$.	

Under this notation, we can define the Toda lattice Lax matrix (originally introduced in \cite{Fla74,Man74}) as follows.

\begin{definition}[Lax matrix]
	
	\label{lt}
	
	Fix $t \in \mathbb{R}_{\ge 0}$. Define the \emph{Lax matrix} $\L(t) = [L_{i j}(t)]_{i, j \in \llbracket N_1, N_2 \rrbracket}$ by 
\begin{align}\label{eqn:Flaschka_vars_inf}
    L_{i i} (t) &= b_i, \quad i\in  \text{for $\llbracket N_1, N_2 \rrbracket$}; \qquad   L_{i,i+1} (t) = L_{i+1,i} (t) = a_i, \quad \text{for $i \in \llbracket N_1, N_2 -1 \rrbracket$}.
\end{align}

\noindent Also set $L_{ij} = 0$ for any $i,j \in \mathbb{Z}$ not of the above form. 

\end{definition}

The use of the Lax matrix is in the following fact, originally due to \cite{Fla74} (see also \cite[Section 2]{Mos75}), indicating that its eigenvalues are conserved. 

\begin{lem}[{\cite{Fla74,Mos75}}]
	
	\label{lem:lax_eig}
	
    The eigenvalues of $\L(t)$ do not depend on $t$, if $N<\infty$.
    
\end{lem}

We will study the Toda lattice under random initial data, in which the Flaschka variables are sampled under a probability measure called thermal equilibrium. In it, the $(a_i, b_i)$ are mutually independent, the $(a_i)$ are chi random variables, and the $(b_i)$ are Gaussian random variables.

\begin{definition}[Thermal equilibrium]\label{def:inf_thermal_eq}

Fix real numbers $\theta, \beta > 0$. \emph{Thermal equilibrium} $\mu = \mu_{\beta,\theta} = \mu_{\beta,\theta;N_1, N_2}$ is the probability measure on $(\mathbf{a}; \mathbf{b}) \in \mathbb{R}_{>0}^{N-1} \times \mathbb{R}^N$, where $\mathbf{a}= (a_i )_{i \in \llbracket N_1, N_2-1 \rrbracket}$ and $\mathbf{b} = (b_i)_{i \in \llbracket N_1, N_2 \rrbracket}$, defined by 
    \begin{equation}\label{eqn:equil}
\mu(d\mathbf{a}; d\mathbf{b})
= \prod_{j\in \llbracket N_1, N_2-1 \rrbracket} \frac{2\beta^{\theta}}{\Gamma(\theta)} a_j^{2\theta - 1} e^{-\beta a_j^2} \, da_j
 \cdot \prod_{j\in \llbracket N_1, N_2 \rrbracket}   (2\pi \beta^{-1})^{-1/2} e^{-\beta b_j^2 / 2} \, db_j.
\end{equation}

\noindent If $\llbracket N_1, N_2 \rrbracket = \mathbb{Z}$, we abbreviate $\mu_{\beta,\theta;N_1, N_2} = \mu_{\beta,\theta;\infty}$. Observe that \eqref{ajbjlimit} holds almost surely under thermal equilibrium $\mu_{\beta,\theta;\infty}$, so the Toda lattice on $\mathbb{Z}$ is well-defined under sampling its initial data $(\mathbf{a}(0); \mathbf{b}(0))$ under this measure.

\end{definition}

Thermal equilibrium $\mu_{\beta,\theta;\infty}$ is invariant \cite[Proposition 2.5]{Agg25a} for the Toda lattice on $\mathbb{Z}$ (and is perhaps its most natural invariant measure). We will be interested in analyzing, under thermal equilibrium, the fluctuations of the Toda lattice's local charges and associated integrated currents. 

To explain them, observe that Lemma \ref{lem:lax_eig} provides a large family of conserved quantities for the Toda lattice, given by the eigenvalues of the Lax matrix. However, these are ``non-local,'' in the sense that they depend on all of the Flaschka variables $(\mathbf{a}(t); \mathbf{b}(t))$, as opposed to only the $(a_i(t), b_i(t))$ for $i$ in some (uniformly) bounded interval. To remedy this, observe that Lemma \ref{lem:lax_eig} also implies that $\Tr \mathbf{L}(t)^m$ is conserved for any integer $m \ge 0$. Since $\mathbf{L}(t)$ is tridiagonal, the diagonal entries of $\mathbf{L}(t)^m$ are local quantities, whose total is preserved; they are called local charges and have associated currents. In what follows, for any matrix $\mathbf{M}$, we let $[\mathbf{M}]_{ij}$ denote the $(i, j)$-entry of $\mathbf{M}$.

\begin{definition}[Local charges]\label{kti00}
	
 Fix an integer $m \geq 0$, index $i \in \llbracket N_1, N_2 \rrbracket$, and real number $t \geq 0$. Define the $m$-th \emph{local charge} $\mathfrak{k}_i^{[m]}(t)$ of $\L(t)$ at $i$ by setting
\begin{flalign} 
\label{mki} 
\mathfrak{k}_i^{[m]}(t) = [\L(t)^m]_{ii}.
\end{flalign}

\noindent Following \cite[Equation (2.6)]{Spo20}, also define a modification of this $m$-th local charge $\mathfrak{K}_i^{[m]} (t)$ by  
\begin{flalign*}
\mathfrak{K}_i^{[0]} (t) = r_i (t) = q_{i+1} (t) - q_i (t), \qquad \text{and} \qquad 
    \mathfrak{K}_i^{[m]} (t) = \mathfrak{k}_i^{[m]} (t), \quad \text{if $m>0$}. 
\end{flalign*}

\end{definition} 

For example, the first local charge $\mathfrak{k}_i^{[1]} (t) = \mathfrak{K}_i^{[1]} (t) = b_i (t) = p_i(t)$ is the momentum of the $i$-th particle in the Toda lattice at time $t$. The zero-th local charge $\mathfrak{k}_i^{[0]} (t)$ is always equal to $1$, so one may view it as counting the $i$-th particle in the Toda lattice. The zero-th modified local charge $\mathfrak{K}_i^{[0]} (t) = q_{i+1}(t) - q_i (t)$ denotes the spacing between the $i$-th and $(i+1)$-th particles at time $t$ (so its inverse is a local prescription for the particle density around $q_i (t)$).

\begin{definition}[Spatial integrated currents]\label{def:integrated_curr_spatial}

Next, for any integer $m \geq 0$, index $i \in \mathbb{Z}$, and real number $t \geq 0$, define the $m$-th \emph{local current} $\mathfrak{j}_i^{[m]}(t)$ by
\begin{flalign} 
\label{jti}
\mathfrak{j}_i^{[m]}(t) =
a_{i-1}(t) \cdot [\L(t)^m]_{i,i-1} .
\end{flalign}

\noindent Further fixing real numbers $q,q' \in \mathbb{R}$, define the $m$-th spatial integrated current $J_t^{[m]}(q, q')$ (across the spatial interval $[q,q']$) as follows. If $\llbracket N_1, N_2 \rrbracket$ is finite, then set 
\begin{flalign}
	\label{currentn} 
	J_t^{[m]}(q, q') =
	\sum_{j : q_j (0) < q}  \mathfrak{k}_j^{[m]}(0) -\sum_{j: q_j(t) < q'}  \mathfrak{k}_j^{[m]}(t).
\end{flalign}

\noindent If $\llbracket N_1, N_2 \rrbracket = \mathbb{Z}$, then both sums above might diverge. So, if $\llbracket N_1, N_2 \rrbracket = \mathbb{Z}$ and there exists an integer $k \in \mathbb{Z}$ such that $\max_{j \le k} (q_{j}(t), q_{j}(0)) < \min(q,q')$, then we set\footnote{If we write $J_t^{[m]} (q,q')$ in an equation when $\llbracket N_1, N_2 \rrbracket = \mathbb{Z}$, then implicit in the equation is the claim that such an integer $k$ exists.} 
\begin{equation}\label{eqn:inf_integrated_curr}
 J_t^{[m]}(q, q') =
 \sum_{j \geq k: q_j(0) < q}  \mathfrak{k}_i^{[m]}(0) -\sum_{j \geq k: q_j(t) < q'}  \mathfrak{k}_i^{[m]}(t) + \int_0^t \mathfrak{j}_{k}^{[m]}(s) ds.
\end{equation}
\end{definition}

\begin{remark} 

\label{km00} 

Let us briefly comment on the equivalence between \eqref{currentn} and \eqref{eqn:inf_integrated_curr}. Using \eqref{eqn:Flaschka_ev}, it can be verified (see \cite[Equation (2.24)]{Spo23} for the derivation) that 
\begin{equation}\label{eqn:current_rel}
\partial_t \mathfrak{k}_k^{[m]}(t) = \mathfrak{j}_{k}^{[m]}(t) -\mathfrak{j}_{k+1}^{[m]}(t) ,\qquad k \in \llbracket N_1, N_2 \rrbracket.
\end{equation}

\noindent In this way, the current $\mathfrak{j}_i^{[m]}$ signifies the ``rate of transfer'' of the total local charge $\mathfrak{k}^{[m]}$ from the right of site $k$ to the left of it (which is why it is referred to as a ``current''). 

It quickly follows from \eqref{eqn:current_rel} that the definitions \eqref{currentn} and \eqref{eqn:inf_integrated_curr} for $J_t^{[m]} (q,q')$ coincide if $\llbracket N_1, N_2 \rrbracket$ is finite and $k \in \llbracket N_1, N_2 \rrbracket$ satisfies $\max_{j \in \llbracket N_1, k \rrbracket} (q_{j}(t), q_{j}(0)) < \min(q,q')$. It also quickly follows from \eqref{eqn:current_rel} that, when $\llbracket N_1, N_2 \rrbracket = \mathbb{Z}$, the definition \eqref{eqn:inf_integrated_curr} is independent of the choice of $k \in \mathbb{Z}$ satisfying $\max_{j \le k} (q_{j}(t), q_{j}(0)) < \min(q,q')$.

\end{remark} 

 For example, the zero-th current $J_t^{[0]} (q,q')$ counts the number of particles in the Toda lattice crossing from the left of $q$ at time $0$ to the right of $q'$ at time $t$. The first current $J_t^{[1]} (q,q')$ may be viewed as the net momentum passing across the spatial interval $[q,q']$ from time $0$ to $t$.

\subsection{Dressing operator and effective velocity}

\label{Operatorv}

In this section we define the dressing operator and effective velocity, which are objects that appear in the limiting current fluctuations of the Toda lattice. We begin by defining a probability density $\varrho : \mathbb{R} \rightarrow \mathbb{R}$; see Remark \ref{lrho} below for its interpretation. The functions $\mathfrak{F}$ and $\varrho_{\beta}$ below originally appeared in the study of high-temperature beta ensembles in \cite[Equation (16)]{IEC} and \cite[Theorem 1.1(ii)]{SMRME}; the function $\varrho$ originally appeared in \cite[Equation (3.5)]{Spo20GGE}.

\begin{definition}[Stretch and density of states]\label{def:alpha_Laxdos}

Fix real numbers $\beta, \theta > 0$. Define the real number, also called the \emph{stretch},
\begin{equation}\label{eq:alpha_def}
    \alpha = \log \beta - \frac{\Gamma'(\theta)}{\Gamma(\theta)}, \qquad \text{and assume that } \alpha \neq 0.
\end{equation}

\noindent For any $x \in \mathbb{R}$, set
\[
\mathfrak{F}(\theta; x)
    = \left( \frac{\theta}{\Gamma(\theta)} \right)^{1/2}
      \int_0^{\infty} y^{\theta - 1} e^{i x y - y^2 / 2} \, dy.
\]

\noindent Define the function $\varrho_\beta : \mathbb{R} \to \mathbb{R}$ and the \emph{density of states} $\varrho : \mathbb{R} \to \mathbb{R}$ by, for any $x \in \mathbb{R}$, setting
\[
\varrho_\beta(x) = \varrho_{\beta; \theta}(x)
    = \left( \frac{\beta}{2\pi} \right)^{1/2}
      \, \bigl| \mathfrak{F}_\theta(\beta^{1/2} x) \bigr|^{-2}
      \, e^{-\beta x^2 / 2};
    \qquad
    \varrho(x) = \partial_\theta \bigl( \theta \, \varrho_{\beta; \theta}(x) \bigr).
\]
\end{definition}

    It is quickly verified from Definition \ref{def:inf_thermal_eq} that, under the thermal equilibrium $\mu_{\beta,\theta;\infty}$, we have 
    \begin{flalign} 
    \label{alphaexpectation} 
    \alpha = \mathbb{E}[q_{i+1} - q_i], 
    \end{flalign} 
    
    \noindent which is why $\alpha$ is called the ``stretch''; see \cite[Lemma 3.11]{Agg25a} for the derivation. Thus, the condition $\alpha \neq 0$ is natural, as it is equivalent to the finiteness of the system's particle density.

\begin{remark} 

\label{lrho} 

 The function $\varrho$ is called the density of states since it prescribes  the limiting empirical eigenvalue density for the (finite $N$) Lax matrix of the Toda lattice under thermal equilibrium, in the following sense. Sample an $N \times N$ Lax matrix $\mathbf{L}_N$ (as in Definition \ref{lt}) under thermal equilibrium $\mu_{\beta,\theta;1,N}$ (as in Definition \ref{def:inf_thermal_eq}), and let $(\lambda_1, \lambda_2, \ldots , \lambda_N)$ denote its eigenvalues. Then, for any bounded continuous function $f: \mathbb{R} \rightarrow \mathbb{R}$, we have \cite[Lemma 4.3]{Maz22} that  
\begin{equation*}
    \lim_{N \rightarrow \infty} \mathbb{E}\left[ \displaystyle\frac{1}{N} \sum_{i=1}^N f(\lambda_i) \right] = \int_{-\infty}^{\infty}f(\lambda) \varrho(\lambda) d\lambda.
\end{equation*}

\end{remark} 

    Our limiting functions will reside in (and operators will act on) the following Hilbert space.

\begin{definition}[Hilbert space $\mathcal{H}$ and operator $\mathbf{T}$]\label{def:Hdef}

Let $\mathcal{H} = L^2(\mathbb{R}, \varrho(x) dx) $ denote the Hilbert space associated with the inner product  
\begin{equation}\label{def:inner_prod}
    \langle f, g \rangle_{\varrho} \coloneqq \int_{-\infty}^{\infty}f(x) \overline{g(x)} \varrho(x) dx.
\end{equation} 
We have by \cite[Lemma 3.1]{Agg25} (see also Lemma \ref{eqn:varrho_tail_bds} below) that $\varsigma_n \in \mathcal{H}$ for any integer $n \ge 0$, where the function $\varsigma_n: \mathbb{R} \rightarrow \mathbb{R}$ is defined by 
\begin{flalign}
    \label{nx} 
	\varsigma_n (x) \coloneqq x^n, \qquad \text{for any $x \in \mathbb{R}$}. 
\end{flalign} 

\noindent For any function $h \in \mathcal{H}$, denote the associated multiplication operator by $\boldsymbol{h}$, defined by setting $\boldsymbol{h} f = hf$ for any $f : \mathbb{R} \rightarrow \mathbb{R}$; if $h \in \mathcal{H}$ is constant (that is, of the form $h = a \varsigma_0$ for some $a \in \mathbb{C}$), we identify $h = \boldsymbol{h}$. Further define the integral operator $\T$ on functions of a real variable by
\begin{equation}
\label{operatort}
    \T f(x) = 2 \int_{-\infty}^{\infty}\log|x-y| f(y) dy,
\end{equation}
whenever $f$ is such that the integral is finite for almost all $x \in \mathbb{R}$.

\end{definition}

By \cite[Lemma 1.4]{Agg25}, the operator $\T   \boldsymbol{\varrho}_{\beta}= \T  \circ\boldsymbol{\varrho}_{\beta} : \mathcal{H} \rightarrow \mathcal{H}$ is bounded. By (this and) \cite[Lemma 1.5]{Agg25}, the operator $\mathrm{Id} -  \theta \T  \boldsymbol{\varrho}_{\beta} : \mathcal{H} \rightarrow \mathcal{H}$ is bounded and invertible. Its inverse is called the dressing operator, which is used to define the effective velocity, as in \cite[Equation (6.20)]{Spo23}.

\begin{definition}[Dressing operator and effective velocity]
	
	\label{def:dress_intro}
    \label{def:v_eff_intro}
    
    The bounded operator $(1- \theta \T  \varrhobetabf )^{-1}:\mathcal{H} \rightarrow \mathcal{H}$ is called the \emph{dressing operator}. For any $f \in \mathcal{H}$, set
    \begin{equation}
    \label{drf}
    f^{\dr} \coloneqq (1- \theta \T \varrhobetabf )^{-1}f \in \mathcal{H}.
    \end{equation}
	
	\noindent Further define the \emph{effective velocity} $\ve: \mathbb{R} \rightarrow \mathbb{R}$ by, for any $x \in \mathbb{R}$, setting
    \begin{equation}\label{eqn:v_eff_intro}
        \ve(x) \coloneqq \varsigma_1^{\dr}(x) \cdot \varsigma_0^{\dr}(x)^{-1}.
    \end{equation}

	\noindent This function is well-defined since, by \cite[Lemma 1.7]{Agg25}, there exists a constant $c>0$ such that $ \varsigma_0^{\dr}(x) \cdot \mathrm{sgn} (\alpha) > c$ for all $x \in \mathbb{R}$.
\end{definition}

\subsection{Current fluctuations}
\label{subsec:main_res}

In this section we state our results on the current fluctuations of the Toda lattice. Throughout, we fix real numbers $\beta, \theta > 0$ and recall the notation from Section \ref{Operatorv}. We begin by defining the white noise that will be used to describe the scaling limits of the quantities from Definitions \ref{kti00} and \ref{def:integrated_curr_spatial}. In what follows, we set
\begin{flalign}
    \label{functionslambdar}
	L^2 (dr \otimes \varrho d \lambda) = \bigg\{  f : \mathbb{R}^2 \rightarrow \mathbb{R} : \displaystyle\int_{-\infty}^{\infty} \displaystyle\int_{-\infty}^{\infty} f(r,\lambda)^2 \varrho (\lambda) dr d \lambda < \infty \bigg\}.
\end{flalign} 

\begin{definition}[White noise]\label{def:white_noise_dress}
    Let $\mathcal{W}$ denote a white noise with respect to $L^2(d r\otimes \varrho d \lambda)$ meaning that, for any functions $f, g  \in L^2(d r\otimes \varrho d \lambda)$, the real-valued random variables $\mathcal{W}(f)$ and $\mathcal{W}(g)$ are centered and jointly Gaussian with covariance $\int_{\mathbb{R}} \int_{\mathbb{R}} f(r, \lambda) g(r, \lambda) \varrho(\lambda) d\lambda dr$.

\end{definition}

 We may now state the precise and general version of Theorem \ref{limitj}, which identifies the joint scaling limit of the Toda lattice currents $J_t^{[m]}(q,q')$ (from Definition \ref{def:integrated_curr_spatial}), under thermal equilibrium $\mu_{\beta,\theta;\infty}$ (from Definition \ref{def:inf_thermal_eq}). Its proof is in Section \ref{subsec:cf_proof} below.

\begin{thm}[Current fluctuations]
\label{thm:current_fluct_intro}

    For any $\beta > 0$, there exists a constant $\theta_0 (\beta)>0$ such that the following holds for any fixed $\theta \in (0, \theta_0)$. Let $(\mathbf{a}(t);\mathbf{b}(t))$ denote the Toda lattice on $\mathbb{Z}$ with initial data $(\mathbf{a}(0);\mathbf{b}(0))$, sampled from the thermal equilibrium $\mu_{\beta,\theta;\infty}$. Further fix an integer $k \ge 1$ and, for each $i \in \llbracket 1, k \rrbracket$, fix an integer $m_i$ and real numbers $\mathfrak{q}_i, \mathfrak{q}_i' \in \mathbb{R}$ and $\tau_i \in \mathbb{R}_{\ge 0}$. As $T$ tends to $\infty$, the $k$-tuple
    \begin{flalign*}
       \Big[ T^{-1/2} \cdot \big(J_{T \tau_i}^{[m_i]}(T \mathfrak{q}_i, T \mathfrak{q}_i') - \mathbb{E}[J_{T \tau_i}^{[m_i]}(T \mathfrak{q}_i, T \mathfrak{q}_i')] \big) \Big]_{i \in \llbracket 1, k \rrbracket}, 
    \end{flalign*} 
    
    \noindent converges in law to the $k$-tuple
  
    \begin{flalign}
        \label{w1} 
 \bigg[     \mathcal{W} \Big(
\varsigma_{m_i}^{\dr}(\lambda) \cdot \big(
\mathbbm{1}\{\mathfrak{q}_i>\alpha r\}
-\mathbbm{1}\{\mathfrak{q}_i'>\alpha r+\tau_i \ve(\lambda)\} \big)
\Big)
\bigg]_{i \in \llbracket 1, k \rrbracket}.
  \end{flalign}

\end{thm}

In this way, the fluctuation scaling limit of the currents for the Toda lattice is given by the Gaussian process obtained by integrating the white noise $\mathcal{W}$ against explicit functions. Let us mention that, while the convergence in Theorem \ref{thm:current_fluct_intro} is stated in the sense of finite-dimensional distributions, it can likely be proven in the sense of uniform convergence on compact subsets. However, we will not pursue this here. 

We next discuss two consequences of (the proof of) Theorem \ref{thm:current_fluct_intro}. First, we have the following corollary, which is the precise formulation of Theorem \ref{q0t00} and is proved in Section \ref{subsec:eulerfluctproof}.

\begin{cor}[Convergence of $q_0(t)$]
\label{cor:q0bm}

  For any $\beta > 0$, there exists a constant $\theta_0 (\beta)>0$ such that the following holds for any fixed $\theta \in (0, \theta_0)$. Let $(\mathbf{a}(t);\mathbf{b}(t))$ denote the Toda lattice on $\mathbb{Z}$ with initial data $(\mathbf{a}(0);\mathbf{b}(0))$, sampled from the thermal equilibrium $\mu_{\beta,\theta;\infty}$. We have the convergence in distribution, in the sense of finite dimensional marginals,
    \begin{flalign}
    	\label{eqn:q0_conv}
    	\begin{aligned} 
    T^{-1/2} q_0( T \tau )\rightarrow \mathcal{B}(\tau) ,
     \end{aligned} 
    \end{flalign}
    where $\mathcal{B}(\tau)$ is a Brownian motion satisfying
    \begin{equation}\label{eqn:q0var}
        \var \mathcal{B}(\tau) =
      \alpha \tau
        \int_{-\infty}^{\infty}
        \bigl|\varsigma_0^{\dr}(\lambda)\bigr|^2\,|\ve(\lambda)|\,\varrho(\lambda)\,d\lambda .
\end{equation}
\end{cor}

To facilitate the statement of the next corollary, we introduce the following normalization of the dressing operator.

 \begin{definition}

\label{operatorf} 

Define the operator $\mathbf{F} : \mathcal{H} \rightarrow \mathcal{H}$ by, for any $h \in \mathcal{H}$, setting
    \begin{flalign*} 
    \mathbf{F} h = h^{\dr} - \langle h, \varsigma_0 \rangle_{\varrho} \cdot \varsigma_0^{\dr}.
    \end{flalign*} 

 \end{definition}

In what follows, we recall the local charges $\mathfrak{K}_i^{[m]}(t)$ from Definition \ref{kti00};\footnote{In Corollary \ref{cor:twopoint_intro} we use the $\mathfrak{K}_j$ local charges since, if we used the $\mathfrak{k}_j$ ones instead, then since $\mathfrak{k}_j^{[0]} = 1$ almost surely, covariances involving the $0$-th charge would be equal to $0$.} and the function $\varsigma_n(x) = x^n$ from \eqref{nx}. Below is the precise and general formulation of Theorem \ref{expectationp}, concerning the \emph{two-point function} $S_{m, n}(j,t) \coloneqq \Cov(\mathfrak{K}_j^{[m]}(t), \mathfrak{K}_0^{[n]}(0))$. It is proved in Section \ref{subsec:tpproof}.

\begin{cor}[Two-point function decay]
\label{cor:twopoint_intro}

  For any $\beta > 0$, there exists a constant $\theta_0 (\beta)>0$ such that the following holds for any fixed $\theta \in (0, \theta_0)$. Let $(\mathbf{a}(t);\mathbf{b}(t))$ denote the Toda lattice on $\mathbb{Z}$ with initial data $(\mathbf{a}(0);\mathbf{b}(0))$, sampled from the thermal equilibrium $\mu_{\beta,\theta;\infty}$. For any integers $m, n \ge 0$; site $j \in \mathbb{Z}$; and real number $t \ge 0$, let
 \begin{equation}\label{eqn:two_pt_def}
	S_{m, n}(j,t) \coloneqq \Cov(\mathfrak{K}_j^{[m]}(t), \mathfrak{K}_0^{[n]}(0)).
\end{equation}

    \noindent Then the following statements hold, for any compactly supported, smooth function $f: \mathbb{R} \rightarrow \mathbb{R}$.
    
    \begin{enumerate} 
    
    \item For any integers $m, n \geq 1$, we have that
 \begin{flalign}
    	\label{eqn:two_point_conv_intro}
    	\begin{aligned} 
     \lim_{T \rightarrow \infty}   \sum_{j \in \mathbb{Z}} f(jT^{-1})  S_{m, n}(j,\tau T) = \int_{-\infty}^{\infty} \mathbf{F} \varsigma_m(\lambda) \cdot \mathbf{F} \varsigma_n(\lambda) \cdot 
     f(\alpha^{-1} \tau \ve(\lambda)) \varrho(\lambda) d \lambda.
     \end{aligned} 
    \end{flalign}
    \item For any integer $n \ge 1$, we have that 
    \begin{flalign}
    	\label{eqn:two_point_conv_intro_n0}
    	\begin{aligned} 
     \lim_{T \rightarrow \infty}   \sum_{j \in \mathbb{Z}} f(jT^{-1})  S_{0, n}(j,\tau T) & = -\alpha \int_{-\infty}^{\infty} \varsigma_0^{\dr}(\lambda) \cdot \mathbf{F} \varsigma_n(\lambda) \cdot 
     f(\alpha^{-1} \tau \ve(\lambda)) \varrho(\lambda) d \lambda \\
     & =  \lim_{T \rightarrow \infty}   \sum_{j \in \mathbb{Z}} f(jT^{-1})  S_{n,0}(j,\tau T).
     \end{aligned} 
    \end{flalign}

    \item We have that
    \begin{flalign}
    	\label{eqn:two_point_conv_intro_both0}
    	\begin{aligned} 
     \lim_{T \rightarrow \infty}  \sum_{j \in \mathbb{Z}} f(jT^{-1})  S_{0, 0}(j,\tau T)= \alpha^2 \int_{-\infty}^{\infty} \varsigma_0^{\dr}(\lambda) \cdot  \varsigma_0^{\dr}(\lambda) \cdot 
     f(\alpha^{-1} \tau \ve(\lambda)) \varrho(\lambda) d \lambda.
     \end{aligned} 
    \end{flalign}
    \end{enumerate} 
\end{cor}

Observe that Corollary \ref{cor:twopoint_intro} can be informally rephrased as follows. Denoting the delta distribution by $\delta(\cdot)$, we have for any real number $x \in \mathbb{R}$ that
    \begin{multline}\label{eqn:two_pt_asymp_intro}
	\displaystyle\lim_{T \rightarrow \infty} T \cdot S_{m, n}( \lfloor xT \rfloor, \tau T) \\
    = 
    \begin{pmatrix}
     \alpha^2  \langle \varsigma_0^{\dr}, \delta(x - \alpha^{-1}\tau \ve(\lambda) ) \cdot  \varsigma_0^{\dr}  \rangle_{\varrho}   &  -\alpha  \langle \mathbf{F} \varsigma_m, \delta(x - \alpha^{-1}\tau \ve(\lambda) ) \cdot  \varsigma_0^{\dr}  \rangle_{\varrho} \\
     -\alpha \langle  \varsigma_0^{\dr}, \delta(x - \alpha^{-1}\tau \ve(\lambda) ) \cdot \mathbf{F} \varsigma_n  \rangle_{\varrho}  & \langle \mathbf{F} \varsigma_m, \delta(x - \alpha^{-1}\tau \ve(\lambda) ) \cdot \mathbf{F} \varsigma_n  \rangle_{\varrho}
    \end{pmatrix},
\end{multline}

\noindent  in the sense of integration against smooth, compactly supported test functions. In the matrix in \eqref{eqn:two_pt_asymp_intro}, the cases $m=0$ and $n=0$ play special roles. Its left column corresponds to the case $m=0$, and the top row corresponds to the case $n =0 $, while the right column and bottom row correspond to the cases when $m \geq 1$ and $n \geq 1$, respectively. This recovers the asymptotic form for the two-point function $S_{m,n}(j,t)$ predicted as \cite[Equation (3.15)]{Spo20}.

\subsection{Quasi-particle fluctuations}

\label{FluctuationsQ}

In this section we first recall the quasi-particles for the Toda lattice. We then state our fluctuation result (Theorem \ref{thm:quasi_fluct_intro}) for them, which is what will eventually enable the proofs of the results from Section \ref{subsec:main_res}. Throughout this section, we recall the notation of Section \ref{subsec:model_and_objs} (including the Lax matrix $\L(t)$ from Definition \ref{lt}), and assume that $N_1$ and $N_2$ are finite (and thus so is $N$). 

The definition of quasi-particles will make use of localization properties of eigenvectors, so we begin by setting notation for this.

\begin{definition}
\label{def:loccent}
    \label{center}
	Let $\mathbf{u} = (u(N_1),\dots, u(N_2)) \in \mathbb{R}^N$ be a unit vector. For any $\zeta \in \mathbb{R}_{\ge 0}$, we say that an index $\varphi \in \llbracket N_1, N_2 \rrbracket$ is a \emph{$\zeta$-localization center} of $\mathbf{u}$ if $|u(\varphi)| \ge \zeta$. 
	
	Next, let $\mathbf{M}$ be an $N \times N$ real symmetric matrix, with rows and columns indexed by $\llbracket N_1, N_2\rrbracket$. If $\lambda$ is an eigenvalue of $\mathbf{M}$, we call $\varphi$ a $\zeta$-localization center for $\lambda$ if $\varphi$ is a $\zeta$-localization center of some unit eigenvector $\mathbf{u}$ of $M$ with eigenvalue $\lambda$. Further let $(\mathbf{u}_1, \mathbf{u}_2, \ldots , \mathbf{u}_N)$ denote an orthonormal eigenbasis for $\mathbf{M}$. We call a bijection $\varphi : \llbracket 1, N \rrbracket \rightarrow \llbracket N_1, N_2 \rrbracket$ a \emph{$\zeta$-localization center bijection} if $\varphi(j)$ is a $\zeta$-localization center for $\mathbf{u}_j$, for each $j= \llbracket 1, N \rrbracket$. By \cite[Lemma 2.7]{Agg25a}, any symmetric $N \times N$ matrix $\mathbf{M}$ admits a $(2N)^{-1}$-localization center bijection.
\end{definition}

Localization center bijections are not always unique. However they, and in fact individual localization centers, are unique up to a small error, if the entries of $\mathbf{M}$ are sufficiently random, such as under thermal equilibrium (see Lemma \ref{lem:loc_cent_dif} below). Due to this approximate uniqueness, our choice to take $\varphi$ to be a bijection is notationally convenient but not essential.

To fix the setup, when we work with Toda lattice on a finite interval $\llbracket N_1, N_2 \rrbracket$, we will typically impose the following assumption.

\begin{ass}\label{ass:NT_assumption}

    As in Definition \ref{def:alpha_Laxdos}, fix real numbers $\beta,\theta>0$, and assume that $\alpha\neq 0$.  
For each real number $t\ge 0$, let $\L(t)=[L_{ij}(t)]_{i, j \in \llbracket N_1, N_2 \rrbracket}$ denote the Lax matrix for the Toda lattice $(\mathbf{a}(t); \mathbf{b}(t))$ on $\llbracket N_1,N_2\rrbracket$ (as in Definition \ref{lt}). Denote the eigenvalues of $\L(t)$ by $\lambda_1 \ge \lambda_2 \ge \cdots \ge \lambda_N$, which do not depend on $t$ by Lemma \ref{lem:lax_eig}.  
At $t=0$, abbreviate $\L=\L(0)$ and $(\mathbf{a};\mathbf{b})=(\mathbf{a}(0); \mathbf{b}(0))$.  
Assume that the initial data $(\mathbf{a}(0);\mathbf{b}(0))$ is sampled under the thermal equilibrium $\mu_{\beta,\theta;N_1,N_2}$ from Definition \ref{def:inf_thermal_eq}.  
Let $(\mathbf{p}(s);\mathbf{q}(s))$, over $s\in\mathbb{R}_{\ge 0}$, denote the Toda state space dynamics associated with $(\mathbf{a}(t);\mathbf{b}(t))$. Further let $T\ge 1$ and $\zeta\ge 0$ be real numbers, and suppose that
\begin{equation}\label{eqn:N1N2}
N_1 \le -N(\log N)^{-7} \le N(\log N)^{-7} \le N_2,
\qquad
N = \floor{T^{100}},
\qquad
\zeta = e^{-100(\log N)^{3/2}} .
\end{equation}
\end{ass}

Let us briefly comment on Assumption \ref{ass:NT_assumption}. First, we set the interval $\llbracket N_1, N_2 \rrbracket$ to have its length $N$ be much larger than the time $T$, so that the behavior of the Toda lattice near its endpoints (which do not exist on $\mathbb{Z}$) does not affect its bulk. While in \eqref{eqn:N1N2} we made the specific choice $N = \lfloor T^{100} \rfloor$, the parameter $T$ could have been taken to grow with $N$ as quickly as $N(\log N)^{-C}$ for $C$ sufficiently large. However, we do not pursue this here, as the setting of Assumption \ref{ass:NT_assumption} will suffice to prove the results stated in Section \ref{subsec:main_res}. Second, the precise choice of $\zeta$ in \eqref{eqn:N1N2} is not essential; we only require it to decay faster than a polynomial, but slower than stretched exponentially, in $N$.

Now, following \cite[Equation (2.9)]{Agg25a}, we can define the quasi-particles mentioned in Section \ref{ParticleQ}. 

\begin{definition}[Quasi-particles] \label{def:quasi_intro}
Adopt Assumption \ref{ass:NT_assumption}. For each $s\in\mathbb{R}$, let $\mathbf{u}_j (s)\in\mathbb{R}^N$ denote a unit eigenvector of $\L(s)$ with eigenvalue $\lambda_j$.  
Under the orthonormal basis $(\mathbf{u}_1(s), \mathbf{u}_2(s),\ldots, \mathbf{u}_N(s))$ of $\L(s)$, let $\varphi_s:\llbracket 1,N\rrbracket \to \llbracket N_1,N_2\rrbracket$ be a $\zeta$-localization center bijection for $\mathbf{L}(s)$, and define the \emph{quasi-particle}
\begin{equation}\label{eqn:quasi_part_intro}
Q_j(s)=q_{\varphi_s(j)}(s),
\qquad
\text{for each }(j,s)\in\llbracket 1,N\rrbracket\times\mathbb{R}_{\ge 0}.
\end{equation}

\noindent For any $s \ge 0$ and $i \in \llbracket N_1, N_2 \rrbracket$, further define (as in \eqref{ziqt})
\begin{equation}
    \label{lambdait}
    \Lambda_i(t) = \lambda_{\varphi_t^{-1}(i)}, 
\end{equation}

\noindent which we view as the eigenvalue with localization center $i$ at time $t$. Also abbreviate $\Lambda_i = \Lambda_i (0)$.

\end{definition}

While the above defines quasi-particles for the Toda lattice on the long but finite interval $\llbracket N_1, N_2 \rrbracket$, it is also possible to define them directly on the infinite line $\mathbb{Z}$, without changing any of the asymptotics results to be stated in this section; this is done in Appendix \ref{ParticlesInfinite} below. However, for the purposes of proving the results from Section \ref{subsec:main_res}, it will in fact be more convenient to analyze the scaling limits of these quasi-particles on the interval $\llbracket N_1, N_2 \rrbracket$ above (which will in any case be used to prove their infinite volume counterparts in Appendix \ref{ParticlesInfinite}, through a limiting procedure). 

As mentioned in \eqref{qjt}, it was shown in \cite{Agg25} that (for $\theta \le \theta_0 (\beta)$ sufficiently small) the quasi-particles $Q_j (t)$ evolve approximately linearly with velocity $\ve (\lambda_j)$ (as in Definition \ref{def:v_eff_intro}), namely, $Q_{i}(t) - Q_i(0) = t \ve(\lambda_i) + O(T^{1/2+o(1)})$. The results from Section \ref{subsec:main_res} will be eventual consequences of an exact computation of the scaling limit for the fluctuations of the quasi-particles.

In order to state this result (that gives the precise form of Theorem \ref{qjt0}), we must define the limiting process. For $f \in L^2 (dr \otimes \varrho d \lambda)$, define the function $(1- \theta \mathbf{T}  \boldsymbol{\varrho}_{\beta})^{-1} f$ to be the one obtained by having the dressing operator $(1- \theta \mathbf{T}  \boldsymbol{\varrho}_{\beta})^{-1}$ act on the second argument $\lambda$ of $f$. Since $(1- \theta \mathbf{T}  \boldsymbol{\varrho}_{\beta})^{-1} : \mathcal{H} \rightarrow \mathcal{H}$ is bounded and $f \in L^2 (dr \otimes \varrho d \lambda)$, we have $(1- \theta \mathbf{T}  \boldsymbol{\varrho}_{\beta})^{-1} f \in L^2 (dr \otimes \varrho d \lambda)$. 

\begin{definition}[Dressed white noise]
\label{def:Wdress}
Let $\mathcal{W}$ denote the white noise from Definition \ref{def:white_noise_dress}. For any $f \in L^2 (dr \otimes \varrho d \lambda)$, define
\begin{equation}
\label{wf00}
   \mathcal{W}^{\dr}(f) = \mathcal{W}\left( (1- \theta \mathbf{T}  \boldsymbol{\varrho}_{\beta})^{-1} f \right).
\end{equation}
\end{definition}

Next, for any $(\Lambda, \kappa, \tau) \in \mathbb{R} \times \mathbb{R} \times \mathbb{R}_{\geq 0}$, define the function $\psi_{\Lambda,\kappa,\tau} : \mathbb{R}^2 \rightarrow \mathbb{R}$ by 
\begin{equation}\label{eqn:logphi_intro}
    \psi_{\Lambda, \kappa,\tau}(r,\lambda) \coloneqq  2 \log|\Lambda - \lambda| \cdot 
\big(
    \mathbbm{1}\{ \alpha (\kappa -  r ) > 0  \} 
    - 
    \mathbbm{1}\{ \alpha (\kappa  -  r )  +   \tau (\ve(\Lambda) - \ve(\lambda)) > 0\}
  \big).
\end{equation}

\noindent We have the following lemma on the function $
    \psi_{\Lambda, \kappa,\tau}$; its proof is in Appendix \ref{app:lim_holder}.
\begin{lem}\label{lem:psi_phi_L2_intro}
    The functions $\psi_{\Lambda, \kappa,\tau}$ from \eqref{eqn:logphi_intro} satisfy the following properties.
    \begin{enumerate}
        \item For any $(\Lambda, \kappa, \tau) \in \mathbb{R} \times \mathbb{R} \times \mathbb{R}_{\geq 0}$, we have 
                \begin{equation}
            \psi_{\Lambda, \kappa,\tau} \in L^2(\mathbb{R}^2, d r\otimes \varrho d \lambda).
        \end{equation}
        \item For any $Q, \tau \in \mathbb{R} \times \mathbb{R}_{\geq 0}$, the random function
        \begin{equation}
            \Lambda \mapsto \mathcal{W}^{\dr}(\psi_{\Lambda, (Q-\tau \ve(\Lambda)) / \alpha,\tau} ) 
        \end{equation}
        is almost surely an element of $\mathcal{H}$.
    \end{enumerate}
\end{lem}

\begin{definition}[Dressed L\'{e}vy--Chentsov field]\label{def:Z_intro}
Define the stochastic process $\mathfrak{X} : \mathbb{R}^2 \times \mathbb{R}_{\ge 0} \rightarrow \mathbb{R}$ by setting 
\begin{equation}\label{eqn:Xdef}
\mathfrak{X}(\Lambda, \kappa, \tau) =  \mathcal{W}^{\dr}(\psi_{\Lambda, \kappa,\tau}),
\end{equation}

\noindent for any $(\Lambda, \kappa, \tau) \in \mathbb{R}^2 \times \mathbb{R}_{\geq 0}$; observe by the first part of Lemma \ref{lem:psi_phi_L2_intro} that $\mathfrak{X}(\Lambda, \kappa, \tau)$ is almost surely finite. By the second part of Lemma \ref{lem:psi_phi_L2_intro}, the (random) function $\Lambda \mapsto \mathfrak{X}(\Lambda, (\mathfrak{q}-\tau \ve(\Lambda)) / \alpha, \tau)$ is almost surely an element of $\mathcal{H}$ for any $\mathfrak{q} \in \mathbb{R}$ and $\tau \in \mathbb{R}_{\ge 0}$. Thus, we may also define the \emph{dressed L\'{e}vy--Chentsov field} $\mathcal{Z}: \mathbb{R}^2 \times \mathbb{R}_{\ge 0} \rightarrow \mathbb{R}$ by setting
\begin{equation}\label{eqn:intro_Z}
\mathcal{Z}(\Lambda, \mathfrak{q}, \tau) \coloneqq \varsigma_0^{\dr}(\Lambda)^{-1} \cdot (1-\theta \mathbf{T} \boldsymbol{\varrho}_{\beta})^{-1} \mathfrak{X} \big(\Lambda, (\mathfrak{q}-\tau \ve(\Lambda)) / \alpha, \tau \big), \end{equation}

\noindent for any $(\Lambda, \mathfrak{q}, \tau) \in \mathbb{R}^2 \times \mathbb{R}_{\geq 0}$. In \eqref{eqn:intro_Z}, the dressing operator $(1-\theta \mathbf{T} \boldsymbol{\varrho}_{\beta})^{-1}$ acts in the variable $\Lambda$. 
\end{definition}

\begin{remark}[$\mathfrak{X}$ as a L\'{e}vy--Chentsov field]
\label{rmk:LCrmk}

 The Gaussian process $\mathfrak{X}(s) = \mathcal{W}^{\dr}(\psi_{\Lambda,\kappa,\tau})$ from \eqref{eqn:Xdef} is an example of what is known in the literature as a \emph{L\'{e}vy--Chentsov field} \cite{PSM,MSP,HRHF}. Recall from Section \ref{ProofFluctuationZ}  that it will govern the fluctuations of the field $\Xi$ from \eqref{xilambda}. Indeed, the form \eqref{eqn:Xdef} of $\mathfrak{X}$ (recalling $\psi$ from \eqref{eqn:logphi_intro}) looks analogous to that of $\Xi$ in \eqref{xilambda}, upon replacing $\chi(x)$ there with $\mathbbm{1}_{x>0}$; the $q_i (0)$ there with $\alpha i$; and the sum over $i$ there with an integral against $\mathcal{W}^{\dr}$. 

Let us briefly comment on the structure of $\mathfrak{X}$. Define the \emph{noise} $\omega_{\Lambda} = 2 \log|\Lambda - \lambda| \cdot \mathcal{W}^{\dr}$, which is a random distribution acting on functions $f \in L^2 (dr \otimes \varrho d\lambda)$ by
    $$\omega_{\Lambda} (f) = \mathcal{W}^{\dr}(2f \log|\Lambda - \lambda|).$$

\noindent For any $\Lambda \in \mathbb{R}$, using the noise $\omega_{\Lambda}$ one can define a function $\eta : \mathbb{R}^3 \rightarrow \mathbb{R}$, sometimes known \cite[Section 2]{HRHF} as a L\'{e}vy-Chentsov field, with respect to $\omega_{\Lambda}$ as follows. First, for any $q,\Lambda, \tau \in \mathbb{R}$, define the functions $L_{+,\Lambda} (q,\tau), L_{-,\Lambda}(q,\tau) : \mathbb{R}^2 \rightarrow \mathbb{R}$ by 
\begin{align*}
    L_{+, \Lambda}(q, \tau) (r,\lambda) \coloneqq \mathbbm{1} \{  \alpha r > q, \alpha r + \tau \ve(\lambda) < q + \tau \ve(\Lambda)\};  \\
    L_{-, \Lambda}(q,\tau) (r,\lambda) \coloneqq \mathbbm{1}\{ \alpha r < q, \alpha r + \tau \ve(\lambda) > q + \tau \ve(\Lambda)\}.
\end{align*}

\noindent One may view these functions as indexing lines crossing the space-time segment connecting $(q,0)$ to $(q+\tau\ve(\Lambda),\tau)$, from right to left and from left to right, respectively. Then, let
\begin{equation*}
    \eta(\Lambda, q, \tau) \coloneqq \omega_{\Lambda}(L_{+, \Lambda}(q, \tau)) - \omega_{\Lambda}(L_{-, \Lambda}(q, \tau)).
\end{equation*}
With this definition, the field $\mathfrak{X}(\Lambda, \kappa, \tau) = \mathcal{W}^{\dr}(\psi_{\Lambda, \kappa, \tau})$ may be identified with $\eta$ via
\begin{equation*}
    \{\mathfrak{X}(\Lambda, \kappa, \tau)\}_{(\Lambda, \kappa, \tau) \in \mathbb{R}^3} \stackrel{d}{=}  \{\eta(\Lambda, \alpha \kappa, \tau)\}_{(\Lambda, \kappa, \tau) \in \mathbb{R}^3},
\end{equation*}

\noindent which is a variant of the L\'{e}vy--Chentsov field considered in \cite[Section 2]{HRHF} (with a different choice of the Hilbert space $\mathcal{H}$ and of the covariance for the white noise $\mathcal{W}$). 

Since $\mathcal{Z}$ is a dressing of $\mathfrak{X}$, we refer to $\mathcal{Z}$ as a dressed L\'{e}vy--Chentsov field.
\end{remark}

The following theorem, proven in Section \ref{ProofParticleQ}, states that this dressed L\'{e}vy--Chentsov field governs the fluctuations for quasi-particles in the Toda lattice. The analog of Theorem \ref{thm:quasi_fluct_intro} can also be shown for the Toda lattice on the infinite line $\mathbb{Z}$; see Corollary \ref{q200} below. 

\begin{thm}[Quasi-particle fluctuations] 
	
	\label{thm:quasi_fluct_intro}
	
	For any $\beta>0$, there exists a constant $\theta_0 (\beta)>0$ such that the following holds. Adopt Assumption \ref{ass:NT_assumption}, and assume that $\theta < \theta_0 (\beta)$. For each $i \in \llbracket N_1, N_2 \rrbracket$ and $t \ge 0$, denote the quasi-particle fluctuations by
	\begin{equation}\label{eqn:quasi_part_fluct_intro}
		Z_{i}^{\mathcal Q}(t) \coloneqq  T^{-1/2} \cdot (Q_{\varphi_0^{-1}(i)}(t) - q_i(0)- t \ve(\Lambda_i)).
	\end{equation}
    
    \noindent Then there exists a constant $\mathfrak{c}>0$ and a coupling between $\L(0)$ and $\mathcal{Z}$ such that the following holds with probability at least $1-\mathfrak{c}^{-1}e^{-\mathfrak{c} (\log N)^2}$. For any integer $k \in \llbracket -T \log N,  T \log N \rrbracket$ and real number $t \in [0, T \log N]$, denoting $\kappa = kT^{-1}$ and $\tau = tT^{-1}$, we have 
    \begin{equation}
        \label{zkqtestimate00}
        \big| Z_{k}^{\mathcal Q}(t) - \mathcal{Z}(\Lambda_k, \alpha \kappa + \tau \ve(\Lambda_k), \tau) \big| \leq T^{-\mathfrak{c}}.
    \end{equation}
    
\end{thm}

Observe that Theorem \ref{thm:quasi_fluct_intro} simultaneously couples every quasi-particle (in an order $T$ neighborhood of the origin) to the dressed L\'{e}vy--Chentsov field. Thus, it provides the limiting correlations between quasi-particles that are both initially very distant (say, of order $T$ apart) and closer (related to Navier--Stokes corrections); Corollary \ref{zqi2} below provides a consequence for closer ones. 

To state the next result, we require some notation. Recall the density of states $\varrho$ of Definition \ref{def:alpha_Laxdos}, the dressing operator of Definition \ref{def:dress_intro}, and the effective velocity $\ve$ defined in Definition \ref{def:v_eff_intro}. Moreover, denote by 
\begin{flalign}
\label{eqn:T_kernel}
T(\lambda, \lambda') \coloneqq 2 \log |\lambda - \lambda'|,
\end{flalign} 

\noindent the kernel of the integral operator $\mathbf{T}$. Also let $T^{\dr}(\lambda,\lambda') = (1 - \theta \mathbf{T} \boldsymbol{\varrho_{\beta}})^{-1} T(\cdot, \lambda')$, where the dressing operator $(1 - \theta \mathbf{T} \boldsymbol{\varrho_{\beta}})^{-1}$ acts on the first argument of $T$. In what follows, we further write $X \stackrel{d}{=} Y$ if $X$ and $Y$ are random variables that have the same law.

Under this notation, we have the below proposition, which will be shown in Section \ref{subsec:tracer_D} (see the content following it for its interpretation).

\begin{prop}[Tracer quasi-particle fluctuations]\label{prop:tracer_fluct_intro}
	
	Fix $\Lambda, \mathfrak{q} \in \mathbb{R}$, and recall the Gaussian process $\mathcal{Z}$ defined in \eqref{eqn:intro_Z}. We have the representation
    \begin{equation}\label{eqn:Z_brownian}
   	\big(  \tau \mapsto \mathcal{Z} (\Lambda, \mathfrak{q}+ \tau \ve (\Lambda), \tau) \big) \stackrel{d}{=}  \big( \tau \mapsto \mathcal{D}(\Lambda)^{1/2} \cdot \mathcal{B}_\tau  \big),
    \end{equation}
 
 	\noindent where $\mathcal{B}_{\tau}$ denotes a standard Brownian motion, and 
    \begin{equation}\label{eqn:diffusivityintro}
        \mathcal{D}(\Lambda) \coloneqq |\alpha|^{-1} (\varsigma_0^{\dr}(\Lambda))^{-2} \int_{-\infty}^{\infty}  |T^{\dr}(\lambda, \Lambda)|^2 |\ve(\lambda) - \ve(\Lambda)|
    \varrho(\lambda) d\lambda.
    \end{equation}
\end{prop}
As we discuss briefly in Remark \ref{qtlambda00}, the process $\tau \mapsto  \mathcal{Z} (\Lambda, \mathfrak{q}+ \tau \ve (\Lambda), \tau)$ can be viewed as describing the fluctuations of a \emph{tracer quasi-particle} with spectral parameter $\Lambda$ and initial position $\mathfrak{q}$, around its limiting linear trajectory $\mathfrak{q} + \tau \ve(\Lambda)$. The formula \eqref{eqn:diffusivityintro} for its variance was predicted for general (quantum) integrable systems in \cite[Equation (6)]{GHKV18} and \cite[Equation (4.29)]{BDD19}.\footnote{We thank Herbert Spohn for translating these predictions in the context of the Toda lattice for us.}

We conclude this section with the following result essentially indicating that the fluctuations of quasi-particles, with nearly equal spectral parameters, are asymptotically completely correlated, as long as their initial separation is $o(T)$. It is a quick consequence of Theorem \ref{thm:quasi_fluct_intro}, together with the continuity of the dressed L\'{e}vy--Chentsov field $\mathcal{Z}$, and it will be shown in Section \ref{subsec:tracer_D}. This result can also be shown for the Toda lattice on the infinite line $\mathbb{Z}$; see Corollary \ref{zqi200} below. 

\begin{cor}

\label{zqi2}

For any $\beta>0$, there exists a constant $\theta_0 (\beta)>0$ such that the following holds. Adopt Assumption \ref{ass:NT_assumption}, and assume that $\theta < \theta_0 (\beta)$. For any real number $\delta > 0$, there is a constant $\mathfrak{c} = \mathfrak{c}(\delta) > 0$ such that the below holds with probability at least $1 - \mathfrak{c}^{-1} e^{-\mathfrak{c} (\log N)^2}$. For any integers $k, k' \in \llbracket -T \log N, T \log N \rrbracket$ satisfying $|k-k'| \le T^{1-\delta}$ and $|\Lambda_k - \Lambda_{k'}| \le T^{-\delta}$, we have that 
\begin{flalign*} 
\displaystyle\sup_{t \in [0, T \log N]} |Z_k^{\mathcal{Q}} (t) - Z_{k'}^{\mathcal{Q}} (t)| \le T^{-\mathfrak{c}}.
\end{flalign*} 

\end{cor}

In the case $\delta \ge 1/2$, the counterpart of Corollary \ref{zqi2} was proven in \cite[Section 3]{FO25} for the hard rods model and as \cite[Theorem 4.11]{SLSBS} for the box ball system (our time parameter $T$ is denoted by $\varepsilon^{-2}$ in \cite{FO25} and by $n^2$ in \cite{SLSBS}).

\subsection{Notation}

Throughout, $c$ and $C$ will denote small and large constants, respectively, which may change between lines. For any real numbers $a<b$, we set $\llbracket a,b \rrbracket = [a,b] \cap \mathbb{Z}$.  For any interval $\mathcal{I} \subseteq \mathbb{R}$ or event $\mathsf{E}$, let $\mathcal{I}^{\complement} = \mathbb{R} \setminus \mathcal{I}$ and $\mathsf{E}^{\complement}$ denote the complements of $\mathcal{I}$ and $\mathsf{E}$, respectively. For any integer $d \ge 1$ and point $v = (v_1, v_2, \ldots , v_d) \in \mathbb{C}^d$, we set $\| v \|_2 = (\sum_{i=1}^d v_i^2)^{1/2}$. For any two subsets $\mathcal{S}_1, \mathcal{S}_2 \subseteq \mathbb{C}^d$, we also let $\dist (\mathcal{S}_1, \mathcal{S}_2) = \inf_{s_1 \in \mathcal{S}_1} \inf_{s_2 \in \mathcal{S}_2} \|s_1 - s_2 \|_2$.

 Let $N \ge 1$ be an integer and $\mathsf{E}_N$ be an event depending on $N$ (and possibly other parameters). We say that $\mathsf{E}_N$ \emph{holds with overwhelming probability}, or is \emph{overwhelmingly probable}, if there exists a constant $c > 0$ such that $\mathbb{P} [\mathsf{E}_N^{\complement}] \le c^{-1} e^{-c(\log N)^2}$. Here, $c$ is an \emph{absolute constant}, meaning that it is independent of $N$ (and all other parameters involved in the definition of $\mathsf{E}_N$, except for possibly $\beta$ and $\theta$ from \eqref{eqn:equil}). Observe that, whenever proving that $\mathsf{E}_N$ is overwhelmingly probable, we may assume that $N \ge N_0$ is sufficiently large; we will often do so implicitly (and without comment). 

We will often study fluctuations of the form $X - \mathbb{E}[X]$, for a random variable given by a fairly intricate sum. So, to ease notation, we may occasionally abbreviate $X - \mathbb{E}[X] = X - \mathbb{E}[\cdots]$ (where it should be understood that $\mathbb{E}[\cdots] = \mathbb{E}[X]$ is chosen to ensure that $X - \mathbb{E}[\cdots]$ has mean $0$).

At the beginning of each (sub)section, we will recall the notation that will be used (without comment) throughout the (sub)section.

\section{Proof of quasi-particle fluctuations}

\label{ProofQ}

As we will see in Sections \ref{sec:current_fluct} and \ref{sec:charge_fluct} below, our results in Section \ref{subsec:main_res} will follow from the quasi-particle fluctuations given by Theorem \ref{thm:quasi_fluct_intro} (together with known results approximating currents and local charges through quasi-particles). So, in this section, we establish the latter, conditional on various intermediate results that will be proven later in this paper. Throughout this section, we recall the notation from Assumption \ref{ass:NT_assumption} and Definition \ref{def:quasi_intro}.

\subsection{Preliminary identities and estimates}

\label{Estimate0} 

In this section we collect several identities and estimates that will be used at various points while proving Theorem \ref{thm:quasi_fluct_intro}. We first recall the following properties of $\varrho$ (of Definition \ref{def:alpha_Laxdos}) due to \cite{GM23} (though stated in the below form in \cite{Agg25}).

\begin{lem}[{\cite[Lemma 3.1]{Agg25}}]\label{lem:varrho_bd}

There exists a constant $C>1$ such that
\begin{equation}\label{eqn:varrho_tail_bds}
\begin{aligned}
\varrho_\beta(x)
&\le C (|x|+1)^{2\theta} e^{-\beta x^2/2};
\qquad
\varrho(x)
&\le C (|x|+1)^{2\theta+1} e^{-\beta x^2/2}.
\end{aligned}
\end{equation}
\end{lem}

The following lemma provides bounds on the dressing operator, the effective velocity, and $\varsigma_0^{\dr}$ (recall Definition \ref{def:dress_intro}). Its first part is due to \cite[Lemma 3.6]{Agg25}, its second to \cite[Corollary 3.8]{Agg25}, and its third to \cite[Lemma 1.7]{Agg25}.

	\begin{lem}[\cite{Agg25}]
		\label{lem:ve_tail}
		
		There exist constants $c>0$ and $C>1$ such that the following three statements hold. 
        
        \begin{enumerate} 
        \item For any function $f \in \mathcal{H}$ and real number $x \in \mathbb{R}$, we have 
		\begin{flalign*} 
			|f^{\dr} (x)| \le C \cdot |f(x)| + C \| f \|_{\mathcal{H}} \cdot \log (|x|+2),
		\end{flalign*} 
and
    \begin{flalign*} 
		\bigl|\partial_x f^{\mathrm{dr}}(x)\bigr|
\le C\cdot\bigl|f'(x)\bigr|
\;+\; C\bigl(\|f'\|_{\mathcal H} + \|f\|_{\mathcal H}\bigr)\,\log\bigl(|x|+2\bigr).
		\end{flalign*} 

        \item For any real number $A \ge 2$, we have
		\begin{flalign}\label{eqn:vebd}
			\displaystyle\sup_{|x| \le A} | \ve (x) | \le CA; \qquad \displaystyle\sup_{|x| \le A} | \partial_x \ve (x) | \le CA \log A. 
		\end{flalign}

        \item For any real number $x \in \mathbb{R}$, we have $\varsigma_0^{\dr} (x) \cdot \sgn (\alpha) \ge c$. 

        \end{enumerate} 
        
	\end{lem}

    We next have the following identities (shown in \cite{Agg25}, upon taking into account the fact that the dressing operator of \cite[Definition 1.6]{Agg25} differs from ours \eqref{drf} by a factor of $\theta$).  
    
	\begin{lem}[{\cite[Corollary 3.4]{Agg25}}]
		\label{rho0}
		
		For any $x \in \mathbb{R}$, we have 
		\begin{flalign} 
			\label{rho2} 
			\varrho (x) = \alpha \theta \cdot \varsigma_0^{\dr} (x) \cdot \varrho_{\beta} (x), \qquad \text{and} \qquad \varsigma_0^{\dr} (x) = \alpha^{-1} \cdot \big( \bm{\mathrm{T}} \varrho (x) + \alpha \big).
		\end{flalign} 
		
	\end{lem}
    
    Next we recall distance bounds between the Toda particles $(q_i(s))$ under thermal equilibrium. As indicated below Definition \ref{def:alpha_Laxdos}, the stretch $\alpha = \mathbb{E}[q_i - q_{i+1}]$ is the average spacing between these particles, one should expect $q_i - q_j  \approx \alpha (i-j)$. This is made precise through the below lemma (which essentially amounts to a Chernoff bound).

\begin{lem}[{\cite[Lemma 7.2]{Agg25a}}]\label{lem:q_spacing}

Adopt Assumption \ref{ass:NT_assumption}, but assume more generally that $T \in [0, N (\log N)^{-7}]$. The following two statements hold with overwhelming probability.

		\begin{enumerate} 
			
			\item For any $s \in [0,T]$ and $i, j \in \llbracket N_1 + T(\log N)^3, N_2 - T (\log N)^3 \rrbracket$, we have
			\begin{flalign}
				\label{qiqjs4}
				\big| q_i (s) - q_j (s) - \alpha (i-j) \big| \le |i-j|^{1/2} (\log N)^2.
			\end{flalign}
			
			\item If $T \geq 1$, then for any $s \in [0,T]$ and $i \in \llbracket N_1, N_2 \rrbracket$ with $|i-j| \ge T (\log N)^5$, we have 
			\begin{flalign}
				\label{qiqjs5}
				\big(  q_i (s) - q_j (s) \big) \cdot \sgn (\alpha i - \alpha j) \ge \displaystyle\frac{|\alpha|}{2} \cdot |i-j|.
			\end{flalign}
			
		\end{enumerate} 

\end{lem}

The next result bounds the entries and eigenvalues of the Lax matrix $\mathbf{L}(t)$ (recall Definition \ref{lt}) under the evolution of the Toda lattice; it requires the below definition.

\begin{definition}

\label{eigenvaluesm}

    If $\mathbf{M}$ is an $N \times N$ symmetric matrix, then we let $\eig \mathbf{M}$ denote the eigenvalues of $\mathbf{M}$, $\eig \mathbf{M} = \{\lambda_1(M) \geq \cdots \geq \lambda_N(M)\}$. For any real number $A \ge 0$, define the event
\begin{equation}
\label{eventbounded}
    \mathsf{BND}_{\mathbf{M}}(A) \coloneqq  \bigcap_{1 \le i,j \le N} \{ |M_{i j}| \leq A\}  \cap   \bigcap_{\lambda \in \eig \mathbf{M}} \{ |\lambda| \leq A\}.
\end{equation}
\end{definition}

\begin{lem}[{\cite[Lemma 3.15]{Agg25}}]
\label{lem:bd_lem}

Under Assumption \ref{ass:NT_assumption}, there exists a constant $c>0$ such that the following holds. For any real number $A \ge 1$, we have
\begin{flalign*}
    \mathbb{P}\Bigg( \bigcap_{t \ge 0} \mathsf{BND}_{\mathbf{L}(t)} (A) \Bigg) \ge 1 - c^{-1} N e^{-cA^2}.
\end{flalign*}

\noindent In particular, $\bigcap_{t \ge 0}\mathsf{BND}_{\L(t)}(\log N)$ holds with overwhelming probability. 

\end{lem} 

The last result we recall here is a type of continuity estimate for the localization centers (recall Definition \ref{center}) of the Lax matrix $\mathbf{L}(t)$, which are not too close to the boundary of $\llbracket N_1, N_2 \rrbracket$. Observe in particular that taking $t=t'$ in \eqref{eqn:loc_cent_dif} below implies that these localization centers are approximately unique, up to a logarithmic error (that will frequently be negligible in our context). 

\begin{lem}[{\cite[Lemmas 3.21 and 3.22]{Agg25}}]\label{lem:loc_cent_dif}

Adopt \Cref{ass:NT_assumption}, but assume more generally that $\zeta \geq e^{-150 (\log N)^{3/2}}$ and $T \leq N (\log N)^{-6}$. The following holds with overwhelming probability. 

\begin{enumerate} 

\item Fix any eigenvalue \(\lambda \in \mathrm{eig}\,\L(0)\) and any \(\zeta\)-localization center
\(\varphi \in \llbracket N_1,N_2 \rrbracket \) of \(\lambda\) with respect to \(\L(0)\). Then, for each real number
\(s \in [0,T]\), there does not exist an index \(m \in \llbracket N_1,N_2 \rrbracket \) satisfying
\[
|m-\varphi| \ge T(\log N)^2,
\]
such that \(m\) is a $\zeta$-localization center for \(\lambda\) with respect to \(\L(s)\). 

\item Fix any real numbers $t, t' \in [0, T]$; an eigenvalue $\lambda \in \eig \mathbf{L}(t)$; and $\zeta$-localization centers $\varphi \in \llbracket N_1, N_2 \rrbracket$ and $\varphi' \in \llbracket N_1, N_2 \rrbracket$ of $\lambda$ with respect to $\mathbf{L}(t)$ and $\mathbf{L}(t')$, respectively. If $N_1 + T (\log N)^4 \le \varphi \le N_2 - T (\log N)^4$, then
\begin{equation}\label{eqn:loc_cent_dif}
|\varphi - \varphi'| \le (|t-t'| + 2) (\log N)^3. 
\end{equation}
\end{enumerate} 
\end{lem}

\subsection{Proof of Theorem \ref{thm:quasi_fluct_intro}}

\label{ProofParticleQ} 

In this section we establish Theorem \ref{thm:quasi_fluct_intro}, assuming Theorems \ref{prop:dressing_Z_approx_outline} and \ref{prop:zcouple_outline} below. Throughout, we recall the quasi-particles $Q_k(t)$ from Definition \ref{def:quasi_intro}. We begin with the following definition for a smoothing $\chi$ of the step function (on some scale~$\mathfrak{M}$) and smoothed logarithm $\mathfrak{l}$.

\begin{definition} \label{def:chi_def} \label{def:smoothed_log}

Adopt Assumption \ref{ass:NT_assumption}; set $\mathfrak{M} = T^{19/20}$; and define $\mathfrak{l} = \mathfrak{l}_N : \mathbb{R} \rightarrow \mathbb{R}$ by 
 \begin{equation}\label{eqn:sl}
 \mathfrak{l}(x) = \frac{1}{2}\log (x^2 + e^{-10 (\log N)^2}).
 \end{equation}
 
 \noindent Also let $\chi : \mathbb{R} \rightarrow \mathbb{R}_{\ge 0}$ be a smooth function, with $\chi'$ even and $0 \le \chi(x) \le 1$ for all $x \in \mathbb{R}$, such that $\chi(x) \equiv 1$ for~$ x > \mathfrak{M}$ and $\chi(x) \equiv 0$ for~$x < - \mathfrak{M}$, and 
\begin{equation}
    \label{estimatesm0} 
\begin{gathered}
\supp \chi' \subseteq [-\mathfrak{M},\mathfrak{M}], \\
\sup_{x \in \mathbb{R}} |\chi'(x)| \leq 10 \mathfrak{M}^{-1}, \qquad
\sup_{x \in \mathbb{R}} |\chi''(x)| \leq 100 \mathfrak{M}^{-2}, \qquad
\sup_{x \in \mathbb{R}} |\chi'''(x)| \leq 1000 \mathfrak{M}^{-3}.
\end{gathered}
\end{equation} 

\end{definition}

     The starting point for the proof of Theorem \ref{thm:quasi_fluct_intro} is contained in the following lemma. Its first part, given by \cite[Proposition 6.2]{Agg25}, provides the smoothened asymptotic scattering relation for the Toda lattice quasi-particles, as described in \eqref{ziqt00}. Its second part, given by \cite[Theorem 1.13]{Agg25}, provides the linear limit shape for the quasi-particles, as described in \eqref{qjt} and \eqref{ziqt0}. In what follows, we recall the effective velocity $\ve$ from Definition \ref{def:v_eff_intro} and the notation of Theorem \ref{thm:quasi_fluct_intro}.

 \begin{lem}[{\cite{Agg25}}]
 \label{lem:asymptotic-scattering}
 \label{thm:Z_apriori}
Adopt Assumption \ref{ass:NT_assumption}, but more generally assume $1 \leq T \leq N (\log N)^{-7}$. The following two statements hold with overwhelming probability, for any integer $k \in \llbracket 1,N \rrbracket$ satisfying
\begin{equation}\label{kn1n2}
N_1 + T (\log N)^6 \le \varphi_0 (k) \le N_2 - T (\log N)^6 .
\end{equation}

\begin{enumerate} 
\item Recalling Definition \ref{def:chi_def}, we have
\begin{flalign}\label{eq:2.15}
\begin{aligned} 
\displaystyle\sup_{t \in [0,T]} \Bigg|
\lambda_k t - Q_k(t) + Q_k(0)
- 2 \sum_{i=1}^{N}
\mathfrak{l}(\lambda_k - \lambda_i) \cdot 
\big(
\chi(Q_k(t) - Q_i(& t))
- \chi( Q_k(0) - Q_i(0))
\big)
\Bigg| \\
& \le \mathfrak{M}^{1/2} (\log N)^{16}.
\end{aligned} 
\end{flalign}

\item Recalling $Z_{\varphi_0 (k)}^{\mathcal{Q}}(t)$ from \eqref{eqn:quasi_part_fluct_intro}, there exists a constant
$\theta_0=\theta_0(\beta)>0$ such that the following holds. If $\theta\in(0,\theta_0)$, then
  \begin{equation}\label{eqn:Z_apriori}
  \displaystyle\sup_{t \in [0,T]} |Z_{\varphi_0 (k)}^{\mathcal Q} (t)| \leq (\log N)^{35}.
  \end{equation}

  \end{enumerate} 
\end{lem}

We next define the following two three-parameter families of functions. The first is denoted by $\Xi$, which was described by \eqref{xilambda}, and provides the fluctuations of the sum on the right side of \eqref{eq:2.15}, upon replacing the $(Q_i(t))$ by their limiting trajectories $(Q_i(0) + t \ve(\lambda_i))$. The second, denoted by $Z$, is a certain dressing of $\Xi$, as mentioned at the end of Section \ref{ProofFluctuationZ}.

\begin{definition} 

\label{xilambdakt} 

Recall Definition \ref{def:chi_def}. For any integer $k \in \llbracket N_1, N_2 \rrbracket$ and real numbers $\Lambda \in \mathbb{R}$ and $t \geq 0$, define
\begin{multline}
    \label{eqn:psiZktau}
 \Xi(\Lambda, k, t) =   \frac{1}{T^{1/2}} \Bigg( 2 \sum_{i = N_1}^{N_2} \mathfrak{l}(\Lambda- \Lambda_i) \cdot 
\Big(  
    \chi(q_{k}(0) - q_i(0)) - 
    \chi \big(q_{k}(0) - q_i(0) + t  \left(\ve(\Lambda) - \ve(\Lambda_i) \right) \big)
   \Big)  \\
    - \mathbb{E} \bigg[  2 \sum_{i = N_1}^{N_2} \mathfrak{l}(\Lambda- \Lambda_i) \cdot 
\Big( 
    \chi(q_{k}(0) - q_i(0)) - 
    \chi \big(q_{k}(0) - q_i(0) + t  \left(\ve(\Lambda) - \ve(\Lambda_i) \right) \big)
    \Big)  \bigg] \Bigg).
 \end{multline}
 When $k$ is not an integer, we extend $ \Xi(\Lambda, k, t)$ by linear interpolation. Also, if $k \leq N_1$, we simply set $\Xi(\Lambda, k, t) = \Xi(\Lambda, N_1, t)$, and if $k \geq N_2$, we set $\Xi(\Lambda, k, t) = \Xi(\Lambda, N_2, t)$. 

 Moreover, for any real numbers $\Lambda, Q \in \mathbb{R}$ and $t \geq 0$, define  (recalling Definition \ref{def:dress_intro})
\begin{equation} \label{eqn:Z_outline}
    Z(\Lambda,Q, t) \coloneqq 
     \varsigma_0^{\dr}(\Lambda)^{-1} \cdot \Big( \Lambda \mapsto \Xi \big(\Lambda, (Q - t  \ve(\Lambda))/\alpha, t \big) \Big)^{\dr}.
\end{equation}

 \end{definition} 

 Theorem \ref{thm:quasi_fluct_intro} will follow quickly from the below two theorems. The first, proven in Section \ref{sec:quasi_ptcl_couple}, states that the quasi-particle fluctuations \eqref{eqn:quasi_part_fluct_intro} approximate the function $Z$ from \eqref{eqn:Z_outline}. The second, proven in Section \ref{CoupleXiX}, states that the latter function $Z$ can with high probability be coupled to approximate the Gaussian process $\mathcal{Z}$ from \eqref{eqn:intro_Z}. In what follows, and in later sections, the range of the index parameter $k$ (or spatial parameters $q$ and $Q$ below) will be taken to be considerably larger than $\llbracket -T, T \rrbracket$. However, the exact form of this range will not be so relevant.

\begin{thm}[Quasi-particle fluctuations approximate $Z$]\label{prop:dressing_Z_approx_outline}
Adopt Assumption \ref{ass:NT_assumption}. There is a constant $\theta_0 (\beta) > 0$ such that the following holds whenever $\theta < \theta_0 (\beta)$. There exists a constant $\mathfrak{c} > 0$ such that, with overwhelming probability, we have 
\begin{equation}\label{eqn:Z_approx_outline}
\displaystyle\sup_{t \in [0, T \log N]} \displaystyle\max_{k \in \llbracket -T^{7/6}, T^{7/6} \rrbracket} |Z_{k}^{\mathcal Q} (t) - Z(\Lambda_k, \alpha k + t  \ve(\Lambda_k), t)| \leq T^{-\mathfrak{c}}.
\end{equation}
\end{thm}

\begin{thm}[Limit for $Z$]\label{prop:zcouple_outline}

Adopt Assumption \ref{ass:NT_assumption}. There is a constant $\theta_0 (\beta) > 0$ such that the following holds whenever $\theta < \theta_0 (\beta)$. There exists a constant $\mathfrak{c} > 0$ and a coupling between $\mathcal{Z}$ with $\L(0)$ such that, with overwhelming probability, we have  
\begin{equation}\label{eqn:tracer_couple_outline}
\displaystyle\sup_{\tau \in [0, \log N]} \displaystyle\sup_{|\Lambda| \le \log N} \displaystyle\sup_{|\mathfrak{q}| \le (\log N)^3} \left| \mathcal{Z}(\Lambda, \mathfrak{q}, \tau) - Z(\Lambda, T \mathfrak{q}, T \tau) \right| \leq T^{-\mathfrak{c}}.
\end{equation}
\end{thm}

\begin{proof}[Proof of Theorem \ref{thm:quasi_fluct_intro}]

By Theorems \ref{prop:dressing_Z_approx_outline} and \ref{prop:zcouple_outline}, to show the theorem, it suffices to verify that, for each integer $k \in \llbracket -T \log N, T \log N \rrbracket$ and real number $t \in [0, \log N]$, the parameters $(\Lambda_k, \alpha \kappa + \tau \ve(\Lambda_k), \tau)$ are likely in the range specified by \eqref{eqn:tracer_couple_outline}. Specifically, denoting $k = \kappa T$, $t = \tau T$, and $\mathfrak{q} = \alpha \kappa + \tau \ve(\Lambda_k)$, we must confirm that $(\Lambda_k, \mathfrak{q}, \tau) \in [-\log N, \log N] \times [-(\log N)^3, (\log N)^3] \times [0, \log N]$ holds, for all $(k,t) \in \llbracket -T \log N, T \log N \rrbracket \times  [0, \log N]$, with overwhelming probability. 

Fix such a pair $(k,t)$. By Lemma \ref{lem:bd_lem}, the event $\mathsf{BND}_{\mathbf{L}(0)} (\log N)$ from Definition \ref{eventbounded} holds with overwhelming probability, so we may restrict to it. Then $|\Lambda_k| \le \log N$, so the second part of Lemma \ref{lem:ve_tail} gives $|\ve(\Lambda_k)| \le C \log N$. Hence, $|\mathfrak{q}| \le (|\alpha| \kappa + C \tau) \log N \le (\log N)^3$, so $(\Lambda_k, \mathfrak{q}, \tau) \in [-\log N,\log N] \times [-(\log N)^3, (\log N)^3] \times [0, \log N]$, which shows the theorem.
\end{proof}

\subsection{Proof outline for Theorem \ref{prop:dressing_Z_approx_outline}}
\label{subsec:heur}

\subsubsection{Linearizing the asymptotic scattering relation}

To prove Theorem \ref{prop:dressing_Z_approx_outline}, we first show the following lemma. It indicates that (upon fixing the eigenvalues $(\Lambda_i)$ of $\mathbf{L}(0)$) the quasi-particle fluctuations $Z_k^{\mathcal Q}$ approximately solve a linear equation in terms of $\Xi$, as described in \eqref{qktqk01}.

\begin{lem}[Linearized asymptotic scattering relation]\label{lem:Z_eq_outline}
Adopt Assumption \ref{ass:NT_assumption}. There is a constant $\theta_0 (\beta) > 0$ such that the following holds whenever $\theta < \theta_0 (\beta)$. There exists a constant $\mathfrak{c}>0$ so that, with overwhelming probability, we have
\begin{flalign}\label{eqn:Zeq_1}
\begin{aligned}
   \displaystyle\sup_{t \in [0,T \log N]} \displaystyle\sup_{k \in \llbracket -T^{3/2}, T^{3/2} \rrbracket} \Bigg| Z_{k}^{\mathcal Q} (& t) + 2 \sum_{i=N_1}^{N_2} ( Z_{k}^{\mathcal Q} (t)- Z_{i}^{\mathcal Q} (t)) \cdot \mathfrak{l}(\Lambda_{k}-\Lambda_i) \\
   &  \times \chi' \big(\alpha( k- i) +  t  (\ve(\Lambda_{k})  - \ve(\Lambda_{i}) ) \big)-  \Xi(\Lambda_{k}, k, t) \Bigg| \leq T^{-\mathfrak{c}}.
   \end{aligned} 
\end{flalign}
\end{lem}

Lemma \ref{lem:Z_eq_outline} essentially follows from Taylor expanding the asymptotic scattering relation \eqref{eq:2.15} (as was also anticipated in \cite[Section 2.3]{Agg25}). The following two lemmas will facilitate this. The first, to be shown in Section \ref{LinearExpectation}, approximately evaluates the mean of the sum involved in the definition \eqref{eqn:psiZktau} of $\Xi$. The second, to be shown in Appendix \ref{ProofQDifference}, bounds the arguments of $\chi$ in the definition \eqref{eqn:psiZktau} of $\Xi$, enabling us to estimate the number of terms in the support of $\chi'$. 

\begin{lem}\label{lem:xi_expectation_outline}
    Adopt Assumption \ref{ass:NT_assumption}. There exist constants $\theta_0 (\beta) > 0$ and $\mathfrak{c}>0$ such that the following holds whenever $\theta < \theta_0 (\beta)$. For any $\Lambda \in [-\log N, \log N]$, $t \in [0,  T \log N]$, and $k \in \llbracket - T^{3/2}, T^{3/2} \rrbracket$, we have 
    \begin{flalign*}
    \begin{aligned} 
        \Bigg|\mathbb{E} \bigg[\sum_{i=N_1}^{N_2}\mathfrak{l}(\Lambda-\Lambda_i) \big( \chi(q_k(0)  -q_i(0) + t \ve(\Lambda)- t \ve(\Lambda_i))- \chi(q_k(0)  -q_i(0)  )\big) \bigg]  -  t( & \Lambda - \ve (\Lambda))  \Bigg| \\
        & 
        \leq   T^{1/2-\mathfrak{c}}.
        \end{aligned} 
    \end{flalign*}
\end{lem}

\begin{lem}\label{lem:second_der_bd}
Adopt Assumption \ref{ass:NT_assumption}. There is a constant $\theta_0 (\beta) > 0$ such that the following holds whenever $\theta < \theta_0 (\beta)$. For any real number $t \in [0, T \log N ]$, and integers $i, k \in \llbracket N_1, N_2 \rrbracket$ with $\varphi_0 (k) \in \llbracket -T^2, T^2 \rrbracket$ with $|\varphi_t(i)-\varphi_t (k)| \ge \mathfrak{M}(\log N)^5$, we have with overwhelming probability that
\begin{flalign*}
     & \min \big( Q_k(t) - Q_i(t), Q_k(0) - Q_i(0) + t(\ve (\lambda_k) - \ve (\lambda_i)) \big) > 10 \mathfrak{M}, \quad \text{if $\varphi_t (k) > \varphi_t (i)$}; \\
     & \max \big( Q_k(t) - Q_i(t), Q_k(0) - Q_i(0) + t(\ve (\lambda_k) - \ve (\lambda_i)) \big) < -10 \mathfrak{M}, \quad \text{if $\varphi_t (k) < \varphi_t (i)$}.
\end{flalign*}

\end{lem}

\begin{proof}[Proof of Lemma \ref{lem:Z_eq_outline}]

Throughout this proof, for any integer $i \in \llbracket N_1, N_2 \rrbracket$, we abbreviate $\hat{i} = \varphi_0^{-1} (i)$. We moreover fix an integer $k \in \llbracket -T^{3/2}, T^{3/2} \rrbracket$, and let $t \in [0, T \log N]$ be any real number. 

First apply Lemma \ref{lem:asymptotic-scattering}, with the $k$ and $T$ there equal to $\hat{k} = \varphi_0^{-1} (k)$ and $N^{4/5} \le N(\log N)^{-7}$ here, respectively. Since \eqref{kn1n2} holds (as $\varphi_0 (\hat{k}) = k \in \llbracket -T^{3/2}, T^{3/2} \rrbracket \subseteq [N_1 + N^{4/5} (\log N)^6, N_2 - N^{4/5} (\log N)^6]$), \eqref{eq:2.15} does as well. Together with the facts that $(Q_{\hat{i}} (0), \lambda_{\hat{i}}) = (q_i (0), \Lambda_i)$ (by Definition \ref{def:quasi_intro}) and that $Q_{\hat{i}} (t) = q_i (0) + t \ve (\Lambda_i) + T^{1/2} \cdot Z_i^{\mathcal{Q}} (t)$ (by \eqref{eqn:quasi_part_fluct_intro}), this yields 
\begin{flalign}
    \label{sumlambdakq} 
    \begin{aligned} 
    \Bigg| 2 \sum_{i=N_1}^{N_2}
\mathfrak{l}(\Lambda_k - \Lambda_i) \cdot 
\big( 
\chi (Q_{\hat{k}} (t ) - Q_{\hat{i}} (t))
- \chi (q_k (0) - q_i (0))
\big)  + t (\ve (\Lambda_k & ) - \Lambda_k) + T^{1/2} \cdot Z_k^{\mathcal{Q}}(t) \Bigg|  \\
&  \le
\mathfrak{M}^{1/2} (\log N)^{16}.
\end{aligned} 
\end{flalign}

We next Taylor expand the summands involving $\chi(Q_{\hat{k}} (t)- Q_{\hat{i}} (t))$ in \eqref{sumlambdakq}. Recalling the definition \eqref{eqn:quasi_part_fluct_intro} of $Z_i^{\mathcal{Q}} (t)$ and abbreviating $Z_i^{\mathcal{Q}} (t) = Z_i^{\mathcal{Q}}$, this yields for each $i \in \llbracket N_1, N_2 \rrbracket$ a real number $\xi_i$ between $Q_{\hat{k}} (t)- Q_{\hat{i}} (t)$ and $q_k(0) - q_i(0) + t(\ve (\Lambda_k) - \ve (\Lambda_i))$ such that 
\begin{flalign}\label{eqn:te_chi}
\begin{aligned} 
    \chi( Q_{\hat{k}} (t)- Q_{\hat{i}} (t)) &  = \chi \big(q_k(0)-q_i(0) + t (\ve(\Lambda_k)-\ve(\Lambda_i)) \big) + \frac{T}{2} \cdot \chi''(\xi_i) \cdot (Z_{k}^{\mathcal{Q}}- Z_{i}^{\mathcal{Q}})^2 \\
    & \qquad + T^{1/2} \cdot \chi' \big(q_k(0)-q_i(0) + t (\ve(\Lambda_k)-\ve(\Lambda_i)) \big) \cdot (Z_{k}^{\mathcal{Q}}- Z_{i}^{\mathcal{Q}}). 
    \end{aligned} 
\end{flalign}

By Lemma \ref{lem:second_der_bd}, with the fact that $\supp \chi \subseteq [-\mathfrak{M}, \mathfrak{M}]$, both sides of \eqref{eqn:te_chi} are equal to $0$ unless $|\varphi_t (\hat{k})- \varphi_t (\hat{i})| \le \mathfrak{M} (\log N)^5$, with overwhelming probability. We restrict to this event in what follows. Further observe by Lemma \ref{lem:loc_cent_dif} that $\varphi_t (\hat{k}) \in \llbracket -2T^{3/2}, 2T^{3/2} \rrbracket$ holds with overwhelming probability, since $k = \varphi_0 (\hat{k}) \in \llbracket -T^{3/2}, T^{3/2} \rrbracket$ and $T = N^{1/100}$ (by \eqref{eqn:N1N2}). On this event, if $|\varphi_t (\hat{k})- \varphi_t (\hat{i})| \le \mathfrak{M} (\log N)^5$, then $\varphi_t (\hat{i}) \in \llbracket N_1 + T^2, N_2 - T^2 \rrbracket$. So, Lemma \ref{thm:Z_apriori} yields with overwhelming probability that $|Z_i^{\mathcal Q}| \leq (\log N)^{35}$, whenever $|\varphi_t (\hat{k})- \varphi_t (\hat{i})| \le \mathfrak{M} (\log N)^5$. We further restrict to this event below. Inserting these estimates, together with \eqref{eqn:te_chi}, into \eqref{sumlambdakq}, then gives 
\begin{flalign*}
\Bigg| &  
2 \sum_{i=N_1}^{N_2}
\mathfrak{l}(\Lambda_k - \Lambda_i) \cdot 
\Big( \chi \big(q_k(0)-q_i(0) + t (\ve(\Lambda_k)-\ve(\Lambda_i)) \big) - \chi( q_k(0) - q_i(0))
\Big) - t(\Lambda_k - \ve (\Lambda_k)) \\ 
& \qquad + T^{1/2} \bigg( Z_k^{\mathcal{Q}} + 2 \displaystyle\sum_{i=N_1}^{N_2} \mathfrak{l}(\Lambda_k - \Lambda_i) \cdot  \chi' \big(q_k(0)-q_i(0) + t (\ve(\Lambda_k)-\ve(\Lambda_i)) \big) \cdot (Z_{k}^{\mathcal{Q}}- Z_{i}^{\mathcal{Q}})  \bigg) 
\Bigg| \\
& \qquad \qquad \le T (\log N)^C \displaystyle\sum_{i=N_1}^{N_2} |\chi'' (\xi_i)| \cdot \mathbbm{1}_{|\varphi_t (\hat{k})- \varphi_t (\hat{i})| \le \mathfrak{M} (\log N)^5} + B \mathfrak{M}^{1/2} (\log N)^{16} \\
& \qquad \qquad \le (\log N)^{2C} ( T \mathfrak{M}^{-1} + \mathfrak{M}^{1/2}) \le B (\log N)^{3C} \mathfrak{M}^{1/2},
\end{flalign*} 

\noindent where in the last inequality we used \eqref{estimatesm0} and the fact that $\mathfrak{M} \ge T^{3/4}$ (recall Definition \ref{def:chi_def}). This, together with the definition \eqref{eqn:psiZktau} of $\Xi$ and Lemma \ref{lem:xi_expectation_outline} (and the facts that $B \mathfrak{M}/T \le T^{-c}$ and that $\log N = 100 \log T$, by \eqref{eqn:N1N2}), yields 
\begin{flalign}
\label{eqn:zeqnbd_out2}
\begin{aligned} 
\Bigg| Z_k^{\mathcal{Q}} + 2 \displaystyle\sum_{i=N_1}^{N_2} \mathfrak{l} (\Lambda_k - \Lambda_i) \cdot \chi' \big(q_k(0)-q_i(0) + t (\ve(\Lambda_k)-\ve(\Lambda_i)) \big) \cdot (& Z_{k}^{\mathcal{Q}}- Z_{i}^{\mathcal{Q}}) \\
& - \Xi ( \Lambda_k, k, t) \Bigg| \le T^{-c}.
\end{aligned} 
\end{flalign} 

To show \eqref{eqn:Zeq_1}, it remains to replace $q_k(0)-q_i(0)$ in the argument of $\chi'$ in \eqref{eqn:zeqnbd_out2} with $\alpha (k - i )$, namely, to show with overwhelming probability that 
\begin{flalign}
    \label{sumalphaq}
    \begin{aligned} 
\displaystyle\sum_{i=N_1}^{N_2} \Big| \chi' \big(q_k(0)-q_i(0) + t (\ve(\Lambda_k)-\ve(\Lambda_i)) \big) -  \chi' & \big( \alpha (k-i) + t (\ve(\Lambda_k)-\ve(\Lambda_i)) \big) \Big|  \\
& \times \mathbbm{1}_{|\varphi_t (\hat{k})-\varphi_t (\hat{i})|\le \mathfrak{M} (\log N)^5} \le T^{-c},
\end{aligned} 
\end{flalign}

\noindent where we have again used the facts (from Lemma \ref{lem:second_der_bd} and the inclusion $\supp \chi \subseteq [-\mathfrak{M}, \mathfrak{M}]$) that the summand in \eqref{sumalphaq} is equal to $0$ unless $|\varphi_t(\hat{k})-\varphi_t (\hat{i})| \le \mathfrak{M} (\log N)^5$, and that $|Z_i^{\mathcal{Q}}| \le (\log N)^C$ and $|\mathfrak{l}(\Lambda_k-\Lambda_i)| \le (\log N)^C$ for such $i$. By a Taylor expansion and the $T=0$ case of Lemma \ref{lem:q_spacing} (together with the fact that, by Lemma \ref{lem:loc_cent_dif}, we have with overwhelming probability that $|k-i| \le T (\log N)^C$ whenever $|\varphi_t (\hat{k}) - \varphi_t (\hat{i})| \le \mathfrak{M} (\log N)^5$), the left side of \eqref{sumalphaq} is at most equal to 
\begin{flalign*}
\displaystyle\sum_{i=N_1}^{N_2} \mathbbm{1}_{|\varphi_t (\hat{k})- \varphi_t (\hat{i})|\le \mathfrak{M} (\log N)^5} \cdot | & q_k(0) - q_i(0) - \alpha(k-i)| \cdot \displaystyle\sup_{x \in \mathbb{R}} | \chi'' (x)| \\
& \le 3\mathfrak{M} (\log N)^5 \cdot 100 \mathfrak{M}^{-2} \cdot T^{1/2} (\log N)^C \le T^{-c},
\end{flalign*}

\noindent with overwhelming probability. This establishes \eqref{sumalphaq} and thus the lemma. 
\end{proof}

\subsubsection{An approximate solution to \eqref{eqn:Zeq_1}}

\label{ProofZ00}

    In view of Lemma \ref{lem:Z_eq_outline}, we seek to solve the linear equation \eqref{eqn:Zeq_1} for $Z_k^{\mathcal Q}$. We will show that, recalling $Z(\Lambda,Q,t)$ from \eqref{eqn:Z_outline}, the function  $Z(\Lambda_{k}, \alpha k + t \ve(\Lambda_{k}), t)$ approximately solves this equation. Indeed, the following proposition states that this holds more generally (as mentioned in Section \ref{ProofFluctuationZ}), upon replacing the parameters $(\Lambda_k, \alpha k + t \ve(\Lambda_{k}))$ by $(\Lambda, Q)$. Observe in particular below that two parameters $(\Lambda, Q)$ are arbitrary, while $(\Lambda_k, \alpha k + t \ve(\Lambda_{k}))$ are both implicitly governed by the single parameter $k$.

\begin{prop}[$Z$ solves \eqref{eqn:Zeq_1}]\label{prop:sumZlambdaQ}

  Adopt Assumption \ref{ass:NT_assumption}, and fix $t \in [0, T \log N]$. There exist constants $\theta_0 (\beta) > 0$ and $\mathfrak{c} = \mathfrak{c}(\theta, \beta)>0$ such that the following holds with overwhelming probability whenever $\theta < \theta_0 (\beta)$. For all $\Lambda \in [- \log N, \log N]$ and $Q \in [-T^{6/5}, T^{6/5}]$, the function $Z(\Lambda, Q) = Z(\Lambda, Q, t)$ (from \eqref{eqn:Z_outline}) satisfies
\begin{flalign}\label{eqn:Zeq_2_prime}
\begin{aligned} 
   \Bigg|  2 \sum_{i=N_1}^{N_2}\big (Z(\Lambda, Q)-Z(\Lambda_{i}, \alpha i + t  \ve(\Lambda_i)) & \big) \cdot \mathfrak{l}(\Lambda -\Lambda_i) \cdot \chi' \big(Q- \alpha i - t  \ve(\Lambda_{i})  \big)   \\
    & + Z(\Lambda, Q)  - \Xi \big(\Lambda, (Q- t \ve(\Lambda))/\alpha,t \big)  \Bigg| \leq T^{-\mathfrak{c}}.
    \end{aligned} 
\end{flalign}
\end{prop}

If $(\Lambda, Q) = (\Lambda_k, \alpha k + t \ve (\Lambda_k))$ then $k = (Q - t \ve (\Lambda_k))/\alpha$, so Proposition \ref{prop:sumZlambdaQ} implies that $Z(\Lambda_{k}, \alpha k + t \ve(\Lambda_{k}), t)$ approximately solves \eqref{eqn:Zeq_1} for $Z_k^{\mathcal Q}$. To show $Z_k^{\mathcal Q} \approx Z(\Lambda_{k}, \alpha k + t \ve(\Lambda_{k}), t)$, it thus suffices to verify that the linear equation \eqref{eqn:Zeq_1} has a unique solution (up to a small error). The latter is essentially stated as Lemma \ref{lem:S_inv}, which we use to prove Theorem \ref{prop:dressing_Z_approx_outline} in Section \ref{sec:quasi_conc_estimates}. 

Proposition \ref{prop:sumZlambdaQ} will follow from the following bound proven in Section \ref{sec:Zsumtoint}.

\begin{prop}\label{prop:sumZconc}
  Adopt Assumption \ref{ass:NT_assumption}, and fix $t \in [0, T \log N]$. There exist constants $\theta_0 (\beta) > 0$ and $\mathfrak{c}= \mathfrak{c}(\theta, \beta)>0$ such that the following holds with overwhelming probability whenever $\theta < \theta_0 (\beta)$. For all $\Lambda \in [- \log N, \log N]$, $Q \in [- T^{6/5}, T^{6/5}]$, the function $Z(\Lambda, Q) = Z(\Lambda, Q, t)$ (from \eqref{eqn:Z_outline}) satisfies
\begin{multline}\label{eqn:Zeq_2_prime_outline}
   \Bigg| \sum_{i=N_1}^{N_2} \big( Z(\Lambda, Q)- Z(\Lambda_{i}, \alpha i + t  \ve(\Lambda_i)) \big) \cdot   
    \mathfrak{l}(\Lambda -\Lambda_i) \cdot \chi'(Q -\alpha i - t \ve(\Lambda_{i}) )  \\
    -\alpha^{-1} \int_{-\infty}^{\infty}\mathfrak{l}(\Lambda - \lambda) \cdot (Z(\Lambda, Q)-Z(\lambda,Q)) \cdot \varrho(\lambda) d\lambda \Bigg| \leq   T^{-\mathfrak{c}}
\end{multline}
\end{prop}

Observe that the integral on the left side of \eqref{eqn:Zeq_2_prime_outline} is obtained from the sum by two modifications. The first is to replace the sum over the $(\Lambda_i)$ by an integral with respect to the measure $\varrho (\lambda) d \lambda$. The second is to omit the function $\chi' (Q-\alpha i - \ve(\Lambda_{i}))$ and replace $\alpha i + \ve(\Lambda_{i})$ by $Q$. Heuristically, the former is plausible since $\varrho$ is the density of the eigenvalues $(\Lambda_i)$ of $\mathbf{L}(0)$ (see Remark \ref{lrho}), and the latter is plausible since $\chi'$ behaves as the $\delta$ distribution on scales larger than $\mathfrak{M} \ll T$. While concentration bounds of a seemingly similar nature were shown as \cite[Proposition 4.3]{Agg25}, the situation here is made more intricate by the fact that $Z$ depends on all eigenvalues of $\mathbf{L}(0)$ (as opposed to on only a single one).

 \begin{proof}[Proof of Proposition \ref{prop:sumZlambdaQ}]

 Throughout this proof, we restrict to the event on which Proposition \ref{prop:sumZconc} and the event $\mathsf{BND}_{\mathbf{L}} (\log N)$ (recall Definition \ref{eigenvaluesm}) both hold; this is overwhelmingly probable by Lemma \ref{lem:bd_lem}. We also let $\Lambda \in [-\log N, \log N]$, $Q \in [-T^{6/5}, T^{6/5}]$ be any real numbers; $t \in [0, T \log N]$ be fixed as in the proposition; and abbreviate $Z(\Lambda,Q) = Z(\Lambda,Q,t)$. By \eqref{eqn:Zeq_2_prime_outline} and the definition \eqref{operatort} of $\mathbf{T}$, we have 
 \begin{flalign}
     \label{sumz01}
     \begin{aligned}
    \Bigg|  2 \sum_{i=N_1}^{N_2}\big (Z( & \Lambda, Q)-Z(\Lambda_{i}, \alpha i + t  \ve( \Lambda_i)) \big) \cdot \mathfrak{l}(\Lambda -\Lambda_i) \cdot \chi'(Q -\alpha i - t \ve(\Lambda_{i}) ) \\
    & \qquad \qquad \qquad \qquad \qquad + Z(\Lambda, Q)  - \big( 1 + \alpha^{-1} (\mathbf{T} \varrho (\Lambda) - \mathbf{T} \bm{\varrho}) \big) Z(\Lambda, Q) \Bigg| \\
    & \leq T^{-c} + C \displaystyle\int_{-\infty}^{\infty} \big| \mathfrak{l} (\Lambda - \lambda) - \log |\Lambda - \lambda| \big| \cdot | Z(\Lambda,Q) - Z(\lambda,Q)| \cdot \varrho (\lambda) d \lambda ,
    \end{aligned}
 \end{flalign}

 \noindent where all operators act in the first argument $\Lambda$ of $Z$. By \eqref{eqn:sl}, we have $|\log |x| - \mathfrak{l}(x)| \le e^{-(\log N)^2}$ if $|x| \ge e^{-(\log N)^2}$ and $|\log |x| - \mathfrak{l}(x)| \le 2|\log |x||$ otherwise. Hence,
 \begin{flalign}
     \label{zsum00}
     \begin{aligned}
    \displaystyle\int_{-\infty}^{\infty} \big| \mathfrak{l} & (\Lambda - \lambda) - \log |\Lambda - \lambda| \big| \cdot | Z(\Lambda,Q) - Z(\lambda,Q)| \cdot \varrho (\lambda) d \lambda \\
    & \le 2\displaystyle\int_{-\infty}^{\infty} |Z(\Lambda,Q)-Z(\lambda,Q)| \cdot (|\log |\Lambda-\lambda|| \cdot \mathbbm{1}_{|\Lambda-\lambda| \le e^{-(\log N)^2}} + e^{-(\log N)^2}) \cdot \varrho (\lambda) d \lambda.
    \end{aligned}
 \end{flalign}

\noindent Moreover, from Lemma \ref{lem:ve_tail} and \eqref{eqn:Z_outline}, we have  
 \begin{flalign}
     \label{zxi}
    |Z(\lambda, Q, t)| \le C \big| \Xi \big(\lambda, (Q-t\ve(\lambda))/\alpha, t \big) \big| + C \big\|  \Xi \big( \cdot , (Q-t\ve(\lambda))/\alpha, t \big) \big\|_{\mathcal{H}} \cdot \log (|\lambda|+2).
 \end{flalign}

 \noindent Also, by Definition \ref{def:smoothed_log}, we have $|\chi| \le 1$ and $|\mathfrak{l}(\lambda - \Lambda_i)| \le (\log N)^3 \cdot \log (|\lambda|+2)$ for all $i \in \llbracket N_1, N_2 \rrbracket$ (as $|\Lambda_i| \le \log N$ for all $i$, by our restriction to $\mathsf{BND}_{\mathbf{L}} (\log N)$). These, with \eqref{eqn:psiZktau}, imply that $|\Xi(\lambda,k,t)| \le N \log (|\lambda|+2)$, so that $\| \Xi (\cdot, k, t) \|_{\mathcal{H}} \le CN$, for all $\lambda, k, t \in \mathbb{R}$. Together with \eqref{zxi}, this gives $|Z(\lambda,Q,t)| \le CN \log (|\lambda|+2)$ for all $\lambda, Q, t \in \mathbb{R}$. Inserting this into \eqref{zsum00}, and using the decay of $\varrho$ provided by Lemma \ref{lem:varrho_bd}, yields
 \begin{flalign}
    \displaystyle\int_{-\infty}^{\infty} \big| \mathfrak{l} & (\Lambda - \lambda) - \log |\Lambda - \lambda| \big| \cdot | Z(\Lambda,Q) - Z(\lambda,Q)| \cdot \varrho (\lambda) d \lambda \le e^{-c(\log N)^2},
 \end{flalign}

 \noindent which, together with \eqref{sumz01}, implies
 \begin{flalign}
     \label{sumz0}
     \begin{aligned}
    \Bigg|  2 \sum_{i=N_1}^{N_2}\big (Z( & \Lambda, Q)-Z(\Lambda_{i}, \alpha i + t  \ve( \Lambda_i)) \big) \cdot \mathfrak{l}(\Lambda -\Lambda_i) \cdot \chi'(Q -\alpha i - t \ve(\Lambda_{i}) )   \\
    & \qquad \qquad \qquad \qquad \qquad + Z(\Lambda, Q)  - \big( 1 + \alpha^{-1} (\mathbf{T} \varrho (\Lambda) - \mathbf{T} \bm{\varrho}) \big) Z(\Lambda, Q) \Bigg| \leq T^{-c}.
    \end{aligned}
 \end{flalign}

  Now observe that 
\begin{flalign}
\label{eqn:operator_identity}
\begin{aligned} 
 \big( 1 + \alpha^{-1}  (\mathbf{T}\varrho (\Lambda) - \mathbf{T} \bm{\varrho} ) \big) Z(\Lambda,Q) & = (\boldsymbol{\varsigma_0^{\dr}}  - \alpha^{-1} \mathbf{T} \bm{\varrho}) Z(\Lambda, Q) \\
 & =  (1 - \theta\mathbf{T}   \boldsymbol{\varrho_{\beta}})  \boldsymbol{\varsigma_0^{\dr}} Z(\Lambda, Q) 
 =  \Xi \big(\Lambda, (Q - t\ve(\Lambda))/\alpha \big),
 \end{aligned} 
 \end{flalign}

 \noindent where the first and second statements follows from the second and first equalities in \eqref{rho2}, respectively, and the third follows from the definition \eqref{eqn:Z_outline} of $Z(\Lambda, Q)$ (and that \eqref{drf} for the dressing operator). Combining \eqref{sumz0} and \eqref{eqn:operator_identity} yields the proposition.
 \end{proof}

\subsection{Proof of Theorem \ref{prop:zcouple_outline}}

\label{ProofZ} 

\subsubsection{Coupling $\Xi$ to $\mathfrak{X}$} 

\label{CoupleXiX} 

Recall that $Z$ and $\mathcal{Z}$ (from \eqref{eqn:Z_outline} and \eqref{eqn:intro_Z}) are obtained, up to a factor of $\varsigma_0^{\dr} (\Lambda)^{-1}$, by applying the dressing operator to the processes $\Xi$ and $\mathfrak{X}$. Thus, to show $Z \approx \mathcal{Z}$, we will show that $\Xi \approx \mathfrak{X}$ (recall \eqref{eqn:Xdef}). This is done by the below theorem, shown in Section \ref{ProofXLambda}.

\begin{thm}[Limit for $\Xi$]\label{thm:couplethm_outline}

Adopt Assumption \ref{ass:NT_assumption}. There are constants $\theta_0 (\beta) > 0$ and $\mathfrak{c}>0$ such that the following holds whenever $\theta < \theta_0 (\beta)$. There is a coupling between $\mathfrak{X}$ and $\Xi$ such that, with overwhelming probability, we have
\begin{equation}\label{eqn:xi_couple_outline} 
\displaystyle\sup_{\tau \in [0, \log N]} 
\displaystyle\sup_{|\Lambda| \le \log N} \displaystyle\sup_{|\kappa| \le (\log N)^4}\left| \Xi(\Lambda, T \kappa, T \tau) - \mathfrak{X}(\Lambda, \kappa, \tau)   \right| \leq T^{-\mathfrak{c}} .
\end{equation}
\end{thm}

Theorem \ref{prop:zcouple_outline} will follow quickly from Theorem \ref{thm:couplethm_outline}. The proof is facilitated by the following estimate on the processes $\Xi$ and $\mathfrak{X}$, shown in Appendix \ref{app:xi_bd_outline}.

\begin{lem}\label{lem:xi_bd_outline}

Adopt Assumption \ref{ass:NT_assumption}. There exist constants $\theta_0 (\beta) > 0$ and $\mathfrak{C}>1$ such that the following holds with overwhelming probability whenever $\theta < \theta_0 (\beta)$. For any real numbers $\Lambda \in  \mathbb{R}$ and $t \in [0,  T \log N]$, we have
\begin{equation}\label{eqn:small_Lambda_bd_outline}
     \displaystyle\max_{k \in \llbracket N_1, N_2 \rrbracket}  | \Xi(\Lambda, k, t) | \leq (\log N)^{\mathfrak{C}} T^{-1/2} (t + \mathfrak{M})^{1/2} (|\Lambda|+1)^2,
    \end{equation}
    
    \noindent and for any $\Lambda \in  \mathbb{R}$,
    \begin{equation}\label{eqn:mathfrakX_bd_outline}
    \displaystyle\sup_{\tau \in [0, \log N]} \displaystyle\sup_{|\mathfrak{q}| \le (\log N)^4} |\mathfrak{X}(\Lambda, (\mathfrak{q} - \tau \ve(\Lambda)) / \alpha, \tau)| \leq (\log N)^{\mathfrak{C}} (|\Lambda|+1)^2.
    \end{equation}
\end{lem}

\begin{proof}[Proof of Theorem \ref{prop:zcouple_outline}]

    Restrict to the event on which Theorem \ref{thm:couplethm_outline} and Lemma \ref{lem:xi_bd_outline} hold. Let $\tau \in [0, \log N]$, $\Lambda \in [-\log N, \log N]$, and $\mathfrak{q} \in [-(\log N)^3, (\log N)^3]$. Setting $t = \tau T$ and $Q = \mathfrak{q} T$, let
    \begin{flalign*} 
    f(\Lambda) = \mathfrak{X}\big(\Lambda, (\mathfrak{q} - \tau \ve(\Lambda)) / \alpha, \tau \big) -   \Xi \big(\Lambda, (Q - t \ve(\Lambda)) / \alpha, t \big).
    \end{flalign*} 

    \noindent Then the definitions \eqref{eqn:intro_Z} and \eqref{eqn:Z_outline} of $\mathcal{Z}$ and $\Xi$, respectively, yield
    \begin{flalign*} 
    \mathcal{Z} (\Lambda, \mathfrak{q}, \tau) - Z(\Lambda, T \mathfrak{q}, T \tau) = \varsigma_0^{\dr} (\Lambda)^{-1} \cdot 
        f^{\dr} (\Lambda). 
    \end{flalign*}
    
    \noindent Together with the first and third statements of Lemma \ref{lem:ve_tail}, this implies that 
    \begin{equation}\label{eqn:fdr_g_rel_outline}
         |\mathcal{Z} (\Lambda, \mathfrak{q}, \tau) - Z(\Lambda, T \mathfrak{q}, T \tau)| \leq C \big( |f(\Lambda)| +  \|f\|_{\mathcal{H}} \cdot \log(|\Lambda|+2) \big).
    \end{equation}

    Since $|\mathfrak{q} - \tau \ve(\Lambda)| / |\alpha| \le C(\log N)^3$ (as $|\mathfrak{q}| \le (\log N)^3$,  $\tau \le \log N$, and $|\ve(\Lambda)| \le C \log N$ by the second part of Lemma \ref{lem:ve_tail}, since $|\Lambda| \le \log N$), \eqref{eqn:xi_couple_outline} applies to yield $|f(\Lambda)| \le T^{-c}$. Hence, it suffices to show that $\| f \|_{\mathcal{H}} \le T^{-c}$, to which end observe that
    \begin{multline}
    \label{frho00}
    \| f \|_{\mathcal{H}} \leq  \Bigg| \int_{-\log N}^{\log N} f(y)^2 \varrho(y) dy \Bigg| + \Bigg| \int_{|x| \ge \log N} f(y)^2 \varrho(y) dy \Bigg| \leq T^{-c} + (\log N)^C e^{-c (\log N)^2},
    \end{multline}
    \noindent where in the second inequality we bounded the first integral using \eqref{eqn:xi_couple_outline}, and we bounded the second using \eqref{eqn:small_Lambda_bd_outline}, \eqref{eqn:mathfrakX_bd_outline}, and the fact from Lemma \ref{lem:varrho_bd} that $\varrho (x) \le C(|x| + 1)^{2\theta+1} e^{-\beta x^2/2}$. 
\end{proof}

\subsubsection{Proof of Theorem \ref{thm:couplethm_outline}} 

\label{ProofXLambda} 

In this section we show Theorem \ref{thm:couplethm_outline} as a direct consequence of the more general coupling statement given by Theorem \ref{thm:couplethm12} below. We begin with the following definition of the process $\Xi_1$, obtained by replacing the quantities $q_i(0)$ in the definition \eqref{eqn:psiZktau} of $\Xi$ with $\alpha i$. Although \eqref{qiqjs4} suggests this is reasonable, it will incur a non-negligible error that must be tracked (see \eqref{eqn:xibr_couple} below).

\begin{definition} 

\label{xi2} 
Adopt Assumption \ref{ass:NT_assumption}. Recalling Definitions \ref{def:quasi_intro} and \ref{def:chi_def}, for any real numbers $k, \Lambda \in \mathbb{R}$ and $t \in \mathbb{R}_{\ge 0}$, denote 
\begin{flalign}\label{eqn:Xi2_outline}
\begin{aligned} 
    \Xi_1(\Lambda, k, t) & = T^{-1/2} \cdot \Bigg( 2 \sum_{i = N_1}^{N_2} \mathfrak{l}(\Lambda- \Lambda_i) \cdot 
    \Big( \chi(\alpha (k - i)) -  
   \chi \big( \alpha (k- i) + t   \left(\ve(\Lambda) - \ve(\Lambda_i) \right) \big) \Big)  \\
  & \quad - \mathbb{E} \bigg[2 \sum_{i = N_1}^{N_2} \mathfrak{l}(\Lambda- \Lambda_i) \cdot 
    \Big( \chi(\alpha (k - i)) -  
   \chi \big( \alpha (k- i) + t   \left(\ve(\Lambda) - \ve(\Lambda_i) \right) \big) \Big) \bigg] \Bigg).
   \end{aligned} 
\end{flalign}

\end{definition} 

We next introduce the following counterparts $\Xi^{[m]}$ and $\Xi_1^{[m]}$ of $\Xi$ and $\Xi_1$ (from Definitions \ref{xilambdakt} and \ref{xi2}), respectively. They involve two free spatial parameters $(q,q')$ (instead of the $q_k(0)$ and $q_k (0) + t \ve(\Lambda)$ in $\Xi$) and multiplication by the power $\lambda_j^m$ (instead of the approximate logarithm $\mathfrak{l}(\Lambda-\Lambda_i)$ in $\Xi$ and $\Xi_1$). For the sole purpose of proving Theorem \ref{thm:couplethm_outline}, these processes are not necessary and can be disregarded. However, we mention them, since to analyze current fluctuations of the Toda lattice later in Section \ref{sec:current_fluct}, it will be necessary to approximate $(\Xi, \Xi_1, \Xi^{[m]}, \Xi_1^{[m]})$ simultaneously, as stated in Theorem \ref{thm:couplethm12} below.

\begin{definition} 

\label{xim} 

Adopt Assumption \ref{ass:NT_assumption}. Recalling Definitions \ref{def:quasi_intro} and \ref{def:chi_def}, for any real numbers $q,q' \in \mathbb{R}$ and $t \in \mathbb{R}_{\ge 0}$, denote
\begin{flalign}\label{eqn:xim_outline}
\begin{aligned} 
\Xi^{[m]}(q,q', t)&  \coloneqq T^{-1/2} \cdot \Bigg( \sum_{j=1}^{N} \big( \lambda_j^m \cdot  \chi(q-Q_j(0)) - \lambda_j^m \cdot \chi(q' -Q_j(0)- t \ve(\lambda_j)) \big)  \\
 & \qquad \qquad - \mathbb{E}\bigg[\sum_{j=1}^{N} \big(\lambda_j^m \cdot \chi(q-Q_j(0)) - \lambda_j^m \cdot \chi(q'-Q_j(0)- t \ve(\lambda_j)) \big) \bigg] \Bigg),
 \end{aligned}
\end{flalign}

and
\begin{flalign}\label{eqn:Xim2_outline}
\begin{aligned} 
    \Xi_1^{[m]}(q,q', t) & = T^{-1/2} \cdot \Bigg( \sum_{i=N_1}^{N_2}   \Lambda_i^m \cdot \big( \chi( q- \alpha i) - \chi(q'- \alpha i - t \ve(\Lambda_i)) \big) \\
  & \qquad \qquad \qquad -\mathbb{E} \bigg[\sum_{i=N_1}^{N_2}   \Lambda_i^m \cdot \big( \chi (q- \alpha i) - \chi(q'- \alpha i  - t \ve(\Lambda_i)) \big) \bigg] \Bigg).
\end{aligned} 
\end{flalign}

\end{definition} 

We next define several Gaussian processes $\mathfrak{X}_0$, $\mathfrak{X}_1$, $\mathfrak{Y}$, $\mathfrak{Y}_0$, and $\mathfrak{Y}_1$ that, as indicated by Theorem \ref{thm:couplethm12}, denote the scaling limits of $\Xi - \Xi_1$, $\Xi_1$, $\Xi^{[m]}$, $\Xi^{[m]} - \Xi_1^{[m]}$, and $\Xi_1^{[m]}$, respectively. The definitions of $\mathfrak{Y}$, $\mathfrak{Y}_0$, and $\mathfrak{Y}_1$ necessitate notation for a new family of test functions. For any $(\mathfrak{q}, \mathfrak{q}', \tau, m) \in \mathbb{R} \times \mathbb{R} \times \mathbb{R}_{\geq 0} \times \mathbb{Z}_{\geq 0}$, define the function $\phi_{\mathfrak{q}, \mathfrak{q}',\tau}^{[m]}: \mathbb{R}^2 \rightarrow \mathbb{R}$ by
\begin{equation}\label{eqn:phim_intro}
\phi_{\mathfrak{q}, \mathfrak{q}',\tau}^{[m]}(r, \lambda) \coloneqq \lambda^m \cdot \big(\mathbbm{1}\{\mathfrak{q} > \alpha r \} - \mathbbm{1}\{\mathfrak{q}' > \alpha r +  \tau \ve(\lambda) \} \big).
\end{equation}
Now, we may define the Gaussian processes.

\begin{definition} 

\label{x2x0y2y0} 
Recalling Definitions \ref{def:Wdress} and \ref{def:Hdef}, define the (Gaussian) processes $\mathfrak{X}_0, \mathfrak{X}_1 : \mathbb{R}^2 \times \mathbb{R}_{\ge 0} \times \mathbb{R}_{\ge 0} \rightarrow \mathbb{R}$ and $\mathfrak{Y}_0, \mathfrak{Y}_1 : \mathbb{R}^2 \times \mathbb{R}_{\ge 0} \times \mathbb{Z}_{\ge 0} \rightarrow \mathbb{R}$ by setting
\begin{align}
   & \mathfrak{X}_0(\Lambda, \kappa, \tau) =  \mathcal{W}^{\dr} \big(\langle \psi_{\Lambda, \kappa,\tau}(r,\cdot), \varsigma_0  \rangle_{\varrho} \cdot \varsigma_0 \big); \label{eqn:X2out} \\
   & \mathfrak{X}_1(\Lambda, \kappa, \tau)  =\mathcal{W}^{\dr} \big(\psi_{\Lambda, \kappa,\tau} - \langle \psi_{\Lambda, \kappa,\tau}(r,\cdot), \varsigma_0  \rangle_{\varrho} \cdot \varsigma_0 \big);  \label{eqn:X1out} \\
   &  \mathfrak{Y}_0(\mathfrak{q}, \mathfrak{q}', \tau, m) =  \mathcal{W}^{\dr}(\langle \phi_{\mathfrak{q}, \mathfrak{q}',\tau}^{[m]}(r, \cdot), \varsigma_0  \rangle_{\varrho} \cdot \varsigma_0); \label{eqn:Y2out} \\
    & \mathfrak{Y}_1(\mathfrak{q}, \mathfrak{q}', \tau, m) =\mathcal{W}^{\dr}(\phi_{\mathfrak{q}, \mathfrak{q}',\tau}^{[m]} - \langle \phi_{\mathfrak{q},\mathfrak{q}',\tau}^{[m]}(r, \cdot), \varsigma_0  \rangle_{\varrho} \cdot \varsigma_0),	\label{eqn:Y1out}  
\end{align} 

\noindent for any $(\Lambda, \kappa, \tau) \in \mathbb{R}^2 \times \mathbb{R}_{\geq 0}$ and $(\mathfrak{q},\mathfrak{q}', \tau, m) \in \mathbb{R}^2 \times \mathbb{R}_{\geq 0} \times \mathbb{Z}_{\geq 0}$. Moreover, define the process $\mathfrak{Y} : \mathbb{R}^2 \times \mathbb{R}_{\ge 0} \times \mathbb{Z}_{\ge 0} \rightarrow \mathbb{R}$ by setting
\begin{equation}
\label{yprocess}
\mathfrak{Y}(\mathfrak{q}, \mathfrak{q}', \tau, m) = \mathcal{W}^{\dr}(\phi_{\mathfrak{q}, \mathfrak{q}', \tau}^{[m]}),
\end{equation}

\noindent for any $(\mathfrak{q}, \mathfrak{q}', \tau, m) \in \mathbb{R} \times \mathbb{R} \times \mathbb{R}_{\geq 0} \times \mathbb{Z}_{\geq 0}$.
 \end{definition} 
 \begin{remark}\label{rmk:Yfin}
 Similarly to Lemma \ref{lem:psi_phi_L2_intro}, for any $(\mathfrak{q}, \mathfrak{q}', \tau, m) \in \mathbb{R} \times \mathbb{R} \times \mathbb{R}_{\geq 0} \times \mathbb{Z}_{\geq 0}$, it is quickly verified that $\phi_{\mathfrak{q}, \mathfrak{q}', \tau}^{[m]} \in L^2(\mathbb{R}^2, d r\otimes \varrho d \lambda)$ (see Lemma \ref{lem:bscL2bds} below). Therefore, the random variable $\mathfrak{Y}(\mathfrak{q}, \mathfrak{q}', \tau, m)$ is almost surely finite (as are $\mathfrak{Y}_0(\mathfrak{q}, \mathfrak{q}', \tau, m) $ and $\mathfrak{Y}_1(\mathfrak{q}, \mathfrak{q}', \tau, m)$). Moreover, note that we have the trivial equalities
$
 \mathfrak{Y}(\mathfrak{q}, \mathfrak{q}', \tau, m) = \mathfrak{Y}_0(\mathfrak{q}, \mathfrak{q}', \tau, m) + \mathfrak{Y}_1(\mathfrak{q}, \mathfrak{q}', \tau, m) 
$ and  (recalling $\mathfrak{X}$ from Definition \ref{def:Z_intro}) $
 \mathfrak{X}(\Lambda, \kappa, \tau)  = \mathfrak{X}_0(\Lambda, \kappa, \tau)  +  \mathfrak{X}_1(\Lambda, \kappa, \tau)  .
$
 \end{remark}

The following theorem, to be shown in Section \ref{subsec:main_coupling_arg}, shows that $(\Xi_1, \Xi - \Xi_1, \Xi_1^{[m]}, \Xi^{[m]}-\Xi_1^{[m]})$ approximate $(\mathfrak{X}_1, \mathfrak{X}_0, \mathfrak{Y}_1, \mathfrak{Y}_0)$.

\begin{thm}[Limits for $(\Xi,\Xi_1,\Xi^{[m]}, \Xi_1^{[m]})$]\label{thm:couplethm12}

Adopt Assumption \ref{ass:NT_assumption}, and recall Definitions \ref{xilambdakt}, \ref{xi2}, \ref{xim}, and \ref{x2x0y2y0}. There is a constant $\theta_0 (\beta) > 0$ such that, whenever $\theta < \theta_0(\beta)$, there exists a constant $\mathfrak{c}>0$ and a coupling of the tuple of processes $(\mathfrak{X}_0,\mathfrak{X}_1,\mathfrak{Y}_1, \mathfrak{Y}_0)$ with $\L(0)$ such that the following holds with overwhelming probability. For any real numbers $\tau \in [0, \log N]$ and $\Lambda \in [-\log N, \log N]$, we have that 
\begin{flalign}\label{eqn:xi2_couple}
& \displaystyle\sup_{|\kappa| \le (\log N)^4} \left| \Xi_1(\Lambda, T \kappa, T \tau) - \mathfrak{X}_1(\Lambda, \kappa, \tau)   \right| \leq T^{-\mathfrak{c}}; \\ 
\label{eqn:xibr_couple}
& \displaystyle\sup_{|\kappa| \le (\log N)^4} \left| \Xi(\Lambda, T \kappa, T \tau)-\Xi_1(\Lambda, T \kappa, T \tau) - \mathfrak{X}_0(\Lambda, \kappa, \tau)   \right| \leq T^{-\mathfrak{c}}; \\
\label{eqn:xim2_couple}
& \displaystyle\max_{m \in \llbracket 0, (\log N)^{1/10} \rrbracket}  \displaystyle\sup_{|\mathfrak{q}|,|\mathfrak{q}'| \le (\log N)^4} \left| \Xi_1^{[m]}(T \mathfrak{q}, T \tau) - \mathfrak{Y}_1(\mathfrak{q},\mathfrak{q}', \tau, m)  \right| \leq T^{-\mathfrak{c}}; \\
\label{eqn:ximbr_couple}
& \displaystyle\max_{m \in \llbracket 0, (\log N)^{1/10} \rrbracket} \displaystyle\sup_{|\mathfrak{q}|,|\mathfrak{q}'| \le (\log N)^4} \Big| \Xi^{[m]}(T \mathfrak{q}, T \mathfrak{q}',  T \tau  )-\Xi_1^{[m]}(T \mathfrak{q}, T \mathfrak{q}', T \tau) - \mathfrak{Y}_0(\mathfrak{q},\mathfrak{q}', \tau, m)  \Big| \leq T^{-\mathfrak{c}}.
\end{flalign}

\end{thm}

\begin{proof}[Proof of Theorem \ref{thm:couplethm_outline}]

This follows from \eqref{eqn:xi2_couple} and \eqref{eqn:xibr_couple}.
\end{proof} 

While we will not establish Theorem \ref{thm:couplethm12} here, let us briefly comment on its proof. We restrict the discussion below to $(\Xi, \Xi_1)$, as similar reasoning will apply to $(\Xi^{[m]}, \Xi_1^{[m]})$. 

1. First, we show that the processes $\Xi_1$ and $\Xi-\Xi_1$ can be approximated by sums of independent random variables; see Proposition \ref{prop:equiv_expr_cor}. To heuristically see this for the difference $\Xi - \Xi_1$, recall that $\Xi$ and $\Xi_1$ differ in that the arguments of $\chi$ involve $q_i (0)$ in the definition \eqref{eqn:psiZktau} of $\Xi$ but involve $\alpha i$ in that \eqref{eqn:Xi2_outline} of $\Xi_1$. Thus, Taylor expanding $\Xi - \Xi_1$ to first order yields a linear combination of the $(q_i (0) - \alpha i)$, with random coefficients (dependent on the eigenvalues $(\Lambda_i)$) that should average to deterministic ones. Since the $(q_i - q_{i-1})$ are independent under thermal equilibrium (by \eqref{abr} and Definition \ref{eqn:equil}), this linear combination is expressible as a sum of independent random variables. 

For $\Xi_1$, the underlying heuristic is different and will be based on the notion that the eigenvalues $\Lambda_i$ of $\mathbf{L}(0)$ are \emph{approximately local}; this means that they essentially only depend on the entries of $\mathbf{L}(0)$ near the $(i,i)$ one (see Lemma \ref{lem:Lax_eig_coupling} for a more precise statement). By the independence of the entries in $\mathbf{L}(0)$, this indicates that $\Lambda_i$ and $\Lambda_j$ are almost independent, unless $|i-j|$ is small. Hence, if one subdivides the sum over $i$ in the definition \eqref{eqn:Xi2_outline} of $\Xi_1$ into small intervals (say, of length $T^c$), the resulting summands should nearly be independent, providing the approximation for $\Xi_1$ as a sum of independent random variables. 

2. Second, we identify the limiting covariances of certain linear statistics of $\mathbf{L}(0)$, with effective  rates of convergence; see Proposition \ref{prop:limvar}. For polynomial test functions, it is known \cite{MM24} that these covariances tend to explicit limits, but without a rate of convergence. In our context, we must allow the test functions to involve the smoothed logarithm $\mathfrak{l}$, and we must ensure that errors are polynomial in $T$ (to offset the factors of $(\log N)^C$ that frequently arise in errors, such as in Lemma \ref{lem:asymptotic-scattering}). This is done through a multi-scale argument, using the above-mentioned approximate locality of eigenvalues of $\mathbf{L}(0)$ and estimates from polynomial approximation theory. 

3. Third, we use the above two statements to couple $(\Xi_1, \Xi - \Xi_1)$ to a Gaussian process. When restricting the parameters $(\Lambda, \kappa, \tau)$ to a mesh of spacing $T^{-c}$, this follows from known high-dimensional Gaussian approximation results \cite{BM06} (see also \cite{EDMCL}). To extend these parameters to their full domains, we require a H\"{o}lder-type continuity estimate of $\Xi$; this is shown as Lemma \ref{lem:xi_xim_cont}, using concentration estimates for functionals of random matrices.

\subsection{Outline}

The remainder of this paper is organized as follows. In Section \ref{sec:bg} we collect preliminary identities and estimates that will be used throughout. In Section \ref{sec:ind_linear} we show how certain functionals (involving both eigenvalues and localization centers) of the random Lax matrix $\mathbf{L}(0)$ can be approximated by sums of independent variables. In Section \ref{sec:effective_conv} we show effective convergence results for the limiting covariances of linear statistics of $\mathbf{L}(0)$. In Section \ref{sec:continuity} we show continuity bounds on $\Xi$ and $\Xi^{[m]}$, which we use to establish the Gaussian coupling result Theorem \ref{thm:couplethm12} in Section 
\ref{sec:levy_chent}. In Section \ref{sec:quasi_coupling} we show Proposition \ref{prop:sumZconc}, which we use in Section \ref{sec:quasi_conc_estimates} to prove Theorem \ref{prop:dressing_Z_approx_outline}. In Section \ref{sec:current_fluct} we apply these results to compute fluctuations of integrated currents, proving Theorem \ref{thm:current_fluct_intro}. Finally, Section \ref{sec:charge_fluct} computes fluctuations of local charges and uses them to analyze associated two-point functions, proving Corollary \ref{cor:twopoint_intro}. The proofs of several auxiliary results are deferred to the appendices.

\section{Miscellaneous preliminaries}
\label{sec:bg}

In this section we collect some preliminary facts that will be used throughout the paper. The statements here are either known or consequences of known results, so we defer their proofs (when we do not have a reference for them) to the Appendices.

\subsection{Properties of limiting quantities} 
\label{subsec:Lax_DOS}

In this section we state the following three lemmas that approximate the effective velocity $\ve$ (Definition \ref{def:v_eff_intro}); estimate the H\"{o}lder continuity of the processes $\mathfrak{X}_0$, $\mathfrak{Y}_0$, $\mathfrak{X}_1$, and $\mathfrak{Y}_1$ (Definition \ref{x2x0y2y0}); and estimates the process $\mathfrak{X}$ (Definition \ref{def:Z_intro}). We prove Lemma \ref{lem:veff_inc} in Appendix  \ref{sec:veff_lemma} and Lemmas \ref{lem:holder_cont} and \ref{lem:Xbd1} in Appendix \ref{app:lim_holder}. 

\begin{lem}\label{lem:veff_inc}

There exist constants $\mathfrak{C}>1$ and $\theta_0 (\beta) > 0$ such that the following two statements hold whenever $\theta < \theta_0 (\beta)$. First, we have 
    \begin{equation}
        \ve(x) = \frac{x + \theta R(x)}{1 + \theta R_0(x)}
    \end{equation}

    \noindent for some real numbers $R(x), R_0 (x) \in \mathbb{R}$ satisfying
    \begin{equation}\label{eqn:veff}
       |R(x)| + |R_0(x)| \leq \mathfrak{C} \log(|x|+2) .
    \end{equation}
    Second, $\ve'(x) > 0$ (that is, $\ve$ is strictly increasing), and for any real number $A \ge 2$ we have 
    \begin{equation}\label{eqn:dveff}
    \ve'(x) \geq \mathfrak{C}^{-1} (1+\theta \log A)^{-1}, \qquad \text{for all $x \in [-A, A]$}.
    \end{equation}
\end{lem}

\begin{lem}\label{lem:holder_cont}

There exists a constant $\mathfrak{c}>0$ such that the following holds. Let $N \ge 1$ be a real number; set $A = \log N$, and let $m \in \llbracket 0, (\log N)^{1/10} \rrbracket$ be an integer. Then, for each $i \in \{ 0, 1 \}$,  
\begin{flalign*}
    & \mathbb{P} \bigg(\sup_{s, s' \in [-A, A]^2 \times [0, A]}|\mathfrak{X}_i(s) - \mathfrak{X}_i(s')| > (\log N)^5 \|s-s'\|_2^{2/5} \bigg) \leq \mathfrak{c}^{-1}  e^{-\mathfrak{c} (\log N)^2}; \\
    & \mathbb{P} \bigg(\sup_{s, s' \in [-A, A]^2 \times [0, A] } |\mathfrak{Y}_i(s, m) - \mathfrak{Y}_i(s', m)| > (\log N)^{m+5} \|s-s'\|_2^{2/5} \bigg) \leq \mathfrak{c}^{-1}  e^{-\mathfrak{c} (\log N)^2}.
    \end{flalign*}
\end{lem}

\begin{lem}\label{lem:Xbd1}

There exists a constant $\mathfrak{c}>0$ such that the following holds. Let $A_1 \geq A_2 \geq 1$ and $\tau \in [0, A_2]$ be real numbers. Then, for any real number $Q \in \mathbb{R}$, we have 
    \begin{equation}\label{eqn:Xnmbd}
       \mathbb{P} \left(  \sup_{|\Lambda| \leq A_1} \big|\mathfrak{X} \big(\Lambda, (Q-\tau \ve(\Lambda))/\alpha, \tau \big) \big| > A_1^2 A_2 \right) \leq \mathfrak{c}^{-1} e^{- \mathfrak{c} A_1^2}.
    \end{equation}
\end{lem}

\subsection{Random Lax matrix estimates} 

In this section we provide several estimates and limiting results for a Lax matrix $\mathbf{L}$ (Definition \ref{lt}), sampled under thermal equilibrium (Definition \ref{def:inf_thermal_eq}). The following definition prescribes the event on which the eigenvalues of a matrix $\mathbf{M}$ are spaced, and the lemma below states that this event is likely for a random Lax matrix.

\begin{definition}

\label{lambda0lambdaevent}

Let $\mathbf{M}$ denote an $N \times N$ symmetric random matrix. For any real number $\delta \ge 0$, define the event
\begin{equation}
    \mathsf{SEP}_{\mathbf{M}}(\delta) \coloneqq \Bigg\{ \min_{\substack{\lambda, \lambda' \in \eig \mathbf{M} \\ \lambda \neq \lambda'}} |\lambda - \lambda'| \geq \delta \Bigg\}.
\end{equation}
\end{definition}

\begin{lem}[{\cite[Lemma 3.18]{Agg25a}}]\label{lem:sep_lem}
Adopt Assumption \ref{ass:NT_assumption}. There exists a constant $\mathfrak{c}>0$ such that, for any $\delta>0$, the Lax matrix $\mathbf{L} = \mathbf{L}(0)$ satisfies
\begin{equation}\label{sep_lem}
\mathbb{P}(\mathsf{SEP}_{\L} (\delta)) \ge 1- \mathfrak{c}^{-1} (\delta N^3 + e^{-\mathfrak{c} N^2}).
\end{equation}
\end{lem}

We next provide a limiting result, essentially due to \cite{MM24,Spo20}, for the limiting variance of polynomial linear statistics of the random Lax matrix $\mathbf{L}$. The formula for this limiting variance requires the below bilinear functional on $\mathcal{H} \times \mathcal{H}$. In what follows, we recall the notation from Definition \ref{def:alpha_Laxdos} and Definition \ref{def:Hdef}.

\begin{definition}\label{def:Cdef}

 For any function $\phi \in \mathcal{H}$, define 
\begin{flalign*} 
\mu_{\phi} \coloneqq \int_{-\infty}^{\infty} \varrho(x) \phi(x) dx = \langle \phi, \varsigma_0 \rangle_{\varrho}.
\end{flalign*} 

\noindent Further define the bilinear functional $\mathscr{C}: \mathcal{H} \times \mathcal{H} \rightarrow \mathbb{C}$ by, for any two functions $\phi_1, \phi_2 \in \mathcal{H}$, setting 
    \begin{equation}\label{eqn:Cdef}
        \mathscr{C}(\phi_1, \phi_2) = \big\langle (1- \theta \T  \boldsymbol{\varrho}_{\beta})^{-1} \big(\phi_1 - (1+\alpha) \mu_{\phi_1} \varsigma_0 \big) , (1-\theta \T  \boldsymbol{\varrho}_{\beta})^{-1} \big(\phi_2 - (1+\alpha) \mu_{\phi_2} \varsigma_0\big) \big\rangle_{\varrho}. 
    \end{equation}

\noindent Additionally, for any $\phi \in \mathcal{H}$, let
\begin{equation}\label{eqn:sig2def}
    \sigma^2(\phi) \coloneqq \mathscr{C}(\phi-\mu_{\phi} \varsigma_0, \phi-\mu_{\phi} \varsigma_0).
\end{equation}
\end{definition}

\begin{lem}\label{lem:lim_var_formula}

Adopt Assumption \ref{ass:NT_assumption}. For any polynomial $p: \mathbb{R} \rightarrow \mathbb{R}$, we have (recalling \eqref{eqn:sig2def}) 
\begin{equation}\label{eqn:var_explicit1}
\displaystyle\lim_{N \rightarrow \infty} 
\frac{1}{N} \var \left(\sum_{i=1}^{N} p(\lambda_i) \right) = \sigma^2(p).
\end{equation}

\noindent Moreover, we have (recalling \eqref{eqn:Cdef})
\begin{equation}\label{eqn:s0p_var}
\displaystyle\lim_{N \rightarrow \infty}
    \frac{1}{N} \Cov\left(\sum_{i=1}^{N} p(\lambda_i), \sum_{i=N_1}^{N_2-1} r_i \right) = \mathscr{C}(\varsigma_0, p - \mu_p \varsigma_0),
\end{equation}

\noindent where we have recalled $r_i = r_i(0)$ for $i \in \llbracket N_1, N_2 -1 \rrbracket$ from \eqref{abr}. Finally, we have 
\begin{equation}\label{eqn:s0_var}
    \var( r_i ) = \mathscr{C}(\varsigma_0, \varsigma_0) .
\end{equation}
\end{lem}

As mentioned above, Lemma \ref{lem:lim_var_formula} is essentially a consequence of results from \cite{MM24} and \cite{Spo20}. Specifically, \cite{MM24} proved convergence of the variances \eqref{eqn:var_explicit1}, \eqref{eqn:s0p_var}, and \eqref{eqn:s0_var}, to certain quantities involving free energy functionals associated with beta random matrix ensembles. In the earlier work \cite{Spo20}, the latter functionals were expressed (using the calculus of variations) in terms of the bilinear form from Definition \ref{def:Cdef}. However, since the statement of Lemma \ref{lem:lim_var_formula} does not directly seem to appear in the literature as stated, we provide its proof in Appendix \ref{app:lin_var}.

We next state an estimate similar to Lemma \ref{lem:second_der_bd}, which will be shown in Appendix \ref{ProofQDifference}.

\begin{lem}\label{lem:bd_num_termsU}
Adopt Assumption \ref{ass:NT_assumption}. There exists a constant $\theta_0 (\beta) > 0$ such that the following holds whenever $\theta < \theta_0 (\beta)$. Fix a real number $U \in [T^{5/6}, T \log N]$. For any real number $t \in [0, T \log N]$, and integers $i, k \in \llbracket 1, N \rrbracket$ with $\varphi_0 (k) \in \llbracket - T^{8/5} \log N, T^{8/5} \log N \rrbracket$ and $|\varphi_t(i)-\varphi_t (k)| \ge U (\log N)^5$, we have with overwhelming probability that
\begin{multline}
     \min \big( Q_k(0) - Q_i(0) + t(\ve (\lambda_k) - \ve (\lambda_i)), \\
     \alpha \varphi_0 (k) - \alpha \varphi_0(i) + t(\ve (\lambda_k) - \ve (\lambda_i)) \big) > 10 U, 
     \quad \text{if $\varphi_t (k) > \varphi_t (i)$};
     \end{multline}
     \begin{multline}
 \max \big( Q_k(0) - Q_i(0) + t(\ve (\lambda_k) - \ve (\lambda_i)), \\
 \alpha \varphi_0 (k) - \alpha \varphi_0(i) + t(\ve (\lambda_k) - \ve (\lambda_i)) \big) < -10 U, 
 \quad \text{if $\varphi_t (k) < \varphi_t (i)$}.
     \end{multline}

\end{lem}

\subsection{Approximate locality of eigenvalues}

\label{Localization0}

In this section we state a result indicating that the eigenvalues of a Lax matrix $\mathbf{L}$ (Definition \ref{lt}), sampled under thermal equilibrium (Definition \ref{def:inf_thermal_eq}), are ``approximately local,'' in that they essentially only depend on a few entries of the matrix. To make this precise, the following assumption begins with a random Lax matrix $\mathbf{L}$ and ``resamples'' entries indexed by some subset $\mathcal{D}$, to form a new matrix $\tilde{\mathbf{L}}$. 

\begin{ass}\label{ass:Lcoupling}
    Adopt Assumption \ref{ass:NT_assumption}. Let $\tilde{\L} = [\tilde{L}_{ij}]$ be another random, tridiagonal symmetric matrix, with rows and columns indexed by $\llbracket N_1, N_2 \rrbracket$, that is obtained from $\mathbf{L}$ in one of the following three ways. Below, we sample $\tilde{\mathbf{a}} = (\tilde{a}_i)_{i \in \llbracket N_1, N_2 - 1 \rrbracket}$ and $\tilde{\mathbf{b}} = (\tilde{b}_i)_{i \in \llbracket N_1, N_2 \rrbracket}$ under thermal equilibrium \eqref{eqn:equil}, independently of $\mathbf{L}$. Moreover, we let  $\mathcal{D} \subseteq \llbracket N_1, N_2 \rrbracket$ denote a subset of indices. 
    \begin{enumerate}  
        \item For each $i \notin \mathcal{D}$, we set $\tilde{L}_{i,i+1} = \tilde{L}_{i+1,i} = L_{i,i+1} = L_{i+1,i} = a_i$ and $ \tilde{L}_{ii} = L_{i i} = b_i$. For each $i \in \mathcal{D}$, we set $\tilde{L}_{i+1,i}  =\tilde{L}_{i,i+1} = \tilde a_i$ and $\tilde{L}_{i i} = \tilde b_i$.
        \item For each $i,j \notin \mathcal{D}$, we set $L_{i j} = \tilde{L}_{ij}$. Otherwise, we set $\tilde{L}_{i j}  = 0$.

        \item Fix a real number $t \in [0, T (\log N)^{10}]$. Let $\tilde{\mathbf{L}}$ have the same law as $\mathbf{L}(t)$, and be coupled with $\mathbf{L}$ such that the following holds. For $\mathcal{D} = [N_1, N_1 + T^2] \cup [N_2-T^2, N_2]$, we have with overwhelming probability that $|L_{i j} - \tilde{L}_{ij}| \leq e^{-T^2/5}$  for all $i, j \notin \mathcal{D}$. The existence of such an $\tilde{\mathbf{L}}$ follows from \cite[Proposition 2.5]{Agg25a}. 
    \end{enumerate}
\end{ass}

The following lemma indicates that the eigenvalues and their associated localization centers (if they are not too close to $\mathcal{D}$) of $\mathbf{L}$ are approximately equal to those of $\tilde{\mathbf{L}}$. Its first and second parts follow from \cite[Corollary 5.5 ]{Agg25a} and \cite[Corollary 5.6]{Agg25a}, respectively, with the $\delta$ there equal to either $0$ (in the first two cases of Assumption \ref{ass:Lcoupling}) or $e^{-T^2/5} < e^{-(\log N)^3}$ (in the third case).
\begin{lem}[{\cite[Corollaries 5.5 and 5.6]{Agg25a}}]\label{lem:Lax_eig_coupling}
    Adopt Assumption~\ref{ass:Lcoupling}, and recall Definition \ref{def:quasi_intro}. There exists a constant $\mathfrak{c}>0$ such that the following holds with overwhelming probability. 
    \begin{enumerate}
    \item For each $i \in \llbracket N_1, N_2 \rrbracket$ with $\text{dist}(i, \mathcal{D}) >  (\log N)^3$, there is an eigenvalue $\tilde{\lambda} \in \eig \tilde{\mathbf{L}}$ such that 
    \begin{equation*}
        |\Lambda_{i} - \tilde{\lambda}| \leq e^{-c (\log N)^3},
        \end{equation*}
        and $i$ is an $N^{-1} \zeta$-localization center for $\tilde{\lambda}$.
   \item   For each $i \in \llbracket N_1, N_2 \rrbracket$ with $\text{dist}(i, \mathcal{D}) >  (\log N)^3$, there is a unique $\psi(i) \in \llbracket N_1, N_2 \rrbracket$ such that 
        \begin{equation*}
        |\Lambda_{\psi(i)} - \tilde{\Lambda}_{i}| \leq e^{-c (\log N)^3}.
        \end{equation*}
    
    \noindent Moreover, this $\psi(i)$ satisfies $|\psi(i) - i| \leq (\log N)^2$.
    \end{enumerate}
\end{lem}

\subsection{Resolvents} 

\label{Resolvent} 

In this section we discuss properties of the resolvent of random Lax matrices. Fix an index set $\mathscr{I} \subset \mathbb{Z}$, and let $n = |\mathscr{I}|$. Let $\mathbf{M} = [M_{ij}]$ denote an $n \times n$ symmetric matrix, whose rows and columns are indexed by $i, j \in \mathscr{I}$. For any $z \in \mathbb{C} \setminus \eig \mathbf{M}$, the \emph{resolvent} of $\mathbf{M}$ is $\mathbf{G}(z) = (\mathbf{M}-z)^{-1} = [G_{ij} (z)]$. Denote $\eig \mathbf{M} = (\lambda_1, \lambda_2, \ldots , \lambda_n)$, and let $(\mathbf{u}_1, \mathbf{u}_2, \ldots , \mathbf{u}_n)$ denote an orthonormal family of eigenvectors of $\mathbf{M}$, such that $\mathbf{u}_j = (\mathbf{u}_j (i))_{i \in \mathscr{I}}$ is an eigenvector of $\mathbf{M}$ with eigenvalue $\lambda_j$ for each $j \in \llbracket 1, n \rrbracket$. Then, for any $i,j \in \mathscr{I}$ observe, denoting $E = \Real z$ and $\eta = \Imaginary z$, that  
\begin{flalign}
\label{gijuij}
    & G_{ij}(z) = \sum_{k=1}^n \displaystyle\frac{u_k (i) u_k (j)}{\lambda_k - z}, \qquad \text{so that} \qquad \Imaginary G_{ii} (z) = \eta \displaystyle\sum_{k=1}^n \frac{u_k (i)^2}{(\lambda_k - E)^2 + \eta^2}.
\end{flalign} 

\noindent Since for each $k \in \llbracket 1, n \rrbracket$ we have 
$\eta\int_{-\infty}^{\infty} ((\lambda_k-E)^2+\eta^2)^{-1} d E = \pi$, it follows for any index $i \in \llbracket 1, \mathscr{I} \rrbracket$, real numbers $E \in \mathbb{R}$ and $\eta > 0$; interval $J \in \mathbb{R}$; and function $\phi: \mathbb{R} \rightarrow \mathbb{R}$ that 
\begin{equation}
\label{integralg2}
\frac{1}{\pi} \int_J \big| \phi(E) \cdot \Imaginary G_{i i}(E+\i\eta) \big| dE \leq \displaystyle\sup_{x \in J} |\phi(x)|.
\end{equation}

The next lemma approximates certain functionals of a random Lax matrix through its resolvent. 

\begin{lem}[{\cite[Lemma D.3]{Agg25}}]\label{lem:res_linst}
Adopt Assumption \ref{ass:NT_assumption}, and denote the resolvent of $\mathbf{L}$ by $\mathbf{G} = [G_{ij}] = (\mathbf{L} - z)^{-1}$, for any complex number $z \in \mathbb{C}$. Let $H : \mathbb{R} \rightarrow \mathbb{R}$ be a continuous function such that $|H(x)| \le A$ for all $x \in \mathbb{R}$, and
\begin{equation}
\label{eqn:Hconds}
|H(x) - H(y)| \le e^{-(\log N)^2},
\qquad
\text{for any } x,y \in [-\log N, \log N]
\text{ with } |x-y| \le e^{-(\log N)^{5/2}}.
\end{equation}
Fix integers $n_1, n_2 \in \llbracket N_1, N_2 \rrbracket$ with $n_2 \ge n_1$; denote $n = n_2 - n_1 + 1$, and assume that $n \ge (\log N)^5$. 
In addition, let $n_1' = n_1 + \floor{(\log N)^5}$ and $n_2' = n_2 - \floor{(\log N)^5}$, and set $\eta = e^{-(\log N)^3}$. Then, with overwhelming probability, we have
\begin{equation}\label{eqn:Hlr}
\left|
\sum_{i=n_1}^{n_2} H(\Lambda_i)
-
\frac{1}{\pi}
\int_{-\log N}^{\log N}
H(E)
\sum_{i=n_1'}^{n_2'}
\Imaginary G_{ii}(E + i\eta)\, dE
\right|
\le 6A (\log N)^5.
\end{equation}
\end{lem}

The next lemma indicates that the resolvents of the matrix $\mathbf{L}$ and its resampled perturbation $\tilde{\mathbf{L}}$, from Assumption \ref{ass:Lcoupling}, are approximately equal (away from entries with an index close to $\mathcal{D}$).

\begin{lem}[{\cite[Lemma 5.4]{Agg25a}}]\label{lem:resolv_coupling}
There exists a constant $\mathfrak{c} > 0$ such that the following holds with overwhelming probability.
Adopt Assumption~\ref{ass:Lcoupling}; let $\delta = 0$ in the first two cases listed there, and let $\delta = e^{-T^2/5}$ the third. For any complex number $z \in \mathbb{C}$, set $\mathbf{G}(z) = [ G_{ij} (z)] = (\mathbf{L}-z)^{-1}$ and $\tilde{\mathbf{G}} (z) = [\tilde{G}_{ij} (z)] = (\tilde{\mathbf{L}}-z)^{-1}$. Let $\eta \in [\delta,1]$ be a real number, and define the set
\begin{flalign}
    \label{omegaset}
\Omega = \{ z \in \mathbb{C} : -N \le \Re z \le N,\ \eta \le \Imaginary z \le 1 \}.
\end{flalign}

\noindent Then, for any integers $i,j \in \llbracket  N_1, N_2  \rrbracket$, we have
\begin{equation} \label{eq:5.7}
\sup_{z \in \Omega} \bigl| G_{ij}(z) - \tilde{G}_{ij}(z) \bigr|
\le e^{(\log N)^2} \eta^{-2} \bigl( \delta^{1/4} + e^{-\mathfrak{c}\,\mathrm{dist}(i,\mathcal{D}) - \mathfrak{c}\,\mathrm{dist}(j,\mathcal{D})} \bigr).
\end{equation}
\end{lem}

Observe that, while Lemma \ref{lem:resolv_coupling} is stated for matrices indexed by the same set $\llbracket N_1, N_2 \rrbracket$, it also directly applies to matrices indexed by different sets $\llbracket N_1, N_2 \rrbracket \subseteq \llbracket \mathfrak{N}_1, \mathfrak{N}_2 \rrbracket$ (upon extending the matrix indexed by $\llbracket N_1, N_2 \rrbracket$ by $0$ to $\llbracket \mathfrak{N}_1, \mathfrak{N}_2 \rrbracket$). However, the factor of $e^{(\log N)^2}$ in \eqref{eq:5.7} would then become $e^{(\log \mathfrak{N})^2}$, where $\mathfrak{N} = \mathfrak{N}_2 - \mathfrak{N}_1 + 1$ is the dimension of the larger matrix. So, if $\mathfrak{N}$ is much larger than $N$, then the right side of \eqref{eq:5.7} might not be small. 

The following lemma addresses this issue under the following notation. Let $\mathfrak{N}_1 \le N_1 < N_2 \le \mathfrak{N}_2$ be integers; denote $N = N_2 - N_1 + 1$ and $\mathfrak{N} = \mathfrak{N}_2 - \mathfrak{N}_1+1$; and let $\mathbf{L} = [L_{ij}]$ and $\bm{\mathfrak{L}} = [\mathfrak{L}_{ij}]$ denote $N \times N$ and $\mathfrak{N} \times \mathfrak{N}$ Lax matrices (Definition \ref{lt}) with rows and columns indexed by $i,j \in \llbracket N_1, N_2 \rrbracket$ and $i,j \in \llbracket \mathfrak{N}_1, \mathfrak{N}_2 \rrbracket$, respectively, sampled under thermal equilibrium (Definition \ref{def:inf_thermal_eq}). Let them be coupled so that $L_{ij} = \mathfrak{L}_{ij}$ for each $i,j \in \llbracket N_1, N_2 \rrbracket$. For any $z \in \mathbb{C}$, denote the resolvents $\mathbf{G}(z) = [G_{ij} (z)] = (\mathbf{L}-z)^{-1}$ and $\bm{\mathfrak{G}} (z) = [\mathfrak{G}_{ij} (z)] = (\bm{\mathfrak{L}} - z)^{-1}$. 
	
	\begin{lem}[{\cite[Lemma D.1]{Agg25}}] 
		
		\label{glambda} 

        Under the notation above, there exists a constant $\mathfrak{c}>0$ such that the following holds with probability at least $1 - \mathfrak{c}^{-1} e^{-\mathfrak{c} (\log N)^2}$. Set $\eta = e^{-(\log N)^3}$, and define $\Omega$ as in \eqref{omegaset}. Then, for any integer $i \in \llbracket N_1 + (\log N)^4, N_2 - (\log N)^4 \rrbracket$, we have  
		\begin{flalign*}
			\displaystyle\sup_{z \in \Omega} | \mathfrak{G}_{ii} (z) - G_{ii} (z) | \le \mathfrak{c}^{-1} e^{-\mathfrak{c} (\log N)^4}.
		\end{flalign*}
	\end{lem}

\subsection{Concentration estimates}
\label{subsec:conc_estimates}

In this section we state several concentration estimates for functionals of a random Lax matrix. Throughout, we recall $\alpha$ and $\varrho$ from Definition \ref{def:alpha_Laxdos}. 

We begin with the below lemma, proven in Appendix \ref{ProofConcentrationK}. It is a minor generalization of \cite[Proposition 4.2]{Agg25} (which coincides with it at $\mathfrak{K} = \log N$). In \eqref{eqn:agg4.2} below, the integral of $F(\lambda) \varrho(\lambda)$ arises since $\lambda$ is the density of states for the eigenvalues $(\lambda_i)$ of $\mathbf{L}$ (see Remark \ref{lrho}); that of $G(\alpha q)$ arises since the average spacing between the $(q_i)$, and thus the $(Q_i)$, is $\alpha$ (see \eqref{alphaexpectation}).

\begin{lem}\label{lem:Ndeltabd}

    Adopt Assumption \ref{ass:NT_assumption} and recall Definition \ref{def:quasi_intro}.  Let $A, B \ge 0$; $\mathfrak{K} \in [\log N, N^{1/200}]$; and $s, S \in [0, T \log N]$ be real numbers, with $S \ge 1$. Also let $F : \mathbb{R} \rightarrow \mathbb{R}$ be a continuous function satisfying the following properties. 
   \begin{enumerate}
       \item We have $\sup_{|x| \le \mathfrak{K}} |F(x)| \leq A$, and $F(x) \leq A e^{|x|^{1/2}}$ for all $x \in \mathbb{R}$.
       \item For all $x,y \in [-\mathfrak{K}, \mathfrak{K}]$ satisfying $|x-y| \leq e^{-\mathfrak{K}^{5/2}}$, we have $ |F(x)-F(y)| \leq A e^{-\mathfrak{K}^2}$. 
   \end{enumerate}
   In addition, let $G : \mathbb{R} \rightarrow \mathbb{R}$ be a function satisfying the following properties.
   \begin{enumerate}
       \item We have $\supp G \subseteq [-S,S]$
       \item For all real numbers $x \ge y$, we have $\lvert G(x) - G(y) \rvert \le B S^{-1} (x - y) + B\,\mathbbm{1}_{\{x \ge 0 \ge y\}}$

   \end{enumerate}
   
   \noindent Then, with probability at least $1-\mathfrak{c}^{-1} e^{-\mathfrak{c} \mathfrak{K}^{2}}$, for any index $j \in \llbracket 1,N \rrbracket$ satisfying $\varphi_s (j) \in \llbracket N_1 + T^4, N_2 - T^4 \rrbracket$ we have  
\begin{equation}\label{eqn:agg4.2}
\left|
\sum_{i=1}^N F(\lambda_i) \cdot G\bigl(Q_i(s) - Q_j(s)\bigr)
\;-\;
\int_{-\infty}^{\infty} F(\lambda)\,\varrho(\lambda)\, d\lambda
\int_{-\infty}^{\infty} G(\alpha q)\, dq
\right|
\;\le\;
A B S^{1/2} \mathfrak{K}^{13}.
\end{equation}

\end{lem}

To state the next estimates, we require the following assumptions on functions $F$ and $G$; they are close to those in Lemma \ref{lem:Ndeltabd} when $\mathfrak{K} = \log N$. 

\begin{ass}[Assumptions on $F$]

\label{ass:F}

Let $A \ge 0$ be a real number and $F : \mathbb{R} \rightarrow \mathbb{R}$ be a continuous function satisfying the following properties.

    \begin{enumerate}
    \item We have 
    \begin{equation}
    \label{estimatef} 
     \sup_{|x| \le \log N} |F(x)| \leq A, \qquad \text{and} \qquad |F(x)| \leq A e^{|x|^{1/2}}, \quad \text{for all $x \in \mathbb{R}$}.
     \end{equation}
    \item   For any $x,y \in [- \log N, \log N] $ with $|x-y| \leq e^{- (\log N)^{5/2}}$, we have 
    \begin{equation}\label{eqn:F_lip}
            |F(x) - F(y)| \leq A e^{-(\log N)^2}.
    \end{equation}
        
    \end{enumerate}

    \end{ass}

\begin{ass}[Assumptions on $G$]\label{ass:G}

Let $B \ge 0$ and $S \ge U \ge 1$ be real numbers and $G : \mathbb{R} \rightarrow \mathbb{R}$ be a continuously differentiable function satisfying the following properties. 
    \begin{enumerate} 
    
    \item We have $\supp G \subseteq [- S, S]$.
    \item For all $x \in \mathbb{R}$, we have $|G(x)| \leq B$ and $|G'(x)| \leq BU^{-1}$.
     \item There exist two intervals $I_1, I_2 \subset \mathbb{R}$, with $|I_j| \leq U$ for each $j \in \{ 1,2 \}$, such that $\supp G' \subseteq  I_1 \cup I_2$. 
\end{enumerate}
    \end{ass}

Under these assumptions, we have the following two lemmas. The first is a concentration bound for Lax matrix functionals\footnote{In these statements, the ranges of the parameters $S$, $U$, and $q$ are not optimal but will be sufficient for our purposes (where $S$ and $U$ will typically be slightly less than $T$, and $q$ will be of order $T$).} similar to (but slightly different from) those in Lemma \ref{lem:Ndeltabd}; the second approximates the expectations of these functionals (where the error there is likely not optimal, but it is less than $T^{1/2-c}$, so it will suffice for our purposes). The first lemma is proved in Appendix \ref{ProofEstimateConcentration}, and the second in Appendix \ref{app:lemma3.18}.

\begin{lem} 

\label{lem:no_q_conc}

Adopt Assumption \ref{ass:NT_assumption} and recall Definition \ref{def:quasi_intro}. There exist constants $\theta_0 (\beta) > 0$ and $\mathfrak{C}>1$ such that the following holds whenever $\theta < \theta_0 (\beta)$. Let $A, B>0$ and $S \geq U \geq T^{7/8}$ be real numbers, and suppose that $F, G : \mathbb{R} \rightarrow \mathbb{R}$ satisfy Assumption \ref{ass:F} and Assumption \ref{ass:G}, respectively. Further fix integers $N_1 \leq n_1 \leq n_2 \leq N_2$ and real numbers $q \in [-T^{8/5}, T^{8/5}]$ and $t \in [0, T \log N]$. 

\begin{enumerate} 

\item With overwhelming probability, we have 
    \begin{equation}\label{eqn:no_q_conc}
        \Bigg|\sum_{i=n_1}^{n_2}F(\Lambda_i) G(q-\alpha i- t \ve(\Lambda_i)) - \mathbb{E} \bigg[ \sum_{i=n_1}^{n_2}F(\Lambda_i) G(q-\alpha i- t \ve(\Lambda_i)) \bigg] \Bigg|\leq A B S^{1/2}(\log N)^{\mathfrak{C}}.
    \end{equation}
\label{item:conc1}
\item  Assume further that $S \ge \mathfrak{M}$, and fix a real number $\Delta \in \mathbb{R}$ with $|\Delta| \leq S$ and such that $t + \Delta \in [0, T \log N]$. For each real number $\Lambda \in \mathbb{R}$ and each integer $i \in \llbracket N_1, N_2 \rrbracket$, define the (random) function $G_i : \mathbb{R} \rightarrow \mathbb{R}$ by setting 
   \begin{equation}
 G_i(x)  = \chi(x) - \chi \big(x + \Delta (\ve(\Lambda) - \ve(\Lambda_i)) \big),
\end{equation}

    \noindent for each $x \in \mathbb{R}$, where $\chi$ is as in Definition \ref{def:chi_def}. The following holds with overwhelming probability. For each $\Lambda \in \mathbb{R}$, we have
    \begin{equation}\label{eqn:no_q_i_conc}
        \Bigg|\sum_{i=n_1}^{n_2}F(\Lambda_i) G_i(q -\alpha i- t \ve(\Lambda_i)) - \mathbb{E} \bigg[ \sum_{i=n_1}^{n_2}F(\Lambda_i) G_i(q-\alpha i- t \ve(\Lambda_i)) \bigg]\Bigg|\leq A (|\Lambda|+1)^2 S^{1/2}(\log N)^{\mathfrak{C}}.
    \end{equation}
\label{item:conc2}
\item Adopt the assumptions and notation of Item \ref{item:conc2} above, and assume further that $\Delta \in [0, S]$. Then the following holds with overwhelming probability. For any integer
$k \in \llbracket N_1, N_2 \rrbracket$ and real number $\Lambda \in \mathbb{R}$, we have
\begin{equation}
\label{eqn:no_q_i_conc2}
\Bigg|\sum_{i=n_1}^{n_2}F(\Lambda_i) G_i(q_k(0)-q_i(0))
- \mathbb{E} \bigg[ \sum_{i=n_1}^{n_2}F(\Lambda_i) G_i(q_k(0)-q_i(0)) \bigg]\Bigg|
\leq A (|\Lambda|+1)^2 S^{1/2}(\log N)^{\mathfrak{C}}.
\end{equation}
\label{item:conc3}
\end{enumerate} 

\end{lem}

\begin{lem}\label{lem:no_q_conc_expect}
  
   Adopt Assumption \ref{ass:NT_assumption} and recall Definition \ref{def:quasi_intro}. There exist constants $\theta_0 (\beta) > 0$ and $\mathfrak{C}>1$ such that the following holds, whenever $\theta < \theta_0 (\beta)$. Let $A, B > 0$, and $S \in [T^{7/8}, T \log N]$ be real numbers. Also let $F, G: \mathbb{R} \rightarrow \mathbb{R}$ be functions satisfying Assumption \ref{ass:F} and Assumption \ref{ass:G} at $U=S$, respectively, and further suppose that $\sup_{x\in \mathbb{R}} |G''(x)| \leq B S^{-2}$. Then for any real numbers $q \in [-T^{8/5}, T^{8/5}]$ and $t \in [0,  T (\log N)^{10}]$, and any integers $n_1, n_2 \in \llbracket N_1, N_2 \rrbracket$ with $n_2 - n_1 + 1 \geq N^{10^{-5}}$, we have 
    \begin{multline}\label{eqn:expectation_approx}
        \Bigg|\mathbb{E} \left[\sum_{i=n_1}^{n_2}F(\Lambda_i) \cdot G( q-\alpha i - t \ve(\Lambda_i))\right] -   \int_{n_1}^{n_2} \int_{-\infty}^{\infty} F(\lambda) \cdot G(q  -\alpha q' - t \ve(\lambda)) \cdot \varrho(\lambda) d \lambda d q'   \Bigg| \\
        \leq AB S^{-1/3} T^{2/3} (\log N)^{\mathfrak{C}}.
    \end{multline}
\end{lem}

\section{Random Lax matrix functionals as independent sums}
\label{sec:ind_linear}

Let us first set the notation and assumptions that will be in place throughout this section, even when not explicitly stated. First, we adopt Assumption \ref{ass:NT_assumption} and recall Definitions \ref{def:alpha_Laxdos}, \ref{def:dress_intro}, \ref{def:quasi_intro} and \ref{def:smoothed_log}. We will additionally assume that $\theta<\theta_0(\beta)$ is sufficiently small, so that all of the results from Sections \ref{ProofQ} and \ref{sec:bg} apply. We moreover fix the following choices for certain mesoscopic large parameters (recall $\mathfrak{M}$ below is as in Definition \ref{def:smoothed_log}) that we adopt throughout.

\begin{ass}\label{ass:R_M_ass}

Define the integer $R \ge 1$ and real number $\mathfrak{a}_0>0$ by 
    \begin{flalign}\label{ra0}
R \coloneqq \floor{T^{1/5}}; \qquad \mathfrak{a}_0 = 10^{-2},
\end{flalign} 

    \noindent so that (as $T = N^{1/100}$ by Assumption \ref{ass:NT_assumption} and $\mathfrak{M} = T^{19/20}$ by Definition \ref{def:smoothed_log})
    \begin{flalign} 
    \label{eqn:R_cond}
    \quad RT\mathfrak{M}^{-1} \le T^{1/2-2\mathfrak{a}_0}; \qquad T^{3/2} (R^{1/2} \mathfrak{M})^{-1} \le T^{1/2-2\mathfrak{a}_0}; \qquad \mathfrak{M} \le T^{1-2\mathfrak{a}_0}.
\end{flalign}

\end{ass}

\noindent Before proceeding, we also set notation for (centered) differences between consecutive $(q_i (0))$.

\begin{definition} 

\label{xi} 

    For each integer $i \in \llbracket N_1, N_2 - 1 \rrbracket$, define the \emph{stretch random variables}
   \begin{equation}\label{eqn:Xidef}
        X_i \coloneqq q_{i+1}(0)-q_i(0)- \alpha, \qquad \text{so that} \quad \mathbb{E}[X_i] = 0, \quad \text{by \eqref{alphaexpectation}}.
   \end{equation}

\end{definition} 

\subsection{Expressions for $\Xi$ and $\Xi^{[m]}$ as independent sums} 

\label{CoupleXH} 

In this section we state a result (Proposition \ref{prop:equiv_expr_cor}) indicating how the functions $(\Xi, \Xi_1, \Xi^{[m]}, \Xi_1^{[m]})$ from Theorem \ref{thm:couplethm12} can be approximated by independent sums. Its statement requires some notation, that involves certain constants $c^{(\Lambda)}$ and $c^{[m]}$ and also certain eigenvalues $(\lambda_j^{(r)})$. 

To that end, for any real numbers $\Lambda, r, k, q, q', t \in \mathbb{R}$ and integer $m \ge 0$, define (recalling $R$ from \eqref{ra0}; $\mathfrak{l}$ and $\chi$ from Definition \ref{def:smoothed_log}; $\alpha$ and $\varrho$ from Definition \ref{def:alpha_Laxdos}; and $\ve$ from Definition \ref{def:dress_intro})
\begin{flalign}
 & c^{(\Lambda)}(r,k,t) = 2 \int_{-\infty}^{\infty} \big( \chi( \alpha k - \alpha r R )
- \chi(\alpha k - \alpha r R + t \ve(\Lambda) - t\ve(\lambda)) \big) \cdot \mathfrak{l}(\Lambda - \lambda)  \varrho(\lambda) d\lambda; \label{eqn:clambda1} \\
& c^{[m]}(r,q,q',t) = \int_{-\infty}^{\infty} \big(\chi(q-\alpha r R )- \chi(q' -\alpha rR - t\ve(\lambda)) \big) \cdot \lambda^m \varrho(\lambda) d \lambda . \label{eqn:cm1}
\end{flalign}

\begin{definition}[Lax submatrices and eigenvalues] \label{def:Lr}

    For each integer $r \in \llbracket N_1/R, N_2/R-1 \rrbracket$, let $\mathbf{L}^{[r]} = [L_{ij}^{[r]}]$ denote the $R \times R$ submatrix of $\mathbf{L}$ (recall Assumption \ref{ass:NT_assumption}) whose rows and columns are indexed by $i, j \in \llbracket Rr, (R+1)r-1 \rrbracket$, namely, set $L_{ij}^{[r]} = L_{ij}$ for each $i, j \in \llbracket Rr, (R+1)r-1 \rrbracket$. Further let $\eig \mathbf{L}^{[r]} = (\lambda_1^{(r)}, \lambda_2^{(r)}, \ldots , \lambda_R^{(r)})$. 
    
\end{definition}

Note that the matrices $\mathbf{L}^{[r]}$ from Definition \ref{def:Lr} are each sampled according to thermal equilibrium \eqref{eqn:equil} (with the $N$ there equal to $R$ here), and they are moreover mutually independent. 

We next state the following proposition, to be shown in Section \ref{SumGeneral}, that estimates the quantities $(\Xi, \Xi_1, \Xi^{[m]}, \Xi_1^{[m]})$ from Definitions \ref{xilambdakt}, \ref{xim}, and \ref{xi2}. Since the $(\lambda_j^{(r)})$ are independent over $r \in \llbracket N_1/R, N_2/R - 1 \rrbracket$ (and as are the $(X_j)$), the below proposition can be interpreted as approximating these functions and their differences by sums of independent random variables.

\begin{prop}[Independent sum approximation]\label{prop:equiv_expr_cor}

There exists a constant $\mathfrak{c}>0$ such that the following holds. Let $t \in [0, T \log N]$ be a real number. Further fixing an integer $k \in \mathbb{Z}$ with $|k| \leq T^{3/2}$ and a real number $\Lambda \in [-\log N, \log N]$, we have with overwhelming probability that 
 \begin{equation}\label{eqn:psiZktauapprox}
 \Bigg| \Xi(\Lambda, k, t) -  \Xi_1(\Lambda, k, t) + \alpha^{-1} T^{-1/2}   \sum_{r = \lceil N_1/R \rceil}^{\lfloor N_2/R \rfloor - 1} \alpha^{-1} c^{(\Lambda)}(r,k,t) \sum_{i = r R+1}^{(r+1) R} X_{i}   \Bigg| 
\leq T^{-\mathfrak{c}},
 \end{equation} 
 and
  \begin{flalign}
  \label{eqn:psi2Zktauapprox}
  \begin{aligned}
\Bigg|  & \Xi_1  (\Lambda, k, t) \\
& -  2T^{-1/2} \Bigg( \sum_{r = \lceil N_1/R \rceil}^{\lfloor N_2/R \rfloor - 1} \sum_{s=1}^{R} \mathfrak{l}(\Lambda -\lambda_s^{(r)}) \big( \chi(\alpha k- \alpha r R)  - \chi (\alpha k - \alpha rR  + t \ve(\Lambda)- t \ve(\lambda_s^{(r)})) \big)   \\
 & -\mathbb{E} \bigg[ \sum_{r = \lceil N_1/R \rceil}^{\lfloor N_2/R \rfloor - 1} \sum_{s=1}^{R} \mathfrak{l}(\Lambda -\lambda_s^{(r)}) \big( \chi(\alpha k- \alpha rR)   - \chi (\alpha k- \alpha rR + t \ve(\Lambda)- t \ve(\lambda_s^{(r)}) ) \big) \bigg]  \Bigg)\Bigg|  \leq T^{-\mathfrak{c}}.
  \end{aligned}
 \end{flalign} 

\noindent Additionally, fixing an integer $m \in \llbracket 0, (\log N)^{1/10} \rrbracket$ and real numbers $q, q' \in [-  (\log N)^{4} T, (\log N)^{4} T]$, we have with overwhelming probability that 
 \begin{multline}\label{eqn:lambdam_approx}
\Bigg| \Xi^{[m]}(q,q', t) - \Xi_1^{[m]}(q,q', t) 
+ \alpha^{-1} T^{-1/2} \sum_{r = \lceil N_1/R \rceil}^{\lfloor N_2/R \rfloor - 1} c^{[m]}(r,q,q',t) \sum_{i = r R}^{(r+1) R-1} X_{i} \Bigg| \leq T^{-\mathfrak{c}}
  \end{multline}
  and
 \begin{flalign}\label{eqn:lambdam2_approx}
 \begin{aligned}
 \Bigg|  \Xi_1^{[m]}( & q,q', t)  - T^{-1/2} \Bigg(  \sum_{r = \lceil N_1/R \rceil}^{\lfloor N_2/R \rfloor - 1} \sum_{s=1}^{R} (\lambda_s^{(r)})^m \big( \chi(q- \alpha r R ) - \chi(q'- \alpha r R  - t \ve(\lambda_s^{(r)}) ) \big) \\
  &   - \mathbb{E}\bigg[ \sum_{r = \lceil N_1/R \rceil}^{\lfloor N_2/R \rfloor - 1} \sum_{s=1}^{R} (\lambda_s^{(r)})^m \big( \chi(q- \alpha r R) - \chi (q' - \alpha r R - t \ve(\lambda_s^{(r)})) \big) \bigg] \Bigg) \Bigg| \leq T^{-\mathfrak{c}}.
  \end{aligned}
 \end{flalign} 
\end{prop}

\subsection{Fluctuations of differences of functionals}\label{subsec:ext_BM}

In this section we approximate differences between functionals of the Lax matrix, such as those appearing in \eqref{eqn:psiZktauapprox} and \eqref{eqn:lambdam_approx}, by linear combinations of the $(X_j)$ (from Definition \ref{xi}). To eventually address both of these estimates simultaneously (and since it will also be useful later in Section \ref{sec:charge_fluct} below), we show a broader result Lemma \ref{lem:HBM_lem} that studies the following functionals $\Xi^{(H)}$ and $\Xi_1^{(H)}$ generalizing those from Definitions \ref{xilambdakt}, \ref{xi2}, and \ref{xim}. Observe that the function $\mathcal{G}_N^{(H)}$ defined by \eqref{eqn:GNH_def} below is similar to the linear combination of the $(X_i)$ appearing in \eqref{eqn:psiZktauapprox} and \eqref{eqn:lambdam_approx}

    \begin{definition}\label{def:Hfunc_def}

    Let $F, G : \mathbb{R} \rightarrow \mathbb{R}$ be functions, and define $H : \mathbb{R}^2 \rightarrow \mathbb{R}$ by setting
\begin{equation}\label{eqn:H_FG_def}
    H(\lambda, q) \coloneqq F(\lambda) \cdot G(q), \qquad \text{for any $(\lambda, q) \in \mathbb{R}^2$}.
\end{equation}

\noindent Given a function $H : \mathbb{R}^2 \rightarrow \mathbb{R}$, define the functions $\Xi^{(H)} : \mathbb{R} \times \mathbb{R}_{\ge 0} \rightarrow \mathbb{R}$ by setting 
   \begin{flalign}\label{eqn:H_functional}
       & \Xi^{(H)}(k, t) \coloneqq \sum_{i = N_1}^{N_2} H(\Lambda_i, q_i(0)-q_k(0) + t \ve(\Lambda_i)); \\
       \label{eqn:H_functional2}
       & \Xi_1^{(H)}(k, t) \coloneqq \sum_{i = N_1}^{N_2} H(\Lambda_i, \alpha (i - k) + t \ve(\Lambda_i)),
    \end{flalign}

    \noindent for any integer $k \in \llbracket N_1, N_2 \rrbracket$ and real number $t \geq 0$. Note that, unlike in \eqref{eqn:psiZktau} and \eqref{eqn:Xi2_outline}, we have not subtracted the means of the sums above. Further define (recalling the $(X_j)$ from \eqref{eqn:Xidef})
\begin{multline}\label{eqn:GNH_def}
  \mathcal{G}_N^{(H)}(k, t) \\
  =  - (\alpha T^{1/2})^{-1} \Bigg( \sum_{r = \lceil k/R \rceil}^{\lfloor N_2/R \rfloor} \left(  \int_{-\infty}^{\infty} \big( H(\lambda, \alpha (r R-k) + t \ve(\lambda) )-H(\lambda, \infty) \big) \varrho(\lambda) d\lambda \right) \sum_{j=r R}^{(r+1) R-1}  X_{j}  \\
  + \sum_{r = \lceil N_1/R \rceil}^{\lfloor k/R-1 \rfloor} \bigg(  \int_{-\infty}^{\infty} \big( H(\lambda, \alpha (r R-k) + t \ve(\lambda) ) - H(\lambda, -\infty) \big) \varrho(\lambda) d\lambda   \bigg)  \sum_{j=r R}^{(r+1) R-1}  X_{j}\Bigg).
\end{multline}
    \end{definition}

The below lemma states that $\Xi^{(H)} - \Xi_1^{(H)} \approx T^{1/2} \cdot \mathcal{G}_N^{(H)}$, under the following assumption on the function $G$ prescribing $H$ in Definition \ref{def:Hfunc_def}.

\begin{ass}\label{ass:G2}

    Fix real numbers $B > 0$ and $U \geq 1$. Let $G : \mathbb{R} \rightarrow \mathbb{R}$ be a thrice continuously differentiable function such that, for some $\mathfrak{Q} \in [-T(\log N)^5, T (\log N)^5]$, we have 
        \begin{equation}
        \label{estimateg}
        \begin{gathered}
           \supp G' \subset [\mathfrak{Q}-U, \mathfrak{Q}+U], \\
           |G(x)| \leq B; \qquad  |G'(x)| \leq BU^{-1}; \qquad |G''(x)| \leq BU^{-2}; \qquad |G^{(3)}(x)| \leq BU^{-3}.
        \end{gathered}
        \end{equation}

\end{ass}

\begin{lem}\label{lem:HBM_lem}

There exists a constant $\mathfrak{c}>0$ such that the following holds. Fix an integer $k $ with $|k| \leq 10 T^{3/2}$, and a real number $t \in [0, T \log N]$. Also let $A, B > 0$ and $U \in [\mathfrak{M}, T \log N]$ be real numbers, and $F, G : \mathbb{R} \rightarrow \mathbb{R}$ be functions satisfying Assumption \ref{ass:F} and Assumption \ref{ass:G2}, respectively. Recalling the notation from Definition \ref{def:Hfunc_def}, we have with overwhelming probability that 
    \begin{equation}\label{eqn:G_Nbd}
       \left| \Xi^{(H)}(k, t)  - \Xi_1^{(H)}(k, t) 
    - T^{1/2} \mathcal{G}_N^{(H)}(k, t) \right| \leq A B T^{1/2-\mathfrak{c}}.
    \end{equation}
    In addition, for such $H, k, t$
    \begin{equation}\label{eqn:G_Nmeanbd}
       \left| \mathbb{E}[\Xi^{(H)}(k, t)]  - \mathbb{E}[\Xi_1^{(H)}(k, t) ]\right| \leq A B T^{1/2-\mathfrak{c}}.
    \end{equation}
\end{lem}

\begin{proof}

Throughout this proof, we abbreviate $q_i = q_i (0)$. Moreover, since $\supp G'$ is compactly supported, the limits $G(-\infty) \coloneqq \lim_{x \rightarrow -\infty} G(x)$ and $G(\infty) \coloneqq \lim_{x \rightarrow \infty} G(x)$ both exist. As such, let $H(\lambda, \pm \infty) \coloneqq F(\lambda) \cdot G(\pm \infty)$ for any $\lambda \in \mathbb{R}$. 

To establish the lemma, we first use a Taylor expansion, which by \eqref{eqn:H_functional} yields
\begin{flalign}\label{eqn:HBMstep1}
\begin{aligned} 
\Xi^{(H)}(k, t)  &  =  \sum_{i=N_1}^{N_2} 
H \big(\Lambda_i, \alpha (i - k) + t \ve(\Lambda_i) \big)  + \displaystyle\frac{1}{2} \sum_{i=N_1}^{N_2}  \partial_q^2 H\left(\Lambda_i, \xi_{k, i}\right) \cdot \left(q_{i}-q_{k} - \alpha(i -k) \right)^2 \\
& \qquad + \sum_{i=N_1}^{N_2} 
\partial_q H \big(\Lambda_i, \alpha (i - k) + t \ve(\Lambda_i) \big) \cdot \left(
q_{i}-q_{k} - \alpha(i -k) \right),
    \end{aligned} 
    \end{flalign}
   
   \noindent where, for each $i \in \llbracket N_1, N_2 \rrbracket$, the parameter $\xi_{k,i}$ is a real number between $q_i - q_k + t \ve(\Lambda_i)$ and $\alpha(i-k) + t \ve(\Lambda_i)$; in particular, it satisfies
   \begin{equation}\label{eqn:xibd}
      | \xi_{k, i} -\alpha (i-k)  - t \ve(\Lambda_i) | \leq |q_{k}-q_{i} - \alpha(k -i)|.
   \end{equation}

     Define the events 
    \begin{flalign*}
        \mathsf{E}_1 &\coloneqq 
        \bigcap_{i, k \in \llbracket N_1 , N_2 \rrbracket} \big\{ |q_{k}-q_{i} - \alpha(k -i)| \leq (\log N)^2 |k -i|^{1/2} \big\}; \qquad \mathsf{E}_2 \coloneqq \mathsf{BND}_{\L(0)}( \log N).
    \end{flalign*}

 \noindent Moreover, for any index $i \in \llbracket N_1, N_2 \rrbracket$, set $i_t \coloneqq \varphi_t(\varphi_0^{-1}(i)) \in \llbracket N_1, N_2 \rrbracket$, and let $\mathsf{E}_3$ denote the event on which the conditions 
\begin{multline}
    \min \big(q_k -q_i + t(\ve (\Lambda_k) - \ve (\Lambda_i)),  \alpha k - \alpha i + t(\ve (\Lambda_k) - \ve (\Lambda_i))  \big) >  10 U, \qquad \\
    \text{for all $k \in  \llbracket - T^{3/2} \log N, T^{3/2} \log N \rrbracket$ with $k_t-i_t >  U (\log N)^5$},
    \end{multline}
    and 
    \begin{multline}
    \max \big(q_k -q_i + t(\ve (\Lambda_k) - \ve (\Lambda_i)),  \alpha k - \alpha i + t(\ve (\Lambda_k) - \ve (\Lambda_i)) \big) <  -10 U, \qquad  \\
    \text{for all $k \in  \llbracket - T^{3/2} \log N, T^{3/2} \log N \rrbracket$ with $k_t-i_t <  - U (\log N)^5$}, 
\end{multline}
both hold. In what follows, we restrict to $\bigcap_{i=1}^3 \mathsf{E}_i$ . By Lemmas \ref{lem:q_spacing}, \ref{lem:bd_lem}, and \ref{lem:bd_num_termsU}, each of these events is overwhelmingly probable. 

 We first claim that, for any $i$ such that $\xi_{k,i} \in \supp G'$, we have $i \in \llbracket k - T (\log N)^8, k + T(\log N)^8 \rrbracket$.  Indeed, assuming to the contrary that this is false, without loss of generality for some $i > k + T (\log N)^8$, then both $\alpha (i - k) > T (\log N)^7$ and $q_i(0) - q_k(0) > T (\log N)^7$ (by our restriction to $\mathsf{E}_1$). So, $\alpha(i - k) + t \ve(\Lambda_i) > T (\log N)^6$ and $q_i(0) - q_k(0) + t \ve(\Lambda_i) > T (\log N)^6$, as $t \le T \log N$ and $|\ve (\Lambda_i)| \le C \log N$ by Lemma \ref{lem:ve_tail} (and our restriction to $\mathsf{E}_2$). Hence, $\xi_{k,i} > T (\log N)^6$, meaning by \eqref{estimateg} that $\xi_{k,i} \notin \supp G'$. This is a contradiction, which verifies the claim.

We next claim that there are at most $2U(\log N)^5$ nonzero terms in the sum of second derivatives in \eqref{eqn:HBMstep1}, namely, that there are at most $2U(\log N)^5$ indices $i \in \llbracket n_1, n_2 \rrbracket$ for which $\partial_q^2 H(\Lambda_i, \xi_{k,i}) \ne 0$. Indeed, let $I \in \llbracket N_1, N_2 \rrbracket$ denote the smallest index such that $\xi_{k, I} \in \supp G'$. By the above, for any other index $i \in \llbracket N_1, N_2 \rrbracket$ with $\xi_{k, i} \in \supp G'$ we must have $I, i \in \llbracket k - T(\log N)^8, k + T(\log N)^8 \rrbracket$, meaning that $I, i \in \llbracket -T^{3/2} \log N, T^{3/2} \log N \rrbracket$ (since $|k| \le 10T^{3/2}$). If $i \in \llbracket -T^{3/2} \log N, T^{3/2} \log N \rrbracket$ satisfies $ i_t-I_t =\varphi_t(\varphi_0^{-1}(i))-\varphi_t(\varphi_0^{-1}(I))  >  U (\log N)^5$, then 
\begin{multline}\label{eqn:xikidf}
\xi_{k,i}-\xi_{k,I} > (\alpha (i-k)  + t \ve(\Lambda_i)) - (\alpha (I-k)  + t \ve(\Lambda_{I})) - T^{5/6} \geq 10 U - T^{5/6} > 9 U,
\end{multline}

\noindent where to obtain the first lower bound in \eqref{eqn:xikidf} we used \eqref{eqn:xibd} and the restriction to $\mathsf{E}_1$; the fact that $k,i,i_0 \in \llbracket -T^{3/2} \log N, T^{3/2} \log N \rrbracket$; and the bound $T^{3/4} (\log N)^3 < T^{5/6}$. To obtain the second lower bound we used our restriction to $\mathsf{E}_3$, and to obtain the third we used the fact $U \ge \mathfrak{M} = T^{19/20}$. However, since $\xi_{k,I}, \xi_{k,i} \in \supp G'$, we have by \eqref{estimateg} that $|\xi_{k,i} - \xi_{k,I}| \le 2U$, which is a contradiction. Thus, for any $i$ satisfying $\partial_q^2 H(\Lambda_i, \xi_{k,i}) \ne 0$, we must have $\xi_{k, i} \in \supp G'$, so we must have $\varphi_t(\varphi_0^{-1}(i))-\varphi_t(\varphi_0^{-1}(I)) \le  U (\log N)^5$ and, by similar reasoning, $\varphi_t(\varphi_0^{-1}(i))-\varphi_t(\varphi_0^{-1}(I)) \ge - U (\log N)^5$. Since each $\varphi_s$ is a bijection, there are at most $2U (\log N)^5$ such indices $i$.

By this estimate and our restriction to $\mathsf{E}_2$,  using \ref{eqn:H_FG_def} with Assumptions \ref{ass:F} and \ref{ass:G2} on $F$ and $G$ to bound $|\partial_q^2 H(\Lambda_i , \xi_{k,i})| \le ABU^{-2}$, we may bound the sum of second derivatives in \eqref{eqn:HBMstep1} by
  \begin{equation}\label{eqn:lastlnbd}
    \left| \displaystyle\frac{1}{2} \sum_{i=N_1}^{N_2} \partial_q^2 H\left(\Lambda_i, \xi_{k, i}
    \right) \cdot \left(
q_{i}-q_{k} - \alpha(i -k)
    \right)^2 \right| \leq AB (\log N)^{C} \cdot TU^{-1}.
  \end{equation}
  
  \noindent Above, we have also used the fact that for $\partial_q^2 H (\Lambda_i, \xi_{k,i}) \ne 0$ we must have $\xi_{k,i} \in \supp G'$ and thus $|i-k| \leq T (\log N)^8$ (as shown above).
    
We next simplify the last line in \eqref{eqn:HBMstep1}. To that end, observe from \eqref{eqn:Xidef} that
   \begin{flalign}\label{eqn:HBMstep2}
   \begin{aligned} 
    & \sum_{i=N_1}^{N_2} 
\partial_q H( \Lambda_i, \alpha (i - k) + t \ve(\Lambda_i) ) \cdot \left(
q_{i}-q_{k} - \alpha(i -k)
    \right) \\
   & \quad = 
   \sum_{i=N_1}^{N_2} 
\partial_q H(\Lambda_i, \alpha (i - k) + t \ve(\Lambda_i) ) \cdot \Bigg(
\mathbbm{1}_{i \geq k} \sum_{j =k}^{i-1}X_{j}
 -  \mathbbm{1}_{i < k} \sum_{j =i}^{k-1}X_{j}  
    \Bigg)    \\
    & \quad =
 \sum_{j =k}^{N_2}X_{j} \sum_{i = j+1}^{N_2} \partial_q H( \Lambda_{i},
    \alpha (i -k) +t  \ve(\Lambda_i) ) 
 - \sum_{j =N_1}^{k-1}X_{j} \sum_{i = N_1 }^{j} \partial_q H( \Lambda_{i},
    \alpha (i -k) +t  \ve(\Lambda_i) )    .
    \end{aligned} 
       \end{flalign}

Next, we enlarge or truncate the inner sums over $i$, so that they always begin and end at integer multiples of $R$. By \eqref{estimatef}, \eqref{estimateg}, and our restriction to $\mathsf{E}_2$, we obtain that 
   \begin{flalign}\label{eqn:Hfloor_err_est1}
   \begin{aligned} 
\left|\sum_{i = R \floor{j/R} }^{N_2} \partial_q H( \Lambda_{i},
    \alpha (i -k) +t  \ve(\Lambda_i) ) -  \sum_{i = j+1}^{N_2} \partial_q H( \Lambda_{i},
    \alpha (i -k) +t  \ve(\Lambda_i) ) \right| 
\leq 
A B R U^{-1}; \\ 
\left|\sum_{i = N_1 }^{R \floor{j/R}} \partial_q H( \Lambda_{i},
    \alpha (i -k) +t  \ve(\Lambda_i) )  -  \sum_{i = N_1 }^{j}  \partial_q H( \Lambda_{i},
    \alpha (i -k) +t  \ve(\Lambda_i) ) \right| 
\leq 
A B R U^{-1} .
\end{aligned} 
       \end{flalign}

    \noindent Thus, 
    \begin{flalign}\label{eqn:BMstep2_equiv}
    \begin{aligned} 
    & \Xi^{(H)}( k, t)  - \Xi_1^{(H)}(k, t)  \\
    & \quad =  \sum_{r = \lceil k/R \rceil}^{\lfloor N_2/R \rfloor - 1}   \sum_{i = R r }^{N_2} \partial_q H( \Lambda_{i},
    \alpha (i -k) +t  \ve(\Lambda_i) ) \sum_{j=r R}^{(r+1) R-1}  X_{j}  \\
    & \qquad  - \sum_{r = \lceil N_1/R \rceil}^{\lfloor k/R \rfloor} \sum_{i = N_1 }^{R r} \partial_q H( \Lambda_{i},
    \alpha (i -k) +t  \ve(\Lambda_i) ) \sum_{j=r R}^{(r+1) R-1}  X_{j}+ O(A B RTU^{-1} (\log N)^C),
    \end{aligned} 
    \end{flalign}

    \noindent where the error satisfies $| O(A B RTU^{-1} (\log N)^C| \leq A B RTU^{-1} (\log N)^C$. To obtain this, we used \eqref{eqn:BMstep2_equiv} on the coefficients of $X_j$, when $j \ge k$ and $j \in \llbracket \lceil k/R \rceil R, \lfloor N_2/R \rfloor R - 1 \rrbracket$, or when $j \le k-1$ and $j \in \llbracket \lceil N_1/R \rceil R, k-1 \rrbracket$. By \eqref{eqn:BMstep2_equiv}, the error incurred by each such coefficient is at most $ABRU^{-1}$. Further observe by our restriction to $\mathsf{E}_1$ that 
    \begin{flalign}
        \label{xj2}
        \displaystyle\max_{j \in \llbracket N_1, N_2-1 \rrbracket} |X_j| \le (\log N)^2,
    \end{flalign}
    
    \noindent and that we have 
    \begin{flalign}
    \label{ikh}
    \partial_q H(\Lambda_i, \alpha (i-k) + t \ve(\Lambda_i)) = 0, \quad \text{unless $|i-k| \le T(\log N)^7$},
    \end{flalign} 
    
    \noindent as otherwise $|\alpha(i-k)+t\ve(\Lambda_i)| \ge T(\log N)^6 - T (\log N)^2 > 2T(\log N)^5$ (by Lemma \ref{lem:ve_tail}), meaning that $\supp G' \subseteq \supp \partial_q H$ cannot contain $\alpha(i-k)+t\ve(\Lambda_i)$ (by \eqref{estimateg}). As such, there are at most $2T(\log N)^7$ indices $j$ for which the coefficient of $X_j$ in either the right side of \eqref{eqn:HBMstep2} or \eqref{eqn:BMstep2_equiv} is nonzero. The error incurred by these terms for the above $j$ is thus at most $A B RTU^{-1} (\log N)^C$. 
    
    There are also at most $4R$ indices $j$ (those satisfying $j-N_1 \le R$, $|j-k| \le R$, or $N_2-j \le R$) not covered by the above cases. We bound the coefficients of these $X_j$ directly. To that end, observe for each $i$ that $|\partial_q H(\Lambda_i, \alpha(i-k) + t \ve(\Lambda_i))| \le ABU^{-1}$, by Assumptions \ref{ass:F} and \ref{ass:G2} (with our restriction to $\mathsf{E}_2$). Again using \eqref{xj2} and \eqref{ikh}, we deduce that the error incurred by these terms is also at most $A B RTU^{-1} (\log N)^C$. This verifies \eqref{eqn:BMstep2_equiv}.

Next we approximate the coefficients of the $X_j$ in \eqref{eqn:BMstep2_equiv}, by applying the first part of Lemma \ref{lem:no_q_conc} and Lemma \ref{lem:no_q_conc_expect}, both with the $(F(\lambda),G(x),q)$ there equal to $(F(\lambda),G'(\mathfrak{Q}-x), \mathfrak{Q} - \alpha k)$ here; observe that $G'$ satisfies Assumption \ref{ass:G} with the $(S,B)$ there equal to $(U, BU^{-1})$ here, and $|G^{(3)}(x)| \leq B U^{-3}$, so $G'(\mathfrak{Q}-x)$ also satisfies the assumptions of Lemma \ref{lem:no_q_conc_expect} with $(S, B) = (U, BU^{-1})$. By these lemmas, we have with overwhelming probability that, for any $j \in \llbracket k, N_2 \rrbracket$, 
\begin{multline}\label{eqn:HBM_bd11}
\Bigg| \sum_{i = j}^{N_2} \partial_q H( \Lambda_{i}, 
    \alpha (i-k)  +t  \ve(\Lambda_i) ) -  \int_{j}^{\infty} \int_{-\infty}^{\infty}
    \partial_q H( \lambda, 
    \alpha (q-k)  +t  \ve(\lambda) ) \cdot \varrho(\lambda) d \lambda d q\Bigg|
\\
    \leq A B T^{1/2} U^{-1} (\log N)^C.
\end{multline}
Similarly, for any $j \in \llbracket N_1, k-1 \rrbracket$,
\begin{multline}\label{eqn:HBM_bd11N1}
\Bigg| \sum_{i =N_1}^{j} \partial_q H( \Lambda_{i}, 
    \alpha (i-k)  +t  \ve(\Lambda_i) ) -  \int_{-\infty}^{j} \int_{-\infty}^{\infty}
    \partial_q H( \lambda, 
    \alpha (q-k)  +t  \ve(\lambda) ) \cdot \varrho(\lambda) d \lambda d q\Bigg|
\\
    \leq A B T^{1/2} U^{-1} (\log N)^C.
\end{multline}

\noindent Inserting \eqref{eqn:HBM_bd11} and \eqref{eqn:HBM_bd11N1} in \eqref{eqn:BMstep2_equiv} and performing the integration over $q$ in both integrals, yields 
    \begin{flalign}
    \label{eqn:HGN_F_0_intermed}
    \begin{aligned} 
\Xi^{(H)} & (k, t)  - \Xi_1^{(H)}(k, t)  \\
& = - \alpha^{-1} \Bigg(  \sum_{r = \lceil k/R \rceil}^{\lfloor N_2/R-1 \rfloor} \int_{-\infty}^{\infty}
    \big( H( \lambda, 
    \alpha (r R-k)  +t  \ve(\lambda)) - H( \lambda, 
    \infty) \big) \varrho(\lambda) d \lambda \sum_{j=r R}^{(r+1) R-1}  X_{j}  \\
  & \qquad \quad + \sum_{r \in \lceil N_1/R \rceil}^{\lfloor k/R \rfloor} \int_{-\infty}^{\infty} \big( H(\lambda, \alpha (r R-k) + t \ve(\lambda) ) - H(\lambda, -\infty) \big) \varrho(\lambda) d\lambda  \sum_{j=r' R}^{(r+1) R-1}  X_{j}\Bigg)\\
 & \qquad \quad  + O \big( A B R TU^{-1} (\log N)^C \big) + O \big(A B R^{1/2} T^{3/2} U^{-1} (\log N)^C \big).
  \end{aligned} 
    \end{flalign}
To obtain the error bound on the last line, we also used the following facts. First, by our restriction to $\mathsf{E}_1$, we have for each $r \in \llbracket N_1/R, N_2/R-1 \rrbracket$ that 
\begin{flalign*}
\Bigg| \sum_{j=r R}^{(r+1) R-1}  X_{j}  \Bigg| \leq (\log N)^2 R^{1/2}.
\end{flalign*}

\noindent Second, by \eqref{ikh}, the coefficient of $\sum_{j=rR}^{(r+1)R-1} X_j$ in \eqref{eqn:BMstep2_equiv} is nonzero for at most $TR^{-1} (\log N)^8$ indices $r$, namely, those satisfying $|rR-k| \le T(\log N)^7$. Third, also for all but at most $TR^{-1} (\log N)^8$ indices $r$, namely, those for which $|rR-k| \ge T(\log N)^7$, the coefficient of $\sum_{j=rR}^{(r+1)R-1} X_j$ in \eqref{eqn:HGN_F_0_intermed} is small, that is, it is bounded by 
\begin{flalign}
    \label{integral0} 
    \begin{aligned} 
 \int_{-\infty}^{\infty}
    \big(& \big| H( \lambda, 
    \alpha (r R-k)  +t  \ve(\lambda)) - H( \lambda, 
    \infty) \big| \cdot \mathbbm{1}_{k \le rR - T (\log N)^7} \\
    & + \big| H(\lambda, \alpha (r R-k) + t \ve(\lambda) ) - H(\lambda, -\infty) \big| \cdot \mathbbm{1}_{k \ge rR + T(\log N)^7} \big) \cdot \varrho(\lambda) d \lambda \le CAB e^{-(\log N)^2}.
    \end{aligned} 
\end{flalign}

\noindent Indeed, we only have $H(\lambda, \alpha (r R-k)  +t  \ve(\lambda)) \ne 0$ and $|k-rR| \ge T(\log N)^7$ if $\ve(\lambda) \ge (\log N)^6$, meaning that $\lambda > (\log N)^5$, by Lemma \ref{lem:ve_tail}. This quickly implies \eqref{integral0} by our assumptions on $F$ and $G$, with the fact from Lemma \ref{lem:varrho_bd} that $|\varrho(x)| \leq c^{-1} e^{-c x^2}$. Combining these three points with \eqref{eqn:BMstep2_equiv} and \eqref{eqn:HBM_bd11} yields an error of $ABT^{1/2} U^{-1} (\log N)^C \cdot (\log N)^C R^{1/2} \cdot T (\log N)^C = ABR^{1/2} T^{3/2} U^{-1} (\log N)^{3C}$ when passing from  \eqref{eqn:HBM_bd11N1} to \eqref{eqn:HGN_F_0_intermed}; this confirms the latter. 

Observing that the main term on the right side of \eqref{eqn:HGN_F_0_intermed} is $T^{1/2} \mathcal{G}_N^{(H)} (k,t)$, and that the error in its last line is at most $A B T^{1/2-c}$ (by \eqref{eqn:N1N2}, \eqref{eqn:R_cond}, and the fact that $U \ge \mathfrak{M}$), \eqref{eqn:HBM_bd11N1} establishes the first part of the lemma, \eqref{eqn:G_Nbd}.

Now we move on to the second statement. Note that $\mathcal{G}_N^{(H)}(k,t)$ has mean zero, since the $X_i$ do (by \eqref{eqn:Xidef}). Thus, if $\mathsf{E} = \bigcap_{i=1}^3 \mathsf{E}_i$ denotes the event on which we have \eqref{eqn:G_Nbd}, in order to show \eqref{eqn:G_Nmeanbd}, it suffices to show $\mathbb{E} [\mathbbm{1}_{\mathsf{E}^{\complement}} \cdot (| \Xi^{(H)}(k, t) | +|\Xi_1^{(H)}(k, t) | ) ] \leq A B c^{-1} e^{-c (\log N)^2}$; this holds since  
     \begin{flalign}
   \label{expectationH}
       \begin{aligned} 
    & \mathbb{E} \Bigg[\mathbbm{1}_{\mathsf{E}^{\complement}} \cdot \Bigg|\sum_{i = N_1}^{N_2} H(\Lambda_i, q_i(0)-q_k(0) + t \ve(\Lambda_i))  \Bigg| \Bigg] 
    + \Bigg[\mathbbm{1}_{\mathsf{E}^{\complement}} \cdot \Bigg|\sum_{i = N_1}^{N_2} H(\Lambda_i, \alpha (i - k) + t \ve(\Lambda_i))  \Bigg| \Bigg] \\
    & \quad \le 2 N \cdot   \mathbb{P} [\mathsf{E}^{\complement}]^{1/2} \cdot \sup_{q \in \mathbb{R}} |G(q)| \cdot \mathbb{E} \bigg[ \displaystyle\max_{i \in \llbracket N_1, N_2 \rrbracket}  |F(\Lambda_i)|^2 \bigg]^{1/2}   \\
    & \quad \le 2 N \cdot c^{-1} e^{-c(\log N)^2/2} \cdot B \cdot \left( \displaystyle\int_{0}^{\infty} A^2 e^{2 |\lambda|^{1/2}} \cdot \mathbb{P} \bigg[ \displaystyle\max_{\Lambda \in \eig \mathbf{L}} |\Lambda| \ge \lambda \bigg] d \lambda \right)^{1/2} \le C AB e^{-c(\log N)^2/4},
    \end{aligned}
    \end{flalign}

    \noindent where to obtain the first bound we used the definition of $H$; to obtain the second we used the fact that $\mathbb{P}(\mathsf{E}) \leq c^{-1} e^{-c (\log N)^2}$ and \eqref{estimatef}; and to obtain the third we used Lemma \ref{lem:bd_lem}. This completes the proof.
\end{proof}

\subsection{Expectations of Lax matrix functionals}

\label{LinearExpectation}

In this section we approximate the expectations of certain functionals of the Lax matrix, through the following lemma.

\begin{lem}\label{lem:xi_expectation}

     There exists a constant $\mathfrak{c}>0$ such that the following holds. For any integer $k \in \llbracket-T^{3/2} \log N, T^{3/2} \log N \rrbracket$, real numbers $\mathfrak{q} \in [-(\log N)^4, (\log N)^4]$ and $t \in [0,  T \log N]$, and function $F: \mathbb{R} \rightarrow \mathbb{R}$ satisfying Assumption \ref{ass:F},  we have 
    \begin{multline}\label{eqn:mean_we_want}
        \Bigg|\mathbb{E} \bigg[\sum_{i=N_1}^{N_2}F(\Lambda_i) \cdot \big( \chi(q_k(0)  -q_i(0) + t \mathfrak{q}- t \ve(\Lambda_i))- \chi(q_k(0)  -q_i(0) )\big)\bigg] \\
        -   t \alpha^{-1} \int_{-\infty}^{\infty} F(\lambda)(\mathfrak{q}-\ve(\lambda)) \varrho(\lambda) d \lambda   \Bigg| 
        \leq  A T^{1/2-\mathfrak{c}}.
    \end{multline}
\end{lem}

\begin{proof}

Throughout this proof, we abbreviate $q_i(0) = q_i$. First observe by Definition \ref{def:chi_def} that the function $G(x) \coloneqq \chi( x + t \mathfrak{q})$ satisfies Assumption \ref{ass:G2}, with the $(B,U)$ there equal to $(10,\mathfrak{M})$ here. Denote $H(\lambda,q) = F(\lambda) \cdot \chi(x + t \mathfrak{q})$, so that under Definition \ref{def:Hfunc_def} we have 
    \begin{flalign}\label{eqn:F_functional2_nomean}
    \begin{aligned} 
      &   \Xi^{(H)}(k, t) = \sum_{i = N_1}^{N_2}  F( \Lambda_i) \cdot  
    \chi (q_k - q_i + t \mathfrak{q} - t \ve(\Lambda_i)); \\
  & \Xi_0(t)  \coloneqq \Xi_1^{(H)}(k, t) =  \sum_{i = N_1}^{N_2}  F( \Lambda_i)  \cdot 
    \chi (\alpha (k - i)  + t \mathfrak{q} - t \ve(\Lambda_i)),
    \end{aligned} 
    \end{flalign}
   \noindent where we introduced the notation $\Xi_0$ to suppress dependence on other variables (as below we will differentiate it in $t$). Under this notation, Lemma \ref{lem:HBM_lem} implies that 
   \begin{flalign}
       \label{sumchif} 
    \begin{aligned} 
    \Bigg| \mathbb{E} \bigg[\sum_{i=N_1}^{N_2}F(\Lambda_i) \cdot \big( \chi(q_k - q_i + t \mathfrak{q}- t \ve(\Lambda_i))- \chi(q_k -q_i )& \big)\bigg] - \big( \mathbb{E}[\Xi_0(t) - \Xi_0(0)] \big) \Bigg| \le ABT^{1/2-c}.
    \end{aligned} 
   \end{flalign}

   \noindent
    Since
   \begin{flalign} 
   \label{xi0stintegral} 
   \mathbb{E}[\Xi_0(t)] - \mathbb{E}[\Xi_0 (0)] = \int_0^t \mathbb{E} [\Xi_0'(s)] ds,
   \end{flalign} 

    \noindent it therefore suffices to approximate $\mathbb{E}[\Xi_0' (s)]$.

    To that end, first observe since 
   \begin{flalign} 
   \label{eqn:xi_deriv}
       \Xi_0'(s) =  \sum_{i = N_1}^{N_2}  F( \Lambda_i) \cdot \left(\mathfrak{q} -\ve(\Lambda_i)\right) \cdot 
    \chi'(\alpha (k - i)  +  s \mathfrak{q} - s \ve(\Lambda_i)),
   \end{flalign}

    \noindent we have by Lemma \ref{lem:no_q_conc_expect}, with the $(F,G)$ there equal to $(F(\lambda) \cdot (\mathfrak{q} - \ve(\lambda)), \chi')$ here ($F$ satisfies Assumption \ref{ass:F} with $A$ there equal to $A(\log N)^5$ here and, by Definition \ref{def:chi_def}, $\chi'$ satisfies the assumptions on $G$ of Lemma \ref{lem:no_q_conc_expect} with the $S=U$ there equal to $\mathfrak{M}$ here), that 
   \begin{flalign}
   \label{derivativexi0}
    \begin{aligned}
    \Bigg|  \mathbb{E}[ \Xi_0'(s)]-  \int_{N_1}^{N_2} \displaystyle\int_{-\infty}^{\infty} F(\lambda) \cdot (\mathfrak{q} -\ve(\lambda))  \varrho(& \lambda) \cdot \chi'( \alpha k + s \mathfrak{q} - \alpha q + s \ve(\lambda)) d \lambda dq  \Bigg| \\
     & \qquad \le A T^{2/3} \mathfrak{M}^{-4/3} (\log N)^C,
    \end{aligned}
   \end{flalign}
    
   \noindent for any $s \in [0, t]$. We also have for any $s \in [0,t]$ that 
   \begin{flalign} \label{eqn:expectationbds}
\begin{aligned}
     \Bigg| \int_{-\infty}^{\infty} & F(\lambda) \cdot (\mathfrak{q} -\ve(\lambda)) \varrho(\lambda) \cdot \bigg( \alpha^{-1} -  \int_{N_1}^{N_2} \chi' ( \alpha k + s \mathfrak{q} - \alpha q + s \ve(\lambda)) d q \bigg) d \lambda  \Bigg| \\
    & =   \Bigg|   \int_{[N_1,N_2]^{\complement}} \displaystyle\int_{-\infty}^{\infty} F(\lambda) \cdot (\mathfrak{q} -\ve(\lambda))  \varrho(\lambda) \cdot \chi'( \alpha k + s \mathfrak{q} - \alpha q + s \ve(\lambda)) d \lambda d q  \Bigg| \\
    & \le \Bigg| \displaystyle\int_{-\infty}^{\infty} F(\lambda) \cdot (\mathfrak{q} -\ve(\lambda))  \varrho(\lambda) \cdot \mathbbm{1}_{|\lambda| \ge T} d \lambda  \Bigg| \leq A B e^{-T},
    \end{aligned}
   \end{flalign}

   \noindent where the first statement holds from the equality $\int_{-\infty}^{\infty} \chi' (-\alpha q)dq = \alpha^{-1}$; the second holds from the fact that, for $q \notin [N_1, N_2]$, we have by \eqref{estimatesm0} that $\chi' (\alpha k + s\mathfrak{q} - \alpha q + s\ve(\lambda)) = 0$ only if $|\ve (\lambda)| \ge t^{-1} (|\alpha q| - \alpha k - t \mathfrak{q} - \mathfrak{M}) \ge T^{10}$ (where the last bound holds as $|q| \ge N (\log N)^{-7} \ge T^{50}$ by \eqref{eqn:N1N2}), meaning by Lemma \ref{lem:ve_tail} that $|\lambda| \ge T$; and the third holds since $|F(\lambda)| \cdot |\mathfrak{q}-\ve(\lambda)| \cdot \varrho(\lambda) \leq A c^{-1} e^{-c \lambda^2}$ holds whenever $|\lambda| > T$, by Lemma \ref{lem:varrho_bd} and \eqref{estimatef} (with Lemma \ref{lem:ve_tail}).
   
   Summing \eqref{derivativexi0} and \eqref{eqn:expectationbds} gives
      \begin{flalign*}
    \Bigg|  \mathbb{E}[ \Xi_0'(s)]-  \alpha^{-1}  \displaystyle\int_{-\infty}^{\infty} F(\lambda) \cdot (\mathfrak{q} -\ve(\lambda))  \varrho( \lambda) d \lambda   \Bigg| \le A T^{2/3} \mathfrak{M}^{-4/3} (\log N)^C,
   \end{flalign*}

   \noindent which upon integrating over $s \in [0, t] \subseteq [0, T \log N]$ yields
   \begin{equation}\label{eqn:final_expect_est}
     \Bigg|  \int_0^t \mathbb{E} [\Xi_0'(s)] ds - t \alpha^{-1}  \int_{-\infty}^{\infty}F(\lambda) \cdot (\mathfrak{q} -\ve(\lambda)) \varrho(\lambda) d \lambda \Bigg| \leq A T^{5/3} \mathfrak{M}^{-4/3} (\log N)^C.
   \end{equation}

   \noindent Recalling from Definition \ref{def:smoothed_log} that $\mathfrak{M} = T^{19/20}$ (and from \eqref{eqn:N1N2} that $\log N \le 100 \log T$), the right side of \eqref{eqn:final_expect_est} is at most $A T^{1/2-c}$. Together with \eqref{xi0stintegral} and \eqref{sumchif}, this establishes the lemma.
\end{proof}

We can now quickly establish Lemma \ref{lem:xi_expectation_outline}.
using Lemma \ref{lem:xi_expectation} and the following identity, originally shown as \cite[Equations (6.20) and (6.21)]{Spo23} (though under our notation it appears in \cite{Agg25}). 

	\begin{lem}[{\cite[Lemma 3.5]{Agg25}}] 
		\label{vt} 
		
		We have $(  \theta^{-1} \cdot \bm{\varsigma_0^{\dr}} - \alpha^{-1} \cdot \bm{\mathrm{T}} \bm{\varrho}) \ve  = \varsigma_1$.
	
	\end{lem} 
    
\begin{proof}[Proof of Lemma \ref{lem:xi_expectation_outline}]

By the $F (\lambda) = \mathfrak{l} (\Lambda - \lambda)$ and $\mathfrak{q} = \ve (\Lambda)$ case of Lemma \ref{lem:xi_expectation}, which is quickly verified to satisfy Assumption \ref{ass:F} at $A = (\log N)^3$ (and $|\mathfrak{q}| \le (\log N)^2$, by the second statement of Lemma \ref{lem:ve_tail}), it suffices to show that 
\begin{flalign}
\label{integralexpectation0}
\Bigg| \alpha^{-1} \displaystyle\int_{-\infty}^{\infty} \mathfrak{l} (\Lambda - \lambda) \cdot (\ve (\Lambda) - \ve (\lambda)) \varrho (\lambda) d \lambda - \big( \Lambda - \ve (\Lambda) \big) \Bigg| \le N^{-1}.
\end{flalign}

\noindent To do so, first observe that 
\begin{flalign}
\label{integralexpectation1}
\begin{aligned} 
\Bigg| \displaystyle\int_{-\infty}^{\infty} &  \mathfrak{l} (\Lambda - \lambda) \cdot (\ve (\Lambda) - \ve (\lambda)) \varrho (\lambda) d \lambda 
- 2 \displaystyle\int_{-\infty}^{\infty} \log |\Lambda - \lambda| \cdot (\ve (\Lambda) - \ve (\lambda)) \varrho (\lambda) d \lambda \Bigg| \\
& \le C \displaystyle\int_{-\infty}^{\infty} \big( e^{-(\log N)^2} + (\log N)^3 \cdot \mathbbm{1}_{|x| \le e^{-(\log N)^2}} \big) \cdot (\log N + \log (|\lambda|+1)) \varrho (\lambda) d \lambda \le N^{-1}.
\end{aligned} 
\end{flalign}

\noindent Here, in the first statement we used the facts that $|\mathfrak{l}(x) -  \log |x|| \le e^{-(\log N)^2} + (\log N)^3 \cdot \mathbbm{1}_{|x| \le e^{-(\log N)^2}}$ (by the definition \eqref{eqn:sl} of $\mathfrak{l}$) and that $|\ve (\Lambda)| + |\ve (\lambda)| \le C \log N + C \log (|\lambda|+1)$ (by the second statement in Lemma \ref{lem:ve_tail}, with the fact that $|\Lambda| \le \log N$). In the second statement, we estimated the integral using the bound on $\varrho$ from \eqref{lem:varrho_bd}. Next, we have 
\begin{flalign*}
 2 \alpha^{-1} \displaystyle\int_{-\infty}^{\infty} \log |\Lambda - \lambda| \cdot (\ve (\Lambda) - \ve (\lambda)) \varrho (\lambda) d \lambda & = \alpha^{-1} (\mathbf{T} \varrho (\Lambda) - \mathbf{T} \bm{\varrho}) \ve (\Lambda) \\
 & = (  \theta^{-1} \cdot \bm{\varsigma_0^{\dr}} - 1 - \alpha^{-1} \cdot \bm{\mathrm{T}} \bm{\varrho}) \ve(\Lambda),
\end{flalign*}

\noindent where the first statement follows from the definition \eqref{operatort} of $\mathbf{T}$ and the second from the second statement in \ref{rho2}. This, together with Lemma \ref{vt} and \eqref{integralexpectation1} yields \eqref{integralexpectation0}, and thus the lemma. 
\end{proof} 

\subsection{Approximations by sums of the $\lambda_j^{(r)}$} \label{subsec:ind_sum}

 In this section, analyze sums of the form
\begin{equation}\label{eqn:H_prelim}
\sum_{i=N_1}^{N_2}F(\Lambda_i) \cdot
   \big( \chi\left(\alpha (k- i)  + q_1\right) - \chi\left(\alpha (k- i)  + q_2- t \ve(\Lambda_i) \big) \right)  -\mathbb{E}[\cdots] .
\end{equation}

\noindent Instead of working with expressions of the form \eqref{eqn:H_prelim}, in order to lighten notation (and since we will require a generalization of these functions in Section \ref{sec:charge_fluct} below), we consider a slightly more general two variable function $H(\lambda, i)$ which satisfies the following assumptions.

\begin{ass}\label{ass:Hass}

Let $A>0$ be a real number, and $H: \mathbb{R} \times \llbracket N_1, N_2 \rrbracket \rightarrow \mathbb{R}$ be a function satisfying the following conditions.
\begin{enumerate}
\item For each $i \in \llbracket N_1, N_2 \rrbracket$, the function $\lambda \mapsto H(\lambda,i)$ satisfies Assumption \ref{ass:F}.\label{item:H1}
\item For all $|\lambda| \le \log N$ and $i, j \in \llbracket N_1, N_2 \rrbracket$, we have $|H(\lambda,i) -H(\lambda,j)| \leq A\mathfrak{M}^{-1} |i-j|$. 
\label{item:H2}
\item There exists some $k_0 \in \llbracket -T^{3/2} , T^{3/2}  \rrbracket$ such that the following holds. For all $|\lambda| \le \log N$ and $i \in \llbracket N_1, N_2 \rrbracket$ with $|i-k_0| >  T (\log N)^{12}$, we have $H(\lambda,i) = 0$.  \label{item:H3}
\end{enumerate}

\end{ass}

We will see that the summand in \eqref{eqn:H_prelim} will serve as an example of an $H$ satisfying Assumption \ref{ass:Hass}. The next lemma estimates sums over $H(\Lambda_i,i)$ through the resolvent of $\mathbf{L}$. Since its proof is similar to that of \cite[Lemma D.3]{Agg25}, we defer it to Appendix \ref{app:res_lem}.

\begin{lem}\label{lem:resolvent_lemma} 

There exists a constant $\mathfrak{c}>0$ such that the following holds. For any $z \in \mathbb{C} \setminus \eig \mathbf{L}$, set $\G(z) = (\L- z \cdot \Id)^{-1}$. Fix a real number $A > 0$ and a function $H : \mathbb{R} \times \llbracket N_1, N_2 \rrbracket$ satisfying Assumption \ref{ass:Hass}. Denoting $\eta = e^{-(\log N)^3}$, we have with overwhelming probability that 
    \begin{equation}\label{eqn:Hsum}
       \Bigg| \sum_{i=N_1}^{N_2} H(\Lambda_i, i) 
        - \frac{1}{\pi} \sum_{r=\lceil N_1/R \rceil}^{ \floor{N_2/R}-1} \int_{-\log N}^{\log N} 
 H(E, r R) \sum_{s=0}^{R-1} \Imaginary G_{rR+s,rR+s}(E + \i \eta) d E \Bigg| \le A T^{1/2-\mathfrak{c}}.
    \end{equation}
    
\end{lem}

The next lemma approximates the fluctuations of the integral appearing in \eqref{eqn:Hsum} by linear combinations of the $(\lambda_j^{(r)})$ from Definition \ref{def:Lr}, analogously to in \eqref{eqn:psi2Zktauapprox} and \eqref{eqn:lambdam2_approx}.

\begin{lem}\label{lem:Psi2_final_approx}
Adopt the notation and assumptions of Lemma \ref{lem:resolvent_lemma}. There exists a constant $\mathfrak{c}>0$ such that, with overwhelming probability, we have (recalling the $(\lambda_i^{(r)})$ from Definition \ref{def:Lr}) that
 \begin{flalign*}
    \frac{1}{\pi} & \sum_{r=\lceil N_1/R \rceil }^{\floor{N_2/R}-1} \int_{-\log N}^{\log N}  H(E, r R) \sum_{s=0}^{R-1} \Imaginary G_{rR+s,rR+s}(E + \i \eta)  d E \\
    & \qquad - \mathbb{E}\Bigg[\frac{1}{\pi} \sum_{r=\lceil N_1/R \rceil}^{\floor{N_2/R}-1}  \int_{-\log N}^{\log N}  H(E, r R) \sum_{s=0}^{R-1} \Imaginary G_{rR+s,rR+s}(E + \i \eta) d E \Bigg] \\
    & \qquad \qquad = \sum_{r=\lceil N_1/R \rceil}^{\floor{N_2/R}-1}  \sum_{s=0}^{R-1} H(  \lambda_s^{(r)}, r R) - \mathbb{E} \Bigg[ \sum_{r=\lceil N_1/R \rceil}^{\floor{N_2/R}-1}  \sum_{s=1}^{R} H(  \lambda_s^{(r)}, r R) \Bigg] + O(A T^{1/2-\mathfrak{c}}),
 \end{flalign*}
 where $O(A T^{1/2-\mathfrak{c}})$ is a real number satisfying $|O(A T^{1/2-\mathfrak{c}})| \le A T^{1/2-\mathfrak{c}}$.
\end{lem}

\begin{proof}

Throughout this proof, set $R' =  R+ 2 \floor{(\log N)^5}$. As in Definition \ref{def:Lr}, for each index $r \in \llbracket N_1/R, N_2/R-1 \rrbracket$ let $\tilde{\mathbf{L}}^{[r]} = [\tilde{L}_{ij}^{[r]}]$ denote the $R' \times R'$ submatrix of $\mathbf{L}$ whose rows and columns are indexed by $i,j \in \llbracket r R - (\log N)^5, (r+1) R + (\log N)^5-1 \rrbracket$ (so that $\tilde{L}_{ij}^{[r]} = L_{ij}$ for each such $(i,j)$). Observe that $\tilde{\mathbf{L}}^{[r]}$ is sampled under the thermal equilibrium \eqref{eqn:equil}, with the $N$ there equal to $R'$ here. For any complex number $z \in \mathbb{C} \setminus \eig \tilde{\mathbf{L}}^{[r]}$, further let $\tilde{\G}^{[r]}(z)= [\tilde G_{i j}^{[r]}] = (\tilde{\mathbf{L}}^{[r]} - z)^{-1}$ denote the resolvent of this matrix. Then, with overwhelming probability, we have
\begin{multline}\label{eqn:psik4equiv}
  \frac{1}{\pi} \sum_{r=\lceil N_1/R \rceil}^{\floor{N_2/R}-1} \int_{-\log N}^{\log N} 
 H(E, r R) \sum_{s=0}^{R-1} \Imaginary G_{rR+s,rR+s}(E + \i \eta) d E
   \\
   = 
   \frac{1}{\pi} \sum_{r=\lceil N_1/R \rceil}^{\floor{N_2/R}-1} \int_{-\log N}^{\log N} 
 H(E, r R) \sum_{s=0}^{R-1} \Imaginary \tilde G_{rR+s,rR+s}^{[r]}(E + \i \eta) d E + O( A T^{1/2-\mathfrak{c}}),
 \end{multline}

\noindent by $\delta=0$ case of Lemma \ref{lem:resolv_coupling} (with the $\mathcal{D}$ there equal to $\llbracket N_1, N_2 \rrbracket \setminus \llbracket rR-(\log N)^5, (r+1)R + (\log N)^5 \rrbracket$ here, observing that $\dist (rR+s, \mathcal{D}) \ge (\log N)^5$ whenever $s \in \llbracket 0, R-1 \rrbracket$, and recalling $\eta = e^{-(\log N)^3}$).

    We next apply Lemma \ref{lem:res_linst}, to which end we must first set some notation. For each $r \in \llbracket N_1/R, N_2/R-1\rrbracket$ and $\lambda \in \mathbb{R}$, set 
    \begin{flalign*} 
    \mathcal{H}_r(\lambda) = H(\lambda, r R) \cdot \mathbbm{1}_{|H(\lambda, r R)| \leq A} + A \cdot \mathbbm{1}_{H(\lambda, r R)> A} - A \cdot \mathbbm{1}_{H(\lambda, r R)<-A}.
    \end{flalign*}

    \noindent For each $r \in \llbracket N_1/R, N_2/R - 1 \rrbracket$, denote $\eig \tilde{\mathbf{L}}^{[r]} = (\tilde{\lambda}_1^{(r)}, \tilde{\lambda}_2^{(r)}, \ldots , \tilde{\lambda}_{R'}^{(r)})$, and let $\tilde{\varphi}^{[r]} : \llbracket 1, R' \rrbracket \rightarrow \llbracket r R -\floor{(\log N)^5}, (r+1) R+\floor{(\log N)^5} - 1 \rrbracket$ denote a $\zeta'$-localization center bijection for $\tilde{\mathbf{L}}^{[r]}$, where $\zeta' = \zeta_{R'} = e^{-100 (\log R')^{3/2}}$. For each $s \in \llbracket r R -\floor{(\log N)^5}, (r+1) R+\floor{(\log N)^5} - 1\rrbracket$, denote by
    \begin{equation}
        \tilde{\Lambda}_{s}^{(r)} \coloneqq \tilde{\lambda}_{(\tilde{\varphi}^{[r]})^{-1}(s)}^{(r)}
    \end{equation}
      the eigenvalue of $\tilde{\mathbf{L}}^{[r]}$ with localization center $s$. The first part of Assumption \ref{ass:Hass} gives $\mathcal{H}_r(\lambda) = H(\lambda, r R)$ for $|\lambda| \leq \log N$, on the overwhelmingly probable event $\mathsf{BND}_{\tilde{\mathbf{L}}^{[r]}}(\log N)$ (by Lemma \ref{lem:bd_lem}), we have $\mathcal{H}_r(\tilde{\Lambda}_{s}^{(r)}) = H(\tilde{\Lambda}_{s}^{(r)}, r R)$ for each $s \in \llbracket r R -\floor{(\log N)^5}, (r+1) R+\floor{(\log N)^5} - 1\rrbracket$. 
      
      Now, let us apply \eqref{eqn:Hlr} from Lemma \ref{lem:res_linst} with the $\mathbf{L}$ there equal to $\tilde{\mathbf{L}}^{[r]}$ here; the $(n_1,n_2)$ there equal to $(r R - \lfloor (\log N)^5 \rfloor, (r+1) R + \lfloor (\log N)^5 \rfloor -1)$; and the $H$ there given by $\mathcal{H}_r$ here. Then denoting
    \begin{flalign}\label{eqn:LHS_resolv_toeig}
    \mathcal{E}_{r,1} = \frac{1}{\pi}   \int_{-\log N}^{\log N} 
 H(E, r R) \sum_{s=0}^{R-1} \Imaginary \tilde G_{rR+s,rR+s}^{(r)}(E + \i \eta) d E  -
     \sum_{s=- \floor{(\log R')^5}}^{R+\floor{(\log R')^5}-1} H(\tilde \Lambda_{rR+s}^{(r)}, r R),
    \end{flalign}
    
     \noindent for each $r \in \llbracket N_1/R, N_2/R-1 \rrbracket$, we have with overwhelming probability that each $\mathcal{E}_{r,1}$ satisfies
    \begin{equation}\label{eqn:ER1bd}
    |\mathcal{E}_{r,1}| \leq 6 A (\log N)^5.
    \end{equation}
      
 Next, in view of \eqref{eqn:psik4equiv} and \eqref{eqn:ER1bd}, we would like to replace the eigenvalues $\tilde{\Lambda}_j^{(r)}$ of $\tilde{\mathbf{L}}^{[r]}$ appearing in \eqref{eqn:LHS_resolv_toeig} with those $\lambda_j^{(r)}$ of $\mathbf{L}^{[r]}$. To that end, define $  \{\lambda_s^{(r)}\}_{s=1,\dots, R}$, as the nonzero eigenvalues of the matrix $\L^{[r]}$. In addition, let $\varphi^{(r)} : \llbracket 1, R \rrbracket \rightarrow \llbracket r R, (r+1)R -1\rrbracket$ denote a localization center bijection for $\L^{[r]}$, and define $\Lambda_{r R+s}^{(r)} = \lambda_{(\varphi^{(r)})^{-1}(rR+s)}^{(r)}$ for each $s \in \llbracket 0, R-1 \rrbracket$.

 By Lemma \ref{lem:Lax_eig_coupling} (applied with the $(\mathbf{L}, \tilde{\mathbf{L}})$ there equal to $(\tilde{\mathbf{L}}^{[r]}, \mathbf{L}^{[r]})$ here), for each $r \in \llbracket N_1/R, N_2/R-1 \rrbracket$ there exists an overwhelmingly probable (by Lemma \ref{lem:Lax_eig_coupling}) event $\mathsf{E}(r)$ on which the following holds. There is an injective map $\psi : \llbracket r R, (r+1)R-1 \rrbracket \rightarrow \llbracket r R- \floor{(\log N)^5}, (r+1)R+ \floor{(\log N)^5}-1 \rrbracket$, such that for any $s \in \llbracket (\log N)^3, R-(\log N)^3 \rrbracket$ we have 
 \begin{flalign} 
 \label{lambdalambda} 
 |\Lambda_{rR +s}^{(r)}-\tilde \Lambda_{\psi(rR +s)}^{(r)}| \leq e^{-c(\log R)^3} \leq e^{-c' (\log N)^3}.
\end{flalign}

\noindent For each $r \in \llbracket N_1/R,N_2/R-1\rrbracket $, denote
\begin{equation}\label{eqn:Er2}
    \mathcal{E}_{r,2} \coloneqq  \sum_{s=- \floor{(\log R')^5}}^{R+\floor{(\log R')^5}-1} H(\tilde \Lambda_{r R+s}^{(r)}, r R)  -
  \sum_{s=0}^{R-1} H(\Lambda_{r R+s}^{(r)}, r R).
\end{equation}

On the event $\mathsf{E}(r) $, we have for each $r \in \llbracket N_1/R,N_2/R-1\rrbracket $,
\begin{multline}\label{eqn:ER2bd}
   \left|\mathcal{E}_{r,2}  \right| 
  \leq \sum_{s=2 \floor{(\log N)^3}}^{R-2 \floor{(\log N)^3}} \left| H(\tilde \Lambda_{\psi(r R + s)}^{(r)}, r R)  -
 H(\Lambda_{r R+s}^{(r)}, r R)   \right| + O(A (\log N)^5) \leq C A (\log N)^5,
\end{multline}

\noindent where to obtain the first inequality we used \eqref{estimatef} and to obtain the second we used \eqref{eqn:F_lip} and \eqref{lambdalambda} (both with Item \ref{item:H1} of Assumption \ref{ass:Hass}).

 Next, for each $r \in \llbracket N_1/R , N_2/R-1 \rrbracket$, define the event 
 \begin{equation}\label{eqn:etildr}
 \tilde{\mathsf{E}}(r) \coloneqq  \mathsf{BND}_{\tilde{\mathbf{L}}^{[r]}}(\log N) \cap \mathsf{BND}_{\mathbf{L}^{[r]}}(\log N) \cap \{|\mathcal{E}_{r,1}| + |\mathcal{E}_{r,2}| \leq   A (\log N)^6 \}.
 \end{equation}
 In what follows, we restrict to the overwhelmingly probable (by Lemma \ref{lem:bd_lem}, \eqref{eqn:ER1bd}, and \eqref{eqn:ER2bd}) event 
\begin{flalign}\label{eqn:bigevent}
\mathsf{E} = \bigcap_{r \in \llbracket N_1/R , N_2/R-1 \rrbracket}   \tilde{\mathsf{E}}(r).
\end{flalign} 

\noindent By definition, we have 
    \begin{flalign}\label{eqn:rsumss}
    \begin{aligned} 
\frac{1}{\pi} \sum_{r=\lceil N_1/R \rceil}^{\floor{N_2/R}-1} \int_{-\log N}^{\log N} 
 H(E, r R &) \sum_{s=0}^{R-1} \Imaginary \tilde G_{rR+s,rR+s}^{[r]}(E + \i \eta) d E  = \sum_{r=\lceil N_1/R \rceil}^{\floor{N_2/R}-1} \left( \sum_{s=0}^{R-1} H( \Lambda_s^{(r)}, r R) + \mathcal{E}_r \right),
 \end{aligned} 
\end{flalign}

\noindent where we have denoted (recalling $\mathcal{E}_{r,1}$ and $\mathcal{E}_{r,2}$ from \eqref{eqn:LHS_resolv_toeig} and \eqref{eqn:Er2}, respectively)
\begin{equation}\label{eqn:ERdf}
    \mathcal{E}_r = \mathcal{E}_{r,1} + \mathcal{E}_{r,2}.
\end{equation}

 \noindent By the definition of the event $\tilde{\mathsf{E}}(r)$ in \eqref{eqn:etildr}, the error $\mathcal{E}_r$ satisfies
    \begin{flalign} \label{eqn:Erbd}
    |\tilde{\mathcal{E}}_r| \leq  A(\log N)^6, \qquad \text{for each $\lceil N_1/R \rceil  \le r \le \floor{N_2/R}-1$}, \qquad \text{where $\tilde{\mathcal{E}}_r=  \mathcal{E}_r  \cdot \mathbbm{1}_{\tilde{\mathsf{E}}(r)}$}.
    \end{flalign} 

   \noindent It is also quickly verified (by closely following the reasoning in \eqref{expectationH} and recalling the event $\tilde{\mathsf{E}}(r)$ is overwhelmingly probable) that 
   \begin{equation}\label{eqn:expec_tb}
    \mathbb{E} [||\mathcal{E}_{r} - \tilde{\mathcal{E}}_r|] = \mathbb{E}[|\mathcal{E}_{r}| \cdot  \mathbbm{1}_{\tilde{\mathsf{E}}(r)^{\complement}}] \leq A R e^{-c (\log N)^2}.
   \end{equation}

    \noindent  By \eqref{eqn:rsumss} and \eqref{eqn:expec_tb}, to establish the lemma, it suffices to show that, with overwhelming probability, 
   \begin{flalign}
   \label{sume2}
       \Bigg| \displaystyle\sum_{r=\lceil N_1/R \rceil}^{\lfloor N_2/R \rfloor - 1} (\tilde{\mathcal{E}}_r - \mathbb{E}[\tilde{\mathcal{E}}_r]) \Bigg| \le AT^{1/2-c}.
   \end{flalign}

    Next recall from \eqref{eqn:bigevent} that on $\mathsf{E}$ we have for all $r \in \llbracket N_1/R, N_2/R-1 \rrbracket$ that $|\lambda_s^{(r)}| \leq \log N$ for each $s \in \llbracket 1, R \rrbracket$ and that $|\tilde{\lambda}_s^{(r)}| \leq \log N$ for each $s \in \llbracket 1, R' \rrbracket$. Thus, whenever $r \in \llbracket N_1/R, N_2/R-1 \rrbracket$ satisfies $|rR-k_0| \le T(\log N)^{12}$, we have from Item \ref{item:H3} of Assumption \ref{ass:Hass} (with \eqref{eqn:LHS_resolv_toeig} and \eqref{eqn:Er2}) that $\mathcal{E}_{r,1} = \mathcal{E}_{r,2} = 0$ and thus that $\tilde{\mathcal{E}}_r = \mathcal{E}_r  = 0$ if $|r R - k_0| > T (\log N)^{12}$. Therefore, to show \eqref{sume2} (and thus the lemma), it suffices to verify that 
    \begin{flalign} 
    \label{sumrek0}
\Bigg| \sum_{r: |rR - k_0| \leq T (\log N)^{12}} (\tilde{\mathcal{E}}_r - \mathbb{E}[\tilde{\mathcal{E}}_r]) \Bigg| \le AT^{1/2-c}.
    \end{flalign}

    Now observe that $\tilde{\mathbf{L}}^{[r]}$ and $\tilde{\mathbf{L}}^{[r']}$ are independent whenever $|r-r'| \ge 2$, as they involve disjoint subsets of the entries of $\mathbf{L}$. Therefore, $\tilde{\mathcal{E}}_r$ and $\tilde{\mathcal{E}}_{r'}$ are independent whenever $|r-r'| \ge 2$, by their definitions \eqref{eqn:Erbd}, \eqref{eqn:ERdf}, \eqref{eqn:LHS_resolv_toeig}, and \eqref{eqn:Er2}. This, together with the bound $|\tilde{\mathcal{E}}_r| \le A (\log N)^C$ and a Chernoff concentration bound, it follows that the sum of $\tilde{\mathcal{E}}_r - \mathbb{E}[\tilde{\mathcal{E}}_r]$ over even indices $r$ is not too large with high probability, namely,
    \begin{equation}
    \label{sumep} 
        \mathbb{P}\Bigg(\bigg|\sum_{p:|2pR-k_0| < T (\log N)^{12}} (\tilde{\mathcal{E}}_{2p} - \mathbb{E}[\tilde{\mathcal{E}}_{2p}])\bigg| > A (\log N)^C (TR^{-1})^{1/2} \Bigg) \leq c^{-1} e^{-c (\log N)^2}.
    \end{equation}
    
    \noindent Similarly bounding with overwhelming probability the sum of $\tilde{\mathcal{E}}_r - \mathbb{E}[\tilde{\mathcal{E}}_r]$ over odd indices $r$, summing with \eqref{sumep}, and recalling from \eqref{ra0} that $R = \lfloor T^{1/5} \rfloor$, we deduce \eqref{sumrek0} and thus the lemma. 
\end{proof}

\subsection{Expressions for Lax matrix functionals as independent sums}

\label{SumGeneral} 

In this section we show the following theorem, which will quickly imply Proposition \ref{prop:equiv_expr_cor} as a special case. Its brief proof combines the earlier lemmas of this section. In the statement of the theorem, we recall the $X_i$ from Definition \ref{xi} and $\{\lambda_s^{(r)}\}_{s=1}^R$ from Definition \ref{def:Lr}.

  \begin{thm}\label{thm:Hindsum_thm}
  
    There exists a constant $\mathfrak{c}>0$ such that the following two statements hold, for any real number $A > 0$ and function $F: \mathbb{R} \rightarrow \mathbb{R}$ satisfying Assumption \ref{ass:F}.  
  
  \begin{enumerate}
      \item  Fix an integer $k \in \llbracket -T (\log N)^{10}, T (\log N)^{10} \rrbracket$ and real numbers $\mathfrak{q}_1, \mathfrak{q}_2 \in [-T(\log N)^4, T(\log N)^4]$ and $t \in [0, T \log N]$. Define the function $\tilde{H}: \mathbb{R}^2 \rightarrow \mathbb{R}$ by setting
    \begin{equation}\label{eqn:Hf}
         \tilde{H}(\lambda, q) \coloneqq F(\lambda) \cdot \big( \chi(\mathfrak{q}_1 - q ) - \chi(\mathfrak{q}_2 - q - t \ve(\lambda)) \big),
     \end{equation}

     \noindent for each $\lambda, q \in \mathbb{R}$. Then, with overwhelming probability, we have  
\begin{flalign}\label{eqn:xiH_decomp0}
\begin{aligned} 
\Bigg|  \sum_{i = N_1}^{N_2} & \tilde{H}(\Lambda_i, q_i(0)-q_k(0) ) - \mathbb{E}\left[ \sum_{i = N_1}^{N_2} \tilde{H}(\Lambda_i, q_i(0)-q_k(0)) \right] \\
 & \qquad - \bigg(\sum_{i = N_1}^{N_2} \tilde{H}(\Lambda_i, \alpha (i- k )) - \mathbb{E} \Big[ \sum_{i = N_1}^{N_2} \tilde{H}(\Lambda_i, \alpha (i- k )) \Big] \bigg) \\
 & \qquad + \alpha^{-1} \sum_{r = \lceil N_1/R \rceil}^{\lfloor N_2/R \rfloor - 1} \int_{-\infty}^{\infty} \tilde{H}(\lambda, \alpha (r R-k)  ) \varrho(\lambda) d\lambda \sum_{i = r R}^{(r+1) R-1} X_{i}  \Bigg| \leq A T^{1/2-\mathfrak{c}},
 \end{aligned} 
 \end{flalign}
and 
\begin{flalign}\label{eqn:xiH_decomp2}
\begin{aligned} 
\Bigg|   \sum_{i = N_1}^{N_2} \tilde{H}(\Lambda_i, \alpha (i- k )) - \mathbb{E}\bigg[ & \sum_{i = N_1}^{N_2} \tilde{H}(\Lambda_i, \alpha (i- k )) \bigg] - \sum_{r = \lceil N_1/R \rceil}^{\lfloor N_2/R  \rfloor - 1}  \sum_{s=1}^{R} \tilde{H}(\lambda_s^{(r)}, \alpha (r R- k) )  \\
   &  + \mathbb{E} \bigg[\sum_{r = \lceil N_1/R \rceil}^{\lfloor N_2/R \rfloor - 1} \sum_{s=1}^{R} \tilde{H}(\lambda_s^{(r)}, \alpha (r R- k) )  \bigg]  \Bigg|\leq A T^{1/2-\mathfrak{c}}.
    \end{aligned} 
\end{flalign}

\label{item:isth1}

\item Fix an integer $k \in \llbracket -T^{3/2}, T^{3/2} \rrbracket$; real numbers $B > 0$ and $U \in \mathbb{R}$ satisfying $\mathfrak{M} \leq U \leq T \log N$, and  $t \in [0, T(\log N)^{10}]$; and a function $G : \mathbb{R} \rightarrow \mathbb{R}$ of compact support, satisfying Assumption \ref{ass:G2}. Set $H(\lambda,q) \coloneqq F(\lambda) \cdot G(q)$ for each $(\lambda, q) \in \mathbb{R}^2$. Then, with overwhelming probability, we have
\begin{flalign}\label{eqn:xiH_decomp00}
\begin{aligned}
 \Bigg| & \sum_{i = N_1}^{N_2} H(\Lambda_i, q_i(0)-q_k(0) + t \ve(\Lambda_i) ) - \mathbb{E}\bigg[ \sum_{i = N_1}^{N_2} H(\Lambda_i, q_i(0)-q_k(0) + t \ve(\Lambda_i)) \bigg]  \\
 &  - \Bigg( \sum_{i = N_1}^{N_2} H(\Lambda_i, \alpha (i -k) + t \ve(\Lambda_i) ) - \mathbb{E}\bigg[ \sum_{i = N_1}^{N_2} H(\Lambda_i, \alpha (i -k) + t \ve(\Lambda_i)) \bigg] \Bigg)
 \\
&  + \alpha^{-1} \sum_{r = \lceil N_1/R \rceil}^{\lfloor N_2/R \rfloor - 1}  \int_{-\infty}^{\infty} H(\lambda, \alpha (r R-k) + t \ve(\lambda) ) \varrho(\lambda) d\lambda \sum_{i = r R}^{(r+1) R-1} X_{i}  \Bigg| \leq A B T^{1/2-\mathfrak{c}},
\end{aligned} 
  \end{flalign}
  and
  \begin{flalign}\label{eqn:xiH_decomp22}
  \begin{aligned}
\Bigg| \sum_{i = N_1}^{N_2} H( & \Lambda_i, \alpha (i -k) + t \ve(\Lambda_i) ) - \mathbb{E}\bigg[ \sum_{i = N_1}^{N_2} H(\Lambda_i, \alpha (i -k) + t \ve(\Lambda_i)) \bigg]  \\ 
 & - \sum_{r = \lceil N_1/R \rceil}^{\lfloor N_2/R \rfloor - 1}  \sum_{s=1}^{R} H \big(\lambda_s^{(r)}, \alpha (r R- k)  + t \ve(\lambda_s^{(r)}) \big) \\
 & + \mathbb{E}\bigg[\sum_{r = \lceil N_1/R \rceil}^{\lfloor N_2/R \rfloor - 1}  \sum_{s=1}^{R} H \big(\lambda_s^{(r)}, \alpha (r R- k) + t \ve(\lambda_s^{(r)}) \big) \bigg] \Bigg| \le  A B T^{1/2-\mathfrak{c}}.
 \end{aligned} 
 \end{flalign}
 \label{item:isth2}
  \end{enumerate}
  
 \end{thm}

 \begin{proof}

We begin by addressing the first case of the theorem. To show \eqref{eqn:xiH_decomp0}, we apply Lemma \ref{lem:HBM_lem} (twice) to approximate 
 \begin{align}
 &\sum_{i = N_1}^{N_2}F(\Lambda_i) \cdot \chi\left(q_k(0)-q_i(0)  + \mathfrak{q}_1\right); \qquad \sum_{i = N_1}^{N_2}F(\Lambda_i) \cdot \chi\left(q_k(0)-q_i(0) + \mathfrak{q}_2- t \ve(\Lambda_i) \right),
 \end{align}

\noindent whose difference constitutes the sums on the left hand side of \eqref{eqn:xiH_decomp0} (so in one application $t$ there is given by $0$, and in another it matches $t$ here). Combining these two applications, we obtain \eqref{eqn:xiH_decomp0}. To show \eqref{eqn:xiH_decomp2}, define for each $(\lambda,i) \in \mathbb{R} \times \mathbb{Z}$ the quantity 
\begin{equation*}
H_0(\lambda, i) \coloneqq F(\lambda)\cdot  \big( \chi ( \mathfrak{q}_1 + \alpha(k- i)) - \chi( \mathfrak{q}_2 + \alpha(k- i) - t \ve(\lambda)) \big).
\end{equation*}

\noindent Observe by Definition \ref{def:chi_def} (and Lemma \ref{lem:ve_tail}) that $H_0$ satisfies the properties of Assumption \ref{ass:Hass} (with the $(A,k_0)$ there equal to $(10A, k)$ here). Since the first two sums in \eqref{eqn:xiH_decomp2} equal $\sum_{i=N_1}^{N_2} H_0(\lambda, i) - \mathbb{E}[\sum_{i=N_1}^{N_2} H_0(\lambda, i)]$, \eqref{eqn:xiH_decomp2} follows from applying Lemmas \ref{lem:resolvent_lemma} and \ref{lem:Psi2_final_approx}.

We next address the second case of the theorem. The bound \eqref{eqn:xiH_decomp00} follows from Lemma \ref{lem:HBM_lem}, together with the fact that $H(\lambda, \infty) = H(\lambda, -\infty) = 0$ (which holds since $G$ has compact support). 

To show \eqref{eqn:xiH_decomp22}, for each $(\lambda,i) \in \mathbb{R} \times \mathbb{Z}$, define $H_1(\lambda, i) \coloneqq H(\lambda, \alpha (i - k)+ t \ve(\lambda))$, which by Assumption \ref{ass:G2} (and Lemma \ref{lem:ve_tail}) satisfies Assumption \ref{ass:Hass} (with the $(A,k_0)$ there given by $(AB,k)$ here). Since the first two sums in \eqref{eqn:xiH_decomp22} equal $\sum_{i=N_1}^{N_2}H_1(\Lambda_i, i) - \mathbb{E}[\sum_{i=N_1}^{N_2}H_1(\Lambda_i, i)]$,
     \eqref{eqn:xiH_decomp22} follows from applying Lemmas \ref{lem:resolvent_lemma} and \ref{lem:Psi2_final_approx}.
 \end{proof}

\begin{proof}[Proof of \Cref{prop:equiv_expr_cor}]
    Both of these families of observables $\Xi$ and $\Xi^{[m]}$ both fit into the framework of (part \ref{item:isth1} of) Theorem \ref{thm:Hindsum_thm}, so long as $\Lambda \in [-\log N, \log N]$ and $m \leq (\log N)^{1/10}$. Indeed, with these assumptions, $\lambda \mapsto 2 \mathfrak{l}(\Lambda- \lambda)$ satisfies Assumption \ref{ass:F} with $A = 4 (\log N)^2$, and $\lambda \mapsto \lambda^m$ satisfies Assumption \ref{ass:F} with $A = (\log N)^m$. Applying the theorem, and noting that in the current proposition we divide all observables by $T^{1/2}$, completes the proof.
\end{proof}

\section{Effective convergence of linear statistics}
\label{sec:effective_conv}

In this section we establish an effective estimate, given by Proposition \ref{ass:phiN}, for the limiting variance of linear statistics of Lax matrices under from thermal equilibrium. Throughout this section, we let $\mathbf{L} = [L_{ij}]$ denote an $N \times N$ Lax matrix (Definition \ref{lt}), with rows and columns indexed by $i,j \in \llbracket 1,N \rrbracket$, sampled from thermal equilibrium (Definition \ref{def:inf_thermal_eq}). Throughout this section, neither Assumption \ref{ass:NT_assumption} nor the localization centers of Definition \ref{def:loccent} will be used; thus, without loss of generality we have shifted set $\llbracket N_1, N_2 \rrbracket$ indexing the rows and columns of $\mathbf{L}$ to $\llbracket 1, N \rrbracket$. Further define the stretch random variables $X_i$ as in Definition \ref{xi}, where the $q_i$ there are determined from the $a_i = L_{i,i+1}$ here by \eqref{abr}. Also let $\eig \mathbf{L} = (\lambda_1, \lambda_2, \ldots , \lambda_N)$. Throughout, we further recall the notation from Definitions \ref{def:dress_intro} and \ref{def:Cdef}, which will be often used. To ease notation, throughout the section we fix the constant 
\begin{flalign}
\label{gamma} 
\gamma = 3 / 5.
\end{flalign}

\subsection{Variance of linear statistic with effective error} \label{subsec:lin_var_err}

Under the above notation, the following proposition, to be shown in Section \ref{subsec:limvar_proof}, states an effective error estimate\footnote{The proof of Proposition \ref{prop:limvar} will indicate that the $\mathfrak{c}$ there can be taken to be any constant less than $1$, but we will not need this here.} for variances of certain linear statistics of the random Lax matrix $\mathbf{L}$.

\begin{prop}[Effective variance approximation]
\label{prop:limvar}

There exists a constant $\mathfrak{c}>0$ such that the following holds. Set $M = e^{\sqrt{\log N}}$, and let $\phi  : \mathbb{R} \rightarrow \mathbb{R}$ be a differentiable function satisfying, for each $x \in \mathbb{R}$,
\begin{flalign} \label{eqn:limvar_phiasmps}
| \phi (x) | \cdot \mathbbm{1}_{|x| \le (\log N)^{\gamma}} \leq M; \qquad |\phi(x)| \leq M e^{ |x|^{3/2} }; \qquad |\phi'(x)| \leq M e^{10^3 (\log N)^2+|x|^{3/2} }.
\end{flalign} 

\noindent Then, we have 
\begin{equation}\label{eqn:limvar_Xi}
\displaystyle\sup_{|a| \le M} \Bigg| \frac{1}{N} \var \bigg(   \sum_{i=1}^{N} \big( \phi(\lambda_i) + a  X_i \big) \bigg)
-\mathscr{C} \big(\phi + (a-\mu_{\phi}) \varsigma_0, \phi + (a-\mu_{\phi}) \varsigma_0 \big) \Bigg| \leq  N^{-\mathfrak{c}}.
\end{equation} 

\end{prop}

\subsection{Estimates on regularized linear statistics using resolvents}
\label{subsec:limvar_prep}

In this subsection we use resolvent estimates to approximate the variances of linear statistics of a Lax matrix of some size $N$ by those of one of some size $\mathfrak{N} \ge N$. Throughout, for any function $\phi : \mathbb{R} \rightarrow \mathbb{R}$ and real numbers $\eta > 0$ and $x \in \mathbb{R}$, we denote (recalling $\gamma = 3/5$ from \eqref{gamma})
\begin{equation}\label{eqn:phietn}
        \phi_{\eta; N}(x) \coloneqq \frac{1}{\pi}\int_{-(\log N)^{\gamma}}^{(\log N)^{\gamma}} \phi(E) \frac{\eta}{(x-E)^2 +\eta^2} d E.
\end{equation}

In what follows, we let $\mathfrak{N} \ge N$ be a positive integer. Let $\mathbf{L}_{\mathfrak{N}} = [L_{ij;\mathfrak{N}}]$ denote an $\mathfrak{N} \times \mathfrak{N}$ Lax matrix (Definition \ref{lt}) sampled under thermal equilibrium \eqref{eqn:equil}, with rows and columns indexed by $i,j \in \llbracket 1, N \rrbracket$. Set $\eig \mathbf{L}_{\mathfrak{N}} = (\lambda_{1;\mathfrak{N}}, \lambda_{2;\mathfrak{N}}, \ldots , \lambda_{\mathfrak{N};\mathfrak{N}})$. We also let $X_{i;\mathfrak{N}}$ denote the quantity $q_{i+1;\mathfrak{N}} - q_{i;\mathfrak{N}} - \alpha$ (as in \eqref{eqn:Xidef}), for $i \in \llbracket \mathfrak{N}_1, \mathfrak{N}_2 -1\rrbracket$; here, the quantities $q_{i+1;\mathfrak{N}} - q_{i;\mathfrak{N}} $ are defined via $a_{i;\mathfrak{N}} = e^{(q_{i;\mathfrak{N}}-q_{i+1;\mathfrak{N}})/2}$ as in \eqref{abr}, where $a_{i;\mathfrak{N}} = L_{i,i+1;\mathfrak{N}}$ are the off-diagonal entries of $\mathbf{L}_{\mathfrak{N}}$.

Throughout this subsection, we will work under the assumption below. Let us mention that the function $\phi$ in the below assumption (and in the results adopting it) will eventually be played not by the function $\phi$ in Proposition \ref{prop:limvar}, but instead by a certain polynomial approximation to it. Thus, the assumptions on $\phi$ below are slightly weaker than those imposed in Proposition \ref{prop:limvar}. 

\begin{ass}\label{ass:phiN}
Let $M \ge 1$ be a real number and $\phi : \mathbb{R} \rightarrow \mathbb{R}$ be a differentiable function such that, for each $x \in \mathbb{R}$, we have (recalling that $\gamma = 3/5$)
\begin{flalign} \label{eqn:phiasmps}
| \phi (x) | \cdot \mathbbm{1}_{|x| \le (\log N)^{\gamma}} \leq M; \qquad |\phi(x)| \leq M e^{10 |x|^{3/2} }; \qquad |\phi'(x)| \leq M e^{10^6 (\log N)^2+10 |x|^{3/2} }.
\end{flalign} 
\end{ass}

\begin{prop}\label{prop:eff_varbound}

There exists a constant $\mathfrak{C}>1$ such that the following holds. Under the above notation, assume that $\mathfrak{N}  \ge N \ge \mathfrak{C}$; let $\eta = e^{-(\log N)^3}$; and suppose that $\mathfrak{N}$ is divisible by $N$. Letting $M \ge 1$ be a real number and $\phi: \mathbb{R} \rightarrow \mathbb{R}$ denote a function satisfying Assumption \ref{ass:phiN}, we have for any $a \in [-M, M]$ that (recalling $\phi_{\eta;N}$ from \eqref{eqn:phietn})
\begin{equation}\label{eqn:stable_witheta} 
\Bigg| \frac{1}{N}\var \bigg(\sum_{i=1}^{N} \big( \phi_{\eta; N}(\lambda_{i;N}) + a X_{i;N} \big) \bigg) - \frac{1}{\mathfrak{N}} \var \bigg(\sum_{i=1}^{\mathfrak{N}} \big( \phi_{\eta;N}(\lambda_{i;\mathfrak{N}}) + a X_{i;\mathfrak{N}} \big) \bigg)  \Bigg| \leq M^2 N^{-1} (\log N)^{11}.
\end{equation}
\end{prop}

The proof essentially proceeds by viewing the $\mathfrak{N} \times \mathfrak{N}$ Lax matrix $\mathbf{L}_{\mathfrak{N}}$ as consisting of $N^{-1}\mathfrak{N}$ submatrices of size $N \times N$, each with the same law as $\mathbf{L}$. We will then show that two different such $N \times N$ Lax submatrices are approximately independent. This will enable us to approximate the variance of a linear statistic for $\mathbf{L}_{\mathfrak{N}}$ by $N^{-1} \mathfrak{N}$ multiplied by the corresponding variance for $\mathbf{L}$. 

Lemma \ref{lem:res_cov_bd} below makes precise the above notion of ``approximate independence,'' by estimating the covariances of distant resolvent entries of $\mathbf{L}_{\mathfrak{N}}$. Before stating it, we state a quick bound on the difference between covariances of random variables that are likely approximately equal.

\begin{lem}\label{lem:cb}

   Let $\varepsilon, \delta, M \geq 0$ be real numbers; $\mathsf{E}$ be an event; and $X_1$, $X_2$, $Y_1$, and $Y_2$ be random variables. Assume for each $i \in \{ 1, 2 \}$ that $\mathbb{E}[X_i^4 + Y_i^4] \leq M^4$ and $\mathbbm{1}_{\mathsf{E}} \cdot |X_i-Y_i| \leq M \varepsilon$, and that $\mathbb{P}(\mathsf{E}^{\complement}) \leq \delta$. Then, 
\begin{equation}\label{eqn:cb}
    |\Cov(X_1, X_2) - \Cov(Y_1, Y_2) |
    \leq 4 M^2 (\varepsilon +\delta^{1/4}) .
\end{equation}

\end{lem}
\begin{proof}

This lemma follows from the estimates
    \begin{flalign}\label{eqn:cov_decomp}
    \begin{aligned} 
    \Cov (X_1, X_2 & ) - \Cov(Y_1, Y_2) \\
    & = \Cov ( \mathbbm{1}_{\mathsf{E}} \cdot (X_1- Y_1), X_2) + \Cov (\mathbbm{1}_{\mathsf{E}^{\complement}} \cdot (X_1 - Y_1), X_2) \\
    & \qquad - \Cov(Y_1, \mathbbm{1}_{\mathsf{E}} \cdot (Y_2 - X_2)) - \Cov (Y_1, \mathbbm{1}_{\mathsf{E}^{\complement}} \cdot (Y_2 - X_2)) \\
    & \le \mathbb{E} [\mathbbm{1}_{\mathsf{E}} \cdot (X_1- Y_1)^2]^{1/2} \cdot \mathbb{E}[X_2^2]^{1/2}  + \mathbb{P} [\mathsf{E}^{\complement}]^{1/4}  \cdot \mathbb{E}[(X_1 - Y_1)^4]^{1/4} \cdot \mathbb{E}[X_2^2]^{1/2} \\
    & \qquad \mathbb{E}[\mathbbm{1}_{\mathsf{E}} \cdot (X_2- Y_2)^2]^{1/2} \cdot \mathbb{E}[Y_1^2]^{1/2}  + \mathbb{P} [\mathsf{E}^{\complement}]^{1/4}  \cdot \mathbb{E}[(X_2 - Y_2)^4]^{1/4} \cdot \mathbb{E}[Y_1^2]^{1/2} \\
    & \le (M \varepsilon  + 2M \delta^{1/4}) \cdot (\mathbb{E}[X_2^2]^{1/2} + \mathbb{E}[Y_1^2]^{1/2}) \le 4M^2 (\varepsilon + \delta^{1/4}).
    \end{aligned} 
\end{flalign}

\noindent Here, in the third statement we used that $\mathbbm{1}\cdot |X_i - Y_i| \le M \varepsilon$, $\mathbb{E}[(X_i - Y_i)^4] \le 16 (\mathbb{E}[X_i^4] + \mathbb{E}[Y_i])^4 \le 16 M^4$, and $\mathbb{P}[\mathsf{E}^{\complement}] \le \delta$; in the fourth we used that $\mathbb{E}[Z^2]^{1/2} \le \mathbb{E}[Z^4]^{1/4} \le M$ for $Z \in \{ X_2, Y_1 \}$.
\end{proof}

 \begin{lem}\label{lem:res_cov_bd}

 There exists a constant $\mathfrak{c}>0$ such that the following holds. Adopt the notation of Proposition \ref{prop:eff_varbound} and, for any $z \in \mathbb{C} \setminus \eig \mathbf{L}_{\mathfrak{N}}$, set $\bm{\mathfrak{G}}(z) = [\mathfrak{G}_{ij} (z)] = (\mathbf{L}_{\mathfrak{N}} - z \cdot \Id)^{-1}$. Let $i, j \in \llbracket 1, \mathfrak{N} \rrbracket$ be integers such that either we have $|i-j| \ge (\log \mathfrak{N})^3$, or we have both $|i-j| > (\log N)^4$ and $\dist ( \{ i, j \}, \{ 1, \mathfrak{N}\}) > (\log \mathfrak{N})^6$. Then,
  \begin{flalign}
 \label{eqn:exp_decay1}
 \begin{aligned} 
\Bigg|\Cov \Bigg( & \int_{-(\log N)^{\gamma}}^{(\log N)^{\gamma} }  \phi(E) \cdot \Imaginary \mathfrak{G}_{i i}(E + \i \eta) d E + a X_{i;\mathfrak{N}}, \\ 
 & \quad \int_{-(\log N)^{\gamma}}^{(\log N)^{\gamma} }  \phi(E) \cdot\Imaginary \mathfrak{G}_{j j}(E + \i \eta) d E + a X_{j;\mathfrak{N}} \Bigg) \Bigg|
 \leq \mathfrak{c}^{-1} M^2 (e^{-\mathfrak{c}|i-j|} + e^{-\mathfrak{c} (\log \mathfrak{N})^2}).
 \end{aligned} 
 \end{flalign}

 \end{lem}
 
\begin{proof}

Throughout this proof, we abbreviate $\mathbf{L}_{\mathfrak{N}} = \bm{\mathfrak{L}} = [\mathfrak{L}_{mn}]$ (which is sampled according to thermal equilibrium) and assume that $i < j$.  Observe by \eqref{integralg2} and \eqref{eqn:Xidef} that, for any $i \in \llbracket 1, \mathfrak{N} \rrbracket$,  
\begin{equation}\label{eqn:res_intbd}
\Bigg| \frac{1}{\pi} \int_{-(\log N)^{\gamma}}^{(\log N)^{\gamma}}  \phi(E) \cdot \Imaginary \mathfrak{G}_{i i}(E+\i\eta)  d E \Bigg| \leq M, \qquad \text{and} \qquad \mathbb{E} [|X_{i;\mathfrak{N}}|^4] \le C.
\end{equation}

\noindent Now fix an index $k \in \llbracket (2i+j)/3, (i+2j)/3 \rrbracket$, and define the $\mathfrak{N} \times \mathfrak{N}$ matrix $\tilde{\bm{\mathfrak{L}}} = [\tilde{\mathfrak{L}}_{mn}]$ by setting $\tilde{\mathfrak{L}}_{mn} = \mathfrak{L}_{mn}$ whenever $(m,n) \notin \{ (k,k), (k+1,k), (k,k+1) \}$ and $\tilde{\mathfrak{L}}_{mn} = 0$ otherwise. For any $z \in \mathbb{C} \setminus \eig \tilde{\bm{\mathfrak{L}}}$, denote $\tilde{\bm{\mathfrak{G}}}(z) = [\tilde{\mathfrak{G}}_{ij} (z)] = (\tilde{\bm{\mathfrak{L}}} - z \cdot \Id)^{-1}$. Observe under this notation that 
\begin{flalign}
\label{giijj}
(\tilde{\mathfrak{G}}_{ii} (z), X_{i;\mathfrak{N}})\quad \text{and} \quad  (\tilde{\mathfrak{G}}_{jj} (z), X_{j;\mathfrak{N}}) \quad \text{are independent}.
\end{flalign} 

By Lemma \ref{lem:resolv_coupling} (recalling that $\eta = e^{-(\log N)^3}$ and that $\dist ( \{ i, j \}, k) \ge (j-i)/3 \ge (\log N)^4/3$), there exists a constant $c_1 \in (0,1)$ and an event $\mathsf{E}_1$ such that $\mathbb{P}[\mathsf{E}_1^{\complement}] \le c_1^{-1} e^{-c_1 (\log \mathfrak{N})^2}$ and 
\begin{equation}\label{eqn:GGtilddf}
\displaystyle\max_{m \in \{ i, j \}} \displaystyle\sup_{|E| \le (\log N)^{\gamma}} |\mathfrak{G}_{mm} (E + \i \eta) - \tilde{\mathfrak{G}}_{mm} (E + \i \eta)|  \cdot \mathbbm{1}_{\mathsf{E}_1} \leq  e^{(\log \mathfrak{N})^2-c_1 (j-i)}.
\end{equation}

\noindent If $|i-j| \ge 2c_1^{-1} (\log \mathfrak{N})^2$, then the right side of \eqref{eqn:GGtilddf} is at most $e^{-c_1 |i-j|/2}$. This, with \eqref{eqn:res_intbd} (and its analog with $\tilde{\mathfrak{G}}_{ii} (z)$ replacing $\mathfrak{G}_{ii} (z)$), \eqref{giijj}, and Lemma \ref{lem:cb} applied with the $(X_1, X_2, Y_1, Y_2)$ there equal to $(X_i, X_j, Y_i, Y_j)$ here, where for each $m \in \{ i, j \}$ the latter are given by 
\begin{align*}
    & X_m = \frac{1}{\pi} \int_{-(\log N)^{\gamma}}^{(\log N)^{\gamma} }  \phi(E) \cdot \Imaginary \mathfrak{G}_{mm} (E + \i \eta) d E + a X_{m;\mathfrak{N}}; \\
    &Y_m = \frac{1}{\pi} \int_{-(\log N)^{\gamma}}^{(\log N)^{\gamma} }  \phi(E) \cdot \Imaginary \tilde{\mathfrak{G}}_{mm}(E + \i \eta) d E + a X_{m;\mathfrak{N}},
\end{align*}

\noindent implies the lemma. So, we may assume in what follows that $j-i < 2c_1^{-1} (\log \mathfrak{N})^2$. Since $j-i \ge (\log N)^4$, this implies that $\log \mathfrak{N} \ge c_1 (\log N)^2/2$. 

In this case, we will first ``shrink'' the $\mathfrak{N} \times \mathfrak{N}$ matrix $\bm{\mathfrak{L}}$ to an $\mathfrak{N}' \times \mathfrak{N}'$ matrix $\bm{\mathfrak{L}}'$, where $\mathfrak{N}'$ satisfies $(\log \mathfrak{N}')^2 < c_1 (j-i)/2$. Then, we will use Lemma \ref{glambda} to compare the resolvents of $\bm{\mathfrak{L}}$ and $\bm{\mathfrak{L}}'$. Next, using Lemma \ref{lem:resolv_coupling} again, we will show that the resolvent of $\bm{\mathfrak{L}}'$ satisfies a bound analogous to \eqref{eqn:GGtilddf}, with the $\mathfrak{N}$ there replaced by $\mathfrak{N}'$ here. In particular the factor of $e^{(\log \mathfrak{N})^2}$ will be replaced by $e^{(\log \mathfrak{N}')^2} \le e^{ c_1 (j-i)/2}$, which when multiplied by $e^{-c_1 |i-j|}$ will be at most $e^{-c_1 |i-j|/2}$.

To implement this, set $\mathfrak{N}' = \lceil e^{c_1 (j-i)^{1/2}/2} \rceil \le \mathfrak{N}^{3/4}$, and let $\mathfrak{N}_1' \in \llbracket (\log \mathfrak{N})^5 + 1, i-(\log \mathfrak{N}')^5 \rrbracket$ and $\mathfrak{N}_2 \in \llbracket j+(\log \mathfrak{N}')^5, \mathfrak{N} - \lceil (\log \mathfrak{N})^5 \rrbracket$ be integers such that $\mathfrak{N}' = \mathfrak{N}_2' - \mathfrak{N}_1' + 1$ (such integers exist since $\dist ( \{ i, j \}, \{ \mathfrak{N}_1, \mathfrak{N}_2 \}) \ge (\log \mathfrak{N})^6$ when $|i-j| < (\log \mathfrak{N})^3$). Then define the $\mathfrak{N}' \times \mathfrak{N}'$ matrix $\bm{\mathfrak{L}}' = [\mathfrak{L}_{mn}']$ with rows and columns indexed by $m,n \in \llbracket \mathfrak{N}_1', \mathfrak{N}_2' \rrbracket$ by setting $\mathfrak{L}_{mn}' = \mathfrak{L}_{mn} \cdot \mathbbm{1}_{m,n \in \llbracket \mathfrak{N}_1', \mathfrak{N}_2' \rrbracket}$ for each such $(m,n)$. Denote the resolvent of this matrix by $\bm{\mathfrak{G}}' (z) = [\mathfrak{G}_{mn}' (z)] = (\bm{\mathfrak{L}}' - z \cdot \Id)^{-1}$. 

Applying Lemma \ref{glambda} (with the $(N,\mathfrak{N})$ there given by $(\mathfrak{N}', \mathfrak{N})$ here, using the fact that $\log \mathfrak{N}' \ge (j-i)^{1/3}$) yields a constant $c_2 \in (0, 1)$ and an event $\mathsf{E}_2$ with $\mathbb{P}[\mathsf{E}_2^{\complement}] \le c_2^{-1} e^{-c_2 (\log \mathfrak{N}')^2}$ such that
\begin{flalign}
\label{2gij2}
\displaystyle\max_{m \in \{ i, j \} } \displaystyle\sup_{|E| \le (\log N)^{\gamma}} | \mathfrak{G}_{mm} (E + \mathrm{i} \eta) - \mathfrak{G}_{mm}' (E + \mathrm{i} \eta)| \cdot \mathbbm{1}_{\mathsf{E}_2}\le c_2^{-1} e^{-c(\log \mathfrak{N}')^4} \le e^{-|i-j|}.
\end{flalign} 

 Next, similarly to in the definition of $\tilde{\bm{\mathfrak{L}}}$, define the $\mathfrak{N}' \times \mathfrak{N}'$ matrix $\tilde{\bm{\mathfrak{L}}}' = [\tilde{\mathfrak{L}}_{mn}']$ by setting $\tilde{\mathfrak{L}}_{mn}' = \mathfrak{L}_{mn}'$ whenever $(m,n) \notin \{ (k,k), (k+1,k), (k,k+1) \}$ and $\tilde{\mathfrak{L}}_{mn}' = 0$ otherwise. For any $z \in \mathbb{C} \setminus \eig \bm{\tilde{\mathfrak{L}}}'$, set $\tilde{\mathfrak{G}}' (z) = [\tilde{\mathfrak{G}}_{mn}' (z)] = (\bm{\tilde{\mathfrak{L}}}' - z \cdot \Id)^{-1}$. Analogously to \eqref{eqn:GGtilddf}, we have that 
\begin{flalign}
\label{giijj2}
(\tilde{\mathfrak{G}}_{ii}' (z), X_{i;\mathfrak{N}})\quad \text{and} \quad  (\tilde{\mathfrak{G}}_{jj}' (z), X_{j;\mathfrak{N}}) \quad \text{are independent}.
\end{flalign}

\noindent By Lemma \ref{lem:resolv_coupling}, there exists an event $\mathsf{E}_3$ such that $\mathbb{P}[\mathsf{E}_3^{\complement}] \le c_1^{-1} e^{-c_1 (\log \mathfrak{N}')^2} \le c_1^{-1} e^{-c_1^3 (j-i)/8}$ (recalling that $\log \mathfrak{N}' \ge c_1 (j-i)^{1/2} / 4$) and 
\begin{equation}\label{2gkk2}
\displaystyle\max_{m \in \{ i, j \}} \displaystyle\sup_{|E| \le (\log N)^{\gamma}} |\mathfrak{G}_{mm}' (E + \i \eta) - \tilde{\mathfrak{G}}_{mm}' (E + \i \eta)|  \cdot \mathbbm{1}_{\mathsf{E}_3} \leq  e^{(\log \mathfrak{N}')^2-c_1 (j-i)} \le e^{-c_1(j-i)/2},
\end{equation}

\noindent where in the last inequality we used the fact that $\log \mathfrak{N}' \le c_1 (j-i)^{1/2}/2$. 

Combining \eqref{eqn:res_intbd}, \eqref{2gij2}, \eqref{giijj2}, and \eqref{2gkk2}, with Lemma \ref{lem:cb} applied with the $(X_1, X_2, Y_1, Y_2)$ there equal to $(X_i, X_j, Y_i, Y_j)$ here, where for each $m \in \{ i, j \}$ the latter are given by 
\begin{align*}
    & X_m = \frac{1}{\pi} \int_{-(\log N)^{\gamma}}^{(\log N)^{\gamma} }  \phi(E) \cdot \Imaginary \mathfrak{G}_{mm} (E + \i \eta) d E + a X_{m;\mathfrak{N}}; \\
    &Y_m = \frac{1}{\pi} \int_{-(\log N)^{\gamma}}^{(\log N)^{\gamma} }  \phi(E) \cdot \Imaginary \tilde{\mathfrak{G}}_{mm}' (E + \i \eta) d E + a X_{m;\mathfrak{N}},
\end{align*}

\noindent then establishes the lemma.
\end{proof}

\begin{proof}[Proof of Proposition \ref{prop:eff_varbound}]

Throughout this proof, we abbreviate $\phi_{\eta; N} = \phi_{\eta}$ and denote the left side of \eqref{eqn:stable_witheta} by $\Delta$. To ease notation, we will also assume throughout that $a = 0$, as the proof for arbitrary $a \in [-M,M]$ is entirely analogous. As in Lemma \ref{lem:res_cov_bd}, for any $z \in \mathbb{C}$, let $\mathbf{G} (z) = [G_{ij} (z)] = (\mathbf{L} - z \cdot \Id)^{-1}$ and $\bm{\mathfrak{G}} = [\mathfrak{G}_{ij} (z)] = (\mathbf{L}_{\mathfrak{N}} - z \cdot \Id)^{-1}$ denote the resolvents of $\mathbf{L}$ and $\mathbf{L}_{\mathfrak{N}}$, respectively. 

Then observe  from \eqref{eqn:phietn} and \eqref{gijuij} (with the orthonormality of the $(\mathbf{u}_1, \mathbf{u}_2, \ldots , \mathbf{u}_N)$ there) that 
\begin{multline}\label{eqn:main_comp}
	N \Delta =  \sum_{i, j=1}^{N} \Cov \left( \frac{1}{\pi} \int_{-(\log N)^{\gamma}}^{(\log N)^{\gamma} } \phi(E) \Imaginary G_{ii}(E + \i \eta)  d E, \frac{1}{\pi} \int_{-(\log N)^{\gamma}}^{(\log N)^{\gamma} } \phi(E) \Imaginary G_{j j}(E + \i \eta)  d E \right) \\
	- \frac{N}{\mathfrak{N}} \sum_{i, j=1}^{\mathfrak{N}}\Cov\left( \frac{1}{\pi} \int_{-(\log N)^{\gamma}}^{(\log N)^{\gamma} } \phi(E) \Imaginary \mathfrak{G}_{ii}(E + \i \eta) dE, \frac{1}{\pi} \int_{-(\log N)^{\gamma}}^{(\log N)^{\gamma} } \phi(E) \Imaginary \mathfrak{G}_{j j}(E + \i \eta) d E   \right)    .
\end{multline}

\noindent Next, for any $i,j \in \llbracket 1, \mathfrak{N} \rrbracket$, if there is an integer $k \in \llbracket 0,  \mathfrak{N}/N- 1 \rrbracket$ such that $kN+1 \leq  i, j \le (k+1)N$, then define $i' =N_1+ i- k N$ and $j' = N_1+ j - k N$. Then let $\mathbf{G}^{(k)} (z) = [G_{ij}^{(k)} (z)]$ be the $N \times N$ matrix, with entries indexed by $i, j \in \llbracket kN+1, (k+1)N \rrbracket$, defined by setting $G_{ij}^{(k)} (z) = G_{i'j'} (z)$ for each such $(i,j)$; in particular, $\mathbf{G}^{(k)}$ is obtained from $\mathbf{G}$ by shifting its indices by $kN$. Then, by writing the first sum in \eqref{eqn:main_comp} as $\mathfrak{N}/N$ copies of it, multiplied by $N/\mathfrak{N}$, we obtain 
\begin{flalign}\label{eqn:dif}
	N \Delta = \displaystyle\sum_{1\le i, j \le \mathfrak{N}} \Delta_{ij},
\end{flalign} 

\noindent where for each $i,j \in \llbracket 1, \mathfrak{N} \rrbracket$ we have denoted (observing that the sum over $k$ below is supported on at most one term)
\begin{flalign}
	\label{deltaij} 
	\begin{aligned}  
	\Delta_{ij} & =  \frac{N}{\mathfrak{N}} \Bigg|   \sum_{k = 0}^{\mathfrak{N}/N-1} \mathbbm{1}_{k N + 1 \leq  i, j \le (k+1) N} \cdot  \Cov \bigg( \frac{1}{\pi} \int_{-(\log N)^{\gamma}}^{(\log N)^{\gamma} } \phi(E) \Imaginary G_{ii}^{(k)}(E + \i \eta)  d E, \\
	& \qquad \qquad \qquad \qquad \qquad \qquad \qquad \qquad \qquad \quad  \frac{1}{\pi} \int_{-(\log N)^{\gamma}}^{(\log N)^{\gamma} } \phi(E) \Imaginary G_{j j}^{(k)}(E + \i \eta)  d E \bigg) \\
	& \qquad \qquad \qquad -\Cov \bigg( \frac{1}{\pi} \int_{-(\log N)^{\gamma}}^{(\log N)^{\gamma} } \Imaginary \mathfrak{G}_{ii}(E + \i \eta) dE, \frac{1}{\pi} \int_{-(\log N)^{\gamma}}^{(\log N)^{\gamma} } \Imaginary \mathfrak{G}_{j j}(E + \i \eta) d E   \bigg) \Bigg|.
	\end{aligned} 
\end{flalign}

Next, we will bound the $\Delta_{ij}$ depending on four cases for $(i,j) \in \llbracket 1, \mathfrak{N} \rrbracket^2$. 

\begin{enumerate}[A.]
	\item Let $\mathcal{S}_{\mathrm{A}}$ denote the set of $(i,j)$ such that $\dist (\{ i, j \}, \{ 1, \mathfrak{N} \}) \le (\log \mathfrak{N})^6$ and $|i-j| \le (\log \mathfrak{N})^4$. 
	\label{item:A}
	\item Let $\mathcal{S}_{\mathrm{B}}$ denote the set of $(i,j) \notin \mathcal{S}_{\mathrm{A}}$ with $|i-j| > (\log N)^4$.  
	\label{item:B}
	\item Let $\mathcal{S}_{\mathrm{C}}$ denote the set of $(i,j) \notin \mathcal{S}_{\mathrm{A}} \cup \mathcal{S}_{\mathrm{B}}$ with $\dist ( \{ i, j \},  kN) \le (\log N)^4$, for some $k \in \mathbb{Z}$.
	\label{item:C}
	\item Let $\mathcal{S}_{\mathrm{D}}$ denote the set of $(i,j) \notin \mathcal{S}_{\mathrm{A}} \cup \mathcal{S}_{\mathrm{B}} \cup \mathcal{S}_{\mathrm{C}}$.
	\label{item:D}
\end{enumerate}

\noindent For any index $\mathrm{X} \in \{ \mathrm{A}, \mathrm{B}, \mathrm{C}, \mathrm{D}\}$, define $\Delta_{\mathrm{X}} = \sum_{(i,j) \in \mathcal{S}_{\mathrm{X}}} \Delta_{ij}$. Then, \eqref{eqn:dif} implies
\begin{flalign}
	\label{ndelta} 
	N \Delta = \Delta_{\mathrm{A}} + \Delta_{\mathrm{B}} + \Delta_{\mathrm{C}} + \Delta_{\mathrm{D}},
\end{flalign} 

\noindent so we will estimate $\Delta_{\mathrm{X}}$ in each of the cases $\mathrm{X} \in \{ \mathrm{A}, \mathrm{B}, \mathrm{C}, \mathrm{D}\}$. Observe from \eqref{integralg2} and \eqref{eqn:Xidef}that 
\begin{flalign}
	\label{deltaijestimate} 
	\displaystyle\max_{i,j \in \llbracket 1, \mathfrak{N} \rrbracket} \Delta_{ij} \le C N \mathfrak{N}^{-1} M^2.
\end{flalign}

\textbf{Case~\ref{item:A}:} By \eqref{deltaijestimate} and the fact that $|\mathcal{S}_{\mathrm{A}}| \le C(\log \mathfrak{N})^{10}$, we have that
\begin{flalign}
	\label{deltaa} 
	\Delta_{\mathrm{A}} \le C M^2 N \mathfrak{N}^{-1} (\log \mathfrak{N})^{10}.
\end{flalign} 

\textbf{Case~\ref{item:B}:} For any $(i,j) \in \mathcal{S}_{\mathrm{B}}$, Lemma \ref{lem:res_cov_bd} applies and yields 
\begin{flalign*}
	\Delta_{ij} \le C M^2 N \mathfrak{N}^{-1} (e^{-c|i-j|} + e^{-c(\log \mathfrak{N})^2}).
\end{flalign*}

\noindent Summing over $(i,j)$ we deduce, since $|i-j| > (\log N)^4$ and $\mathfrak{N} \ge N$, that
\begin{flalign}
	\label{deltab} 
	\mathcal{S}_{\mathrm{B}} \le CM^2 N \mathfrak{N}^{-1} \displaystyle\sum_{\substack{i,j \in \llbracket 1, \mathfrak{N} \rrbracket \\ |i-j| > (\log N)^4}} (e^{-c|i-j|} + e^{-c(\log \mathfrak{N})^2}) \le CM^2 N \mathfrak{N}^{-1} \cdot \mathfrak{N} e^{-c(\log N)^2/2} \le CM^2.
\end{flalign}

\textbf{Case \ref{item:C}:} Observe that $|\mathcal{S}_{\mathrm{C}}| \le C N^{-1} \mathfrak{N} (\log N)^8$. Indeed, fixing $(i,j) \in \mathcal{S}_{\mathrm{C}}$, we have for some $h \in \{ i, j \}$ and $k \in \mathbb{Z}$ that $|h-kN| \le (\log N)^4$. There are at most $N^{-1} \mathfrak{N}$ choices for $k$; at most $C(\log N)^4$ choices for $h$ given $k$; and at most $C(\log N)^4$ choices for the element $\{ i, j \} \setminus \{ h \}$ given $h$ (since $|i-j| \le (\log N)^4$, as $(i,j) \notin \mathcal{S}_{\mathrm{R}} \cup \mathcal{S}_{\mathrm{B}}$). Together with \eqref{deltaijestimate}, this yields 
\begin{flalign}
	\label{deltac}
	\Delta_{\mathrm{C}} \le C M^2 (\log N)^8.
\end{flalign}

\textbf{Case \ref{item:D}:} Let $(i,j) \in \mathcal{S}_{\mathrm{D}}$. Then, we have $|i-j| \le (\log N)^4$; we have $\dist (\{i, j\}, \{ 1, \mathfrak{N} \}) > (\log \mathfrak{N})^6$; and we have $\dist (\{ i, j\}, kN) > (\log N)^4$ for all $k \in \mathbb{Z}$. In particular, there exists some integer $k \in \llbracket 0, \mathfrak{N}/N-1 \rrbracket$ such that $i,j \in \llbracket kN+1, (k+1)N \rrbracket$, and $\dist (\{i,j\}, \{ kN+1, (k+1)N\}) \ge (\log N)^4$. Then, by Lemma \ref{glambda}, there exists an overwhelmingly probable event $\mathsf{E}$ for which 
\begin{flalign}\label{eqn:resdfC}
	\begin{aligned}
		\mathbbm{1}_{\mathsf{E}} \cdot \big( |\mathfrak{G}_{ii}(E + \i \eta)-G_{ii}^{(k)}(E + \i \eta)|  |\mathfrak{G}_{j j}(E + \i \eta)- G_{jj}^{(k)}(E + \i \eta)| \big) &\leq c^{-1}e^{-c (\log N)^4} .
	\end{aligned}
\end{flalign}
Therefore, we may apply Lemma \ref{lem:cb}, with the $\delta$ there equal to $c^{-1} e^{-c(\log N)^4}$ here; the $M$ there equal to $M$ here (by \eqref{eqn:res_intbd}) the $(X_1, X_2, Y_1, Y_2; \mathsf{E})$ there equal to $(X_i, X_j, Y_i, Y_j; \mathsf{E})$ here, where for each $m \in \{ i, j \}$ the latter are given by
\begin{flalign}\label{eqn:x1x2y1y2C}
	\begin{aligned}
		X_m &=  \frac{1}{\pi} \int_{-(\log N)^{\gamma}}^{(\log N)^{\gamma} }   \phi(E) \Imaginary \mathfrak{G}_{mm}(E + \i \eta)dE;  \quad Y_m =  \frac{1}{\pi} \int_{-(\log N)^{\gamma}}^{(\log N)^{\gamma} }   \phi(E) \Imaginary G_{mm}^{(k)}(E + \i \eta)dE.
	\end{aligned}
\end{flalign}

\noindent Thus, under the above notation, we have from \eqref{deltaij} that
\begin{equation}
	\Delta_{ij} = N \mathfrak{M}^{-1}  |\Cov(X_i, X_j) - \Cov(Y_i, Y_j) | \leq C M^2 N \mathfrak{N}^{-1} e^{-c (\log N)^2}.
\end{equation}

\noindent Since $|\mathcal{S}_{\mathrm{D}}| \le C \mathfrak{N} (\log N)^4$ (as $i-j| \le (\log N)^4$), it follows upon summing the above estimate over $(i,j) \in \mathcal{S}_{\mathrm{D}}$ that  
\begin{flalign}
	\label{deltad} 
	\Delta_{\mathrm{D}} \le CM^2.
\end{flalign}

\noindent Summing \eqref{deltaa}, \eqref{deltab}, \eqref{deltac}, and \eqref{deltad} yields $\Delta \le M^2 (\log N)^{11}$. Together with \eqref{ndelta}, this establishes the proposition.
\end{proof}

\subsection{Changing between different scales}
\label{subsec:multi_scale}

Throughout this section, for any integer $\mathfrak{N} \ge N$, we the notation $\mathbf{L}_{\mathfrak{N}}$, $\lambda_{i;\mathfrak{N}}$, and $X_{i;\mathfrak{N}}$ defined above Proposition \ref{prop:eff_varbound} (and recall that $ \lambda_i = \lambda_{i;N}$ and $X_{i;N} = X_i$ for each $i$). The next proposition bounds how the variance of the linear statistic $\sum_{i=1}^{N}\phi_{\eta;N}(\lambda_i) $ changes upon increasing $N$ and decreasing $\eta$.

\begin{prop}\label{prop:dyadic}

There exists constants $\mathfrak{c}>0$ and $\mathfrak{C}>1$ such that the following holds. Let $\eta = e^{-(\log N)^3}$ and $\eta' = e^{-(\log \mathfrak{N})^3}$. Let $M \geq 1$ be a real number, $\phi : \mathbb{R} \rightarrow \mathbb{R}$ be a function satisfying Assumption \ref{ass:phiN}, and $a \in \mathbb{R}$ be a real number satisfying $|a| \leq M$. Then, for any $\mathfrak{N} \geq N \ge \mathfrak{C}$ of the form $\mathfrak{N} = N^{2^m}$ for some integer $m \geq 0$, we have 
\begin{equation}\label{eqn:bd1}
\Bigg|\frac{1}{N} \var \bigg(   \sum_{i=1}^{N}\phi_{\eta;N}(\lambda_i) + a  X_i\bigg)  - \frac{1}{\mathfrak{N}} \var \bigg(   \sum_{i=1}^{\mathfrak{N}}\phi_{\eta'; \mathfrak{N}}(\lambda_{i;\mathfrak{N}}) + a X_{i;\mathfrak{N}}  \bigg) \Bigg| \leq  M^2 N^{-\mathfrak{c}}. 
\end{equation}

 \noindent Additionally, 
\begin{equation}\label{eqn:bd2}
\Bigg|\frac{1}{\mathfrak{N}} \var \bigg(   \sum_{i=1}^{\mathfrak{N}}\phi(\lambda_{i;\mathfrak{N}}) + a X_{i;\mathfrak{N}}\bigg)  - \frac{1}{\mathfrak{N}} \var \bigg(   \sum_{i=1}^{\mathfrak{N}}\phi_{\eta';\mathfrak{N}}(\lambda_{i;\mathfrak{N}}) + a X_{i;\mathfrak{N}} \bigg) \Bigg| \leq  \displaystyle\frac{M^2}{e^{\mathfrak{c} (\log \mathfrak{N})^{6/5}}}, 
\end{equation}

\end{prop}
\begin{proof}
Throughout the proof, for any integer $K \ge N$, we let (recalling below the notation for the $K \times K$ Lax matrix $\mathbf{L}_K$ sampled under thermal equilibrium)
\begin{flalign}
    \label{lkf} 
\mathcal{L}_K (f) = K^{-1/2} \sum_{i=1}^K \big( f(\lambda_{i;K}) + a X_{i;K} \big); \qquad \bar{\mathcal{L}}_K (f) \coloneqq \mathcal{L}_K (f)  - \mathbb{E}[\mathcal{L}_K (f)].
\end{flalign} 

\noindent Observe for any $K_1,K_2 \ge N$ and $\eta'' \le \eta$ that  
\begin{equation}\label{eqn:LKsizebd}
|\mathcal{L}_{K_1} (\phi_{\eta'';K_2})| \le K_1^{1/2} \displaystyle\sup_{|E| \le (\log K_2)^{\gamma}} |\phi(E)| \le M K_1^{1/2} e^{10 (\log K_2)^{3\gamma/2}},
\end{equation}

\noindent where the first inequality holds by \eqref{lkf}, \eqref{integralg2}, and \eqref{eqn:phietn}, and the second by Assumption \ref{ass:phiN}.

Let us set $N_j \coloneqq N^{2^j}$ for each $j \in \llbracket 0, m \rrbracket$, where $N_0 = N$ and $m$ is as in the statement, so $N_m = \mathfrak{N}$; in addition, let $\eta_j = \eta^{8^j}$, so that $\eta_j = e^{-(\log N_j)^3}$ for each $j \in \llbracket 0, m \rrbracket$. We will bound the error in $\Var \mathcal{L}_{N_j} (\phi_{\eta_j;N_j})$ incurred by increasing $N_j$ to $N_{j+1}$; so, let $j \in \llbracket 0, m-1 \rrbracket$. Under the above notation, we have by \eqref{eqn:stable_witheta} (with $(N,\mathfrak{N}; M)$ there equal to $(N_j,N_{j+1}; M e^{10 (\log N_j)^{3\gamma/2}}) $ here) that 
\begin{equation}\label{eqn:step1}
|  \Var \bar{\mathcal{L}}_{N_j} (\phi_{\eta_j; N_j}) -  \Var \bar{\mathcal{L}}_{N_{j+1}} (\phi_{\eta_{j}; N_{j}})| \le C M^2 e^{10 (\log N_j)^{3\gamma/2}} N_j^{-1} (\log N_j)^{11}.
\end{equation}

To show \eqref{eqn:bd1} we will we will estimate $|\Var \bar{\mathcal{L}}_{N_{j+1}} (\phi_{\eta_{j}; N_{j}}) - \Var \bar{\mathcal{L}}_{N_{j+1}} (\phi_{\eta_{j+1}; N_{j+1}})|$. Together with \eqref{eqn:step1}, this would bound $|\Var \bar{\mathcal{L}}_{N_{j}} (\phi_{\eta_{j}; N_{j}}) - \Var \bar{\mathcal{L}}_{N_{j+1}} (\phi_{\eta_{j+1}; N_{j+1}})|$ for any $j$, enabling us (by summing over $j \in \llbracket 0, m-1 \rrbracket$) to then bound $|\bar{\mathcal{L}}_N (\phi_{\eta; N}) - \Var \bar{\mathcal{L}}_{\mathfrak{N}} (\phi_{\eta'; \mathfrak{N}})|$. To that end, observe that
\begin{flalign}\label{eqn:step2}
\begin{aligned} 
|& \Var \bar{\mathcal{L}}_{N_j^2} (\phi_{\eta_j; N_j}) - \Var \bar{\mathcal{L}}_{N_j^2} (\phi_{\eta_j^8; N_j^2})| \\
& \leq   \big( \var \big( \mathcal{L}_{N_j^2}(\phi_{\eta_j^8;N_j^2}) + \mathcal{L}_{N_j^2}(\phi_{\eta_j;N_j}) \big) \big)^{1/2}  \cdot \big(\var \big( \mathcal{L}_{N_j^2}(\phi_{\eta_j^8;N_j^2}) -  \mathcal{L}_{N_j^2}(\phi_{\eta_j;N_j}) \big) \big)^{1/2}  \\
&  \leq 4 M e^{10 (2 \log N_j)^{3 \gamma /2}} N_j \cdot     \big(\Var \mathcal{L}_{N_j^2}(\phi_{\eta_j^8;N_j^2} - \phi_{\eta_j;N_j}) \big)^{1/2},
  \end{aligned}
\end{flalign}

\noindent  by \eqref{eqn:LKsizebd}. To estimate the right side of \eqref{eqn:step2}, observe for any index $n \in \{ j, j + 1 \}$ and real number $x \in [-(\log N_n)^{\gamma}/2, (\log N_n)^{\gamma}/2]$ that  
\begin{flalign}
    \label{netan0}
    \begin{aligned} 
    |\phi_{\eta_n; N_n}(x) - \phi(x)| & \le \frac{1}{\pi} \int_{-(\log N_n)^{\gamma}}^{(\log N_n)^{\gamma}}  \frac{\eta_n  |\phi(E)-\phi(x)|}{(x-E)^2 +\eta_n^2} d E + \frac{1}{\pi} \int_{|E| > (\log N_n)^{\gamma}}  \frac{\eta_n |\phi(x)| }{(x-E)^2 +\eta_n^2} d E  \\
& \leq  \displaystyle\sup_{\substack{|E|, |x| \le (\log N_n)^{\gamma} \\ |E-x| \le \eta_n^{1/3}}}   |\phi(E)-\phi(x)|  \\
& \qquad + 3M e^{10 ( \log N_n)^{3 \gamma /2}} \int_{|E-x| > \eta_n^{1/3}} \frac{\eta_n}{(x-E)^2 +\eta_n^2} d E \\
& \leq C \eta_n^{1/3} M e^{C (\log N_n)^2} + C M e^{10 ( \log N_n)^{3 \gamma /2}} \eta_n^{1/3}  \leq M e^{-c (\log N_n)^3},
\end{aligned} 
\end{flalign}

\noindent where in the first bound we used \eqref{eqn:phietn} with the fact that $\pi^{-1} \eta_n \int_{-\infty}^{\infty} ((x-E)^2 + \eta_n^2)^{-1} dE = 1$; in the second we decomposed the first integral into the regions where $|x-E| \le \eta_n^{1/3}$ and where $|x-E| > \eta_n^{1/3}$ and used Assumption \ref{ass:phiN} (also with the above integral identity); in the third we again used Assumption \ref{ass:phiN}; and in the fourth we used the fact that $\eta_n = e^{-(\log N_n)^3}$.

Thus, denoting the event $\mathsf{E} = \mathsf{BND}_{\mathbf{L}_{N_n}} ( (\log N_n)^{\gamma} / 2)$ (from Definition \ref{eventbounded}), \eqref{lkf} and \eqref{netan0} together yield $\mathbbm{1}_{\mathsf{E}} \cdot \mathcal{L}_{N_n} (\phi_{\eta_n;N_n} - \phi) \le M e^{-c(\log N_n)^3}$. Since, by Lemma \ref{lem:bd_lem}, the event $\mathsf{E}$ holds with probability at least $1 - c^{-1} e^{-c(\log N_n)^{2\gamma}}$ and $\gamma = 3/5$, it quickly follows (see, for example, the derivation of \eqref{expectationH}) that 
\begin{flalign} 
\label{nln} 
\big( \Var \mathcal{L}_{N_n} (\phi_{\eta_n; N_n} - \phi) \big)^{1/2}  & \leq  C M e^{-c (\log N_n)^{6/5}}.
\end{flalign}

\noindent Inserting \eqref{nln} at $n \in \{ j, j+1 \}$ into \eqref{eqn:step2} (and recalling that $(N_{j+1}, \eta_{j+1}) = (N_j^2, \eta_j^8)$ and $\gamma = 3/5$) yields  
\begin{flalign}
\label{nj1l2}
|\Var \bar{\mathcal{L}}_{N_{j+1}} (\phi_{\eta_j;N_j}) - \Var \mathcal{L}_{N_{j+1}} (\eta_{j+1}; N_{j+1})| \le C M e^{-c(\log N_j)^{6/5}}.
\end{flalign} 

Combining the bounds \eqref{eqn:step1} and \eqref{nj1l2} (realizing that the left side of the former is larger), and recalling that $3\gamma/2 = 9/10$, we obtain for any $j \in \llbracket 0, m-1 \rrbracket$ that
\begin{equation}\label{eqn:onestep}
\left|\Var \bar{\mathcal{L}}_{N_{j}} (\phi_{\eta_{j}; N_{j}}) - \Var \bar{\mathcal{L}}_{N_{j+1}} (\phi_{\eta_{j+1}; N_{j+1}})\right|
\leq C M^2 e^{20 (\log N_j)^{9/10}} N_j^{-1} (\log N_j)^{12} \le CM^2 N_j^{-c},
\end{equation}

\noindent where $c>0$ may be any constant less than $1$. Summing \eqref{eqn:onestep}  over $j \in \llbracket 0, m-1 \rrbracket$ (and using the fact that $(\mathfrak{N}, \eta') = (N_m, \eta_m)$), we obtain
\begin{flalign*}
|\Var \bar{\mathcal{L}}_{\mathfrak{N}} (\phi_{\eta';\mathfrak{N}}) - \Var \bar{\mathcal{L}}_N (\phi_{\eta;N})| \le CM^2 \sum_{j =0}^{m-1} N_j^{-c} \le CM^2 N^{-c}, 
\end{flalign*}

 \noindent which yields \eqref{eqn:bd1} (at any constant $\mathfrak{c}<1$). 

 The last statement \eqref{eqn:bd2} of the lemma follows from the estimates 
\begin{flalign}\label{eqn:bd2intermed}
\begin{aligned} 
|\var & \left(   \mathcal{L}_{\mathfrak{N}}(\phi_{\eta';\mathfrak{N}})\right) - \var \left(  \mathcal{L}_{\mathfrak{N}}(\phi) \right)| \\
&  \leq \big( \var( \mathcal{L}_{\mathfrak{N}}(\phi_{\eta';\mathfrak{N}}) )^{1/2}+ \var( \mathcal{L}_{\mathfrak{N}}(\phi))^{1/2}\big) \cdot \var \left( \bar{\mathcal{L}}_{\mathfrak{N}}(\phi_{\eta';\mathfrak{N}})-  \bar{\mathcal{L}}_{\mathfrak{N}}(\phi) \right)^{1/2} \\
& \le C\mathfrak{N}^{1/2} M e^{10  (\log \mathfrak{N})^{3 \gamma /2   }}  \cdot \var \left( \bar{\mathcal{L}}_{\mathfrak{N}}(\phi_{\eta';\mathfrak{N}})-  \bar{\mathcal{L}}_{\mathfrak{N}}(\phi) \right)^{1/2} \le CM^2 e^{-c(\log \mathfrak{N})^{6/5}}.
\end{aligned} 
\end{flalign}

\noindent Here, the first and second bounds follow from  \eqref{eqn:LKsizebd} and the fact by Assumption \ref{ass:phiN} that we have $|\mathcal{L}_{\mathfrak{N}} (\phi)| \le \mathfrak{N}^{1/2} M e^{10 (\log \mathfrak{N})^{3\gamma/2}}$  on the event $\mathsf{BND}_{\mathbf{L}_{\mathfrak{N}}} ((\log \mathfrak{N})^{\gamma}/2)$ (recall Definition \ref{eventbounded}) that holds with probability at least $1 - c^{-1} e^{-c (\log \mathfrak{N})^{6/5}}$ by Lemma \ref{lem:bd_lem}. Moreover, the third bound follows from \eqref{nln} (at $j=m$). This establishes the lemma. 
\end{proof}

As a corollary of the previous proposition, we obtain a quantitative error bound on the difference between $N^{-1} \var ( \sum_{i=1}^{N}\phi(\lambda_i)) $ and its limit as $N$ tends to $\infty$ (keeping $\phi$ fixed) along a particular subsequence. Let us mention that, with further effort, it can be shown that the latter limit can be taken without passing to a subsequence, but we will not require this here.

\begin{cor}\label{cor:stable_noeta}

There exist constants $\mathfrak{c}>0$ and $\mathfrak{C}>1$ such that the following holds. Let $M \ge 1$ and $a \in [-M,M]$ be real numbers and $\phi: \mathbb{R} \rightarrow \mathbb{R}$ be a function satisfying Assumption \ref{ass:phiN}. Then, for any integer $\mathfrak{N} \geq N \ge \mathfrak{C}$ of the form $\mathfrak{N} = N_m = N^{2^{m}}$ for some integer $m \ge 0$, we have 
\begin{equation}\label{eqn:stable_noeta}
\Bigg|\frac{1}{N} \var \bigg(   \sum_{i=1}^{N} \big( \phi(\lambda_i) + a X_i \big) \bigg)  - \frac{1}{\mathfrak{N}} \var \bigg(   \sum_{i=1}^{\mathfrak{N}} \big( \phi(\lambda_{i;\mathfrak{N}}) + a X_{i;\mathfrak{N}} \big) \bigg) \Bigg| \leq M^2 N^{-\mathfrak{c}}.
\end{equation}

\noindent Moreover, the limit $\lim_{k \rightarrow \infty} N_k^{-1}  \var (\sum_{i=1}^{N_k} (\phi(\lambda_{i;N_k}) + aX_{i;N_k}) ) $ exists, and 
\begin{equation}\label{eqn:stable_noetalim}
\Bigg|\frac{1}{N} \var \bigg(   \sum_{i=1}^{N} \big( \phi(\lambda_i) + a X_i \big) \bigg)  - \lim_{k \rightarrow \infty} \frac{1}{N_k} \var \bigg(   \sum_{i=1}^{N_k} \big( \phi(\lambda_{i;N_k}) + aX_{i;N_k} \big) \bigg) \Bigg| \leq M^2 N^{-\mathfrak{c}}.
\end{equation}
\end{cor}

\begin{proof}

The first statement \eqref{eqn:stable_noeta} of the corollary follows from applying \eqref{eqn:bd2} (with the $\mathfrak{N}$ equal to $N$ here); then applying \eqref{eqn:bd1}; and next applying \eqref{eqn:bd2} again (now with the $\mathfrak{N}$ there equal to $\mathfrak{N}$ here). This bound \eqref{eqn:stable_noeta} (applied with the $(N, N_k)$ there equal to any $(N_j, N_k)$ here) implies that $(N_k^{-1} \Var (\sum_{i=1}^{N_k} (\phi(\lambda_{i;N_k})+aX_{i;N_k}))_{k \ge 0}$ forms a Cauchy sequence, meaning that it has a limit as $k$ tends to $\infty$. This, together with \eqref{eqn:stable_noeta} yields \eqref{eqn:stable_noetalim}.
\end{proof}

\subsection{Proof of Proposition \ref{prop:limvar}} 
\label{subsec:limvar_proof}

In this section we prove Proposition \ref{prop:limvar}; we adopt the notation of that proposition throughout. 

We will proceed by first replacing the function $\phi$ with a proxy $\psi$, which is analytic in a strip of the complex plane around the real line. Due to its analyticity, $\psi$ will be closely approximable by a polynomial $p$. After verifying that $p$ satisfies Assumption \ref{ass:phiN}, we will apply Corollary \ref{cor:stable_noeta} to $p$, where the limiting variance there has the explicit form given by Lemma \ref{lem:lim_var_formula}. Since $p \approx \psi \approx \phi$, we will then deduce that the same formulas hold with $p$ replaced by $\phi$, confirming the proposition. 

To implement this, set $\eta = e^{-(\log N)^3}$, and let $\tilde{\phi}: \mathbb{R} \rightarrow \mathbb{R}$ be any smooth function such that
\begin{flalign}
    \label{eqn:tildphibd}
    \begin{aligned} 
& \tilde{\phi} (x) = \phi(x), \quad \text{for all $x \in [-(\log N)^{\gamma}, (\log N)^{\gamma}]$}; \\
& \displaystyle\sup_{x \in \mathbb{R}} |\tilde{\phi}(x)| \le \displaystyle\frac{11M}{10}; \qquad \displaystyle\sup_{x \in \mathbb{R}} |\tilde{\phi}'(x)| \le Me^{2000(\log N)^{3/2}},
\end{aligned} 
\end{flalign}

\noindent which exists by \eqref{eqn:limvar_phiasmps}. Next, denoting the strip $\mathfrak{S} = \{ z \in \mathbb{C} : |\Imaginary z| \le \eta/2 \}$, define the analytic function $\psi : \mathfrak{S} \rightarrow \mathbb{R}$ by setting
\begin{equation}\label{eqn:psietadf}
    \psi(x) =  (2\pi)^{-1/2} \eta^{-1} \int_{-\infty}^{\infty}  e^{-(x-E)^2/(2\eta^2)} \cdot \tilde{\phi}(E) dE,
\end{equation}

\noindent for any $x \in \mathfrak{S}$. Since \eqref{eqn:tildphibd} implies that $|e^{-(x+\i y)^2/(2 \eta^2)} \cdot \tilde{\phi}(E)| \leq e^{1/8-x^2/(2 \eta^2)} \cdot 11M/10 \le e^{-x^2/(2\eta^2)} \cdot 2M$ for any $x + \mathrm{i} y \in \mathfrak{S}$ and $E \in \mathbb{R}$, from \eqref{eqn:psietadf} that
\begin{align}
\label{psiz} 
    |\psi(z)| &\leq 2M \cdot (2\pi)^{-1/2} \eta^{-1} \int_{-\infty}^\infty e^{-(\Real z-E)^2/(2\eta^2)} dE  \leq 2 M,
\end{align}

\noindent for any $z \in \mathfrak{S}$. In particular,
\begin{equation}\label{eqn:Mphieta}
    \mathcal{M}(\psi) \coloneqq 2 \sup_{z \in \mathfrak{S}} |\psi(z)| \leq 4 M .
\end{equation}

The following lemma from \cite{Tre13} (rescaled by $A$, specialized to the function $\psi$ above, and taken with the $\rho$ there equal to $1 + \eta/(2A)$ here) yields a polynomial approximation for $\psi$ through its Chebyshev projections. In what follows, we let $T_n(x) = ((x+\sqrt{x^2-1})^n + (x-\sqrt{x^2-1})^n)/2$ denote the $n$-th Chebyshev polynomial of the first kind, for any integer $n \ge 0$. 

\begin{lem}[{\cite[Theorems 8.1 and 8.2]{Tre13}}] 

\label{etap}

Let $A \ge 1$ be a real number and $D \ge 1$ be an integer. Define the polynomial $p_{A,D}: \mathbb{R} \rightarrow \mathbb{R}$ by setting
\begin{equation}\label{eqn:pexpans}
p_{A,D} (x)=\sum_{k=0}^{D} a_k \cdot T_k(xA^{-1}), \qquad \text{where} \qquad 
a_k=\frac{2-\mathbbm{1}_{k=0}}{\pi}\int_{-1}^{1}
\frac{\psi(Ax)\,T_k(x)}{(1-x^2)^{1/2}}\,dx .
\end{equation}

\noindent Then, for any integer $k \ge 0$, we have 
\begin{equation}\label{eqn:akbd}
    |a_k|\le \mathcal{M}(\psi) \cdot \bigg( 1 + \displaystyle\frac{\eta}{2 A} \bigg)^{-k}; \qquad 
\sup_{|x|\le A} |p_{A,D} (x)-\psi(x)| \le \displaystyle\frac{4A}{\eta} \cdot \mathcal{M}(\psi) \cdot \bigg( 1 + \displaystyle\frac{\eta}{2A} \bigg)^{-D}.
\end{equation}

\end{lem} 

In what follows, under the notation of Lemma \ref{etap}, we set 
\begin{flalign}
\label{0p0}
    A = \eta^{-10}; \qquad D = \lfloor \eta^{-12} \rfloor; \qquad p = p_{A,D}. 
\end{flalign}

\noindent Then observe from \eqref{eqn:akbd}, \eqref{eqn:Mphieta}, and \eqref{0p0} (with the fact that $M = e^{\sqrt{\log N}}$) that 
\begin{equation}\label{eqn:peps}
\sup_{|x|\le A} |p(x)-\psi(x)| \le  e^{-1/(20\eta)}.
\end{equation}

\begin{lem} 

\label{p0} 

The function $p: \mathbb{R} \rightarrow \mathbb{R}$ from \eqref{0p0} (and Lemma \ref{etap}) satisfies Assumption \ref{ass:phiN}, with the $M$ there equal to $4e^{\sqrt{\log N}}$ here. 
\end{lem} 

\begin{proof} 

By \eqref{eqn:peps} and \eqref{psiz}, we have that $|p(x)| \le 4e^{\sqrt{\log N}}$ for all $|x| \le A$; since $(\log N)^{\gamma} < A$, this confirms the first statement in \eqref{eqn:phiasmps}, as well as the second for $|x| \le A$. For the second bound in \eqref{eqn:phiasmps} for $|x| > A =\eta^{-10}$, we estimate 
\begin{flalign}\label{eqn:ptail}
|p(x)| \leq (D+1) \cdot  4 e^{\sqrt{\log N}} \cdot (2 \eta^{10} |x|)^{D} \le  e^{ |x|^{3/2}} ,
\end{flalign}

\noindent where the first statement holds by the first equality in \eqref{eqn:pexpans}, the first bound in \eqref{eqn:akbd}, \eqref{eqn:Mphieta}, and the fact that $|T_k(xA^{-1})| \le (2A^{-1}|x|)^k$ for $x \ge 1$; the second follows from individually considering the cases when $A < |x| \le D$ and $|x| \ge D$. 

It remains to verify the last inequality in \eqref{eqn:phiasmps} for $p$. For $|x|>A$, we have
\begin{flalign*} 
|p'(x)| \leq |x|^{-1} \cdot D^2 \sup_{|y| \leq |x|} |p(y)| \le D^2 e^{|x|^{3/2}} \leq  e^{2|x|^{3/2}},
\end{flalign*}

\noindent where the first inequality follows from the Markov brothers inequality; the second from \eqref{eqn:ptail}; and the third from the bounds $|x|>A = \eta^{-10}$ and $D \le \eta^{-12}$. This addresses \eqref{eqn:phiasmps} when $|x|>A$. 

To address the last bound in \eqref{eqn:phiasmps} when $|x| \le A$, fix $x \in [-A, A]$ and observe by expanding $\psi$ over the basis of Chebyshev polynomials $(T_k)$ that $\psi(x) = \sum_{k=0}^{\infty} a_k \cdot T_k (xA^{-1})$. Hence, 
\begin{flalign*}
    |\psi'(x) - p'(x)| \le \displaystyle\sum_{k=D+1}^{\infty} |a_k| \cdot |T_k'(xA^{-1})| & \le M(\psi) \displaystyle\sum_{k=D+1}^{\infty}  |T_k'(xA^{-1})| \cdot (1+\eta^{11})^{-k}  \\
    & \le 4e^{\sqrt{\log N}} \displaystyle\sum_{k=D+1}^{\infty} k^2 (1+\eta^{11})^{-k} < 1,
\end{flalign*}

\noindent where in the first estimate we used the above expansion of $\psi$; in the second we used \eqref{eqn:Mphieta} and \eqref{eqn:akbd} (with the fact that $A = \eta^{-10}$); in the third we applied the Markov brothers inequality to $T_k$ (using the fact that $|T_k(yA^{-1})| \le 1$ for $|y| \le A$); and in the fourth we used the fact that $D = \lfloor \eta^{-12} \rfloor$. We then deduce \eqref{eqn:phiasmps} from this with the bound 
\begin{equation*}
    |\psi'(x)| \leq (2 \pi)^{-1/2} \eta^{-1} \displaystyle\int_{-\infty}^{\infty} e^{-(x-E)^2/(2\eta^2)} \cdot |\tilde{\phi}'(E)| dE \le 4Me^{2000(\log N)^{3/2}} \le M e^{1000 (\log N)^2},
\end{equation*}

\noindent which holds from \eqref{eqn:psietadf} and the last bound in \eqref{eqn:tildphibd}.
\end{proof}

\begin{proof}[Proof of Proposition \ref{prop:limvar}]

To ease notation, we only establish the lemma when $a = 0$, as the proof for arbitrary $a \in [-M, M]$ entirely analogous. Denoting 
\begin{flalign*} 
\mathcal{L}_N(f) = N^{-1/2} \sum_{i=1}^N f(\lambda_i),
\end{flalign*} 

\noindent for any  function $f: \mathbb{R} \rightarrow \mathbb{R}$ (as in \eqref{lkf}), we have
\begin{multline}\label{eqn:varNbd}
 \big|\var (   \mathcal{L}_{N}(\phi)) - \var (  \mathcal{L}_{N}(p)) \big|  \\
\leq \big ( \var( \mathcal{L}_{N}(\psi) )^{1/2}+ \var( \mathcal{L}_{N}(p))^{1/2}  \big) \cdot \var \big(   \mathcal{L}_{N}(\phi_{\eta;N})-   \mathcal{L}_{N}(p) \big)^{1/2} \\
+ \big( \var( \mathcal{L}_{N}(\phi) )^{1/2}+ \var( \mathcal{L}_{N}(\psi))^{1/2}  \big) \cdot \var \big(   \mathcal{L}_{N}(\phi)-   \mathcal{L}_{N}(\psi) \big)^{1/2} 
\end{multline}

Observe that
\begin{flalign}\label{eqn:vbds}
\begin{aligned}
    \var( \mathcal{L}_{N}(\psi) )  \leq 4 N M^2; \qquad \var( \mathcal{L}_{N}(p)) \leq 20 N M^2; \qquad    \var( \mathcal{L}_{N}(\phi)) \leq 20 N M^2.
\end{aligned}
\end{flalign}

\noindent  Indeed, the first bound in \eqref{eqn:vbds} follows from \eqref{psiz}. The second and third bounds there follow from the facts that $|\phi(\lambda)| \le M$ and $|p(\lambda)| \le 4M$ for $|\lambda| \le (\log N)^{\gamma}$ (by \eqref{eqn:limvar_phiasmps} and Lemma \ref{p0}); that, denoting the event $\mathsf{E} = \bigcap_{i=1}^N \{ |\lambda_i| \le (\log N)^{\gamma} \}$, we have $\mathbb{P}[\mathsf{E}] \ge 1 - c^{-1} e^{-c(\log N)^{6/5}}$ by Lemma \ref{lem:bd_lem} (and the fact that $\gamma = 3/5$); and the fact that the contribution to the variances in \eqref{eqn:vbds} off of this event is at most $c^{-1} e^{-c(\log N)^{6/5}}$ (see, for example, the derivation of \eqref{expectationH}). 
 
  Next, to bound \eqref{eqn:varNbd}, we claim that
 \begin{flalign}\label{eqn:Lphipdif}
 \begin{aligned}
 \var \left(   \mathcal{L}_{N}(\psi)-   \mathcal{L}_{N}(p) \right) \leq c^{-1} e^{-c(\log N)^{6/5}}; \qquad \var \left(   \mathcal{L}_{N}(\psi)-   \mathcal{L}_{N}(\phi) \right) \leq  c^{-1} e^{-c (\log N)^{6/5}}. 
 \end{aligned}
 \end{flalign}
 
 \noindent Indeed, the first bound follows from the fact by \eqref{eqn:peps} that $|\mathcal{L}_N (\psi) - \mathcal{L}_N (\phi)| \le Ne^{-1/(20\eta)} \le Ce^{-c(\log N)^3}$ on the above event $\mathsf{E} = \bigcap_{i=1}^N \{ |\lambda_i| \le (\log N)^{\gamma} \}$ that holds with probability at least $1 - c^{-1} e^{-c(\log N)^{6/5}}$. To obtain the second bound, first observe for $|x| \le (\log N)^{\gamma}$ that
 \begin{flalign}\label{eqn:psieta_phi}
 \begin{aligned} 
    |\psi(x) - \phi(x)| & \leq (2 \pi)^{-1/2} \eta^{-1} \displaystyle\int_{-\infty}^{\infty} e^{-(x-E)^2/(2\eta^2)} \cdot | \tilde{\phi}(E) - \phi(x)| dE  \\
    & \le  
    3M \eta^{1/2}  + (2 \pi)^{-1/2} \eta^{-1} \int_{|E-x|>\eta^{1/2}} e^{-(x-E)^2/(2 \eta^2)} \cdot (|\tilde{\phi}(E)| + |\phi(x)|) dE  \\
    &   \leq 3M\eta^{1/2} + M \eta^{-1} \int_{|E-x|>\eta^{1/2}}  e^{-(x-E)^2/(2 \eta^2)} dE  \leq C e^{-c (\log N)^3},
      \end{aligned} 
 \end{flalign}

 \noindent where the first estimate follows from \eqref{eqn:psietadf} with the identity $(2 \pi)^{-1/2} \eta^{-1} \int_{-\infty}^{\infty} e^{-(x-E)^2/(2\eta^2)} dE = 1$; the second follows from decomposing the integral into the regions where $|E-x| \le \eta^{1/2}$ and where $|E-x| > \eta^{1/2}$, and using the above identity, with the facts that $|\tilde{\phi}(E)| \le 11M/10$ by \eqref{eqn:tildphibd} and $|\phi(x)| \le M$ by \eqref{eqn:limvar_phiasmps} and the fact that $|x| \le (\log N)^{\gamma}$; the third follows from the above bounds on $|\tilde{\phi}|$ and $|\phi|$; and the fourth and fifth follow from performing the integral and using the definitions $\eta = e^{-(\log N)^3}$ and $M = e^{\sqrt{\log N}}$. Thus, $|\mathcal{L}_N (\psi) - \mathcal{L}_N (\phi)| \le C e^{-c(\log N)^3}$ on $\mathsf{E}$, which holds with probability at least $1 - c^{-1} e^{-c(\log N)^3}$; so, we deduce the second bound in \eqref{eqn:Lphipdif}.

 Next, by Corollary \ref{cor:stable_noeta} (and Lemma \ref{p0}), we obtain 
$$ \Big|\lim_{k \rightarrow \infty} \var( \mathcal{L}_{N^{2^k}}(p) ) - \var( \mathcal{L}_N(p) ) \Big| \leq C M^2 N^{-c}.$$

\noindent Since  Lemma \ref{lem:lim_var_formula} implies that $\lim_{\mathfrak{N} \rightarrow \infty} \var \left(  \mathcal{L}_{\mathfrak{N}}(p) \right) = \sigma^2(p)$, this with \eqref{eqn:Lphipdif} (and the fact that $M = e^{\sqrt{\log N}}$) yields 
\begin{equation}\label{eqn:varbound2}
|\var \left(   \mathcal{L}_{N}(\phi)\right) - \sigma^2(p) | \leq N^{-c}.
\end{equation}

\noindent Thus, to prove the proposition it suffices to show the two upper bounds 
\begin{align}
|\sigma^2(p) - \sigma^2(\psi)| \leq  N^{-\mathfrak{c}}; \qquad 
|\sigma^2(\psi) - \sigma^2(\phi)| \leq N^{-\mathfrak{c}} . \label{eqn:lim_equiv1}
\end{align}

We only show the second bound in \eqref{eqn:lim_equiv1}, as the proof of the first is very similar. Recalling Definition \ref{def:Cdef} (and its notation), we have 
\begin{equation}\label{eqn:var_explicit}
\sigma^2(\phi) = \big\langle (1- \theta \T \boldsymbol{\varrho}_{\beta})^{-1} (\phi - \mu_{\phi}\varsigma_0) , (1- \theta \T \boldsymbol{\varrho}_{\beta})^{-1} (\phi - \mu_{\phi} \varsigma_0) \rangle_{\varrho} = \langle \mathbf{F}  \phi, \mathbf{F} \phi \big\rangle_{\varrho}.
\end{equation}

\noindent where we have recalled the (bounded) operator $\mathbf{F} : \mathcal{H} \rightarrow \mathcal{H}$ from Definition \ref{operatorf}. In this way, we have that 
\begin{flalign}
\label{sigmapsi2}
|\sigma^2(\phi) - \sigma^2(\psi)|  
&\leq  \big| \langle \mathbf{F} \psi, \mathbf{F}(   \phi -  \psi) \rangle_{\varrho} \big|  + \big| \langle \mathbf{F}  \phi , \mathbf{F}(   \phi -  \psi) \rangle_{\varrho} \big| \le CM \cdot \| \phi - \psi \|_{\mathcal{H}},
\end{flalign} 

\noindent where in the last bound we used the facts that $\| \phi \|_{\mathcal{H}} + \| \psi \|_{\mathcal{H}} \le CM$ (by the growth conditions on $\phi$ and $\psi$ implied by \eqref{eqn:limvar_phiasmps} and Lemma \ref{p0}). Moreover, we have that 
\begin{flalign*}
    \| \phi - \psi \|_{\mathcal{H}}^2 & \le  \displaystyle\int_{|x| \le (\log N)^{\gamma}} |\phi(x) - \psi(x)| \cdot \varrho (x) dx + \displaystyle\int_{|x| > (\log N)^{\gamma}} \big( |\phi(x)| + | \psi(x)| \big) \cdot \varrho (x) \\
    & \le C e^{-c(\log N)^{6/5}}, 
\end{flalign*}

\noindent where in the first bound we used \eqref{def:inner_prod}; in the second we used \eqref{eqn:psieta_phi}, which holds for $|x| \le (\log N)^{\gamma}$, and \eqref{eqn:psieta_phi} with the decay on $\varrho$ provided by Lemma \ref{lem:varrho_bd} (and \eqref{eqn:limvar_phiasmps} and Lemma \ref{p0}). Inserting this into \eqref{sigmapsi2} gives second bound in \eqref{eqn:lim_equiv1}; as mentioned above, the proof of the first is entirely analogous. This establishes the proposition. 
\end{proof}

\section{Continuity of prelimit processes}
\label{sec:continuity}

Throughout this section, we adopt Assumptions \ref{ass:NT_assumption} and \ref{ass:R_M_ass}, and we recall the notation in Definitions \ref{def:chi_def}, \ref{xilambdakt}, \ref{xi2}, \ref{xim}, and \ref{xi}. Moreover, we assume throughout (even when not explicitly stated) that $\theta<\theta_0(\beta)$ is sufficiently small, so that all of the results from Sections \ref{ProofQ}, \ref{sec:bg}, and \ref{sec:ind_linear} apply.

\subsection{Continuity bounds}
\label{subsec:lscont}

In this section we prove the below continuity bound for the functional $\Xi$ from  Definition \ref{xilambdakt}.

\begin{prop}[Continuity of $\Xi$]\label{cor:xicont}

Let $S \in [\mathfrak{M},T]$ be a real number. There exist constants $\mathfrak{c}>0$ and $\mathfrak{C}>1$ such that the following holds with overwhelming probability. For any integers $k_1, k_2 \in \llbracket - T^{3/2},  T^{3/2} \rrbracket$ with $|k_1-k_2| \leq S$; real numbers $t_1, t_2 \in [0,  T \log N]$ with $|t_1 - t_2| \leq S$; and real numbers $\Lambda_1, \Lambda_2 \in [- \log N, \log N]$ with $|\Lambda_1 - \Lambda_2| \leq ST^{-1}$, we have 
\begin{equation*}
|\Xi(\Lambda_1, k_1, t_1)- \Xi(\Lambda_2, k_2, t_2)| \leq  ( S^{1/4} T^{-1/4} +T^{-\mathfrak{c}}) (\log N)^{\mathfrak{C}}.
\end{equation*}
\end{prop}

To show Proposition \ref{cor:xicont}, we require several lemmas that (with some of the intermediate statements involved in their proofs) will also be used in Section \ref{sec:levy_chent} below. Here, we recall the functionals $\Xi_1$, $\Xi^{[m]}$, and $\Xi_1^{[m]}$ from  Definitions \ref{xi2} and \ref{xim}. The following lemma, to be shown in Section \ref{subsec:xi2xim2_contproofs}, provides H\"older type continuity bounds for $\Xi_1$ and $\Xi_1^{[m]}$.

\begin{lem}[Continuity of $\Xi_1$ and $\Xi_1^{[m]}$]\label{lem:xi_xim_cont}

Let $S \in [\mathfrak{M}, T]$ be a real number. There exist constants $\mathfrak{c} > 0$ and $\mathfrak{C}>1$ such that the following two statements hold with overwhelming probability, for all real numbers $t_1, t_2 \in [0, T \log N]$ with $|t_1 - t_2| \le S$. 

\begin{enumerate} 

\item For any integers $k_1, k_2 \in \llbracket - T^{3/2},  T^{3/2} \rrbracket$ with $|k_1-k_2| \leq S$ and real numbers $\Lambda_1, \Lambda_2 \in [- \log N, \log N]$ with $|\Lambda_1 - \Lambda_2| \leq ST^{-1}$, we have 
\begin{equation}\label{eqn:xi_2_bound1}
|\Xi_1(\Lambda_1, k_1, t_1)- \Xi_1(\Lambda_2, k_2, t_2)| \leq ( S^{1/4} T^{-1/4} +T^{-\mathfrak{c}}) (\log N)^{\mathfrak{C}}.
\end{equation}

\item For any integer $m \in \llbracket 0, (\log N)^{1/10} \rrbracket$ and real numbers $q_1, q_2,q_1',q_2' \in [- (\log N)^8 T, (\log N)^8 T]$ with $|q_1-q_2| \le S$ and $|q_1'-q_2'| \le S$, we have 
\begin{equation}\label{eqn:xim_2_cont}
   | \Xi_1^{[m]}(q_1,q_1', t_1) - \Xi_1^{[m]}(q_2,q_2', t_2)| \leq S ^{1/2} T^{-1/2}(\log N)^{m+\mathfrak{C}}.
\end{equation}

\end{enumerate} 
\end{lem}

For the next lemma, define the function $\mathcal{G}_N : \mathbb{R}^3 \times \mathbb{R}$ implicitly appearing in \eqref{eqn:psiZktauapprox}, by setting 
\begin{equation}\label{eqn:GN_def}
    \mathcal{G}_N(\Lambda, k, t) \coloneqq - (\alpha T^{1/2})^{-1} \sum_{r = \lceil N_1/R \rceil}^{\lfloor N_2/R \rfloor-1} c^{(\Lambda)}(r,k,t) \sum_{i=r R}^{(r+1)R-1} X_i,
\end{equation}

\noindent for all $\Lambda, k, t \in \mathbb{R}$, where we recall the coefficients $ c^{(\Lambda)}(r,k,t)$ from \eqref{eqn:clambda1} (and the integer $R$ from Assumption \ref{ass:R_M_ass} and the stretch random variables $X_i$ from Definition \ref{xi}). The following lemma, to be shown in Section \ref{ProofContinuousG}, states a continuity bound for $\mathcal{G}_N$.

\begin{lem}[Continuity of $\mathcal{G}_N$]\label{lem:Xi_brownian_cont_allparam_pl}

    Let $S \in [\mathfrak{M}, T]$ be a real number. There exists a constant $\mathfrak{C}>1$ such that following holds with overwhelming probability. For any integers $k_1, k_2 \in \llbracket - T^{3/2},  T^{3/2} \rrbracket $ with $|k_1-k_2| \leq S $; real numbers $t_1, t_2 \in [0,  T \log N]$ with $|t_1 - t_2| \leq S$, and real numbers $\Lambda_1, \Lambda_2 \in [- \log N, \log N]$ such that $|\Lambda_1 - \Lambda_2| \leq \mathcal{S}/T$, we have
    \begin{equation}\label{eqn:xibrownian_cont_bd_all}
         |\mathcal{G}_N(|\Lambda,k_1,t_1) -  \mathcal{G}_N(\Lambda_2,k_2,t_2)| \leq   S^{1/4} T^{-1/4} (\log N)^{\mathfrak{C}}.
    \end{equation}
\end{lem}

For the next lemma, for any function $F: \mathbb{R} \rightarrow \mathbb{R}$, define the function $\mathcal{G}_N^{(F)}: \mathbb{R}^3 \rightarrow \mathbb{R}$ by setting 
\begin{equation}\label{eqn:GNF_gen}
    \mathcal{G}_N^{(F)}(q,q',t) \coloneqq - (\alpha T)^{-1/2} \sum_{r = \lceil N_1/R \rceil}^{\lfloor N_2/R \rfloor - 1} c^{(F)}(r,q,q',t) \sum_{i=r R}^{(r+1)R-1} X_i,
\end{equation}

\noindent for all $q,q',t \in \mathbb{R}$ where (for any $r \in \mathbb{R}$) the coefficients $c^{(F)} (r,q,q',t)$ are defined by 
\begin{equation}\label{eqn:cmF}
c^{(F)}(r,q,q',t) = \int_{-\infty}^{\infty}F(\lambda) \cdot \big(\chi(q-\alpha r R)- \chi(q' -\alpha rR - t\ve(\lambda)) \big) \varrho(\lambda) d \lambda .
\end{equation}

\noindent Observe in particular that for all $r,q,q',t \in \mathbb{R}$ and $m \in \mathbb{Z}_{\ge 0}$ we have $c^{(\varsigma_m)} (r,q,q',t) = c^{[m]} (r,q,q',t)$, where $c^{[m]}$ is given by \eqref{eqn:cm1} (and $\varsigma_m$ by \eqref{nx}), and that $G_N^{(\varsigma_m)}$ appears implicitly in \eqref{eqn:lambdam_approx}.

The following lemma, to be shown in Section \ref{subsec:brownian_cont_proofs} (and which will be relevant for the proof of Theorem \ref{thm:couplethm12}), provides a continuity bound for $\mathcal{G}_N^{(F)}$ from \eqref{eqn:GNF_gen}.

\begin{lem}[Continuity of $\mathcal{G}_N^{(F)}$]\label{lem:gen_brownian_cont_allparam_prelim}

    There exists a constant $\mathfrak{C}>1$ such that the following holds. Let $A> 0$ and $S \in [\mathfrak{M},T]$ be real numbers, and let $F: \mathbb{R} \rightarrow \mathbb{R}$ be a function satisfying Assumption \ref{ass:F}. With overwhelming probability, we have that
    \begin{equation}\label{eqn:Fbrownian_cont_bd_all}
         |\mathcal{G}_N^{(F)}(q_1,q_1',t_1) -  \mathcal{G}_N^{(F)}(q_2,q_2',t_2)| \leq A S^{1/4} T^{-1/4} (\log N)^{\mathfrak{C}},
    \end{equation}

    \noindent  for all real numbers $t_1, t_2 \in [0,  T \log N]$ with $|t_1 - t_2| \leq S$, and all real numbers $q_1, q_2,q_1',q_2' \in [- T(\log N)^8,  T (\log N)^8] $ with $|q_1-q_2| \leq S $ and $|q_1'-q_2'| \leq S $.
\end{lem}

The next lemma, to be shown in Section \ref{ProofContinuousXi}, states that the bounds \eqref{eqn:psiZktauapprox} and \eqref{eqn:lambdam_approx} likely hold uniformly in the parameters $(\Lambda, k, t)$ and $(q, q',t)$ (as opposed to only at a single point). 

\begin{lem}\label{lem:xi_decomp_lem}
    There exists a constant $\mathfrak{c}>0$ such that, with overwhelming probability, we have for all $t \in [0, T \log N]$ that
    \begin{flalign}
     &  \displaystyle\sup_{|\Lambda| \le \log N} \displaystyle\sup_{|k| \le T (\log N)^4} \left|  \Xi(\Lambda, k, t) - \Xi_1(\Lambda, k, t) - \mathcal{G}_N(\Lambda, k,t ) \right| \leq 
       T^{-\mathfrak{c}}; \label{eqn:xi12_dif} \\
      & \displaystyle\max_{m \in \llbracket 0, (\log N)^{1/10} \rrbracket}  \displaystyle\sup_{|q|, |q'| \le T(\log N)^5} \big|  \Xi^{[m]}(q, q',t) - \Xi_1^{[m]}(q, q',t) - \mathcal{G}_N^{(\varsigma_m)}(q,q',t) \big| \leq 
       T^{-\mathfrak{c}}.\label{eqn:xim12_dif}
    \end{flalign}
\end{lem}

\begin{proof}[Proof of Proposition \ref{cor:xicont}]
    This follows from Lemmas \ref{lem:xi_xim_cont}, \ref{lem:Xi_brownian_cont_allparam_pl}, and \ref{lem:xi_decomp_lem}, 
\end{proof}

\subsection{Proofs of Lemma \ref{lem:gen_brownian_cont_allparam_prelim}}
\label{subsec:brownian_cont_proofs}

To show Lemma \ref{lem:gen_brownian_cont_allparam_prelim}, we  begin with the several intermediate lemmas.

\begin{lem}\label{lem:gen_brownian_cont_prelim_close}
Let $A > 0$ be a real number and $F: \mathbb{R} \rightarrow \mathbb{R}$ be a function satisfying Assumption \ref{ass:F}. Let $S \ge 0$ be a real number; $t_1, t_2 \in \mathbb{R}$ be real numbers with $|t_1 - t_2| \leq S$, and $q_1, q_2,q_1',q_2' \in \mathbb{R}$ be real numbers with $|q_1-q_2| \leq S $ and $|q_1'-q_2'| \leq S $. Then, for any real number $r \in \mathbb{R}$, we have 
  \begin{align}
  \label{c2f} 
      |c^{(F)}(r,q_1,q_1',t_1) - c^{(F)}(r,q_2,q_2',t_2)| &\leq A S \mathfrak{M}^{-1}.
        \end{align}
\end{lem}

\begin{proof}
    This follows from the definition \eqref{eqn:cmF} of $c^{(F)}$, together with the facts that $|\chi'| \le 10 \mathfrak{M}^{-1}$ from Definition \ref{def:smoothed_log}; that $|F(\lambda)| \cdot \varrho(x) \le C Ae^{-cx^2}$ from Lemma \ref{lem:varrho_bd} and Assumption \ref{estimatef}; and that $|\ve(\lambda)| \le C \log(|\lambda|+1)$ from Lemma \ref{lem:ve_tail}. 
\end{proof}

To state the next lemma, for any real numbers $q \in \mathbb{R}$ and $S \geq 1$ define (recalling the parameter $R$ from \eqref{ra0}) the interval $\mathcal{B}_q = \mathcal{B}_q (S)$, on which certain estimates will not hold, by
\begin{equation}\label{eqn:Bqdef}
    \mathcal{B}_q \coloneqq [q(\alpha R)^{-1} -(\log N)^5 (ST)^{1/2} R^{-1}, q(\alpha R)^{-1}+(\log N)^5 (ST)^{1/2} R^{-1}].
\end{equation}

\begin{lem}\label{lem:gen_brownian_cont}

There exists a constant $\mathfrak{C}>1$ such that the following holds. Let $A> 0$ be a real number and $F : \mathbb{R} \rightarrow \mathbb{R}$ be a function satisfying Assumption \ref{ass:F}. In addition, let $S \in [\mathfrak{M}, T]$ be a real number; let $t_1, t_2 \in [0,  T \log N]$ be real numbers with $|t_1 - t_2| \leq S$; and let $q_1, q_2,q_1',q_2' \in [- T(\log N)^8,  T (\log N)^8] $ be real numbers with $|q_1-q_2| \leq S $ and $|q_1'-q_2'| \leq S $.
 
 \begin{enumerate} 
 \item For any real number $r \in \mathbb{R} \setminus (\mathcal{B}_{q_1} \cup \mathcal{B}_{q_1'})$, we have
  \begin{align}
     |c^{(F)}(r,q_1,q_1',t_1) - c^{(F)}(r,q_2,q_2',t_2)| &\leq A  S^{1/2} T^{-1/2} (\log N)^{\mathfrak{C}}. 
     \label{eqn:cFq_bd}
        \end{align}
    \item With overwhelming probability, we have 
    \begin{equation}\label{eqn:Fbrownian_cont_bd}
         |\mathcal{G}_N^{(F)}(q_1,q_1',t_1) -  \mathcal{G}_N^{(F)}(q_2,q_2',t_2)| \leq A  S^{1/4} T^{-1/4} (\log N)^{\mathfrak{C}}.
    \end{equation}
    \end{enumerate} 
\end{lem}

\begin{proof}

    To show \eqref{eqn:cFq_bd}, it suffices by \eqref{eqn:cmF} to verify the following two bounds. First, for any real numbers $t \in [0,  T \log N]$; $q_1,q_2 \in [- T(\log N)^8,  T (\log N)^8] $ with $|q_1-q_2| \leq S $; and $r \notin \mathcal{B}_{q_1}$, we have 
    \begin{equation}\label{eqn:q_cont_toshow}
     \left|   \int_{-\infty}^{\infty}F(\lambda) \big(\chi(q_1-\alpha rR - t\ve(\lambda) )- \chi(q_2 -\alpha rR - t\ve(\lambda)) \big) \varrho(\lambda) d \lambda \right| \leq  A S^{1/2}T^{-1/2} (\log N)^C. 
    \end{equation}
    
     \noindent Second, for any real numbers $t_1,t_2 \in [0,  T \log N]$ with $|t_1-t_2| \leq S $; $q \in  [- T(\log N)^8,  T (\log N)^8] $; and $r \notin \mathcal{B}_q$, we have 
     \begin{equation}\label{eqn:t_cont_toshow}
     \left|   \int_{-\infty}^{\infty}F(\lambda) \big(\chi(q-\alpha rR - t_1\ve(\lambda))- \chi(q -\alpha rR - t_2\ve(\lambda) ) \big) \varrho(\lambda) d \lambda \right| \leq  A S^{1/2}T^{-1/2} (\log N)^C .
    \end{equation}
   
     To show \eqref{eqn:q_cont_toshow}, first suppose that $t \geq (T S)^{1/2}$. Then denoting the interval
$$\mathcal{I} \coloneqq [-\log N, \log N] \cap \ve^{-1}\left( [ (q_1 - \alpha r R - 10 S) t^{-1}, (q_1 - \alpha r R + 10 S) t^{-1}] \right),  $$

\noindent observe that its size is bounded by 
\begin{equation}\label{eqn:FmathcalIbd}
    |\mathcal{I}| \leq C \log N \cdot St^{-1} \leq (\log N)^2 S^{1/2} T^{-1/2},
\end{equation}

\noindent since $\inf_{|\lambda| \le \log N} \ve'(\lambda) \geq c (\log N)^{-1}$ (by Lemma \ref{lem:veff_inc}). Thus, since $|\chi| \le 1$ and $\supp \chi' \subseteq [-\mathfrak{M}, \mathfrak{M}] \subseteq [-S, S]$ by Definition \ref{def:smoothed_log}, we may upper bound the left hand side of \eqref{eqn:q_cont_toshow} as
\begin{flalign}\label{eqn:Fnew_decomp}
    2 \int_{\mathcal{I}} | F(\lambda) | \cdot \varrho(\lambda) d \lambda +  2\int_{|\lambda| > \log N} |F(\lambda) |  \varrho(\lambda) d \lambda \le CA S^{1/2} T^{-1/2} (\log N)^2, 
\end{flalign}

\noindent due to \eqref{eqn:FmathcalIbd} and the fact that $|F(\lambda)| \cdot \varrho (\lambda) \le CA e^{-c|\lambda|^2}$ by Lemma \ref{lem:varrho_bd} and \eqref{estimatef}.

This verifies \eqref{eqn:q_cont_toshow} when $t \ge (TS)^{1/2}$, so suppose instead that $t \leq (T S)^{1/2}$. Then, observe that $\chi(q_1-\alpha rR - t\ve(\lambda) )- \chi(q_2 -\alpha rR - t\ve(\lambda)) = 0$ if $\lambda \in [-\log N, \log N]$. Indeed, since $r \notin \mathcal{B}_{q_1}$ and $\sup_{|\lambda| \leq \log N}|\ve(\lambda)| \leq C \log N$ (by Lemma \ref{lem:ve_tail}), we have  
\begin{multline}
|q_1-\alpha rR - t\ve(\lambda)| \geq |q_1-\alpha rR| - t |\ve(\lambda)| \geq  (\log N)^5 (T S)^{1/2} - C (T S)^{1/2} (\log N)^2 \\
> \frac{1}{2}(\log N)^5 (T S)^{1/2} >10 S.
\end{multline}

\noindent Since $|q_1-q_2| \le S$ and $\supp \chi' \subseteq [-\mathfrak{M}, \mathfrak{M}] \subseteq [-S, S]$, it follows that $\chi(q_1-\alpha rR - t\ve(\lambda) )- \chi(q_2 -\alpha rR - t\ve(\lambda)) = 0$. Thus, the left side of 
\eqref{eqn:q_cont_toshow} is bounded above by 
\begin{flalign}
    \label{integralfrho0lambda}
2 \int_{|\lambda| > \log N} |F(\lambda) |  \varrho(\lambda) d \lambda \le CAe^{-c (\log N)^2} \le CA S^{1/2} T^{-1/2} (\log N)^2,
\end{flalign}

\noindent again using the facts that $|\chi| \le 1$ and $|F(\lambda)| \cdot \varrho (\lambda) \le CA e^{-c|\lambda|^2}$. This confirms \eqref{eqn:q_cont_toshow}.

Now we must show \eqref{eqn:t_cont_toshow}. Assume without loss of generality that $t_1 < t_2$ and that $q \ge \alpha rR$. Set $\Delta \coloneqq q - \alpha r R$. Observe that $|\Delta| \ge (ST)^{1/2}$, since $r \notin \mathcal{B}_q$. We separately consider the cases when $\Delta (\log N)^{-2} \leq t_2 $ and when $t_2 \leq \Delta (\log N)^{-2}$. 

Assume first that $\Delta (\log N)^{-2} \leq t_2$. Then, since $r \notin \mathcal{B}_q$, we have $S \leq \Delta (\log N)^{-4}$, which implies that $t_1 \ge t_2 - S \ge t_2/2 \ge \Delta (\log N)^{-2}/2$. Thus,
    \begin{flalign}
        \label{delta0} 
     \Delta (t_1^{-1} - t_2^{-1}) \leq  S \Delta t_1^{-2} \leq 4 S\Delta^{-1} (\log N)^2; \qquad S t_1^{-1} \le 2St_2^{-1} \le 2S \Delta^{-1} (\log N)^2 .
    \end{flalign}
    Therefore, since $\supp \chi' \subseteq [-S,S]$, we only have $\chi(q-\alpha r R - t_1 \ve(\lambda)) - \chi(q - \alpha r R - t_2 \ve (\lambda)) \ne 0$ in \eqref{eqn:t_cont_toshow} if 
    $$\ve(\lambda) \in I \coloneqq [(q - \alpha r R) t_2^{-1} - St_1^{-1},  (q - \alpha r R) t_1^{-1} + St_1^{-1}].$$
    
     \noindent By \eqref{delta0}, we have $|I| \leq 8 S\Delta^{-1} (\log N)^2 \le 8S^{1/2} T^{-1/2} (\log N)^2$ (as $\Delta \ge (ST)^{1/2}$). Thus, denoting 
    $$\mathcal{I}' \coloneqq \ve^{-1}(I) \cap [-\log N, \log N] ,$$ 
    we have $|\mathcal{I}'| \le CS^{1/2} T^{-1/2} (\log N)^3$, where we have again used the bound $\inf_{|\lambda| \le \log N} \ve'(\lambda) \ge c(\log N)^{-1}$ (from Lemma \ref{lem:veff_inc}). Hence, as in \eqref{eqn:Fnew_decomp}, the left hand side of \eqref{eqn:t_cont_toshow} is at most  
    \begin{flalign*}
    2 \int_{\mathcal{I}'} | F(\lambda) | \cdot \varrho(\lambda) d \lambda +  2\int_{|\lambda| > \log N} |F(\lambda) |  \varrho(\lambda) d \lambda \le CA S^{1/2} T^{-1/2} (\log N)^3.
    \end{flalign*}

    This shows \eqref{eqn:t_cont_toshow} if $\Delta (\log N)^{-2} \le t_2$, so assume instead that $\Delta (\log N)^{-2} \ge t_2$. In this case, since $t_1 \le t_2$; $\Delta \ge (ST)^{1/2}$; and $\sup_{|\lambda| \le \log N} |\ve(\lambda)| \le C \log N$ (by Lemma \ref{lem:ve_tail}), we have that $q - \alpha r R - t_i \ve (\lambda) \ge \Delta - C\Delta (\log N)^{-1} \ge (ST)^{1/2}/2 \ge 10\mathfrak{M}$, for each $i \in \{ 1, 2 \}$. Together with the fact that $\supp \chi' \subseteq [-\mathfrak{M}, \mathfrak{M}]$, this implies that the integrand in \eqref{eqn:t_cont_toshow} vanishes for $\lambda \in [- \log N, \log N]$. Hence (since $|\chi| \le 1$), the left hand side of \eqref{eqn:t_cont_toshow} is  bounded above by \eqref{integralfrho0lambda}, which again confirms \eqref{eqn:t_cont_toshow}. This establishes \eqref{eqn:cFq_bd}. 

    It remains to verify \eqref{eqn:brownian_cont_bd_lam}. By \eqref{eqn:GNF_gen}, it suffices to show that 
    \begin{flalign}
    \label{sumr0c}
        \Bigg| \sum_{r = \lceil N_1/R \rceil}^{\lfloor N_2/R \rfloor - 1} \big(c^{(F)}(r,q_1,q_1',t_1)-c^{(F)}(r,q_2,q_2',t_2) \big)  \sum_{i=r R}^{(r+1)R-1} X_i \Bigg| \le A (ST)^{1/4} (\log N)^C.
    \end{flalign}

    \noindent Let us separately estimate the contribution to the left side of \eqref{sumr0c} over $r \in \mathcal{B}_{q_1} \cup \mathcal{B}_{q_1'}$ and over the remaining $r \in \llbracket N_1/R, N_2/R-1\rrbracket$. For the first, since $|\chi| \le 1$ and $|F(\lambda)| \cdot \varrho (\lambda) \le CA e^{-c|\lambda|^2}$, we find from \eqref{eqn:cmF} that $|c^{(F)} (r,q,q',t)| \le CA$, for all $r,q,q',t \in \mathbb{R}$. Thus, since $R \cdot |(\mathcal{B}_{q_1} \cup \mathcal{B}_{q_2}) \cap \mathbb{Z}| \le 5(ST)^{1/2} (\log N)^5$ (by \eqref{eqn:Bqdef}), using a Chernoff bound and $\mathbb{E}[X_i] = 0$ (by \eqref{xi}), we deduce with overwhelming probability that 
    \begin{flalign}\label{eqn:dom_err}
    \Bigg|\sum_{r \in \mathcal{B}_{q_1} \cup \mathcal{B}_{q_1'}} \sum_{i=r R}^{(r+1)R-1} X_i \cdot \big (c^{(F)}(r,q_1,q_1',t_1)-c^{(F)}(r,q_2,q_2',t_2) \big)  \Bigg| & \le A (ST)^{1/4} (\log N)^C. 
    \end{flalign}

    To estimate the remaining summands in \eqref{sumr0c}, set $\mathcal{J} = \llbracket- T R^{-1} (\log N)^{10}, TR^{-1} (\log N)^5 \rrbracket$. Observe that, for any index $r \notin \setminus \mathcal{J}$ and real number $q \in \{ q_1, q_2, q_1', q_2' \}$ that $|q - \alpha rR| \ge \alpha T(\log N)^{10} - T(\log N)^8 \ge T(\log N)^9$. Hence, by Lemma \ref{lem:ve_tail} and the fact that $\supp \chi' \subseteq [-S,S]$, there exist some $r \notin \mathcal{J}$; $Q_1, Q_2 \in \{ q_1, q_2, q_1', q_2' \}$; and $s_1, s_2 \in [0, T \log N]$, for which $\chi (Q_1 - \alpha rR - s_1 \ve(\lambda)) - \chi (Q_2 -  \alpha rR - s_2 \ve(\lambda)) \ne 0$, only if $\lambda \ge \log N$. Thus, by \eqref{eqn:cmF} (and the fact that $|\chi| \le 1$), we have that $|c^{(F)} (r,q_1,q_1',t_1) - c^{(F)} (r,q_2, q_2', t_2)|$ is bounded above by \eqref{integralfrho0lambda} and is therefore at most $CA e^{-c(\log N)^2}$. As such, since $|X_i| \le (\log N)^2$ holds with overwhelming probability for all $i \in \llbracket N_1, N_2-1 \rrbracket$ by Lemma \ref{lem:q_spacing}, we have with overwhelming probability that
    \begin{flalign}
        \label{sumr0} 
         \sum_{\substack{r \in \llbracket N_1/R, N_2/R \rrbracket \\ r \notin \mathcal{J}}} \sum_{i=r R}^{(r+1)R-1} |X_i| \cdot \big| c^{(F)}(r,q_1,q_1',t_1)-c^{(F)}(r,q_2,q_2',t_2) \big| \le C e^{-c(\log N)^2}.
    \end{flalign}

    \noindent For the remaining terms, we have that 
    \begin{flalign}
        \label{sumr1}
        \begin{aligned}
        & \sum_{\substack{r \in  \mathcal{J} \\ r \notin \mathcal{B}_{q_1} \cup \mathcal{B}_{q_1'}}} \sum_{i=r R}^{(r+1)R-1} X_i \cdot \big( c^{(F)}(r,q_1,q_1',t_1)-c^{(F)}(r,q_2,q_2',t_2) \big) \\
        & \qquad \qquad \qquad \le AS^{1/2} T^{-1/2} (\log N)^C \cdot T^{1/2} (\log N)^C \le A(ST)^{1/4} (\log N)^{2C},
        \end{aligned}
    \end{flalign}

    \noindent where in the first estimate we applied \eqref{eqn:cFq_bd} and a Chernoff bound (using the fact from \eqref{xi} that $\mathbb{E}[X_i] = 0$), together with the bound $R \cdot |\mathcal{J}| \le 3T (\log N)^{10}$; in the second, we used the fact that $S \le T$. We then deduce \eqref{sumr0c}, and thus the lemma, by combining \eqref{eqn:dom_err}, \eqref{sumr0}, and \eqref{sumr1}.
\end{proof}

Now, we are in a position to prove Lemma \ref{lem:gen_brownian_cont_allparam_prelim}.

\begin{proof}[Proof of Lemma \ref{lem:gen_brownian_cont_allparam_prelim}]
    Set $\mathcal{U} =  [- T(\log N)^8,  T (\log N)^8] \cap (N^{-20} \cdot \mathbb{Z})$ and $\mathcal{V} = [0, T \log N] \cap (N^{-20} \cdot \mathbb{Z})$. By Lemma \ref{lem:gen_brownian_cont} and a union bound, there exists an event $\mathsf{E}_1$ on which \eqref{eqn:Fbrownian_cont_bd_all} holds for all $q_1, q_2, q_1',q_2' \in \mathcal{U}$ with $|q_1-q_2| \leq S$ and $|q_1'-q_2'| \leq S$, and all $t_1,t_2 \in \mathcal{V}$ with $|t_1-t_2| \leq S$. In addition, by Lemma \ref{lem:q_spacing}, there exists an overwhelmingly probable event $\mathsf{E}_2$ on which $|X_i| \leq (\log N)^2$ holds $i \in \llbracket N_1, N_2 - 1 \rrbracket$. We therefore restrict to $\mathsf{E} = \mathsf{E}_1 \cap \mathsf{E}_2$ in what follows. Lemma \ref{lem:gen_brownian_cont_prelim_close} and \eqref{eqn:GNF_gen} then imply on $\mathsf{E}$ that, for any $q, q' \in [-T(\log N)^8, T(\log N)^8]$ and $t \in [0, T \log N]$, we have $|\mathcal{G}_N (q,q',t) - \mathcal{G}_N (\tilde{q}, \tilde{q}', \tilde{t})| \le AN^{-18}$, where $\tilde{q}, \tilde{q}' \in \mathcal{U}$ and $\tilde{t} \in \mathcal{V}$ are the closest elements to $q$, $q'$, and $t'$, respectively. Combining these estimates yields \eqref{eqn:Fbrownian_cont_bd_all} in general, thereby establishing the lemma.
\end{proof}

\subsection{Proof of Lemma \ref{lem:Xi_brownian_cont_allparam_pl} and Lemma \ref{lem:xi_decomp_lem}}

\label{ProofContinuousG}
\label{ProofContinuousXi}

To prove Lemma \ref{lem:Xi_brownian_cont_allparam_pl} and Lemma \ref{lem:xi_decomp_lem}, we first require several intermediate lemmas.

\begin{lem} \label{lem:lc}

There exists $\mathfrak{C}>1$ so that, for any real number $\xi > 0$ there is a constant $\mathfrak{c} = \mathfrak{c}(\xi)>0$, for which the following holds with probability at least $1 - \mathfrak{c}^{-1} e^{-\mathfrak{c} (\log N)^2}$. For any integer $k \in \llbracket - T^{3/2},  T^{3/2} \rrbracket $, and real numbers $t \in [0,  T \log N]$ and $\Lambda, \Lambda' \in [-\log N, \log N]$ with $|\Lambda-\Lambda'| \leq e^{-\xi (\log N)^2}$, we have
   \begin{flalign}
      & |\Xi(\Lambda,k,t) - \Xi(\Lambda',k,t)| \leq  T^{-1/2} (\log N)^{\mathfrak{C}}; \label{eqn:lambda_weak_simul0} \\
       &  |\Xi_1(\Lambda,k,t) - \Xi_1(\Lambda',k,t)| \leq  T^{-1/2} (\log N)^{\mathfrak{C}};  \label{eqn:lambda2_weak_simul} \\
        & |\mathcal{G}_N(\Lambda,k,t) - \mathcal{G}_N(\Lambda',k,t)| \leq  T^{-1/2} (\log N)^{\mathfrak{C}}. \label{eqn:GNlambda_weak_simul}
   \end{flalign}
\end{lem}
\begin{proof}

Throughout this proof, we will let the constants $c>0$ and $C>1$ (which might change between lines) depend on $\xi$.
    To show \eqref{eqn:lambda_weak_simul0} and \eqref{eqn:lambda2_weak_simul}, define the event $\mathsf{E} = \mathsf{BND}_{\mathbf{L}} (\log N) \cap \mathrm{SEP}_{\L}(e^{-\xi (\log N)^2/2})$ (recall Definition \ref{eigenvaluesm} and Definition \ref{lambda0lambdaevent}), which by Lemma \ref{lem:bd_lem}
 and Lemma \ref{lem:sep_lem} satisfies $\mathbb{P}[\mathsf{E}] \ge 1- c^{-1} e^{-c(\log N)^2}$. We therefore restrict to $\mathsf{E}$ in what follows. 
 
    On $\mathsf{E}$, for any $\Lambda \in [- \log N, \log N]$, there can be at most one eigenvalue $\Lambda_i \in \eig \mathbf{L}$ satisfying $|\Lambda_i-\Lambda| \leq e^{-\xi (\log N)^2/2}/2$. In addition, if $|\Lambda-\Lambda'| \leq e^{-\xi (\log N)^2}$ and $|\Lambda_i-\Lambda| > e^{-\xi (\log N)^2/2}/2$, then $|\Lambda_i-\Lambda' | > e^{-\xi (\log N)^2/2}/4$. Since $|\mathfrak{l}'(\lambda) | \leq |\lambda|^{-1}$ by \eqref{eqn:sl}, this yields
\begin{equation}\label{eqn:ldbd0}
    |\mathfrak{l}(\Lambda-\Lambda_i) - \mathfrak{l}(\Lambda'-\Lambda_i)| \leq 4 e^{\xi (\log N)^2/2} \cdot |\Lambda-\Lambda'| \leq C e^{-\xi/2 (\log N)^2}.
\end{equation}
Moreover, since $|\chi'| \le 1$ (by Definition \ref{def:chi_def}) and $\sup_{|\lambda| \le \log N} |\ve'(\lambda)| \le C \log N$ (by Lemma \ref{lem:ve_tail}), if $|\Lambda-\Lambda'| \leq e^{-\xi (\log N)^2}$, then we have
\begin{equation}\label{eqn:chidiffbd0}
\big| \chi \big(\alpha (k - i)+ t  (\ve(\Lambda) - \ve(\Lambda_i) ) \big)-\chi\big(\alpha (k - i)+ t  (\ve(\Lambda') - \ve(\Lambda_i) ) \big) \big| \leq C e^{-\xi (\log N)^2/2}.
    \end{equation}

Therefore, for any real numbers $t \in [0,  T \log N]$ and $\Lambda, \Lambda' \in [-\log N, \log N]$ with $|\Lambda- \Lambda'| \leq e^{-\xi (\log N)^2}$, and any integer $k \in \llbracket - T^{3/2},  T^{3/2} \rrbracket $, we have on $\mathsf{E}$ that 
\begin{equation*}
     |\Xi_1(\Lambda,k,t) - \Xi_1(\Lambda',k,t)| \leq Ce^{-c (\log N)^2} + T^{-1/2} (\log N)^C.
\end{equation*}

\noindent Indeed, the first term on the right hand side above arises from applying \eqref{eqn:ldbd0} and \eqref{eqn:chidiffbd0} to bound the $i$-th summand from \eqref{eqn:Xi2_outline} in the difference $\Xi_1(\Lambda,k,t) - \Xi_1(\Lambda',k,t)$ whenever $|\Lambda_i-\Lambda| > e^{-\xi (\log N)^2/2}/4$; the second term corresponds to the at most one value of $i$ for which $|\Lambda_i-\Lambda| \leq e^{-\xi (\log N)^2/2}/2$, using the fact from \eqref{eqn:sl} that $|\mathfrak{l}(\Lambda_i - \Lambda)| + |\mathfrak{l}(\Lambda_i - \Lambda')| \le (\log N)^3$. This confirms \eqref{eqn:lambda_weak_simul0}; the proof that \eqref{eqn:lambda2_weak_simul} holds on $\mathsf{E}$ is entirely analogous and is therefore omitted. 

It remains to verify \eqref{eqn:GNlambda_weak_simul}; we claim it holds on the event $\mathsf{E}' = \bigcap_{i=N_1}^{N_2-1} \{ |X_i| \le (\log N)^2 \}$, which satisfies $\mathbb{P}[\mathsf{E}'] \ge 1 - c^{-1} e^{-c(\log N)^2}$, by Lemma \ref{lem:q_spacing}; we restrict to $\mathsf{E}'$ below. To show this, it suffices by \eqref{eqn:GN_def} to verify for any $\Lambda, \Lambda' \in [-\log N, \log N]$ with $|\Lambda - \Lambda'| \leq e^{-\xi (\log N)^2}$ that 
\begin{flalign*}
   \displaystyle\sup_{r,k,t \in \mathbb{R}} |c^{(\Lambda)} (r,k,t) - c^{(\Lambda')} (r,k,t)| \le c^{-1} e^{-c(\log N)^2}.
\end{flalign*}

\noindent To that end, observe if $|\lambda - \Lambda|> e^{-\xi (\log N)^2/2}$, then $|\lambda - \Lambda'|> e^{-\xi(\log N)^2/2}/2$. Thus, since $|\mathfrak{l}'(\lambda) |  \leq |\lambda|^{-1}$ and $|\chi'| \le 10\mathfrak{M}^{-1}$, we have for any $r, k \in \mathbb{R}$ and $t \in [0, T \log N]$ that 
\begin{multline}
    \Bigg| 2 \mathfrak{l}(\Lambda - \lambda) \cdot \big( \chi(\alpha k -\alpha r R)
- \chi \big(\alpha k - \alpha r R + t (\ve(\Lambda) - \ve(\lambda)) \big) \big) \\
-  2 \mathfrak{l}(\Lambda' - \lambda) \cdot \big( \chi(\alpha k -\alpha r R)
- \chi \big(\alpha k - \alpha r R + t (\ve(\Lambda) - \ve(\lambda)) \big) \big) \Bigg| \leq e^{-c (\log N)^2},
\end{multline}

\noindent whenever $\lambda \in [-\log N, \log N]$ satisfies $|\lambda - \Lambda|> e^{-\xi (\log N)^2 /2}$. By \eqref{eqn:clambda1} and the fact that $|\chi| \le 1$, it therefore suffices to show for $\Lambda_0 \in \{ \Lambda, \Lambda' \}$ that 
\begin{flalign*}
    \displaystyle\int_{-\infty}^{\infty} (\mathbbm{1}_{|\lambda| > \log N} + \mathbbm{1}_{|\lambda - \Lambda| \le  e^{-\xi (\log N)^2}}) \cdot |\mathfrak{l}(\Lambda_0 - \lambda)| \cdot \varrho (\lambda) d \lambda \le c^{-1} e^{-c(\log N)^2},
\end{flalign*}

\noindent which holds from the definition \eqref{eqn:sl} of $\mathfrak{l}$ and the fact that $\varrho(\lambda) \le c^{-1} e^{-c(\log N)^2}$ (by Lemma \ref{lem:varrho_bd}).
\end{proof}

For the next lemma, we recall the set $\mathcal{B}_q$ from \eqref{eqn:Bqdef}.

\begin{lem}\label{lem:Xi_brownian_cont}
  There exists a constant $\mathfrak{C}>1$ such that the following holds. Let $ S \in [T^{1/2}, T]$ be a real number. Further let $k_1, k_2 \in \llbracket - T^{3/2},  T^{3/2} \rrbracket $ be integers with $|k_1-k_2| \leq S $, and let $t_1, t_2 \in [0,  T \log N]$ and $\Lambda_1, \Lambda_2 \in [- \log N, \log N]$ be real numbers with $|t_1 - t_2| \leq S$ and $|\Lambda_1 - \Lambda_2| \leq ST^{-1}$. 
  
  \begin{enumerate} 
  \item For any real number $r \in \mathbb{R} \setminus (\mathcal{B}_{\alpha k_1}\cup  \mathcal{B}_{\alpha k_1 + t_1 \ve(\Lambda_1)})$, we have 
  \begin{align}
      |c^{(\Lambda_1)}(r,k_1,t_1) - c^{(\Lambda_2)}(r,k_2,t_2)| &\leq S^{1/2} T^{-1/2} (\log N)^{\mathfrak{C}}. \label{eqn:clam_bd}
  \end{align}
  \item With overwhelming probability, we have 
  \begin{equation}\label{eqn:brownian_cont_bd_lam}
     | \mathcal{G}_N(\Lambda_1, k_1, t_1) -  \mathcal{G}_N(\Lambda_2, k_2, t_2)| \leq S^{1/4} T^{-1/4} (\log N)^{\mathfrak{C}}; 
  \end{equation}
  \end{enumerate} 
\end{lem}

\begin{proof}

    The proof of \eqref{eqn:brownian_cont_bd_lam} given \eqref{eqn:clam_bd} is entirely analogous to that of \eqref{eqn:Fbrownian_cont_bd} given \eqref{eqn:cFq_bd} and is therefore omitted. To show \eqref{eqn:clam_bd}, first observe by \eqref{eqn:sl} that, if $|\Lambda| \le \log N$, then $\lambda \mapsto \mathfrak{l}(\Lambda - \lambda)$ satisfies Assumption \ref{ass:F} with $A = (\log N)^3$. Thus by (the $(q_i, q_i') = (\alpha k_i, \alpha k_i + t_i \ve(\Lambda_i))$ and $F(\lambda) = \mathfrak{l}(\Lambda - \lambda)$ case of) \eqref{eqn:cFq_bd}, and the definitions \eqref{eqn:clambda1} and \eqref{eqn:cmF} of $c^{(\Lambda)}$ and $c^{(F)}$, to show \eqref{eqn:clam_bd}, it suffices to assume that $k_1=k_2$ and $t_1=t_2$; below, we let $k=k_1=k_2$ and $t=t_1=t_2$. 
    
    To that end, note since $|\chi| \le 1$ that 
    \begin{flalign}\label{eqn:lamdif_tbd}
    \begin{aligned}
        |& c^{(\Lambda_1)}(r,k,t) - c^{(\Lambda_2)}(r,k,t)| \\
        & \quad  \leq  2 \int_{-\infty}^{\infty}  \big|  \chi \big(\alpha k - \alpha r R + t (\ve(\Lambda_2) - \ve(\lambda)) \big) - \chi \big(\alpha k - \alpha r R + t (\ve(\Lambda_1) - \ve(\lambda)) \big) \big|   \\
        & \qquad \qquad \times |\mathfrak{l}(\Lambda_2 -\lambda)| \cdot \varrho(\lambda) d\lambda + 8 \int_{-\infty}^{\infty} |\mathfrak{l}(\Lambda_1 -\lambda) - \mathfrak{l}(\Lambda_2 -\lambda)| \cdot \varrho(\lambda) d\lambda.
        \end{aligned}
    \end{flalign}
    
    \noindent Next, observe that $|\ve(\Lambda_1) - \ve(\Lambda_2)| \leq ST^{-1} (\log N)^2$, by Lemma \ref{lem:ve_tail} and the facts that $\Lambda_1, \Lambda_2 \in [-\log N, \log N]$ and $|\Lambda_1 - \Lambda_2| \le S$.  Hence, the $(q_1, q_2) = (\alpha k + t \ve (\Lambda_1), \alpha k + t \ve (\Lambda_2))$ case of \eqref{eqn:cFq_bd} (with the $S$ equal to $S (\log N)^3$ here, as $t \le T \log N$) bounds the first integral in \eqref{eqn:lamdif_tbd} by 
    \begin{flalign}
    \label{integral00} 
    \begin{aligned} 
    2 \int_{-\infty}^{\infty}  \big|  \chi \big(\alpha k - \alpha r R + t (\ve(\Lambda_2) - \ve(& \lambda)) \big) - \chi \big( \alpha k - \alpha r R + t ( \ve(\Lambda_1) - \ve(\lambda)) \big) \big| \\
    & \qquad \times |\mathfrak{l}(\Lambda_2 -\lambda)| \cdot \varrho(\lambda) d\lambda  \le S^{1/2} T^{-1/2} (\log N)^C.
    \end{aligned} 
    \end{flalign} 

    To complete a proof of \eqref{eqn:clam_bd}, it remains to bound the integral on the second line of \eqref{eqn:lamdif_tbd}. To do so, observe that 
    \begin{flalign*}
        \displaystyle\int_{-\infty}^{\infty} |\mathfrak{l}(\Lambda_1-\lambda) - \mathfrak{l}(\Lambda_2-\lambda)| \cdot \varrho (\lambda) d \lambda & \le |\Lambda_1 - \Lambda_2| \displaystyle\int_{-\infty}^{\infty} \big( |\mathfrak{l}' (\Lambda_1 - \lambda)| + |\mathfrak{l}' (\Lambda_2-\lambda)| \big) \cdot \varrho (\lambda) d \lambda \\
        & \le CST^{-1} (\log N)^2,
    \end{flalign*}

    \noindent where in the last bound we used the estimates $|\mathfrak{l}' (\lambda)| \le C \min \{ |\lambda|^{-1}, e^{-10(\log N)^2} \}$ (by \eqref{eqn:sl}), $\varrho (\lambda) \le c^{-1} e^{-c\lambda^2}$ (by Lemma \ref{lem:varrho_bd}), and $|\Lambda_1-\Lambda_2| \le ST^{-1}$. This, together with \eqref{integral00}, verifies \eqref{eqn:clam_bd}. 
\end{proof}

Now we can show Lemma \ref{lem:Xi_brownian_cont_allparam_pl} and Lemma \ref{lem:xi_decomp_lem}.

\begin{proof}[Proof of Lemma \ref{lem:Xi_brownian_cont_allparam_pl}]
The proof will be similar to that of Lemma \ref{lem:gen_brownian_cont_allparam_prelim}. Lemma \ref{lem:Xi_brownian_cont} states \eqref{eqn:xibrownian_cont_bd_all}, except only at a fixed choice of parameters $(\Lambda,k,t)$. Specifically, Lemma \ref{lem:Xi_brownian_cont} and a union bound together imply the following, for a sufficiently small constant $\xi>0$. Let $\mathcal{U} = [-\log N, \log N] \cap (e^{-\xi (\log N)^2} \cdot \mathbb{Z})$ and $\mathcal{V} = [0, T \log N] \cap (e^{-\xi (\log N)^2} \cdot \mathbb{Z})$, and let $\mathsf{E}_0$ denote the event on which  \eqref{eqn:xibrownian_cont_bd_all} holds for all $(\Lambda,k,t) \in \mathcal{U} \times \llbracket -T(\log N)^4, T(\log N)^4 \rrbracket \times \mathcal{V}$. Then $\mathsf{E}_0$ holds with overwhelming probability. Further letting $\mathsf{E}_1$ denote the event on which \eqref{eqn:GNlambda_weak_simul} holds, we restrict to the event $\mathsf{E} = \mathsf{E}_0 \cap \mathsf{F} \cap \bigcap_{i=N_1}^{N_2-1} \{ |X_i| \le (\log N)^2 \}$ in what follows. By Lemmas \ref{lem:q_spacing} and \ref{lem:lc}, $\mathsf{E}$ is overwhelmingly probable. 

 We must then extend \eqref{eqn:xibrownian_cont_bd_all} to hold for parameters outside of the above mesh, which we do through a (weak) continuity argument. Specifically, for any parameters $\Lambda \in [-\log N, \log N]$; $k \in \llbracket -T(\log N)^4, T(\log N)^4 \rrbracket$; and $t \in [0, T \log N]$, observe  that 
 \begin{flalign} 
 \label{xixi201} 
 |\mathcal{G}_N (\Lambda, k, t) - \mathcal{G}_N (\Lambda, k, \tilde{t})| \le T^{-1/2} (\log N)^C; \qquad  |\mathcal{G}_N (\Lambda, k, t) - \mathcal{G}_N (\tilde{\Lambda}, k, t)| \le T^{-1/2} (\log N)^C,
 \end{flalign} 
 
 \noindent where $\tilde{\Lambda} \in \mathcal{U}$ and $\tilde{t} \in \mathcal{V}$ are the closest elements to $\Lambda$ and $t$, respectively. Indeed, the first bound follows from \eqref{c2f} (with the $(q_1, q_1', F(\lambda))$ there equal to $(\alpha k, \alpha k + t \ve (\Lambda), \mathfrak{l}(\Lambda-\lambda))$ here), \eqref{eqn:clambda1}, and \eqref{eqn:GN_def} (with the fact that each $|X_i| \le (\log N)^2$), and the second from \eqref{eqn:GNlambda_weak_simul}. Combining \eqref{xixi201} with the fact that \eqref{eqn:xibrownian_cont_bd_all} holds for all $(\Lambda, k, t) \in \mathcal{U} \times \llbracket -T(\log N)^4, T(\log N)^4 \rrbracket \times \mathcal{V}$, we obtain \eqref{eqn:xibrownian_cont_bd_all} for all $(\Lambda, k, t) \in [-\log N, \log N]  \times \llbracket -T(\log N)^4, T(\log N)^4 \rrbracket \times [0, T \log N]$, establishing the lemma. 
\end{proof}

\begin{proof}[Proof of Lemma \ref{lem:xi_decomp_lem}]

The proofs of \eqref{eqn:xi12_dif} and \eqref{eqn:xim12_dif} are entirely analogous to those of Lemma \ref{lem:Xi_brownian_cont_allparam_pl} (and Lemma \ref{lem:gen_brownian_cont_allparam_prelim}), replacing the use of Lemma \ref{lem:Xi_brownian_cont} (with \eqref{eqn:GNlambda_weak_simul}) in the latter by Lemma \ref{lem:HBM_lem} (with \eqref{eqn:Fbrownian_cont_bd}, Lemma \ref{lem:lc}, and Lemma \ref{lem:Xi_brownian_cont_allparam_pl}) here. We omit further details.
\end{proof}

\subsection{Proof of Lemma \ref{lem:xi_xim_cont}}
\label{subsec:xi2xim2_contproofs}

To show Lemma \ref{lem:xi_xim_cont}, we begin with the following continuity bound.

\begin{lem}\label{lem:initial_cont_lem}

    There exists a constant $\mathfrak{C}>1$ so that the following holds. Let $A>0$ and $S \in [\mathfrak{M}, T]$ be real numbers and $F: \mathbb{R} \rightarrow \mathbb{R}$ be a function satisfying Assumption \ref{ass:F}. For any real numbers $q,q' \in \mathbb{R}$ and $t \geq 0$, define  
    \begin{multline}
        \Xi_1^{(F)}(q,q', t) \coloneqq T^{-1/2} \cdot \Bigg(\sum_{i=N_1}^{N_2} F(\Lambda_i) \cdot 
        \big( \chi\left(q - \alpha i   \right) - \chi\left(q'- \alpha i  - t \ve(\Lambda_i) \right) \big)  \\
        -\mathbb{E} \bigg[\sum_{i=N_1}^{N_2} F(\Lambda_i) \cdot
        \big( \chi\left(q- \alpha i  \right) - \chi\left(q'- \alpha i - t \ve(\Lambda_i) \right) \big) \bigg] \Bigg).
    \end{multline}
 
    \noindent  Then, the following estimate holds overwhelming probability. For any real numbers $q_1,q_2,q_1',q_2' \in [-10 \alpha T^{3/2}, 10 \alpha T^{3/2}]$ with $|q_1-q_2| \leq S$ and $|q_1'-q_2'| \leq S$, and real numbers $t_1,t_2 \in [0, T \log N]$ with $|t_1-t_2| \leq S$, we have 
    \begin{equation}\label{eqn:xi2F_dif_bd}
        \big| \Xi_1^{(F)}(q_1,q_1', t_1) - \Xi_1^{(F)}(q_2,q_2', t_2) \big| \leq A S^{1/2} T^{-1/2} (\log N)^{\mathfrak{C}}.
    \end{equation}
\end{lem}

\begin{proof}

 We only confirm that \eqref{eqn:xi2F_dif_bd} holds with overwhelming probability for a fixed choice of parameters $(q_1, q_1', t_1; q_2, q_2', t_2)$. Its overwhelmingly probable extension to all parameters simultaneously then follows from a union bound over a lattice mesh, with a continuity estimate. We omit further details on the latter, as it is entirely analogous to the proofs of Lemmas \ref{lem:Xi_brownian_cont_allparam_pl} and \ref{lem:gen_brownian_cont_allparam_prelim}. 

So, fix $(q_1, q_1', t_1; q_2, q_2', t_2)$ as in the conditions of the lemma. Let
    \begin{equation*}
        G(q) \coloneqq \chi(q) - \chi(q + q_2-q_1); \qquad \tilde{G}(q) \coloneqq \chi(q) - \chi(q + q_2'-q_1').
    \end{equation*}
    Note that for any $t \in [0, T \log N]$, we may write
    \begin{flalign}\label{eqn:xiF_bdnot}
    \begin{aligned}
        T^{1/2} \cdot \big( & \Xi_1^{(F)}(q_1,q_1', t) -  \Xi_1^{(F)}(q_2,q_2', t) \big) \\
        & = \sum_{i=N_1}^{N_2}F(\Lambda_i) \cdot G(q_1 -\alpha i) - \mathbb{E}\left[ \sum_{i=N_1}^{N_2}F(\Lambda_i) \cdot G(q_1 -\alpha i)\right]
        \\
        & \qquad -  \sum_{i=N_1}^{N_2}F(\Lambda_i) \cdot \tilde{G}(q_1'-\alpha i - t\ve(\Lambda_i)) - \mathbb{E} \left[\sum_{i=N_1}^{N_2}F(\Lambda_i) \cdot \tilde{G} \big(q_1' -\alpha i - t\ve(\Lambda_i) \big)\right].
        \end{aligned}
    \end{flalign}
    
    By \eqref{estimatesm0}, both $G$ and $\tilde{G}$ satisfy Assumption \ref{ass:G} with $(U,S)$ there equal to $(\mathfrak{M}, 2S)$ here. Thus, the first part of Lemma \ref{lem:no_q_conc} yields a concentration estimate for both lines on the right hand side of \eqref{eqn:xiF_bdnot} above. We apply this lemma for $t = 0$ (for the first line), and for both $t = t_1$ and $t= t_2$ (for the second line). For each $t \in \{t_1, t_2\}$, this yields 
\begin{equation}\label{eqn:xi2fixedtbd}
         \big| \Xi_1^{(F)}(q_1,q_1', t) -  \Xi_1^{(F)}(q_2,q_2', t) \big| \leq 
    A S^{1/2} T^{-1/2} (\log N)^C.
\end{equation}

Moreover, defining 
\begin{equation*}
 G_i(x)  = \chi(x) - \chi \big(x    - (t_2-t_1) \cdot \ve(\Lambda_i) \big),
\end{equation*}
we may apply the second part of Lemma \ref{lem:no_q_conc} (at $\Lambda = \Lambda_0$ chosen to satisfy $\ve(\Lambda_0) = 0$, which exists by Lemma \ref{lem:veff_inc} and satisfies $|\Lambda_0| \leq \log N$ for $N$ sufficiently large), with $G_i$ above, to obtain the bound 
\begin{equation}\label{eqn:xi2tdifbd}
   \big| \Xi_1^{(F)}(q,q', t_1) - \Xi_1^{(F)}(q,q', t_2) \big| \leq A S^{1/2} T^{-1/2} (\log N)^C,
   \end{equation}

\noindent for any $(q,q') \in \{ (q_1,q_1'), (q_2, q_2') \}$. The lemma then follows from \eqref{eqn:xi2fixedtbd} and \eqref{eqn:xi2tdifbd}.
\end{proof}

\begin{proof}[Proof of Lemma \ref{lem:xi_xim_cont}]
 
Since $m \leq (\log N)^{1/10}$, the function $\lambda \mapsto \varsigma_m(\lambda) =\lambda^m$ satisfies Assumption \ref{ass:F} with $A = (\log N)^m$. Thus, Lemma \ref{lem:initial_cont_lem} implies \eqref{eqn:xim_2_cont}. It therefore remains to verify \eqref{eqn:xi_2_bound1}. We only show that this bound holds with overwhelming probability for a fixed choice of parameters $(\Lambda_1, k_1, t_1; \Lambda_2, k_2, t_2)$, as then its overwhelmingly probable extension to all parameters simultaneously is proven very similarly to the proof of Lemma \ref{lem:Xi_brownian_cont_allparam_pl}; we omit further details. 

To that end, observe for $\Lambda \in [- \log N, \log N]$ that function $\lambda \mapsto \mathfrak{l}(\Lambda-\lambda)$ satisfies Assumption \ref{ass:F} with $A = (\log N)^3$. Thus, by Lemma \ref{lem:initial_cont_lem}, we may set $(k_1, t_1) = (k_2, t_2)$. More specifically, it suffices to show, for any integer $k \in \llbracket - T^{3/2},  T^{3/2} \rrbracket $; real number $t \in [0,  T \log N]$; and real numbers $\Lambda_1, \Lambda_2 \in [- \log N, \log N]$ with $|\Lambda_1-\Lambda_2| \leq ST^{-1}$, that with overwhelming probability we have 
   \begin{equation}\label{eqn:lambda_fixed_S}
   |\Xi_1(\Lambda_1,k,t) - \Xi_1(\Lambda_2,k,t)| \leq A S^{1/2}T^{-1/2} (\log N)^C.
   \end{equation}

 To do so, fix such parameters $(\Lambda_1,\Lambda_2,k,t)$. Then, defining $\Delta^{(I)}$ and $\Delta^{(II)}$ by 
\begin{flalign}\label{eqn:decomp1}
\begin{aligned} 
\Delta^{(I)} & = \Bigg|  \sum_{i = N_1}^{N_2} (\mathfrak{l}(\Lambda_1 - \Lambda_i) -  \mathfrak{l}(\Lambda_2 - \Lambda_i)) \\
& \qquad \qquad \times \big( 
    \chi(\alpha (k - i)) - 
    \chi\left(\alpha (k - i)+ t  \left(\ve(\Lambda_1) - \ve(\Lambda_i) \right) \right) \big)  - \mathbb{E}[\cdots]   \Bigg|; \\
    \Delta^{(II)} & = \Bigg| \sum_{i = N_1}^{N_2} \mathfrak{l}(\Lambda_2 - \Lambda_i) \cdot 
  \Big( \chi \big(\alpha (k - i)+ t  \left(\ve(\Lambda_1) - \ve(\Lambda_i) \right) \big) \\ 
  & \qquad \qquad \qquad \qquad \qquad -\chi \big(\alpha (k - i)+ t  \left(\ve(\Lambda_2) - \ve(\Lambda_i) \right)  \big) \Big)
    - \mathbb{E}[\cdots]  \Bigg|,
    \end{aligned}
\end{flalign} 

\noindent we must bound $\Delta^{(I)}$ and $\Delta^{(II)}$. For the latter, observe by Lemma \ref{lem:ve_tail} (since $t \le T \log N$ and $|\Lambda_1-|\Lambda_2| \le ST^{-1}$) that $t \cdot |\ve(\Lambda_1) - \ve(\Lambda_2)| \leq S (\log N)^3$. Therefore, Lemma \ref{lem:initial_cont_lem} applies at $q_1=q_2 = 0$; $q_1' = \alpha k + t \ve(\Lambda_1)$; $q_2' = \alpha k + t \ve(\Lambda_2)$; $F(\lambda) = \mathfrak{l}(\Lambda_2-\lambda)$. This yields with overwhelming probability that
\begin{equation}\label{eqn:delta2bd}
    \Delta^{(II)} \leq (\log N)^C S^{1/2}.
\end{equation}

It remains to bound $\Delta^{(I)}$. To that end, set $U \coloneqq 10 S^{1/2} T^{-1/2}$, and let $b : \mathbb{R} \rightarrow \mathbb{R}$ denote a smooth bump function such that $0 \le b(x) \le 1$ for all $x \in \mathbb{R}$; such that $b(x) = 1$ for all $x \in [-U/2, U/2]$; and such that
\begin{flalign*}
    \supp b \subseteq [-U,U]; \qquad \displaystyle\sup_{x \in \mathbb{R}} |b'(x)| \le 20U^{-1}.
\end{flalign*}

Now observe that $\mathfrak{l}(\lambda) = F_{<}(\lambda) + F_{>}(\lambda)$, where
\begin{align}
    \label{2f} 
F_{<}(\lambda)   &= b(\Lambda_1 - \lambda) \cdot \mathfrak{l}(\lambda); \qquad F_{>}(\lambda)   =  (1-b(\Lambda_1 - \lambda)) \cdot \mathfrak{l}(\lambda).
\end{align}

\noindent In this way, we have 
\begin{flalign}
    \label{delta3}
    \Delta^{(I)} \le \Delta^{(I, <)} + \Delta^{(I, >)},
\end{flalign}

\noindent where 
\begin{flalign*}
\Delta^{(I, <)} & = \Bigg|  \sum_{i = N_1}^{N_2} \big( F_{<}(\Lambda_1-\Lambda_i) -  F_{<}(\Lambda_2-\Lambda_i) \big) \\
& \qquad \qquad \quad \times 
\big(  \chi (\alpha k - \alpha i) - 
    \chi(\alpha k - \alpha i+ t \ve(\Lambda_1) - t\ve(\Lambda_i) ) \big)  - \mathbb{E}[\cdots]   \Bigg|;  \\
    \Delta^{(I, >)} & = \Bigg|  \sum_{i = N_1}^{N_2} \big(F_{>}(\Lambda_1-\Lambda_i) -  F_{>}(\Lambda_2-\Lambda_i) \big) \\
    & \qquad \qquad \quad \times \big( 
    \chi(\alpha k - \alpha i) - 
    \chi(\alpha k - \alpha i + t \ve(\Lambda_1) - t \ve(\Lambda_i)) \big)  - \mathbb{E}[\cdots]   \Bigg|.
\end{flalign*}

We must then bound $\Delta^{(I,<)}$ and $\Delta^{(I,>)}$. For the latter, we claim (using Lemma \ref{lem:no_q_conc}) that 
\begin{equation}\label{eqn:Delta>bd}
    \Delta^{(I, >)} \leq S^{1/2} (\log N)^C.
\end{equation}

\noindent To verify this, first observe by \eqref{eqn:sl} that $|\mathfrak{l}(x)| \le (\log N)^C$ for $|x| \le \log N$ and that $|\mathfrak{l}'(x)| \leq  |x|^{-1}$ for all $x \in \mathbb{R}$. Together with the facts that $|\Lambda_1-\Lambda_2| \le ST^{-1}$, that $|b'| \le C(T/S)^{1/2}$, and that $b(x) = 1$ for $|x| \le 5(S/T)^{1/2}$, these imply by \eqref{2f} that for $|\lambda| \le  \log N$ we have $|F_{>}(\Lambda_1-\lambda) -  F_{>}(\Lambda_2-\lambda)| \leq  (S/T)^{1/2} (\log N)^C$. Applying Item \ref{item:conc2} of Lemma \ref{lem:no_q_conc}, with $(A, t, \Lambda, q, \Delta, S)$ there given by $((S/T)^{1/2} (\log N)^C, 0, \Lambda_1, \alpha k, t, T \log N) $ (and $F(\lambda) = F_{>} (\Lambda_1-\lambda) - F_{<} (\Lambda_2-\lambda)$) here, yields 
 \begin{equation}
     \label{delta00}
 \Delta^{(I, >)} \leq   (\log N)^C S^{1/2} .
 \end{equation}
Here, since $|\Lambda_1| \leq \log N$, we have absorbed the factor $(|\Lambda_1|+1)^2$ into the quantity $(\log N)^C$. 

To bound $\Delta^{(I,<)}$, observe that  $\Delta^{(I, <)} = |\sum_{i=N_1}^{N_2} H(\lambda,i)|$, where 
\begin{equation}\label{eqn:H_lamdif}
    H(\lambda, i) = \big( F_{<}(\Lambda_1-\lambda) -  F_{<}(\Lambda_2-\lambda) \big) \cdot
\big( 
    \chi(  \alpha k- \alpha i ) - 
    \chi (\alpha k- \alpha i+ t \ve(\Lambda_1) - t \ve(\lambda)) \big)
    \big).
\end{equation}

\noindent Under this notation, it is quickly verified that $H$ satisfies Assumption \ref{ass:Hass} (with the $A$ there equal to $1$ here). Therefore, Lemmas \ref{lem:resolvent_lemma} and \ref{lem:Psi2_final_approx} apply and yield, with $R$ as in Assumption \ref{ass:R_M_ass},
\begin{flalign}\label{eqn:F<_ind_decomp}
 \Delta^{(I,<)}    = \Bigg| \sum_{r = \lceil N_1/R \rceil}^{\lfloor N_2/R-1 \rfloor} \sum_{s=1}^{R} \phi_r (\lambda_s^{(r)}) - \mathbb{E}[\cdots] \Bigg| + O(T^{1/2-\mathfrak{c}}),
\end{flalign}

\noindent where we recall from Definition \ref{def:Lr} that the $(\lambda_i^{[r]})$ are the eigenvalues of the $R \times R$ submatrices $\mathbf{L}^{[r]}$ of $\mathbf{L}$ with rows and columns indexed by $\llbracket Rr, (R+1)r-1\rrbracket$. Moreover, in \eqref{eqn:F<_ind_decomp}, for each $r \in \mathbb{Z}$, we have denoted the function $\phi_r : \mathbb{R} \rightarrow \mathbb{R}$ by setting $\phi_r (\lambda) = H(\lambda, R r)$, for each $\lambda \in \mathbb{R}$. 

Then, Proposition \ref{prop:limvar} yields for any $r \in \llbracket N_1/R, N_2/R-1\rrbracket$ the approximation
\begin{equation}\label{eqn:variance_r}
  R^{-1} \cdot \var Y_r = \sigma^2(\phi_r) + O(T^{-c}), \qquad \text{where} \qquad Y_r = \displaystyle\sum_{s=1}^R \phi_r(\lambda_s^{(r)}),
\end{equation}
 
 \noindent and $\sigma^2(\phi_r)$ is given by  Definition \ref{def:Cdef}. Observe that 
\begin{equation}\label{eqn:varphirbd}
    \sigma^2(\phi_r) \leq C \| \phi_r \|_{\mathcal{H}}^2 \leq (\log N)^3 
    \int_{-\infty}^{\infty} b(\Lambda_1 - \lambda)^2  \cdot \varrho(\lambda) d \lambda \leq  10 S^{1/2} T^{-1/2} (\log N)^3,
\end{equation}

\noindent where the first bound follows from \eqref{eqn:Cdef} and the fact that $(1 - \theta \mathbf{T} \bm{\varrho}_{\beta})^{-1}$ is bounded, and the second and third follow from \eqref{def:inner_prod} and the fact that $\supp b \subseteq [-5(S/T)^{1/2},5(S/T)^{1/2}] \subseteq [-5,5]$ (with the bound $|\mathfrak{l}(\lambda)| \le 10(\log N)^2$ for $|\lambda| \le \log N$).

 Now observe that, with overwhelming probability, there at most $TR^{-1} (\log N)^C$ indices $r \in \mathbb{Z}$ for which $Y_r \ne 0$, as $\supp \chi \subseteq [-\mathfrak{M},\mathfrak{M}]$ and $|t\ve (\lambda_s^{(r)})| \le T (\log N)^3$ holds for all $s$ and $r$ with overwhelming probability (by Lemmas \ref{lem:ve_tail} and \ref{lem:bd_lem}). Moreover, the $(Y_r)$ are mutually independent and satisfy $|Y_r| \le R (\log N)^C$ with overwhelming probability (by Lemma \ref{lem:bd_lem}). These with  \eqref{eqn:variance_r}, \eqref{eqn:varphirbd}, and the Bernstein estimate yield, with overwhelming probability,
 \begin{flalign*}
 \Delta^{(I,<)}    = \Bigg| \sum_{r = \lceil N_1/R \rceil}^{\lfloor N_2/R-1 \rfloor} Y_r - \mathbb{E}[\cdots] \Bigg| + O(T^{1/2-c}) \le (S^{1/4} T^{1/4} + T^{1/2-c}) (\log N)^C.
 \end{flalign*} 

\noindent Combining this with \eqref{delta00}, \eqref{delta3}, \eqref{eqn:delta2bd}, and \eqref{eqn:decomp1} establishes \eqref{eqn:lambda_fixed_S} and thus the lemma.
\end{proof}

\section{Coupling with a Gaussian process}
\label{sec:levy_chent}

In this section we prove Theorem \ref{thm:couplethm12}. Throughout this section, we adopt Assumptions \ref{ass:NT_assumption} and we recall the notation from \eqref{functionslambdar} and Definitions \ref{def:white_noise_dress}, \ref{def:Wdress}, \ref{def:Z_intro}, \ref{def:chi_def}, \ref{xilambdakt}, \ref{xim}, \ref{x2x0y2y0}, and \ref{def:Cdef}. Moreover, we assume throughout (even when not explicitly stated) that $\theta < \theta_0(\beta)$ is sufficiently small (so that all of the results from Sections \ref{ProofQ}, \ref{sec:bg}, \ref{sec:ind_linear}, and \ref{sec:continuity} apply). 

\subsection{Comparing prelimiting and limiting variances}
\label{subsec:couplepre}

To motivate the following preliminary lemmas, recall the Gaussian process $\mathcal{W}^{\dr}$ (Definition \ref{def:Wdress}), whose covariance structure of $\mathcal{W}^{\dr}$ is related to the bilinear form $\mathscr{C}$ (Definition \ref{def:Cdef}). Specifically, recall that $\mathcal{W}^{\dr}$ is a Gaussian process indexed by functions $f \in L^2(dr \otimes \varrho d \lambda)$ (from \eqref{functionslambdar}). In what follows, for any $f \in L^2(dr \otimes \varrho d \lambda)$ and $r \in \mathbb{R}$, define $f_r \in \mathcal{H}$ by setting $f_r (\lambda) \coloneqq f(r, \lambda)$ for each $\lambda \in \mathbb{R}$; we also set $\mu_f : \mathbb{R} \rightarrow \mathbb{R}$ by setting $\mu_f (r) = \mu_{f_r}$ (Definition \ref{def:Cdef}) for each $r \in \mathbb{R}$. The Gaussian process $\mathcal{W}^{\dr}$ has the property that for any two functions $f, g \in L^2(dr \otimes \varrho d \lambda)$, the pair $  (\mathcal{W}^{\dr}(f- \mu_{f} \varsigma_0), \mathcal{W}^{\dr}(g- \mu_{g} \varsigma_0))$ is (by \eqref{eqn:Cdef}) jointly Gaussian with covariance
\begin{equation}\label{eqn:wdr_cov}
    \Cov_{\mathcal{W}^{\dr}}(f- \mu_{f} \varsigma_0, g - \mu_{g} \varsigma_0) = \int_{-\infty}^{\infty} \mathscr{C}(f_r-\mu_{f_r} \varsigma_0, g_r-\mu_{g_r} \varsigma_0) dr .  
\end{equation}

Before stating the main lemma needed for the proof of Theorem \ref{thm:couplethm12}, we define (and recall) some more notation. In particular, next we will introduce three different families of prelimiting and corresponding limiting functions.

\begin{definition} 

\label{ff} 

 We define three families of \emph{prelimiting functions}, each of which will be associated with a \emph{corresponding limiting function}. 

  \textbf{Family 1}: Let $\Lambda \in [-\log N, \log N]$; $k \in [- T(\log N)^4, T (\log N)^4]$; and $t \in [0, T \log N]$ be real numbers. For any $r \in \mathbb{R}$, define the function $\psi_{\Lambda,k,t;r}^{(N)}: \mathbb{R} \rightarrow \mathbb{R}$ by 
\begin{equation}\label{eqn:phi_lc_reg}
\psi_{\Lambda, k,t; r}^{(N)}(\lambda) =  2 \mathfrak{l}(\Lambda - \lambda)
\cdot \big(
    \chi(\alpha k - \alpha r R) 
    - 
    \chi\left( \alpha k - \alpha r R + t \ve(\Lambda) - t \ve(\lambda) \big)
   \right).
\end{equation}

\noindent The \emph{prelimiting function} $f_{r}^{(N)}: \mathbb{R} \rightarrow \mathbb{R}$ is given by 
\begin{equation}\label{eqn:psi_regf}
   f_{r}^{(N)}(\lambda) \coloneqq  \psi_{\Lambda, k,t; r}^{(N)}(\lambda) .
   \end{equation}

 \noindent For any $r \in \mathbb{R}$, the \emph{corresponding limiting function} $f_r: \mathbb{R} \rightarrow \mathbb{R}$ is given by (recalling $\psi_{\Lambda, \kappa,\tau} : \mathbb{R}^2 \rightarrow \mathbb{R}$ from \eqref{eqn:logphi_intro})
\begin{equation}\label{eqn:psi_r_2f}
    f_{r}(\lambda) \coloneqq  \psi_{\Lambda, k/T,t/T}( r, \lambda) .
   \end{equation}

 \textbf{Family 2}: Let $m \in \llbracket 0, (\log N)^{1/10} \rrbracket$ be an integer and $q,q' \in [- T(\log N)^4, T (\log N)^4]$ and $t \in [0, T \log N]$ be real numbers. For any $r \in \mathbb{R}$, define the function $\phi_{q,q',t; r}^{[m, N]}(\lambda): \mathbb{R} \rightarrow \mathbb{R}$ by 
\begin{equation}\label{eqn:phi_m_reg}
\phi_{q,q',t; r}^{[m, N]}(\lambda) \coloneqq \lambda^m \cdot \big(\chi(q-\alpha r R)  - \chi(q'-\alpha r R - t \ve(\lambda)) \big).
\end{equation}

\noindent The \emph{prelimiting function} $f_{r}^{(N)}: \mathbb{R} \rightarrow \mathbb{R}$ is given by 
\begin{equation}\label{eqn:psi_m_regf}
    f_{r}^{(N)}(\lambda) \coloneqq  \phi_{q,q',t; r}^{[m, N]}(\lambda) .
   \end{equation}

 \noindent For any $r \in \mathbb{R}$, the \emph{corresponding limiting function} $f_r: \mathbb{R} \rightarrow \mathbb{R}$ is given by (recalling $\phi_{\mathfrak{q}, \mathfrak{q}',\tau}^{[m]}:\mathbb{R}^2 \rightarrow \mathbb{R}$ from \eqref{eqn:phim_intro}) 
  \begin{equation}\label{eqn:phi_m_noreg_2f}
    f_{r}(\lambda) \coloneqq  \phi_{q/T, q'/T,t/T}^{[m]}(r, \lambda).
   \end{equation}

 \textbf{Family 3}: Let $A,B \in [1, (\log N)^3]$; $U \in [\mathfrak{M},T \log N]$; and $t \in [0, T \log N]$ be real numbers, and let $F : \mathbb{R} \rightarrow \mathbb{R}$ and $G : \mathbb{R} \rightarrow \mathbb{R}$ be functions. Suppose that $F$ satisfies Assumption \ref{ass:F}, and that $G$ has compact support and satisfies Assumption \ref{ass:G2}. For any $r \in \mathbb{R}$, the \emph{prelimiting function} $f_r^{(N)}: \mathbb{R} \rightarrow \mathbb{R}$ is given by 
 \begin{equation}\label{eqn:fFG_prelim}
     f_r^{(N)}(\lambda) \coloneqq F(\lambda) \cdot G(\alpha r R + t \ve(\lambda)).
 \end{equation} 

 \noindent For any $r \in \mathbb{R}$, the \emph{corresponding limiting function} $f_r : \mathbb{R} \rightarrow \mathbb{R}$ is given by 
\begin{equation}\label{eqn:fFG_lim}
    f_{r}(\lambda) \coloneqq F(\lambda) \cdot  G( T \alpha r  + t \ve(\lambda) ).
   \end{equation}
\end{definition}

With this notation, we have the following lemma, to be shown in Section \ref{subsec:couplepreproofs}. It provides the covariances between the terms in the sums approximating $(\Xi-\Xi_1, \Xi_1, \Xi^{[m]} - \Xi_1^{[m]}, \Xi_1^{[m]})$ from Proposition \ref{prop:equiv_expr_cor}. Below, we recall the $R \times R$ matrix $\mathbf{L}^{[r]}$ and its eigenvalues $(\lambda_i^{(r)})$ from Definition \ref{def:Lr}. We also recall the stretch random variables $(X_i)$ from Definition \ref{xi}.

\begin{lem}[Covariances between summands in Proposition \ref{prop:equiv_expr_cor}]\label{lem:cov_error_estimate_specific}

There exists a constant $\mathfrak{c}>0$ such that the following holds. Let $f_r^{(N)}, g_r^{(N)} : \mathbb{R} \rightarrow \mathbb{R}$ denote any of the prelimiting functions from Definition \ref{ff}, and let $f_{r}, g_{r} : \mathbb{R} \rightarrow \mathbb{R}$ denote the corresponding limiting functions, respectively. Then (recalling Definition \ref{def:Lr}) we have 
\begin{flalign}
  & \left| T^{-1} \displaystyle\sum_{r = \lceil N_1/R \rceil}^{\lfloor N_2/R \rfloor - 1}  \Cov\left( \sum_{i=1}^R f_r^{(N)}(\lambda_i^{(r)}),  \sum_{i=1}^R g_r^{(N)}(\lambda_i^{(r)}) \right) - \int_{-\infty}^{\infty} \mathscr{C}(f_{r}-\mu_{f_{r}} \varsigma_0, g_{r}-\mu_{g_{r}} \varsigma_0) d r \right| \leq T^{-\mathfrak{c}}; \label{eqn:covariance_bd_touse} \\
 &  \Bigg| T^{-1} \displaystyle\sum_{r = \lceil N_1/R \rceil}^{\lfloor N_2/R \rfloor - 1}  \Cov\left( \sum_{i=1}^R f_r^{(N)}(\lambda_i^{(r)}), \langle g_{r}^{(N)}, \varsigma_0 \rangle_{\varrho}  \sum_{i=1}^R X_i^{(r)} \right) - \int_{-\infty}^{\infty} \mu_{g_{r}} \mathscr{C}(f_{r}-\mu_{f_{r}} \varsigma_0, \varsigma_0) d r \Bigg| \leq T^{-\mathfrak{c}};\label{eqn:cov_f_br_bd_to_use} \\
 &  \Bigg| T^{-1} \displaystyle\sum_{r = \lceil N_1/R \rceil}^{\lfloor N_2/R \rfloor - 1}  \Cov\left(  \langle f_{r}^{(N)}, \varsigma_0 \rangle_{\varrho}  \sum_{i=1}^R X_i^{(r)}, \langle g_{r}^{(N)}, \varsigma_0 \rangle_{\varrho}  \sum_{i=1}^R X_i^{(r)} \right)  - \int_{-\infty}^{\infty} \mu_{f_{r}} \mu_{g_{r}}  \mathscr{C}(\varsigma_0, \varsigma_0) d r \Bigg| \leq T^{-\mathfrak{c}}. \label{eqn:cov_br_br_bd_to_use}
\end{flalign}

\end{lem}

\subsection{Gaussian coupling lemmas}

\label{subsec:gauss_couple_prelim}

In this section we state two lemmas that will be used in the proof of Theorem \ref{thm:couplethm12}. The first provides, in the special case we need, a multi-dimensional version of the Kolm\'{o}s--Major--Tusn\'{a}dy (KMT) coupling. The second states that multivariate Gaussian random variables, with approximately equal covariance matrices, can be coupled to likely be close to each other. 

Below, for any integer $d \ge 0$ and $d$-dimensional random variable $v = (v_1, v_2, \ldots , v_d)$, the  \emph{covariance matrix} of $v$ is the $d \times d$ matrix, whose $(i,j)$-entry is equal to $\Cov (v_i, v_j)$ for each $i, j \in \llbracket 1, d \rrbracket$.

\begin{lem}[{\cite[Equation (1.1)]{BM06}}]\label{lem:final_couple}

There exist constants $\mathfrak{c}>0$ and $\mathfrak{C}>1$ such that the following holds. Let $n, d \ge 1$ be integers; $B>0$ be a real number; and $Y_1, \ldots, Y_n \in \mathbb{R}^d$ be mutually independent, centered ($d$-dimensional) random variables, such that $\max_{1 \le i \le n} \|Y_i\|_2 \le B$. Then, for any real number $\delta>0$, there exists a coupling between $(Y_1,\dots, Y_n)$ and a family $(Z_1, \ldots, Z_n)$ of mutually independent, centered ($d$-dimensional) Gaussian random variables, such that $Y_i$ and $Z_i$ have the same covariance matrix for each $i \in \llbracket 1, n \rrbracket$, and
\begin{equation}
\label{eqn:coupling}
\mathbb{P}\Bigg[
\bigg| \sum_{i=1}^n (Y_i - Z_i) \bigg| > \delta
\Bigg]
\le
\mathfrak{C} d^2 \exp\!\left( -\frac{\mathfrak{c} \delta}{Bd^2} \right).
\end{equation}
\end{lem}

The next lemma (which is not optimal but sufficient for our purposes) is proved in Appendix \ref{app:gauss_couple}.

\begin{lem}\label{lem:gauss_couple}
   Let $d \ge 1$ be an integer; $\gamma \ge 0$ be a real number; and $W = (W_1, W_2, \ldots , W_d) \in \mathbb{R}^d$ and $\tilde{W} = (\tilde{W}_1, \tilde{W}_2, \ldots , \tilde{W}_d) \in \mathbb{R}^d$ be two centered ($d$-dimensional) Gaussian random variables vectors, with covariance matrices $\Sigma$ and $\tilde \Sigma$, respectively. Assume they satisfy 
   \begin{equation}\label{eqn:cov_bd}
 \sup_{1\leq i, j \leq d}  |\Sigma_{i j} - \tilde \Sigma_{ij}| \leq \gamma.
   \end{equation}
  
   \noindent Then, there exists a coupling between $W$ and $\tilde W$ such that 
\begin{equation}
\label{eq:coupling}
\mathbb{P} \bigg[
\max_{i \in \llbracket 1, d \rrbracket} | W_i - \tilde W_i |  > \gamma^{1/2} d \log d
\bigg] 
\le
(2\pi)^{-1/2} e^{-(\log d)^2 / 2}.
\end{equation}
\end{lem}

\subsection{Proof of Theorem \ref{thm:couplethm12}}
\label{subsec:main_coupling_arg}

In this section we prove Theorem \ref{thm:couplethm12}. First, we define the function $B: \mathbb{R} \rightarrow \mathbb{R}$ by setting (where below we view $\mathbbm{1}_{x \le r \le 0} \cdot \varsigma_0 (\lambda) - \mathbbm{1}_{0 \le r \le x} \cdot \varsigma_0 (\lambda)$ as an $x$-dependent function in $L^2 (dr \otimes \varrho d \lambda)$)
\begin{equation}\label{eqn:Brownian}
B(x) \coloneqq \alpha \cdot \mathcal{W}^{\dr} \big( \mathbbm{1}_{x \le r \le 0} \cdot \varsigma_0 (\lambda) - \mathbbm{1}_{0 \le r \le x} \cdot \varsigma_0 (\lambda)  \big),
\end{equation} 

\noindent for any $x \in \mathbb{R}$. By 
Definition \ref{def:Wdress}, we have $\Var (B(x) - B(y)) =  (x-y) \cdot \alpha^2 \| \varsigma_0^{\dr} \|_{\mathcal{H}}^2$ (and thus only depends on $x-y$) for any $x \ge y$. Moreover, for any $x_1 \le x_2 \le \cdots \le x_n$, the $(B(x_{i+1})-B(x_i))$ are mutually independent Gaussian random variables, by \eqref{eqn:wdr_cov} (as the supports in $r$ of the functions $\mathbbm{1}_{x_i \le r \le x_{i+1}})$ are disjoint over $i$). Hence, $B(x)$ is a two-sided Brownian motion of variance $\alpha^2 \| \varsigma_0^{\dr} \|_{\mathcal{H}}^2$.

It is moreover quickly verified that, if we define the coefficients  
\begin{flalign}
\mathfrak{c}(\Lambda, \kappa, \tau; r)  \coloneqq  \int_{-\infty}^{\infty}\psi_{\Lambda, \kappa,\tau}( r, \lambda) \varrho(\lambda) d\lambda 
= \langle \psi_{\Lambda, \kappa,\tau}( r,\cdot) , \varsigma_0 \rangle_{\varrho}; \label{eqn:mathfrakc1} \\
\mathfrak{c}^{[m]}(\mathfrak{q}, \mathfrak{q}',\tau; r)  \coloneqq  \int_{-\infty}^{\infty} \phi_{\mathfrak{q}, \mathfrak{q}',\tau}^{[m]}(r, \lambda) \varrho(\lambda) d \lambda = \langle \phi_{\mathfrak{q}, \mathfrak{q}',\tau}^{[m]}(r, \cdot), \varsigma_0\rangle_{\varrho},\label{eqn:mathfrakc2}
\end{flalign}
then (recalling the processes $\mathfrak{X}_1, \mathfrak{X}_0, \mathfrak{Y}_1, \mathfrak{Y}_0$ from Definition \ref{x2x0y2y0}) we have the representations 
\begin{align} 
			\mathfrak{X}_0(\Lambda, \kappa, \tau) &=-\alpha^{-1} \int_{-\infty}^{\infty} \mathfrak{c}(\Lambda, \kappa, \tau; r) dB(r);
            \label{eqn:X0eq} 
			      \\
           \mathfrak{Y}_0(\Lambda, \kappa, \tau) &=-\alpha^{-1} \int_{-\infty}^{\infty} \mathfrak{c}^{[m]}(\mathfrak{q}, \mathfrak{q}',\tau; r) dB(r)  . \label{eqn:Y0eq} 
\end{align} 

\noindent Observe that \eqref{eqn:X0eq}  and \eqref{eqn:Y0eq} mirror the definitions \eqref{eqn:GN_def} of $\mathcal{G}_N$ and \eqref{eqn:GNF_gen} of $\mathcal{G}_N^{(\varsigma_n)}$. Indeed, we will see that the Brownian motion $B(x)$ will represent the fluctuation scaling limit for the stretch random variables $X_i$ from Definition \ref{xi} (in that functionals of the form $\sum_{i=N_1}^{N_2} f(iT^{-1}) X_i$ will converge in distribution to $\int_{-\infty}^{\infty} f(x) dB(x)$ for certain test functions $f: \mathbb{R} \rightarrow \mathbb{R}$).

\begin{remark} 

\label{mufr}

Observe that, if the prelimiting function $f_r^{(N)}$ from Definition \ref{ff} is $f_r^{(N)} = \psi_{\Lambda,k,t;r}^{(N)}$ as in \eqref{eqn:psi_regf}, then $\langle f_r^{(N)}, \varsigma_0 \rangle_{\rho} = c^{(\Lambda)} (r,k,t)$ from \eqref{eqn:clambda1}. Similarly, if $f_r^{(N)} = \phi_{q,q',t;r}^{[m,N]}$ as in \eqref{eqn:phi_m_noreg_2f}, then $\langle f_r^{(N)}, \varsigma_0 \rangle_{\varrho} = c^{[m]} (r,q,q',t)$ from \eqref{eqn:cm1}. Correspondingly, if the limiting function $f_r$ from Definition \ref{ff} is $f_r = \psi_{\Lambda,\kappa,\tau} (r, \cdot)$ as in \eqref{eqn:psi_r_2f}, then $\mu_{f_r} = \langle f_r, \varsigma_0 \rangle_{\varrho} = \mathfrak{c} (\Lambda, \kappa, \tau; r)$ from \eqref{eqn:mathfrakc1}. Analogously, if it satisfies $f_r = \phi_{q,q',\tau}^{[m]} (r, \cdot)$ as in \eqref{eqn:phi_m_noreg_2f}, then $\mu_{f_r} = \langle f_r, \varsigma_0 \rangle_{\varrho} = \mathfrak{c}^{[m]} (\mathfrak{q}, \mathfrak{q}', \tau; r)$ from \eqref{eqn:mathfrakc2}. 

\end{remark} 

Recalling $\mathfrak{X}$ and $\mathfrak{Y}$ from Definitions \ref{def:Z_intro} and \ref{x2x0y2y0}, Lemma \ref{lem:holder_cont} states that they are almost surely H\"{o}lder continuous in their parameters. By Definition \ref{x2x0y2y0}, we almost surely have 
\begin{equation}
    \mathfrak{X}(\Lambda, \kappa, \tau) = \mathfrak{X}_1(\Lambda, \kappa, \tau) + \mathfrak{X}_0(\Lambda, \kappa, \tau); \qquad 
    \mathfrak{Y}(\mathfrak{q}, \mathfrak{q}', \tau, m) = \mathfrak{Y}_1(\mathfrak{q}, \mathfrak{q}', \tau, m) + \mathfrak{Y}_0(\mathfrak{q}, \mathfrak{q}', \tau, m).
\end{equation}

With this notation, we can now establish of Theorem \ref{thm:couplethm12}.

\begin{proof}[Proof of Theorem \ref{thm:couplethm12}]

 Let us outline how we will proceed. First, we will use Proposition \ref{prop:equiv_expr_cor} and Lemma \ref{lem:xi_decomp_lem} to show that $\Xi_1 (\Lambda,k,t)$ and $\Xi(\Lambda, k, t) - \Xi_1(\Lambda,k,t)$ (or $\Xi_1^{[m]} (q,q',t)$ and $\Xi^{[m]} (q,q', t) - \Xi_1^{[m]} (q,q',t)$) are  with overwhelming probability approximately equal to sums of independent random variables, in their full ranges of parameters $(\Lambda,k,t)$ (or $(q,q',t)$, respectively). Second, using the variance estimates from Lemma \ref{lem:cov_error_estimate_specific}, we apply multivariate Gaussian approximation bounds (Lemmas \ref{lem:final_couple} and \ref{lem:gauss_couple}) to couple $(\Xi_1,\Xi-\Xi_1)$ (or $(\Xi_1^{[m]}, \Xi^{[m]} - \Xi_1^{[m]})$) to approximate the Gaussian random variables $(\mathfrak{X}_1, \mathfrak{X}_0)$ (or $(\mathfrak{Y}_1, \mathfrak{Y}_0)$, respectively), on a certain lattice mesh of parameter values. Third, we extend from this mesh to the full parameter range using continuity bounds for the prelimiting and limiting processes, as shown in Section \ref{sec:continuity} (and Lemma \ref{lem:holder_cont}).

We begin with the first step above. By Lemma \ref{lem:xi_decomp_lem}, we have with overwhelming probability that 
\begin{flalign}
    \label{xigxi}
    |\Xi(\Lambda, k, t)-\Xi_1(\Lambda, k, t) - \mathcal{G}_N(\Lambda, k, t)| + |\Xi^{[m]}(q, q',t) - \Xi_1^{[m]}(q, q',t) -\mathcal{G}_N^{(\varsigma_m)}(q,q',t)| \le T^{-c},
\end{flalign}

\noindent for all integers $m \in \llbracket 0, (\log N)^{1/10}$ and real numbers $|\Lambda| \le \log N$; $T \in (\log N)^4$; $|k| \le T(\log N)^4$; and $|q|, |q'| \in [-T(\log N)^4, T(\log N)^4]$. As such, we will in what follows replace $\Xi(\Lambda, k, t)-\Xi_1(\Lambda, k, t)$ by $\mathcal{G}_N(\Lambda, k, t)$, and $\Xi^{[m]}(q, q',t) - \Xi_1^{[m]}(q, q',t)$  by $ \mathcal{G}_N^{(\varsigma_m)}(q,q',t)$. 

To implement the second, we define a mesh covering the relevant parameter ranges. Let $\delta_0 = \delta_0(N) = T^{-\nu}$ be a small parameter, where $\nu$ is a small constant to be chosen later, and set 
$$M_0 \coloneqq (\log N)^{1/10}.$$
Then define the lattices 
\begin{align}
\mathfrak{L}_{\Lambda} &= \{ i \delta_0 :i \in \mathbb{Z}, i \delta_0  \in [-\log N, \log N] \} \label{eqn:lat_1}; \\
\mathfrak{L}_{\kappa} &= \{ k   \delta_0  : k \in \mathbb{Z}, k  \delta_0  \in [- (\log N)^4,  (\log N)^4]  \}   \label{eqn:lat_2}; \\
\mathfrak{L}_{\mathfrak{q}} &= \{q \delta_0  : q \in \mathbb{Z}, q \delta_0 \in [-  (\log N)^5,   (\log N)^5] \}   \label{eqn:lat_3}; \\
\mathfrak{L}_{\tau} &= \{t \delta_0 : t \in \mathbb{Z}, t \delta_0 \in [0, \log N] \}; \label{eqn:lat_4} \\
\mathfrak{L}_1 &\coloneqq  \mathfrak{L}_{\Lambda}  \times \mathfrak{L}_{\kappa} \times  \mathfrak{L}_{\tau} \label{eqn:lat_5}; \\
\mathfrak{L}_2 &\coloneqq  \mathfrak{L}_{\mathfrak{q}}^2 \times \mathfrak{L}_{\tau}.  \label{eqn:lat_6} 
\end{align} 
We wish to couple the set of random variables
\begin{flalign*}
    & \{\mathcal{G}_N(\Lambda, T \kappa,T\tau)\}_{(\Lambda,  \kappa, \tau) \in \mathfrak{L}_1} \cup \{\Xi_1(\Lambda, T \kappa, T \tau)\}_{(\Lambda,  \kappa, \tau) \in \mathfrak{L}_1} \\
    & \qquad \cup \{ \mathcal{G}_N^{(\varsigma_m)}(\mathfrak{q},\mathfrak{q}', \tau) \}_{(\mathfrak{q},\mathfrak{q}', \tau) \in \mathfrak{L}_2,  0 \leq m \leq M_0} \cup \{ \Xi_1^{[m]}(T \mathfrak{q}, T \mathfrak{q}', T \tau) \}_{(\mathfrak{q},\mathfrak{q}', \tau) \in \mathfrak{L}_2,  0 \leq m \leq M_0},
\end{flalign*}
with their limiting counterparts (which are jointly Gaussian). 

To do so, recall the $(X_i)$ from Definition \ref{xi}; the $(\lambda_i^{(r)})$ from Definition \ref{def:Lr}; and the $c^{(\Lambda)}$ and $c^{[m]}$ from \eqref{eqn:clambda1} and \eqref{eqn:cm1}. For any $r \in \mathbb{Z}$ and $m \in \mathbb{Z}_{\ge 0}$, and $\mathfrak{q},\mathfrak{q}' \in \mathbb{R}$ and $\tau \in \mathbb{R}_{\ge 0}$, define  
\begin{flalign}
\label{y001}  
\begin{aligned} 
Y_{\mathfrak{q}, \tau; r;0}^{[m], (N)} &= -  \alpha^{-1} R^{-1/2} \cdot c^{[m]}(r,q,t) \sum_{i = r R}^{(r+1) R-1} X_{i}; \\
 Y_{\mathfrak{q}, \tau; r;1}^{[m], (N)} &=R^{-1/2} \sum_{s=1}^{R} (\lambda_s^{(r)})^m \cdot \big( \chi (q- \alpha r R ) - \chi (q- \alpha r R- t \ve(\lambda_s^{(r)})) \big)  -\mathbb{E}[\cdots].
 \end{aligned} 
  \end{flalign}

  \noindent Moreover,  for any integer $r \in \mathbb{Z}$ and real numbers $\Lambda, \kappa \in \mathbb{R}$ and $\tau \geq 0$, define 
  \begin{flalign}
  \label{x001} 
  \begin{aligned} 
  X_{\Lambda, \kappa, \tau; r;0}^{(N)} &=    - \alpha^{-1} R^{-1/2} \cdot c^{(\Lambda)}(r,k,t) \sum_{i = r R}^{(r+1) R-1} X_{i}; \\
  X_{\Lambda, \kappa, \tau; r;1}^{(N)} &= R^{-1/2} \sum_{s=1}^{R} \mathfrak{l}(\Lambda -\lambda_s^{(r)})  \cdot \big( \chi (\alpha k- \alpha r R)   - \chi(\alpha k- \alpha r R  + t \ve(\Lambda)- t \ve(\lambda_s^{(r)})) \big)  -\mathbb{E}[\cdots].
  \end{aligned} 
  \end{flalign}

 Proposition \ref{prop:equiv_expr_cor}, \eqref{eqn:GN_def}, \eqref{eqn:GNF_gen}, and a union bound imply, with overwhelming probability, that for all $(\Lambda,\kappa, \tau) \in \mathfrak{L}_1$, $(\mathfrak{q},\mathfrak{q}',\tau) \in \mathfrak{L}_2$, and $m \in \llbracket 1, M_0 \rrbracket$, we have 
\begin{align}
\label{xisum} 
\begin{aligned} 
\Xi_1(\Lambda, T \kappa, T \tau) &= R^{1/2} T^{-1/2} \sum_{r=\lceil N_1/R \rceil}^{\lfloor N_2/R \rfloor - 1} X_{\Lambda, \kappa, \tau; r;1}^{(N)}   + O(T^{-\mathfrak{c}}); \\
\mathcal{G}_N(\Lambda, T \kappa, T \tau) &= R^{1/2} T^{-1/2} \sum_{r=\lceil N_1/R \rceil}^{\lfloor N_2/R \rfloor - 1} X_{\Lambda, \kappa, \tau; r;0}^{(N)};   \\
\Xi_1^{[m]}(T \mathfrak{q}, T \tau) &= R^{1/2} T^{-1/2} \sum_{r=\lceil N_1/R \rceil}^{\lfloor N_2/R \rfloor - 1} Y_{\mathfrak{q}, \tau; r;1}^{[m], (N)}  +  O(T^{-\mathfrak{c}}); \\
\mathcal{G}_N^{(\varsigma_m)}(\mathfrak{q},\mathfrak{q}', \tau)  &= R^{1/2} T^{-1/2} \sum_{r=\lceil N_1/R \rceil}^{\lfloor N_2/R \rfloor - 1} Y_{\mathfrak{q}, \tau; r;0}^{[m], (N)}.
\end{aligned} 
\end{align}

\noindent Moreover, applying the second part of Lemma \ref{lem:no_q_conc} (with the $N$ there equal to $R$ here) to $X_{\Lambda, \kappa, \tau; r;1}^{(N)}$ and $Y_{\mathfrak{q}, \tau; r;1}^{[m], (N)}$ and Lemma \ref{lem:q_spacing} to bound $X_{\Lambda, \kappa, \tau; r;0}^{(N)}$ and $Y_{\mathfrak{q}, \tau; r;0}^{[m], (N)}$, we find with overwhelming probability that, for all values of parameters as above, and for all $r \in \llbracket N_1/R, N_2/R-1 \rrbracket$, we have for each $i \in \{ 0, 1 \}$ that 
\begin{align*}
   |X_{\Lambda, \kappa, \tau; r;i}^{(N)}| &\leq (\log N)^{10}; \qquad  |Y_{\mathfrak{q}, \tau; r;i}^{[m], (N)} | \leq (\log N)^{m+10}. 
\end{align*}

Now, observe (from Definition \ref{def:Lr}) that the $(X_{\Lambda,\kappa,\tau;r;0}^{(N)}, X_{\Lambda,\kappa,\tau;r;1}^{(N)}, Y_{\mathfrak{q},\tau;r;0}^{(N)}, Y_{\mathfrak{q},\tau;r;1}^{(N)})$ are mutually independent over the index $r$. We will ultimately choose $\delta_0 = T^{-\nu}$ in the definition of \eqref{eqn:lat_1}, \eqref{eqn:lat_2}, \eqref{eqn:lat_3}, and \eqref{eqn:lat_4} so that the total number of parameters satisfies, for some small absolute constant $\mathfrak{c}_0 > 0$, 
\begin{equation}\label{eqn:ddef}
    d \coloneqq 2 \left( |\mathfrak{L}_1| + M_0 |\mathfrak{L}_2| \right) \leq T^{\mathfrak{c}_0}.
\end{equation}

\noindent Hence, we may apply Lemma \ref{lem:final_couple} to the sums of independent random variables, multiplied by $T^{1/2} R^{-1/2}$, in  \eqref{xisum}; we apply it with $B = d^{1/2} (\log N)^{m+10}$ and $\delta = T^{1/2-\mathfrak{c}_1} R^{-1/2}$, for some constant $\mathfrak{c}_1>0$. The probability in \eqref{eqn:coupling} will be less than $c^{-1} e^{-c (\log N)^2}$, if $\mathfrak{c}_0$ and $\mathfrak{c}_1$ are sufficiently small so that 
\begin{flalign}
    \label{c0c1}
T^{1/2-3\mathfrak{c}_0 - \mathfrak{c}_1} R^{-1/2} \ge 1.
\end{flalign} 

Thus, in what follows, we fix the constants $\mathfrak{c}_0$ and $\mathfrak{c}_1$ to satisfy \eqref{c0c1}. Then  Lemma \ref{lem:final_couple} (with \eqref{xisum}) yields a coupling between the random Lax matrix $\L(0)$ and a Gaussian random vector 
\begin{equation}
\label{eqn:big_N_dep_Gaus}
(X_{l_1;1},X_{l_1;0})_{l_1 \in \mathfrak{L}_1} \cup  (Y_{l_2;1}^{[m]} , Y_{l_2;0}^{[m]})_{l_2 \in \mathfrak{L}_2, 0 \leq m \leq M_0},
\end{equation}
such that with overwhelming probability we have
\begin{align}
    \label{xxigxgy} 
    \begin{aligned} 
\displaystyle\max_{l_1 \in \mathfrak{L}_1} \left|  \Xi_1(l_1) - X_{l_1;1} \right| &\leq T^{-\mathfrak{c}_1}; \qquad 
\displaystyle\max_{m \in \llbracket 0, M_0 \rrbracket} \displaystyle\sup_{l_2 \in \mathfrak{L}_2}  \left|  \Xi_1^{[m]}(l_2) - Y_{l_2;1} \right| \leq T^{-\mathfrak{c}_1}; \\
  \displaystyle\max_{l_1 \in \mathfrak{L}_1} \left|  \mathcal{G}_N(l_1) - X_{l_1;0} \right| &\leq  T^{-\mathfrak{c}_1}; \qquad 
 \displaystyle\max_{m \in \llbracket 0, M_0 \rrbracket} \displaystyle\sup_{l_2 \in \mathfrak{L}_2}   \left|  \mathcal{G}_N^{(\varsigma_m)}(l_2) - Y_{l_2;0} \right| \leq T^{-\mathfrak{c}_1},
 \end{aligned} 
\end{align}

\noindent and such that the covariances between the entries of \eqref{eqn:big_N_dep_Gaus} coincide with those between  
\begin{equation} 
\label{xi2g2}
(\Xi_1 (l_1), \mathcal{G}_N (l_1) )_{l_1 \in \mathfrak{L}_1} \cup  (\Xi_1^{[m]} (l_2), \mathcal{G}_N^{(\varsigma_m)} (l_2) \big)_{l_2 \in \mathfrak{L}_2, 0 \leq m \leq M_0},
\end{equation}

Denote the Gaussian vector \eqref{eqn:big_N_dep_Gaus} by $W = (W_1, W_2, \ldots , W_d)$. We seek to couple $W$ with the Gaussian vector $\tilde{W} = (\tilde{W}_1, \tilde{W}_2, \ldots , \tilde{W}_d)$ given by 
\begin{equation}
\label{x0x2y2y0}
(\mathfrak{X}_1 (l_1), \mathfrak{X}_0 (l_1))_{l_1 \in \mathfrak{L}_1} \cup  (\mathfrak{Y}_1^{[m]} (l_2), \mathfrak{Y}_0^{[m]} (l_2))_{l_2 \in \mathfrak{L}_2, 0 \leq m \leq M_0},
\end{equation}

\noindent While the vector $W$ describing \eqref{eqn:big_N_dep_Gaus} and the vector $\tilde{W}$ describing \eqref{x0x2y2y0} are Gaussian, the covariance of the former is dependent on $N$ and therefore does not exactly equal that of the latter. To show that $W$ and $\tilde{W}$ can be coupled to approximately coincide, we must estimate the difference between the covariance matrices of $W$ and $\tilde{W}$. We denote these covariance matrices by $\Sigma = (\Sigma_{ij})$ and $\tilde{\Sigma} = (\tilde{\Sigma}_{ij})$, respectively, so that $\Sigma_{ij} = \Cov (W_i, W_j)$ and $\tilde{\Sigma}_{ij} = \Cov (\tilde{W}_i, \tilde{W}_j)$ for all $i, j \in \llbracket 1, d \rrbracket$.

 Then by Lemma \ref{lem:cov_error_estimate_specific} and Remark \ref{mufr} (together with \eqref{y001}, \eqref{x001}, and \eqref{xisum} for the entries of \eqref{xi2g2} that prescribe $\Sigma$, and \eqref{eqn:X1out}, \eqref{eqn:Y1out}, \eqref{eqn:X0eq}, and \eqref{eqn:Y0eq} for the entries \eqref{x0x2y2y0} of $\tilde{W}$ that with \eqref{eqn:wdr_cov} prescribe $\tilde{\Sigma}$), we deduce for some constant $\mathfrak{c}_2 > 0$ that 
\begin{equation}\label{eqn:cov_mat_dif}
\displaystyle\max_{i,j \in \llbracket 1, d \rrbracket} |\Sigma_{i j} - \tilde \Sigma_{i j} | \leq T^{-\mathfrak{c}_2}.
\end{equation}

\noindent Fix the constants $\mathfrak{c}_0$ and $\mathfrak{c}_1$ to be sufficiently small so that 
 \begin{equation}\label{eqn:delt0}
 T^{-\mathfrak{c}_0 - \mathfrak{c}_1} \geq \delta_0 \geq  T^{-\mathfrak{c}_2/20},
 \end{equation}

\noindent holds (in addition to \eqref{c0c1}). Since the dimension of $W$ is $d \leq T^{\mathfrak{c}_2/20}$,  by Lemma \ref{lem:gauss_couple} and the bound \eqref{eqn:cov_mat_dif}, we can couple $W$ with $\tilde{W}$,  such that, with overwhelming probability,
 \begin{equation}\label{eqn:lim_gaussian2}
\displaystyle\max_{i,j \in \llbracket 1, d \rrbracket} |\tilde W_i - W_i | \leq T^{-\mathfrak{c}_2/4}.
 \end{equation}

\noindent Combining \eqref{eqn:lim_gaussian2} and \eqref{xxigxgy} yields a coupling between the random variables  
\begin{multline}\label{eqn:first_joint_gauss_approx2}
    \{\mathcal{G}_N(\Lambda, T \kappa,T\tau)\}_{(\Lambda,  \kappa, \tau) \in \mathfrak{L}_1} 
    \cup  \{\Xi_1(\Lambda, T \kappa, T \tau)\}_{(\Lambda,  \kappa, \tau) \in \mathfrak{L}_1} \\ 
    \cup\{ \mathcal{G}_N^{(\varsigma_m)}(T \mathfrak{q},T\mathfrak{q}', \tau) \}_{(\mathfrak{q},\mathfrak{q}', \tau) \in \mathfrak{L}_2,  0 \leq m \leq M_0}
    \cup \{ \Xi_1^{[m]}(T \mathfrak{q}, T \mathfrak{q}', T \tau) \}_{(\mathfrak{q},\mathfrak{q}', \tau) \in \mathfrak{L}_2,  0 \leq m \leq M_0}.
\end{multline}

\noindent and $(\mathfrak{X}_0, \mathfrak{X}_1, \mathfrak{Y}_0, \mathfrak{Y}_1)$ such that, for some constant $\mathfrak{c}_4>0$, with overwhelming probability we have 
 \begin{equation}\label{eqn:lim_gaussian}
 \begin{aligned} 
 \displaystyle\max_{(\Lambda,\kappa,\tau) \in \mathfrak{L}_1} \big( |\mathfrak{X}_0(\Lambda, \kappa,\tau) - \mathcal{G}_N(\Lambda, T \kappa,T\tau) | + |\mathfrak{X}_1 (\Lambda, \kappa,\tau) - \Xi_1 (\Lambda, T \kappa,T\tau) | \big) \leq T^{-\mathfrak{c}_4}; \\
 \displaystyle\max_{\substack{(\mathfrak{q}, \mathfrak{q}',\tau) \in \mathfrak{L}_2 \\ m \in \llbracket 0, M_0 \rrbracket}} \big( |\mathfrak{Y}_0^{[m]} (\mathfrak{q}, \mathfrak{q}',\tau) - \mathcal{G}_N^{(\varsigma_m)} (T\mathfrak{q}, T\mathfrak{q}', T\tau) | + |\mathfrak{Y}_1^{[m]} (\mathfrak{q}, \mathfrak{q}',\tau) - \Xi_1^{[m]} (T\mathfrak{q}, T\mathfrak{q}', T\tau) | \big) \leq T^{-\mathfrak{c}_4}.
  \end{aligned} 
 \end{equation}

Now we extend from $\mathfrak{L}_1, \mathfrak{L}_2$ to all parameter values using H\"older bounds for the prelimiting processes (shown in Section \ref{sec:continuity}) and the limiting ones (shown as Lemma \ref{lem:holder_cont}). Specifically, from the definitions \eqref{eqn:lat_1}, \eqref{eqn:lat_2}, \eqref{eqn:lat_3}, and \eqref{eqn:lat_4}, we see that $\delta_0 \leq T^{-\mathfrak{c}_3}$ implies that any point $(\Lambda, \kappa, \tau) \in [- \log N, \log N] \times [-(\log N)^4, (\log N)^4] \times [0, \log N]$ is at of distance at most $T^{1-c}$ away from some point $(\tilde{\Lambda}, \tilde{\kappa}, \tilde{\tau}) \in \mathfrak{L}_1$. Similarly, any point $(\mathfrak{q}, \mathfrak{q}', \tau) \in [-(\log N)^5, (\log N)^5] \times [-(\log N)^5, (\log N)^5] \times [0, \log N]$ is of distance at most $T^{-c}$ from some point $(\tilde{\mathfrak{q}}, \tilde{\mathfrak{q}}', \tilde{\tau}) \in \mathfrak{L}_2$. Therefore, by Lemma \ref{lem:xi_xim_cont}; Lemma \ref{lem:gen_brownian_cont_allparam_prelim} with $F = \varsigma_m$; and Lemma \ref{lem:Xi_brownian_cont_allparam_pl}, all applied with the $S$ there equal to $T^{1-c}$ here, we deduce with overwhelming probability that the, under the above notation, 
 \begin{flalign}\label{xilambdalambda}
 \begin{aligned} 
 |\mathcal{G}_N (\Lambda, T \kappa, T\tau) - \mathcal{G}_N( \tilde{\Lambda}, T \tilde{\kappa},T \tilde{\tau}) | + |\Xi_1 (\Lambda, T \kappa, T\tau) - \Xi_1 ( \tilde{\Lambda}, T \tilde{\kappa},T \tilde{\tau}) | \le T^{-c}; \\
  |\mathcal{G}_N^{(\varsigma_m)} (T\mathfrak{q}, T\mathfrak{q}', T\tau) - \mathcal{G}_N^{(\varsigma_m)} (T \tilde{\mathfrak{q}}, T \tilde{\mathfrak{q}}', T \tilde{\tau}) | + |\Xi_1^{[m]} ( T\mathfrak{q}, T\mathfrak{q}', T\tau) - \Xi_1^{[m]} (T \tilde{\mathfrak{q}}, T \tilde{\mathfrak{q}}', T \tilde{\tau})|  \leq T^{-c}.
  \end{aligned} 
 \end{flalign}

\noindent Analogously, the H\"{o}lder continuity of $(\mathfrak{X}_0, \mathfrak{X}_1, \mathfrak{Y}_0, \mathfrak{Y}_1)$ provided by Corollary \ref{lem:holder_cont} implies that  
  \begin{flalign}\label{x02y02}
 \begin{aligned} \ |\mathfrak{X}_0(\Lambda, \kappa,\tau) - \mathfrak{X}_0(\tilde{\Lambda}, \tilde{\kappa}, \tilde{\tau}) | + |\mathfrak{X}_1 (\Lambda, \kappa,\tau) - \mathfrak{X}_1 (\tilde{\Lambda}, \tilde{\kappa}, \tilde{\tau}) |  \leq T^{-c}; \\
 |\mathfrak{Y}_0^{[m]} (\mathfrak{q}, \mathfrak{q}',\tau) - \mathfrak{Y}_0^{[m]} (\tilde{\mathfrak{q}}, \tilde{\mathfrak{q}}', \tilde{\tau}) | + |\mathfrak{Y}_1^{[m]} ( \mathfrak{q}, \mathfrak{q}',\tau) - \mathfrak{Y}_1^{[m]} (\tilde{\mathfrak{q}}, \tilde{\mathfrak{q}}', \tilde{\tau}) |  \leq T^{-c}.
  \end{aligned} 
 \end{flalign}

\noindent The theorem then follows from combining \eqref{xigxi}, \eqref{eqn:lim_gaussian}, \eqref{xilambdalambda}, and\eqref{x02y02}.
\end{proof}

Next, we provide the following generalization of Theorem \ref{thm:couplethm12} (with essentially the same proof), which will be useful to us in Section \ref{sec:current_fluct} below.

\begin{thm}\label{thm:couplethmFG}

There exists a constant $\mathfrak{c}>0$ such that the following holds. Let $A, B \in [1, (\log N)^3]$ be real numbers, and let $K \in \llbracket 1, \log N \rrbracket$ be an integer. For each $i \in \llbracket 1, K \rrbracket$, let $\tau_i \in [0, \log N]$ and $U_i \in [\mathfrak{M}, T \log N]$  be real numbers, and let $F_i, G_i : \mathbb{R} \rightarrow \mathbb{R}$ be functions, such that each $F_i$ satisfies Assumption \ref{ass:F}, and each $G_i$ has compact support and satisfies Assumption \ref{ass:G2} with this choice of $B$ and with $U = U_i$. For each $i \in \llbracket 1, K \rrbracket$, also set $H_i(\lambda, q) \coloneqq F(\lambda) \cdot G(q)$ and, defining $\Xi^{(H_i)} : \mathbb{R} \times \mathbb{R}_{\ge 0} \rightarrow \mathbb{R}$ as in \eqref{eqn:H_functional}, let
\begin{equation}
    \bar{\Xi}^{(i)}(t) \coloneqq T^{-1/2} \cdot \big(\Xi^{(H_i)}(0,t) - \mathbb{E}[\Xi^{(H_i)}(0,t) ] \big).
\end{equation}

\noindent For any index $i \in \llbracket 1, K \rrbracket$, and real numbers $\tau, r, \lambda \in \mathbb{R}$, further define 
\begin{equation}\label{eqn:phi_i_lim}
    \phi_{\tau}^{(i)}(r, \lambda) \coloneqq G( T \alpha r  + T \tau \ve(\lambda) ) \cdot F(\lambda); \qquad 
    \mathcal{G}^{(i)}(\tau) \coloneqq \mathcal{W}^{\dr}(\phi_{\tau}^{(i)}) .
\end{equation}

\noindent  Then there exists a coupling between $\mathcal{W}^{\dr}$ and $\L(0)$ such that, with overwhelming probability,  
    \begin{equation}\label{eqn:couplei}
        \sup_{i \in \llbracket 1, K \rrbracket} | \bar{\Xi}^{(i)}(T \tau_i) - \mathcal{G}^{(i)}(\tau_i)  | \leq T^{-\mathfrak{c}} \qquad \text{and} \quad \text{\eqref{eqn:xi_couple_outline}} \qquad \text{simultaneously hold}. 
    \end{equation}

\end{thm}

\begin{proof}[Proof (Outline)]

The proof of this theorem is similar to that of Theorem \ref{thm:couplethm12}, so we only briefly describe the differences. Recall that in the proof of  Theorem \ref{thm:couplethm12} that we approximately expressed $(\Xi_1, \mathcal{G}_N, \Xi_1^{[m]}, \mathcal{G}_N^{(\varsigma_m)})$ as sums of independent random variables, through \eqref{xisum} (with \eqref{y001} and \eqref{x001}). Throughout, we restrict to the overwhelmingly probable event that \eqref{xisum} holds. 

The counterparts of such expressions for $\bar{\Xi}^{(i)}$ is as follows. For any indices $r \in \llbracket N_1/R, N_2/R \rrbracket$ and $i \in \llbracket 1, K \rrbracket$, set  
    \begin{flalign*}
  Z_{i;r}^{(N)} & =    R^{-1/2} \sum_{s=1}^{R}  H \big(\lambda_i^{(r)}, \alpha r R  + T \tau_i \ve(\lambda_i^{(r)}) \big) - \mathbb{E}[\cdots] \\
  & \qquad \qquad \qquad - \bigg( (\alpha R^{1/2})^{-1} \cdot c^{(H_i)}(T \tau_i; r)  \sum_{j=rR}^{(r+1)R-1} X_j - \mathbb{E}[\cdots] \bigg),
  \end{flalign*}
  
  \noindent where, for any $t, r \in \mathbb{R}$, we have denoted
  \begin{flalign*} 
  c^{(H_i)}(t; r) \coloneqq \int_{-\infty}^{\infty} H_i(\lambda, \alpha r R + t \ve(\lambda) ) \cdot \varrho(\lambda) d\lambda.
    \end{flalign*} 

\noindent We then further restrict to the overwhelmingly probable event provided by Item \ref{item:isth2} of Theorem \ref{thm:Hindsum_thm}, on which we have for each $i \in \llbracket 1, K \rrbracket$ that 
    \begin{equation}\label{eqn:Xii_dcmp}
\bar{\Xi}^{(i)}(T \tau_i) = R^{1/2} T^{-1/2} \sum_{r = \lceil N_1/R \rceil}^{\lfloor N_2/R \rfloor - 1}  Z_{i; r}^{(N)}  + O(T^{-\mathfrak{c}}).
\end{equation}

By Definitions \ref{xi} and \ref{def:Lr}, the $(Z_{i;r}^{(N)})$ are mutually independent over the index $i$. Moreover, by the first part of Lemma \ref{lem:no_q_conc} and Lemma \ref{lem:q_spacing}, we have with overwhelming probability that $|Z_{i;r}^{(N)}| \le (\log N)^C$ for each $i \in \llbracket 1, K \rrbracket$ and $r \in \llbracket N_1/R,  N_2/R-1 \rrbracket$. Thus, applying Lemma \ref{lem:final_couple} (with the $B$ there equal to $2 d^{1/2} (\log N)^{C}$, with $d$ as in \eqref{eqn:ddef}, and with the $\delta$ there equal to $T^{1/2} R^{-1/2} T^{-\mathfrak{c}_1}$ for a sufficiently small constant $\mathfrak{c}_1>0$ satisfying \eqref{c0c1}), we obtain a coupling between the Lax matrix $\mathbf{L}(0)$ and a Gaussian random vector 
\begin{flalign}
    \label{zxxyy}
    (Z_i)_{i \in \llbracket 1, K \rrbracket} \cup (X_{l_1;1},X_{l_1;0})_{l_1 \in \mathfrak{L}_1} \cup  (Y_{l_2;1}^{[m]} , Y_{l_2;0}^{[m]})_{l_2 \in \mathfrak{L}_2, 0 \leq m \leq M_0},
\end{flalign}

\noindent such that we have \eqref{xxigxgy}; $|Z_i - \bar{\Xi}^{(i)} (T\tau_i)| \le T^{-\mathfrak{c}_1}$; and the covariances between the entries of \eqref{zxxyy} coincide with those between 
\begin{flalign}
\label{w2}
(\bar{\Xi}^{(i)} (T\tau_i))_{i \in \llbracket 1, K \rrbracket} \cup (\Xi_1 (l_1), \mathcal{G}_N (l_1) )_{l_1 \in \mathfrak{L}_1} \cup  ( \Xi_1^{[m]} (l_2), \mathcal{G}_N^{(\varsigma_m)} (l_2) \big)_{l_2 \in \mathfrak{L}_2, 0 \leq m \leq M_0}.
\end{flalign} 

Denote the Gaussian vector \eqref{zxxyy} by $W = (W_1, W_2, \ldots , W_D)$. We then seek to couple $W$ with the Gaussian vector $\tilde{W} = (\tilde{W}_1, \tilde{W}_2, \ldots , \tilde{W}_D)$ given by 
\begin{flalign}
    \label{w3} 
(\mathcal{G}^{(i)} (\tau_i))_{i \in \llbracket 1, K \rrbracket} \cup (\mathfrak{X}_1 (l_1), \mathfrak{X}_0 (l_1))_{l_1 \in \mathfrak{L}_1} \cup  (\mathfrak{Y}_1^{[m]} (l_2), \mathfrak{Y}_0^{[m]} (l_2))_{l_2 \in \mathfrak{L}_2, 0 \leq m \leq M_0},
\end{flalign} 

\noindent Denoting the covariances matrices of these vectors by $\Sigma = [\Sigma_{ij}]$ and $\tilde{\Sigma} = [\tilde{\Sigma}_{ij}]$, where $\Sigma_{ij} = \Cov (W_i, W_j)$ and $\tilde{\Sigma}_{ij} = \Cov (\tilde{W}_i, \tilde{W}_j)$, it follows quickly from Lemma \ref{lem:cov_error_estimate_specific}, Remark \ref{mufr}, \eqref{eqn:phi_i_lim}, and \eqref{eqn:wdr_cov} (as in the proof of Theorem \ref{thm:couplethm12}) that $|\Sigma_{ij} - \tilde{\Sigma}_{ij}| < T^{-\mathfrak{c}_3}$ for some constant $\mathfrak{c}_3 > 0$. By Lemma \ref{lem:gauss_couple}, it follows that we can couple $W$ and $\tilde{W}$ such that $|W_i - \tilde{W}_i| \le T^{-c}$ for all $i$. Hence, there exists a coupling between the random variables \eqref{w2} and \eqref{w3} such that \eqref{eqn:lim_gaussian} holds, and we have $|\bar{\Xi}^{(i)} (T\tau_i) - \mathcal{G}^{(i)} (\tau_i)| < T^{-c}$ for each $i \in \llbracket 1, K \rrbracket$. This establishes the first statement in \eqref{eqn:couplei}; the proof of the second is entirely analogous to that in the proof of Theorem \ref{thm:couplethm12} (and follows from \eqref{xilambdalambda} and \eqref{x02y02}, with \eqref{eqn:lim_gaussian}). 
\end{proof}

\subsection{Proof of Lemma \ref{lem:cov_error_estimate_specific}}
\label{subsec:couplepreproofs}

To show Lemma \ref{lem:cov_error_estimate_specific}, we begin with the following lemma.

\begin{lem}\label{lem:var_error_estimate}

There exists a constant $\mathfrak{c}>0$ such that the following holds. Below,  $f_{r}^{(N)}$ denotes any prelimiting function, as in Definition \ref{ff}, and $f_{r}$ denotes the corresponding limiting function. 

\begin{enumerate}
    \item For any real numbers $r, r' \in \mathbb{R}$ with $r' \in [r, r+1]$, 
\begin{equation}\label{eqn:rdif_Hnorm}
    \|f_{r}^{(N)} - f_{r'}^{(N)} \|_{\mathcal{H}} \leq T^{-\mathfrak{c}} .
\end{equation} \label{item:limvar1}
\item For any integer $r \in \llbracket N_1/R, N_2/R \rrbracket$ with $|r| > TR^{-1} (\log N)^7$,  
\begin{equation}\label{eqn:variance_tailbd}
  \frac{1}{R} \var \left(   \sum_{i=1}^{R}f_r^{(N)}(\lambda_i^{(r)}) \right)  \leq \mathfrak{c}^{-1} e^{-\mathfrak{c} (\log N)^2}.
\end{equation} \label{item:limvar2}
    \item For any real number $r \in \mathbb{R}$ with $|r - q/(T\alpha)| > (\mathfrak{M}/T)^{1/2}$ and $|r - q'/(T\alpha)| > (\mathfrak{M}/T)^{1/2}$, 
   \begin{equation}\label{eqn:chi_1_dif}
       \| f_{r T / R}^{(N)} - f_{r} \|_{\mathcal{H}} \leq T^{-\mathfrak{c}}
   \end{equation}\label{item:limvar3}
   \item For any real number $r \in \mathbb{R}$ with $|r| \geq (\log N)^7$,  
   \begin{equation}\label{eqn:tilder_tailbd}
   \| f_{r} \|_{\mathcal{H}} \leq \mathfrak{c}^{-1} e^{-\mathfrak{c} r^2}.
\end{equation} \label{item:limvar4}
\end{enumerate}

\end{lem}

Before proving Lemma \ref{lem:var_error_estimate}, let us establish Lemma \ref{lem:cov_error_estimate_specific} assuming it.

\begin{proof}[Proof of Lemma \ref{lem:cov_error_estimate_specific}]

We only establish \eqref{eqn:covariance_bd_touse} in detail, as the proofs of \eqref{eqn:cov_f_br_bd_to_use} and \eqref{eqn:cov_br_br_bd_to_use} are entirely analogous. Throughout this discussion, for any two functions $\phi_1, \phi_2 \in \mathcal{H}$, we abbreviate  
\begin{equation}\label{eqn:tildmthscrC}
    \tilde{\mathscr{C}}(\phi_1, \phi_2) \coloneqq \mathscr{C}(\phi_1 - \mu_{\phi_1} \varsigma_0, \phi_2- \mu_{\phi_2} \varsigma_0).
\end{equation}

\noindent Observe for any $\phi_1, \phi_2 \in \mathcal{H}$ that 
\begin{flalign}
    \label{estimatec0}
    \begin{aligned} 
    \tilde{\mathscr{C}} (\phi_1, \phi_2) & = \| (1 - \theta \mathbf{T} \bm{\varrho}_{\beta})^{-1} (\phi_1 - \mu_{\phi_1} \varsigma_0) \|_{\mathcal{H}} \cdot \| (1 - \theta \mathbf{T} \bm{\varrho}_{\beta}) (\phi_1 - \mu_{\phi_1} \varsigma_0) \|_{\mathcal{H}} \\
    & \le C \cdot \| \phi_1 - \mu_{\phi_1} \varsigma_0 \|_{\mathcal{H}} \cdot \| \phi_2 - \mu_{\phi_2} \varsigma_0 \|_{\mathcal{H}} \le 4C \cdot \| \phi_1 \|_{\mathcal{H}} \cdot \| \phi_2 \|_{\mathcal{H}},
    \end{aligned} 
\end{flalign}

\noindent where the first statement follows from \eqref{eqn:tildmthscrC} and \eqref{eqn:sig2def}; the second from the boundedness of the operator $(1 - \theta \mathbf{T} \bm{\varrho}_{\beta})^{-1} : \mathcal{H} \rightarrow \mathcal{H}$; and the third from the fact that $|\mu_{\phi_i}| = |\langle \phi_i, \varsigma_0 \rangle_{\varrho}| \le \| \phi_i \|_{\mathcal{H}}$.

First, we argue that for $r \in \mathbb{R} $, and any $r' \in [r, r+1]$,
    \begin{equation}\label{eqn:limvar_fR}
\Bigg| \frac{1}{R} \Cov \left(   \sum_{i=1}^{R}f_r^{(N)}(\lambda_i^{(R)}),  \sum_{i=1}^{R}g_r^{(N)}(\lambda_i^{(R)}) \right) -\tilde{\mathscr{C}}(f_{r'}^{(N)}, g_{r'}^{(N)} ) \Bigg| \leq    R^{-c}.
\end{equation} 

\noindent To that end, first observe by three applications of Proposition \ref{prop:limvar} (with the $\phi$ there equal to each of $(f_r^{(N)}, g_r^{(N)}, f_r^{(N)} + g_r^{(N)})$) that
   \begin{equation}\label{eqn:limvar_R_0}
\Bigg| \frac{1}{R} \Cov \left(   \sum_{i=1}^{R}f_r^{(N)}(\lambda_i^{(R)}),  \sum_{i=1}^{R}g_r^{(N)}(\lambda_i^{(R)}) \right) -\tilde{\mathscr{C}}(f_{r}^{(N)}, g_{r}^{(N)}) \Bigg| \leq    R^{-c}.
\end{equation} 

\noindent To deduce \eqref{eqn:limvar_fR}, it suffices to show for $|r-r'| \leq 1$ that
\begin{equation}\label{eqn:limvar_R_1}
\big|\tilde{\mathscr{C}}(f_{r}^{(N)}, g_{r}^{(N)}) -\tilde{\mathscr{C}}(f_{r'}^{(N)}, g_{r'}^{(N)}) \big| \leq    R^{-c}.
\end{equation} 

 \noindent By \eqref{estimatec0}, we have 
\begin{flalign}\label{eqn:sigma2_dif}
\begin{aligned} 
    \big| \tilde{\mathscr{C}}  (f_{r}^{(N)}, g_{r}^{(N)}) -\tilde{\mathscr{C}}(f_{r'}^{(N)}, g_{r'}^{(N)}) \big|  & \leq  \big| \tilde{\mathscr{C}}(f_{r}^{(N)},g_{r}^{(N)} - g_{r'}^{(N)} ) \big|  + \big| \tilde{\mathscr{C}} (g_{r'}^{(N)}, f_{r}^{(N)}-f_{r'}^{(N)}) \big|  \\
   &  \leq  C \big( \| f_{r}^{(N)}\|_{\mathcal{H}} \cdot \|g_{r}^{(N)} - g_{r'}^{(N)} \|_{\mathcal{H}} + \|g_{r'}^{(N)}\|_{\mathcal{H}}  \cdot \|f_{r}^{(N)} - f_{r'}^{(N)} \|_{\mathcal{H}}  \big) .
    \end{aligned} 
\end{flalign}

\noindent Also, by \eqref{eqn:phi_lc_reg}, \eqref{eqn:phi_m_reg}, and \eqref{eqn:fFG_prelim} (together with Definition \ref{eqn:sl} and \eqref{estimatef}), we have that $|f_r^{(N)} (\lambda)| \le C (\log  N)^{10} \lambda^m e^{|\lambda|^{3/2}}$ for all $\lambda \in \mathbb{R}$. Together with \eqref{def:inner_prod} and Lemma \ref{lem:varrho_bd}, this implies that 
\begin{equation}\label{eqn:Hn_triv}
\| f_{r}^{(N)}\|_{\mathcal{H}} + \|g_{r}^{(N)}\|_{\mathcal{H}} \leq C (\log N)^{m+10}.
\end{equation}

\noindent Then \eqref{eqn:rdif_Hnorm} \eqref{eqn:Hn_triv}, and \eqref{eqn:sigma2_dif} together yield \eqref{eqn:limvar_R_1} and thus \eqref{eqn:limvar_fR}; the latter implies that
   \begin{multline}\label{eqn:almostlastbd_var_fNint}
     \Bigg| \frac{R}{T} \sum_{r \in \llbracket - T(\log N)^7 /R, T(\log N)^7 /R\rrbracket  } \frac{1}{R} \Cov\left( \sum_{i=1}^R f_r^{(N)}(\lambda_i^{(r)}), \sum_{i=1}^R g_r^{(N)}(\lambda_i^{(r)}) \right) \\
     - \frac{R}{T}\int_{|r|\leq T(\log N)^7 /R} \tilde{\mathscr{C}}(f_{r}^{(N)},g_{r}^{(N)}) dr \Bigg| \leq T^{-c}.
   \end{multline}

\noindent Together with  \eqref{eqn:variance_tailbd}, this gives 
   \begin{multline}\label{eqn:secondtolastbd_var_fNint}
     \Bigg| \frac{R}{T} \sum_{r \in \llbracket N_1/R, N_2/R-1\rrbracket  } \frac{1}{R} \Cov\left( \sum_{i=1}^R f_r^{(N)}(\lambda_i^{(r)}), \sum_{i=1}^R g_r^{(N)}(\lambda_i^{(r)}) \right) \\
     - \frac{R}{T}\int_{|r|\leq T (\log N)^7 / R} \tilde{\mathscr{C}}(f_{r}^{(N)}, g_{r}^{(N)}) dr \Bigg| \leq T^{-c}.
   \end{multline}

Next, we will show that 
   \begin{equation}\label{eqn:limfNint_cutoff}
      \left|  \int_{|r| \leq (\log N)^7} \tilde{\mathscr{C}}(f_{r}, g_{r})d r - \frac{R}{T}\int_{|r|\leq T (\log N)^7 / R} \tilde{\mathscr{C}}(f_{r}^{(N)}, g_{r}^{(N)}) dr \right| < T^{-c}.
   \end{equation}

   \noindent Let us first obtain a weak bound on the integrands $\tilde{\mathscr{C}} (f_{r}, g_{r})$ and $\tilde{\mathscr{C}} (f_{rT/R}^{(N)}, g_{rT/R}^{(N)})$ in \eqref{eqn:limfNint_cutoff} for all $r \in \mathbb{R}$. Following the derivation of \eqref{eqn:Hn_triv}, observe by \eqref{eqn:psi_r_2f}, \eqref{eqn:psi_m_regf}, and \eqref{eqn:fFG_lim}(together with Definition \ref{eqn:sl} and \eqref{estimatef}), we have that $|f_r^{(N)} (\lambda)| \le C (\log  N)^{10} \lambda^m e^{|\lambda|^{3/2}}$ for all $\lambda \in \mathbb{R}$. Together with \eqref{def:inner_prod} and Lemma \ref{lem:varrho_bd} (and \eqref{eqn:Hn_triv}), this implies that 
\begin{equation*}
\| f_{rT/R}^{(N)}\|_{\mathcal{H}} + \| f_{r} \|_{\mathcal{H}} + \|g_{rT/R}^{(N)}\|_{\mathcal{H}} + \| g_{r} \|_{\mathcal{H}} \leq C (\log N)^{m+10}.
\end{equation*} 

\noindent Together with \eqref{estimatec0}, it follows that 
\begin{flalign}
\label{cfngncfg}
\big| \tilde{\mathscr{C}} (f_{rT/R}^{(N)}, g_{rT/R}^{(N)}) \big| \le C(\log N)^{2m+20}; \qquad \| \tilde{\mathscr{C}} (f_{r}, g_{r}) \big| \le C(\log N)^{2m+20}. 
\end{flalign} 

 Next fix a real number $r \in \mathbb{R}$ with $|r| \le (\log N)^7$, such that $|r - q/(T\alpha)| > (\mathfrak{M}/T)^{1/2}$ and $|r - q'/(T\alpha)| > (\mathfrak{M}/T)^{1/2}$. We claim that 
\begin{equation}\label{eqn:chi_1s2_dif}
       \big| \tilde{\mathscr{C}}(f_{rT/R}^{(N)}, g_{rT/R}^{(N)})-\tilde{\mathscr{C}}(f_{r}, g_{r}) \big| \leq T^{-c}.
   \end{equation}
   
\noindent Indeed, for such $r$, \eqref{eqn:chi_1_dif} yields 
   \begin{equation}\label{eqn:dif_again}
       \| f_{rT/R}^{(N)} - f_{r} \|_{\mathcal{H}}^2 \leq T^{-\mathfrak{c}},
   \end{equation}

 \noindent  so we quickly deduce \eqref{eqn:chi_1s2_dif} by (following the reasoning in) \eqref{eqn:sigma2_dif}. Integrating \eqref{eqn:chi_1s2_dif} over $r \in [-(\log N)^7, (\log N)^7]$ satisfying  $|r - q/(T\alpha)| > (\mathfrak{M}/T)^{1/2}$ and $|r - q'/(T\alpha)| > (\mathfrak{M}/T)^{1/2}$, and integrating \eqref{cfngncfg} over the remaining such $r \in [-(\log N)^7, (\log N)^7]$ (which has measure at most $4(\mathfrak{M}/T)^{1/2} \le T^{-c}$), we obtain \eqref{eqn:limfNint_cutoff}. Together with \eqref{eqn:secondtolastbd_var_fNint}, this gives 

 \begin{multline*}
     \Bigg| T^{-1} \sum_{r = \lceil N_1/R \rceil}^{\lfloor N_2/R \rfloor -1}  \Cov\left( \sum_{i=1}^R f_r^{(N)}(\lambda_i^{(r)}), \sum_{i=1}^R g_r^{(N)}(\lambda_i^{(r)}) \right)   - \int_{|r| \leq (\log N)^7} \tilde{\mathscr{C}}(f_{r}, g_{r}) d r \Bigg| \leq T^{-c}.
   \end{multline*}

\noindent It therefore remains to verify that 
   \begin{equation}
       \int_{|r| > (\log N)^7} \tilde{\mathscr{C}}(f_{r}, g_{r}) d r < Ce^{-(\log N)^2}.
   \end{equation}

   \noindent which holds by \eqref{estimatec0} and \eqref{eqn:tilder_tailbd}. 
   \end{proof}

\begin{proof}[Proof of Lemma \ref{lem:var_error_estimate}]

We only establish the lemma for functions $f_{r}^{(N)}$ from the second family of in Definition \ref{ff}, so we assume throughout that $f_r^{(N)}$ is of the form \eqref{eqn:psi_m_regf}). Indeed, the proofs when $f_r^{(N)}$ is either in the first and third families there are entirely analogous (for the third, we replace our estimates on $\chi$ below with ones on $G$, using Assumption \ref{ass:G2}) and thus are omitted. 

 To prove Item \ref{item:limvar1} of the lemma, fix $r, r' \in \mathbb{R}$ with $|r-r'| \le 1$. Since $|\chi| \leq 10\mathfrak{M}^{-1}$ (by Definition \ref{def:smoothed_log}), we have from \eqref{eqn:phi_m_reg} that 
\begin{flalign*} 
|f_{r}^{(N)}(\lambda) - f_{r'}^{(N)}(\lambda)| \leq 20 RM^{-1} |\lambda|^m,
\end{flalign*} 

\noindent for all $\lambda \in \mathbb{R}$. Together with \eqref{def:inner_prod}, this gives
\begin{equation}\label{eqn:frrp_hbd}
     \|f_{r}^{(N)} - f_{r'}^{(N)} \|_{\mathcal{H}}^2 \leq 100 (R \mathfrak{M}^{-1})^2 \int_{-\infty}^{\infty} |\lambda|^{2m} \varrho(\lambda) d \lambda \le T^{-c},
\end{equation}

\noindent where in the last inequality used the facts that $m \le (\log N)^{1/10}$ and that $\varrho (\lambda) \le C e^{-c|\lambda|}$ by Lemma \ref{lem:varrho_bd} (with the bound $R\mathfrak{M}^{-1} \le T^{-1/2}$ by Definition \ref{def:smoothed_log} and \eqref{ra0}). This confirms \eqref{eqn:rdif_Hnorm}.

 To show \eqref{eqn:variance_tailbd}, first observe  that, if $|r| > TR^{-1} (\log N)^7$ then, whenever $|q|, |q'|, |t| \le T(\log N)^4$,
\begin{equation}\label{eqn:chidifzero}
\chi(q - \alpha  r R)  - 
    \chi\left( q' - \alpha r R - t \ve(\lambda) \right) = 0, \qquad \text{for all $\lambda \in [-\log N, \log N]$}. 
    \end{equation}
    Indeed, since $|\ve(\lambda) | \leq C \log N$ (by Lemma \ref{lem:ve_tail}), we have for such $q$, $q'$, and $t$ that $|q - \alpha r R | > T (\log N)^6 > 10 \mathfrak{M}$, that $|q' - \alpha r R - t \ve(\lambda) | > 10 \mathfrak{M}$, and that $q - r R $ and $q' - r R - t \ve(\lambda) $ have the same sign. Together with the fact that $\supp \chi' \subseteq [-\mathfrak{M}, \mathfrak{M}]$, this implies \eqref{eqn:chidifzero}. So, on the event $\mathsf{E} = \bigcap_{i=1}^R \{ |\lambda_i^{(r)}| \le \log N \}$ that is overwhelmingly probable by Lemma \ref{lem:bd_lem}, we have 
$$\chi(q - \alpha r R)  - 
    \chi ( q' - \alpha r R - t \ve(\lambda_i^{(R)}) ) = 0,$$

    \noindent for all $i \in \llbracket 1, R \rrbracket$. Therefore, by \eqref{eqn:phi_m_reg}, we have $f_r^{(N)}(\lambda_i^{(R)}) \cdot \mathbbm{1}_{\mathsf{E}} = 0$ for all $i$. Thus,
    \begin{flalign*}
        \Var \Bigg( \displaystyle\sum_{i=1}^R f_r^{(N)} (\lambda_i^{(R)}) \Bigg) = \Var \Bigg( \mathbbm{1}_{\mathsf{E}^{\complement}} \cdot \displaystyle\sum_{i=1}^R f_r^{(N)} (\lambda_i^{(R)}) \Bigg).
    \end{flalign*}
    
    \noindent Since $\mathbb{P}[\mathsf{E}^{\complement}] < Ce^{-c(\log N)^2}$, it is quickly verified (for example, following the derivation of \eqref{expectationH}) that the right side of the above inequality is at most $C e^{-c(\log N)^2}$, which confirms \eqref{eqn:variance_tailbd}.

 Next we show Item \ref{item:limvar3} of the lemma, so fix $r \in \mathbb{R}$ satisfying $|r - q/(T\alpha)| > (\mathfrak{M}/T)^{1/2}$ and $|r - q'/(T\alpha)| > (\mathfrak{M}/T)^{1/2}$. First, observe that 
 \begin{flalign}
     \label{frnfr}
     \begin{aligned} 
    \| f_r^{(N)} - f_{r} \|_{\mathcal{H}}^2 & \le \displaystyle\int_{-\infty}^{\infty}  \big| \chi(q-\alpha r T)-\chi(q'-\alpha rT - t \ve(\lambda)) \\
    & \qquad \qquad \qquad   - \mathbbm{1}_{q/T-\alpha r >0} + \mathbbm{1}_{q'/T-\alpha r  - t \ve(\lambda)/T}) \big|\cdot |\lambda|^{2m} \varrho (\lambda) d \lambda \\
    & \le \displaystyle\int_{-\infty}^{\infty} (\mathbbm{1}_{|q-\alpha r T| \le \mathfrak{M}} + \mathbbm{1}_{|q-\alpha r T - t \ve(\lambda)| \le \mathfrak{M}}) \cdot  |\lambda|^{2m} \varrho (\lambda) d \lambda \\
    & = \displaystyle\int_{-\infty}^{\infty} \mathbbm{1}_{|q-\alpha r T - t \ve(\lambda)| \le \mathfrak{M}} \cdot  |\lambda|^{2m} \varrho (\lambda) d \lambda  \\
    & \le \displaystyle\int_{-\log N}^{\log N}  \mathbbm{1}_{|q-\alpha r T - t \ve(\lambda)| \le \mathfrak{M}} \cdot  |\lambda|^{2m} \varrho (\lambda) d \lambda + C e^{-(\log N)^2},
    \end{aligned} 
 \end{flalign}
 
 \noindent where in the first bound we used \eqref{eqn:phi_m_reg}, \eqref{eqn:phim_intro}; in the second we used the facts that $\supp \chi \subseteq [-\mathfrak{M}, \mathfrak{M}]$ and $|\chi| \le 1$; in the third we used the fact that $|q - \alpha r T| > \alpha (\mathfrak{M}T)^{1/2} > \mathfrak{M}$; and in the fourth we used the facts that $\varrho(\lambda) \le Ce^{-|\lambda|^2}$ (by Lemma \ref{lem:varrho_bd}) and $m \le (\log N)^{1/10}$.

 Now let us understand when the indicator function on the right side of \eqref{frnfr} is nonzero. Since $|q - \alpha r T| \ge (\mathfrak{M}T)^{1/2}$ and $|\ve (\lambda)| \le C\log N$ for $|\lambda| \le \log N$ (by Lemma \ref{lem:ve_tail}), we have that $|q-\alpha r T - t \ve(\lambda)| \le \mathfrak{M}$ can only hold if $t \ge  (\mathfrak{M} T)^{1/2} (\log N)^{-2}$. In this case, by lower bound on the derivative of $\ve$ in Lemma \ref{lem:veff_inc}, the measure of the set $\{ \lambda \in \mathbb{R} : |q-\alpha r T - t \ve(\lambda)| \le \mathfrak{M} \}$ is at most $t^{-1} \mathfrak{M} \le (\mathfrak{M}/T)^{1/2} (\log N)^{-2} \le T^{-c}$. Inserting this into \eqref{frnfr} (and again using Lemma \ref{lem:varrho_bd} and the bound $m \le (\log N)^{1/10}$) yields \eqref{eqn:chi_1_dif}.

   It remains to prove Item \ref{item:limvar4} of the lemma, so fix $r \in \mathbb{R}$ with $|r| > (\log N)^7$. Then observe, whenever $|\mathfrak{q}|, |\mathfrak{q}'|, |\tau| \le (\log N)^4$, we have that $|q - \alpha r| > (\log N)^6$; that $|\mathfrak{q}' - \alpha r - \tau \ve(\lambda)| > (\log N)^6 \cdot \mathbbm{1}_{|\lambda| \ge \mathfrak{c} |r|}$ for some constant $\mathfrak{c}>0$ (by Lemma \ref{lem:ve_tail}); and that $q - \alpha r$ and $\mathfrak{q}' - \alpha r - \tau \ve(\lambda)$ are of the same sign. Therefore, by \eqref{eqn:psi_m_regf} and \eqref{eqn:phim_intro}, we deduce  $|f_{r}(\lambda)| \le 2 |\lambda|^m \cdot \mathbbm{1}_{|\lambda| \ge \mathfrak{c} |r|}$, and so 
   \begin{flalign*}
       \| f_{r} \|_{\mathcal{H}}^2 = \displaystyle\int_{-\infty}^{\infty} f_{r} (\lambda) \cdot \varrho (\lambda) d \lambda \le 2 \displaystyle\int_{|\lambda| > \mathfrak{c}|r|} |\lambda|^m \cdot \varrho (\lambda) d \lambda \le C e^{-c r^2},
   \end{flalign*}

   \noindent where in the last inequality we applied Lemma \ref{lem:varrho_bd} (with the facts that $|r| \ge (\log N)^7$ and that $m \le (\log N)^{1/10}$). This establishes \eqref{eqn:tilder_tailbd}. 
\end{proof}

\section{Concentration estimate for sums of \texorpdfstring{$Z$}{}}
\label{sec:quasi_coupling}

In this section we establish Proposition \ref{prop:sumZconc}, which provides a concentration estimate for linear combinations of the function $Z$ from \eqref{eqn:Z_outline}. 

Throughout this section, we adopt Assumption \ref{ass:NT_assumption} and recall the notation from Definitions \ref{def:quasi_intro} and \ref{xilambdakt}. Moreover, we assume throughout (even when not explicitly stated) that $\theta < \theta_0(\beta)$ is sufficiently small (so that all of the results from the previous sections apply).

\subsection{Bounds on $Z$}
\label{subsec:prep_lem_Z}

In this section we show several estimates for $Z(\Lambda, Q) = Z(\Lambda, Q, t)$ defined by \eqref{eqn:Z_outline} (below we also recall the process $\mathcal{Z}$ from Definition \ref{def:Z_intro}).

\begin{lem}\label{lem:Z_bd}\label{lem:mathcalZ_bd}
   There exists a constant $\mathfrak{C}>1$ such that the following holds with overwhelming probability. For any real number $\Lambda \in \mathbb{R}$, we have 
    \begin{equation}\label{eqn:small_Lambda_Zbd}
       \displaystyle\sup_{t \in [0, T \log N]} \displaystyle\sup_{Q \in \mathbb{R}} \left| Z(\Lambda, Q, t) \right| \leq  (t + \mathfrak{M})^{1/2} T^{-1/2} (|\Lambda|+1)^2 (\log N)^{\mathfrak{C}},
    \end{equation}

    \noindent and \begin{equation}\label{eqn:small_Lambda_mcZbd}
       \displaystyle\sup_{t \in [0, T \log N]} \displaystyle\sup_{|\mathfrak{q}| \leq (\log N)^4} \left| \mathcal{Z}(\Lambda, \mathfrak{q}, t) \right| \leq  (|\Lambda|+1)^2 (\log N)^{\mathfrak{C}}.
    \end{equation}

\end{lem}

\begin{proof}
Throughout this proof, we restrict to the overwhelmingly probable event $\mathsf{E}$ on which Lemma \ref{lem:xi_bd_outline} holds. Define $f: \mathbb{R} \rightarrow \mathbb{R}$ by setting 
$f(\Lambda) = \Xi(\Lambda, (Q - t \ve(\Lambda))/\alpha, t)$ for each $\Lambda \in \mathbb{R}$. Then Lemma \ref{lem:ve_tail} yields 
   \begin{equation}\label{eqn:bd_fh_stuff}
      |f^{\dr}(\Lambda)| \leq C (|f(\Lambda)| +  \|f\|_{\mathcal{H}} \log(|\Lambda|+2)).
   \end{equation}

   \noindent By \eqref{eqn:Z_outline}, it suffices to prove that the right hand side of \eqref{eqn:bd_fh_stuff} is upper bounded by the right hand side of \eqref{eqn:small_Lambda_Zbd}, since the third item of Lemma \ref{lem:ve_tail} indicates that $|\varsigma_0^{\dr}(\Lambda)| \ge c$ for all $\Lambda \in \mathbb{R}$. The first term satisfies $|f(\Lambda)| \leq (t + \mathfrak{M})^{1/2}(|\Lambda|+1)^2 T^{-1/2}(\log N)^C$ by \eqref{eqn:small_Lambda_bd_outline}, which implies by \eqref{def:inner_prod} and Lemma \ref{lem:varrho_bd} that $\|f\|_{\mathcal{H}} \leq  (t + \mathfrak{M})^{1/2}T^{-1/2} (\log N)^C$. Combining these bounds, and using \eqref{eqn:bd_fh_stuff}, proves \eqref{eqn:small_Lambda_Zbd}. The proof of \eqref{eqn:bd_fh_stuff} is entirely analogous, by using \eqref{eqn:mathfrakX_bd_outline} in place of \eqref{eqn:small_Lambda_bd_outline}. 
\end{proof}

The next lemma provides a weak H\"older continuity bound for $Z(\Lambda, Q, t)$ in the variable $\Lambda$.

\begin{lem}\label{lem:Z_cont_lambda}

    There exist constants $\mathfrak{c}>0$ and $\mathfrak{C}>1$ such that the following holds. Let $S \in [\mathfrak{M}, T]$ be a real number. With overwhelming probability we have, for any real numbers $\Lambda, \Lambda' \in [-\log N, \log N]$ with $|\Lambda - \Lambda'| \leq ST^{-1}$, that
   \begin{equation}\label{eqn:Lambda_Zdif}
   \displaystyle\sup_{t \in [0, T \log N]} \displaystyle\sup_{|Q| \le T^{4/3}} | Z(\Lambda, Q, t) - Z(\Lambda', Q, t) | \leq ( S^{1/4}T^{-1/4} +T^{-\mathfrak{c}}) (\log N)^{\mathfrak{C}} .
    \end{equation}
    In particular, if $S \in [T^{1-\mathfrak{c}}, T]$, then 
    \begin{equation}\label{eqn:Lambda_Zdif_hold}
       \displaystyle\sup_{t \in [0, T \log N]} \displaystyle\sup_{|Q| \le T^{4/3}} | Z(\Lambda, Q, t) - Z(\Lambda', Q, t) | \leq 2S^{1/4}T^{-1/4} (\log N)^{\mathfrak{C}} .
    \end{equation}
\end{lem}

\begin{proof}
Since \eqref{eqn:Lambda_Zdif} implies \eqref{eqn:Lambda_Zdif_hold}, it suffices to prove the former. Throughout the proof, we restrict to the event $\mathsf{E}$ defined as the intersection of the event $\mathsf{BND}_{\L}(\log N)$ from Definition \ref{eventbounded} (which has overwhelming probability by Lemma \ref{lem:bd_lem}); the overwhelmingly probable event of Lemma \ref{lem:Z_bd} on which \eqref{eqn:small_Lambda_Zbd} holds; and the overwhelmingly probable event provided by Proposition \ref{cor:xicont}, with $S$ there given by $(\log N)^4 S$ here, on which  
\begin{equation}\label{eqn:xi_dif_kkp}
    |\Xi(\Lambda, k_1, t) -\Xi(\Lambda', k_2, t) | \leq (\log N)^C( S^{1/4} T^{-1/4} +T^{-c}).
    \end{equation}
for all $\Lambda, \Lambda' \in [-\log N, \log N]$ with $|\Lambda-\Lambda'| \leq S T^{-1} (\log N)^4 $, all $k_1,k_2 \in \llbracket - T^{3/2},  T^{3/2} \rrbracket$ with $|k_1-k_2| \leq S (\log N)^4$, and all $t \in [0, T \log N]$.

Let $f(\Lambda) = \Xi(\Lambda, (Q - t \ve(\Lambda))/\alpha, t)$. We will start from the general equation which holds (by \eqref{operatort} and \eqref{drf}) for any $f \in \mathcal{H}$ and $g = f^{\dr}$,
    \begin{multline}\label{eqn:firstgeq}
        g(x) =  f(x) + 2 \theta \int_{-\infty}^{\infty}\log |x-y| g(y) \varrho_{\beta}(y) dy  = f(x) + 2 \alpha^{-1} \int_{-\infty}^{\infty}\log |x-y| g(y) \varsigma_0^{\dr}(y)^{-1}  \varrho(y) dy,
    \end{multline}
where the last equality is due to the first equation in \eqref{rho2}. Let $h(y) \coloneqq g(y) \varsigma_0^{\dr}(y)^{-1}$. Equation \eqref{eqn:firstgeq} implies (using Cauchy-Schwarz)
\begin{multline}\label{eqn:Delta_g_first}
    |g(\Lambda)-g(\Lambda')| \leq  |f(\Lambda)-f(\Lambda')| \\
    + 2\alpha^{-1} \left( \int_{-\infty}^{\infty}(\log |\Lambda-\lambda|- \log|\Lambda'-\lambda|)^2 \varrho(\lambda) d\lambda \right)^{1/2} \left(\int_{-\infty}^{\infty} h(\lambda)^2 \varrho(\lambda) d\lambda \right)^{1/2}.
\end{multline}
By \eqref{eqn:Z_outline}, since $f(\Lambda) = \Xi(\Lambda, (Q - t \ve(\Lambda))/\alpha, t)$, we have $g(\Lambda) = \varsigma_0^{\dr}(\Lambda) Z(\Lambda, Q, t)$ and $h(\Lambda) = Z(\Lambda, Q, t)$. Also note that (by Lemma \ref{lem:varrho_bd}),  for any $x, x' \in \mathbb{R}$ with $|x-x'| \leq 1$, we have
\begin{equation}\label{eqn:Delta_g_second}
    \left( \int_{-\infty}^{\infty}(\log |x-y|- \log|x'-y|)^2 \varrho(y) dy\right)^{1/2} \leq C |x-x'|^{1/2} .
\end{equation}

Now, we prove that \eqref{eqn:Lambda_Zdif} holds on the event $\mathsf{E}$ defined above. Let $k_1 = (Q - t \ve(\Lambda'))/\alpha$ and $k_2 =(Q - t \ve(\Lambda))/\alpha$. Observe by the first bound in \eqref{eqn:vebd} that $k_1,k_2 \in [-\floor{T^{3/2}}, \floor{T^{3/2}}]$. Note that, using the second bound in \eqref{eqn:vebd},
\begin{equation}\label{eqn:kdif}
    |k_1-k_2| = \alpha^{-1} t | \ve(\Lambda) - \ve(\Lambda')| \leq T (\log N)^4 |\Lambda-\Lambda'| \leq S (\log N)^4.
\end{equation}
Note that the bound \eqref{eqn:xi_dif_kkp} continues to hold even though $k_1,k_2$ are not necessarily integers because we extend $\Xi(\Lambda, k, t)$ to non-integer values of $k $ by linear interpolation (see Definition \ref{xilambdakt}). Thus, on $\mathsf{E}$, by \eqref{eqn:xi_dif_kkp},
    \begin{equation}\label{eqn:xi_dif_two}
    \Bigg|\Xi\left(\Lambda, \frac{Q - t \ve(\Lambda)}{\alpha}, t\right) -\Xi\left(\Lambda', \frac{Q - t \ve(\Lambda')}{\alpha}, t\right) \Bigg| \leq (\log N)^C ( S^{1/4} T^{-1/4} +T^{-\mathfrak{c}}).
    \end{equation}

    Now, we bound the first and second term on the right hand side of \eqref{eqn:Delta_g_first}. Recalling that $f(\Lambda) = \Xi(\Lambda, (Q - t \ve(\Lambda))/\alpha, t)$, we may use \eqref{eqn:xi_dif_two} to bound the first term on the right hand side of \eqref{eqn:Delta_g_first}.  Next, recall also that $g(\Lambda) = \varsigma_0^{\dr}(\Lambda) Z(\Lambda, Q, t)$ and $h(\Lambda) = Z(\Lambda, Q, t)$. By \eqref{eqn:small_Lambda_Zbd} (which is valid on $\mathsf{E}$) and the bound on $\varrho$ from Lemma \ref{lem:varrho_bd}, we have
 \begin{equation}\label{eqn:hnbd}
 \int_{-\infty}^{\infty} h(\lambda)^2 \varrho(\lambda) d\lambda <  (\log N)^C.
 \end{equation}
    To bound the second term on the right hand side of \eqref{eqn:Delta_g_first}, we use \eqref{eqn:Delta_g_second} with $(x,x') = (\Lambda, \Lambda')$ together with \eqref{eqn:hnbd} to obtain a bound of $(\log N)^C S^{1/2}T^{-1/2}$. Adding the obtained bounds on the two terms in the right hand side of \eqref{eqn:Delta_g_first} gives 
     \begin{equation}\label{eqn:g_dif_two}
    |g(\Lambda) -g(\Lambda') | \leq (\log N)^C ( S^{1/4} T^{-1/4} +T^{-c}).
    \end{equation}
    
    Finally, we thus obtain
    \begin{flalign}\label{eqn:SLZdif_two}
    \begin{aligned}
\left| Z(\Lambda, Q, t) - Z(\Lambda', Q, t) \right|  
&\leq  |\varsigma_0(\Lambda)^{-1}- \varsigma_0(\Lambda')^{-1}| \cdot |g(\Lambda)|  + |\varsigma_0(\Lambda')^{-1}|\cdot |g(\Lambda)-g(\Lambda')|
\\
&\leq C |\Lambda-\Lambda'| (\log N)^C  + C |g(\Lambda)-g(\Lambda')|  \\
&\leq (\log N)^C( S^{1/4} T^{-1/4} +T^{-c}).
    \end{aligned}
    \end{flalign}
   Above, we have used the fact that for all $x \in \mathbb{R}$, $|\varsigma_0^{\dr}(x)^{-1}| \leq C$ and $|(\varsigma_0^{\dr})'(x)| \leq C \log (|x|+2)$ (which holds by the third item of Lemma \ref{lem:ve_tail} and the first item of the same lemma applied with $f = \varsigma_0$); the bound \eqref{eqn:small_Lambda_Zbd}; and the bound \eqref{eqn:g_dif_two}. This completes the proof.
\end{proof}

The next lemma provides a weak continuity of $Z(\Lambda, Q, t)$ in the variables $Q$ and $t$.

\begin{lem}\label{lem:Z_cont}
There exists a constant $\mathfrak{c}>0$ such that the following holds with overwhelming probability. For all $Q, Q' \in [- T^{4/3},  T^{4/3}]$ with $|Q - Q'|\leq  \mathfrak{M}$; all $t,t' \in [0, T \log N]$ with $|t-t'| \leq \mathfrak{M}$; and all $\Lambda \in [-\log N, \log N]$, we have 
\begin{equation}\label{eqn:Z_dif_bound}
|Z(\Lambda, Q,t)- Z(\Lambda, Q',t')| \leq   T^{-\mathfrak{c}}.
\end{equation}
\end{lem}

\begin{proof}
  By Corollary \ref{cor:xicont} with $S \coloneqq \mathfrak{M} (\log N)^3$, there exists an event $\mathsf{E}$ of overwhelming probability, on which for all $k, k' \in \llbracket-T^{3/2}, T^{3/2} \rrbracket$ with $|k - k'| \leq \mathfrak{M} (\log N)^3$; all $\Lambda \in [-\log N, \log N]$; and all $t,t' \in  [0,  T \log N]$ satisfying $|t-t'| \leq \mathfrak{M} (\log N)^3$, we have 
  \begin{equation}\label{eqn:kcont}
  |\Xi(\Lambda, k',t') - \Xi(\Lambda, k,t)|  \leq (\log N)^C \mathfrak{M}^{1/2} T^{-1/2}.
  \end{equation}
  
  \noindent From \eqref{eqn:kcont}, we obtain on $\mathsf{E}$, for all $Q,Q'\in [- T^{4/3},  T^{4/3}]$ with $|Q - Q'| \leq  \mathfrak{M}$ and all $t,t' \in [0, T \log N]$ with $|t-t'| \leq \mathfrak{M} $, 
  \begin{equation}\label{eqn:kcont2}
\big | \Xi \big(\Lambda, (Q - t  \ve(\Lambda)) / \alpha,t \big) -   \Xi \big(\Lambda, (Q' - t'  \ve(\Lambda)) / \alpha,t' \big) \big| \leq (\log N)^C ( \mathfrak{M}^{1/4} T^{-1/4} +T^{-c}).
 \end{equation}
Above, we used the bound $|\ve(\Lambda)|\leq C \log N$ (by Lemma \ref{lem:ve_tail}) to conclude that $(Q - t  \ve(\Lambda))/ \alpha \in [-T^{3/2}, T^{3/2}]$, and that $|t-t'| |\ve(\Lambda)| \leq |\alpha| \mathfrak{M} (\log N)^2$.

Now set $f(\Lambda) =\Xi(\Lambda, (Q- t \ve(\Lambda))/\alpha,t) - \Xi(\Lambda, (Q'- t' \ve(\Lambda))/\alpha, t')$. By Lemma \ref{lem:ve_tail}, we have 
\begin{equation}\label{eqn:fdr_bd}
    |f^{\dr}(x)| \leq C (|f(x)| +  \|f\|_{\mathcal{H}} \log(|x|+2)).
\end{equation}
Using Lemma \ref{lem:xi_bd_outline} to obtain, after further restricting to an event of overwhelming probability, the bounds $|\Xi(\Lambda, (Q- t \ve(\Lambda))/{\alpha},t)| \leq (\log N)^C (|\Lambda| +1)^2$ and $|\Xi(\Lambda, (Q'- t' \ve(\Lambda))/\alpha, t')| \leq (\log N)^C (|\Lambda| +1)^2$, we can bound 
\begin{equation}\label{eqn:fhnbd}
    \|f\|_{\mathcal{H}} \leq  T^{-c} +  c^{-1} e^{-c (\log N)^2} \le 2 T^{-c}.
\end{equation}
Indeed, the first term comes from the bound \eqref{eqn:kcont2} for $|x| \leq \log N$, and the second term comes from bounding $\int_{|x|> \log N} |f(x)|^2\varrho(x) dx  $ using Lemma \ref{lem:varrho_bd} and the above tail bounds on $\Xi$.

Inserting \eqref{eqn:kcont2} and \eqref{eqn:fhnbd} into \eqref{eqn:fdr_bd} (with \eqref{eqn:Z_outline} and the third part of Lemma \ref{lem:ve_tail}) yields \eqref{eqn:Z_dif_bound} and thus the lemma.
\end{proof}

\subsection{Approximating indices}
\label{subsec:prep_conc}
In this section we show the following lemma approximating random Lax matrix functionals of the form $\sum_{i=1}^{N}F(\lambda_i) \cdot G(Q  -Q_i(t))$ with those of the form $\sum_{i=1}^{N}F(\lambda_i) \cdot G(Q_k(t)  -Q_i(t))$, for a suitable index $k$ (which is the, possibly random, one for which $Q_k (t) \approx Q$). This will more directly enable the use of certain concentration estimates (such as Lemma \ref{lem:Ndeltabd}) later.

\begin{lem}\label{lem:conc_estimate_fixedq}
   Let $A, B > 0$ and $S \in [1, T \log N]$ be real numbers. There exists a constant $\mathfrak{C}>0$ such that the following holds with overwhelming probability. Suppose $F, G : \mathbb{R} \rightarrow \mathbb{R}$ are any functions with $\sup_{|x| \leq \log N} |F(x)| \leq A$, satisfying the below conditions. 
    \begin{enumerate}
        \item For all real numbers $x \ge y$, we have \begin{equation}\label{eqn:AmolG_ass1}
        |G(x) - G(y)| \leq BS^{-1} (x-y) + B \cdot \mathbbm{1}_{x\geq 0 \geq y}.
        \end{equation}
        \item We have 
        \begin{equation}\label{eqn:AmolG_ass2}
        \supp G \subset [-S, S].
        \end{equation}
    \end{enumerate}
    
     \noindent Then, for all $Q \in [- T^{4/3},  T^{4/3}]$ and $t \in [0, T \log N]$, there exists an index $k = k(Q)$ satisfying $\varphi_t(k) \in \llbracket -T^{3/2}, T^{3/2} \rrbracket$, such that $|Q_k(t) - Q| \leq (\log N)^3$ and
    \begin{equation}\label{eqn:FGconc_q}
       \left| \sum_{i=1}^{N}F(\lambda_i) \cdot G(Q  -Q_i(t)) - \sum_{i=1}^{N}F(\lambda_i) \cdot G(Q_k(t)  -Q_i(t)) \right| \leq A B (\log N)^{\mathfrak{C}} .
    \end{equation}

\end{lem}

\begin{proof}
    Let us restrict to the event $\mathsf{E}_1$ of overwhelming probability defined in Lemma \ref{lem:q_spacing} with $T$ there given by $T \log N$ here; to the event $\mathsf{E}_2 = \bigcap_{t \geq 0} \mathsf{BND}_{\L(t)}(\log N)$ (which we may by Lemma \ref{lem:bd_lem}). Using the restriction to $\bigcap_{t \geq 0} \mathsf{BND}_{\L(t)}(\log N)$, with the fact that $q_j(t)-q_0(t) > \alpha j /2 $ for $j> T (\log N)^6$ (by Lemma \ref{lem:q_spacing}), we observe that $j > T^{3/2}$ implies that $q_j(t)-Q >(\log N)^{-1} T^{3/2} -T^{4/3} >  T \log N \geq 10 S$. Similarly, if $j < - T^{3/2}$, then $q_j(t)-Q  < - 10 S$. Using this, we claim that on $\mathsf{E} = \mathsf{E}_1 \cap \mathsf{E}_2$, there is some (random) index $k$ such that  $\varphi_t(k) \in \llbracket N_1 + T^2, N_1-T^2 \rrbracket$, and satisfying $|Q_k(t) - Q| \leq (\log N)^3$. Indeed, if this were not the case, then for some $k \in \llbracket N_1 + T^2, N_2 - T^2 \rrbracket$ we would have $|q_k(t) - q_{k+1}(t)| >(\log N)^2  $ which is impossible on $\mathsf{E}$ (by Lemma \ref{lem:q_spacing}). In addition, by our previous observation, we know that any such $k = k(Q)$ satisfies $\varphi_t(k) \in \llbracket -T^{3/2}, T^{3/2} \rrbracket$.

    As a consequence, on this event, choosing this random index $k = k(Q)$,
\begin{flalign}\label{eqn:q_k_diff}
\begin{aligned}
     &  \left| \sum_{i=N_1}^{N_2}F(\Lambda_i(t)) G(Q  -q_i(t)) - \sum_{i=N_1}^{N_2}F(\Lambda_i(t)) G(q_k(t)  -q_i(t)) \right|  \\
       &\leq \sum_{i=N_1}^{N_2}\left|F(\Lambda_i(t))\right| \cdot \left|G(Q  -q_i(t)) -  G(q_k(t)  -q_i(t)) \right| \mathbbm{1}_{|i-k| < 2 S}\\
     &\leq  5 A B (\log N)^3 +  2 A B \sum_{i=N_1}^{N_2}  \left( \mathbbm{1}_{q_k(t) - q_i(t) \ge 0 \ge Q - q_i(t)}
       +\mathbbm{1}_{Q - q_i(t) \ge 0 \ge q_k(t) - q_i(t)} \right) 
       \leq A B (\log N)^C .
    \end{aligned}
    \end{flalign}
To obtain \eqref{eqn:q_k_diff}, we have used the assumptions \eqref{eqn:AmolG_ass1} and \eqref{eqn:AmolG_ass2} on $G$, together with the following fact: On $\mathsf{E}$ and with $k$ as above, there at most $C (\log N)^5$ indices $i$ for which
$q_k(t) - q_i(t) \ge 0 \ge Q - q_i(t)$ or $q - q_i(t) \ge 0 \ge q_k(t) - q_i(t)$. Indeed, if on the contrary, say, $q_k(t) \geq q_{i}(t) \geq q$ for both $i = i_0$ and $i = i_1$, and moreover $|i_1-i_0| > (\log N)^5$, then $|q_{i_1}(t) - q_{i_0}(t)| \leq  (\log N)^3$. However, this contracts that (for $N$ large enough) we must have $|q_{i_1}(t) - q_{i_0}(t)| > c (\log N)^5 >  (\log N)^3$ on $\mathsf{E}$. In addition, we have used the restriction to $\bigcap_{t \geq 0} \mathsf{BND}_{\L(t)}(\log N)$ in order to bound $|F(\Lambda_i(t))| \leq A$.
\end{proof}

\subsection{Uniform concentration estimate}
\label{subsec:key_conc}

 In this section we show the following lemma, which is formally obtained by replacing $G(Q_i(s) - Q_j(s))$ in Lemma \ref{lem:Ndeltabd} by $\chi'(Q_j(t) - Q_i(t) )$ and $F(\lambda)$ there by $F(\lambda) \cdot f(\Lambda - \lambda) \cdot Z(\lambda, Q)$ here (and then using Lemma \ref{lem:conc_estimate_fixedq}). However, Lemma \ref{lem:Ndeltabd} is not directly applicable to such a choice of $F$, because $Z(\lambda, Q)$ is a random function, as it depends on the Lax matrix $\L(0)$ through the family of observables $\Xi(\Lambda, k, t)$ (which in turn depend on the full spectrum of $\mathbf{L}(0)$).

\begin{lem}\label{lem:Z_sum_conc}

    There exists a constant $\mathfrak{c}>0$ such that the following holds. Let $t \in [0, T \log N]$ and $A \in [1, N]$ be real numbers; let $f : \mathbb{R} \rightarrow \mathbb{R}$ be a function either equal to $f \equiv 1$ or $f = \mathfrak{l}$; and let $F : \mathbb{R} \rightarrow \mathbb{R}$ be a function satisfying Assumption \ref{ass:F}. Further assume that, for all real numbers $\lambda, \lambda' \in [-\log N, \log N]$ satisfying $|\lambda-\lambda'| \leq e^{-(\log N)^{1/2}}$, the function $F$ satisfies 
    \begin{equation}\label{eqn:Fhold}
        |F(\lambda) - F(\lambda')| \leq A |\lambda-\lambda'|^{1/4}.
    \end{equation}
    
    \noindent  Then, with overwhelming probability, we have (recalling $Z(\Lambda, Q) = Z(\Lambda, Q, t)$ from \eqref{eqn:Z_outline}), 
   \begin{multline}\label{eqn:main_approx_intermed}
      \displaystyle\sup_{|\Lambda| \le \log N} \displaystyle\sup_{|Q| \le T^{4/3}} \Bigg|  \sum_{i=1}^{N} F(\lambda_i) \cdot f(\Lambda - \lambda_i) \cdot Z(\lambda_i, Q) \cdot  \chi'(Q - Q_i(t) ) \\
      -\alpha^{-1} \int_{-\infty}^{\infty}F(\lambda) \cdot f(\Lambda - \lambda) \cdot Z(\lambda,Q) \cdot \varrho(\lambda) d\lambda \Bigg| \leq  A T^{-\mathfrak{c}},
    \end{multline}
    and 
    \begin{multline}\label{eqn:main_approx_intermed2}
      \displaystyle\sup_{|\Lambda| \le \log N} \displaystyle\sup_{|Q| \le T^{4/3}} \Bigg|  \sum_{i=1}^{N} F(\lambda_i) \cdot f(\Lambda - \lambda_i)  \cdot  \chi'(Q - Q_i(t) ) \\
      -\alpha^{-1} \int_{-\infty}^{\infty}F(\lambda) \cdot f(\Lambda - \lambda)  \cdot \varrho(\lambda) d\lambda \Bigg| \leq  A T^{-\mathfrak{c}}. 
    \end{multline}
\end{lem}

\begin{proof}

We only prove \eqref{eqn:main_approx_intermed}, as the proof of \eqref{eqn:main_approx_intermed2} is very similar. Moreover, we will assume that $f = \mathfrak{l}$, as the proof is entirely analogous if instead $f = 1$. We will proceed by defining a net of functions, fine enough so that the function $\lambda \mapsto F(\lambda) \mathfrak{l}(\Lambda-\lambda) \cdot Z(\lambda, Q)$ can be approximated by some function in this net, but small enough so that a union bound applies to ensure that Lemma \ref{lem:Ndeltabd} likely holds simultaneously for each function in this net. 

To begin, let $\delta > 0$ be a sufficiently small constant (namely, smaller than $\mathfrak{c}/50$, where $\mathfrak{c}$ is the constant in Lemma \ref{lem:Z_cont_lambda}).  Let $\gamma= \delta/10$. Define $\mathcal{L}_{N;  \gamma} = \mathcal{L}_{N} \coloneqq [-\log N, \log N] \cap (T^{-\gamma} \cdot \mathbb{Z})$, and let $x_{\min} = \min \mathcal{L}_N$ and $x_{\max} = \max \mathcal{L}_N$. We call $x,y \in \mathcal{L}_N$ \emph{adjacent} if $|x-y| = T^{-\gamma}$.

 Let us rename $Q$ in the second argument of $Z$ to $Q_0$, and define, for $Q_0 \in \mathbb{R}$,
 \begin{equation}\label{eqn:ZQ0}
 Z_{Q_0}(\lambda) \coloneqq Z(\lambda, Q_0,t).
 \end{equation}
 By Lemma \ref{lem:Z_cont_lambda}, we have with overwhelming probability that, for any $Q_0 \in [-T^{4/3}, T^{4/3}]$ the function $\lambda \mapsto Z_{Q_0}(\lambda) $ satisfies, for $x,y \in [-\log N, \log N]$ with $|x-y| \leq T^{-\gamma}$, the estimate
 \begin{equation}\label{eqn:ZDB}
     |Z_{Q_0}(x) -  Z_{Q_0}(y)| \leq (\log N)^C T^{-\gamma/4}.
 \end{equation}

 \noindent Moreover, by Lemma \ref{lem:Z_bd}, we have with overwhelming probability that 
\begin{equation}\label{eqn:ZBD2}
\sup_{|Q| \leq T^{4/3}} |Z(\Lambda, Q,t)| \leq (\log N)^C (1+|\Lambda|)^2, \qquad \text{for all } \Lambda \in \mathbb{R}
\end{equation}

\noindent Let $\mathsf{E}$ denote the event on which Lemma \ref{lem:q_spacing}, $\mathsf{BND}_{\mathbf{L}} (\log N)$ (recall Definition \ref{eigenvaluesm} and Lemma \ref{lem:bd_lem}), Lemma \ref{lem:conc_estimate_fixedq}, \eqref{eqn:ZDB}, and \eqref{eqn:ZBD2} all hold. We restrict to $\mathsf{E}$ in what follows. 

Let us define a piecewise linear function $H_{Q_0}: \mathbb{R} \rightarrow \mathbb{R}$ as follows. For $x \in \mathcal{L}_N$, set
 \begin{equation}\label{eqn:barZQ0}
     H_{Q_0}(x) \coloneqq \floor{T^{\gamma} F(x) Z_{Q_0}(x)} T^{-\gamma} .
 \end{equation}
Then, for $x \in [x_{\min}, x_{\max}] \setminus \mathcal{L}_N$, we extend $H$ linearly between adjacent lattice points. Finally, for $x \leq x_{\min}$ or $x \geq x_{\max}$, we extend $H_{Q_0}(x) = H_{Q_0}(x_{\min})$ or $H_{Q_0}(x) = H_{Q_0}(x_{\max})$ by a constant out to $-\infty$ or $+\infty$, respectively. Let us show several bounds for $H_{Q_0}$. 

If $x,y \in \mathcal{L}_N$ are adjacent lattice points, then (by \eqref{eqn:Fhold}, \eqref{eqn:ZDB}, \eqref{estimatef}, and \eqref{eqn:ZBD2})
\begin{flalign}\label{eqn:FZdfbd}
\begin{aligned}
|H_{Q_0}(x) - H_{Q_0}(y)|  
&\leq |F(x) Z_{Q_0}(x) - F(y) Z_{Q_0}(y)| + 2 T^{-\gamma} \\
     &\leq  |F(x)-F(y)| \cdot |Z_{Q_0}(x)| + |F(y)| \cdot |Z_{Q_0}(x)- Z_{Q_0}(y)|  + 2 T^{-\gamma}\\
     &\leq A (\log N)^C T^{-\gamma/4} .
\end{aligned}
\end{flalign}
    Moreover, for all $x \in [- \log N, \log N]$, 
    \begin{equation}\label{eqn:ZFbd}
 |H_{Q_0}(x)| \leq |Z_{Q_0}(x) F(x) | + 2 T^{-\gamma} \leq A (\log N)^C.
    \end{equation}

    \noindent Define $C_0$ to be equal to the constant $C$ for which \eqref{eqn:FZdfbd} and \eqref{eqn:ZFbd} both hold. Now we further shrink $\mathsf{E}$ so that the following events hold on $\mathsf{E}$: The outcome of Lemma \ref{lem:conc_estimate_fixedq} with $(A, B, S)= (A (\log N)^{C_0}, 10, \mathfrak{M})$; the event of Lemma \ref{lem:q_spacing} with $T $ there given by $T \log N$ here; and the event $\mathsf{BND}_{\L(0)}(\log N)$.  
    
  Next, define the collection of functions $\mathcal{F}_N$ as follows. A function $F_0  : \mathbb{R} \rightarrow \mathbb{R}$ is in $\mathcal{F}_N$ if and only if it satisfies the following properties:
    \begin{enumerate}
        \item For all $x \in \mathcal{L}_{N}$, we have $F_0(x) \in [-A (\log N)^{C_0}, A (\log N)^{C_0}] \cap T^{-\gamma}\cdot  \mathbb{Z}$.
        \item For all adjacent $x,y \in \mathcal{L}_N$, we have $|F_0(x)-F_0(y)| \leq A (\log N)^{C_0} |x-y|^{1/4} $.
        \item For all adjacent $x, y \in \mathcal{L}_N$ with $x < y$, and $z \in [x,y]$, we have $F_0(z)= t F_0(y) + (1-t) F_0(x)$.
        \item For all $z \leq x_{\min}$, we have $F_0(z) = F_0(x_{\min})$ and for all $z \geq x_{\max}$, we have $F_0(z) = F_0(x_{\max})$.
    \end{enumerate}

\noindent By the definition of $C_0$, \eqref{eqn:barZQ0} and \eqref{eqn:ZFbd} imply the first item above, and \eqref{eqn:FZdfbd} implies the second, for $F_0 = H_{Q_0}$. The third and fourth items above hold by the definition of $H_{Q_0}$. Thus,
\begin{equation}\label{eqn:H0inFN}
    H_{Q_0} \in \mathcal{F}_N, \qquad \text{ for any } Q_0 \in [ -T^{4/3}, T^{4/3} ].
\end{equation}

  Consider the finite collection of functions $\mathcal{F}_N$ defined above, and denote its cardinality by $|\mathcal{F}_N|$. Note that 
  \begin{equation}\label{eqn:FNCbd}
  |\mathcal{F}_N| \leq 2 A (\log N)^{C_0}  T^{\gamma} \cdot e^{T^{\gamma} (\log N)^2}  \leq A  e^{T^{\gamma} (\log N)^3} .
  \end{equation}
  Here we have used that the number of points in $\mathcal{L}_N$ is $\leq T^{\gamma} \cdot 2 \log N$; the number of possible values of $F(x_{\min})$ is $2  A (\log N)^{C_0} T^{\gamma} $; and, if we choose values of $F$ sequentially at lattice points $x_i \in \mathcal{L}_N$ starting from $x_{\min}$, number of possible choices for each subsequent $F(x_i)$, given the value $F(x_{i-1})$ at the adjacent point $x_{i-1}$ to the left of $x_i$, is upper bounded by $T^{\gamma}$, by the H\"older condition. 

    Let us define another mesh $\mathcal{L}_{\eig} \coloneqq [- \log N, \log N] \cap (e^{-(\log N)^3} \cdot \mathbb{Z})$. Note 
    \begin{equation}\label{eqn:FNLCbd}
        |\mathcal{F}_N\times \mathcal{L}_{\eig}| \leq A  e^{T^{\gamma} (\log N)^3} \times e^{2(\log N)^3} \leq A  e^{T^{\gamma} (\log N)^4},
    \end{equation}
    by \eqref{eqn:FNCbd} and the fact that there are at most $e^{2 (\log N)^3}$ choices of $\Lambda \in \mathcal{L}_{\eig}$. 

  For a function $F$ and an index $k \in \llbracket 1, N \rrbracket$, denote the (random) functional $\Sigma_k(F, \Lambda)$ by
  $$
  \Sigma_k(F, \Lambda) \coloneqq  \sum_{i=1}^{N} F(\lambda_i) \mathfrak{l}(\Lambda-\lambda_i) \chi'(Q_k(t) - Q_i(t) ) -\alpha^{-1} \int_{-\infty}^{\infty}F(\lambda) \mathfrak{l}(\Lambda-\lambda_i) \varrho(\lambda) d\lambda  .
  $$
  By Lemma \ref{lem:Ndeltabd} with $(F, G; A, B; \mathfrak{K}) = (F_0\cdot \mathfrak{l}(\Lambda-\cdot),\chi'; A (\log N)^{C_0+3}, 10;T^{\delta})$, and a union bound (with \eqref{eqn:FNLCbd}), for some constant $\mathfrak{c}_1>0$
  \begin{multline}\label{eqn:FN_prob_bd}
      \mathbb{P}\left[ \exists F_0 \in \mathcal{F}_N, \Lambda \in \mathcal{L}_{\eig}, k \in \varphi_t^{-1}\llbracket -T^2, T^2 \rrbracket  : \left|\Sigma_k(F_0,\Lambda)\right| >  10 A (\log N)^{C_0+3} T^{13 \delta} \mathfrak{M}^{-1/2} \right] \\
      \leq 
     |\mathcal{F}_N \times \mathcal{L}_{\eig}| \cdot (2 T^2) \cdot \mathfrak{c}_0^{-1}e^{-\mathfrak{c}_0 T^{2\delta}} \leq  \mathfrak{c}_1^{-1} e^{-\mathfrak{c}_1 T^{2\delta}} .
  \end{multline}

Define, for any real numbers $Q, \Lambda \in \mathbb{R}$ and function $H : \mathbb{R} \rightarrow \mathbb{R}$,
\begin{equation}\label{eqn:SigmaQ}
     \Sigma(H, \Lambda, Q) \coloneqq  \sum_{i=1}^{N} H(\lambda_i) \mathfrak{l}(\Lambda-\lambda_i) \chi'(Q - Q_i(t) ) -\alpha^{-1} \int_{-\infty}^{\infty}
     H(\lambda) \mathfrak{l}(\Lambda-\lambda_i) \varrho(\lambda) d\lambda .
\end{equation}

Note that by Lemma \ref{lem:conc_estimate_fixedq} (which recall holds on $\mathsf{E}$) with $(F, G) = (H_{Q_0}, \chi')$, there exists some $k_0 \in \llbracket 1, N \rrbracket$ satisfying $\varphi_t(k_0) \in \llbracket -T^{-3/2}, T^{-3/2} \rrbracket$, $|Q_{k_0}(t)-Q| \leq (\log N)^3$, and
\begin{equation}\label{eqn:sigmaQsigmak}
 \left| \Sigma_{k_0}(H_{Q_{0}}, \Lambda) - \Sigma(H_{Q_{0}}, \Lambda, Q) \right|  \leq A (\log N)^{C_0+C} \mathfrak{M}^{-1}, \qquad \text{for all } Q_0 \in [-T^{4/3}, T^{4/3}] .
 \end{equation}
By \eqref{eqn:FN_prob_bd}, \eqref{eqn:H0inFN}, and \eqref{eqn:sigmaQsigmak} with $Q_0 = Q$, we have
 \begin{equation}\label{eqn:FNQ_prob_bd}
      \mathbb{P}\left[ \exists  \Lambda \in \mathcal{L}_{\eig}, Q \in [ -T^{4/3}, T^{4/3}] : \big|\Sigma(H_{Q},\Lambda,Q)\big| > A (\log N)^{C_0+C} T^{13 \delta} \mathfrak{M}^{-1/2} \right] \leq    c^{-1} e^{-c (\log N)^2} .
  \end{equation}

Note that, by following the reasoning in \eqref{eqn:FZdfbd}, we have
  \begin{equation}\label{eqn:ZFZFbbd}
     \sup_{x \in [- \log N, \log N]} |H_{Q_0}(x) -  Z_{Q_0}(x) F(x)| \leq (\log N)^C A T^{-\gamma/4}  .
 \end{equation}
In addition, since the outcome of Lemma \ref{lem:q_spacing} and the event $\mathsf{BND}_{\L(0)}(\log N)$ both hold on $\mathsf{E}$, using the fact that $|Q_{k_0}(t)-Q| \leq (\log N)^3$, we obtain
\begin{flalign}\label{eqn:sig_lhsdif}
 \begin{aligned}
    \Big|  \sum_{i=1}^{N} & \mathfrak{l}(\Lambda - \lambda_i) \cdot F(\lambda_i) Z_{Q}(\lambda_i) \cdot \chi'(Q - Q_i(t) )
   -   \sum_{i=1}^{N}  \mathfrak{l}(\Lambda - \lambda_i) \cdot H_{Q}(\lambda_i) \cdot  \chi'(Q - Q_i(t) ) \Big| \\
   &\leq (\log N)^3\sum_{i=1}^{N} \left|F(\lambda_i) Z_{Q}(\lambda_i) - H_{Q}(\lambda_i) \right| \chi'(Q - Q_i(t) ) \mathbbm{1}_{|\varphi_t(i)-\varphi_t(k_0)| \leq \mathfrak{M} (\log N)^5} \\
   & \leq A (\log N)^C T^{-\gamma/4}.
 \end{aligned}
 \end{flalign}
To obtain the first inequality, we used the fact that the summand in the second line vanishes unless $|\varphi_t(i)-\varphi_t(k_0)| \leq \mathfrak{M} (\log N)^5$ (by Lemma \ref{lem:q_spacing}). To obtain the final inequality we have used \eqref{eqn:ZFZFbbd}. Moreover, 
\begin{flalign}\label{eqn:sig_rhsdif}
 \begin{aligned}
   \Big|  \int_{-\infty}^{\infty} & \mathfrak{l}(\Lambda-\lambda) \cdot F(\lambda) Z_{Q}(\lambda) \cdot \varrho(\lambda) d\lambda 
    -  \int_{-\infty}^{\infty} \mathfrak{l}(\Lambda-\lambda) \cdot H_{Q}(\lambda) \cdot \varrho(\lambda) d\lambda \Big|  \\ 
 &\leq  2  A (\log N)^C T^{-\gamma/4} \int_{-\log N}^{\log N} \varrho(\lambda) d\lambda \\
    & \qquad +  (\log N)^C \int_{[-\log N, \log N]^{\complement}}  A e^{2 |\lambda|^{1/2}}  (1+|\lambda|)^2  |\mathfrak{l}(\Lambda-\lambda)|  \varrho(\lambda) d\lambda  \\
    &\leq   A (\log N)^C T^{-\gamma/4} +c^{-1} A e^{-c(\log N)^2 }.
 \end{aligned}
 \end{flalign}
 The first term in the middle display is due to the fact that $\sup_{|x|\leq 2 \log N} |\mathfrak{l}(x)| \leq (\log N)^3$, which implies $\sup_{|\lambda| \leq \log N} |F(\lambda)| |\mathfrak{l}(\Lambda-\lambda)| \leq A (\log N)^3$, together with \eqref{eqn:ZFZFbbd}. The second term is from the fact that $|F(\lambda)| |\mathfrak{l}(\Lambda-\lambda)| \leq A (\log N)^3 e^{2 |\lambda|^{1/2}}$ (by \eqref{estimatef}); and the bound \eqref{eqn:ZBD2} which holds on $\mathsf{E}$. The third bound holds since $\varrho(\lambda) \leq c^{-1} e^{-c \lambda^2}$ (by Lemma \ref{lem:varrho_bd}).

 As a consequence of \eqref{eqn:sig_lhsdif} and \eqref{eqn:sig_rhsdif}, on $\mathsf{E}$, for all $Q \in [-T^{4/3}, T^{4/3}]$ and $\Lambda \in \mathcal{L}_{\eig}$,
 \begin{equation}\label{eqn:sigma_actual_sigmaFN}
 \Big|  \Sigma(F \cdot Z_{Q}, \Lambda,Q) - \Sigma(H_Q, \Lambda,Q) \Big| \leq  A (\log N)^C T^{-\gamma/4} .
 \end{equation}
By \eqref{eqn:sigma_actual_sigmaFN} and \eqref{eqn:FNQ_prob_bd}, we have 
  \begin{equation}
  \label{eventlambda0}
      \mathbb{P}\left[ \exists \Lambda \in \mathcal{L}_{\eig}, Q \in [ -T^{4/3}, T^{4/3}] : \left|\Sigma(F \cdot Z_{Q},\Lambda,Q)\right| > A T^{-\gamma/5} \right] 
      \leq c^{-1}e^{-c (\log N)^2} .
 \end{equation}
Now, because the spacing of the mesh $\mathcal{L}_{\eig}$ is $e^{-(\log N)^3}$, a quick continuity argument shows that, on the event $\tilde{\mathsf{E}}$ defined as the complement of the event whose probability is upper bounded in \eqref{eventlambda0}, the desired bound \eqref{eqn:main_approx_intermed} holds for all $\Lambda \in [-\log N, \log N]$. This completes the proof.
\end{proof}

\subsection{Proof of Proposition \ref{prop:sumZconc}} 
\label{sec:Zsumtoint}

Before proving Proposition \ref{prop:sumZconc}, we require a lemma that modifies Lemma \ref{lem:Z_sum_conc} by replacing $Q_i(t)$ there by its limiting trajectory $Q_i(0) + t \ve(\lambda_i)$ (or $ \alpha \varphi_0(i) + t \ve(\lambda_i)$). Again, recall $Z(\Lambda, Q) = Z(\Lambda, Q, t)$ from \eqref{eqn:Z_outline}.

\begin{lem}\label{lem:z_sumconc_withF}
 Let $A \in [1, N]$ and $t \in [0, T \log N]$ be real numbers, let $F, f : \mathbb{R} \rightarrow \mathbb{R}$ be two functions satisfying the assumptions of Lemma \ref{lem:Z_sum_conc}. There exist a constant $\mathfrak{c}> 0$ such that the following holds with overwhelming probability. For all $\Lambda \in [- \log N, \log N]$ and $Q \in [- T^{6/5},  T^{6/5}]$, we have 
\begin{align}
   &  \Bigg| \sum_{i=N_1}^{N_2} Z(\Lambda_{i}, \alpha i + t  \ve(\Lambda_i)) \cdot   
    F(\Lambda_i) f(\Lambda -\Lambda_i) \cdot \chi'(Q -\alpha i   - t \ve(\Lambda_{i}) )  \notag \\
   & \qquad \qquad \qquad \qquad  -\alpha^{-1} \int_{-\infty}^{\infty}F(\lambda)  f(\Lambda - \lambda) \cdot Z(\lambda,Q) \cdot \varrho(\lambda) d\lambda \Bigg| \leq  A T^{-\mathfrak{c}}; \label{eqn:Zeq_2_Ff} \\
    &  \Bigg| Z(\Lambda, Q) \sum_{i=N_1}^{N_2}     
    F(\Lambda_i) f(\Lambda -\Lambda_i) \cdot \chi'(Q -\alpha i   - t \ve(\Lambda_{i}) )  \notag \\
   & \qquad \qquad \qquad  \qquad  -\alpha^{-1} Z(\Lambda, Q)  \int_{-\infty}^{\infty}F(\lambda)  f(\Lambda - \lambda)  \cdot \varrho(\lambda) d\lambda \Bigg| \leq  A T^{-\mathfrak{c}}; \label{eqn:Zeq_2_Ff2} \\
    &  \Bigg| \sum_{i=N_1}^{N_2} Z(\Lambda_{i}, \alpha i + t  \ve(\Lambda_i)) \cdot   
    F(\Lambda_i) f(\Lambda -\Lambda_i) \cdot \chi'(Q -q_i(0)   - t \ve(\Lambda_{i}) )  \notag \\
   & \qquad \qquad \qquad \qquad  -\alpha^{-1} \int_{-\infty}^{\infty}F(\lambda)  f(\Lambda - \lambda) \cdot Z(\lambda,Q) \cdot \varrho(\lambda) d\lambda \Bigg| \leq  A T^{-\mathfrak{c}}. \label{eqn:Zeq_2_Ff3}
    \end{align}
    
\end{lem}

\begin{proof}
    We only show \eqref{eqn:Zeq_2_Ff} and \eqref{eqn:Zeq_2_Ff3}, since the proof of \eqref{eqn:Zeq_2_Ff2} is very similar. 

First observe that Lemma \ref{lem:Z_sum_conc} yields an event $\mathsf{E}_1$ of overwhelming probability, on which for all $\Lambda \in [-\log N, \log N]$ and $Q \in [- T^{6/5}, T^{6/5}]$ we have 
\begin{multline}\label{eqn:Zeq_approx}
  \Big|  
   \sum_{i=1}^{N} Z(\lambda_i, Q, t) 
  \cdot F(\lambda_i) f(\Lambda -\lambda_i) \cdot \chi'(Q-  Q_i(t)) )\\
-\alpha^{-1} \int_{-\infty}^{\infty} F(\lambda) f(\Lambda - \lambda) \cdot Z(\lambda,Q) \cdot \varrho(\lambda) d\lambda \Big| \leq A T^{-c}.
\end{multline}

\noindent In what follows, we restrict to the event $\mathsf{E}_1$ on which \eqref{eqn:Zeq_approx} holds. We must then approximate the sums in \eqref{eqn:Zeq_2_Ff} and \eqref{eqn:Zeq_2_Ff3} by that in \eqref{eqn:Zeq_approx} on $\mathsf{E}_1 \cap \mathsf{E}_2$, for some overwhelmingly probable event $\mathsf{E}_2$.

To do so, we restrict to the event $\mathsf{E}_{2,1}$ on which the conclusion of Lemma \ref{lem:Z_cont} holds; the event $\mathsf{E}_{2,2}$ on which the conclusion in the second part of Lemma \ref{lem:asymptotic-scattering} holds and the event given in Lemma \ref{lem:loc_cent_dif} holds; the event $\mathsf{E}_{2,3}$ on which both the conclusion of Lemma \ref{lem:bd_num_termsU} with $U = \mathfrak{M}$, the conclusion of Lemma \ref{lem:second_der_bd} holds, and the conclusion of Lemma \ref{lem:q_spacing} with $T=0$ holds; the event $\mathsf{E}_{2,4}$ on which the conclusion \eqref{eqn:small_Lambda_Zbd} of Lemma \ref{lem:Z_bd} holds; the event $\mathsf{E}_{2,5} \coloneqq \bigcap_{i=1}^N \{|\lambda_i| \leq \log N \}$, which is overwhelmingly probable by Lemma \ref{lem:bd_lem}. We restrict to the event $\mathsf{E}_2 \coloneqq \bigcap_{i=1}^5 \mathsf{E}_{2, i}$ in what follows, and let $\Lambda \in [- \log N, \log N]$ and $Q \in [- T^{6/5}, T^{6/5}]$ be real numbers.

Define $i_t \coloneqq \varphi_t(\varphi_0^{-1}(i))$ for any $i \in \llbracket N_1, N_2 \rrbracket$. Moreover, let $\mathfrak{S} \coloneqq \{i \in \llbracket N_1, N_2 \rrbracket : |\alpha i + t  \ve(\Lambda_i) - Q| \leq 2 \mathfrak{M} \}$, and $k \coloneqq \min \mathfrak{S}$, that is, let $k$ be the minimal index such that $|\alpha k + t  \ve(\Lambda_k) - Q| \leq 2 \mathfrak{M}$.  By Lemma \ref{lem:ve_tail}, since $Q \in [-T^{6/5}, T^{6/5}]$, we have $k \in \llbracket - T^{6/5} (\log N)^{10}, T^{6/5} (\log N)^{10} \rrbracket$ on $\mathsf{E}_{2,5}$. With this notation, on $\mathsf{E}_{2,1} \cap \mathsf{E}_{2,2}\cap \mathsf{E}_{2,3}$, we have 
\begin{flalign}\label{eqn:ai_Q_and_conc}
\begin{aligned}
 & \mathbbm{1}_{|\alpha i + t  \ve(\Lambda_i) - Q| \leq \mathfrak{M}}   | Z(\Lambda_i, \alpha i + t  \ve(\Lambda_i), t) -  Z(\Lambda_i, Q, t)| \leq  T^{-c}; \\
&  |\alpha i + t  \ve(\Lambda_i) - Q| > \mathfrak{M} \quad \text{holds whenever} \quad |i_t - k_t| \geq  \mathfrak{M} (\log N)^5 .
\end{aligned}
\end{flalign}
Thus, on $\mathsf{E}_{2,1} \cap \mathsf{E}_{2,2} \cap \mathsf{E}_{2,3} \cap \mathsf{E}_{2,5}$, using \eqref{eqn:ai_Q_and_conc}, the fact that $|F(\Lambda_i) f(\Lambda - \Lambda_i)| \le A (\log N)^C$, and the fact that $|\chi'(x)| \leq 10\mathfrak{M}^{-1}$ for all $x$, we have 
\begin{multline}\label{eqn:Zeq_approx01}
  2 \sum_{i=N_1}^{N_2}  \big|  Z(\Lambda_i, \alpha i + t  \ve(\Lambda_i), t) -  Z(\Lambda_i, Q, t)\big| \cdot 
 |  F(\Lambda_i) f(\Lambda-\Lambda_i) | \cdot \chi'(Q- \alpha i - t   
 \ve(\Lambda_{i} )) \leq A T^{-c},
\end{multline}

\noindent since there are at most $3 \mathfrak{M} (\log N)^5 $ nonzero terms in the sum. By similar reasoning,
\begin{multline}\label{eqn:Zeq_approx02}
  2 \sum_{i=N_1}^{N_2}  \big|  Z(\Lambda_i, \alpha i + t  \ve(\Lambda_i), t) -  Z(\Lambda_i, Q, t)\big| \cdot
 |  F(\Lambda_i) f(\Lambda-\Lambda_i) | \cdot \chi'(Q- q_i(0) - t   
 \ve(\Lambda_{i} )) \leq A T^{-c}.
\end{multline}

Next, we have the following on $\mathsf{E}_{2}$. First, by Lemma \ref{lem:loc_cent_dif} and Lemma \ref{lem:q_spacing}, we have for all $i \in \llbracket N_1, N_2 \rrbracket$ that 
$$
\mathbbm{1}_{|i_t-k_t| \leq  \mathfrak{M} (\log N)^5} \cdot | \alpha i -q_i(0) | \leq |i|^{1/2} (\log N)^2 \le T^{3/5}  (\log N)^C.
$$
Using this; the mean value theorem to rewrite $\chi'(Q- \alpha i - t   
 \ve(\Lambda_{i} )) - \chi'(Q-  q_i(0) - t   
 \ve(\Lambda_{i})) = (q_i(0) - \alpha i) \chi''(\xi_i ) $ for some $\xi_i$ between $Q- \alpha i - t   
 \ve(\Lambda_{i})$ and $Q-  q_i(0) - t   
 \ve(\Lambda_{i})$; and the fact that $\chi''(\xi_i) \neq 0$ implies that $|i_t-k_t| \leq \mathfrak{M} (\log N)^5$, we obtain 
\begin{flalign}\label{eqn:chidifs1}
\begin{aligned}
   \sum_{i=N_1}^{N_2} & | Z(\Lambda_i, Q, t)| \cdot  \big|
  F(\Lambda_i) f(\Lambda-\Lambda_i) \big| \cdot  \big|  \chi'(Q- \alpha i - t   
 \ve(\Lambda_{i} )) - \chi'(Q-  q_i(0) - t   
 \ve(\Lambda_{i}) )\big| \\
 &\leq  A (\log N)^C \sum_{i=N_1}^{N_2}   \mathbbm{1}_{|k_t - i_t| \leq  \mathfrak{M} (\log N)^5} \cdot  | \alpha i -q_i(0)|  \cdot |  \chi''(\xi_i )  | \leq A  (\log N)^C T^{3/5} \mathfrak{M}^{-1} .
\end{aligned}
\end{flalign}
In the second inequality, we bounded each $|Z(\Lambda_i, Q, t)| \leq (\log N)^C$ by our restriction to $\mathsf{E}_{2,4} \cap \mathsf{E}_{2,5}$; bounded each  $|\mathfrak{l}(\Lambda-\Lambda_i)| \leq (\log N)^3$; and bounded $\sup_{|\lambda| \leq \log N } |F(\lambda)| \leq A$. In the last inequality we used the fact from Definition \ref{def:chi_def} that $|\chi''| \le C\mathfrak{M}^{-2}$.

Note that by \eqref{eqn:chidifs1}, \eqref{eqn:Zeq_approx02}, and \eqref{eqn:Zeq_approx01}, to complete the proof of \eqref{eqn:Zeq_2_Ff} it now suffices to show \eqref{eqn:Zeq_2_Ff3}. To do the latter, it suffices by \eqref{eqn:Zeq_approx} and \eqref{eqn:Zeq_approx02} to show 
\begin{multline}\label{eqn:Zeq_approx3}
  \Big|  
   \sum_{i=N_1}^{N_2}  Z(\Lambda_i, Q, t) \cdot
 F(\Lambda_i) f(\Lambda-\Lambda_i ) \cdot\chi'(Q- q_i(0) - t   
 \ve(\Lambda_{i} ) ) \\
-  \sum_{i=1}^{N}  Z(\lambda_i, Q, t)  \cdot F(\lambda_i)
  f(\Lambda-\lambda_i)  \cdot \chi'(Q- Q_i(t) )\Big| \leq T^{-c}.
\end{multline}
By \eqref{eqn:Z_apriori} (which holds on $\mathsf{E}_{2,2}$) and the fact that $k \in \llbracket - T^{6/5} (\log N)^{10}, T^{6/5} (\log N)^{10} \rrbracket$, we have
$$
\mathbbm{1}_{|i_t-k_t| \leq  \mathfrak{M} (\log N)^5} \cdot \big| q_i(0) + t \ve(\Lambda_i) - Q_{\varphi_0^{-1}(i)}(t) \big| \leq T^{1/2} (\log N)^C .
$$

\noindent Thus, by entirely analogous reasoning as used to deduce \eqref{eqn:chidifs1}, we obtain 
 \begin{flalign}\label{eqn:chidifs2}
\begin{aligned}
   \sum_{i=N_1}^{N_2}  | Z( & \Lambda_i, Q, t)| \cdot \big|  F(\Lambda_i) f(\Lambda-\Lambda_i) \big| \cdot  \big|  \chi'(Q-q_i(0) - t   
 \ve(\Lambda_{i}) ) - \chi'(Q- Q_{\varphi_0^{-1}(i)}(t)) \big| \\
 &\leq A  (\log N)^C T^{1/2} \mathfrak{M}^{-1} .
\end{aligned}
\end{flalign}

\noindent The bound \eqref{eqn:chidifs2} implies \eqref{eqn:Zeq_approx3} (upon changing the summation index $i$ to $\varphi_0^{-1}(i)$ in the second sum there); thus, \eqref{eqn:Zeq_2_Ff3} holds on $\mathsf{E}_1 \cap \mathsf{E}_2$. As explained above \eqref{eqn:Zeq_approx3}, this yields \eqref{eqn:Zeq_2_Ff} and \eqref{eqn:Zeq_2_Ff3}.
\end{proof}

\begin{proof}[Proof of Proposition \ref{prop:sumZconc}]

 The bound \eqref{eqn:Zeq_2_prime_outline} follows from subtracting the two estimates \eqref{eqn:Zeq_2_Ff} and \eqref{eqn:Zeq_2_Ff2} from Lemma \ref{lem:z_sumconc_withF}, applied at  $F \equiv 1 = A$, and $f = \mathfrak{l}$. This establishes the proposition. 
\end{proof}

\section{Comparisons and limits for quasi-particle fluctuations}
\label{sec:quasi_conc_estimates}
Throughout this subsection we adopt Assumption \ref{ass:NT_assumption}. In addition, we use the notation $\Xi(\Lambda, k, t)$ as defined in Definition \ref{xilambdakt}, $Z(\Lambda, Q, t)$ as defined in \eqref{eqn:Z_outline}, and the notation from Definition \ref{def:quasi_intro}. 

The goal of this section is to prove Theorem \ref{prop:dressing_Z_approx_outline}, which computes the scaling limit of the quasi-particle fluctuations $Z_k^{\mathcal{Q}} (t)$ from \eqref{eqn:quasi_part_fluct_intro}; it may be useful to recall the discussion around Proposition \ref{prop:sumZlambdaQ} in Section \ref{subsec:heur}.

\subsection{Proof of Theorem \ref{prop:dressing_Z_approx_outline}} 
\label{sec:quasi_ptcl_couple}

To show Theorem \ref{prop:dressing_Z_approx_outline}, we use the following lemma, essentially indicating that the matrix describing the (approximate) linear system \eqref{eqn:Zeq_1} satisfied by the $(Z_k^{\mathcal{Q}}(t))$ is (quantitatively) invertible. Its proof is similar to those of \cite[Lemmas 6.5 and 6.6]{Agg25}, so it is outlined in Appendix \ref{app:Sinv_proofs}.

\begin{lem} \label{lem:both_S_inv}
\label{lem:S_inv}

There exists a constant $\theta_0 (\beta)>0$ such that the below holds if $\theta < \theta_0$. Let $t \in [0, T \log N]$, and let $l, m \in \mathbb{Z}$ be integers satisfying $ -T^{4/3}   \leq  l \leq - T (\log N)^{10} < 0 < T ( \log N)^{10} \leq m \leq T^{4/3}/2  $. Consider the random $(2m-2l+1)  \times (2m-2l+1)$ matrix $\mathbf{S} = \mathbf{S}^{[2l,2m]} = [S_{k i}]_{k, i}$, with entries indexed by $k, i \in \llbracket 2 l, 2 m \rrbracket$, defined by  
\begin{multline}\label{eqn:S0def}
S_{k i} = - 2 \mathfrak{l}(\Lambda_{k} - \Lambda_i) \cdot \chi' \big(\alpha (k - i) + t (\ve(\Lambda_{k}) - \ve(\Lambda_{i}) ) \big)  \\
+ \left(1+ 2 \sum_{j = N_1}^{N_2} \mathfrak{l}(\Lambda_{k} - \Lambda_j) \cdot \chi' \big(\alpha (k - j) + t  (\ve(\Lambda_{k}) - \ve(\Lambda_{j}) ) \big)  \right) \cdot \mathbbm{1}_{i = k}.
\end{multline}
Define $\mathfrak{Q}_j \coloneqq \alpha j + t \ve(\Lambda_j)$, for $j \in \llbracket N_1, N_2 \rrbracket$. The following holds with overwhelming probability. 

\begin{enumerate} 

\item The matrix $\mathbf{S}$ is invertible. Moreover, for any index $j \in \llbracket 2 l, 2 m \rrbracket$, we have 
\begin{equation}\label{eqn:Sbd}
|S_{jj}| \geq (2 + (\log N)^{-1}) \sum_{i = N_1}^{N_2} \left| \mathfrak{l}(\Lambda_j - \Lambda_i) \cdot \chi'\left( \mathfrak{Q}_j - \mathfrak{Q}_i \right) \right| + (\log N)^{-1}.
\end{equation}

\item Let $v = (v_{2 l }, \dots, v_{2 m}) \in \mathbb{R}^{2m-2l+1}$ be a vector such that $|v_i| \leq \delta$ for each $i \in \llbracket 3 l/2, 3 m/2 \rrbracket $. Suppose in addition that, for some real number $\mathfrak{C} \ge 0$, we have $|(\S^{-1} v)_k| \leq  (\log N)^{\mathfrak{C}}$ for each $k \in \llbracket 2 l, 2 m \rrbracket$. Then, for any index $k' \in \llbracket l, m \rrbracket$, we have  
\begin{equation}\label{eqn:sinvbd}
|(\S^{-1} v)_{k'} | \leq \delta \log N + (\log N)^{\mathfrak{C}} e^{- (\log N)^3   / 8 }.
\end{equation}
\end{enumerate} 

\end{lem}

\begin{proof}[Proof of Theorem \ref{prop:dressing_Z_approx_outline}]

   Throughout the proof, let $l = m = \floor{T^{7/6}}$, choices which satisfy the conditions $ -T^{4/3}   \leq  l \leq - T (\log N)^{10} < 0 < T ( \log N)^{10} \leq m \leq T^{4/3} /2$ from Lemma \ref{lem:both_S_inv}. Let $\mathbf{S} $ be the $(2m-2l+1) \times (2m-2l+1)$ matrix defined in Lemma \ref{lem:both_S_inv}, and for $k \in \llbracket N_1, N_2 \rrbracket$, abbreviate $\Xi_k  \coloneqq \Xi(\Lambda_k, k, t)$ (recall Definition \ref{xilambdakt}). In addition, define 
\begin{equation}
    \tilde{Z}_k \coloneqq Z(\Lambda_k, \alpha k  + t \ve(\Lambda_k), t),
\end{equation}

\noindent for each $k \in \llbracket N_1, N_2 \rrbracket$, and define the vector $\tilde{Z}^{\llbracket 2 l , 2 m \rrbracket} = (\tilde{Z}_k)_{k \in \llbracket 2 l , 2 m \rrbracket} \in \mathbb{R}^{2m-2l+1}$.

For any index $k \in \llbracket 2 l, 2 m \rrbracket$, define
\begin{equation}\label{eqn:time-0-linear}
\mathcal E^{(0)}_k = (\mathbf{S} Z^{\mathcal{Q}, \llbracket 2l , 2 m \rrbracket})_k -\Xi_k ,
\end{equation}
where $\mathbf{S} Z^{\mathcal{Q}, \llbracket 2l , 2 m \rrbracket} \in \mathbb{R}^{2m-2l+1}$ denotes the matrix multiplication of $\mathbf{S}$ by the vector $Z^{\mathcal{Q}, \llbracket 2l , 2 m \rrbracket} \coloneqq (Z_k^{\mathcal{Q}})_{k \in  \llbracket 2l , 2 m \rrbracket} \in \mathbb{R}^{2m-2l+1}$. In addition, for any $k \in \llbracket 2 l, 2 m \rrbracket$, define 
\begin{equation}\label{eqn:time-0-linear1}
\mathcal E^{(1)}_k = (\mathbf{S} \tilde{Z}^{ \llbracket 2l , 2 m \rrbracket})_k -\Xi_k .
\end{equation}

Fix $t \in [0, T \log N]$ as in the statement in the theorem. We claim that, with overwhelming probability, we have for each $i \in \{ 0, 1 \}$ that 
\begin{equation}\label{eqn:E0-bound}
\max_{ k \in \llbracket 3 l/2, 3m/2 \rrbracket } |\mathcal E^{(i)} _k|
\leq 
  T^{-c}.
\end{equation}
We will then use \eqref{eqn:E0-bound} to show that, for the given (fixed) $t \in [0, T \log N]$, we have  with overwhelming probability that 
\begin{equation}\label{eqn:Zapproxwts_state0}
\max_{k \in \llbracket l, m \rrbracket} |Z_k^{\mathcal{Q}} -\tilde{Z}_k | \leq T^{-\mathfrak{c}_1}.
\end{equation}
Then we will extend to all $t$ by union bounding over a mesh and using a continuity argument.

The first bound, \eqref{eqn:E0-bound} with $i = 0$ follows from \eqref{eqn:Zeq_1}. Indeed, the difference between \eqref{eqn:time-0-linear} and the left side of \eqref{eqn:Zeq_1} is that the sum over $i$ in \eqref{eqn:Zeq_1}  is over $\llbracket N_1, N_2 \rrbracket$, whereas by definition the matrix multiplication by $\mathbf{S}$ in \eqref{eqn:time-0-linear} is a sum over $\llbracket 2 l, 2 m \rrbracket$. However, on the event $\mathsf{BND}_{\L}(\log N)$ (which has overwhelming probability by Lemma \ref{lem:bd_lem}), if $k \in \llbracket 3 l/2, 3m/2 \rrbracket$ and $i \notin \llbracket 2 l, 2 m \rrbracket$, then we have $|\alpha (k - i) +t (\ve(\Lambda_k) - \ve(\Lambda_i))| > \alpha|i-k| - T \log N \cdot 2 C \log N > c T^{7/6} > \mathfrak{M}$ (by Lemma \ref{lem:ve_tail}), meaning since $\supp \chi' \subseteq [-\mathfrak{M}, \mathfrak{M}]$ that 
\begin{flalign} 
\label{chi0} 
\chi'(\alpha (k - i) +t (\ve(\Lambda_k) - \ve(\Lambda_i))) = 0.
\end{flalign} 

\noindent Therefore, the sum can be truncated to $\llbracket 2 l, 2 m \rrbracket$ without any error. Thus, with $\mathcal{E}_k^{(0)}$ as in \eqref{eqn:time-0-linear}, we have verified \eqref{eqn:E0-bound} at $i = 0$.

To confirm \eqref{eqn:E0-bound} with overwhelming probability when $i = 1$, we restrict to $\mathsf{E}_1 \coloneqq \mathsf{BND}_{\L}(\log N)$, and invoke Proposition \ref{prop:sumZlambdaQ} to obtain an event $\mathsf{E}_2$ of overwhelming probability, on which 
\begin{multline}\label{eqn:Zeq_approx00}
 \sup_{|k| \le 2 T^{7/6}} \Bigg|   
  2 \sum_{i=N_1}^{N_2} \big(  Z(\Lambda_k, \alpha k + t  \ve(\Lambda_k), t) -  Z(\Lambda_i, \alpha i + t  \ve(\Lambda_i), t) \big) \cdot 
  \mathfrak{l}(\Lambda_{k}-\Lambda_i) \\ 
  \times     \chi' \big(\alpha( k- i) +  t  (\ve(\Lambda_{k}) 
- \ve(\Lambda_{i}) ) \big)
+ Z(\Lambda_k, \alpha k + t  \ve(\Lambda_k), t) -\Xi(\Lambda_k, k, t)\Bigg| \leq T^{-\mathfrak{c}_1}.
\end{multline}
To see that the above holds, let us consider the event $\mathsf{E}_2$ provided by Proposition \ref{prop:sumZlambdaQ}. On this event, \eqref{eqn:Zeq_2_prime} simultaneously holds for all $\Lambda \in [- \log N, \log N]$ and $Q \in [-T^{6/5}, T^{6/5}]$, so that we may substitute $\Lambda = \Lambda_{k} \in [-\log N, \log N ]$ (where the last inclusion is by our restriction to $\mathsf{E}_1$), and $Q= \alpha k + t  \ve(\Lambda_k) \in [-T^{6/5}, T^{6/5}]$ (where the last inclusion is by  Lemma \ref{lem:ve_tail}, since $|k| \leq 2 T^{7/6}$ and $t \leq T \log N$), and the bound still holds.

Next, we argue that \eqref{eqn:Zeq_approx00} implies that \eqref{eqn:E0-bound} holds for $i = 1$. Indeed, if $k \in \llbracket 3 l/2, 3m/2 \rrbracket$, then 
\begin{multline}\label{eqn:e1trunc}
      \sum_{i=N_1}^{N_2}Z(\Lambda_i, \alpha i + t  \ve(\Lambda_i), t) \cdot \mathfrak{l}(\Lambda_{k}-\Lambda_i) \cdot \chi'(\alpha( k- i) +  t  (\ve(\Lambda_{k}) 
- \ve(\Lambda_{i}) )) \\
= \sum_{i=2 l }^{2 m}Z(\Lambda_i, \alpha i + t  \ve(\Lambda_i), t) \cdot \mathfrak{l}(\Lambda_{k}-\Lambda_i) \cdot  \chi'(\alpha( k- i) +  t  (\ve(\Lambda_{k}) 
- \ve(\Lambda_{i}) )),
\end{multline}

\noindent since if $i \notin \llbracket 2l , 2 m \rrbracket$, then $\chi(\alpha( k- i) +  t  (\ve(\Lambda_{k}) 
- \ve(\Lambda_{i}) )) = 0$ by \eqref{chi0}. Combining \eqref{eqn:e1trunc} and \eqref{eqn:Zeq_approx00} implies \eqref{eqn:E0-bound} for $i = 1$.

Next, we argue that \eqref{eqn:E0-bound} implies that \eqref{eqn:Zapproxwts_state0} holds (both with the fixed value of $t \in [0, T \log N]$). We restrict to $\mathsf{BND}_{\L}(\log N)$; to the event on which \eqref{eqn:small_Lambda_Zbd} holds; to the event where both outcomes \eqref{eqn:Sbd} and \eqref{eqn:sinvbd} of Lemma \ref{lem:S_inv} hold; and to the event where Lemma \ref{thm:Z_apriori} holds. We claim that upon such restrictions, \eqref{eqn:Zapproxwts_state0} holds.

Indeed, by the statement in the second part of Lemma \ref{thm:Z_apriori}, $|Z_{k}^{\mathcal{Q}}| \leq (\log N)^C$ for all $k \in \llbracket 2 l, 2 m \rrbracket$. By \eqref{eqn:small_Lambda_Zbd} and the restriction to $\mathsf{BND}_{\L(0)}(\log N)$, we have that $|\tilde{Z}_k|\leq (\log N)^C $ for all $k \in \llbracket 2 l, 2 m \rrbracket$. Combining these two bounds yields
\begin{equation}
    \sup_{k \in \llbracket 2 l , 2 m \rrbracket } |Z_k^{\mathcal{Q}} - \tilde{Z}_k | \leq (\log N)^C.
\end{equation}

\noindent Thus, \eqref{eqn:Zapproxwts_state0} follows from taking the difference between \eqref{eqn:E0-bound} at $i = 0$ and \eqref{eqn:E0-bound} at $i = 1$, and then applying Lemma \ref{lem:S_inv} with $v_k = \mathcal E^{(1)}_k - \mathcal E^{(0)}_k$ and $\delta = T^{-c}$.

The bound \eqref{eqn:Zapproxwts_state0} completes the proof, up to showing that this estimate holds for all $t \in [0, T \log N]$ simultaneously. To do this, we define the event $\mathsf{E}'$ on which \eqref{eqn:Zapproxwts_state0} holds for all $t$ in an $N^{-10}$-mesh $\mathcal{T}$ of $[0, T \log N]$. We further restrict to the event on which $\bigcap_{t \geq 0} \mathsf{BND}_{\L(t)}(\log N)$ holds, as well as do Lemmas \ref{lem:Z_cont} and \ref{lem:loc_cent_dif} hold. By \eqref{eqn:Z_dif_bound} (and Lemma \ref{lem:ve_tail}), for any $k \in \llbracket 2l, 2 m \rrbracket$, $|Z(\Lambda_k, \alpha k + t \ve(\Lambda_k), t) - Z(\Lambda_k, \alpha k + t' \ve(\Lambda_k), t')| \le T^{-c}$  if $|t-t'| \le N^{-10}$. Moreover, if $|t-t'|\leq N^{-10}$, then $\max_{k \in \llbracket 2l, 2 m \rrbracket} |Z_k^{\mathcal Q}(t) -Z_k^{\mathcal Q}(t')| \leq (\log N)^C T^{-1/2}$ by \eqref{eqn:loc_cent_dif} and the fact that $|q_j'(t)| = |L_{jj}(t)| \leq \log N$ for all indices $j$. Combining these estimates with the fact that \eqref{eqn:Zapproxwts_state0} holds for all $t \in \mathcal{T}$, we deduce it holds for all $t \in [0, T \log N]$, thereby establishing the theorem. 
\end{proof}

\subsection{Diffusivity of a tracer quasi-particle} 
\label{subsec:tracer_D}

In this section we show Proposition \ref{prop:tracer_fluct_intro} and Corollary \ref{zqi2}, The first identifies the process $\tau \mapsto \mathcal{Z}(\Lambda, \mathfrak{q} + \tau \ve(\Lambda), \tau)$ (recall Definition \ref{def:Z_intro}) that represents the limiting fluctuations for a tracer quasi-particle (see Remark \ref{qtlambda00}). The second indicates that the fluctuations of the quasi-particle trajectories $Z_k^\mathcal{Q}$ and $Z_{k'}^{\mathcal{Q}}$ essentially coincide if $|k-k'| \ll T$ and $|\Lambda_k - \Lambda_{k'}| \ll 1$. We begin with the second. 

\begin{proof}[Proof of Corollary \ref{zqi2}]

By Theorem \ref{thm:quasi_fluct_intro}, there exists a coupling between $\mathbf{L}(0)$ and the process $\mathcal{Z}$ so that $|Z_i^{\mathcal{Q}} (T\tau ) - \mathcal{Z} (\Lambda_i, \alpha i T^{-1} + \tau \ve (\Lambda_i), \tau)| \le T^{-c}$ holds with overwhelming probability, for each $i \in \{ k, k' \}$ and $\tau \in [0, \log N]$. Thus, denoting $\kappa = kT^{-1}$ and $\kappa' = k'T^{-1}$, it we must confirm that, with overwhelming probability,
\begin{flalign}
    \label{zlambdak00estimate} 
\displaystyle\sup_{\tau \in [0, \log N]} \big| \mathcal{Z} (\Lambda_k, \alpha \kappa + \tau \ve (\Lambda_k), \tau) - \mathcal{Z} (\Lambda_{k'}, \alpha \kappa' + \tau \ve (\Lambda_{k'}), \tau) \big| \le T^{-c}.
\end{flalign} 

\noindent By the definition \eqref{def:Z_intro} of $\mathcal{Z}$, the first and third statements of Lemma \ref{lem:ve_tail}, and the fact that $|\Lambda_k - \Lambda_{k'}| \le T^{-\delta}$, it thus suffices to verify that, with overwhelming probability,
\begin{flalign*} 
\displaystyle\sup_{\tau \in [0, \log N]} |\mathfrak{X} (\Lambda_k, \kappa, \tau) - \mathfrak{X} (\Lambda_{k'}, \kappa' , \tau)| \le T^{-c}.
\end{flalign*} 

\noindent This follows from Lemma \ref{lem:holder_cont}, with the facts that $|\kappa - \kappa'| \le T^{-\delta}$ and $|\Lambda_k - \Lambda_{k'}| \le T^{-\delta}$ (and that $|\Lambda_k|, |\Lambda_{k'}| \le \log N$ holds with overwhelming probability, by Lemma \ref{lem:bd_lem}).
\end{proof} 

To show Proposition \ref{prop:tracer_fluct_intro}, we first require the following lemma. Below, we recall the function $T(\lambda, \lambda') = 2 \log |\lambda - \lambda'|$ from \eqref{eqn:T_kernel}, and its dressing $T^{\dr}(\lambda, \lambda') \coloneqq [(1- \theta \T \varrhobetabf )^{-1}T(\cdot, \lambda')](\lambda)$ in its first argument. 

\begin{lem}\label{lem:Tdr_sym}

For all $\lambda, \mu \in \mathbb{R}$, we have $T^{\dr}(\lambda, \mu) = T^{\dr}(\mu, \lambda)$. 
\end{lem}

\begin{proof}

Observe from \eqref{drf} that, for any $f \in \mathcal{H}$,
\begin{equation}\label{eqn:dr_defining}
    f^{\dr} = f + \theta \T \varrhobetabf f^{\dr}.
\end{equation}

\noindent Also observe for any $f, g \in \mathcal{H}$ that
\begin{equation}\label{eqn:dr_adj}
    \theta \int_{-\infty}^{\infty}f(\lambda) g^{\dr}(\lambda) \varrho_{\beta}(\lambda) d \lambda = \theta \int_{-\infty}^{\infty}f^{\dr}(\lambda) g(\lambda)  \varrho_{\beta}(\lambda) d \lambda,
\end{equation}

\noindent as $\mathbf{T} \bm{\varrho}_{\beta}$ (and thus the dressing operator $(1 - \theta \mathbf{T} \bm{\varrho}_{\beta})^{-1}$) is self-adjoint with respect to $L^2 (\varrho_{\beta} d\lambda)$. So, 
    \begin{flalign*}
        T^{\dr}(\lambda, \mu) & = T(\lambda, \mu) + \theta \int_{-\infty}^{\infty}T(\lambda, \lambda') T^{\dr}(\lambda', \mu)  \varrho_{\beta}(\lambda') d \lambda' \\ 
         & =  T(\mu, \lambda) +  \theta \int_{-\infty}^{\infty}T(\lambda', \lambda) T^{\dr}(\lambda', \mu)  \varrho_{\beta}(\lambda') d \lambda' 
         \\ 
         & =  T(\mu, \lambda) + \theta \int_{-\infty}^{\infty}T^{\dr}(\lambda', \lambda) T(\lambda', \mu)  \varrho_{\beta}(\lambda') d \lambda' 
         =  T^{\dr}(\mu, \lambda).
\end{flalign*}

\noindent where in the first statement we used \eqref{eqn:dr_defining}; in the second we used the fact that $T(\lambda, \mu) = T(\mu, \lambda)$; in the third we used \eqref{eqn:dr_adj}; and in the fourth we again used \eqref{eqn:dr_defining}. This verifies the lemma.
\end{proof}

\begin{proof}[Proof of Proposition \ref{prop:tracer_fluct_intro}]

Since $\mathcal{Z}$ is a Gaussian process, it suffices to compute the covariance of the map on the left side of \eqref{eqn:Z_brownian}, for different time values $\tau_1$ and $\tau_2$. So, fix real numbers $\tau_1 \le \tau_2$. In what follows, we will let $\mathfrak{q}_1 \coloneqq \mathfrak{q} + \tau_1 \ve(\Lambda)$ and $\mathfrak{q}_2 = \mathfrak{q} + \tau_2 \ve(\Lambda)$, and let $f_{\Lambda_1} \coloneqq \psi_{\Lambda_1, (\mathfrak{q}_1-\tau_1 \ve(\Lambda_1))/\alpha, \tau_1}$, and $f_{\Lambda_2} \coloneqq \psi_{\Lambda_2, (\mathfrak{q}_2-\tau_2 \ve(\Lambda_2))/\alpha, \tau_2}$, where $\psi_{\Lambda, \kappa, \tau}$ is defined in \eqref{eqn:logphi_intro}. 

In addition, below we denote by $\mathbf{D} = (1 - \theta \mathbf{T} \bm{\varrho}_{\beta})^{-1}$ the dressing operator (Definition \ref{def:dress_intro}) and, if $f(r,\lambda)$ is a function of two variables, the notation $\mathbf{D}_{\lambda} f (r,\lambda)$ means that we apply $\mathbf{D}$ to the function $\lambda \mapsto f(r, \lambda)$, and then evaluated the result at $\lambda$. Finally, define the operator $\mathbf{O}  = \boldsymbol{(\varsigma_0^{\dr})^{-1}} (1 - \theta \mathbf{T} \boldsymbol{\varrho}_{\beta})^{-1}$ as the composition of the dressing operator with multiplication by $\varsigma_0^{\dr}(\Lambda)^{-1}$; moreover, as with $\mathbf{D}_{\lambda}$, the subscript in $\mathbf{O}_{\Lambda}$ indicates the variable in which the operator acts. 

Under this notation, we compute the limiting covariance 
    \begin{flalign}\label{eqn:Zcov_calc}
    \begin{aligned}
       \Cov & \left(\mathcal{Z}(\Lambda, \mathfrak{q}_1, \tau_1),  \mathcal{Z}(\Lambda, \mathfrak{q}_2, \tau_2) \right) \\
        &=\mathbf{O}_{\Lambda_1} \mathbf{O}_{\Lambda_2}  \Cov\left( \mathcal{W}^{\dr}( f_{\Lambda_1} ), \mathcal{W}^{\dr}( f_{\Lambda_2} ) \right) |_{\Lambda_1=\Lambda_2 = \Lambda} \\
        &= \mathbf{O}_{\Lambda_1} \mathbf{O}_{\Lambda_2} \int_{-\infty}^{\infty}\int_{-\infty}^{\infty}\mathbf{D}_{\lambda} f_{\Lambda_1} (r, \lambda) \mathbf{D}_{\lambda}  f_{\Lambda_2} (r, \lambda)  d r \varrho(\lambda) d \lambda |_{\Lambda_1=\Lambda_2 = \Lambda} \\
        &=   (\varsigma_0^{\dr}(\Lambda))^{-2} \int_{-\infty}^{\infty}\int_{-\infty}^{\infty}\mathbf{D}_{\Lambda_1} \mathbf{D}_{\lambda} \big( T(\lambda, \Lambda_1) \cdot 
\left(
    \mathbbm{1}_{ \mathfrak{q}_1-  \alpha r -\tau_1 \ve(\Lambda_1)   > 0  } 
    - 
    \mathbbm{1}_{ \mathfrak{q}_1   -  \alpha r  -   \tau_1  \ve(\lambda)   > 0}
   \right) \big) \\
   & \quad   \times \mathbf{D}_{\Lambda_2} \mathbf{D}_{\lambda} \big(
T(\lambda,\Lambda_2) \cdot 
\left(
    \mathbbm{1}_{\mathfrak{q}_2 -  \alpha r -\tau_2 \ve(\Lambda_2)    > 0 } 
    - 
    \mathbbm{1}_{ \mathfrak{q}_2  -  \alpha r  -   \tau_2  \ve(\lambda)  > 0}
   \right) 
   \big) d r \varrho(\lambda) d\lambda |_{\Lambda_1=\Lambda_2 = \Lambda} .
    \end{aligned}
    \end{flalign}
    In the display above, we have used \eqref{eqn:intro_Z} and \eqref{eqn:Xdef} (and recalled that $T(\lambda, \lambda') = 2 \log |\lambda-\lambda'|$). 

     We claim that for any $\tau \geq 0$ and any $\Lambda, \lambda, \mathfrak{q}_0,r \in \mathbb{R}$,
    \begin{multline}\label{eqn:key_Tdr_id}
        \mathbf{D}_{\Lambda} \mathbf{D}_{\lambda} T(\lambda, \Lambda) \cdot 
\left(
    \mathbbm{1}_{ \mathfrak{q}_0-  \alpha r -\tau \ve(\Lambda)   > 0  } 
    - 
    \mathbbm{1}_{ \mathfrak{q}_0   -  \alpha r  -   \tau  \ve(\lambda)   > 0}
   \right) \\
        = 
        \big(
    \mathbbm{1}\{ \mathfrak{q}_0-  \alpha r -\tau \ve(\Lambda)   > 0  \} 
    - 
    \mathbbm{1}\{ \mathfrak{q}_0   -  \alpha r  -   \tau \ve(\lambda)   > 0\}
   \big) \cdot T^{\dr}(\lambda, \Lambda).
    \end{multline}

To verify this, we first note that 
    \begin{multline}\label{eqn:Zvar_i}
\mathbf{D}_{\Lambda_1} \mathbf{D}_{\lambda_1}  \big( T(\lambda_1, \Lambda_1) \cdot 
\left( \mathbbm{1}\{\mathfrak{q}_0 - \tau v(\Lambda_1) - \alpha r > 0 \}  - \mathbbm{1}\{\mathfrak{q} - \tau v(\lambda_1) - \alpha r>0\}\right) \big) |_{\Lambda_1 =\Lambda, \lambda_1 = \lambda} \\
 = \mathbf{D}_{\Lambda} \big( T^{\dr}(\lambda, \Lambda) \cdot  \mathbbm{1}\{\mathfrak{q}_0- \tau \ve(\Lambda) - \alpha r > 0 \}  \big)\\
 - \mathbf{D}_{\lambda}   \big( T^{\dr}(\Lambda, \lambda) \cdot \mathbbm{1}\{\mathfrak{q}_0 - \tau \ve(\lambda) - \alpha r > 0 \}  \big).
\end{multline}
Let 
$$F(\lambda, \Lambda) \coloneqq  T^{\dr}(\Lambda, \lambda) \cdot \mathbbm{1}\{\mathfrak{q}_0 - \tau \ve(\lambda) - \alpha r > 0 \}.$$ 
Denote $F^{\dr}(\lambda, \Lambda) = \mathbf{D}_{\lambda} F(\lambda, \Lambda)$ the dressing of $F$ in its first argument $\lambda$. We then obtain for the final expression in \eqref{eqn:Zvar_i}, 
\begin{flalign}\label{eqn:Zvar_i2}
\begin{aligned} 
    F(\Lambda, \lambda) - F(& \lambda, \Lambda)  + \theta \int_{-\infty}^{\infty} \big(T(\Lambda, \mu)F^{\dr}(\mu, \lambda) - T(\lambda, \mu)F^{\dr}(\mu, \Lambda) \big) \cdot  \varrho_{\beta}(\mu) d\mu \\
    & =   \left(
    \mathbbm{1}\{ \mathfrak{q}_0-  \alpha r -\tau \ve(\Lambda)   > 0  \} 
    - 
    \mathbbm{1}\{ \mathfrak{q}_0   -  \alpha r  -   \tau  \ve(\lambda)   > 0\}
   \right) \cdot T^{\dr}(\lambda, \Lambda) \\
  & \qquad  + \theta \int_{-\infty}^{\infty}(T^{\dr}(\mu, \Lambda)F(\mu, \lambda) - T^{\dr}(\mu, \lambda)F(\mu, \Lambda)) \cdot \varrho_{\beta}(\mu) d\mu   \\
&  = \left(
    \mathbbm{1}\{ \mathfrak{q}_0-  \alpha r -\tau \ve(\Lambda)   > 0  \} 
    - 
    \mathbbm{1}\{ \mathfrak{q}_0  -  \alpha r  -   \tau  \ve(\lambda)   > 0\}
   \right) \cdot T^{\dr}(\lambda, \Lambda) .
   \end{aligned} 
\end{flalign}
The first equality in \eqref{eqn:Zvar_i2} is obtained by the definition of $F$ (with Lemma \ref{lem:Tdr_sym}) and the self-adjointness of the dressing operator with respect to the inner product associated with the measure $L^2 (\varrho_{\beta} d \lambda)$ (see also \eqref{eqn:dr_defining} and \eqref{eqn:dr_adj}). The final equality is due to the fact that the term with the integral vanishes, by the symmetry of $T^{\dr}$ (Lemma \ref{lem:Tdr_sym}). Together with \eqref{eqn:Zvar_i}, this establishes \eqref{eqn:key_Tdr_id}.

Thus, we may compute the last two lines of \eqref{eqn:Zcov_calc};
\begin{flalign*}
    & \int_{-\infty}^{\infty}\int_{-\infty}^{\infty}\mathbf{D}_{\Lambda_1} \mathbf{D}_{\lambda} \big( T(\lambda, \Lambda_1) \cdot 
\left(
    \mathbbm{1}\{ \mathfrak{q}_1-  \alpha r  -\tau_1 \ve(\Lambda_1)   > 0  \} 
    - 
    \mathbbm{1}\{ \mathfrak{q}_1   -  \alpha r  -   \tau_1  \ve(\lambda)   > 0\}
   \right) \big) \\
   & \qquad \qquad  \times \mathbf{D}_{\Lambda_2} \mathbf{D}_{\lambda} \left(
T(\lambda, \Lambda_2) \cdot 
\big(
    \mathbbm{1}\{ \mathfrak{q}_2 -  \alpha r -\tau_2 \ve(\Lambda_2)    > 0  \} 
    - 
    \mathbbm{1}\{ \mathfrak{q}_2  -  \alpha r  -   \tau_2  \ve(\lambda)  > 0\}
   \right) 
   \big) \\
   & \qquad \qquad  \times  \varrho(\lambda) dr d\lambda |_{\Lambda_1=\Lambda_2 = \Lambda} \\ 
   & \quad =  \int_{-\infty}^{\infty}|T^{\dr}(\lambda, \Lambda)|^2 \int_{-\infty}^{\infty} 
\left(
    \mathbbm{1}\{ \mathfrak{q}_1-  \alpha r -\tau_1 \ve(\Lambda)   > 0  \} 
    - 
    \mathbbm{1}\{ \mathfrak{q}_1   -  \alpha r  -   \tau_1  \ve(\lambda)   > 0\}
   \right)^2 \varrho(\lambda) dr d\lambda \\
   & \quad = \tau_1 |\alpha|^{-1} \int_{-\infty}^{\infty}|T^{\dr}(\lambda, \Lambda)|^2 |\ve(\lambda) - \ve(\Lambda)|
    \varrho(\lambda) d\lambda.
\end{flalign*}
\noindent where the first equality holds by \eqref{eqn:key_Tdr_id} and the relations $   \mathfrak{q}_1-\tau_1 \ve(\Lambda) =\mathfrak{q}= \mathfrak{q}_2-\tau_2 \ve(\Lambda)$, and the second by performing the integral over $r$. Substituting this back into \eqref{eqn:Zcov_calc}, and recalling \eqref{eqn:diffusivityintro}, we obtain that $\Cov (\mathcal{Z}(\Lambda, \mathfrak{q}_1, \tau_1), \mathcal{Z}(\Lambda, \mathfrak{q}_2 ,\tau_2) = \tau_1 \cdot \mathcal{D}(\Lambda)$. This confirms that the covariance of the Gaussian processes on the left and right sides of \eqref{eqn:Z_brownian} coincide, which establishes the proposition. 
\end{proof}

 \begin{remark}

 \label{qtlambda00}
 
    The random variable $\tau \mapsto \mathcal{Z}(\Lambda, \mathfrak{q} + \tau \ve(\Lambda), \tau)$ represents the scaling limit of the fluctuations of a tracer quasi-particle with spectral parameter $\Lambda$ and macroscopic position $\mathfrak{q}$, at time $\tau$. There does not appear to be a unique microscopic definition of such a ``tracer quasi-particle'' in the literature, although there are multiple reasonable notions of it; let us briefly discuss two here (recalling Definition \ref{def:quasi_intro} below). One way to define a tracer quasi-particle, with spectral parameter $\Lambda$ and initial position $q_k(0) \approx \alpha k$ (for some index $k$), is to set it equal to $Q_k (t)$, conditional on the event that $\{\Lambda_k \approx \Lambda \}$. A second way is to impose that its trajectory $Q_{\Lambda}(t)$ satisfies the asymptotic scattering relation describing quasi-particles, but with a fixed value of its spectral parameter; specifically, it can be set to an approximate solution of
\begin{multline}
\label{eqn:tracer_quasi}
 Q_{\Lambda}(t) \approx  Q_{\Lambda}(0) + \Lambda t 
+ 2\,\mathrm{sgn}(\alpha) \Bigg( 
\displaystyle\sum_{i:\, Q_i(0) < Q_{\Lambda}(0)} \log \lvert \Lambda - \lambda_i \rvert -
\displaystyle\sum_{i:\, Q_i(t) < Q_{\Lambda}(t)} \log \lvert \Lambda - \lambda_i \rvert \Bigg).
\end{multline}

\noindent In the scaling limit, one should have
\begin{flalign}
    \label{qlambdat0} 
Q_{\Lambda}(t) \coloneqq Q_{\Lambda}(0) + t \ve(\Lambda) + T^{1/2} \cdot \mathcal{Z}(\Lambda, \mathfrak{q}, \tau).
\end{flalign}

\noindent Indeed, for $\Lambda = \Lambda_k$ we have by Theorem \ref{thm:quasi_fluct_intro} that these describe the fluctuations of a quasi-particle with spectral parameter $\Lambda_k$, which is in heuristic alignment with the first description above (especially since the coupling in Theorem \ref{thm:quasi_fluct_intro} holds with overwhelming probability, suggesting that it should persist under conditioning on events such as $\{ \Lambda_k \approx \Lambda \}$). Moreover, \eqref{qlambdat0} can be shown (though we will not pursue this in detail here) to be an approximate solution to \eqref{eqn:tracer_quasi}, in alignment with the second description above.

 \end{remark}

\section{Fluctuations of integrated currents}
\label{sec:current_fluct}
In this Section we prove Theorem \ref{thm:current_fluct_intro}. Throughout, we work under Assumption \ref{ass:NT_assumption} (when dealing with the Toda lattice on $\llbracket N_1, N_2 \rrbracket$) and assume implicitly that $\theta<\theta_0(\beta)$ is sufficiently small. We recall the notation $Z(\Lambda, Q, t)$ from \eqref{eqn:Z_outline} and $Z_j(t)$ from \eqref{eqn:quasi_part_fluct_intro}; moreover, we will use the notation from Definitions \ref{kti00}, \ref{def:integrated_curr_spatial}, \ref{def:Wdress}, \ref{def:Z_intro}, \ref{def:quasi_intro}, and \ref{def:smoothed_log}.

\subsection{Preliminaries}

In this section we provide several preliminary statements that will facilitate the proof of Theorem \ref{thm:current_fluct_intro}. Throughout, we recall the Toda lattice Flaschka variables from \eqref{abr} and thermal equilibrium from Definition \ref{def:inf_thermal_eq}. The below lemma indicates that the Toda lattice current fluctuations are closely approximated by functionals of quasi-particles; this will enable us to use Theorem \ref{thm:quasi_fluct_intro} to show Theorem \ref{thm:current_fluct_intro}. Below we recall the integrated currents $J_t^{[m]}(q,q') $ from Definition \ref{def:integrated_curr_spatial}. 

\begin{lem}[{\cite[Proposition 2.10]{Agg25a}}]\label{prop:spatial_int_current}

The following holds with overwhelming probability, with respect to the Toda lattice on $\llbracket N_1, N_2 \rrbracket$, initialized from thermal equilibrium (as in Assumption \ref{ass:NT_assumption}). For any integer $m \in \llbracket 0, (\log N)^{1/10} \rrbracket$ and real numbers $t \in [0, T \log N]$ and $q, q' \in [-\alpha T^2, \alpha T^2]$, 
\begin{equation}\label{eqn:our_version_of_Pro210}
\Bigg| 
 J_t^{[m]}(q,q')    - 
\bigg( \sum_{j : Q_j(0)  < q} \lambda_{j}^{m}-\sum_{j : Q_j(t)  < q'} \lambda_{j}^{m}  \bigg) \Bigg| 
\leq (3 \log N)^{m+10}.
\end{equation}
\end{lem}

Lemma \ref{prop:spatial_int_current} concerns the Toda lattice on $\llbracket N_1, N_2 \rrbracket$, while Theorem \ref{thm:current_fluct_intro} concerns this model on $\mathbb{Z}$. The following lemma indicates that one may couple these dynamics, so that they nearly coincide away from the boundaries of $\llbracket N_1, N_2 \rrbracket$. 

\begin{lem}[{\cite[Proposition 2.5]{Agg25a}}]\label{lem:fin_inf_t_coupling}

Let $(\mathbf{a}^{\llbracket N_1, N_2 \rrbracket}(t); \mathbf{b}^{\llbracket N_1, N_2 \rrbracket}(t))$ denote the Flaschka variables of the Toda lattice on $\llbracket N_1, N_2 \rrbracket$, initialized under thermal equilibrium, where $\mathbf{a}^{\llbracket N_1, N_2 \rrbracket}(t) = ((a_i^{\llbracket N_1, N_2 \rrbracket}(t))$ and $\mathbf{b}^{\llbracket N_1, N_2 \rrbracket}(t) = (b_i^{\llbracket N_1, N_2 \rrbracket}(t))$, over $i \in \llbracket N_1, N_2 \rrbracket$. Further let $(\mathbf{a}(t); \mathbf{b}(t))$ denote the Flaschka variables of the Toda lattice on $\mathbb{Z}$, initialized under thermal equilibrium, where $\mathbf{a}(t)= (a_i(t))$ and $\mathbf{b}(t) =  (b_i(t))$, over $i \in  \mathbb{Z}$. Suppose $(\mathbf{a}(t); \mathbf{b}(t))$ and $(\mathbf{a}^{\llbracket N_1, N_2 \rrbracket}(t); \mathbf{b}^{\llbracket N_1, N_2 \rrbracket}(t))$ share the same initial data on $\llbracket N_1, N_2 \rrbracket$, namely, $a_i(0) = a_i^{\llbracket N_1, N_2 \rrbracket}(0)$ for all $i \in \llbracket N_1, N_2-1\rrbracket$ and $b_i(0) = b_i^{\llbracket N_1, N_2 \rrbracket}(0)$ for all $i \in \llbracket N_1, N_2 \rrbracket$. Then, the following two statements hold.
\begin{enumerate}
    \item For all $t \geq 0$, $(\mathbf{a}(t); \mathbf{b}(t))$ has the same law as $(\mathbf{a}(0); \mathbf{b}(0))$.
    \item  With overwhelming probability, we have
    \begin{equation}\label{eqn:fin_inf_dif}
   \sup_{t \in [0, (\log N)^{10} T]} \max_{i \in \llbracket N_1+T^2, N_2-T^2 \rrbracket} \big( |a_{i}^{\llbracket N_1, N_2 \rrbracket}(t) - a_{i}(t)| +|b_{i}^{\llbracket N_1, N_2 \rrbracket}(t) - b_{i}(t)| \big) \leq e^{-T^2/5}.
\end{equation}
\end{enumerate}
\end{lem}

The next lemma, to be shown in Appendix \ref{sec:inf_vol_bd_proof}, bounds the locations and variables in the infinite volume Toda lattice. It in particular implies that the integrated currents \eqref{eqn:inf_integrated_curr} are well-defined for the Toda lattice on $\mathbb{Z}$ (in that the integer $k$ there exists almost surely).

\begin{lem}\label{lem:inf_vol_bds}

There exists a constant $\mathfrak{c}>0$ such that the following holds. 
    Let $(a_i(t), b_i(t))_{i \in \mathbb{Z}}$, denote the Toda lattice on $\mathbb{Z}$, initialized under thermal equilibrium. 
    Let $T > 1$ and $K \geq T $ be real numbers with $T \leq K^{9/10}$. With probability at least $1- \mathfrak{c}^{-1} e^{-\mathfrak{c} (\log K)^2} $, we have for all $t \in [0, T]$ that
    \begin{flalign}\label{eqn:q_i_inf_bd}
    \begin{aligned}
    \displaystyle\sup_{i \le -K} q_i(t) < - T (\log T)^2; \qquad \displaystyle\inf_{i \ge K}  q_i(t) >  T (\log T)^2,
    \end{aligned}
    \end{flalign}
    and 
    \begin{equation}\label{eqn:ab_i_inf_bd}
       \displaystyle\max_{i \in \llbracket -K^2/2, K^2/2\rrbracket} |a_i(t)| \leq   3 \log K; \qquad \displaystyle\max_{i \in \llbracket - K^2/2, K^2/2 \rrbracket}  |b_i(t)| \leq 3 \log K.
    \end{equation}
\end{lem}

Next, we state a concentration estimate, which  bounds the quantity $\sum_{j=N_1}^{N_2} (\mathbbm{1}\{ Q_j(t) < q\} - \chi( q-Q_j(t)))$ at any time $t$. This will allow us to replace the indicator in the expression for the current provided by \eqref{eqn:our_version_of_Pro210} by its smoothed version, $\chi$.

\begin{lem}\label{lem:lem_chireplacement}
Fix $t \in [0, T \log N]$ and $q \in [-T^{4/3}, T^{4/3}]$. With overwhelming probability, we have for each integer $m \in \llbracket 0, (\log N)^{1/10} \rrbracket$ that 
\begin{equation}\label{eqn:m_current_reg}
        \left|\sum_{ j =1}^N \lambda_j^m \cdot \big( \mathbbm{1}\{ Q_j(t) < q\}  -  \chi(q-Q_j(t)) \big) \right| 
        \leq  \mathfrak{M}^{1/2}(\log N)^{m+15}.
\end{equation}
\end{lem}
\begin{proof}

Let $G(x) = \mathbbm{1}\{ x < 0\} - \chi(-x)$. This satisfies the assumptions on $G$, with $(B,S) = (10,\mathfrak{M})$ (recall Definition \ref{def:chi_def}), of the $\mathfrak{K} = \log N$ case of Lemma \ref{lem:Ndeltabd}. Moreover, $F(\lambda) \coloneqq \lambda^m$ satisfies Assumption \ref{ass:F} with $A = (\log N)^m$, as $m \leq (\log N)^{1/10}$. Thus Lemma \ref{lem:Ndeltabd} and Lemma \ref{lem:conc_estimate_fixedq} together give
\begin{equation}\label{eqn:m_current_conc_2}
        \left|\sum_{ j =1}^N \lambda_j^m \cdot \big( \mathbbm{1}\{ Q_j(t) < q\}  -  \chi(q-Q_j(t)) \big) \right| 
        \leq \mathfrak{M}^{1/2}(\log N)^{m+15}.
\end{equation}

\noindent This, together with a union bound over $m \in \llbracket 0, (\log N)^{1/10} \rrbracket$, yields \eqref{lem:lem_chireplacement}.
\end{proof}

The following concentration bound constitutes a minor generalization of Lemma \ref{lem:z_sumconc_withF}.
\begin{lem}\label{lem:Z_sum_conc_general}
    There exists a constant $\mathfrak{c}>0$ such that the following holds. Let $t \in [0, T \log N]$ be a real number, and denote $\tau = tT^{-1}$. Also let $A \geq 1$ be a real number and $F : \mathbb{R} \rightarrow \mathbb{R}$ be a function satisfying Assumption \ref{ass:F} such that, for all real numbers $\lambda, \lambda' \in [-\log N, \log N]$ with $|\lambda-\lambda'| \leq e^{-(\log N)^{1/2}}$, we have
    \begin{equation}\label{eqn:Fhold2}
        |F(\lambda) - F(\lambda')| \leq A |\lambda-\lambda'|^{1/4}.
    \end{equation}
    
     \noindent Suppose $\mathcal{W}^{\dr}$ and $\L(0)$ are coupled so that \eqref{eqn:tracer_couple_outline} holds. Then we have the below estimates.
    \begin{enumerate}
        \item With overwhelming probability,  \begin{multline}\label{eqn:Z_sum_conc_general1}
    \sup_{|q| \leq T (\log N)^3}  \Bigg| \sum_{i=N_1}^{N_2} \chi' \big(q-(q_i(0) + t \ve(\Lambda_i)) \big) \cdot F(\Lambda_i) \cdot Z_i^{\mathcal Q}(t)    \\
      - \alpha^{-1} \int_{-\infty}^{\infty} F(\lambda) \mathcal{Z}(\lambda, qT^{-1}, \tau) \varrho(\lambda)  d \lambda  \Bigg| \leq A T^{-\mathfrak{c}}.
    \end{multline}
    \item Let $B \in [1, \log N]$, $U  \in [ \mathfrak{M}, T \log N]$, and $\mathfrak{Q} \in [-T(\log N)^2, T (\log N)^2]$ be real numbers, and $G:  \mathbb{R} \rightarrow \mathbb{R}$ be a function satisfying Assumption \ref{ass:G2}. Then, with overwhelming probability,  \begin{multline}\label{eqn:Z_sum_conc_general2}
      \Bigg| \sum_{i=N_1}^{N_2} G'(q_i(0) + t \ve(\Lambda_i) ) \cdot F(\Lambda_i) \cdot  Z_i^{\mathcal{Q}}(t)    \\
      - \alpha^{-1} \int_{-\infty}^{\infty}\int_{-\infty}^{\infty}G'(q) F(\lambda) \mathcal{Z}(\lambda, qT^{-1}, \tau) \varrho(\lambda)  d \lambda dq \Bigg| \leq A  T^{-\mathfrak{c}}.
    \end{multline}
    \end{enumerate}

\end{lem}

\begin{proof}

    We first restrict to $\mathsf{E}_1$, the event provided by Theorem \ref{prop:dressing_Z_approx_outline}; to $\mathsf{E}_2$, the event on which the outcomes of Lemma \ref{lem:second_der_bd} and of Lemma \ref{lem:q_spacing} with $T = 0$ both hold; to $\mathsf{E}_3 = \mathsf{BND}_{\L(0)}(\log N)$. By \eqref{estimatesm0}, \eqref{eqn:Z_approx_outline} (which holds on $\mathsf{E}_1$) and the fact (from our restriction to $\mathsf{E}_2 \cap \mathsf{E}_3$, together with \eqref{eqn:vebd} from Lemma \ref{lem:ve_tail}) that for any $|q|\leq T (\log N)^3$ there are at most $\mathfrak{M} (\log N)^6$ nonzero terms in the sum below, each with $\varphi_0(j) \in \llbracket -T^{7/6}, T^{7/6} \rrbracket$, we have
\begin{multline}
    \label{fz00} 
   \Bigg| \sum_{j=1}^N F(\lambda_j) \cdot Z_{\varphi_0(j)}^{\mathcal{Q}}(t) \cdot \chi'(q-Q_j(0)- t \ve(\lambda_j))  \\
    -\sum_{j=1}^N F(\lambda_j) \cdot Z(\lambda_j, \alpha \varphi_0(j) + t \ve(\lambda_j), t) \cdot \chi'(q-Q_j(0)- t \ve(\lambda_j))   \Bigg|
    \leq A T^{-c}.
\end{multline}

Next, we restrict further to the overwhelmingly probable event $\mathsf{E}_4$ provided by Lemma \ref{lem:z_sumconc_withF} (with $f = 1$), which gives
\begin{multline}\label{eqn:main_est_2}
   \Bigg| \sum_{j=1}^N F(\lambda_j) \cdot Z(\lambda_j, \alpha \varphi_0(j) + t \ve(\lambda_j), t) \cdot \chi'(q-Q_j(0)- t \ve(\lambda_j))  \\ 
    -\alpha^{-1} \int_{-\infty}^{\infty}  F(\lambda) Z(\lambda,q,t)   \varrho(\lambda) d\lambda \Bigg| 
    \leq A T^{-c}.
\end{multline}

\noindent Further restrict to the overwhelmingly probable event $\mathsf{E}_5$, on which \eqref{eqn:small_Lambda_Zbd} and \eqref{eqn:small_Lambda_mcZbd} both hold. By these two estimates (and \eqref{estimatef}), we have 
\begin{flalign} 
\label{fz0} 
\mathbbm{1}_{|\Lambda| > \log N} \cdot |F(\Lambda)| \cdot |\mathcal{Z} (\Lambda, qT^{-1}, tT^{-1}) - Z (\Lambda, q, t)| \leq  A e^{|\Lambda|^{1/2}} \cdot (\log N)^{C}(|\Lambda|+1)^2.
\end{flalign} 

\noindent Hence, for any $|q| \leq T (\log N)^3$, we have
\begin{equation}\label{eqn:main_est_3}
 \left| \int_{-\infty}^{\infty}  F(\lambda) Z(\lambda,q,t)   \varrho(\lambda) d\lambda - \int_{-\infty}^{\infty}  F(\lambda) \mathcal{Z}(\lambda,qT^{-1},tT^{-1})   \varrho(\lambda) d\lambda \right|\leq A T^{-c},
\end{equation}

\noindent which is obtained by using \eqref{eqn:tracer_couple_outline} to estimate the integral over $\lambda \in [- \log N, \log N]$, and using \eqref{fz0} with Lemma \ref{lem:varrho_bd} to bound the integral over $[- \log N, \log N]^{\complement}$ (see \eqref{frho00} for a very similar argument). Combining \eqref{fz00}, \eqref{eqn:main_est_2}, and \eqref{eqn:main_est_3} verifies \eqref{eqn:Z_sum_conc_general1}.

Next, we show \eqref{eqn:Z_sum_conc_general2}. For this, by the bound on the support and derivative of $G$, and by the uniformity over $q$ in \eqref{eqn:Z_sum_conc_general1}, it suffices to show that 
\begin{equation}\label{eqn:G_chi_conv}
    \sup_{x \in \mathbb{R}} \left| G'(x)-\int_{-T (\log N)^3}^{T (\log N)^3} G'(q)\chi'(q-x) dq \right| \leq  T^{-c}.
\end{equation}

\noindent We compute
\begin{align*}
 \left| G'(x)-\int_{-T (\log N)^3}^{T (\log N)^3} G'(q)\chi'(q-x) dq \right| &=  \left| G'(x)-\int_{-\infty}^{\infty} G'(q)\chi'(q-x) dq \right| \\
 &\leq \int_{-\infty}^{\infty}  \left|G'(x)-G'(q)\right|\chi'(q-x) dq  \\
 &\leq 2 \mathfrak{M} \cdot \frac{B}{U^2} \mathfrak{M} \cdot \frac{10}{\mathfrak{M}} \leq 20 B \mathfrak{M}^{-1} \leq  T^{-c},
\end{align*}

\noindent where in the first statement we used the fact that $\supp G' \subseteq [-T(\log N)^3, T(\log N)^3]$ (by \eqref{estimateg}); in the third, we used the facts that $|G''| \le BU^{-2}$, that $\supp \chi' \subseteq [-\mathfrak{M}, \mathfrak{M}]$, and that $|\chi'| \le 10 \mathfrak{M}^{-1}$ (by \eqref{estimatesm0} and \eqref{estimateg}). This gives \eqref{eqn:G_chi_conv}, thereby completing the proof of the lemma. 
\end{proof}

Finally, we have the following lemma involving only limiting objects. The lemma will allow us to uncover a simple form for the multivariate Gaussian describing the joint limit of the integrated currents (in Theorem \ref{thm:current_fluct_intro}), and it will also be useful in the proof of Corollary \ref{cor:twopoint_intro} in Section \ref{sec:charge_fluct}. 
\begin{lem}\label{lem:ttfe}
Let  $\tau \geq 0$ be a real number, and let $F  : \mathbb{R} \rightarrow \mathbb{R}$ be a continuous function satisfying $\sup_{x \in \mathbb{R}} e^{\sqrt{|x|}}  |F(x)| < \infty$; in addition, let $G : \mathbb{R} \rightarrow \mathbb{R}$ be a compactly supported, continuously differentiable function. Define $\mathcal{W}$ and $\mathcal{W}^{\dr}$ as in Definitions~\ref{def:white_noise_dress} and \ref{def:Wdress}, and define $\tilde G:\mathbb{R}^2\to\mathbb{R}$ by setting $\tilde G(r,\lambda):=G (\alpha r+\tau \ve(\lambda))$. Then,  almost surely we have
\begin{flalign} 
\label{eqn:ttfe}
\mathcal{W}^{\dr}( \tilde{G} \cdot F ) +\alpha^{-1} \int_{-\infty}^{\infty} \int_{-\infty}^{\infty} G'(\mathfrak{q}) \mathcal{Z}(\lambda, \mathfrak{q} ,\tau) F(\lambda) \varrho (\lambda) d \lambda d \mathfrak{q}= \mathcal{W}( \tilde{G} \cdot F^{\dr}  ).
\end{flalign} 
Moreover, recalling $\mathfrak{Y}$ from \eqref{yprocess}, for any $(\mathfrak{q}, \mathfrak{q}', \tau, m) \in \mathbb{R} \times \mathbb{R} \times \mathbb{R}_{\geq 0} \times \mathbb{Z}_{\geq 0}$, almost surely we have
\begin{multline}
\label{eqn:ttfe2}
\mathfrak{Y}(\mathfrak{q},\mathfrak{q}',\tau,m)
+ \alpha^{-1} \int_{-\infty}^{\infty} \lambda^m \mathcal{Z}(\lambda,\mathfrak{q}',\tau)\varrho(\lambda)\,d\lambda \\
=
\mathcal{W}\!\left(
\varsigma_m^{\dr}(\lambda) \cdot \big(
\mathbbm{1}\{\mathfrak{q}>\alpha r\}
-\mathbbm{1}\{\mathfrak{q}'>\alpha r+\tau \ve(\lambda)\}
\big)
\right).
\end{multline}
\end{lem}

\begin{proof}
We will show \eqref{eqn:ttfe}, since the proof of \eqref{eqn:ttfe2} is very similar (and formally follows from \eqref{eqn:ttfe} if one could set $G(x) = \mathbbm{1}\{x < \mathfrak{q}'\}$ and $F= \varsigma_m$ there). To that end, we abbreviate $\mathcal{Z}(\lambda, \mathfrak{q}) = \mathcal{Z}(\lambda, \mathfrak{q}, \tau)$ in what follows. By \eqref{eqn:Z_X_relation}; the fact that $\varrho (x) = \alpha \theta \cdot \varsigma_0^{\dr} (x) \cdot \varrho_{\beta} (x)$ (by Lemma \ref{rho0}); and the fact that $\mathbf{D} \coloneqq (1 - \theta \mathbf{T} \bm{\varrho}_{\beta})^{-1}$  is self-adjoint on $L^2 (\varrho_{\beta} d\lambda)$, we have for any $\mathfrak{q} \in \mathbb{R}$ that 
\begin{flalign}\label{eqn:Z_contribution_simplified}
\begin{aligned}
      \int_{-\infty}^{\infty} F(\lambda) \mathcal{Z}(\lambda, \mathfrak{q}) \varrho(\lambda)  d \lambda 
    &=  \alpha \theta \int_{-\infty}^{\infty} F(\lambda)  \cdot \mathbf{D}_{\Lambda}  \mathfrak{X} \big(\Lambda, (\mathfrak{q}- \tau \ve(\Lambda)) / \alpha \big)|_{\Lambda=\lambda} \cdot \varrho_{\beta}(\lambda)  d \lambda \\
  &=  \alpha \theta \int_{-\infty}^{\infty} F^{\dr}(\lambda) \cdot \mathfrak{X} \big(\lambda, (\mathfrak{q}- \tau \ve(\lambda))/ \alpha \big)  \cdot \varrho_{\beta} (\lambda) d \lambda  
\end{aligned}
\end{flalign}
Next, we simplify \eqref{eqn:Z_contribution_simplified} using the equality $\mathfrak{X}(\Lambda, \kappa, \tau)=\mathcal{W}^{\dr}(\psi_{\Lambda, \kappa,\tau})$ (Definition \ref{def:Z_intro}). Letting $g(\lambda) =F^{\dr}(\lambda)$, and using the definition \eqref{eqn:logphi_intro} of $\psi_{\Lambda, \kappa, \tau}$ and the fact that $\mathcal{W}^{\dr}$ is linear, we obtain
\begin{align*}
    \int_{-\infty}^{\infty} & g(\lambda) \cdot \mathcal{W}^{\dr}(\psi_{\lambda, (\mathfrak{q}- \tau \ve(\lambda))/\alpha, \tau}) \cdot \varrho_{\beta}(\lambda)  d \lambda \\
    &= \mathcal{W}^{\dr}\left( \int_{-\infty}^{\infty} g(\lambda) \cdot  \psi_{\lambda, (\mathfrak{q}- \tau \ve(\lambda))/\alpha, \tau} \cdot  \varrho_{\beta}(\lambda) d \lambda \right)\\
     &= \mathcal{W}_{\tilde r, \tilde{\lambda}}^{\dr}\bigg( 2 \int_{-\infty}^{\infty} g(\lambda) \cdot \varrho_{\beta} (\lambda)  \log|\lambda - \tilde{\lambda}| \\
     & \qquad \qquad \quad \times \big( \mathbbm{1}\{ \mathfrak{q}  -\alpha  \tilde r - \tau \ve(\lambda)  > 0  \}    -   \mathbbm{1}\{  \mathfrak{q} - \alpha \tilde r   -  \tau \ve(\tilde{\lambda})  > 0\} \big) d \lambda \bigg),
\end{align*}

\noindent where, in the third statement, we have indicated the variables of the function space acted on by $\mathcal{W}^{\dr}$ as subscripts; the argument of $\mathcal{W}_{\tilde r, \tilde{\lambda}}^{\dr}$ is a function of $(\tilde{r}, \tilde{\lambda})$. Inserting this into \eqref{eqn:Z_contribution_simplified}; and integrating over $\mathfrak{q}$, we obtain 
\begin{flalign}\label{eqn:sts_Z_term_final}
\begin{aligned}
\alpha^{-1} & \int_{-\infty}^{\infty} \int_{-\infty}^{\infty} G'(\mathfrak{q} ) \mathcal{Z}(\lambda, \mathfrak{q} ,\tau) F(\lambda) \varrho (\lambda) d \lambda d \mathfrak{q}  =  \\
& =   \mathcal{W}_{\tilde r, \tilde{\lambda}}^{\dr}\bigg( 2 \theta \int_{-\infty}^{\infty}\int_{-\infty}^{\infty}G'( \mathfrak{q}) \cdot  g(\lambda) \cdot \varrho_{\beta}(\lambda) \log|\lambda - \tilde{\lambda}|\\
    & \qquad \qquad \qquad \times \big( \mathbbm{1}\{ \mathfrak{q}  -\alpha  \tilde r - \tau \ve(\lambda)  > 0  \} 
    - 
    \mathbbm{1}\{  \mathfrak{q} - \alpha \tilde r   -  \tau \ve(\tilde{\lambda})  > 0\} \big) d \lambda d \mathfrak{q} \bigg) \\
    &= \mathcal{W}_{\tilde r, \tilde{\lambda}}^{\dr}\bigg( 2 \theta \int_{-\infty}^{\infty} g(\lambda) \cdot  \varrho_{\beta}(\lambda)  \log|\lambda - \tilde{\lambda}| \cdot \big(
    G( \alpha \tilde r   +  \tau \ve(\tilde{\lambda}))  - G( \alpha  \tilde r + \tau \ve(\lambda))\big) d \lambda \bigg) \\
   &=  \mathcal{W}_{\tilde r, \tilde{\lambda}}^{\dr}\Big(  \theta G( \alpha  \tilde r   +  \theta \tau \ve(\tilde{\lambda})) \cdot  \mathbf{T}  \boldsymbol{\varrho}_{\beta} g
   -\theta \mathbf{T} \boldsymbol{\varrho}_{\beta} \bm{g}  G( \alpha  \tilde r + \tau \ve(\cdot))
    \Big) .
\end{aligned}
\end{flalign}

 \noindent Further observe that, since $g = F^{\dr}$, we have  $\mathcal{W}^{\dr}( \tilde{G} \cdot F ) = \mathcal{W}^{\dr}( \tilde{G} \cdot (1- \theta \mathbf{T}\boldsymbol{\varrho}_{\beta} ) g  )$. Together with \eqref{eqn:sts_Z_term_final} (and the fact that $\tilde{G}(\tilde{r},\tilde{\lambda}) = G(\alpha  \tilde{r} +  \tau \ve(\tilde{\lambda}))$), this implies that the left side of \eqref{eqn:ttfe} is equal to
\begin{align*}
\mathcal{W}^{\dr}( \tilde{G} \cdot F ) &+ \mathcal{W}_{\tilde r, \tilde{\lambda}}^{\dr}\Big(  G( \alpha  \tilde r   +   \tau \ve(\tilde{\lambda})) \cdot \theta \mathbf{T} \boldsymbol{\varrho}_{\beta} [g](\tilde{\lambda})
   -\theta \mathbf{T}\boldsymbol{\varrho}_{\beta}\left[  g \cdot  G( \alpha  \tilde r + \tau \ve(\cdot)) \right](\tilde{\lambda})
    \Big) \\
    &= \mathcal{W}^{\dr}\big( (1-\theta \mathbf{T} \boldsymbol{\varrho}_{\beta})[g \cdot \tilde{G}]\big) = \mathcal{W}(\tilde{G} \cdot F^{\dr}).
    \end{align*}

    \noindent which confirms \eqref{eqn:time_t_finalexp} and thus the lemma. 
\end{proof}

\subsection{Current fluctuations for the Toda lattice on $\llbracket N_1, N_2 \rrbracket$} 
\label{subsec:curr_fluct}

In this section we show the following theorem, which is the counterpart of Theorem \ref{thm:current_fluct_intro} for the Toda lattice on $\llbracket N_1, N_2 \rrbracket$, as it provides the scaling limit for the integrated currents $J_t^{[m]}(q,q')$ of this Toda lattice on a finite (but large) interval.

\begin{thm}[Current fluctuations on $\llbracket N_1, N_2 \rrbracket$] \label{thm:current_fluct}

There is a constant $\mathfrak{c}>0$ so that the following holds. Let $k \in \llbracket 1, N \rrbracket$ be an integer; for each $i \in \llbracket 1, k \rrbracket$, fix an integer $m_i \in \llbracket 0, (\log N)^{1/10} \rrbracket$ and real numbers $\mathfrak{q}_i, \mathfrak{q}_i' \in [-\log N, \log N]$ and $\tau_i \in [0, \log N]$. There exists a coupling between the white noise $\mathcal{W}$ and $\mathbf{L}(0)$ such that, with overwhelming probability, we have 
    \begin{multline}\label{eqn:fin_N_curr_thm}
       \displaystyle\max_{i \in \llbracket 1, k \rrbracket} \Bigg| T^{-1/2} \big(J_{T \tau_i}^{[m_i]}(T \mathfrak{q}_i, T \mathfrak{q}_i') - \mathbb{E} \big[J_{T \tau_i}^{[m_i]}(T \mathfrak{q}_i, T \mathfrak{q}_i') \big] \big) \\
       -  \mathcal{W}\!\Big(
\varsigma_{m_i}^{\dr}(\lambda) \cdot \big(
\mathbbm{1}\{\mathfrak{q}_i>\alpha r\}
-\mathbbm{1}\{\mathfrak{q}_i'>\alpha r+\tau_i \ve(\lambda)\} \big)
\Big) \Bigg| \leq T^{-\mathfrak{c}}.
    \end{multline}

\end{thm}

\begin{proof}

First note by Lemma \ref{prop:spatial_int_current} that, with overwhelming probability, for all $i \in \llbracket 1, k \rrbracket$ we have 
\begin{equation}\label{eqn:fin_N_curr_first_couple}
\Bigg| J_{T \tau_i}^{[m_i]}(T \mathfrak{q}_i,T \mathfrak{q}_i') -   \bigg( -\sum_{j : Q_j(T \tau_i)  < T \mathfrak{q}_i'} \lambda_{j}^{m_i} + \sum_{j : Q_j(0)  < T \mathfrak{q}_i} \lambda_{j}^{m_i} \bigg) \Bigg| \leq  (3 \log N)^{m_i+10}.
\end{equation}

\noindent Restrict to the event on which \eqref{eqn:fin_N_curr_first_couple} and Lemma \ref{lem:lem_chireplacement} hold for each choice of $(q,q',t,m)= (T \mathfrak{q}_i,T \mathfrak{q}_i', T \tau_i, m_i )$ with $i \in \llbracket 1, k \rrbracket$. Then it suffices to couple $(\mathfrak{Y}, \mathfrak{Z})$ and $\mathbf{L}(0)$ so that, with overwhelming probability, \eqref{eqn:fin_N_curr_thm} holds, with $T^{-1/2} (J_{T\tau_i}^{[m_i]} (T\mathfrak{q}_i, T\mathfrak{q}_i') - \mathbb{E}[\cdots])$ there replaced here by 
\begin{equation}\label{eqn:tocompute}
\sum_{ j =1}^N \lambda_j^m \cdot \big( \chi(q- Q_j(0)) - \chi(q'-Q_j(t)) \big) - \mathbb{E}[\cdots ].
\end{equation}

Further restrict to the event $\mathsf{E}_1 = \bigcap_{t \geq 0} \mathsf{BND}_{\L(t)}(\log N)$, which has overwhelming probability by Lemma \ref{lem:bd_lem}. Also restrict to the event $\mathsf{E}_2$ on which Lemma \ref{lem:q_spacing} (with $T = 0$), Lemma \ref{lem:second_der_bd}, and Lemma \ref{lem:loc_cent_dif} all hold. Define, for $j \in \llbracket 1, N \rrbracket$,
\begin{equation}\label{eqn:Zjdf}
    Z_j(t) \coloneqq T^{-1/2} \cdot \big( Q_j(t)-Q_0(t) - t \ve(\lambda_j)\big),
\end{equation}

\noindent observing the slightly different indexing between this definition and the one in \eqref{eqn:quasi_part_fluct_intro}, where instead $j \in \llbracket N_1, N_2 \rrbracket$. Additionally restrict to the event $\mathsf{E}_3$ on which $|Z_j(t)| \leq (\log N)^C$ holds, for all $j$ with $N_1 + T (\log N)^7 \le \varphi_0 (j) \le N_2 - T (\log N)^7$; this is overwhelmingly probable by Lemma \ref{thm:Z_apriori}.

Let $i \in \llbracket 1, k \rrbracket$, and set $(q,q',t,m)= (T \mathfrak{q}_i,T \mathfrak{q}_i', T \tau_i, m_i )$. By Taylor expanding, the sum in \eqref{eqn:tocompute} can be rewritten as
\begin{flalign}\label{eqn:main_est_0}
\begin{aligned}
\sum_{j=1}^{N} & \lambda_j^m \cdot \big( \chi(q-Q_j(0)) -  \chi(q'-Q_j(t)) \big) \\
&=   \sum_{j=1}^{N} \lambda_j^m \cdot \big( \chi(q-Q_j(0)) - \chi(q'-Q_j(0)- t \ve(\lambda_j) - T^{1/2} Z_j(t)) \big) \\
&=   \sum_{j=1}^{N} \lambda_j^m \cdot \bigg(  \chi(q-Q_j(0)) -  \chi(q'-Q_j(0)- t \ve(\lambda_j)) \\
 &\qquad \qquad \qquad +  T^{1/2} \cdot Z_j(t) \cdot \chi'(q'-Q_j(0)- t \ve(\lambda_j)) - \frac{1}{2} \cdot (T^{1/2} Z_j(t))^2 \cdot \chi''(\xi_j) \bigg) ,
 \end{aligned}
\end{flalign}
 where $\xi_j$ is a point between $q'-Q_j(0)-t\,\ve(\lambda_j)$ and $q'-Q_j(0)-t\,\ve(\lambda_j)-T^{1/2}Z_j(t)$, for each $j$. In addition, 
$$
    \sum_{j=1}^{N}\Big| \lambda_j^m (T^{1/2} Z_j(t))^2\chi''(\xi_j) \Big| \leq (\log N)^{m+C} T \mathfrak{M}^{-1},
$$
since $|Z_j| \leq (\log N)^C$ for each $j$ in the support of the sum (by our restriction to $\mathsf{E}_3$); there are at most $\mathfrak{M} (\log N)^6$ nonzero terms in the sum involving second derivatives (using our restriction to $\mathsf{E}_2$, entirely analogously to argument above \eqref{eqn:chidifs1}); $|\lambda_j|^m \leq (\log N)^m$ for all $j \in \llbracket 1, N \rrbracket$ (due to our restriction to $\mathsf{E}_1$);  and $|\chi''|\leq 10 \mathfrak{M}^{-2}$ (by \eqref{estimatesm0}). Together with \eqref{eqn:main_est_0}, this implies
\begin{multline}\label{eqn:main_est_1}
  T^{-1/2} \cdot \Bigg|   \sum_{j=1}^{N} \lambda_j^m \cdot \big( \chi(q-Q_j(0)) - \chi(q'-Q_j(t)) \big) -   \sum_{j=1}^{N} \Bigg( \lambda_j^m \cdot \big(\chi(q-Q_j(0)) -  \chi(q'-Q_j(0)- t \ve(\lambda_j)) \big) \\
+ \lambda_j^m \cdot T^{1/2} \cdot  Z_j(t) \cdot \chi'(q'-Q_j(0)- t \ve(\lambda_j)) \Bigg) \Bigg|  \leq (\log N)^m T^{-c}.
\end{multline}

\noindent We further restrict to $\mathsf{E}_4$, the event provided by Lemma \ref{lem:Z_sum_conc_general} (with the $(F(\lambda), A)$ there equal to $(\lambda^m, (\log N)^m)$ here). By the bound \eqref{eqn:Z_sum_conc_general1}, we have
\begin{multline}\label{eqn:main_est_1b}
   \Bigg| \sum_{j=1}^N \lambda_j^m \cdot Z_j(t) \cdot  \chi'(q'-Q_j(0)- t \ve(\lambda_j))  
    - \alpha^{-1} \int_{-\infty}^{\infty}  \lambda^m \mathcal{Z}(\lambda,q'/T,t/T)   \varrho(\lambda) d\lambda    \Bigg| \leq T^{-c}.
\end{multline}

Combining \eqref{eqn:main_est_1} and \eqref{eqn:main_est_1b}, we deduce
\begin{multline}\label{eqn:main_est_final}
  T^{-1/2} \cdot \Bigg|   \sum_{j=1}^{N} \lambda_j^m \cdot \big( \chi(q-Q_j(0)) - \chi(q'-Q_j(t)) \big)  - \alpha^{-1} \int_{-\infty}^{\infty}  \lambda^m \mathcal{Z}(\lambda,q'/T,t/T)   \varrho(\lambda) d\lambda  \\ 
  - \sum_{j=1}^{N} \lambda_j^m \cdot \big(\chi(q-Q_j(0)) -  \chi(q'-Q_j(0)- t \ve(\lambda_j)) \big) \Bigg| 
   \leq T^{-c}.
\end{multline}

\noindent Taking expectations in \eqref{eqn:main_est_final} and using the fact that $\mathbb{E}[\mathcal{Z}(\lambda,q'/T,t/T)] = 0$ (by \eqref{eqn:intro_Z}, since $\mathfrak{X}$ is centered by \eqref{eqn:Xdef}), we obtain (entirely analogously to in the derivation of \eqref{expectationH}, using Lemma \ref{lem:bd_lem}, as well as \eqref{eqn:intro_Z}, Lemma \ref{lem:Xbd1}, and Lemma \ref{lem:ve_tail} to bound $\mathcal{Z}$ off of the overwhelmingly probable event on which \eqref{eqn:main_est_final} holds) that
\begin{multline}\label{eqn:main_est_finalexp}
\Bigg| \mathbb{E}  \bigg[   \sum_{j=1}^{N} \lambda_j^m \cdot \big( \chi(q-Q_j(0)) -  \chi(q'-Q_j(t)) \big) \\
   -  \sum_{j=1}^{N} \lambda_j^m \cdot \big(\chi(q-Q_j(0)) -  \chi(q'-Q_j(0)- t \ve(\lambda_j)) \big)   \bigg] \Bigg|
   \leq T^{-c}.
\end{multline}

By \eqref{eqn:fin_N_curr_first_couple}, Lemma \ref{lem:lem_chireplacement}, \eqref{eqn:main_est_final}, \eqref{eqn:main_est_finalexp}, and \eqref{eqn:ttfe2} from Lemma  \ref{lem:ttfe} it suffices to couple $(\mathfrak{Y},\mathcal{Z})$ and $\mathbf{L}(0)$ so that, with overwhelming probability, we have (recalling $\Xi^{[m]}$ from Definition \ref{eqn:xim_outline}) that
 \begin{flalign*}
       \displaystyle\max_{i \in \llbracket 1, k \rrbracket} \big| \Xi^{[m_i]} (T\mathfrak{q}_i, T\mathfrak{q}_i', T\tau_i)   -  \mathfrak{Y}(\mathfrak{q}_i,\mathfrak{q}_i', \tau_i, m_i)  \big| \leq T^{-c}.
    \end{flalign*}

    \noindent This follows from Theorem \ref{thm:couplethm12}; the estimates \eqref{eqn:xim2_couple} and \eqref{eqn:ximbr_couple}; and the definitions \eqref{eqn:Y2out}, \eqref{eqn:Y1out}, and \eqref{yprocess}. 
\end{proof}

\subsection{Proof of Theorem \ref{thm:current_fluct_intro}}
\label{subsec:cf_proof}

 In this section we prove Theorem \ref{thm:current_fluct_intro}. Using Theorem \ref{thm:current_fluct}, this will essentially follow from the coupling between the Toda lattices on $\llbracket N_1, N_2 \rrbracket$ and $\mathbb{Z}$, provided by Lemma \ref{lem:inf_vol_bds}.

\begin{proof}[Proof of Theorem \ref{thm:current_fluct_intro}]

First, we set some notation. Denote $\mathbf{a}(t) = (a_i(t))_{i \in \mathbb{Z}}$ and $\mathbf{b}(t) = (b_i(t))_{i \in \mathbb{Z}}$, and recall that $(\mathbf{a}(0);\mathbf{b}(0))$ is sampled under thermal equilibrium $\mu_{\beta,\theta;\infty}$ from \eqref{eqn:equil}. Let $T \ge 1$ be a real number; set $N=\floor{T^{100}}$; and let $N_1 \le N_2$ be integers satisfying Assumption \ref{ass:NT_assumption}. Let $(\tilde{\mathbf{a}}(t);\tilde{\mathbf{b}}(t))$ denote the Toda lattice on $\llbracket N_1, N_2 \rrbracket$, with initial data $(\tilde{\mathbf{a}}(0); \tilde{\mathbf{b}}(0))$ sampled under thermal equilibrium $\mu_{\beta,\theta;N_1,N_2}$, let $(\tilde{\mathbf{p}}(t);\tilde{\mathbf{q}}(t))$ denote the associated Toda state space dynamics. Set $\tilde{\mathbf{a}}(t) = (\tilde{a}_i(t))$, $\tilde{\mathbf{b}}(t) = (\tilde{b}_i(t))$, $\tilde{\mathbf{p}}(t) = (\tilde{p}_i(t))$, and $\tilde{\mathbf{q}}(t) = (\tilde{q}_i (t))$ over $i \in \llbracket N_1, N_2 \rrbracket$. As in Lemma \ref{lem:fin_inf_t_coupling}, couple $(\tilde{\mathbf{a}}(0);\tilde{\mathbf{b}}(0))$ with $(\mathbf{a}(0); \mathbf{b}(0))$ so that $\tilde{a}_i(0) = a_i(0)$ for all $i \in \llbracket N_1, N_2 -1 \rrbracket$, and $\tilde{b}_i (0) = b_i(0)$ for all $i \in \llbracket N_1, N_2 \rrbracket$. We assume throughout that $T$ is sufficiently large so that $k \le \log N$; $\tau_i \in [0, \log N]$; and $\mathfrak{q}_i, \mathfrak{q}_i' \in [- \log N, \log N]$, for each $i \in \llbracket 1, k \rrbracket$. 

    Let $\mathfrak{k}_i^{[m]} (t)$ and $J_t^{[m]} (q,q')$ denote the local charges and spatial integrated currents (Definitions \ref{kti00} and \ref{def:integrated_curr_spatial}) associated with the Toda lattice $(\mathbf{a}(t); \mathbf{b}(t))$, respectively. Similarly let $\tilde{\mathfrak{k}}_i^{[m]} (t)$ and $\tilde{J}_t^{[m]} (q,q')$ denote the local charges and spatial integrated currents associated with the Toda lattice $(\tilde{\mathbf{a}}(t); \tilde{\mathbf{b}}(t))$ on $\llbracket N_1, N_2 \rrbracket$, respectively. Also denote by $\tilde{\L}(t)$ the Lax matrix (Definition \ref{lt}) associated with the latter Toda lattice $(\tilde{\mathbf{a}}(t); \tilde{\mathbf{b}}(t))$. Below, we couple the Gaussian processes $(\mathfrak{Y}, \mathcal{Z})$ (Definitions \ref{def:Z_intro} and \ref{x2x0y2y0}) with $(\tilde{\mathbf{a}}(0);\tilde{\mathbf{b}}(0))$ so that Theorem \ref{thm:current_fluct} holds, with the $J$ there replaced by $\tilde{J}$ here. 

Let $\mathsf{E}_1 \coloneqq \bigcap_{t \geq 0} \mathsf{BND}_{\tilde{\L}(t)}( (\log N)/2)$, which overwhelmingly probable by Lemma \ref{lem:bd_lem}. Also let $\mathsf{E}_2$ denote the event on which Lemma \ref{lem:q_spacing} holds (for $(\tilde{\mathbf{a}}(t); \tilde{\mathbf{b}}(t))$), with the $T$ there both equal to $0$ and $T \log N$ here. Further let $\mathsf{E}_3$ denote the event on which the second part of Lemma \ref{lem:fin_inf_t_coupling} holds; $\mathsf{E}_4$ denote that on which Theorem \ref{thm:current_fluct} holds; and $\mathsf{E}_5$ denote that on which Lemma \ref{lem:inf_vol_bds} holds, with $(T, K)$ there given by $(T \log N, T^2)$ here. Restrict to $\bigcap_{i=1}^5 \mathsf{E}_i $, which is overwhelmingly probable. 

Denote $I_0 = -\lfloor 2 T^2 \rfloor$. Observe by Definition \ref{def:integrated_curr_spatial} (with the $k$ there equal to $I_0$ herer) and Lemma \ref{lem:inf_vol_bds} that, for any $(q,q', t, m) \in [-T \log N, T \log N]^2 \times  [0, T \log N] \times \llbracket 0, (\log N)^{1/10} \rrbracket$, the sum \eqref{eqn:inf_integrated_curr} defining the integrated current in the Toda lattice $(\mathbf{a}(t);\mathbf{b}(t))$ on $\mathbb{Z}$ is given by 
\begin{equation}\label{eqn:Jtilde_approx}
         J_t^{[m]}(q,q')=  \sum_{ i \geq I_0  : q_i(0) < q}  \mathfrak{k}_i^{[m]}(0) -\sum_{ i \geq I_0 : q_i(t) < q}  \mathfrak{k}_i^{[m]}(t)+ \int_0^t \mathfrak{j}_{I_0}^{[m]}(s) ds.
\end{equation}

\noindent Moreover, by Remark \ref{km00}, these currents for the Toda lattice $(\tilde{\mathbf{a}}(t);\tilde{\mathbf{b}}(t))$ on $\llbracket N_1, N_2 \rrbracket$ are given by 
\begin{equation}\label{eqn:J_approx}
        \tilde{J}_t^{[m]}(q,q')=  \sum_{ i \ge I_0  : q_i(0) < q}  \tilde{\mathfrak{k}}_i^{[m]}(0) -\sum_{ i > I_0 : q_i(t) < q'}  \tilde{\mathfrak{k}}_i^{[m]}(t) + \int_0^t \tilde{\mathfrak{j}}_{I_0}^{[m]}(s) ds.
\end{equation}

\noindent Comparing \eqref{eqn:Jtilde_approx} and \eqref{eqn:J_approx}; using \eqref{mki}, \eqref{jti}, \eqref{eqn:fin_inf_dif}, and Lemma \ref{lem:inf_vol_bds} (the latter two of which hold on $\mathsf{E}_3 \cap \mathsf{E}_5$); and using the fact that $\mathfrak{k}_i^{[m]}(s)$ and $\mathfrak{j}_i^{[m]}(s)$ (or $\tilde{\mathfrak{k}}_i^{[m]} (s)$ and $\tilde{\mathfrak{j}}_i^{[m]} (s)$) both only depend on the random variables $(a_j(s), b_j(s))$ (or on the $(\tilde{a}_j(s), \tilde{b}_j(s))$, respectively) for $j \in \llbracket i-m-1, i+m+1 \rrbracket$, we have for any $i \in \llbracket 1, k \rrbracket$ that  
\begin{multline}\label{eqn:finite_inf_couple}
   \left| \tilde J_{T \tau_i}^{[m_i]}(T \mathfrak{q}_i, T \mathfrak{q}_i') -  J_{T \tau_i}^{[m_i]}(T \mathfrak{q}_i, T \mathfrak{q}_i') \right|   \leq  2 (4 T^2+1 + T \log N) \cdot m_i (3 \log N)^{m_i} e^{-T^2/5} \leq e^{-T^2/10}.
\end{multline}

By \eqref{eqn:fin_N_curr_thm} (which holds on $\mathsf{E}_4$), to establish the theorem it suffices to show that, with overwhelming probability,
\begin{multline}\label{eqn:inf_vol_curr_coupling}
   \displaystyle\max_{i \in \llbracket 1, k \rrbracket} \Bigg| T^{-1/2} \left(\tilde J_{T \tau_i}^{[m_i]}(T \mathfrak{q}_i, T \mathfrak{q}_i') - \mathbb{E}[\tilde J_{T \tau_i}^{[m_i]}(T \mathfrak{q}_i, T \mathfrak{q}_i')] \right)  \\
    - T^{-1/2} \left(J_{T \tau_i}^{[m_i]}(T \mathfrak{q}_i, T \mathfrak{q}_i') - \mathbb{E}[J_{T \tau_i}^{[m_i]}(T \mathfrak{q}_i, T \mathfrak{q}_i')] \right)  \Bigg| \leq T^{-c},
\end{multline}

\noindent to which end it suffices by \eqref{eqn:finite_inf_couple} to show that, with overwhelming probability, 
\begin{equation}\label{eqn:expect_bd}
   \displaystyle\max_{i \in \llbracket 1, k \rrbracket} T^{-1/2} \cdot \Big|\mathbb{E}\big[\tilde J_{T \tau_i}^{[m_i]}(T \mathfrak{q}_i, T \mathfrak{q}_i') \big] 
    - \mathbb{E} \big[J_{T \tau_i}^{[m_i]}(T \mathfrak{q}_i, T \mathfrak{q}_i') \big]\Big| \leq T^{-c}.
\end{equation}

The proof of \eqref{eqn:expect_bd} will be similar to the derivation of \eqref{expectationH}. Specifically, since \eqref{eqn:finite_inf_couple} holds with overwhelming probability, it suffices to show that 
\begin{equation}\label{eqn:Echarge_N}
\mathbb{E} \big[ \big|\tilde{\mathfrak{k}}_h^{[m_i]}(T \tau_i) \big|^2 \big] +\big|\tilde{\mathfrak{j}}_h^{[m_i]}(T \tau_i) \big|^2 \big] \leq T^{10}; \qquad  \mathbb{E}\big[J_{T \tau_i}^{[m_i]}(T \mathfrak{q}_i, T \mathfrak{q}_i')^2\big] \leq T^{10},
\end{equation}

\noindent for all $h \in \llbracket N_1, N_2 \rrbracket$ and $i \in \llbracket 1, k \rrbracket$. The first bound in \eqref{eqn:Echarge_N} is proven entirely analogously to in \eqref{expectationH} (using \eqref{mki}, \eqref{jti}, and Lemma \ref{lem:bd_lem}), so let us only show the second bound there. 

To do so, for each integer $r \geq 1$, let $\mathsf{A}_r$ denote the event on which Lemma \ref{lem:inf_vol_bds} holds, with the $(T, K)$ there given by $(T \log N, K_r)$ here, where $K_r = T^{2r}$; also let $\mathsf{A}_0 = \emptyset$. By the Borel--Cantelli lemma, we have that $\mathbb{P}[\bigcup_{r =1}^{\infty} \mathsf{A}_r] = 1$. Moreover, for each $s \in [0, T \log N]$ and $m \geq 0$, we have by \eqref{mki}, \eqref{jti} that
\begin{flalign}\label{eqn:Echarge_inf}
\begin{aligned}
\mathbbm{1}_{\mathsf{A}_r} \cdot \big|\mathfrak{k}_j^{[m]}(s) \big| &\leq  (C \log (T^{2^{r}}))^{m}, \qquad \quad \text{for all } j \in \llbracket -T^2, T^2 \rrbracket, \\
\qquad \mathbbm{1}_{\mathsf{A}_r} \cdot \big| \mathfrak{j}_j^{[m]}(s) \big| &\leq  (C \log (T^{2^{r}}))^{m+1}, \qquad \text{for all } j \in \llbracket -T^2, T^2 \rrbracket.
\end{aligned}
\end{flalign}

\noindent Thus,
\begin{align*}
   \mathbb{E}\left[J_{T \tau_i}^{[m_i]}(T \mathfrak{q}_i, T \mathfrak{q}_i')^2\right] &\leq \mathbb{E}\Bigg[\sum_{r \geq 1} \mathbbm{1}_{\mathsf{A}_r \cap \mathsf{A}_{r-1}^{\complement}} \cdot J_{T \tau_i}^{[m_i]}(T \mathfrak{q}_i, T \mathfrak{q}_i')^2\Bigg] \\
   &\leq C \sum_{r=1}^{\infty}  \big( K_r^2 \cdot   (C \log K_r)^{2 m_i} + (T \log T)^2 (C \log K_r) )^{2 m_i+2} \big)\cdot \mathbb{P} \big[\mathsf{A}_{r-1}^{\complement} \big] \\
   &\leq 2C  \big( K_1^2 \cdot   (C \log K_1)^{2 m_i} + (T \log T)^2 (C \log K_1) )^{2 m_i+2} \big) \le T^{10},
\end{align*}

\noindent where we used that $\mathbb{P}[\mathsf{A}_{r-1}^{\complement}] \le e^{-c(\log K_{r})^2}$ by Lemma \ref{lem:inf_vol_bds} (and that $K_1 = T^2$ and $m_i \le (\log N)^{1/10}$). This shows \eqref{eqn:expect_bd}, verifying \eqref{eqn:inf_vol_curr_coupling}, which (as mentioned above) establishes the theorem.
\end{proof}

\section{Fluctuations of local charges}
\label{sec:charge_fluct}

In this section we prove Corollary~\ref{cor:twopoint_intro}. To show Corollary \ref{cor:twopoint_intro}, we will first prove Theorem \ref{thm:euler_flucts} below, which proves a joint scaling limit for certain linear functionals of local charges at fixed times. Throughout we assume $\theta<\theta_0(\beta)$ is sufficiently small and we adopt Assumption~\ref{ass:NT_assumption} when applicable. We use the local charges $\mathfrak{k}_j^{[m]}(t)$ and $\mathfrak{K}_i^{[m]}(t)$ from Definition~\ref{kti00} and the two-point function $S_{m,n}(j,t)$ from \eqref{eqn:two_pt_def}. We also use the quantities $Z(\Lambda,Q,t)$ from \eqref{eqn:Z_outline} and $Z_j^{\mathcal Q}(t)$ from \eqref{eqn:quasi_part_fluct_intro}, the white noise $\mathcal{W}$ and dressed white noise $\mathcal{W}^{\dr}$ from Definitions~\ref{def:white_noise_dress} and \ref{def:Wdress}, and the processes $\mathfrak{X}$ and $\mathcal{Z}$ from Definition~\ref{def:Z_intro}. Finally, we use the notation of Definition~\ref{def:quasi_intro}. We denote the dressing operator (Definition~\ref{def:dress_intro}) by
\[
\mathbf{D}:=(1-\theta \mathbf{T}\boldsymbol{\varrho}_{\beta})^{-1}:\mathcal{H}\to\mathcal{H},
\]
and write $\mathbf{D}_\Lambda F(\Lambda,\lambda,q)$ for the action of $\mathbf{D}$ on the function $\Lambda\mapsto F(\Lambda,\lambda,q)$ (when defined).

\subsection{Functionals of charges through eigenvalues}
\label{subsec:fininfcouple}

The next theorem identifies the limit of the Toda lattice local charges, on the Euler scale, under thermal equilibrium. Its proof is in Section \ref{subsec:eulerfluctproof} below. In what follows, we recall the operator $\mathbf{F}$ from Definition \ref{operatorf}; and the function $\varsigma_n(x) = x^n$ from \eqref{nx}. 

\begin{theorem}[Euler scale fluctuations]\label{thm:euler_flucts}
	For any $\beta > 0$, there exists a constant $\theta_0 (\beta)>0$ such that the following holds for any fixed $\theta \in (0, \theta_0)$. Let $(\mathbf{a}(t);\mathbf{b}(t))$ denote the Toda lattice with initial data $(\mathbf{a}(0);\mathbf{b}(0))$, sampled from the thermal equilibrium $\mu_{\beta,\theta;\infty}$. Further fix an integer $k \ge 1$ and, for each $i \in \llbracket 1, k \rrbracket$, fix an integer $m_i \ge 0$, a real number $\tau_i \ge 0$, and a smooth function $f_i : \mathbb{R} \rightarrow \mathbb{R}$ of compact support. For $i \in \llbracket 1, k \rrbracket$, set
\begin{align*}
   \mathfrak{g}_i =\begin{cases}
       \mathbf{F} \varsigma_{m_i} &  \text{ if } m_i > 0 \\
      - \alpha \cdot \varsigma_{0}^{\dr}  &  \text{ if } m_i = 0 .
   \end{cases}  
\end{align*}
    
   \noindent Then, as $T$ tends to $\infty$, the $k$-tuple
\begin{flalign}
    \label{sumk} 
  \bigg[ T^{-1/2}  \sum_{j \in \mathbb{Z} } f_i(jT^{-1}) \cdot \big(  \mathfrak{K}_j^{[m_i]}(T \tau_i) - \mathbb{E}[ \mathfrak{K}_j^{[m_i]}(T \tau_i) ] \big) \bigg]_{i \in \llbracket 1, k \rrbracket}  
\end{flalign} 

    \noindent converges in law to the $k$-tuple
    \begin{flalign*} 
     \left[ \mathcal{W}\left(  f_i(r + \alpha^{-1} \tau_i \ve(\lambda)) \cdot \mathfrak{g}_i (\lambda) \right) \right]_{i \in \llbracket 1, k \rrbracket}.
\end{flalign*}
\end{theorem}

We will proceed by coupling the Toda lattice under thermal equilibrium $\mu_{\beta,\theta;\infty}$ to the Toda lattice on $\llbracket N_1, N_2 \rrbracket$. Therefore, to prove this theorem, we require the following lemma approximating charges of the Toda lattice on the interval $\llbracket N_1, N_2 \rrbracket$ by the eigenvalues of the associated Lax matrix (indexed as in \eqref{lambdait}). This will later enable us to use our understanding of quasi-particles to analyze these charges.

\begin{lem}[Approximating charges by eigenvalues]\label{lem:euler_to_linstat}
 There exists a constant $\mathfrak{c}>0$ such that the following holds. Let $K \in \llbracket 1, \log N \rrbracket$ be an integer. For each $i\in\llbracket 1,K\rrbracket$, let $m_i \in \llbracket 0, (\log N)^{1/10} \rrbracket$ be an integer; $\tau_i \in [0, \log N]$, $B_i \ge 1$, and $U_i \in [\mathfrak{M}, T \log N ]$ be real numbers; and $G_i : \mathbb{R} \rightarrow \mathbb{R}$ be a compactly supported function satisfying Assumption~\ref{ass:G2}, with the $(B,U,\mathfrak{Q})$ there equal to $( B_i,U_i,0)$ here. Adopt Assumption~\ref{ass:NT_assumption}; denote the associated Toda lattice here by $(\mathbf{a}^{\llbracket N_1, N_2 \rrbracket}(t); \mathbf{b}^{\llbracket N_1, N_2 \rrbracket}(t))$; and denote the local charges (recall \eqref{mki}) of this Toda lattice by $\mathfrak{k}_j^{[m]}$. Then (recalling \eqref{lambdait}), the following holds with overwhelming probability, for every $i\in\llbracket 1,K\rrbracket$. If $m_i > 0$, then 
\begin{multline}\label{eqn:etl}
\Bigg|
T^{-1/2} \sum_{j=N_1}^{N_2} G_i(j) \cdot \big(\mathfrak{k}_j^{[m_i]}(T\tau_i)-\mathbb{E}\big[\mathfrak{k}_j^{[m_i]}(T\tau_i) \big]\big)
\\
-T^{-1/2} \bigg(
\sum_{j=N_1}^{N_2} G_i(j)  \Lambda_j(T \tau_i)^{m_i}
-\mathbb{E}\Bigl[\sum_{j=N_1}^{N_2} G_i(j)\,\Lambda_j(T \tau_i)^{m_i}\Bigr]
\bigg)
\Bigg|
\le  B_i T^{-\mathfrak{c}}.
\end{multline}

 \noindent If instead $m_i = 0$, then
\begin{multline}\label{eqn:etl0}
\Bigg|
T^{-1/2} \sum_{j=N_1}^{N_2} G_i(j) \cdot \big(r_j(T\tau_i)-\mathbb{E} [ r_j(T\tau_i) ]\big)
\\
+\alpha T^{-1/2} \bigg(
\sum_{j=N_1}^{N_2}  G_i(\alpha^{-1} q_j(T \tau_i)) 
-\mathbb{E} \bigg[\sum_{j=N_1}^{N_2} G_i(\alpha^{-1} q_j(T \tau_i)) \bigg]
\bigg)
\Bigg|
\le  B_i T^{-\mathfrak{c}}.
\end{multline}

\end{lem}

\begin{proof}
We first couple $\L(t)$ to a copy of $\L(0)$ so that they are close for indices $i$ in the bulk. Define $(\mathbf{a}(t); \mathbf{b}(t)) = (a_i(t); b_i(0))_{i \in \mathbb{Z}}$ as the Toda lattice on $\mathbb{Z}$ started at thermal equilibrium and coupled to the Toda lattice $(\mathbf{a}^{\llbracket N_1, N_2 \rrbracket}(t); \mathbf{b}^{\llbracket N_1, N_2 \rrbracket}(t))= ((a_i^{\llbracket N_1, N_2 \rrbracket}(t), b_i^{\llbracket N_1, N_2 \rrbracket}(t)))_{i \in  \llbracket N_1, N_2 \rrbracket} $ on $\llbracket N_1, N_2 \rrbracket$ as in Lemma \ref{lem:fin_inf_t_coupling}. As in Definition \ref{lt}, for each $t \in [0, T \log N]$, define the $N \times N$ matrix $\M(t) = [M_{i j}(t)]_{i, j \in \llbracket N_1, N_2 \rrbracket} $ by setting $M_{i i }(t) = b_i(t)$ for each $i \in \llbracket N_1, N_2 \rrbracket$; setting $M_{i,i+1}(t) =  M_{i+1,i}(t) = a_i(t)$ for each $i \in \llbracket N_1, N_2 -1 \rrbracket$; and setting all $M_{ij} (t) = 0$ for $(i,j)$ not of the above form. By Lemma \ref{lem:fin_inf_t_coupling}, $\M(t)$ is equal in law to $\L(0)$ for any $t \geq 0$ and, with overwhelming probability, 
 \begin{equation}\label{eqn:MLijbd}
   \sup_{t \in [0, T \log N]} \max_{i,j \in \llbracket N_1+T^2, N_2-T^2 \rrbracket}  |M_{i j}(t)- L_{i j}(t)| \leq e^{-T^2/5}.
\end{equation}
In what follows, we restrict to the event $\mathsf{E}_1$ on which \eqref{eqn:MLijbd} holds.
 
For any $t > 0$, let $\eig \mathbf{M}(t) = (\lambda_{1;\mathbf{M}(t)}, \lambda_{2;\mathbf{M}(t)}, \ldots , \lambda_{N;\mathbf{M}(t)})$ denote the eigenvalues of $\M(t)$. Recalling from \eqref{eqn:N1N2} that $\zeta \ge e^{-100 (\log N)^{3/2}}$, also let $\varphi_{\M(t)} : \llbracket 1, N \rrbracket \rightarrow \llbracket N_1, N_2 \rrbracket$ denote a $\zeta$-localization center bijection of $\M(t)$, and (in analogy with in Definition \ref{def:quasi_intro}) define 
\begin{equation}\label{eqn:LamiM}
  \Lambda_{i;\M(t)} =  \lambda_{\varphi_{\M(t)}^{-1}(i); \mathbf{M}(t)}.
\end{equation}

 \noindent In addition, let $(p_{i;\mathbf{M}}(t); q_{i;\M}(t))_{i \in \mathbb{Z}}$ denote the Toda state space dynamics associated with the Toda lattice $(\mathbf{a}(t); \mathbf{b}(t))$ on $\mathbb{Z}$. Also set $r_{j;\M}(t) \coloneqq q_{i+1;\M}(t)-q_{i;\M}(t) = -2 \log a_i(t)$, for each $i \in \mathbb{Z}$.

Fix $i \in \llbracket 1, K \rrbracket$, and abbreviate $G = G_i$, $t = T \tau_i$, $m = m_i$, $B = B_i$, $U = U_i$, and $\M = \M(t)$. Also denote $\mathfrak{k}_{j;\M}^{[m]} \coloneqq \left(\mathbf{M}^m\right)_{j j}$ if $m \ge 1$. We will reduce \eqref{eqn:etl} and \eqref{eqn:etl0} to statements about $\M = \M(t)$, instead of $\L(t)$. To that end, we claim that, with overwhelming probability, if $m \ge 1$ we have 
    \begin{align}\label{eqn:cplLwts}
      & \Bigg|    \sum_{j =N_1 }^{N_2} G(j) \cdot \big( \mathfrak{k}_j^{[m]}(t)  - \mathbb{E}[ \mathfrak{k}_j^{[m]}(t)  ] \big)
         -  \sum_{j =N_1 }^{N_2} G(j) \cdot  \big( \mathfrak{k}_{j;\M}^{[m]}  - \mathbb{E}[ \mathfrak{k}_{j;\M}^{[m]}  ] \big) \Bigg|  
         \leq B T^{1/2-c}; \\
         & \Bigg|   \sum_{j =N_1 }^{N_2} G(j)  \Lambda_{j}(t)^{m} - \mathbb{E}\bigg[ \sum_{j =N_1 }^{N_2} G(j)  \Lambda_{j}(t)^{m}\bigg]  \notag \\ 
         &  \qquad \qquad \qquad \qquad- \bigg(  \sum_{j =N_1 }^{N_2} G(j) \Lambda_{j;\M}^{m} - \mathbb{E}\bigg[ \sum_{j =N_1 }^{N_2} G(j)  \Lambda_{j;\M}^{m}\bigg] \bigg) \Bigg|  \leq B T^{1/2-c} , \label{eqn:cplLwts2}
    \end{align}
    and if $m=0$ we have 
    \begin{align}\label{eqn:cplLwts0}
      & \Bigg|    \sum_{j =N_1 }^{N_2} G(j) \cdot \bigl(r_j(t)-\mathbb{E}[r_j(t)]\bigr)
         -  \sum_{j =N_1 }^{N_2} G(j)  \cdot \big( r_{j;\M}(t)  - \mathbb{E}[ r_{j;\M}(t)  ] \big) \Bigg|  \leq B T^{1/2-c}; \\
        & \Bigg|   
            \sum_{j=N_1}^{N_2} G_i(\alpha^{-1} q_j(t)) 
            -\mathbb{E}\biggl[\sum_{j=N_1}^{N_2} G_i(\alpha^{-1} q_j(t)) \biggr]
             \notag \\ 
        &  \qquad \qquad \qquad \qquad - \biggl(
\sum_{j=N_1}^{N_2} G_i(\alpha^{-1} q_{j;\M}(t)) 
    -\mathbb{E}\biggl[\sum_{j=N_1}^{N_2} G_i(\alpha^{-1} q_{j;\M}(t)) \biggr] \biggr)  \Bigg| \leq B T^{1/2-c} . \label{eqn:cplLwts02}
    \end{align}
We will only verify \eqref{eqn:cplLwts} and \eqref{eqn:cplLwts2}, as the proofs of  \eqref{eqn:cplLwts0} and \eqref{eqn:cplLwts02} are similar (and are in fact more direct consequences of \eqref{eqn:MLijbd}).

Let $\mathsf{E}_2 = \mathsf{BND}_{\L(t)}(\log N) \cap \mathsf{BND}_{\M(t)}(\log N) \cap \mathsf{SEP}_{\M(t)}(e^{-(\log N)^2}) $. Let $\mathsf{E}_3  $ be the event given by Lemma \ref{lem:Lax_eig_coupling}, with the $(\L, \tilde{\L}, \delta, \mathcal{D})$ there equal to $(\M(t), \L(t), e^{-T^2/5}), \llbracket N_1 + T^2, N_2 - T^2 \rrbracket)$ here (using \eqref{eqn:MLijbd}). Let $\mathsf{E}_4$ denote the event of Lemma \ref{lem:loc_cent_dif} (for $\L$). Also let $\mathsf{E}_5 = \bigcap_{j=N_1}^{N_2} \{ |r_{j;\M}| \leq (\log N)^C \}$. Below, we restrict to $\bigcap_{i=1}^5 \mathsf{E}_i$, which is overwhelmingly probable by Lemma \ref{lem:bd_lem}, \ref{lem:sep_lem}, and \ref{lem:inf_vol_bds}.

We denote $\psi : \llbracket N_1, N_2 \rrbracket \rightarrow \llbracket N_1, N_2 \rrbracket$ as the injective function $\psi$ given by Lemma \ref{lem:Lax_eig_coupling} (extended from $\llbracket N_1+T^2, N_2-T^2 \rrbracket$ to $\llbracket N_1, N_2 \rrbracket$ arbitrarily to a bijection). Then, we have
    \begin{flalign}\label{eqn:cplML1}
    \begin{aligned}
      \Bigg|    \sum_{j =N_1 }^{N_2} G(j)  \Lambda_j^{m}(t) 
         -  \sum_{j =N_1 }^{N_2} G(j)  \Lambda_{j;\M}^{m} \Bigg|  
     &=   \Bigg|    \sum_{j =N_1}^{N_2} G(j)  \Lambda_j^{m}(t) 
         -  \sum_{j =N_1 }^{N_2} G(\psi(j))  \Lambda_{\psi(j);\M}^{m} \Bigg|\\
         &\leq  \Bigg|    \sum_{j =N_1+\floor{T^3} }^{N_2-\floor{T^3}} G(j)  \Lambda_j^{m}(t) 
         -  \sum_{j =N_1+\floor{T^3} }^{N_2-\floor{T^3}} G(\psi(j))  \Lambda_{\psi(j);\M}^{m} \Bigg|\\
         &\leq   N B e^{-c(\log N)^3} + (\log N)^m  \sum_{j =N_1+\floor{T^3} }^{N_2-\floor{T^3}} |G(j) - G(\psi(j))|   \\
         &\leq  N B e^{-c(\log N)^3}  + B (\log N)^{m+4} \leq B T^{1/2-c},
    \end{aligned}
    \end{flalign}
    
    \noindent where the first statement holds since $\psi$ is a bijection; the second holds by Lemma \ref{lem:Lax_eig_coupling} and the fact that $\psi(\llbracket N_1+\floor{T^3}, N_2-\floor{T^3} \rrbracket) $ contains the support of the second sum by our restriction to $\mathsf{E}_4$ and the assumptions on $G$; the third holds by the bound in the second part of Lemma \ref{lem:Lax_eig_coupling}; the fourth holds by the fact that $|\psi(j) - j|\leq (\log N)^2$, that $|G'(x)| \leq B U^{-1}$ and that fact $\supp G' \subseteq [-U, U]$; and the fifth holds since $m \leq (\log N)^{1/10}$. In addition, we have from \eqref{eqn:MLijbd} that 
     \begin{equation}\label{eqn:cplML2}
      \Bigg|    \sum_{j =N_1 }^{N_2} G(j)  \mathfrak{k}_j^{[m]}(t) 
         -  \sum_{j =N_1 }^{N_2} G(j)  \mathfrak{k}_{j;\M}^{[m]} \Bigg| 
         \leq B N (3 \log N)^m e^{-T^2/10} \leq B e^{-T^2/20} .
    \end{equation}
Denote by $\mathsf{E}$ the overwhelmingly probable event on which \eqref{eqn:cplML1} and \eqref{eqn:cplML2} hold, we obtain using Lemma \ref{lem:bd_lem} (very similarly to \eqref{expectationH}) that
\begin{align*}
   \mathbb{E}\left[ \mathbbm{1}_{\mathsf{E}^{\complement}} \cdot \Bigg|    \sum_{j =N_1 }^{N_2} G(j)  \mathfrak{k}_j^{[m]}(t) 
         -  \sum_{j =N_1 }^{N_2} G(j)  \mathfrak{k}_{j;\M}^{[m]} \Bigg| \right] & \leq B e^{-c_1 (\log N)^2 } \\
         \mathbb{E}\left[ \mathbbm{1}_{\mathsf{E}^{\complement}} \cdot  \Bigg|    \sum_{j =N_1 }^{N_2} G(j)  \Lambda_j^{m}(t) 
         -  \sum_{j =N_1 }^{N_2} G(j)  \Lambda_{j;\M}^{m} \Bigg|   \right] & \leq B e^{-c_1 (\log N)^2 }.
\end{align*}
The display above, together with \eqref{eqn:cplML1} and \eqref{eqn:cplML2}, proves \eqref{eqn:cplLwts} and \eqref{eqn:cplLwts2}.
    
By \eqref{eqn:cplLwts} and \eqref{eqn:cplLwts2}, to prove \eqref{eqn:etl} it suffices to show 
\begin{equation}\label{eqn:cplMwts}
      \Bigg|    \sum_{j =N_1 }^{N_2} G(j)  \mathfrak{k}_{j;\M}^{[m]} - \mathbb{E}[\cdots]
         -  \bigg( \sum_{j =N_1 }^{N_2} G_i(j)  \Lambda_{j;\M}^{m} -\mathbb{E}[\cdots] \bigg)\Bigg|  
         \leq B T^{1/2-c},
    \end{equation}
    and similarly for \eqref{eqn:etl0}. Since $\M(t)$ has the same law as $\L(0)$, we have reduced \eqref{eqn:etl} and \eqref{eqn:etl0}  to showing that the following holds with overwhelming probability. If $m >0$, then
\begin{equation}\label{eqn:cplL0wts}
      \Bigg|    \sum_{j =N_1 }^{N_2} G(j) \cdot \big(\mathfrak{k}_j^{[m]}(0) -\mathbb{E}[\mathfrak{k}_j^{[m]}(0) ]\big)
         -  \bigg( \sum_{j =N_1 }^{N_2} G(j)  \Lambda_{j}^{m}  - \mathbb{E}\bigg[ \sum_{j =N_1 }^{N_2} G(j)  \Lambda_{j}^{m} \bigg] \bigg) \Bigg|  
         \leq B T^{1/2-c},
    \end{equation}
    and if $m = 0$, then (recalling $X_j$ from Definition \ref{eqn:Xidef})
\begin{equation}\label{eqn:cplL0wtsm0}
     \Bigg| T^{-1/2} \sum_{j=N_1}^{N_2} G_i(j)X_j
+\alpha T^{-1/2}\biggl(
\sum_{j=N_1}^{N_2} G_i(\alpha^{-1} q_j(0)) 
-\mathbb{E}\biggl[\sum_{j=N_1}^{N_2} G_i(\alpha^{-1} q_j(0)) \biggr]
\biggr)
\Bigg|
\le B T^{-\mathfrak{c}}.
    \end{equation}

    \noindent We omit a detailed proof of \eqref{eqn:cplL0wtsm0}, which can be deduced quickly from Lemma \ref{lem:HBM_lem} with $t = k = 0$, and $F = 1$. So, we will focus on verifying \eqref{eqn:cplL0wts}.
    
Abbreviate $\L = \L(0)$. To show \eqref{eqn:cplL0wts}, we first use \eqref{eqn:xiH_decomp22} from Theorem \ref{thm:Hindsum_thm}, applied at $t=k=0$, with the $(F(\lambda), G(q))$ there equal to $(\lambda^m, G(\alpha^{-1} q))$ here. We use the notation of the theorem, recalling in particular $R>0$ from Assumption \ref{ass:R_M_ass}, and the submatrices $\L^{[r]}$ and their eigenvalues $(\lambda_i^{(r)})$ from  Definition \ref{def:Lr}. This yields an overwhelmingly probable event $\tilde{\mathsf{E}}_1$ on which
\begin{multline}\label{eqn:GLm_indsum}
   \Bigg|  \sum_{j =N_1 }^{N_2} G(j)  \Lambda_{j}^{m}  - \mathbb{E}\bigg[ \sum_{j =N_1 }^{N_2} G(j)  \Lambda_{j}^{m} \bigg] 
     - \sum_{r = \lceil N_1/R \rceil}^{\lfloor N_2/R \rfloor - 1} G( r R)   \sum_{s=1}^{R} (\lambda_s^{(r)})^m  \\
  + \mathbb{E}\bigg[\sum_{r = \lceil N_1/R \rceil}^{\lfloor N_2/R \rfloor - 1} G( r R) \sum_{s=1}^{R} (\lambda_s^{(r)})^m   \bigg] \Bigg| \le  A B T^{1/2-\mathfrak{c}},
\end{multline}
Since each inner summation in \eqref{eqn:GLm_indsum} is a trace of $(\L^{[r]})^m$, it follows that
\begin{multline}\label{eqn:GLm_indsum2}
   \Bigg|  \sum_{j =N_1 }^{N_2} G(j)  \Lambda_{j}^{m}  - \mathbb{E}\left[ \sum_{j =N_1 }^{N_2} G(j)  \Lambda_{j}^{m} \right] 
     - \sum_{r = \lceil N_1/R \rceil}^{\lfloor N_2/R \rfloor - 1} G( r R)   \sum_{s=0}^{R-1} (\L^{[r]})_{r R+s, r R+s}^m  \\
  + \mathbb{E}\bigg[\sum_{r = \lceil N_1/R \rceil}^{\lfloor N_2/R \rfloor - 1} G( r R) \sum_{s=0}^{R-1}  ((\L^{[r]})^m)_{r R+s, r R+s} \bigg] \Bigg| \le  A B T^{1/2-\mathfrak{c}}.
\end{multline}
Moreover, by the definition of $\L^{[r]}$ and since $(\L^{[r]})_{s s}^m$ only depends on $(\L^{[r]})_{i j}^m$ for $|i-s| \leq m$ and $|j-i| \leq m$, we have 
\begin{equation}\label{eqn:bulkbd}
  \mathbbm{1}_{s \in \llbracket 2m+1, R-2m-1 \rrbracket} \cdot \big( ((\L^{[r]})^m)_{rR+s,rR+s} - (\L^m)_{rR+s,rR+s} \big) = 0.
\end{equation}
In addition, by the mean value theorem, the properties of $|G'|$ and $ \supp G'$ from \eqref{estimateg}, and the fact that $| \mathfrak{k}_i^{[m]}(0) | \leq 3^{m} \max_{i \in \llbracket N_1, N_2 \rrbracket} \{|a_i|^{m}, |b_i|^{m}\}$, on $\tilde{\mathsf{E}}_2 = \mathsf{BND}_{\L(0)}(\log N) $ we have
\begin{equation}\label{eqn:Grtoibd}
 \Bigg|   \sum_{r = \lceil N_1/R \rceil}^{\lfloor N_2/R \rfloor - 1} G( r R)   \sum_{s=0}^{R-1} \mathfrak{k}_{r R +s}^{[m]}(0)  - 
  \sum_{i=N_1}^{N_2} G(i) \mathfrak{k}_i^{[m]}(0) \Bigg| \leq B R (3 \log N)^{C+ m}.
\end{equation}

\noindent Putting \eqref{eqn:bulkbd}, \eqref{eqn:Grtoibd}, and \eqref{eqn:GLm_indsum2} together, we obtain on $\tilde{\mathsf{E}}_1 \cap \tilde{\mathsf{E}}_2 $,
\begin{multline}\label{eqn:cplL0wts1}
      \Bigg|    \sum_{j =N_1 }^{N_2} G(j) \cdot \big(\mathfrak{k}_j^{[m]}(0) -\mathbb{E}[\mathfrak{k}_j^{[m]}(0) ]\big)
         - \sum_{j=N_1}^{N_2} G( j) \Lambda_j^m  
  + \mathbb{E}\bigg[\sum_{j=N_1}^{N_2} G( j) \Lambda_j^m  \bigg] \Bigg|  \\
         \leq \Bigg|  \sum_{r = \lceil N_1/R \rceil}^{\lfloor N_2/R \rfloor - 1} \mathcal{E}_r - \mathbb{E} \bigg[\sum_{r = \lceil N_1/R \rceil}^{\lfloor N_2/R \rfloor - 1} \mathcal{E}_r \bigg]\Bigg|  + B T^{1/2-c},
    \end{multline}
    where, for any $r \in \llbracket N_1/R, N_2/R-1\rrbracket$, we have denoted
    \begin{equation}
      \mathcal{E}_r \coloneqq  \sum_{|s| \leq 2 m } C_s \big( ((\L^{[r]})^m)_{rR+s,r R+s} -(\L^m)_{r R+s, r R+s}\big) ,
    \end{equation}
    where $C_s = G(r R) \cdot \mathbbm{1}_{s \ge 0} + G((r-1) R) \cdot \mathbbm{1}_{s<0}$. To estimate the right side of \eqref{eqn:cplL0wts1}, we proceed as in the end of the proof of Lemma \ref{lem:Psi2_final_approx} (though the setup here will slightly simplify the analysis).
    
    Specifically, observe that the $(\mathcal{E}_r)$ over $r \in \llbracket N_1/R \rceil, N_2/R -1 \rrbracket$ are mutually independent random variables; moreover, on $\tilde{\mathsf{E}}_2$, they satisfy $|\mathcal{E}_r| \leq 5 m (3 \log N)^m$. In addition, since $\supp G \subseteq [-U, U]$, if $|r| R > T (\log N)^6 $, then $\mathcal{E}_r= 0$. Hence, a Chernoff bound, yields with overwhelming probability,
    \begin{equation}\label{eqn:ersumbd}
        \Bigg|  \sum_{r = \lceil N_1/R \rceil}^{\lfloor N_2/R \rfloor - 1} \mathcal{E}_r - \mathbb{E}\bigg[\sum_{r = \lceil N_1/R \rceil}^{\lfloor N_2/R \rfloor - 1} \mathcal{E}_r \bigg]\Bigg|  \leq B T^{1/2} R^{-1/2} (3 \log N)^{m+C}.
    \end{equation}
    Combining \eqref{eqn:ersumbd} with \eqref{eqn:cplL0wts1}, we obtain the desired bound \eqref{eqn:cplL0wts}. This completes the proof.
\end{proof}

\subsection{Charge functionals for the Toda lattice on $\llbracket N_1, N_2 \rrbracket$}
\label{subsec:euler_and_twopoint}

The main result of this subsection is the next lemma, shows the counterpart of Theorem \ref{thm:euler_flucts} for the Toda lattice on $\llbracket N_1, N_2 \rrbracket$. 

\begin{lem}[Euler fluctuations on $\llbracket N_1, N_2 \rrbracket$] \label{thm:time_t_flucts}

There exists a constant $\mathfrak{c}>0$ such that the following holds. Let $K \in \llbracket 1, \log N \rrbracket$ be an integer and $A \ge 1$ be a real number. For each $i \in \llbracket 1, K \rrbracket$, let $\tau_i \in [0, \log N]$, $B_i > 1$, and $U_i \in [\mathfrak{M}, T \log N]$ be real numbers; let $F_i: \mathbb{R} \rightarrow \mathbb{R}$ be a function satisfying Assumption~\ref{ass:F} and \eqref{eqn:Fhold2}; and let $G_i : \mathbb{R} \rightarrow \mathbb{R}$ be a compactly supported function satisfying Assumption ~\ref{ass:G2}, with the $(B,U,\mathfrak{Q})$ there equal to $(B_i, U_i, 0)$ here. For each $i \in \llbracket 1, K \rrbracket$, set $t_i = T \tau_i$ and define $\tilde G_i:\mathbb{R}^2\to\mathbb{R}$ by setting
\[
\tilde G_i(r,\lambda):=G_i (\alpha T r+T \tau_i\ve(\lambda)),
\]  

\noindent for all $r, \lambda \in \mathbb{R}$. 
Then there exists a coupling between $\L(0)$ and $\mathcal{W}$ such that, with overwhelming probability, we have for each $j \in \llbracket 1, K \rrbracket$ that 
    \begin{equation}\label{eqn:time_t_conv}
      \left| T^{-1/2} \left( \sum_{i=N_1}^{N_2} F_j(\Lambda_i(t_j)) \cdot G_j(q_i(t_j)) - \mathbb{E}\bigg[ \sum_{i=N_1}^{N_2} F_j(\Lambda_i(t_j)) \cdot G_j(q_i(t_j)) \bigg]\right) - \mathcal{W}( \tilde{G}_j \cdot F^{\dr}  ) \right| \leq  A B_j  T^{-\mathfrak{c}}  .
    \end{equation}

\end{lem}

\begin{proof}
    
Throughout this proof we recall the Gaussian field $\mathcal{W}^{\dr}$ from Definition \ref{def:Wdress} and the processes $\mathfrak{X}$ and $\mathcal{Z}$ from Definition \ref{def:Z_intro}, which satisfy 
\begin{equation}\label{eqn:Z_X_relation}
   \mathcal{Z}(\lambda, \mathfrak{q}, \tau) = \mathbf{O}_{\Lambda} \mathfrak{X}\big(\Lambda, (\mathfrak{q}-\tau \ve(\Lambda))/\alpha, \tau \big),
\end{equation}
where we have denoted the operator $\mathbf{O}_{\Lambda} = \bm{(\varsigma_0^{\dr})^{-1}} (1 - \theta \mathbf{T} \boldsymbol{\varrho}_{\beta})^{-1}$ acting in first argument $\Lambda$ of $\mathfrak{X}$. 

We verify \eqref{eqn:time_t_conv} under the coupling provided by Theorem \ref{thm:couplethmFG}. By scaling each $G_j$ if necessary, we may assume that $B_j=1$ below. Let $\mathsf{E}_1$ denote the event provided by Theorem \ref{thm:couplethmFG} (with $(A,B, K, \{F_i, G_i, U_i\}_{i=1}^K)$ there as here); $\mathsf{E}_2$ denote the event on which Lemma \ref{lem:Z_sum_conc_general} holds; $\mathsf{E}_3= \bigcap_{t \geq 0}\mathsf{BND}_{\L(t)}(\log N)$ (which is overwhelmingly probable by Lemma \ref{lem:bd_lem}); $\mathsf{E}_4$ denote the event on which Lemma \ref{lem:q_spacing} holds, with $T$ there equal to both $0$ and $T \log N$ here; $\mathsf{E}_5$ denote the event on which Lemma \ref{thm:Z_apriori} holds; and $\mathsf{E}_6$ denote the event on which Lemma \ref{lem:bd_num_termsU} holds for each $U_j$ with $j \in \llbracket 1, K \rrbracket$. We restrict to $\bigcap_{i=1}^6 \mathcal{E}_i$ in what follows.

Now, we fix some $j \in \llbracket 1, K \rrbracket$ and will confirm \eqref{eqn:time_t_conv}. Abbreviate $F = F_j$, $G = G_j$, $t = t_j$, $A= A_j$, and $U = U_j$. By the definition \eqref{eqn:quasi_part_fluct_intro} of $Z_{i}^{\mathcal Q}(t)$, we may rewrite the sum on the left hand side of \eqref{eqn:time_t_conv} as (by changing the index $i$ there to $\varphi_t (\varphi_0^{-1}(i))$ here)
\begin{equation}\label{eqn:statistic_tostudy}
     \sum_{i=N_1}^{N_2} G \big(q_i(0) + t \ve(\Lambda_i) + T^{1/2} Z_{i}^{\mathcal Q}(t) \big) \cdot F(\Lambda_i).
\end{equation}

\noindent By Taylor expanding, for each $i$ there exists $\xi_i$ between $q_i(0)+t\,\ve(\Lambda_i)$ and $q_i(0)+t\,\ve(\Lambda_i)+T^{1/2}Z_i^{\mathcal Q}(t)$ such that
\begin{multline}
    G\bigl(q_i(0)+t\,\ve(\Lambda_i)+T^{1/2}Z_i^{\mathcal Q}(t)\bigr)
=
G\bigl(q_i(0)+t\,\ve(\Lambda_i)\bigr)
+T^{1/2}Z_i^{\mathcal Q}(t) \cdot G'\bigl(q_i(0)+t\,\ve(\Lambda_i)\bigr) \\
+\frac12 \cdot T(Z_i^{\mathcal Q}(t))^2 \cdot G''(\xi_i).
\end{multline}
By our restriction to $\mathsf{E}_3 \cap \mathsf{E}_4 \cap \mathsf{E}_6$, we have $G''(\xi_i) \ne 0$ for at most $U (\log N)^6$ indices $i \in \llbracket N_1, N_2 \rrbracket$, all of which are in $\llbracket -T^{3/2}, T^{3/2} \rrbracket$. Therefore, abbreviating $Z_i^{\mathcal{Q}} = Z_i^{\mathcal{Q}}(t)$, we obtain (using the bounds $|Z_i^{\mathcal Q}| \le (\log N)^C$ from Lemma \ref{eqn:Z_apriori} and $|G''| \le U^{-2} \log N$ from Assumption \ref{ass:G2}) that
\begin{multline}\label{eqn:sts2}
  \Bigg|  \eqref{eqn:statistic_tostudy} - \left(\sum_{i=N_1}^{N_2}  G(q_i(0) + t \ve(\Lambda_i) ) \cdot F(\Lambda_i) 
    + T^{1/2} \sum_{i=N_1}^{N_2}  G'(q_i(0) + t \ve(\Lambda_i) ) \cdot Z_{i}^{\mathcal Q} \cdot F(\Lambda_i) \right)  \Bigg| \\
\leq A U (\log N)^6 \cdot T (\log N)^{70} \cdot U^{-2} \log N \leq  A T^{1/2-c} .
\end{multline}

Next, we estimate the second summation in the first line of \eqref{eqn:sts2}. Our restriction to $\mathsf{E}_2$ (see the first display in Lemma \ref{lem:Z_sum_conc_general}) implies that 
\begin{multline}\label{eqn:sts_bm_term}
\Bigg|   \sum_{i=N_1}^{N_2}  G'(q_i(0) + t \ve(\Lambda_i) ) \cdot Z_{i}^{\mathcal Q} \cdot F(\Lambda_i) -\alpha^{-1} \int_{-\infty}^{\infty} \int_{-\infty}^{\infty} G'(q) \mathcal{Z}(\lambda, qT^{-1} ,\tau) F(\lambda) \varrho(\lambda)d \lambda dq \Bigg| \leq A T^{-c} . 
\end{multline}

\noindent Using \eqref{eqn:sts2}, \eqref{eqn:sts_bm_term}, a tail bound (following the derivation of \eqref{expectationH}), and the fact that $\mathbb{E}[\mathcal{Z}(\lambda,qT^{-1},\tau)] = 0$ we also have (see the derivation of \eqref{eqn:main_est_finalexp} for entirely analogous reasoning) that 
 \begin{equation}\label{eqn:sts_E}
 T^{-1/2} \cdot \left| \mathbb{E}\bigg[\sum_{i=N_1}^{N_2} G(q_i(t)) F(\Lambda_i(t)) \bigg] 
    -   \mathbb{E}\bigg[\sum_{i=N_1}^{N_2} G(q_i(0) + t \ve(\Lambda_i) ) F(\Lambda_i)\bigg] \right|
   \leq  AT^{-c}.
\end{equation}

\noindent Combining \eqref{eqn:sts_E}, \eqref{eqn:sts_bm_term}, and \eqref{eqn:sts2}, we obtain, on the events defined above,
 \begin{multline}\label{eqn:sts_1_term}
  \Bigg| T^{-1/2} \cdot \bigg(  \sum_{i=N_1}^{N_2} G(q_i(t)) F(\Lambda_i(t)) - \mathbb{E}[\cdots]\bigg)  T^{-1/2} \cdot \bigg(\sum_{i=N_1}^{N_2} G(q_i(0) + t \ve(\Lambda_i) ) F(\Lambda_i) - \mathbb{E}[\cdots] \bigg) \\
   -\alpha^{-1} \int_{-\infty}^{\infty} \int_{-\infty}^{\infty} G'(q) \mathcal{Z}(\lambda, qT^{-1} ,\tau) F(\lambda) \varrho(\lambda)d \lambda dq \Bigg|
   \leq T^{1/2-c}.
\end{multline}

\noindent Finally by Theorem \ref{thm:couplethmFG}, we have, denoting $\tilde{G}(r, \lambda) \coloneqq G(\alpha T r + T \tau \ve(\lambda))$, that
\begin{equation*}
   \left| T^{-1/2} \bigg(\sum_{i=N_1}^{N_2} G(q_i(0) + t \ve(\Lambda_i) ) F(\Lambda_i) - \mathbb{E}[\cdots] \bigg)
   - \mathcal{W}^{\dr}( \tilde{G} \cdot F ) \right| \leq T^{-c}.
\end{equation*}

\noindent By this, together with  \eqref{eqn:sts_1_term}, to establish \eqref{eqn:time_t_conv} it suffices to show that 
\begin{flalign} 
\label{eqn:time_t_finalexp}
\mathcal{W}^{\dr}( \tilde{G} \cdot F ) +\alpha^{-1} \int_{-\infty}^{\infty} \int_{-\infty}^{\infty} G'(q) \mathcal{Z}(\lambda, qT^{-1} ,\tau) F(\lambda) \varrho (\lambda) d \lambda dq= \mathcal{W}( \tilde{G} \cdot F^{\dr}  ).
\end{flalign} 
This follows from \eqref{eqn:ttfe} in Lemma \ref{lem:ttfe}, with $G(\mathfrak{q})$ there given by $G(T \mathfrak{q})$ here.
\end{proof}

Next, recall the operator $\mathbf{F} : \mathcal{H} \rightarrow \mathcal{H}$ (given by  $\mathbf{F} (\phi) = \mathbf{D} \phi - \langle \phi, \varsigma_0 \rangle_{\varrho} \cdot \varsigma_0^{\dr}$) from Definition \ref{operatorf}. The following corollary of (the proof of) Lemma \ref{thm:time_t_flucts} the argument of $G_j$ there.
\begin{cor}\label{cor:lattice_fluct_cor}
There exists a constant $\mathfrak{c}>0$ such that the following holds.
     Adopt the notation and assumptions of Lemma \ref{thm:time_t_flucts}, and assume that $\L(0)$ and $\mathcal{W}$ are coupled as in that theorem. For each $j \in \llbracket 1, K \rrbracket$, denote by $\bar{G}_{j} : \mathbb{R}^2 \rightarrow \mathbb{R}$ the function defined by $\bar{G}_j(r, \lambda) = G(T r + \alpha^{-1} T \tau_j \ve(\lambda))$. Then, with overwhelming probability, we have for each $j \in \llbracket 1, K \rrbracket$ that 
  \begin{equation}\label{eqn:time_t_conv_noQ}
      \left| T^{-1/2}\bigg( \sum_{i=N_1}^{N_2} F_j(\Lambda_i(t_j)) \cdot G_j(i) - \mathbb{E}\bigg[ \sum_{i=N_1}^{N_2} F_j(\Lambda_i(t_j)) \cdot G_j(i) \bigg]\bigg) - \mathcal{W}( \bar{G}_j \cdot \mathbf{F} (F_j) ) \right| \leq  A B_j T^{-\mathfrak{c}}  ,
    \end{equation}
    and (recalling $r_i(t)$ from \eqref{abr}) that
    \begin{equation}\label{eqn:time_t_conv_0}
      \left| T^{-1/2} \bigg( \sum_{i=N_1}^{N_2} G_j(i) \cdot  r_i(t_j) - \mathbb{E}\bigg[ \sum_{i=N_1}^{N_2}  G_j(i) \cdot r_i(t_j) \bigg]\bigg) +\alpha \mathcal{W}\left( \bar{G}_j \cdot \varsigma_0^{\dr}  \right) \right| \leq  B_j T^{-\mathfrak{c}}  .
    \end{equation}
    
\end{cor}
\begin{proof}

    Fix $j \in \llbracket 1, K \rrbracket$, and abbreviate $t = t_j$; set $J \coloneqq J_t^{[0]}(0,0)$. By scaling, we may assume that $B_j=1$ below. Observe that, on an overwhelmingly probable event $\mathsf{E}_1$, we have 
  \begin{flalign}\label{eqn:J00bd}
    \begin{aligned}
        |J| &\leq T^{1/2} (\log N)^{20}; \qquad \left| \mathbb{E}[J] \right| &\leq 1.
    \end{aligned}
    \end{flalign}
  
    \noindent Indeed, the first inequality in \eqref{eqn:J00bd} is due to Theorem \ref{thm:current_fluct}, together with the second inequality in \eqref{eqn:J00bd}. To verify the latter, observe under the symmetry $q_j\mapsto -q_{-j}$ of the Toda lattice that the integrated current $J$ is negated, up to the contribution of particle $q_0$; hence $|\mathbb{E}[J] - \mathbb{E}[-J]| \le 1$, which yields the second estimate in \eqref{eqn:J00bd}. Let $\mathsf{E}_2$ denote the event on which Lemma \ref{lem:q_spacing} holds, with $T$ there given by $0$ and by $T \log N $ here, and let $\mathsf{E}_3 \coloneqq  \bigcap_{t \geq 0} \mathsf{BND}_{\L(t)}(\log N)$. Note that on $\mathsf{E}_1 \cap \mathsf{E}_2 \cap \mathsf{E}_3 $, the following holds. Denote by $m_0 \in \llbracket N_1, -J \rrbracket$ the smallest index such that $q_{m_0}(t) \geq 0$. By the restriction to $\mathsf{E}_2 \cap \mathsf{E}_3$, we have $|q_{m_0}(t)| \leq (\log N)^6$, since otherwise it would contradict Lemma \ref{lem:q_spacing}. On the other hand, $i < m_0$ implies $q_{i}(t)< 0$, while (on $\mathsf{E}_2$) the bound $i > m_0 + (\log N)^6 $ implies $q_{i}(t) > q_{m_0}(t) \geq 0$ (by Lemma \ref{lem:q_spacing}), and thus $|-m_0 -J|\leq (\log N)^7$. Therefore, again using the restriction to $\mathsf{E}_2$ and \eqref{eqn:J00bd}, we have
    \begin{equation}\label{eqn:q0_J_bd}
        |q_0(t) - \alpha J| \leq |q_{0}(t) - q_{m_0}(t) - \alpha J| + (\log N)^8 \leq (|m_0|^{1/2}+1) (\log N)^9 \leq T^{1/4} (\log N)^{C}.
    \end{equation}
   
    \noindent As a consequence, on $\mathsf{E}_1 \cap \mathsf{E}_2 \cap \mathsf{E}_3$, we have 
    \begin{equation}\label{eqn:q0_bd}
        |q_0(t) | \leq  T^{1/2} (\log N)^{C}.
    \end{equation}
   
   Let $S = \floor{T^{7/8}}$. Let $\mathsf{E}_4$ be the event given in Lemma \ref{thm:time_t_flucts}, with $(F_j(\lambda), G_j(q), \tau_j)_{j \in \llbracket 1, K \rrbracket}$ there both given by $ (F_j(\lambda), G_j(\alpha^{-1} q), \tau_j)_{j \in \llbracket 1, K \rrbracket}$ and by $ [(1, G_j(\alpha^{-1} q), \tau_j)]_{j \in \llbracket 1, K \rrbracket}$ here. Also let $\mathsf{E}_5$ denote the event on which, for each $r \in \llbracket -T S^{-1} (\log N)^5, T S^{-1} (\log N)^5 \rrbracket$, setting $G = G_j$, $F = F_j$, $A = A_j$, we have 
   \begin{equation}\label{eqn:FGSbd}
     \left|  \sum_{i = r S}^{(r+1) S-1}F(\Lambda_i(t)) \cdot G'(i) - \int_{r S}^{(r+1) S-1} \int_{-\infty}^{\infty}F(\lambda) G'(q) \varrho(\lambda) d \lambda d q
     \right| \leq A T^{1/2} U^{-1} (\log N)^C.
   \end{equation}

    We claim that $\mathsf{E}_5$ is overwhelmingly probable. To see this, first let $(\tilde{\mathbf{a}}(s); \tilde{\mathbf{b}}(s))$ denote the Toda lattice on $\mathbb{Z}$, coupled to the Toda lattice $(\mathbf{a}(s); \mathbf{b}(s))$ on $\llbracket N_1, N_2 \rrbracket$ so that they share the same initial data, namely, so that $(\tilde{a}_i(0), \tilde{b}_i (0)) = (a_i (0), b_i(0))$ for each $i \in \llbracket N_1, N_2 \rrbracket$ (as in Lemma \ref{lem:fin_inf_t_coupling}), where $\tilde{\mathbf{a}}(s) = (\tilde{a}_i(s))_{i \in \mathbb{Z}}$ and $\tilde{\mathbf{b}} (s) = (\tilde{b}_i(s))_{i \in \mathbb{Z}}$. Then, for any $s \ge 0$, let $\tilde{\L}(s) = [\tilde{L}_{ij} (s)]_{i,j\in \llbracket N_1, N_2 \rrbracket}$ denote the (symmetric, tridiagonal) matrix obtained by setting $\tilde{L}_{i,i+1} (s) = \tilde{L}_{i+1,i} (s) = \tilde{a}_i(s)$ for $i \in \llbracket N_1, N_2-1 \rrbracket$, and $\tilde{L}_{i i} (s) = \tilde{b}_i(s)$ for $i \in \llbracket N_1, N_2 \rrbracket$ (and $\tilde{L}_{ij} (s)= 0$ for all other $(i,j)$). Letting $\eig \tilde{\mathbf{L}} (s) = (\tilde{\lambda}_1 (s), \tilde{\lambda}_2 (s), \ldots , \tilde{\lambda}_N (s))$ and $\tilde{\varphi}_s : \llbracket N_1, N_2 \rrbracket \rightarrow \llbracket 1, N \rrbracket$ denote a $\zeta$-localization center bijection for $\tilde{\mathbf{L}}(s)$, set  $\tilde{\Lambda}_i (s) = \lambda_{\tilde{\varphi}_s^{-1}(i)} (s)$ for each $i \in \llbracket 1, N \rrbracket$ (as in Definition \ref{def:loccent}). Then, it follows quickly from Lemma \ref{lem:Lax_eig_coupling} (see also the derivation of \eqref{eqn:cplML1}) that
    \begin{flalign}
        \label{sumsr}
    \Bigg|  \sum_{i = r S}^{(r+1) S-1}(F(\Lambda_i(t))-F(\tilde{\Lambda}_i (t)) ) \cdot G'(i) \Bigg| \leq A U^{-1} (\log N)^C. 
    \end{flalign} 
    
    \noindent Moreover, using the first part of Lemma \ref{lem:no_q_conc} and Lemma \ref{lem:no_q_conc_expect} (at $t=q=0$, with $(A,B,S,U)$ there given by $(A,  U^{-1} , U, U)$ here up to $\alpha$-dependent constants, the $G(x)$ there given by $G'(-\alpha^{-1} x)$ here, we deduce that, with overwhelming probability, 
    \begin{flalign*} 
    \Bigg| \sum_{i = r S}^{(r+1) S-1} F(\tilde{\Lambda}_i (0)) \cdot  G'(i) - \int_{r S}^{(r+1) S-1} \int_{-\infty}^{\infty}F(\lambda) G'(q) \varrho(\lambda) d \lambda d q
     \Bigg| \leq A  T^{1/2} U^{-1} (\log N)^C.
    \end{flalign*} 
    
    \noindent Combining this with \eqref{sumsr}, and the fact that $(\tilde{\Lambda}_i(t))$ has the same law as $(\tilde{\Lambda}_i(0))$ due to the invariance of $(\tilde{\mathbf{a}} (s); \tilde{\mathbf{b}} (s))$ (by the first part of Lemma \ref{lem:fin_inf_t_coupling}), yields \eqref{eqn:FGSbd}. 
    
    Next, let $\mathsf{E}_6$ be the event on which Lemma \ref{lem:euler_to_linstat} holds. We restrict to $\bigcap_{i=1}^6 \mathsf{E}_i$ in what follows. We claim that for any $j \in \llbracket 1, K \rrbracket$, letting $F = F_j$, $G = G_j$, and $t = t_j$, 
  \begin{flalign}\label{eqn:FGw1}
    \begin{aligned}
        &\sum_{i=N_1}^{N_2}  F(\Lambda_i(t)) \cdot (G(i)-G(\alpha^{-1} q_i(t))) \\
        &= -\sum_{i=-\floor{T (\log N)^5}}^{\floor{T (\log N)^5}} F(\Lambda_i(t)) \cdot G'(i) \cdot (\alpha^{-1} q_i(t) - i) + O(A T^{1/2-c}) \\
        &= -\sum_{r=-\floor{T (\log N)^5/S}}^{\floor{T (\log N)^5/S}-1}  (\alpha^{-1} q_{r S}(t)-  r S)  \sum_{i = r S}^{(r+1) S-1}F(\Lambda_i(t)) \cdot G'(i)  + O(A T^{1/2-c})\\
        &= -\sum_{r=-\floor{T (\log N)^5/S}}^{\floor{T (\log N)^5/S}-1}  (\alpha^{-1} q_{r S}(t)-  r S) \int_{r S}^{(r+1) S-1} \int_{-\infty}^{\infty}F(\lambda) G'(q) \varrho(\lambda) d \lambda d q \\
        & \qquad + O(T S^{-1} (\log N)^C \cdot T^{1/2} (\log N)^C \cdot A T^{1/2} U^{-1} (\log N)^C)  + O(A T^{1/2-c}) \\
        &= -\sum_{r=-\floor{T (\log N)^5/S}}^{\floor{T (\log N)^5/S}-1}  (\alpha^{-1} q_{r S}(t)-  r S) \int_{r S}^{(r+1) S-1} \int_{-\infty}^{\infty}F(\lambda) G'(q) \varrho(\lambda) d \lambda +  O(A T^{1/2-c}),
        \end{aligned}
    \end{flalign}

    \noindent where for any $\Gamma \ge 0$ we have let $O(\Gamma)$ be some real number that can vary between lines satisfying $|O(\Gamma)| \le \Gamma$. Indeed, the first equality follows from first using $\supp G \subseteq [-U, U]$ and the combination of \eqref{eqn:q0_bd} and Lemma \ref{lem:q_spacing} to first truncate the sum, and then Taylor expanding and bounding the second derivative contribution by $A \cdot 3 T (\log N)^5 \cdot T (\log N)^C \cdot U^{-2} \log N \leq AT^{1/4}$ (which again uses \eqref{eqn:q0_bd} and Lemma \ref{lem:q_spacing}, as well as the second derivative bound \eqref{estimateg} satisfied by $G$). The second is from the bound, which by Lemma \ref{lem:q_spacing} holds for each $i \in \llbracket r S, (r+1) S-1 \rrbracket$, 
$$
\left| (\alpha^{-1} q_{r S}(t)-  r S) - (\alpha^{-1} q_{i}(t)-  i)  \right| \leq S^{1/2} (\log N)^C.
$$
The third is by \eqref{eqn:FGSbd} combined with the bound $\left|\alpha^{-1} q_{r S}(t)-  r S \right| \leq T^{1/2} (\log N)^C$ (which in turn follows from Lemma \ref{lem:q_spacing} and \eqref{eqn:q0_bd}). The fourth is simply the observation that the error in the third equality can be absorbed into the  $O(A T^{1/2-c})$ one.

By similar reasoning, the display \eqref{eqn:FGw1} remains valid with overwhelming probability, upon replacing $F$ by $\langle F, \varsigma_0 \rangle_{\varrho} \cdot \varsigma_0$. Upon such a replacement, the final line of \eqref{eqn:FGw1} remains the same, so it follows that with overwhelming probability we have 
\begin{equation}\label{eqn:FGw2}
  \left|   \sum_{i=N_1}^{N_2}  F(\Lambda_i(t)) \cdot \big(G(i)-G(\alpha^{-1} q_i(t)) \big) - \langle F, \varsigma_0 \rangle_{\varrho} \sum_{i=N_1}^{N_2} \big(G(i)-G(\alpha^{-1} q_i(t)) \big) \right| \leq A T^{1/2-c}.
\end{equation}
We further restrict to the event on which \eqref{eqn:FGw2} holds. Now, $\sum_{i=N_1}^{N_2} (G(i)-G(\alpha^{-1} q_i(t))) - \mathbb{E}[\cdots]$ approximates an explicit Gaussian limit, by Lemma \ref{thm:time_t_flucts}. Specifically, letting $\bar{G}_j(r, \lambda) = G(T r + \alpha^{-1} T \tau_j \ve(\lambda)))$, that theorem couples $\L(0)$ and $\mathcal{W}$ so that 
\begin{equation}\label{eqn:GF1fluct}
\left| T^{-1/2} \bigg( \sum_{i=N_1}^{N_2} (G(i)-G(\alpha^{-1} q_i(t))) - \mathbb{E}[\cdots] \bigg)
+ \mathcal{W}\left( \bar{G}_j \cdot \varsigma_0^{\dr} \right) \right| \leq T^{-c}.
\end{equation}
Combining this with \eqref{eqn:FGw2} and the fact by Lemma \ref{thm:time_t_flucts} that 
$$\left| T^{-1/2} \bigg( \sum_{i=N_1}^{N_2}F(\Lambda_i(t)) \cdot G(\alpha^{-1} q_i(t)) - \mathbb{E}[\cdots] \bigg)
- \mathcal{W}\left( \bar{G}\cdot F^{\dr} \right) \right| \leq T^{-c},
$$ 

 \noindent completes the proof of \eqref{eqn:time_t_conv_noQ}. Finally, \eqref{eqn:time_t_conv_0} follows from combining \eqref{eqn:etl0} and \eqref{eqn:GF1fluct}.
\end{proof}

\subsection{Proof of Theorem \ref{thm:euler_flucts} and Corollary \ref{cor:q0bm}}
\label{subsec:eulerfluctproof}

\begin{proof}[Proof of Theorem \ref{thm:euler_flucts}]
Assume throughout that $T\ge T_0$, for some sufficiently large number $T_0=T_0(k)$ (dependent on $k$, $f_1,\dots,f_k$, and $(t_1,m_1),\dots,(t_k,m_k)$). We couple the Toda lattice on $\llbracket N_1,N_2\rrbracket$ (with $(N_1,N_2,N,T)$ satisfying Assumption~\ref{ass:NT_assumption}) to the Toda lattice on $\mathbb{Z}$, both started from thermal equilibrium \eqref{eqn:equil}, as in Lemma~\ref{lem:fin_inf_t_coupling} (see also the proof of Theorem~\ref{thm:current_fluct_intro}). We also use the coupling between $(a_i(0), b_i(0))_{i \in \llbracket N_1, N_2 \rrbracket}$, the initial data of the Toda lattice, and $\mathcal{W}^{\dr}$ constructed in Corollary~\ref{cor:lattice_fluct_cor} with $G_i(x)=f_i(xT^{-1})$. Throughout, we adopt the notation in the first two paragraphs of the proof of Theorem \ref{thm:current_fluct_intro}. In particular, $\tilde{\mathfrak{K}}_j^{[m]}(t)$ is a local charge for the Toda lattice on $\llbracket N_1, N_2 \rrbracket$, while and $\mathfrak{K}_j^{[m]}(t)$ is for the Toda lattice on $\mathbb{Z}$.

In this proof, we will show a slightly stronger statement. We will show that if $G_i$ satisfies Assumption \ref{ass:G2} with the $(B,U)$ there equal to $(B_i, U_i)$ here, real numbers $B_i \geq 1$ and $U_i \in [\mathfrak{M}, T \log N]$, then under the coupling between $(a_i(0), b_i(0))_{i \in \mathbb{Z}}$ and $\mathcal{W}$ adopted above, with probability at least $1-c^{-1} e^{-c (\log T)^2}$, for each $i \in \llbracket 1, k \rrbracket$,
\begin{multline}\label{eqn:general_bd}
T^{-1/2}\Bigg|
\sum_{j\in\mathbb{Z}} G_i(j) \cdot \big(\mathfrak{K}_j^{[m_i]}(T\tau_i)-\mathbb{E}\big[\mathfrak{K}_j^{[m_i]}(T\tau_i)\big]\big)
\\
- \mathcal{W}\left(  G_i\big(T(r + \alpha^{-1} \tau_i \ve(\lambda))\big) \cdot \mathfrak{g}_i (\lambda) \right)
\Bigg|
\le B_i T^{-c},
\end{multline}
where $\mathfrak{g}_i$ is as in the theorem statement, and $c>0$ is a constant depending only on $m_1,\dots, m_k$.

The theorem follows from \eqref{eqn:general_bd} by taking $G_i(x)=f_i(x T^{-1})$; taking $U_i \in [\mathfrak{M}, C T ]$ large enough that $\supp f_i\subset[-U_i T^{-1},U_i T^{-1}]$; and setting each $B_i$ equal to 
\[
B:=\max_{i\in\llbracket 1,k\rrbracket} \Big\{\sup_{x\in\mathbb{R}}|f_i(x)|, U_i T^{-1} \sup_{x\in\mathbb{R}}|f_i'(x)|, U_i^2 T^{-2} \sup_{x\in\mathbb{R}}|f_i''(x)|, U_i^3 T^{-3} \sup_{x\in\mathbb{R}}|f_i^{(3)}(x)|\Big\}.
\]

Let $\mathsf{E}_1:=\bigcap_{t\ge 0}\mathsf{BND}_{\L(t)}(\log N)$ and let $\mathsf{E}_2$ be the event from Lemma~\ref{lem:fin_inf_t_coupling}. Let $\mathsf{E}_3$ be the event from Corollary~\ref{cor:lattice_fluct_cor} on which \eqref{eqn:time_t_conv_noQ} and \eqref{eqn:time_t_conv_0} hold, and let $\mathsf{E}_4$ be the event from Lemma~\ref{lem:euler_to_linstat} on which \eqref{eqn:etl} holds; here, we use Corollary~\ref{cor:lattice_fluct_cor} and Lemma~\ref{lem:euler_to_linstat} with $K = k$, with $F_i = \varsigma_i$, and with $G_i$ and parameters $(B_i,U_i)$ as in \eqref{eqn:general_bd}. For $T$ sufficiently large we have $K \leq \log N$, $U_i \le T (\log N)^{1/2}$, $t_i\le T\log N$, $m_i\le(\log N)^{1/10}$, and $B\le \log N$. We work on the intersection $\mathsf{E}:=\mathsf{E}_1\cap\mathsf{E}_2\cap\mathsf{E}_3\cap\mathsf{E}_4$, which has overwhelming probability.

By the combination of \eqref{eqn:time_t_conv_noQ} with \eqref{eqn:etl} for $m_i >0$, or by \eqref{eqn:time_t_conv_0} for $m_i = 0$, it suffices to show that for each $i\in\llbracket 1,k\rrbracket$,
\begin{multline}\label{eqn:bd_we_want}
T^{-1/2}\Bigg|
\sum_{j\in\mathbb{Z}} G_i(j) \cdot \big(\mathfrak{K}_j^{[m_i]}(T\tau_i)-\mathbb{E}\big[\mathfrak{K}_j^{[m_i]}(T\tau_i)\big]\big)
\\
-\sum_{j =N_1}^{N_2} G_i(j ) \cdot \big(\tilde{\mathfrak{K}}_j^{[m_i]}(T\tau_i)-\mathbb{E}\big[\tilde{\mathfrak{K}}_j^{[m_i]}(T\tau_i)\big]\big)
\Bigg|
\le B_i T^{-c}.
\end{multline}

\noindent By our restriction to $\mathsf{E}_1$, \eqref{eqn:fin_inf_dif} implies that $| \mathfrak{k}_j^{[m_i]}(T \tau_i) -  \tilde{\mathfrak{k}}_j^{[m_i]}(T \tau_i) | \leq m_i (3 \log N)^{m_i}  e^{-T^2/5}$, for $j \in \llbracket N_1+T^2, N_2 - T^2 \rrbracket$. Together with the facts that $\supp G_i\subseteq [-U_i ,U_i ]$, that $\sup_{x\in \mathbb{R}}|G_i(x)| \leq B_i$, and that $m_i \leq (\log N)^{1/10}$, this yields \eqref{eqn:bd_we_want} for $m_i > 0$. If $m_i = 0$, the argument is similar; the only difference here is, to bound $|r_j(t) - \tilde r_j(t)|$ for $|i| \le T (\log N)^{10}$ using \eqref{eqn:fin_inf_dif}, we further restrict to the event $\bigcap_{i=1}^{k} \bigcap_{j=N_1}^{N_2} \{ \tilde a_j(t_i) \ge e^{-(\log N)^2} \}$, which is overwhelmingly probable by the explicit density for $(\tilde{a}_j)$ (from \eqref{eqn:equil} and invariance of the Toda lattice on $\mathbb{Z}$ under thermal equilibrium). This establishes the corollary.
\end{proof}

Next, we reuse the argument in the first paragraph of the previous proof to give a proof of Corollary \ref{cor:q0bm}. 

\begin{proof}[Proof of Corollary \ref{cor:q0bm}]
Adopt the same coupling between the Toda lattice $(\bm{p}(t);\bm{q}(t))$ on $\mathbb{Z}$ and that $(\tilde{\bm{p}}(t); \tilde{\bm{q}}(t))$ on $\llbracket N_1, N_2 \rrbracket$ as used in the proof of Theorem \ref{thm:euler_flucts}. By Lemma \ref{lem:fin_inf_t_coupling}, it suffices to prove convergence (in the sense of finite-dimensional distributions) of $T^{-1/2} \tilde{q}_0 (T \tau)$ to $\mathcal{B}(\tau)$. Recall from \eqref{eqn:q0_J_bd} (and Lemma \ref{lem:fin_inf_t_coupling}) that $|\tilde{q}_0 (T\tau) - \alpha J_{T\tau}^{[0]} (0,0)| \le T^{1/4} (\log N)^C$ holds with overwhelming probability, for any $\tau \in [0, \log N]$. As such, recalling $\mathcal{W}$ from Definition \ref{def:white_noise_dress}, the $m_i = \mathfrak{q}_i = \mathfrak{q}_i' = 0$ case of Theorem \ref{thm:current_fluct_intro} implies that $\tilde{q}_0 (T\tau)$ converges to the Gaussian process 
\begin{flalign}
    \label{b2} 
\tilde{\mathcal{B}} (\tau) \coloneqq \alpha \cdot \mathcal{W} (\varsigma_0^{\dr} (\lambda) \cdot  \big(
\mathbbm{1}\{\alpha r<0\}
-\mathbbm{1}\{\alpha r+\tau \ve(\lambda)<0\}).
\end{flalign} 

It remains to identify $\tilde{\mathcal{B}}(\tau)$ with $\mathcal{B}(\tau)$, to which end we must compute the covariance structure of the process $\tilde{\mathcal{B}}(\tau)$. By \eqref{b2} and Definition \ref{def:white_noise_dress} for $\mathcal{W}$, we obtain for any $0 \le \tau_1 \le \tau_2$ that 
\begin{align}
    \label{b3} 
\begin{aligned} 
\Cov(\tilde{\mathcal{B}}(\tau_1), \tilde{\mathcal{B}}(\tau_2))
&= \alpha^2 \int_{-\infty}^{\infty}\int_{-\infty}^{\infty}|\varsigma_0^{\dr}(\lambda)|^2 \cdot \bigl(\mathbbm{1}\{\alpha r < 0\}-\mathbbm{1}\{\alpha r+\tau_1 \ve(\lambda) < 0\}\bigr) \\
&\qquad\qquad \qquad  \times \bigl(\mathbbm{1}\{\alpha r < 0\}-\mathbbm{1}\{\alpha r+\tau_2 \ve(\lambda) < 0\}\bigr)\varrho(\lambda)\,dr\,d\lambda \\
&= \alpha \tau_1\int_{-\infty}^{\infty}|\varsigma_0^{\dr}(\lambda)|^2 |\ve(\lambda)|\varrho(\lambda)\,d\lambda,
\end{aligned} 
\end{align}

\noindent where in the last equality we used the fact that, for any $\lambda \in \mathbb{R}$, 
\begin{flalign*} 
\int_{-\infty}^{\infty} \bigl(\mathbbm{1}\{\alpha r < 0\}-\mathbbm{1}\{ \alpha r+\tau_1 \ve(\lambda) < 0\} \bigr) \cdot \bigl(\mathbbm{1}\{ & \alpha r < 0\}-\mathbbm{1}\{\alpha r+\tau_2 \ve(\lambda) < 0\}\bigr) dr \\
& = \alpha^{-1} \min \{ \tau_1, \tau_2 \} \cdot |\ve (\lambda)| = \alpha^{-1} \tau_1 |\ve(\lambda)|.  
\end{flalign*}

\noindent Since \eqref{b3} coincides with the covariance structure of the Brownian motion $\mathcal{B}$ of variance \eqref{eqn:q0var}, this yields the corollary. 
\end{proof}

\subsection{Proof of Corollary \ref{cor:twopoint_intro}}
\label{subsec:tpproof}

\begin{proof}[Proof of Corollary \ref{cor:twopoint_intro}]

Let us assume $m, n \geq 1$, as the proof if either $m = 0$ or $n = 0$ is similar. Let $U \coloneqq T^{1-\mathfrak{c}/10}$, where $\mathfrak{c} = \mathfrak{c}(m, n) > 0$ is the constant $c$ appearing in the right hand side of the bound \eqref{eqn:general_bd}. Further define $G:\mathbb{R}\to\mathbb{R}$ by setting $G = \chi_U'$, where $\chi_U$ is as in Definition~\ref{def:chi_def}, but with the $\mathfrak{M}$ there replaced by
$U$ here. Specifically, it is a smooth function such that $\chi_U'$ is even; $0 \le \chi_U (x) \le 1$ for all $x \in \mathbb{R}$; for $ x > U$ we have $\chi_U (x) \equiv 1$; for $x < - U$ we have $\chi_U (x) \equiv 0$; and
\begin{equation}
    \label{estimatesm00u} 
\begin{gathered}
\supp \chi_U' \subseteq [-U,U], \\
\sup_{x \in \mathbb{R}} |\chi_U'(x)| \leq 10U^{-1}; \qquad \sup_{x \in \mathbb{R}} |\chi_U''(x)| \leq 100U^{-2}; \qquad \sup_{x \in \mathbb{R}} |\chi_U'''(x)| \leq 1000U^{-3}.
\end{gathered}
\end{equation}

We use the coupling between the Toda lattice on $\mathbb{Z}$ and $\mathcal{W}^{\dr}$ constructed in the proof of Theorem~\ref{thm:euler_flucts} (via the Toda lattice on $\llbracket N_1,N_2\rrbracket$ with $(N_1,N_2,N,T)$ satisfying Assumption~\ref{ass:NT_assumption}). In the notation of that proof we take $k=2$, let $f$ play the role of $f_1$, $t = t_1$, $m = m_1$, and let $x\mapsto T G(Tx)$ play the role of $f_2$. Denote $f_\tau(r,\lambda):=f(r+\alpha^{-1}\tau\ve(\lambda))$ and $G_T(r):=T G(Tr)$.

Note that Theorem \ref{thm:euler_flucts} does not directly apply because now $f_2$ depends on the large parameter $T$. Nevertheless, \eqref{eqn:general_bd} from the proof still applies in this setting, and in particular by that display, with probability at least $1-\mathfrak{c}^{-1}e^{-\mathfrak{c} (\log T)^2} $, we have
\begin{equation}\label{eqn:fcoupleerr}
     \left| T^{-1/2}  \sum_{j \in \mathbb{Z} } f(j T^{-1}) \cdot \big(  \mathfrak{k}_j^{[m]}(t) - \mathbb{E}[ \mathfrak{k}_j^{[m]}(t) ] \big) -   \mathcal{W}\left(  f_{\tau} \cdot \mathbf{F} \varsigma_{m} \right)\right| \leq T^{-\mathfrak{c}},
\end{equation}
and, by the definition of $U$ with \eqref{estimatesm00u}, 
    \begin{equation}\label{eqn:Gcoupleerr}
  \left| T^{-1/2}  \sum_{j \in \mathbb{Z} }T G(j) \cdot \big(  \mathfrak{k}_j^{[n]}(0) - \mathbb{E}[ \mathfrak{k}_j^{[n]}(0) ] \big) -   \mathcal{W}\left( G_T \cdot \mathbf{F} \varsigma_{n}  \right)\right| \leq  T^{-\mathfrak{c}/2}.
  \end{equation}

By \eqref{eqn:fcoupleerr} and \eqref{eqn:Gcoupleerr}, together with Lemma \ref{lem:cb}, which implies that 
     \begin{equation}\label{eqn:rembd3}
      \left| \Cov\left( \mathcal{W}\left( G_T \cdot \mathbf{F} \varsigma_{n} (\lambda) \right) , \mathcal{W}\left(  f_{\tau} \cdot \mathbf{F} \varsigma_{m} (\lambda) \right)\right)  - \sum_{j_1, j_2 \in \mathbb{Z}}G(j_1) f(j_2 T^{-1}) \Cov\left(  \mathfrak{k}_{j_1}^{[n]}(0), \mathfrak{k}_{j_2}^{[m]}(t) \right) \right| \leq  T^{-c},
    \end{equation}
in order to prove the corollary it suffices to show the following two bounds. First, 
    \begin{equation}\label{eqn:rembd1}
      \left| \sum_{j_1, j_2 \in \mathbb{Z}}G(j_1) f(j_2 T^{-1}) \Cov\left(  \mathfrak{k}_{j_1}^{[n]}(0), \mathfrak{k}_{j_2}^{[m]}(t) \right)  -\sum_{j \in \mathbb{Z}}f(j T^{-1}) \Cov\left(  \mathfrak{k}_{0}^{[n]}(0), \mathfrak{k}_{j}^{[m]}(t) \right)  \right| \leq  T^{-c},
    \end{equation}
    and second (recalling $m, n \geq 1$),
    \begin{equation}\label{eqn:rembd2}
      \left| \Cov\left( \mathcal{W}\left( G_T \cdot \mathbf{F} \varsigma_{n} (\lambda) \right) , \mathcal{W}\left(  f_{\tau} \cdot \mathbf{F} \varsigma_{m} (\lambda) \right)\right)  -\int_{-\infty}^{\infty} \mathbf{F} \varsigma_m(\lambda) \cdot \mathbf{F} \varsigma_n(\lambda) \cdot 
     f(\alpha^{-1} \tau \ve(\lambda)) \varrho(\lambda) d \lambda \right| \leq  T^{-c}.
    \end{equation}

    First we show \eqref{eqn:rembd2}. By the definition of the white noise $\mathcal{W}$ (Definition \ref{def:white_noise_dress}), we obtain 
\begin{multline}\label{eqn:r0}
\Cov\left( \mathcal{W}\left( G_T \cdot \mathbf{F} \varsigma_{n} (\lambda) \right) , \mathcal{W}\left(  f_{\tau} \cdot \mathbf{F} \varsigma_{m}(\lambda) \right)\right) \\
= T \int_{-\infty}^{\infty} \int_{-\infty}^{\infty} G(T r) f(r + \alpha^{-1}\tau \ve(\lambda)) \mathbf{F} \varsigma_{n} (\lambda) \mathbf{F} \varsigma_{m} (\lambda) \varrho(\lambda) dr d \lambda.
\end{multline}

\noindent Moreover, for any $\lambda \in \mathbb{R}$, we have since $\supp G \subseteq [-U, U]$ and $f$ is smooth that 
 \begin{multline}\label{eqn:r1}
  T  \int_{-\infty}^{\infty} G(T r) \cdot \big| f(r + \alpha^{-1}\tau \ve(\lambda)) - f(\alpha^{-1}\tau \ve(\lambda)) \big| d r   \leq  C U \int_{-\infty}^{\infty} G(T r) d r 
      = C U T^{-1} \leq T^{-c}.
    \end{multline}

 \noindent By \eqref{eqn:r0} and \eqref{eqn:r1} (with the boundedness of $\mathbf{F}$), we deduce \eqref{eqn:rembd2}.

 To show \eqref{eqn:rembd1}, observe that upon setting $M = \lceil 20\mathfrak{c}^{-1} \rceil$ we obtain 
 \begin{align}
   &  \Bigg| \sum_{j_1, j_2 \in \mathbb{Z}}G(j_1) f(j_2 T^{-1}) \Cov\left(  \mathfrak{k}_{j_1}^{[n]}(0), \mathfrak{k}_{j_2}^{[m]}(t) \right)  \notag\\
     & \qquad \qquad - \sum_{ j_1 \in \mathbb{Z}}\sum_{j_2 \in \mathbb{Z}} G(j_1) f((j_2-j_1) T^{-1}) \Cov\left(  \mathfrak{k}_{0}^{[n]}(0), \mathfrak{k}_{j_2-j_1}^{[m]}(t) \right)  \Bigg|   \notag \\
 &\leq   \left| \sum_{j_1, j_2 \in \mathbb{Z}} G(j_1) \left( \sum_{p=1}^M f^{(p)}(j_2 T^{-1})  \frac{j_1^p}{T^p p!} +f^{(M+1)}(\xi_{j_1,j_2} T^{-1}) \frac{(j_1 T^{-1})^{M+1}}{(M+1)!}  \right) \Cov\left(  \mathfrak{k}_{j_1}^{[n]}(0), \mathfrak{k}_{j_2}^{[m]}(t) \right)    \right|  \notag\\
 &\leq  \sum_{p=1}^M  \left|  \Cov\left(\sum_{j_1 \in \mathbb{Z}}  G(j_1) \frac{(j_1 T^{-1})^p}{p!}  \mathfrak{k}_{j_1}^{[n]}(0), \sum_{j_2 \in \mathbb{Z}} f^{(p)}(j_2 T^{-1})  \mathfrak{k}_{j_2}^{[m]}(t) \right)    \right| \notag\\
 &\qquad + T \cdot C  \frac{(U T^{-1})^{M+1}}{(M+1)!} \cdot (C \log N)^{m+n}.  \label{eqn:s1}
    \end{align}

     \noindent where $\xi_{j_1,j_2}$ is a real number between $j_2-j_1$ and $j_2$. Here, to obtain the first statement we used the mean value form of the remainder for the Taylor expansion. To obtain the second we used the facts that $|f^{(M+1)}| \leq C$, that $\supp G \subseteq [-U, U]$, that $|G| \le 10U^{-1}$; that $|\Cov(\mathfrak{k}_{0}^{[n]}(0),\mathfrak{k}_{j_2-j_1}^{[m]}(t))| \le (C \log T)^{m+n}$ for $j_1,j_2 \in \mathbb{Z}$ by Lemma \ref{lem:bd_lem}; and that  $\supp f \subseteq [-C,C]$. Since $M \ge 20\mathfrak{c}^{-1}$ and $UT^{-1} = T^{-\mathfrak{c}/10}$, the right side of \eqref{eqn:s1} is at most $T^{-c}$.

Next, we again invoke \eqref{eqn:general_bd} and Lemma \ref{lem:cb} to obtain $M$ different analogs of \eqref{eqn:rembd3}; for each $p \in \llbracket 1, M\rrbracket$, \eqref{eqn:rembd3} holds with $G(x)$ there replaced by $G(x) (x T^{-1})^p/p!$ and $f$ there replaced by $f^{(p)}$. Thus, for each $p \in \llbracket 1, M \rrbracket$, we have by \eqref{estimatesm00u} that
\begin{align}\label{eqn:s1b}
& \left|  \Cov\left(\sum_{j_1 \in \mathbb{Z}}  G(j_1) \frac{(j_1 T^{-1})^p}{p!}  \mathfrak{k}_{j_1}^{[n]}(0), \sum_{j_2 \in \mathbb{Z}} f^{(p)}(j_2 T^{-1})  \mathfrak{k}_{j_2}^{[m]}(t) \right)    \right|  \notag \\
& \leq \left|   \Cov\big( \mathcal{W}( G_T \varsigma_p \mathbf{F} \varsigma_n) , \mathcal{W}( f_{\tau}^{(p)} \mathbf{F} \varsigma_m)  \big)    \right| +T^{-c} \notag \\
&= \int_{-\infty}^{\infty}\int_{-\infty}^{\infty} T G(T r) \cdot r^p \cdot f^{(p)}(r + \alpha^{-1}\tau \ve(\lambda)) \mathbf{F} \varsigma_n(\lambda)  \mathbf{F} \varsigma_m(\lambda)  \varrho(\lambda) d r d \lambda +T^{-c} \notag \\
&\leq 10 TU^{-1} \cdot (U T^{-1})^p \cdot UT^{-1} \displaystyle\max_{|r| \le UT^{-1}} \int_{-\infty}^{\infty} f^{(p)}(r + \alpha^{-1}\tau \ve(\lambda)) \mathbf{F} \varsigma_n(\lambda)  \mathbf{F} \varsigma_m(\lambda)  \varrho(\lambda) d \lambda +T^{-c} \notag \\
&\leq C (U T^{-1})^p +T^{-c} \leq 2C T^{-c}.
\end{align}

    \noindent In addition, note that 
    \begin{flalign}\label{eqn:s2}
    \begin{aligned}
    & \sum_{ j_1 \in \mathbb{Z}}\sum_{j_2 \in \mathbb{Z}} G(j_1) f((j_2-j_1) T^{-1}) \Cov\left(  \mathfrak{k}_{0}^{[n]}(0), \mathfrak{k}_{j_2-j_1}^{[m]}(t) \right) \\
     &\qquad = \sum_{j \in \mathbb{Z}}  f(j T^{-1}) \Cov\left(  \mathfrak{k}_{0}^{[n]}(0), \mathfrak{k}_{j}^{[m]}(t) \right) \sum_{ j_1 \in \mathbb{Z}}G(j_1) \\
     &\qquad =\left(1 + O(U^{-1}) \right) \sum_{j \in \mathbb{Z}}  f(j T^{-1}) \Cov\left(  \mathfrak{k}_{0}^{[n]}(0), \mathfrak{k}_{j}^{[m]}(t) \right) 
      \end{aligned}
      \end{flalign}

      \noindent where for any $\Gamma \ge 0$ we have let $O(\Gamma)$ denote a real number satisfying $|O(\Gamma)| \le C\Gamma$ that might vary between lines. Above, we have used the fact that $\Cov(\mathfrak{k}_{j_1}^{[n]}(0),\mathfrak{k}_{j_2}^{[m]}(t)) =\Cov(\mathfrak{k}_{0}^{[n]}(0),\mathfrak{k}_{j_2-j_1}^{[m]}(t))$ (by the translation-invariance of thermal equilibrium) and the bound   $| \sum_{ j \in \mathbb{Z}} G(j) - 1 | \leq C U^{-1}$ (as $G = \chi_U'$). Thus \eqref{eqn:s1}, \eqref{eqn:s1b}, and \eqref{eqn:s2} combine to give 
      \begin{multline}\label{eqn:s3}
        \sum_{j \in \mathbb{Z}}  f(j T^{-1}) \Cov\left(  \mathfrak{k}_{0}^{[n]}(0), \mathfrak{k}_{j}^{[m]}(t) \right) \\
        =\left( 1 +  O(U^{-1})\right) \sum_{j_1, j_2 \in \mathbb{Z}}G(j_1) f(j_2 T^{-1}) \Cov\left(  \mathfrak{k}_{j_1}^{[n]}(0), \mathfrak{k}_{j_2}^{[m]}(t) \right) + O(T^{-c}).
      \end{multline}

      \noindent By \eqref{eqn:s3} and \eqref{eqn:rembd3} (which with \eqref{eqn:rembd2} shows that the sum in the right hand side of \eqref{eqn:s3} is at most $C$), implies \eqref{eqn:rembd1} and thus the corollary. 
    \end{proof}

\appendix

\section{Proofs of several discrete estimates} 

\label{ProofQT}

\subsection{Proofs of lemmas about $\ve$}
\label{sec:veff_lemma}
\label{app:Tdr_sym}

\begin{proof}[Proof of Lemma \ref{lem:veff_inc}]
We begin with the observation, explained in \cite[Remark A.5]{Agg25}, that the constant $c>0$ appearing in the bound $\mathbf{T} \varrho (x) + \alpha >c$, proved in \cite[Lemma 3.3]{Agg25}, can be taken independent of $\theta$, so long as $\theta \leq 1$.

Throughout this proof we use the present paper's convention for dressing, which is off from that of \cite{Agg25} by a factor of $\theta$:
\begin{equation}\label{eqn:Amol_dress_app}
f^{\dr} = (1-\theta\,\mathbf{T}\boldsymbol{\varrho}_{\beta})^{-1}f .
\end{equation}
Equivalently, $f^{\dr}$ is the unique solution to the integral equation
\begin{equation}\label{eqn:dress_IE_app}
f^{\dr}(x)=f(x)+2\theta\int_{-\infty}^{\infty}\log|x-y|\,\varrho_{\beta}(y)\,f^{\dr}(y)\,dy .
\end{equation}

We apply this with $f(x)=x$ and $f(x)=1$. Writing $\varsigma_1^{\dr}$ for the dressing of $x$ and $\varsigma_0^{\dr}$ for the dressing of $1$, \eqref{eqn:dress_IE_app} gives
\begin{align}
\varsigma_1^{\dr}(x)
&=x+2\theta\int_{-\infty}^{\infty}\log|x-y|\,\varrho_{\beta}(y)\,\varsigma_1^{\dr}(y)\,dy, \label{eqn:s1_rel_app}\\
\varsigma_0^{\dr}(x)
&=1+2\theta\int_{-\infty}^{\infty}\log|x-y|\,\varrho_{\beta}(y)\,\varsigma_0^{\dr}(y)\,dy. \label{eqn:s0_rel_app}
\end{align}

Denoting
\begin{equation}\label{eqn:s0dr_app} 
\varsigma_1^{\dr}(x)=x+\theta R(x), \qquad
\varsigma_0^{\dr}(x)=1+\theta R_0(x),
\end{equation}
the proof of \cite[Lemma 3.6]{Agg25}, specifically (A.12) and (A.14) there, yields the bound
\[
|R(x)|+|R_0(x)|\le C\log(|x|+2),
\]
where the constant $C$ can be chosen independent of $\theta$ (by the observation in the first paragraph above, which implies that the constants in \cite[Lemma 3.6]{Agg25} can be chosen in a way that is uniform in $\theta$). This gives the first statement of the lemma.

Moreover, plugging the second in \eqref{eqn:s0dr_app} into the right hand side of \eqref{eqn:s0_rel_app} and using that
$2\int_{-\infty}^{\infty} \log|x-y|\,\varrho_\beta(y)\,dy=2 \log|x|+O(1)$ as $|x|\to\infty$,
we obtain as $|x| \rightarrow \infty$,
\begin{equation}\label{eqn:tildeR0_app}
R_0(x)=2\log|x|+O(1) + \theta \cdot O\left( \log(|x|+2) \right),
\end{equation}
where $O(1)$ is uniformly bounded in $x$ and $\theta$, and the constant in $O\left( \log(|x|+2) \right)$ does not depend on $\theta$.

Now let $x>0$ be sufficiently large. Splitting up the integral in \eqref{eqn:s1_rel_app}, differentiating, and bounding the boundary terms using Lemma \ref{lem:varrho_bd}, yields
\begin{flalign}\label{eqn:s1db}
\begin{aligned}
(\varsigma_1^{\dr})'(x) &=1+2\theta \int_{-\infty}^{x/10}\frac{1}{x-y}\,\varrho_{\beta}(y)\,\varsigma_1^{\dr}(y)\,dy
+ \theta O(e^{-c x^2}) \\
&+ 2\theta \int_{x/10}^{\infty} \log|x-y| \left( \varrho_{\beta}'(y)\varsigma_1^{\dr}(y) + \varrho_{\beta}(y)\,(\varsigma_1^{\dr})'(y) \right) dy \\
&= 1 +  \theta O( |x|^{-1} ) + \theta O(e^{-c x^2}) =  1 +  \theta O( |x|^{-1} ) .
\end{aligned}
\end{flalign}
To see the bound above, we have used the bound $|\varrho_{\beta}'(x)| \leq C (1+|x|) \varrho_{\beta}(x)$ (and the bound on $\varrho_{\beta}$) from (3.1) in \cite[Lemma 3.1]{Agg25}; the bounds $|\varsigma_1^{\dr}(y)| \leq C |y| + C \log(|y|+2)$ and $|(\varsigma_1^{\dr})'(y)| \leq C \log(|y|+2)$ obtained from Lemmas 3.5 and 3.6 of \cite{Agg25}. The same bound \eqref{eqn:s1db} can be obtained if $x < 0$ with $|x| $ large enough. Similarly, we obtain 
\begin{equation}\label{eqn:s0db}
|(\varsigma_0^{\dr})'(x)|\leq C \theta |x|^{-1}
\end{equation}
for $|x| $ large enough.

Finally, recalling $\ve=\varsigma_1^{\dr}/\varsigma_0^{\dr}$, increasing the implicit constants so that the bounds \eqref{eqn:s1db} and \eqref{eqn:s0db} may be written with $(|x|+1)^{-1}$ in place of $|x|^{-1}$, and then further increasing constants so the bounds hold for all $x \in \mathbb{R}$, we compute (using \eqref{eqn:s0db}, \eqref{eqn:s1db}, and \eqref{eqn:veff})
\begin{align}
\ve'(x)
&=\frac{(\varsigma_1^{\dr})'(x)\varsigma_0^{\dr}(x)-\varsigma_1^{\dr}(x)(\varsigma_0^{\dr})'(x)}{\varsigma_0^{\dr}(x)^2} \notag\\
&=\frac{\bigl(1+\theta O((|x|+1)^{-1})\bigr)\bigl(1+\theta  R_0(x)\bigr)
-\bigl(x+\theta  R(x)\bigr)\,\theta O((|x|+1)^{-1})}{(1+\theta  R_0(x))^2} \notag\\
&\geq \frac{1+\theta  R_0(x)- C \theta }{(1+\theta  R_0(x))^2}. \label{eqn:veprime_app}
\end{align}
Note that the last expression above satisfies 
$$\frac{1+\theta  R_0(x)- C \theta }{(1+\theta  R_0(x))^2} \geq \frac{1+2 \theta  \log|x|  - C' \theta -C' \theta^2 \log(|x|+2)}{(1+\theta  R_0(x))^2}$$
for $|x|$ large enough by  \eqref{eqn:tildeR0_app}. Therefore, choosing $\theta_0$ sufficiently small at the outset, the expression \eqref{eqn:veprime_app} implies $\ve'(x)>0$ for all $x$. Moreover, \eqref{eqn:veff} and \eqref{eqn:veprime_app} yield the stated lower bound on $\ve'$.
\end{proof}

\subsection{Proof of Lemma \ref{lem:second_der_bd}}
\label{ProofQDifference}

\begin{proof}[Proof of Lemma \ref{lem:second_der_bd}]

In what follows, we only address the case when $\varphi_t (k)> \varphi_t (i)$, as the proof when $\varphi_t (k) < \varphi_t (i)$ is entirely analogous. For $\theta < \theta_0 (\beta)$ sufficiently small, we have by \eqref{eq:alpha_def} that $\alpha>0$. Throughout this proof, we restrict to the event on which Lemmas \ref{lem:q_spacing}, \ref{lem:bd_lem} (with $A = \log N$), \ref{lem:loc_cent_dif}, and \ref{thm:Z_apriori}  all hold.

 Since $\varphi_0 (k) \in \llbracket -T^2, T^2 \rrbracket$ and $t \le T \log N$, Lemma \ref{lem:loc_cent_dif} yields $\varphi_t (k) \in \llbracket -2T^2, 2T^2 \rrbracket$, so we either have $\varphi_t (i), \varphi_t (k) \in \llbracket N_1 + T(\log N)^4, N_2 - T(\log N)^4 \rrbracket$ or $\varphi_t (k) - \varphi_t (i) \ge T (\log N)^5$. In both cases,  
\begin{flalign}
    \label{qkit0} 
    \begin{aligned} 
Q_k (t) - Q_i (t) = q_{\varphi_t (k)} - q_{\varphi_t (i)} (t) & \ge \displaystyle\frac{\alpha}{2} \cdot (\varphi_t (k)-\varphi_t (i)) - | \varphi_t (k)- \varphi_t (i)|^{1/2} \cdot (\log N)^2 \\
& \ge \displaystyle\frac{\alpha}{4} \cdot ( \varphi_t (k)- \varphi_t (i)) \ge \mathfrak{M} (\log N)^4,
\end{aligned} 
\end{flalign} 

\noindent where the first inequality holds by  Lemma \ref{lem:q_spacing}, and the second and third by the fact that $\varphi_t (k) - \varphi_t (i) \ge \mathfrak{M} (\log N)^5$. Recalling that $Q_i (t) = q_{\varphi_0 (i)} (t)$, it therefore remains to show that 
\begin{flalign}
    \label{qkit1} 
 q_{\varphi_0 (k)} (0) - q_{\varphi_0 (i)} + t(\ve (\lambda_k) - \ve (\lambda_i))  \ge \mathfrak{M} \log N.
\end{flalign} 

 Whenever $\varphi_0 (i) \in \llbracket N_1 + T(\log N)^6, N_2 - T(\log N)^6 \rrbracket$, we have 
\begin{flalign*}
    \big| Q_k(t) - Q_i (t) - \big( Q_k (0& ) - Q_i (0) + t(\ve (\lambda_k) - \ve (\lambda_i)) \big) \big| \\
    & = T^{1/2} \cdot |Z_{\varphi_0 (k)}^{\mathcal{Q}} (t) - Z_{\varphi_0 (i)}^{\mathcal{Q}} (t)| \le 2T^{1/2} (\log N)^{35} \le \mathfrak{M} \log N,
\end{flalign*}

\noindent where the first statement follows from the definition \eqref{eqn:quasi_part_fluct_intro} of the $Z_j^{\mathcal{Q}}$; the second from the second statement of Lemma \ref{thm:Z_apriori}; and the third from the fact that $\mathfrak{M} \ge T^{3/4}$ (recall Definition \ref{def:chi_def}). Together with \eqref{qkit0}, this implies \eqref{qkit1} if $\varphi_0 (i) \in \llbracket N_1 + T(\log N)^6, N_2 - T(\log N)^6 \rrbracket$. 

Otherwise, since $\varphi_0 (k) \ge \varphi_0 (i)$ and $\varphi_0 (k) \in \llbracket -T^2, T^2 \rrbracket$, we have (as $T = N^{1/100}$) that $\varphi_0 (k) - \varphi_0 (i) \ge T^2$. In this case, we obtain
\begin{flalign*}
 q_{\varphi_0 (k)} (0) - q_{\varphi_0 (i)}(0) + t(\ve (\lambda_k) - \ve (\lambda_i)) & \ge \displaystyle\frac{\alpha}{2} \cdot (\varphi_0 (k) - \varphi_0 (i)) - T \log N \cdot \big( |\ve (\lambda_k)| + |\ve (\lambda_i)| \big) \\
 & \ge \displaystyle\frac{\alpha T^2}{2} - CT (\log N)^2 \ge  \mathfrak{M} \log N,
\end{flalign*}

\noindent where the first statement holds by the second part of Lemma \ref{lem:q_spacing}, with the fact that $t \le T \log N$; the second by the facts that $\varphi_0 (k) - \varphi_0 (i) \ge T^2$ ; and the third by the fact that $|\lambda_j| \le \log N$ for each $j \in \llbracket 1, N \rrbracket$ (by Lemma \ref{lem:bd_lem}), together with the second part of Lemma \ref{lem:ve_tail}. This again implies \eqref{qkit1} and thus the lemma.
\end{proof}

\begin{proof}[Proof of Lemma \ref{lem:bd_num_termsU} (Outline)]
    The proof is very similar to that of Lemma \ref{lem:second_der_bd}, so we only outline the modifications.  Assuming without loss of generality that $\varphi_t(k)>\varphi_t(i)$, the proof of Lemma \ref{lem:second_der_bd} can be followed verbatim up to and including the display \eqref{qkit0}, with $U$ replacing $\mathfrak{M}$ and $\llbracket - T^{8/5} \log N, T^{8/5} \log N \rrbracket$ replacing $\llbracket -T^2, T^2 \rrbracket$. So, we obtain 
    \begin{flalign}
    \label{eqn:qkit0U} 
    \begin{aligned} 
Q_k (t) - Q_i (t) = q_{\varphi_t (k)} - q_{\varphi_t (i)} (t) & \ge \displaystyle\frac{\alpha}{2} \cdot (\varphi_t (k)-\varphi_t (i)) - | \varphi_t (k)- \varphi_t (i)|^{1/2} \cdot (\log N)^2 \\
& \ge \displaystyle\frac{\alpha}{4} \cdot ( \varphi_t (k)- \varphi_t (i)) \ge U (\log N)^4.
\end{aligned} 
\end{flalign} 

    Therefore, it remains to show  
\begin{equation}
    \label{eqn:qkit1} 
 \alpha \varphi_0 (k) - \alpha \varphi_0 (i)+ t(\ve (\lambda_k) - \ve (\lambda_i))  \ge U \log N.
\end{equation} 
We restrict to the event of Lemma \ref{lem:q_spacing} with $T = 0$, and to $\mathsf{BND}_{\L(0)}(\log N)$. In addition, we further restrict to the event on which the conclusions of Lemmas \ref{thm:Z_apriori} and \ref{lem:loc_cent_dif} hold. Set $k' = \varphi_0(k), i'= \varphi_0(i)$. If we also have $i' \in \llbracket - T^{8/5} (\log N)^2, T^{8/5} (\log N)^2 \rrbracket$, then $|q_{i'}(0)-\alpha i'| \leq T^{5/6}$ as well as $|q_{k'}(0)-\alpha k'| \leq T^{5/6}$. In this case, $\varphi_0 (i) \in \llbracket N_1 + T(\log N)^{20}, N_2 - T(\log N)^{20} \rrbracket$, so 
\begin{flalign*}
    \big|  Q_k(t) - Q_i (t)- \big(  Q_k(0) & -  Q_i(0) + t(\ve (\lambda_k) - \ve (\lambda_i)) \big) \big| \\
    & = T^{1/2} \cdot |Z_{\varphi_0 (k)}^{\mathcal{Q}} (t) - Z_{\varphi_0 (i)}^{\mathcal{Q}} (t)|  \le  T^{5/6} \le U \log N.
\end{flalign*}
and by the bounds on $|q_{k'}(0)-\alpha k'|$ and $|q_{i'}(0)-\alpha i'|$ from Lemma \ref{lem:q_spacing},
\begin{flalign*}
    \big|  Q_k(t) - Q_i (t)- \big( \alpha \varphi_0(k) & - \alpha \varphi_0(i) + t(\ve (\lambda_k) - \ve (\lambda_i)) \big) \big| \\
    & = T^{1/2} \cdot |Z_{\varphi_0 (k)}^{\mathcal{Q}} (t) - Z_{\varphi_0 (i)}^{\mathcal{Q}} (t)| + 2  T^{5/6} \le 3 T^{5/6} \le U \log N.
\end{flalign*}
The two displays above combined with \eqref{eqn:qkit0U} give the desired bound in this case. 

Otherwise we have $i'-k' \geq T^{8/5} \log N$. In this case, by Lemma \ref{lem:ve_tail}, 
\begin{flalign*}
\alpha \varphi_0 (k) - \alpha \varphi_0 (i) + t(\ve (\lambda_k) - \ve (\lambda_i)) & \ge \displaystyle \alpha \cdot (\varphi_0 (k) - \varphi_0 (i)) - T \log N \cdot \big( |\ve (\lambda_k)| + |\ve (\lambda_i)| \big) \\
 & \ge \displaystyle T^{8/5}  - CT (\log N)^2 \ge  U \log N.
\end{flalign*}
 This completes the proof.
\end{proof}

\subsection{Proof of Lemma \ref{lem:resolvent_lemma}}
\label{app:res_lem}

\begin{proof}
Let $K = 2 (\log N)^{12}$ for a moment. We begin by noting that by restricting to the event $\mathsf{E}_1 \coloneqq \mathsf{BND}_{\L(0)}(\log N)$,
which has overwhelming probability by Lemma \ref{lem:bd_lem}. On this event, the first summation in the left-hand side of \eqref{eqn:Hsum} is equal to  
$$\sum_{i=N_1}^{N_2} H(\Lambda_i, i) =\sum_{i=k_0-\floor{K T}}^{k_0+\floor{K T}} H(\Lambda_i, i),$$
since other summands are $0$ on this event by item \ref{item:H3} of Assumption \ref{ass:Hass}.

In all that follows, let us denote 
\begin{equation}\label{eqn:r0r1d}
    r_0 = \floor{(k_0-KT)/R}, \qquad r_1 = \floor{(k_0+KT)/R}-1.
\end{equation}

Next, we claim that on $\mathsf{E}_1$ defined above,
\begin{equation}\label{eqn:Hproof_1stbd}
 \left| \sum_{i=k_0-\floor{K T}}^{k_0+\floor{K T}} H(\Lambda_i, i) -  \sum_{r = r_0}^{r_1} \sum_{s=0}^{R-1} H(\Lambda_{r R+s}, rR )  \right| \leq 3 A  K RT/\mathfrak{M}.
\end{equation}
To verify \eqref{eqn:Hproof_1stbd}, first observe that the set of indices $i$ and $r R + s$ in the two sums, respectively, are in bijection, except for possibly a few boundary terms (e.g. if the interval $\llbracket r_0 R, (r_1+1) R-1 \rrbracket$ does not exactly equal the interval $\llbracket k_0 -\floor{ K T}, k_0 +\floor{ K T} \rrbracket$), which vanish anyways on $\mathsf{E}_1$ by item \ref{item:H3} of Assumption \ref{ass:Hass}. Now, \eqref{eqn:Hproof_1stbd} follows from the bound 
$$|H(\Lambda_{r R+s}, rR ) - H(\Lambda_{r R+s}, rR +s)| \leq A  R/\mathfrak{M},$$
which follows from item \ref{item:H2} of Assumption \ref{ass:Hass}. 

Thus, using \eqref{eqn:R_cond} together with \eqref{eqn:Hproof_1stbd},
\begin{equation}\label{eqn:reduc1}
\Bigg| \sum_{i=N_1}^{N_2} H(\Lambda_i, i)  - \sum_{r=r_0}^{r_1} \sum_{s=0}^{R-1} H(\Lambda_{r R+s}, rR ) \Bigg| \leq A T^{1/2-\mathfrak{c}}.
\end{equation}

\noindent By the spectral decomposition \eqref{gijuij} of the resolvent and Item \ref{item:H3} of Assumption \ref{ass:Hass}, we have
\begin{multline}\label{eqn:RHS1}
    \frac{1}{\pi} \int_{-\log N}^{\log N} \sum_{r = \lceil N_1/R \rceil}^{\lfloor N_2/R \rfloor -1} \sum_{s=0}^{R-1} H(E, rR )\Imaginary G_{r R+ s,r R+ s}(E+ \i \eta) dE \\
    = \frac{1}{\pi} \int_{-\log N}^{\log N} \sum_{k=N_1}^{N_2} \sum_{r=r_0}^{r_1}  \sum_{s=0}^{R-1} H(E, r R ) u_{\varphi_0^{-1}(k)}^2(r R+ s)\Imaginary (\Lambda_k- E-  \i \eta)^{-1} dE.
\end{multline}

Next we use eigenvector localization in the form of \cite[Lemma 5.1]{Agg25a} to see that on an event $\mathsf{E}_2$ of overwhelming probability, which we now restrict to,
\begin{equation}\label{eqn:eig_loc}
|u_{\varphi_0^{-1}(k)}^2(i)| \leq e^{-c (\log N)^5} \qquad \text{if } |i-k| \geq (\log N)^5,
\end{equation}
for all $k \in \llbracket N_1, N_2 \rrbracket$. In what follows, let 
$$I_0 \coloneqq \llbracket k_0-\floor{T K \log N}, k_0+\floor{T K \log N} \rrbracket.$$
Let $r(k) \coloneqq \floor{\frac{k}{R}} R$. We claim that, with overwhelming probability, we can with small error restrict the sum on the right-hand side of \eqref{eqn:RHS1} to $(r,s)$ satisfying $|Rr+s-k|<(\log N)^5$ and to $k \in I_0$, namely, we have
\begin{multline}\label{eqn:Important_HG_bd}
    \Big| \frac{1}{\pi} \int_{-\log N}^{\log N} \sum_{k=N_1}^{N_2}\sum_{r=r_0}^{r_1}  \sum_{s=0}^{R-1} H(E, r R ) u_{\varphi_0^{-1}(k)}^2(r R+ s)\Imaginary (\Lambda_k- E-  \i \eta)^{-1} dE\\
     - \frac{1}{\pi} \int_{-\log N}^{\log N} \sum_{k \in I_0} H(E, r(k) R ) \sum_{r,s : |R r + s -k| < (\log N)^5}  u_{\varphi_0^{-1}(k)}^2(r R+ s)\Imaginary (\Lambda_k- E-  \i \eta)^{-1} dE \Big| \\
     \leq A  T (\log N)^{15} R /\mathfrak{M} + 2 A N^4 e^{-c (\log N)^5}.
\end{multline}
To obtain this, we used \eqref{eqn:eig_loc} together with the bound
\begin{flalign}\label{eqn:togetloc}
\begin{aligned}
     & \Big| \frac{1}{\pi} \int_{-\log N}^{\log N} \sum_{k=N_1}^{N_2} \sum_{r,s : |R r + s -k| < (\log N)^5} H(E, r R ) u_{\varphi_0^{-1}(k)}^2(r R+ s)\Imaginary (\Lambda_k- E-  \i \eta)^{-1} dE\\
    & \quad - \frac{1}{\pi} \int_{-\log N}^{\log N} \sum_{k\in I_0} H(E, r(k) R ) \sum_{r,s : |R r + s -k| < (\log N)^5}   u_{\varphi_0^{-1}(k)}^2(r R+ s)\Imaginary (\Lambda_k- E-  \i \eta)^{-1} dE \Big| \\
    & \leq \frac{A R}{\mathfrak{M}} 
     \frac{1}{\pi}\int_{-\log N}^{\log N} \sum_{k\in I_0} \sum_{r,s : |R r + s -k| < (\log N)^5} u_{\varphi_0^{-1}(k)}^2(r R+ s) \Imaginary (\Lambda_k- E-  \i \eta)^{-1} dE \\
    & \quad  + \frac{R}{\mathfrak{M}} 
     \frac{1}{\pi}\int_{-\log N}^{\log N} \sum_{k\in I_0^{\complement}} \sum_{r,s : |R r + s -k| < (\log N)^5} |H(E, r R )| u_{\varphi_0^{-1}(k)}^2(r R+ s) \Imaginary (\Lambda_k- E-  \i \eta)^{-1} dE \\
     &\leq \frac{A R}{\mathfrak{M}} T (\log N)^{14}.
\end{aligned}
\end{flalign}
To see \eqref{eqn:togetloc}, first note that by Item \ref{item:H2} of Assumption \ref{ass:Hass}, for any $k \in \llbracket N_1, N_2 \rrbracket$ and $E \in [-\log N, \log N]$, we have
\begin{multline}
    \left|H(E, r(k) R ) \sum_{r,s : |R r + s -k| < (\log N)^5}  u_{\varphi_0^{-1}(k)}^2(r R+ s)- \sum_{r,s : |R r + s -k| < (\log N)^5} H(E, r R )  u_{\varphi_0^{-1}(k)}^2(r R+ s) \right| \\
    \leq \frac{A R}{\mathfrak{M}} \sum_{r,s : |R r + s -k| < (\log N)^5} u_{\varphi_0^{-1}(k)}^2(r R+ s).
\end{multline}
This provides the first inequality in \eqref{eqn:togetloc}. For the second inequality, we use
\[
\frac{1}{\pi}\int_{-\infty}^{\infty} \Imaginary (\Lambda_k-E-\i\eta)^{-1} \,dE = 1,
\]
the bound $\sum_{r,s : |R r + s -k| < (\log N)^5} u_{\varphi_0^{-1}(k)}^2(r R+ s)\le 1$, that $|I_0| \le 3TK (\log N) \le T(\log N)^{14}$, and Item \ref{item:H3} of Assumption \ref{ass:Hass}. 

Using \eqref{eqn:eig_loc} again, we then obtain
\begin{flalign}
\label{integral001}
\begin{aligned}
\frac{1}{\pi} & \int_{-\log N}^{\log N} \sum_{k\in I_0} H(E,  R r(k) ) \sum_{r,s : |R r + s -k| < (\log N)^5}   u_{\varphi_0^{-1}(k)}^2(r R+ s)\Imaginary (\Lambda_k- E-  \i \eta)^{-1} dE \\
&=
    \frac{1}{\pi} \int_{-\log N}^{\log N} \sum_{k\in I_0}  H(E, R r(k)) \sum_{i=N_1}^{N_2} u_{\varphi_0^{-1}(k)}^2(i)\Imaginary (\Lambda_k- E-  \i \eta)^{-1} dE + O(A e^{-c (\log N)^5}) \\
    &= \frac{1}{\pi} \int_{-\log N}^{\log N} \sum_{k\in I_0}   H(E, R r(k)) \Imaginary (\Lambda_k- E-  \i \eta)^{-1} dE +  O(A e^{-c (\log N)^5}) .
\end{aligned}
\end{flalign}
The error above satisfies $|O(A e^{-c (\log N)^5})| \leq A e^{-c (\log N)^5}$.

Next, a standard argument (which we omit here but see, for example, the end of the proof of \cite[Lemma D.3]{Agg25} in a very similar setting)  using the bound $|H(z,i)-H(z',i)| \leq A e^{-(\log N)^2}$, whenever $z,z' \in [-\log N,\log N]$ satisfy $|z-z'| \leq e^{-(\log N)^{5/2}}$, yields 
\begin{equation}
\label{sumintegral001}
  \left|  \frac{1}{\pi} \int_{-\log N}^{\log N} \sum_{k\in I_0}  \big(H(\Lambda_k,  R r(k)) - H(E, R r(k)) \big)\Imaginary (\Lambda_k- E-  \i \eta)^{-1} dE \right| \leq A e^{-c (\log N)^2}.
\end{equation}

\noindent Finally, on $\mathsf{E}_1$, 
$$\sum_{k \in I_0} H(\Lambda_k,  R r(k)) = \sum_{r = r_0}^{r_1} \sum_{s=0}^{R-1} H(\Lambda_{r R+s}, R r )$$
and hence we deduce the lemma by \eqref{eqn:reduc1}, \eqref{eqn:RHS1}, \eqref{eqn:Important_HG_bd}, \eqref{integral001}, \eqref{sumintegral001}.
\end{proof}

\subsection{Proof of Lemma \ref{lem:both_S_inv}}\label{app:Sinv_proofs}

We prove the two statements in Lemma \ref{lem:both_S_inv} separately.

		\begin{proof}[Proof of \eqref{eqn:Sbd} in Lemma \ref{lem:both_S_inv} (Outline)]

        The bound  \eqref{eqn:Sbd} implies that $\mathbf{S}$ is strictly diagonally dominant and thus invertible. Thus, it suffices to show  \eqref{eqn:Sbd}, whose proof is quite similar to (and in fact simpler than) the proof of \cite[Lemma 6.6]{Agg25}, so we will only outline it.

           It suffices to show that, on an overwhelmingly probable event $\mathsf{E}$, the following two inequalities hold for all $j \in \llbracket 2 l , 2 m  \rrbracket$. First,  
			\begin{flalign}
				\label{sjjlower} 
				\begin{aligned} 
					|S_{jj}| 
					&  \ge \Bigg| 2 \displaystyle\sum_{k=N_1}^{N_2} \mathfrak{l} (\Lambda_j - \Lambda_k) \cdot \chi' \big( Q_{\varphi_0^{-1} (j)} (t) - Q_{\varphi_0^{-1}(k)} (t) \big) + 1 \Bigg|  \\
					& \qquad \qquad  - 10 \mathfrak{M}^{-2} (\log N)^2 \cdot 2 T^{4/5} \cdot 2 \mathfrak{M} (\log N)^5  \ge 2 \alpha^{-1} \displaystyle\int_{-\infty}^{\infty} |\mathfrak{l} (\Lambda_j - x)| \varrho(x) d x + \displaystyle\frac{1}{4}.
				\end{aligned} 
			\end{flalign}

            \noindent Second,
            \begin{flalign}
				\label{integral000} 
				\begin{aligned} 
				2 \alpha^{-1} \displaystyle\int_{-\infty}^{\infty} |\mathfrak{l}(\Lambda_j-x)| \varrho (x) dx + \displaystyle\frac{1}{4}  \ge (2 + (\log N)^{-1}) \displaystyle\sum_{k=N_1}^{N_2} |\mathfrak{l}(\Lambda_j-\Lambda_k)| \cdot \chi' (\mathfrak{Q}_j-\mathfrak{Q}_j) + (\log N)^{-1}.
				\end{aligned}
			\end{flalign} 
			Clearly the two displays above imply \eqref{eqn:Sbd}.

			We will now define $\mathsf{E}$, the event on which we will be able to verify \eqref{sjjlower} and \eqref{integral000}. First, define events $\mathsf{E}_i$, $i =1,\dots, 6$ exactly as in the proof of \cite[Lemma 6.6]{Agg25}, but with $T$ there given by $T \log N$ here, with $s $ there given by $t$ here, with $\ell$ there given by $3 l $ here, and $m$ there by $ 3 m $ here. Second, let $\mathsf{E}_7$ denote the event on which the outcome of Lemma \ref{thm:Z_apriori} holds, but with $T$ there given by $T \log N$. We define $\mathsf{E} = \bigcap_{i=1}^7 \mathsf{E}_i$.
		
			Now we begin to verify \eqref{sjjlower}. The first lower bound in \eqref{sjjlower} is implied by the following facts. On $\mathsf{E}$, 
        \begin{align}
             \chi' (Q_{\varphi_0^{-1}(j)}(t)- Q_{\varphi_0^{-1}(k)}(t))  &= 0,   \qquad k \notin \llbracket 3 l , 3 m \rrbracket  \label{eqn:qqjqk2} \\
            \chi' (\mathfrak{Q}_j-\mathfrak{Q}_k) &= 0, \qquad k \notin \llbracket 3 l , 3 m \rrbracket  \label{eqn:mfqqjqk2} ,
        \end{align}
        and thus, on $\mathsf{E}$, for each $k$ in the support of either sum defining $S_{j j}$, we have $|\mathfrak{Q}_j-\mathfrak{Q}_k - (Q_{\varphi_0^{-1}(j)}(t)- Q_{\varphi_0^{-1}(k)}(t))| \leq 2 T^{4/5}$. This reasoning replaces the assumption in \cite[Lemma 6.6]{Agg25} bounding $|Q_{\varphi_0^{-1}(j)}(t) - \mathfrak{Q}_j|$; therefore, we may proceed in a way very similar to \cite[Equation (6.13)]{Agg25}  to obtain \eqref{sjjlower}. Moreover, as pointed out above \cite[Equation (6.14)]{Agg25}, \eqref{integral000} can be obtained by essentially reversing the reasoning which led to \eqref{sjjlower}  (and using our restriction to $\mathsf{E}_3$ in place of the restriction to $\mathsf{E}_2$).			 
		\end{proof}

\begin{proof}[Proof of \eqref{eqn:sinvbd} in Lemma \ref{lem:both_S_inv}] Both the statement of \eqref{eqn:sinvbd} and its proof are similar to \cite[Lemma 6.5]{Agg25}.
Throughout the proof, we assume we have restricted to the overwhelmingly probable  (by the first part of Lemma \ref{lem:both_S_inv}) event on which \eqref{eqn:Sbd} holds for all $j \in \llbracket 2 l , 2 m \rrbracket$.

Set $w = \mathbf{S}^{-1} v$.  Let $i_0  \in \llbracket l, m \rrbracket $ be the maximizing index of $|w_i|$ over the set $\llbracket l, m \rrbracket $. We may write, using \eqref{eqn:Sbd} and restricting to a smaller event of overwhelming probability, 
\begin{flalign}\label{eqn:v_bound}
\begin{aligned} 
    |v_{i_0}| & = \left| w_{i_0} S_{i_0 i_0} + \sum_{i \neq i_0} S_{i_0 i} w_i \right| \\
    & \geq  (\log N)^{-1} |w_{i_0}| + 2 \sum_{i=l- \floor{t (\log N)^4}}^{m + \floor{t (\log N)^4}} |\mathfrak{l}(\Lambda_{i_0} - \Lambda_i) \chi'(\mathfrak{Q}_{i_0}  - \mathfrak{Q}_i) | \big( (1 + (2\log N)^{-1})|w_{i_0}| - |w_{i}| \big).
    \end{aligned} 
\end{flalign}
Above we have used the fact that on the event $\bigcap_{i=1}^N \{ |\lambda_i| \leq \log N \} \subset \mathsf{BND}_{\L(0)}(\log N)$ (which we may restrict to by Lemma \ref{lem:bd_lem}), $|\ve(\Lambda_i)| \leq C \log N$ for all $i$ (by Lemma \ref{lem:ve_tail}), so that (since $i_0 \in \llbracket l, m \rrbracket$) there are no nonzero terms in either the sum defining $S_{i_0 i_0}$ or the sum $\sum_{i \neq i_0} S_{i_0 i} w_i$ outside of $\llbracket l -t (\log N)^4,  m + t (\log N)^4 \rrbracket \subset \llbracket 2 l, 2 m \rrbracket$. 

Now one of the following two possibilities can occur. 
\begin{enumerate}
    \item For all indices $i \in \llbracket l -t (\log N)^4,  m + t (\log N)^4 \rrbracket $ in the support of the sum, 
    $$|w_{i}| \leq (1 + (2\log N)^{-1}) |w_{i_0}|.$$
    In this case the bound \eqref{eqn:sinvbd} follows since the second term in the right hand side of \eqref{eqn:v_bound} has only nonnegative summands.
    \item There is an index $i_1 \in \llbracket l- t (\log N)^4, l-1 \rrbracket \cup \llbracket m+1 , m + t (\log N)^4 \rrbracket $ such that $|w_{i_1}| > (1 + (2\log N)^{-1}) w_{i_0}$. If so, we consider the index $i_1$ for which $|w_{i_1}|$ is maximal.
\end{enumerate}

We may now repeat the argument, examining \eqref{eqn:v_bound} with $i_1$ in place of $i_0$, and with $\llbracket l- t (\log N)^4, m + t (\log N)^4 \rrbracket$ in place of $\llbracket l, m \rrbracket$. Note that if we return to first case above, since now $|w_{i_1}| > |w_{i_0}|$, we obtain the bound \eqref{eqn:sinvbd}. Otherwise, we iterate the argument and obtain a new index $i_2$ with $w_{i_2} > (1 + (2\log N)^{-1}) w_{i_1} > (1 + (2 \log N)^{-1})^2 w_{i_0}$. Let $U = (\log N)^4$. If we end up repeating this at least $\lfloor U \rfloor +1$ times, then we obtain an index $i_{\lfloor U \rfloor +1} \in \llbracket 2 l , 2 m \rrbracket$ such that $|w_{i_0}| < (1 + (2\log N)^{-1})^{-U} |w_{i_{\lfloor U \rfloor +1}}| \leq e^{-(\log N)^3/8} |w_{i_{\lfloor U \rfloor +1}}| $. In this case the assumed bound $\max_{i \in \llbracket 2l, 2m \rrbracket} |w_i| \le (\log N)^{\mathfrak{C}}$ then implies the result. 
\end{proof}

\subsection{Proof of Lemma \ref{lem:inf_vol_bds}}
\label{sec:inf_vol_bd_proof}

\begin{proof}[Proof of Lemma \ref{lem:inf_vol_bds}]
Throughout the proof, suppose without loss of generality that $\alpha > 0$. We will provide an event on which the first bound in \eqref{eqn:q_i_inf_bd} holds, as the second bound will then follow by a symmetric argument (or by the $q_j \mapsto - q_{-j}$ symmetry of the Toda lattice).

   For any integer $J \geq 0$, let $(a_i^{(J)}(t), b_i^{(J)}(t))$, $i \in \llbracket -J, J \rrbracket $, denote the Toda lattice on $\llbracket -J, J \rrbracket$ with initial data $(a_i^{(J)}(0), b_i^{(J)}(0)) = (a_i(0), b_i(0))$ for all $i \in \llbracket -J, J \rrbracket$. So for any $J \geq 0$, we have a Toda lattice of length $2 J + 1$ which is coupled to the Toda lattice on $\mathbb{Z}$ given by $(a_i(t), b_i(t))$, $i \in \mathbb{Z}$. 

    Let us restrict to the following events. Let $N_0 =  K < N_1 = K^2 < N_2 = K^{2^2} <\cdots$. Let us denote the events
\begin{equation}\label{eqn:bd_allk}
  \mathsf{E}_1 \coloneqq \bigcap_{k \geq 0}  \bigcap_{t \geq 0} \mathsf{BND}_{\L_{N_k}(t)} ( \log N_k),
\end{equation}
    and 
    \begin{equation}\label{eqn:delta_q_bd_N0}
\mathsf{E}_2 \coloneqq \bigcap_{k \geq 0} \bigcap_{i, j \in \llbracket - N_k , N_k  \rrbracket} \{ |q_i^{(N_k)}(0) - q_j^{(N_k)}(0) - \alpha (i-j)| \leq \sqrt{|i-j|} (\log N_k)^2  \},
   \end{equation}
 and let $\mathsf{E}_3$ denote the event that for each $k \geq 0$,
 \begin{equation}\label{eqn:delta_abbds}
\sup_{t \in [0, T]} \max_{j \in \llbracket -N_k + 2 T \log N_k, N_k- 2 T \log N_k \rrbracket } \big(|a_j(t)-a_j^{(N_k)}(t)| +|b_j(t)-b_j^{(N_k)}(t)| \big) \leq e^{-T/5} 
   \end{equation}
 Note that on $\mathsf{E}_3$, for all $k \geq 0$, by \eqref{eqn:delta_abbds}, the coupling of initial conditions, and the fact that $\partial_s q_j(s) = b_j(s)$ and $\partial_s q_j^{(N_k)}(s) = b_j^{(N_k)}(s)$,
 \begin{equation}\label{eqn:fin_inf_q_bd_Nk}
 |q_i^{(N_k)}(s) - q_i(s)| \leq e^{-T/10} \qquad \text{for all $i \in \llbracket -N_k +2 T \log N_k, N_k- 2 T \log N_k\rrbracket$ and $s \in [0, T]$}.
   \end{equation}

By  Lemma \ref{lem:q_spacing}, Lemma \ref{lem:bd_lem}, and \cite[Proposition 2.5]{Agg25a} (with $\mathfrak{p} = 19/20$, $R = \log (3 N_k)$, $T$ there given by $T$ here, $N_1, N_2, N$ there given by $-N_k, N_k, 2 N_k + 1$ here, and $K$ there given by $2 T \log N_k$), together with a union bound, 
 \begin{equation}
     \mathbb{P}(\mathsf{E}_1 \cap \mathsf{E}_2 \cap \mathsf{E}_3) \geq 1- c^{-1} e^{-c (\log K)^2}.
 \end{equation}
 
 Now, suppose $i \leq - K$. Then, fix the smallest $k \geq 1$ such that $i \in \llbracket - N_k +  2 T \log N_k, N_k -  2 T \log  N_k \rrbracket$. Note that we necessarily have $|i| > N_{k-1} -   2 T \log N_{k-1}> N_{k-1} - 2 N_{k-1}^{9/10} \log N_{k-1} > N_{k-1}/2$. Therefore, for such $i$, $q_i^{(N_k)}(0) \leq - \alpha N_{k-1}/4$ by \eqref{eqn:delta_q_bd_N0}. Note also that $ T \log N_k = 2 T \log N_{k-1} \leq 2 N_{k-1}^{9/10} \log N_{k-1}$ by our assumption on $T$. 
 
 As a consequence of the above bounds together with \eqref{eqn:fin_inf_q_bd_Nk}, and because $|b_i^{(N_k)}(s)| \leq \log N_k$ for all $i \in \llbracket -N_k, N_k \rrbracket$ on $\mathsf{E}_1$, we have for $t \in [0, T]$ 
 \begin{equation}\label{eqn:finalbd0}
     q_i(t)  \leq   - \alpha N_{k-1}/4 + T \log N_k +e^{-T} \leq  - \alpha N_{k-1}/10 +e^{-T} < -T (\log T)^2.
 \end{equation}
 Since $i \leq - K$ was arbitrary, this completes the proof of the first two bounds. By a similar argument, in particular using \eqref{eqn:bd_allk} and \eqref{eqn:delta_q_bd_N0} for $k = 1$, for all $i \in \llbracket -K^2/2, K^2/2 \rrbracket$, we have $|a_i(t)|, |b_i(t)| \leq \log N_1 = 2 \log K + e^{-T^2/5} \leq 3 \log K$ for any $t \in [0, T]$. This proves the last statement.
\end{proof}

\section{Proofs of estimates for Gaussian limits}
\label{app:gauss_and219}

\subsection{Proof of coupling}
\label{app:gauss_couple}
\begin{proof}[Proof of Lemma \ref{lem:gauss_couple}]

By \eqref{eqn:cov_bd}, we have $(\sum_{i=1}^d |\Sigma_{i, j} - \tilde \Sigma_{i j}|^2 )^{1/2} \leq \gamma d^{1/2}$, for all $j \in \llbracket 1, d \rrbracket$. Therefore, denoting by $\| \cdot \|_{\op}$ the matrix operator norm, we have 
\begin{equation}\label{eqn:op_bd}
     \|\Sigma-\tilde \Sigma\|_{\op} \leq d \gamma.
\end{equation}

\noindent Letting $Z \in \mathbb{R}^d$ denote a standard Gaussian random variable (with covariance matrix given by the identity), define $W = \Sigma^{1/2} Z$ and $\tilde{W} = \tilde \Sigma^{1/2} Z$, which couples $(W, \tilde{W})$. We have 
\begin{equation}
    \|\Sigma^{1/2}-\tilde \Sigma^{1/2}\|_{\op}
    \leq  \|\Sigma-\tilde \Sigma\|_{\op}^{1/2},
\end{equation}

\noindent by \cite[Theorem X.1.1]{Bha97}. This, together with \eqref{eqn:op_bd}, yields
\begin{equation}
\|(\Sigma^{1/2}-\tilde \Sigma^{1/2}) Z\|_2 \leq \|(\Sigma^{1/2}-\tilde \Sigma^{1/2})\|_{\op} \cdot \|Z\|_{2} 
\leq  \|\Sigma -\tilde \Sigma\|_{\op}^{1/2} \cdot \|Z\|_{2} 
\leq (d \gamma)^{1/2} \cdot  \|Z\|_{2}.
\end{equation}

\noindent Hence, if $\|Z\|_2 \leq d^{1/2} \log d$, then 
\begin{equation}
\displaystyle\max_{i \in \llbracket 1, d \rrbracket} |W_i - \tilde{W}_i| \leq \|(\Sigma^{1/2}-\tilde \Sigma^{1/2}) Z\|_2  \leq \gamma^{1/2} d \log d.
\end{equation}
Since $\| Z \|_2$ is a Gaussian random variable with variance $d$, we have 
\begin{equation}
  \mathbb{P} (\|Z\|_2 > d^{1/2} \log d) \leq (2\pi)^{-1/2} e^{-(\log d)^2 / 2},
\end{equation}

\noindent from which we obtain the lemma. 
\end{proof}

\subsection{H\"older continuity and properties of limit processes}
\label{app:lim_holder}

Throughout this section, we recall the functions $\psi_{\Lambda,\kappa,\tau}$ from \eqref{eqn:logphi_intro} and $\phi_{\mathfrak{q}, \mathfrak{q}',\tau}^{[m]}$ from \eqref{eqn:phim_intro}. We seek to understand continuity properties of the Gaussian processes 
$\mathbb{R}^2 \times \mathbb{R}_{\ge 0} \ni s \mapsto \mathcal{W}^{\dr}(\psi_s)$ and $\mathbb{R}^2 \times \mathbb{R}_{\ge 0} \ni s \mapsto \mathcal{W}^{\dr}(\phi_s^{[m]})$. We begin with the below lemma bounding the $\mathcal{H}$-norms of $\psi$ and $\phi$.
\begin{lem}\label{lem:bscL2bds}
    For any $s = (\Lambda, \kappa, \tau) \in \mathbb{R}^2 \times \mathbb{R}_{\geq 0}$, we have 
    \begin{equation}\label{eqn:psHlem}
       \int_{-\infty}^{\infty}\|\psi_{s}(r,\cdot)\|_{\mathcal{H}}^2 d r < \infty.
    \end{equation}
    \noindent More specifically, there exists an absolute constant $\mathfrak{C}= \mathfrak{C}(\theta, \beta)> 0$ such that
    \begin{equation}\label{eqn:ps0Hlem}
     \sup_{\kappa \in \mathbb{R}}  \int_{-\infty}^{\infty}\|\psi_{(\Lambda, \kappa,\tau)}(r,\cdot)\|_{\mathcal{H}}^2 d r  \leq \mathfrak{C} \log(|\Lambda|+2)^2\left( |\Lambda| + 1\right) \tau  .
    \end{equation}
    In addition, for any $s = (\mathfrak{q}, \mathfrak{q}', \tau) \in \mathbb{R}^2 \times \mathbb{R}_{\geq 0}$ and $m \in \mathbb{Z}_{\geq 0}$, we have
    \begin{equation}\label{eqn:phHlem}
       \int_{-\infty}^{\infty}\|\phi_{s}^{[m]}(r, \cdot)\|_{\mathcal{H}}^2 d r < \infty.
    \end{equation}
\end{lem}
\begin{proof}
We have, by definition
\begin{flalign}\label{eqn:psi_sl2}
    \begin{aligned}
        \int_{-\infty}^{\infty}\|\psi_{s}(r,\cdot)\|_{\mathcal{H}}^2  dr  
        &= \int_{-\infty}^{\infty}\int_{-\infty}^{\infty} \left(
    \mathbbm{1}\{ \alpha (\kappa -   r ) > 0  \} 
    - 
    \mathbbm{1}\{ \alpha (\kappa  -   r )  +   \tau \left(\ve(\Lambda) - \ve(\lambda) \right) > 0\}
   \right)^2 \\
   & \qquad \qquad \times \left(\log |\Lambda -\lambda|\right)^2 \varrho(\lambda) dr d \lambda  \\
   &\leq  \tau \alpha^{-1} \int_{-\infty}^{\infty}\left(\log |\Lambda -\lambda|\right)^2 \left|  \ve(\Lambda) - \ve(\lambda) \right| 
   \varrho(\lambda) d \lambda .
    \end{aligned}
    \end{flalign}
Now, because the singularity of $(\log|\Lambda -\lambda|)^2$ is integrable, together with the bounds $|\ve(\lambda)| \leq C |\lambda|$, $\varrho(\lambda) \leq c^{-1} e^{-c \lambda^2}$ (from Lemma \ref{lem:ve_tail} and Lemma \ref{lem:varrho_bd}, respectively), \eqref{eqn:psi_sl2} is finite. This proves \eqref{eqn:psHlem}; the proof of \eqref{eqn:phHlem} is entirely analogous and thus omitted.

In order to show \eqref{eqn:ps0Hlem}, we begin with the right hand side of \eqref{eqn:psi_sl2}. Using Lemma \ref{lem:ve_tail} to bound $|\ve(\lambda)| \leq C (|\lambda|+1)$, we have
\begin{equation}
  \eqref{eqn:psi_sl2}  \leq  C \tau \alpha^{-1} \int_{-\infty}^{\infty}\left(\log |\Lambda-\lambda|\right)^2  (|\Lambda|+|\lambda|+2)  \varrho(\lambda) d \lambda,
\end{equation}

\noindent which is at most $C\log(|\Lambda|+2)^2 (|\Lambda|+1) \tau$, thereby confirming \eqref{eqn:ps0Hlem}.
\end{proof}

The next lemma will enable us to define a $(1/2)^-$ H\"older continuous modification of the process $\mathfrak{X}$ from Definition \ref{def:Z_intro}.
\begin{lem}\label{lem:var_bd}
There exists a constant $\mathfrak{C}>0$ such that the following holds. Consider the random process indexed by $s=(\Lambda,\kappa,\tau)$ defined by
$s\mapsto \mathcal{W}^{\dr}(\psi_{\Lambda,\kappa,\tau})$. Fix $\Lambda_0=\kappa_0=\tau_0 \geq 1$.
For any
$s,s'\in K:=[-\Lambda_0,\Lambda_0]\times[-\kappa_0,\kappa_0]\times[0,\tau_0]$,
\begin{flalign}
\int_{-\infty}^{\infty}\|\psi_s(r,\cdot)-\psi_{s'}(r,\cdot)\|_{\mathcal H}^2\,dr
& \le  \mathfrak{C}_0\,\Lambda_0 \kappa_0 \tau_0 \,\|s-s'\|_2; \label{eqn:psi_Hbd}\\ 
\mathbb{E}\Bigl[(\mathcal{W}^{\dr}(\psi_s)-\mathcal{W}^{\dr}(\psi_{s'}))^2\Bigr]
& \le  \mathfrak{C}\,\Lambda_0 \kappa_0 \tau_0\,\|s-s'\|_2. \label{eqn:psi_vbd}
\end{flalign}
\end{lem}

\begin{proof}
Since, denoting $\mathbf{D}$ as the dressing operator, $\|\mathbf{D} f(r,\cdot)\|_{\mathcal H}^2 \le C\|f(r,\cdot)\|_{\mathcal H}^2$, we have (c.f. Definition \ref{def:Wdress})
\begin{equation}\label{eqn:var_diff}
\operatorname{Var}\bigl(\mathcal{W}^{\dr}(\psi_s)-\mathcal{W}^{\dr}(\psi_{s'})\bigr)
\le C\int_{-\infty}^{\infty}\|\psi_s(r,\cdot)-\psi_{s'}(r,\cdot)\|_{\mathcal H}^2\,dr.
\end{equation}
Thus it suffices to bound the right-hand side of \eqref{eqn:var_diff}.

Write $s=(\Lambda_1,\kappa_1,\tau_1)$ and $s'=(\Lambda_2,\kappa_2,\tau_2)$. Then
\begin{align}\label{eqn:var_diff2}
\psi_{s'}(r,\lambda)-\psi_s(r,\lambda)
&=\bigl(\psi_{\Lambda_2,\kappa_1,\tau_1}(r,\lambda)-\psi_{\Lambda_1,\kappa_1,\tau_1}(r,\lambda)\bigr) \notag\\
&\quad +\bigl(\psi_{\Lambda_2,\kappa_2,\tau_1}(r,\lambda)-\psi_{\Lambda_2,\kappa_1,\tau_1}(r,\lambda)\bigr) \notag\\
&\quad +\bigl(\psi_{\Lambda_2,\kappa_2,\tau_2}(r,\lambda)-\psi_{\Lambda_2,\kappa_2,\tau_1}(r,\lambda)\bigr) \notag\\
&=:\Delta_\Lambda\psi_{\Lambda_1,\kappa_1,\tau_1}(r,\lambda)
+\Delta_\kappa\psi_{\Lambda_2,\kappa_1,\tau_1}(r,\lambda)
+\Delta_\tau\psi_{\Lambda_2,\kappa_2,\tau_1}(r,\lambda).
\end{align}
Consequently, using $(a+b+c)^2\le 3(a^2+b^2+c^2)$,
\begin{align}\label{eqn:var_diff3}
\int_{-\infty}^{\infty}\|\psi_s(r,\cdot)-\psi_{s'}(r,\cdot)\|_{\mathcal H}^2\,dr
&\le 3\int_{-\infty}^{\infty}\Bigl(
\|\Delta_\Lambda\psi_{\Lambda_1,\kappa_1,\tau_1}(r,\cdot)\|_{\mathcal H}^2 \notag\\
&\hspace{3.8em}+\|\Delta_\kappa\psi_{\Lambda_2,\kappa_1,\tau_1}(r,\cdot)\|_{\mathcal H}^2
+\|\Delta_\tau\psi_{\Lambda_2,\kappa_2,\tau_1}(r,\cdot)\|_{\mathcal H}^2
\Bigr)\,dr.
\end{align}

\textbf{Step 1: the $\kappa$-increment.}
From the definition \eqref{eqn:logphi_intro} of $\psi_{\Lambda,\kappa,\tau}$,
\begin{align}\label{eqn:var_dif4}
&\int_{-\infty}^{\infty}\|\Delta_{\kappa}\psi_{\Lambda_2,\kappa_1,\tau_1}(r,\cdot)\|_{\mathcal H}^2\,dr \notag\\
&\qquad =4 \int_{-\infty}^{\infty}\int_{-\infty}^{\infty}
\Bigl(
\Delta_\kappa \mathbbm{1}\{\alpha(\kappa_1-r)>0\}
-\Delta_\kappa \mathbbm{1}\{\alpha(\kappa_1-r)+\tau_1(\ve(\Lambda_2)-\ve(\lambda))>0\}
\Bigr)^2 \notag\\
&\hspace{9em}\times (\log|\Lambda_2-\lambda|)^2\,\varrho(\lambda)\,d\lambda\,dr \notag\\
&\qquad \le 8\int_{-\infty}^{\infty}\int_{-\infty}^{\infty}
\Bigl( \big( \Delta_\kappa \mathbbm{1}\{\alpha(\kappa_1-r)>0\}\big)^2
+\big(\Delta_\kappa \mathbbm{1}\{\alpha(\kappa_1-r)+\tau_1(\ve(\Lambda_2)-\ve(\lambda))>0\} \big)\Bigr)^2 \notag\\
&\hspace{9em}\times (\log|\Lambda_2-\lambda|)^2\,\varrho(\lambda)\,d\lambda\,dr \notag\\
&\qquad \le C|\kappa_2-\kappa_1|\int_{-\infty}^{\infty}(\log|\Lambda_2-\lambda|)^2\,\varrho(\lambda)\,d\lambda  \le C(\log(|\Lambda_2|+2))^2\,|\kappa_2-\kappa_1|.
\end{align}
Here we used that $\int_{-\infty}^{\infty} (\Delta_\kappa \mathbbm{1}\{\alpha(\kappa_1-r)>0\})^2 dr = \int_{-\infty}^{\infty} (\Delta_\kappa \mathbbm{1}\{\alpha(\kappa_1-r)+\tau_1(\ve(\Lambda_2)-\ve(\lambda))>0\})^2  dr = |\kappa_1 - \kappa_2|$ for any $\tau_1,\lambda \in \mathbb{R}$.

\textbf{Step 2: the $\Lambda$-increment.}
Similarly,
\begin{align}\label{eqn:var_diff5}
&\int_{-\infty}^{\infty}\|\Delta_{\Lambda}\psi_{\Lambda_1,\kappa_1,\tau_1}(r,\cdot)\|_{\mathcal H}^2\,dr \notag\\
&\qquad \le 8\int_{-\infty}^{\infty}\int_{-\infty}^{\infty}
\bigl(\log|\Lambda_2-\lambda|-\log|\Lambda_1-\lambda|\bigr)^2 \notag\\
&\hspace{6em}\times
\Bigl(
\mathbbm{1}\{\alpha(\kappa_1-r)>0\}
-\mathbbm{1}\{\alpha(\kappa_1-r)+\tau_1(\ve(\Lambda_2)-\ve(\lambda))>0\}
\Bigr)^2
\,\varrho(\lambda)\,d\lambda\,dr \notag\\
&\qquad\quad +8\int_{-\infty}^{\infty}\int_{-\infty}^{\infty}
\bigl(\log|\Lambda_1-\lambda|\bigr)^2
\Bigl(\Delta_\Lambda \mathbbm{1}\{\alpha(\kappa_1-r)+\tau_1(\ve(\Lambda_1)-\ve(\lambda))>0\}\Bigr)^2
\,\varrho(\lambda)\,d\lambda\,dr.
\end{align}
The second term is bounded by
\begin{equation}\label{eqn:var_diff5b}
C\,\tau_1\,|\Lambda_0|\,\log(|\Lambda_0|+2) \cdot |\Lambda_2-\Lambda_1| \cdot (\log(|\Lambda_0|+2))^2 = C \tau_1 |\Lambda_0| |\Lambda_2-\Lambda_1| (\log (|\Lambda_0|+2))^3,
\end{equation}
using $|\Lambda_i|\le \Lambda_0$ and the bound on $\ve'$ on $[-\Lambda_0,\Lambda_0]$ from Lemma \ref{lem:ve_tail}.

For the first term in \eqref{eqn:var_diff5}, we compute (again using Lemmas \ref{lem:varrho_bd} and \ref{lem:ve_tail})
\begin{align}\label{eqn:var_diff6}
&\int_{-\infty}^{\infty}\int_{-\infty}^{\infty}
\bigl(\log|\Lambda_2-\lambda|-\log|\Lambda_1-\lambda|\bigr)^2 \notag \\
&\qquad \times
\Bigl(
\mathbbm{1}\{\alpha(\kappa_1-r)>0\}
-\mathbbm{1}\{\alpha(\kappa_1-r)+\tau_1(\ve(\Lambda_2)-\ve(\lambda))>0\}
\Bigr)^2
\,\varrho(\lambda)\,d\lambda\,dr \notag\\
&\quad = C\,\tau_1\int_{-\infty}^{\infty}
\bigl(\log|\Lambda_2-\lambda|-\log|\Lambda_1-\lambda|\bigr)^2\,|\ve(\Lambda_2)-\ve(\lambda)|\,\varrho(\lambda)\,d\lambda  \le C\,\tau_1\,|\Lambda_0| \,|\Lambda_2-\Lambda_1| .
\end{align}
Combining \eqref{eqn:var_diff5}--\eqref{eqn:var_diff6} yields
\begin{equation}\label{eqn:delta_lam_bd}
\int_{-\infty}^{\infty}\|\Delta_{\Lambda}\psi_{\Lambda_1,\kappa_1,\tau_1}(r,\cdot)\|_{\mathcal H}^2\,dr
\le C\,\tau_1\,|\Lambda_0|\,\left(\log(|\Lambda_0|+2)\right)^3\,|\Lambda_2-\Lambda_1| .
\end{equation}

\textbf{Step 3: the $\tau$-increment.}
Finally, using Lemmas \ref{lem:varrho_bd} and \ref{lem:ve_tail}, we obtain 
\begin{align}\label{eqn:var_dif7}
\int_{-\infty}^{\infty} & \|\Delta_{\tau}\psi_{\Lambda_2,\kappa_2,\tau_1}(r,\cdot)\|_{\mathcal H}^2\,dr \notag\\
& =
\int_{-\infty}^{\infty}\int_{-\infty}^{\infty}
\Bigl(
\mathbbm{1}\{\alpha(\kappa_2-r)+\tau_2(\ve(\Lambda_2)-\ve(\lambda))>0\}
 \notag\\
& \qquad \qquad -\mathbbm{1}\{\alpha(\kappa_2-r)+\tau_1(\ve(\Lambda_2)-\ve(\lambda))>0\}
\Bigr)^2 \cdot (\log|\Lambda_2-\lambda|)^2\,\varrho(\lambda)\,d\lambda\,dr \notag\\
& \le C|\tau_2-\tau_1|\int_{-\infty}^{\infty}
|\ve(\Lambda_2)-\ve(\lambda)|\,(\log|\Lambda_2-\lambda|)^2\,\varrho(\lambda)\,d\lambda \notag\\
& \le C\,(\log(|\Lambda_0|+2))^2\,|\Lambda_0|\,|\tau_2-\tau_1|.
\end{align}

Putting \eqref{eqn:var_dif4}, \eqref{eqn:delta_lam_bd}, and \eqref{eqn:var_dif7} into \eqref{eqn:var_diff3}, we obtain
\begin{align}\label{eqn:psi_holder_pre}
\int_{-\infty}^{\infty}\|\psi_s(r,\cdot)-\psi_{s'}(r,\cdot)\|_{\mathcal H}^2\,dr
\le C \Lambda_0 \kappa_0 \tau_0 \Bigl(
|\Lambda_2-\Lambda_1|
+|\kappa_2-\kappa_1|
+|\tau_2-\tau_1|
\Bigr).
\end{align}
Each coordinate difference above is at most $\|s-s'\|_2$.
Thus \eqref{eqn:psi_holder_pre} implies \eqref{eqn:psi_Hbd}, and \eqref{eqn:psi_vbd} follows from \eqref{eqn:var_diff}.
\end{proof}

 Using the above results (with the Kolmogorov--Chentsov theorem), the following lemma defines a $\gamma$-H\"older continuous modification of the process $\mathfrak{X}$, for any $\gamma < 1/2$. For the choice $\gamma = 2/5$, it also gives quantitative probabilistic bounds on the $\gamma$-H\"older norm for $\mathfrak{X}$. 

\begin{lem}\label{lem:holder_norm}
  There exists a modification $\mathfrak{X}(\Lambda, \kappa, \tau)$ of the process $(\Lambda, \kappa, \tau) \mapsto \mathcal{W}^{\dr}(\psi_{\Lambda, \kappa, \tau})$, over $(\Lambda, \kappa, \tau) \in \mathbb{R}^2 \times \mathbb{R}_{\geq 0}$, which is almost surely locally $\gamma$-H\"older continuous  for any $\gamma \in (0, 1/2)$.
   Moreover, there exists a constant $\mathfrak{c} > 0$ such that the following holds. For any real number $A \ge 1$, the $2/5$-H\"older norm of the process $s \mapsto \mathfrak{X}(s) $, over $s \in K_0=[-A,A]^2 \times [0, A]$, has Gaussian tails. Specifically, for any real number $u \ge 0$, we have 
    $$\mathbb{P} \bigg(\sup_{s, s' \in K_0} |\mathfrak{X}(s) - \mathfrak{X}(s')| > u \|s-s'\|_2^{2/5} \bigg) \leq  e^{- \mathfrak{c} u^2 / A^8} .$$

\end{lem}

\begin{proof}

Fixing $\gamma \in (0,1/2)$, we use \cite[Theorem 2.53]{Str11} with the $X$ there equal to $\mathfrak{X}$ here, and the $(N,\alpha, r, C)$ there equal to $(3, \gamma, 1/2-\delta_0 - 3p^{-1}, C A^3 p^{1/2})$ here, for any large and small real numbers $p$ and $\delta$, respectively. By Lemma \ref{lem:var_bd}, these parameters satisfy the assumptions of \cite[Theorem 2.53]{Str11} (since $\mathbb{E}[\mathcal{N}^p] \le (C \sigma p^{1/2})^p$ for a centered normal random variable $\mathcal{N}$ of variance $\sigma$). Taking $p$ sufficiently large and $\delta_0$ sufficiently small, \cite[Theorem 2.53]{Str11}  yields of a locally $\gamma$-H\"older continuous modification $\mathfrak{X}(\Lambda, \kappa, \tau)$ of the process $(\Lambda, \kappa, \tau) \mapsto \mathcal{W}^{\dr}(\psi_{\Lambda, \kappa, \tau})$, verifying the first statement of the lemma. 

To show the second statement of the lemma, we apply the last bound in \cite[Theorem 2.53]{Str11}, with the same parameters as previously, except that we now set $\alpha = 2/5$. Observe that the factor $K$ appearing in \cite[Theorem 2.53]{Str11} is provided by the last line of the proof of that theorem, namely, denoting $\delta = \delta_0 + 3p^{-1}$, it is given by 
    \begin{flalign*} 
    K = \frac{2^{3(1+\frac{1}{p})+\frac{1}{2}-\delta -\frac{2}{5}}}{2^{\frac{1}{2}-\delta-\frac{2}{5}}-1} \left(\frac{2^{\frac{3}{2}-\delta}}{2^{\frac{1}{2}-\delta}-1}+3^{\frac{1}{2}}\right)\le C.
    \end{flalign*} 

\noindent Hence, a Markov estimate and then Lemma \ref{lem:var_bd} with \cite[Theorem 2.53]{Str11} and together yield
\begin{flalign}\label{eqn:markov_ineq}
\begin{aligned}
   \mathbb{P} \bigg(\sup_{s\neq s' \in K_0} \frac{|\mathfrak{X}(s)-\mathfrak{X}(s')|}{\|s-s'\|_2^{2/5}} > u \bigg) & \le u^{-p} \cdot \mathbb{E} \Bigg[ \bigg| \sup_{s\neq s' \in K_0} \frac{|\mathfrak{X}(s)-\mathfrak{X}(s')|}{\|s-s'\|_2^{2/5}} \bigg|^p \Bigg] \\  
   &\leq u^{-p}
  \left( C A^3 p^{1/2} \cdot (2 A)^{1/2-2/5-\delta_0} \right)^p \\
  &\leq e^{c^{-1} p  + 4 p \log A + (p \log p)/2 - p \log u }.
\end{aligned}
\end{flalign}

\noindent Setting $p = e^{-2/c-1} A^{-8} u^2$, the right side of \eqref{eqn:markov_ineq} is at most $e^{- c u^2 / A^6}$, which shows the lemma.
\end{proof}

We have an analogous lemma for the process $\mathfrak{Y}(\mathfrak{q}, \mathfrak{q}', \tau, m)= \mathcal{W}^{\dr}( \phi_{\mathfrak{q},\mathfrak{q}',\tau}^{[m]} ) $ given in Definition \ref{x2x0y2y0}. For each fixed $m$ we will redefine this as a continuous modification of the process $(\mathfrak{q}, \mathfrak{q}', \tau) \mapsto   \mathcal{W}^{\dr}( \phi_{\mathfrak{q},\mathfrak{q}',\tau}^{[m]} )$, where the family of functions $\phi_{\mathfrak{q},\mathfrak{q}',\tau}^{[m]}$ is defined in \eqref{eqn:phim_intro}.

\begin{lem}
   There exists a constant $\mathfrak{C}>0$ such that the following holds. For any real number $A > 1$, integer $m \in \llbracket 0, A^{1/10} \rrbracket$, and triples of real numbers $s=(\mathfrak{q}_1, \mathfrak{q}_1', \tau_1) \in \mathbb{R}^2 \times \mathbb{R}_{\geq 0}$ and $s' = (\mathfrak{q}_2, \mathfrak{q}_2', \tau_2),  \in \mathbb{R}^2 \times \mathbb{R}_{\geq 0}$, we have
\begin{flalign}
\int_{-\infty}^{\infty}\|\phi_{s}^{[m]}(r,\cdot ) - \phi_{s'}^{[m]} (r, \cdot)\|_{\mathcal{H}}^2 dr & \leq \mathfrak{C} A^{m} \|s-s'\|_2; \label{lem:phi_m_Hbd} \\
\mathbb{E} [( \mathcal{W}^{\dr}( \phi_{s}^{[m]} ) -\mathcal{W}^{\dr}( \phi_{s'}^{[m]} ))^2] & \leq \mathfrak{C} A^{m} \|s-s'\|_2. \label{lem:phi_m_varbd}
\end{flalign}
 
\end{lem}
\begin{proof}

    The estimate \eqref{lem:phi_m_varbd} follows from \eqref{lem:phi_m_Hbd}, as in the proof of Lemma \ref{lem:var_bd} (see \eqref{eqn:var_diff}), so it suffices to show \eqref{lem:phi_m_Hbd}. To that end, observe that 
    \begin{flalign}\label{eqn:phimvd}
\begin{aligned}
         \int_{-\infty}^{\infty}\| & \phi_{s}^{[m]}(r,\cdot ) - \phi_{s'}^{[m]} (r, \cdot)\|_{\mathcal{H}}^2 dr \\
         &\leq   \int_{-\infty}^{\infty}\int_{-\infty}^{\infty} \lambda^{2 m} \Bigg( \left( \mathbbm{1}\{-\alpha r + \mathfrak{q}_1 > 0\} -  \mathbbm{1}\{-\alpha r + \mathfrak{q}_2 > 0\} \right) \\
         & \qquad \qquad - \left( \mathbbm{1}\{-\alpha r  - \tau_1 \ve(\lambda) +  \mathfrak{q}_1' > 0\}-\mathbbm{1}\{-\alpha r  - \tau_2 \ve(\lambda) +  \mathfrak{q}_2' > 0\}\right) \Bigg)^2 d r  \varrho(\lambda) d \lambda \\
         &\leq C \int_{-\infty}^{\infty}\lambda^{2 m} \left( |\mathfrak{q}_1-\mathfrak{q}_2|+|\mathfrak{q}_1'-\mathfrak{q}_2'| + |\tau_1-\tau_2| |\ve(\lambda)|  \right)   \varrho(\lambda) d \lambda.
\end{aligned}
\end{flalign}
    Since $m \leq A^{1/10}$, we deduce from Lemmas \ref{lem:varrho_bd} and \ref{lem:ve_tail} the bounds
    $$\int_{-\infty}^{\infty}\lambda^{2 m}   \varrho(\lambda) d \lambda \leq C A^{m}, \qquad \int_{-\infty}^{\infty}\lambda^{2 m}  |\ve(\lambda)|  \varrho(\lambda) d \lambda \leq C A^{m}, $$
    and so the lemma follows from inserting these in \eqref{eqn:phimvd}.
\end{proof}

Next, we have the following lemma, which is the analog of Lemma \ref{lem:holder_norm} for the process $\mathfrak{Y}$. Its proof omitted, as it is entirely analogous to that of Lemma \ref{lem:holder_norm}, by replacing the use of Lemma \ref{lem:var_bd} there with Lemma \ref{lem:phi_m_varbd} here.

\begin{lem}\label{lem:holder_norm_2}

    There exists a constant $\mathfrak{c}>0$ such that the following holds. For any real number $A \ge 1$ and integer $m \in \llbracket 0, A^{1/10} \rrbracket$, the $2/5$-H\"older norm of the process $s \mapsto \mathfrak{Y}(s) $, over $s \in K=[-A,A] \times [0, A]$, has Gaussian tails. Specifically, for any real number $u \ge 0$, we have 
    $$\mathbb{P} \bigg(\sup_{s, s' \in K, s \neq s'} |\mathfrak{Y}(s, m) - \mathfrak{Y}(s', m)| > u \|s-s'\|_2^{2/5} \bigg) \leq \mathfrak{c}^{-1}e^{-\mathfrak{c} u^2/ A^{2 m+8}} .$$
\end{lem}

As a consequence of Lemmas \ref{lem:holder_norm} and \ref{lem:holder_norm_2}, we are able to prove Lemma \ref{lem:holder_cont}.

\begin{proof}[Proof of Lemma \ref{lem:holder_cont}]
    The first part of the lemma follows from Lemma \ref{lem:holder_norm} at $A = \log N$, upon taking the $u$ there to be $(\log N)^5$ here. The second follows from Lemma \ref{lem:holder_norm_2} at $A = \log N$, upon taking the $u$ there to by $(\log N)^{m+5}$ here. 
\end{proof}

\begin{proof}[Proof of Lemma \ref{lem:Xbd1}] 

Observe that \eqref{eqn:Xdef}, \eqref{eqn:ps0Hlem} and the boundedness of the operator $(1-\theta \mathbf{T} \varrhobetabf)^{-1}: \mathcal{H} \rightarrow \mathcal{H}$ together yield (recalling that $\tau \in [0, A_2]$)
$$ \displaystyle\sup_{|\Lambda| \le A_1} \var \mathfrak{X} \big(\Lambda, (Q-\tau \ve(\Lambda))/\alpha, \tau \big)  \leq C A_1^2 A_2. $$
     
     \noindent Hence, defining the event $\mathsf{E}_{\Lambda} \coloneqq \{ |\mathfrak{X}(\Lambda, (Q-\tau \ve(\Lambda))/\alpha, \tau)| >  A_1^2 A_2 / 2 \}$, we have for $|\Lambda| \leq  A_1$ that
    $$\mathbb{P}(\mathsf{E}_{\Lambda}) \leq c^{-1} e^{-c A_1^2 A_2}. $$

Let $\mathcal{M} \coloneqq [-A_1, A_1] \cap (A_1^{-100} \mathbb{Z})$ be a fine mesh approximating $[-A_1, A_1] $, and define $\mathsf{E}_0 = \bigcap_{\Lambda \in \mathcal{M}} \mathsf{E}_{\Lambda}$, which satisfies $\mathbb{P}[\mathsf{E}_0] \le C e^{-cA_1^2 A_2}$. Then, for any $\Lambda \in [-A_1, A_1]$, there is some $\Lambda' \in \mathcal{M}$ such that $|\Lambda-\Lambda'| \leq A_1^{-100}$; moreover, by Lemma \ref{lem:ve_tail}, $\tau |\ve(\Lambda) - \ve(\Lambda')| \leq C A_1^{-90}$. Using these bounds, together with Lemma \ref{lem:holder_norm} (with the $(A,u)$ there equal to $(A_1, A_1^5)$ here), we obtain an event $\mathsf{E}_1$ of probability at least $1-c^{-1} e^{- c A_1^2}$ on which

\begin{equation}\label{eqn:Xdfbd}
\big|\mathfrak{X} \big(\Lambda, (Q-\tau \ve(\Lambda))/\alpha, \tau \big) -
    \mathfrak{X} \big(\Lambda', (Q-\tau \ve(\Lambda'))/\alpha, \tau \big)  \big| \leq A_1^{-30}.
\end{equation} 

\noindent Hence,
\begin{equation*}
  \sup_{|\Lambda| \leq A_1} \big|\mathfrak{X} \big(\Lambda, (Q-\tau \ve(\Lambda))/\alpha, \tau \big) \big| \leq \displaystyle\frac{1}{2} \cdot A_1^2 A_2 + A_1^{-30} \le A_1^2 A_2, \qquad  \text{on }  \mathsf{E}_0 \cap \mathsf{E}_1.
\end{equation*}

 \noindent This, together with the fact that $\mathbb{P}[\mathsf{E}_0 \cap \mathsf{E}_1] \ge 1 - C e^{-cA_1^2}$, completes the proof.
\end{proof}

\begin{proof}[Proof of Lemma \ref{lem:psi_phi_L2_intro}]
  The first part of the lemma follows from \eqref{eqn:psHlem} in Lemma \ref{lem:bscL2bds}. To show the second, set $A_2 = \max(\tau , 1)$ and observe, by Lemma \ref{lem:Xbd1}, 
    \begin{equation}
        \sum_{A \in \mathbb{Z}_{>0} ,A \geq A_2 } \mathbb{P}\left(  \sup_{|\Lambda| \leq A}  \big|\mathfrak{X} \big(\Lambda, (Q-\tau \ve(\Lambda))/\alpha, \tau \big) \big| > \tau A^2  \right) \leq  \sum_{A \in \mathbb{Z}_{>0}} C e^{-c A^2} < \infty.
    \end{equation}

    \noindent The lemma then follows from the Borel--Cantelli lemma, which implies that almost surely, for large enough $|\Lambda| $, we have $|\mathfrak{X}(\Lambda, (Q-\tau \ve(\Lambda))/\alpha, \tau)| \leq 10 \tau \Lambda^2$; this in turn implies that $\Lambda \mapsto \mathfrak{X}(\Lambda, (Q-\tau \ve(\Lambda))/\alpha, \tau)$ is an element of $\mathcal{H}$. 
\end{proof}

\section{Proofs of concentration estimates}
\label{app:mcDiarmid_and_proofs}

\subsection{Minami estimate}
\label{app:minami}

In this section, we record a form of the Minami estimate \cite{Min96}, as well as a consequence of it that will be useful for us.

\begin{lem}\label{lem:min}

There exists a constant $\mathfrak{C}>1$ such that the folllowing holds. Let the integers $(N_1, N_2; N; \zeta)$ and the eigenvalues $(\Lambda_{N_1}, \Lambda_{N_1+1}, \ldots , \Lambda_{N_2})$ of the Lax matrix $\L(0)$ be as in Assumption \ref{ass:NT_assumption} and Definition \ref{def:quasi_intro}. Let $\mathfrak{J} = \{i_0,i_0+1,\dots, i_1\} \subseteq \llbracket N_1, N_2 \rrbracket$ be a set of consecutive indices, and let $I \subset \mathbb{R}$ be an interval. Let $\L^{(\mathfrak{J})} = [L_{i,j}(0)]_{i,j\in \mathfrak{J}}$ denote the principal submatrix of $\L(0)$ indexed by $\mathfrak{J}$, and denote its eigenvalues by $\lambda_1^{(\mathfrak{J})} \ge \lambda_2^{(\mathfrak{J})} \ge \cdots \ge \lambda_{|\mathfrak{J}|}^{(\mathfrak{J})}$.
Then,
	    \begin{equation}\label{eqn:minami1}
	    \mathbb{E} \big[\# \{  j \in \llbracket 1, |\mathfrak{J}| \rrbracket : \lambda_j^{(\mathfrak{J})} \in I  \} \big] \leq \mathfrak{C} |\mathfrak{J}| |I| 
	    \end{equation}
	    and 
\begin{equation}\label{eqn:minami2}
	    \mathbb{E} \big[\# \{ (j,k) \in \llbracket 1, |\mathfrak{J}| \rrbracket^2: \lambda_j^{(\mathfrak{J})}, \lambda_k^{(\mathfrak{J})} \in I; j \ne k  \} \big] \leq \mathfrak{C} (|\mathfrak{J}| |I|)^2.
	    \end{equation}
	    
	\noindent Moreover, suppose $R \in [1, N ]$ is a real number, $\mathcal{C} \subseteq \llbracket N_1, N_2 \rrbracket$ is any consecutive set of indices of size $|\mathcal{C}| \leq  R$, and $I$ satisfies $|I| \leq (\log N)^4 R^{-1}$. Recalling Definition \ref{def:quasi_intro}, set
	\[
	N_{I, \mathcal{C}} \coloneqq \# \{ i \in \mathcal{C} : \Lambda_i \in I \},
	\]
	so $N_{I,\mathcal{C}}$ counts the eigenvalues of $\L(0)$ lying in $I$ whose $\zeta$-localization centers lie in $\mathcal{C}$. Then, we have with overwhelming probability that $N_{I, \mathcal{C}} \le   (\log N)^7$.

\end{lem}

\begin{proof}

    We only prove \eqref{eqn:minami2} in detail, as the proof of \eqref{eqn:minami1} is very similar (obtained by replacing the below use of \cite[Equation (2.65)]{Min96} by \cite[Equation (2.22)]{Min96}).

    Suppose $I = [a, b]$, denote $|a-b| = 2\delta$, let $m = (a+b)/2$, and let $z = m + \i \delta = a + \delta + \i \delta$. Denote by $\G^{(\mathfrak{J})}(z) = (\L^{(\mathfrak{J})} - z)^{-1}$ the resolvent of the matrix $\L^{(\mathfrak{J})}$. The estimate \eqref{eqn:minami2} can be proved by first noting that with our choice of $z$, for $j \neq k$,
    $$\mathbbm{1}_{\lambda_j^{(\mathfrak{J})} \in I} \cdot \mathbbm{1}_{\lambda_k^{(\mathfrak{J})}\in I} \leq 4\delta^2 \cdot \Imaginary \left[ \big( \lambda_j^{(\mathfrak{J})} - z\big)^{-1} \right] \cdot \Imaginary \left[\big( \lambda_k^{(\mathfrak{J})} - z \big)^{-1} \right]. $$ 
    Then, using this, the left hand side of \eqref{eqn:minami2} can be upper bounded by 
    \begin{flalign} 
    \label{eqn:2ptMin}
    \begin{aligned} 
      4 \delta^2 \cdot \mathbb{E} \Bigg[ & \sum_{\substack{j, k \in \llbracket 1, |\mathfrak{J}| \rrbracket \\ j\neq k}}  \Imaginary \big( \lambda_j^{(\mathfrak{J})} - z\big)^{-1} \Imaginary \big( \lambda_k^{(\mathfrak{J})} - z \big)^{-1}  \Bigg] \\
    &= 4 \delta^2 \cdot \mathbb{E} \left[ \big( \Tr \Imaginary \G^{(\mathfrak{J})}(z) \big)^2  - \Tr \big( \Imaginary \G^{(\mathfrak{J})}(z) \big)^2 \right] \\
      & =   4  \delta^2  \sum_{i, k \in \mathfrak{J}; i\neq k } \mathbb{E} \left[
  \det\!
  \begin{bmatrix}
    \Imaginary \G^{(\mathfrak{J})}_{ii}(z) & \Imaginary \G^{(\mathfrak{J})}_{ik}(z) \\
    \Imaginary \G^{(\mathfrak{J})}_{ik}(z) & \Imaginary \G^{(\mathfrak{J})}_{kk}(z) 
  \end{bmatrix}
\right] \leq C \delta^2 |\mathfrak{J}|^2 = \displaystyle\frac{C}{4} \cdot |I|^2 |\mathfrak{J}|^2,
\end{aligned} 
     \end{flalign}

     \noindent where to obtain the final inequality above, we used \cite[Equation (2.65)]{Min96}\footnote{The results there are stated for a different random matrix, but it is quickly verified that its proofs apply for the one we study.} to bound the expectation of each $2 \times 2$ determinant. This shows \eqref{eqn:minami2}.

     It remains to show that $N_{I,\mathcal{C}} \le (\log N)^7$ with overwhelming probability. This holds if $R \leq (\log N)^7$, so assume instead that $R > (\log N)^7$. We decompose $\mathcal{C} = \bigsqcup_{m=1}^M \mathfrak{J}_m$ as a disjoint union of pairwise disjoint consecutive index intervals $\mathfrak{J}_m = \llbracket i_{0,m}, i_{1,m} \rrbracket$, ordered from left to right, with $(\log N)^5/3 \leq |\mathfrak{J}_m| \leq (\log N)^5$. Since $R > (\log N)^7$, we must have $M > (\log N)^2$.

Consider the enlargement $\mathfrak{J}_m'$ of these chunks by $3 (\log N)^3$ on both ends; specifically, if $\mathfrak{J}_m = \llbracket i_{0, m}, i_{1,m} \rrbracket$ then $\mathfrak{J}_m' = \llbracket i_{0, m}-  3 (\log N)^3, i_{1,m} +3 (\log N)^3 \rrbracket$. Let $\L^{(\mathfrak{J}_m')}$ denote the principal submatrix of $\L(0)$ indexed by $\mathfrak{J}_m'$, and let $\widehat{\L}^{(m)}$ denote the $N\times N$ matrix obtained from $\L(0)$ by setting to zero entries outside of $\mathfrak{J}_m'$; $\widehat{L}_{i,j}^{(m)} = L_{i,j}(0)$ if $i, j \in \mathfrak{J}_m'$, and otherwise $\widehat{L}_{i,j}^{(m)} = 0$. Also set
\[
N_{I', \mathfrak{J}_m'}^{(m)} \coloneqq \#\{ j \in \llbracket 1, |\mathfrak{J}_m'| \rrbracket : \lambda_j^{(\mathfrak{J}_m')} \in I' \},
\qquad
I' = [a-N^{-1}, b+N^{-1}].
\]
First, if we restrict to the event that $N_{I, \mathcal{C}} >  (\log N)^7 $; to the event of overwhelming probability coming from the intersection of the $M$ events furnished by Lemma \ref{lem:Lax_eig_coupling} (with $\tilde{\L}$ given by $\widehat{\L}^{(m)}$, $m=1,\dots, M$); and to $\mathsf{SEP}_{\L(0)}(e^{-(\log N)^2})$, then for at least one of $p \in \{ 0, 1 \}$ we have 
\begin{equation}\label{eqn:event_to_bd}
\sum_{r=1}^{\lfloor (M-p)/2 \rfloor} N_{I', \mathfrak{J}_{2r+p}'}^{(2r+p)}  > \displaystyle\frac{1}{2} \cdot (\log N)^{7}.
\end{equation}

Indeed, let us restrict to the events described above, and let $m = 2 r + p$ for some $r \in \llbracket 1,(M-p)/2 \rrbracket $. By Lemma \ref{lem:Lax_eig_coupling}, for each $\Lambda_i \in I$ with $i \in \mathfrak{J}_{m}$ there is a corresponding eigenvalue $\tilde{\lambda} \in \eig(\widehat{\L}^{(m)})$ such that $|\tilde{\lambda}-\Lambda_i| \le e^{-c(\log N)^3}$, so in particular $\tilde{\lambda} \in I'$, and $i$ is an $\zeta N^{-1}$-localization center of $\tilde{\lambda}$. Since $\widehat{\L}^{(m)}$ vanishes outside of $\mathfrak{J}_m'$, every such eigenvalue $\tilde{\lambda}$ is an eigenvalue of the principal block $\L^{(\mathfrak{J}_m')}$, and therefore contributes to $N_{I', \mathfrak{J}_m'}^{(m)}$. Moreover, by the restriction to $\mathsf{SEP}_{\L(0)}(e^{-(\log N)^2})$ and Lemma \ref{lem:sep_lem}, two different  $\Lambda_i, \Lambda_j$ cannot correspond to the same eigenvalue of $\widehat{\L}^{(m)}$. As a consequence, for each eigenvalue of $\L(0)$ whose $\zeta$-localization center is in $\mathfrak{J}_{2 r + p}$ for some $r \in \llbracket 1,(M-p)/2 \rrbracket$, and which contributes to the count $N_{I, \mathcal{C}}$, there is at least one unique eigenvalue contributing to $\sum_{r=1}^{\lfloor (M-p)/2 \rfloor} N_{I', \mathfrak{J}_{2r+p}'}^{(2r+p)}$. This implies \eqref{eqn:event_to_bd}.

Therefore, it suffices to show that for $p \in \{ 0, 1 \}$ the above event \eqref{eqn:event_to_bd} occurs with probability at most $c^{-1}e^{-c(\log N)^2}$. We prove this below in the case $p = 0$, as the proof is identical for $p=1$. 

Now, note that the enlarged intervals $\{\mathfrak{J}_{2r}'\}_{r=1}^{\lfloor M/2 \rfloor}$ are pairwise disjoint, since the original consecutive intervals have length at least $(\log N)^5/3 > 6(\log N)^3$. Hence the random variables $\{N_{I', \mathfrak{J}_{2r}'}^{(2 r)}\}_{r=1}^{\lfloor M/2 \rfloor}$ are independent. By \eqref{eqn:minami1}, $\mathbb{E} [N_{I', \mathfrak{J}_{2r}'}^{(2 r)}] \leq C |I'| |\mathfrak{J}_{2 r}'|$, and in particular 
\begin{flalign*} 
\mu \coloneqq \mathbb{E} \Bigg[ \sum_{r=1}^{\lfloor M/2 \rfloor} N_{I', \mathfrak{J}_{2r}'}^{(2 r)} \Bigg] \leq C |I'| \sum_{r=1}^{\lfloor M/2 \rfloor} |\mathfrak{J}_{2 r}'|  \leq 2C (\log N)^4 .
\end{flalign*} 

\noindent Moreover, if we set $X_r \coloneqq N_{I', \mathfrak{J}_{2r}'}^{(2 r)}$, then
\[
\var(X_r)=\mathbb{E}[X_r^2]-\mathbb{E}[X_r]^2\le \mathbb{E}[X_r^2]
=\mathbb{E}[X_r(X_r-1)] + \mathbb{E}[X_r].
\]
By \eqref{eqn:minami1} and \eqref{eqn:minami2}, this yields
\[
\var(X_r)\le C |I'| |\mathfrak{J}_{2r}'| + C (|I'| |\mathfrak{J}_{2r}'|)^2.
\]
Hence
\[
\sigma^2 \coloneqq \sum_{r=1}^{\lfloor M/2 \rfloor} \var N_{I', \mathfrak{J}_{2r}'}^{(2 r)}
\leq C |I'| \sum_{r=1}^{\lfloor M/2 \rfloor} |\mathfrak{J}_{2r}'|
+ C |I'|^2 \sum_{r=1}^{\lfloor M/2 \rfloor} |\mathfrak{J}_{2r}'|^2.
\]
Since $\sum_{r=1}^{\lfloor M/2 \rfloor} |\mathfrak{J}_{2r}'| \le 2R$, $\max_r |\mathfrak{J}_{2r}'| \le 2(\log N)^5$, and $|I'|\le 2(\log N)^4 R^{-1}$, we further obtain
\[
\sigma^2 \le C(\log N)^4 + C |I'|^2 \max_r |\mathfrak{J}_{2r}'| \sum_{r=1}^{\lfloor M/2 \rfloor} |\mathfrak{J}_{2r}'|
\le C(\log N)^4 + 16 CR^{-1} (\log N)^{13}
\le 20C(\log N)^6,
\]
where in the last step we used $R>(\log N)^7$.
In addition, we have the deterministic bound $|X_r| \leq 2 (\log N)^5$, and thus $|X_r-\mathbb{E}[X_r]|\le C(\log N)^5$. Since $\mu \le 2C(\log N)^4$, for $N$ large enough we also have
\[
\frac{1}{2}(\log N)^7-\mu \ge c(\log N)^7.
\]
So, by Bernstein's inequality,

\begin{align*}
  \mathbb{P} \bigg( \sum_{r=1}^{\floor{M/2}}  N_{I', \mathfrak{J}_{2r}'}^{(2r)} >  \displaystyle\frac{1}{2} \cdot (\log N)^{7} \bigg)
  &\leq \mathbb{P} \bigg( \bigg|\sum_{r=1}^{\floor{M/2}} \Big( N_{I', \mathfrak{J}_{2r}'}^{(2r)} - \mathbb{E}\big[N_{I', \mathfrak{J}_{2r}'}^{(2r)}\big] \Big) \bigg| > c (\log N)^7 \bigg) \\
  &\leq  C \exp \Bigg( -c \min \bigg( \frac{(\log N)^{14}}{\sigma^2}, \frac{(\log N)^7}{(\log N)^5} \bigg) \Bigg) \leq C e^{-c  (\log N)^2},
\end{align*}

\noindent which shows the lemma.
\end{proof}

\subsection{Concentration bounds}

\label{ProofEstimateConcentration}

Throughout this appendix, we will make use of the following concentration estimate for functions of independent random variables, which is a variant of McDiarmid's inequality.
\begin{lem}[Lemma 5.2 of \cite{Agg25}] \label{lem:McD_bd}
    Let $N, m \ge 1$ be integers and $F : \mathbb{R}^N \rightarrow \mathbb{R}$ be a function. Let $x_1,\dots, x_N$ be mutually independent random variables, and let $\mathcal{J}_1\sqcup \mathcal{J}_2 \sqcup \cdots \sqcup \mathcal{J}_m$ be a partition of $\llbracket 1, N \rrbracket$ into a disjoint union of $m$ subsets. Denote by $(\tilde{x}_1,\dots, \tilde{x}_N)$ an independent copy of $(x_1,\dots, x_N)$. Further set $x(\mathcal{J}_k) = (x_j)_{j \in \mathcal{J}_k}$, $\tilde x (\mathcal{J}_k) = (\tilde{x}_j)_{j \in \mathcal{J}_k}$, and $x(\mathcal{J}_k^{\complement}) = (x_j)_{j \notin \mathcal{J}_k}$. Let $p \in (0, 1)$ and $A_1, A_2, \ldots , A_m \ge 0$ be real numbers, and suppose for each $k \in \llbracket 1, m \rrbracket$ that replacing $x(\mathcal{J}_k) \rightarrow \tilde{x} (\mathcal{J}_k)$ in $F$ incurs an error of at most $A_k$ with probability at least $1-p$, namely,
    \begin{equation}
    \mathbb{P} \big( \big|F(x(\mathcal{J}_k^{\complement}), x(\mathcal{J}_k)) - F(x(\mathcal{J}_k^{\complement}), \tilde{x}(\mathcal{J}_k)) \big| \geq A_k \big) \leq p.
    \end{equation}
    
    \noindent Then, denoting  
    \begin{equation}\label{eqn:mcd_S_U}
    S = \sum_{k=1}^m A_k^2 \qquad U = \mathbb{E}[F(x)^2]^{1/2}.
    \end{equation}
    
    \noindent we have for any real number $R \geq 0$ that 
    \begin{equation}
    \mathbb{P} \big(\big| F(x) - \mathbb{E}[F(x)]\big| \geq  R S^{1/2} + 2 m^{1/4} p^{1/2} U \big) <  2 m p^{1/2} + 2 e^{-R^2/4} .
    \end{equation}
  
\end{lem}

Now, we can prove Lemma \ref{lem:no_q_conc}. Since the proofs of Items \ref{item:conc2} and \ref{item:conc3} are very similar to the proof of Item \ref{item:conc1}, we omit their proofs. 

\begin{proof}[Proof of Lemma \ref{lem:no_q_conc}]

We only prove the first part of Lemma \ref{lem:no_q_conc}, as the proof of the second and third are entirely analogous. We further assume throughout that $1 \leq S \leq t$, as otherwise the proof is similar but in fact simpler (and closely follows that of \cite[Lemma 4.5]{Agg25}).

 We will apply the concentration bound, Lemma \ref{lem:McD_bd}. Define the integer $P = \lfloor t/S \rfloor  \geq 1$. We partition indices into $M$ chunks of size $P$, together with one additional chunk containing all remaining indices. Define
\begin{equation}\label{eqn:Cr}
\mathcal{C}_r:=\llbracket rP+1,(r+1)P\rrbracket.
\end{equation}
Let
\[
I_q:=\llbracket  \alpha^{-1} q-t(\log N)^4, \alpha^{-1} q +t(\log N)^4 \rrbracket.
\]
Cover $I_q$ by consecutive intervals of the form $\mathcal{C}_r$: There exist integers $r_0\le r_1$ such that
\begin{equation}\label{eqn:C_cover}
I_q\subset \mathcal{C}:=\bigcup_{r=r_0}^{r_1}\mathcal{C}_r,
\qquad
M:=r_1-r_0+1\le 100\,S(\log N)^4,
\end{equation}
where we used $P\geq t/(2S)$ and $|I_q| \leq 3 t(\log N)^4$. Finally, set
\begin{equation}\label{eqn:Cbig}
\mathcal{C}_{\mathrm{big}}:=\llbracket N_1,N_2\rrbracket\setminus \mathcal{C}.
\end{equation}
 
 Recall we work under Assumption~\ref{ass:NT_assumption}, and that $\Lambda_i$ denotes the eigenvalue of $\L$ with $\zeta$-localization center $i$ (Definition~\ref{def:quasi_intro}). Define 
\begin{equation}\label{eqn:D_def}
D:=\sum_{i=N_1}^{N_2} F(\Lambda_i)\,G(q-\alpha i-t\ve(\Lambda_i)).
\end{equation}

\noindent Below we will define a Lax matrix $\tilde\L$, with rows and columns indexed by $\llbracket N_1, N_2 \rrbracket$, obtained from $\L$ by resampling a specified subset (either $\mathcal{C}_{\mathrm{big}}$ or one of the $\mathcal{C}_r$) of the variables $(a_i,b_i)$ (and replacing them with $(\tilde a_i, \tilde b_i)$ which are independent and identically distributed, and independent of $(a_i, b_i)_{i \in \llbracket N_1, N_2 \rrbracket}$). Letting $\tilde{\varphi} : \llbracket 1, N \rrbracket \rightarrow \llbracket N_1, N_2 \rrbracket$ denote a $\zeta$-localization center bijection for $\tilde{\mathbf{L}}$, and denoting $\eig \tilde{\mathbf{L}} = (\tilde{\lambda}_1, \tilde{\lambda}_2, \ldots , \tilde{\lambda}_N)$, for each $i \in \llbracket N_1, N_2 \rrbracket$ set $\tilde\Lambda_i = \tilde{\lambda}_{\varphi^{-1} (i)}$ (as in Definition~\ref{def:quasi_intro}). Then define $\tilde D$ by 
\begin{equation}\label{eqn:Dtilde_def}
\tilde D:=\sum_{i=N_1}^{N_2} F(\tilde\Lambda_i)\,G(q-\alpha i-t\ve(\tilde\Lambda_i)).
\end{equation}
Below, we will bound $|D-\tilde D|$.

\medskip

\noindent\textbf{Resampling the chunk $\mathcal{C}_{\mathrm{big}}$.}
We first treat the case in which we resample all indices in $\mathcal{C}_{\mathrm{big}}$. Namely, to obtain $\tilde\L$, for each $i\in\mathcal{C}_{\mathrm{big}}$, we replace the $(i,i)$ and $(i,i+1),(i+1,i)$ entries of $\L$ with independent random variables $\tilde b_i$ and $\tilde a_i$ with the same law as $b_i$ and $a_i$, and which are independent of $(\mathbf a;\mathbf b)$. Restrict to the overwhelmingly probable event $\mathsf{E}_1$ on which $\mathsf{SEP}_{\L} (e^{-(\log N)^2}) \cap \mathsf{SEP}_{\tilde{\L}} (e^{-(\log N)^2})$ holds and 
\[
|\Lambda_i|\le \log N\quad\text{and}\quad |\tilde\Lambda_i|\le \log N
\qquad\text{for all } i\in\llbracket N_1,N_2\rrbracket,
\]
which holds with overwhelming probability by Lemmas \ref{lem:bd_lem} and \ref{lem:sep_lem}. 

Then Lemma \ref{lem:ve_tail} gives $|\ve(\Lambda_i)|\le C\log N$ for each $i \in \llbracket 1, N \rrbracket$. Hence, indices within distance $(\log N)^3$ of $\mathcal{C}_{\mathrm{big}}$ do not contribute to $D-\tilde D$, namely,
\begin{equation}\label{eqn:Cbig_far}
\sum_{d(i,\mathcal{C}_{\mathrm{big}})\le (\log N)^3} F(\Lambda_i)\,G(q-\alpha i-t\ve(\Lambda_i))
-\sum_{d(i,\mathcal{C}_{\mathrm{big}})\le (\log N)^3} F(\tilde\Lambda_i)\,G(q-\alpha i-t\ve(\tilde\Lambda_i))
=0.
\end{equation}
Indeed, if $d(i,\mathcal{C}_{\mathrm{big}})\le (\log N)^3$, then $|i-\alpha^{-1} q| \le t (\log N)^4 / 2$, so 
\begin{align}
|q-\alpha i-t\ve(\Lambda_i)|
&\ge \displaystyle\frac{t}{2} \cdot (\log N)^4 - Ct \log N > 4S, \label{eqn:arg_outside_support}
\end{align}

\noindent and similar reasoning applies upon replacing $\Lambda_i$ with $\tilde{\Lambda}_i$. Hence, 
\begin{equation}\label{eqn:argument_dif}
|q-\alpha i-t\ve(\Lambda_i)|>4S; \qquad |q-\alpha i-t\ve(\tilde\Lambda_i)|>4S.
\end{equation}

\noindent As $\supp G \subseteq [-S, S]$ this gives $G(q-\alpha i-t\ve(\Lambda_i)) = 0 = G(q-\alpha i-t\ve(\tilde{\Lambda}_i))$, verifying \eqref{eqn:Cbig_far}.

Thus, we may restrict the sums defining $D$ and $\tilde{D}$ to $\{i:d(i,\mathcal{C}_{\mathrm{big}})>(\log N)^3\}$. Next, restrict to the overwhelmingly probable event $\mathsf{E}_2$ from Lemma~\ref{lem:Lax_eig_coupling}, and let $\psi$ be the injective map provided by that lemma. Let us also restrict to the overwhelmingly probable event $\mathsf{E}_3$ from Lemma~\ref{lem:bd_num_termsU} with $U = S$. Further abbreviate
\[
G_i:=G(q-\alpha i-t\ve(\Lambda_{i})),\qquad
\tilde G_i:=G(q-\alpha i-t\ve(\tilde\Lambda_i)).
\] 
On $\mathsf{E}_1\cap\mathsf{E}_2\cap\mathsf{E}_3$ we have
\begin{align}\label{eqn:finalbd}
|D-\tilde D|
&=\left|\sum_{d(i,\mathcal{C}_{\mathrm{big}})>(\log N)^3}
(F(\Lambda_i)\,G_i - F(\tilde\Lambda_i)\,\tilde G_i)\right| \notag\\
&\le
\left|\sum_{d(i,\mathcal{C}_{\mathrm{big}})>(\log N)^3}
(F(\Lambda_{\psi(i)})\,G_{\psi(i)} - F(\tilde\Lambda_i) \tilde{G}_i)\right|  + C A B (\log N)^3 \notag\\
&\leq B \sum_{d(i,\mathcal{C}_{\mathrm{big}})>(\log N)^3}
|F(\Lambda_{\psi(i)})-F(\tilde \Lambda_{i})| +A \sum_{d(i,\mathcal{C}_{\mathrm{big}})>(\log N)^3} |\tilde G_i -  G_{\psi(i)}| \notag\\
&\qquad + C A B (\log N)^3 \le C A B (\log N)^7, 
\end{align}
where we have used $|\tilde\Lambda_i-\Lambda_{\psi(i)}|\le e^{-c(\log N)^3}$; the bound \eqref{eqn:F_lip} on $F$; the fact that $|G'(x)| \leq B U^{-1}$; and the fact that $|\psi(i)-i|\le (\log N)^2$ (by Lemma~\ref{lem:Lax_eig_coupling}). More explicitly, for the first sum on the second-to-last line, each summand is at most $B \cdot A e^{-c(\log N)^3}$, while the restriction to $\mathsf{E}_3$ implies that there are at most $C U (\log N)^5$ indices with $G_i\neq 0$ or $\tilde G_i\neq 0$; hence this contribution is at most
\[
C A B U(\log N)^5 e^{-c(\log N)^3}\le C A B.
\]
For the second sum, apply the mean value theorem: For each index $i$ in the summation, for some $|\xi_i|\le 2 (\log N)^2$, we may rewrite
\[
G_{\psi(i)}-\tilde G_i
= \big(-\alpha(\psi(i)-i) - t(\ve(\Lambda_{\psi(i)})-\ve(\tilde \Lambda_{i}) ) \big) \cdot G'(q-\alpha \psi(i)-t\ve(\Lambda_{\psi(i)})+\xi_i),
\]
and therefore
\[
|G_{\psi(i)}-\tilde G_i|
\le \frac{B}{U}\Big(\alpha |\psi(i)-i| + t |\ve(\Lambda_{\psi(i)})-\ve(\tilde \Lambda_i)|\Big).
\]
Since $|\tilde\Lambda_i-\Lambda_{\psi(i)}|\le e^{-c(\log N)^3}$, both eigenvalues lie in $[-\log N,\log N]$ on $\mathsf{E}_1$, and since \eqref{eqn:vebd} gives $\sup_{|x|\le \log N} |\ve'(x)|\le C\log N \log\log N$, we have
\[
t |\ve(\Lambda_{\psi(i)})-\ve(\tilde \Lambda_i)| \le C T(\log N)^2 \log\log N \, e^{-c(\log N)^3}\le 1.
\]
Thus since $|\psi(i)-i|\le (\log N)^2$, each summand is at most $CBU^{-1}(\log N)^2$. Finally, $\supp G'$ has length at most $2U$, and our restriction to $\mathsf{E}_3$ implies that there are at most $CU(\log N)^5$ indices for which the argument of $G'$ lies in $\supp G'$. Hence the second sum is at most
\[
A \cdot C U(\log N)^5 \cdot C B U^{-1}(\log N)^2
\le C A B (\log N)^7.
\]
This yields the final bound in \eqref{eqn:finalbd}.

\medskip

\noindent\textbf{Resampling a chunk $\mathcal{C}_r$.}
Fix $r\in\llbracket r_0,r_1\rrbracket$, and let $\tilde\L$ be obtained by resampling $(a_i,b_i)$ for $i\in\mathcal{C}_r$. Let
\[
\mathcal{C}_r':=\llbracket rP-\lfloor(\log N)^4\rfloor,\ (r+1)P+\lfloor(\log N)^4\rfloor\rrbracket.
\]
Let us again restrict to the overwhelmingly probable event $\mathsf{E}_1$ from above, and to $\mathsf{E}_2$ and $\mathsf{E}_3$ from Lemmas~\ref{lem:Lax_eig_coupling} and \ref{lem:bd_num_termsU}, respectively. Then
\begin{align}\label{eqn:Gsum_expand}
&\left|\sum_{i\in\llbracket N_1,N_2\rrbracket\setminus \mathcal{C}_r'}
\bigl(F(\Lambda_i)\,G_i - F(\tilde\Lambda_i)\,\tilde G_i\bigr)\right| \notag\\
& \le
\left|\sum_{i\in\llbracket N_1,N_2\rrbracket\setminus \mathcal{C}_r'}
\bigl(F(\Lambda_{\psi(i)})\,G_{\psi(i)} - F(\tilde\Lambda_i)\,\tilde G_i\bigr)\right|
+ C A B (\log N)^3  \le C A B (\log N)^7,
\end{align}
using again $|\psi(i)-i|\le(\log N)^2$ and the same mean-value-theorem and counting argument as in the case of $\mathcal{C}_{\mathrm{big}}$.

By \eqref{eqn:Gsum_expand}, it remains to bound the contribution from indices in $\mathcal{C}_r'$. We claim that there is an event $\mathsf{E}_r$ of overwhelming probability such that on $\mathsf{E}_r$,
\begin{equation}\label{eqn:Gsumbd1}
\sum_{i\in \mathcal{C}_r'}
A\Bigl(|G(q-\alpha i-t\ve(\Lambda_i))|+|G(q-\alpha i-t\ve(\tilde\Lambda_i))|\Bigr)
\le A B (\log N)^8.
\end{equation}

By Lemma~\ref{lem:veff_inc}, $|\ve'(\lambda)|\ge c(\log N)^{-1}$ for $\lambda\in[-\log N,\log N]$. Hence, defining
\[
I:=\{\lambda \in [-\log N, \log N]:\ q-\alpha rP-t\ve(\lambda)\in[-2S,2S]\},
\]
\noindent we have 
\[
|I|\le CS t^{-1} \log N.
\]
Applying Lemma~\ref{lem:min} with $R=tS^{-1}$, we may restrict to the overwhelmingly probable event
\[
\mathsf{E}_r:=\{N_{I,\mathcal{C}_r'}\le (\log N)^7\}\cap \{\tilde N_{I,\mathcal{C}_r'}\le (\log N)^7\},
\]
where $\tilde N$ is the count from Lemma~\ref{lem:min} for $\tilde\L$. On $\mathsf{E}_r$, Assumption \ref{ass:G} implies \eqref{eqn:Gsumbd1}.

Now we use the bounds derived above to conclude the proof. We compute the quantity $\sum_{k=1}^m A_k^2$ defined in \eqref{eqn:mcd_S_U} in Lemma \ref{lem:McD_bd} using the computations above and recall that $m = M+1 \leq S (\log N)^5$: We obtain
\begin{equation}
  \sum_{k=1}^m A_k^2 \leq  (A B)^2 S (\log N)^{C} .
\end{equation}
Moreover, for the quantity $p$ there we may take $p = c^{-1} e^{-c (\log N)^2}$ for some absolute constant $c> 0$. By \cite[Lemma 5.4]{Agg25}, we also have $U \le 4ABN$. As a consequence, we may conclude by invoking Lemma \ref{lem:McD_bd}.
\end{proof}

\subsection{Proof of Lemma \ref{lem:no_q_conc_expect}}
\label{app:lemma3.18}

In this section we show Lemma \ref{lem:no_q_conc_expect} (omitting details which are either straightforward computations or which follow from arguments similar to ones already used several times throughout this section).

\begin{proof}[Proof of Lemma \ref{lem:no_q_conc_expect}]

For $\lambda \in \mathbb{R}$, let
\begin{equation}\label{eqn:Fbounded}
    \tilde{F}(\lambda) \coloneqq F(\lambda)\cdot \mathbbm{1}_{|F(\lambda)| \leq A} +A\cdot \mathbbm{1}_{F(\lambda) > A} -A\cdot \mathbbm{1}_{F(\lambda) < -A}   .
\end{equation}
Let $\mathsf{E}_0 = \mathsf{BND}_{\L(0)}(\log N)$, which is overwhelmingly probable by Lemma \ref{lem:bd_lem}, and note that 
$$
\mathbbm{1}_{\mathsf{E}_0}\Bigg| \sum_{i=n_1}^{n_2}F(\Lambda_i) \cdot G( q-\alpha i - t \ve(\Lambda_i)) -\sum_{i=n_1}^{n_2}\tilde{F}(\Lambda_i) \cdot G( q-\alpha i - t \ve(\Lambda_i)) \Bigg| = 0.
$$
Moreover, by an argument similar to the one given for \eqref{expectationH}, 
\begin{multline}\label{eqn:FGbounded_dif}
        \mathbb{E} \left[ \mathbbm{1}_{\mathsf{E}_0^{\complement}} \Bigg| \sum_{i=n_1}^{n_2}F(\Lambda_i) \cdot G( q-\alpha i - t \ve(\Lambda_i)) -\sum_{i=n_1}^{n_2}\tilde{F}(\Lambda_i) \cdot G( q-\alpha i - t \ve(\Lambda_i)) \Bigg| \right]   \\
        \leq N\cdot  e^{-c(\log N)^2}  \cdot C AB.
\end{multline}
Therefore, it suffices to prove \eqref{eqn:expectation_approx} with $F$ replaced by $\tilde{F}$.

Let 
\begin{equation}\label{eqn:H_FGbounded}
    H(\lambda, q) \coloneqq G(-q) \tilde{F}(\lambda).
\end{equation}
Notice, $H$ satisfies the following:
\begin{enumerate}
\item For each fixed $q \in \mathbb{R}$, both $\lambda \mapsto H(\lambda, q)$ and $\lambda \mapsto \partial_q H(\lambda, q)$ satisfy \eqref{eqn:F_lip} with $A$ there given by $A B$ here, and with $A $ there given by $A B / S$ here, respectively. \label{item:H10}
\item For each $q,z \in \mathbb{R}$, 
\begin{equation}\label{eqn:Hlip_bd}
|H(z,q)| \leq A B \qquad |\partial_q H(z,q) | \leq A B S^{-1} \qquad 
|\partial_q^2 H(z,q) | \leq A B S^{-2}.
\end{equation}
Moreover, for each $z \in \mathbb{R}$, 
\begin{equation}\label{eqn:Hsuppbd}
\mathrm{Supp}\left(q \mapsto  H(z,q) \right) \subset [-S, S].
\end{equation}
\label{item:H20}
\end{enumerate}

We will break up $\llbracket n_1, n_2 \rrbracket$ into $m$ chunks $\llbracket n_{1,j}, n_{2, j} \rrbracket$, $j=1,\dots, m$, with $n_{1,j+1} = n_{2,j}+1$ for all $j$, $n_{1,1} = n_1$, and $n_{2,m} = n_2$. Let $n \coloneqq n_{2,j}-n_{1,j}+1$; we will assume (without loss of generality) that $n_2-n_1+1$ is divisible by $n$, and choose $n$ (and thus also $m$) such that $(\log N)^{50} \leq n \leq S$, and we will specify all bounds in terms of $m$ and $n$ under this assumption until the end of the proof. Let $\mathsf{E}_0$ be the event on which the event $\mathsf{BND}_{\L(0)}(\log N)$ holds, and let $\mathsf{E}_1$ be the overwhelmingly probable event on which the conclusion of Lemma \ref{lem:bd_num_termsU} holds (with $U $ there given by $S$ here); let us define $\mathsf{E}$ as in the intersection of these events.

Define $n_{av,j} \coloneqq (n_{2,j}+n_{1,j})/2$. Observe by a Taylor expansion that
\begin{flalign}\label{eqn:first_expansion}
\begin{aligned}
&\left|  \sum_{i=n_1}^{n_{2}}H(\Lambda_i,\alpha i -q + t \ve(\Lambda_i))  - \sum_{j=1}^m \sum_{i=n_{1,j}}^{n_{2,j}} H(\Lambda_i,\alpha n_{av,j} -q + t \ve(\Lambda_i))  \right| \\
&\leq  \sum_{j=1}^m \left| \sum_{i=n_{1,j}}^{n_{2, j}} \left( H(\Lambda_i,\alpha i - q + t \ve(\Lambda_i))  - H(\Lambda_i,\alpha n_{av,j} -q  + t \ve(\Lambda_i)) \right) \right| \\
&\leq  \alpha \sum_{j=1}^m \left| \sum_{i=n_{1,j}}^{n_{2, j}}  \partial_q H(\Lambda_i, \alpha n_{av,j} -q +  t \ve(\Lambda_i))  (n_{av,j}-i) \right|\\
& \qquad + \alpha^2 \sum_{j=1}^m  \sum_{i=n_{1,j}}^{n_{2, j}} \frac{1}{2} \left| \partial_q^2 H(\Lambda_i, \alpha i-q+  t \ve(\Lambda_i)+ \xi_i)  (n_{av,j}-i)^2   \right| 
\end{aligned}
\end{flalign}
where on the last line $|\xi_i| \leq n \leq S$ for each $i \in \llbracket N_1, N_2 \rrbracket$. On $\mathsf{E}$, the expression on the final line can be upper bounded by $ C AB n^2 S^{-1} (\log N)^6$ by the last bound in \eqref{eqn:Hlip_bd}, since there are at most $S (\log N)^6$ terms in the summation (as is quickly deduced by, for example, following the argument in the paragraph containing the display \eqref{eqn:xikidf}). In addition, since $\mathsf{E}$ is overwhelmingly probable and the sum is bounded by $\alpha^2 A B N^2 $ deterministically, the expectation of the summation on the last line of \eqref{eqn:first_expansion} can be upper bounded by $ABn^2 S^{-1} (\log N)^7$. 

For fixed, $j$, we may rewrite the inner sum in the second to last line of \eqref{eqn:first_expansion} as
\begin{multline}\label{eqn:second_expansion}
   \sum_{i=n_{1,j}}^{n_{2, j}} \partial_q H(\Lambda_i,\alpha  n_{av,j}-q+ t \ve(\Lambda_i))\left(-\mathbbm{1}_{i > n_{av,j}} \sum_{n_{av,j} < p \leq i} 1 + \mathbbm{1}_{i < n_{av,j}} \sum_{i \leq p < n_{av,j}} 1 \right) \\
     =  \sum_{p=n_{1, j}}^{n_{2,j}} \Big(-\mathbbm{1}_{p > n_{av,j}} \sum_{i=p}^{n_{2, j}} \partial_q H(\Lambda_i,\alpha n_{av,j} - q+ t \ve(\Lambda_i))  \\
     + \mathbbm{1}_{p < n_{av,j}} \sum_{i=n_{1,j}}^{p} \partial_q H(\Lambda_i,\alpha  n_{av,j} - q+ t \ve(\Lambda_i)) \Big) .
\end{multline}

\noindent If $p < n_{2, j}-(\log N)^{10}$ in the second to last line, or $p > n_{1,j} + (\log N)^{10}$ in the last line, then \cite[Lemma 5.7]{Agg25} with $H$ there given by $H(\lambda,\alpha n_{av,j} - q+ t \ve(\lambda))$ here (and $A$ there given by $ABS^{-1}$ here by Items \ref{item:H10} and \ref{item:H20} above), we obtain 
\begin{multline}\label{eqn:lem57_appl}
    \left|  \mathbb{E}\left[ \sum_{i=p}^{n_{2, j}} \partial_q H(\Lambda_i,\alpha n_{av,j} - q+ t \ve(\Lambda_i)) \right] - (n_{2,j}-p)\int_{-\infty}^{\infty}\partial_q H(\lambda,\alpha n_{av,j} - q+ t \ve(\lambda)) \varrho(\lambda) d \lambda \right|\\
    \leq A B S^{-1} (\log N)^{6} ,
\end{multline}
for $p < n_{2, j} - (\log N)^{10}$; similarly, 
\begin{multline}\label{eqn:lem57_appl2}
    \left|  \mathbb{E}\left[ \sum_{i=n_{1,j}}^{p} \partial_q H(\Lambda_i,\alpha n_{av,j} - q+ t \ve(\Lambda_i))\right] - (p-n_{1,j})\int_{-\infty}^{\infty}\partial_q H(\lambda,\alpha n_{av,j} - q+ t \ve(\lambda)) \varrho(\lambda) d \lambda \right| \\
    \leq  A B S^{-1} (\log N)^{6} ,
\end{multline}

 \noindent for $p > n_{1,j} + (\log N)^{10}$.

Substituting \eqref{eqn:lem57_appl} and \eqref{eqn:lem57_appl2} back into \eqref{eqn:second_expansion} we see (upon taking the expectation) that at leading order, the sums over $i$ cancel in pairs, corresponding to when $p = n_{2,j} - k$ in the first sum and $p = n_{1,j} + k$ in the second. Adding up the errors from the right hand sides of \eqref{eqn:lem57_appl} and \eqref{eqn:lem57_appl2} there, and from the boundary terms with $|p-n_{1,j}| + |p-n_{2,j}| \le 2(\log N)^{10}$, we observe that the absolute value of the expectation of \eqref{eqn:second_expansion} is upper bounded by $A B n S^{-1} (\log N)^{7}$. Recalling our goal of bounding the summation in the second to last line of \eqref{eqn:first_expansion}, we sum this error over $j$, and obtain an overall error upper bounded by $A B n m S^{-1}  (\log N)^{8}$. Therefore, combining our bounds on the expectations of the last two lines of \eqref{eqn:first_expansion},
\begin{multline}\label{eqn:first_expansion2}
\mathbb{E}\left[ \left|  \sum_{i=n_1}^{n_{2}}H(\Lambda_i,\alpha i -q + t \ve(\Lambda_i))  - \sum_{j=1}^m \sum_{i=n_{1,j}}^{n_{2,j}} H(\Lambda_i,\alpha n_{av,j} -q + t \ve(\Lambda_i))  \right| \right]\\
\leq A B n m S^{-1} (\log N)^{8}  + AB n^2 S^{-1} (\log N)^7.
\end{multline}

Moreover, for any fixed $j$, again by \cite[Lemma 5.7]{Agg25}
\begin{multline}\label{eqn:E_bd1}
\left|\mathbb{E} \left[\sum_{i=n_{1,j}}^{n_{2,j}}   H(\Lambda_i, \alpha (n_{av,j}-k)+t \ve(\Lambda_i)) \right] - n \int_{-\infty}^{\infty}H(\lambda, \alpha n_{av,j}-q+t \ve(\lambda))\varrho(\lambda) d\lambda \right| \\
\leq  A B  (\log N)^6.
\end{multline}
Thus,
\begin{multline}\label{eqn:Hexpect_finalbd}
\left|\mathbb{E} \left[ \sum_{j=1}^m \sum_{i=n_{1,j}}^{n_{2,j}}   H(\Lambda_i, \alpha n_{av,j}-q+t \ve(\Lambda_i)) \right] -  \int_{n_1}^{n_2} \int_{-\infty}^{\infty}H(\lambda, \alpha q'-q+t \ve(\lambda)) d q' \varrho(\lambda) d\lambda \right| \\
\leq AB (\log N)^8 ( m + n^2 S^{-1} ) .
\end{multline}
In the last line we have used \eqref{eqn:E_bd1}, and also the fact that for any $\lambda,q \in \mathbb{R}$, we have 
$$\left|\sum_{j=1}^m \int_{n_{1,j}}^{n_{1,j+1}} \left( H(\lambda, \alpha n_{av,j} - q + t \ve(\lambda)) 
 -  H(\lambda, \alpha q'-q+ t \ve(\lambda)) \right)  dq'  \right| \leq AB n^2 S^{-1} (\log N).$$
  To obtain the above inequality, we are Taylor expanding to second order in $q'$, using the second derivative bound on $H$ from \eqref{eqn:Hlip_bd}, and using the fact that the sum is supported on at most $CSn^{-1}$ indices $j$ (as $H(\lambda,Q) = 0$ for $|Q| > S$).
 
Set $n = \floor{T^{1/2 + \delta}}$. To conclude, it remains to justify that we may also assume
\[
n_2-n_1+1 \le C T(\log N)^{12}.
\]
Define
\[
\hat n_1 \coloneqq \max\Bigl(n_1,\floor{q/\alpha-T(\log N)^{12}}\Bigr),\qquad
\hat n_2 \coloneqq \min\Bigl(n_2,\lceil q/\alpha+T(\log N)^{12}\rceil\Bigr).
\]
Then $\hat n_2-\hat n_1+1 \le CT(\log N)^{12}$. We claim that replacing $(n_1,n_2)$ by $(\hat n_1,\hat n_2)$ changes the left side of \eqref{eqn:expectation_approx} (with $F$ replaced by $\tilde{F}$) by at most $CABe^{-c(\log N)^2}$. Indeed, on $\mathsf{E}$, we have $|\Lambda_i|\le \log N$ for all $i$, and therefore \eqref{eqn:vebd} gives $|\ve(\Lambda_i)|\le C\log N$. Hence, if $i\in \llbracket n_1, n_2 \rrbracket \setminus \llbracket \hat n_1,\hat n_2\rrbracket$, then
\[
|q-\alpha i-t\ve(\Lambda_i)|
\ge \alpha |i-\alpha^{-1} q|-Ct\log N
\ge cT(\log N)^{12}-CT(\log N)^{11}>2S,
\]
since $S\le T\log N$ and $t\le T(\log N)^{10}$. By \eqref{eqn:Hsuppbd}, the corresponding summand vanishes. Thus the sum in \eqref{eqn:expectation_approx} is unchanged on $\mathsf{E}$ after replacing $(n_1,n_2)$ by $(\hat n_1,\hat n_2)$. Since $\mathsf{E}$ is overwhelmingly probable and each summand is bounded by $AB$, the contribution of $\mathsf{E}^{\complement}$ to the expectation is at most $CABe^{-c(\log N)^2}$ after adjusting constants.

The same truncation applies to the integral term. For $q'\notin [\hat n_1,\hat n_2]$ and $|\lambda|\le \log N$, the same estimate gives
\[
|q-\alpha q'-t\ve(\lambda)|>2S,
\]
so $H(\lambda,\alpha q'-q+t\ve(\lambda))=0$ by \eqref{eqn:Hsuppbd}. For $|\lambda|>\log N$, we use $|H(\lambda,\cdot)|\le AB$ together with the Gaussian tail bound $\varrho (x) \le C e^{-cx^2}$ from Lemma \ref{lem:varrho_bd} to see that the contribution of this region is again at most $CABe^{-c(\log N)^2}$. 

Therefore we may replace $(n_1,n_2)$ with $(\hat n_1,\hat n_2)$. Upon this replacement, $m$ is replaced by $\hat{m} = (\hat{n}_2-\hat{n}_1+1)/n \le CT^{1/2-\delta}(\log N)^{12}$. The bounds \eqref{eqn:first_expansion2} and \eqref{eqn:Hexpect_finalbd} (with $(\hat{n}_1,\hat{n}_2,\hat{m})$ here replacing $(n_1,n_2,m)$ there, and recalling that $n = \lfloor T^{1/2+\delta} \rfloor$), and minimizing in $\delta$ (by taking $T^{\delta} = S^{1/3} T^{-1/6}$), then complete the proof.
\end{proof}

\subsection{Proof of Lemma \ref{lem:Ndeltabd}}

\label{ProofConcentrationK}

\begin{proof}[Proof of Lemma \ref{lem:Ndeltabd}]
   The proof is to first reduce to the case $\mathfrak{K} = \log N$ by a coupling argument using Lemma \ref{lem:Lax_eig_coupling}, and then invoke \cite[Proposition 4.2]{Agg25} (which is the statement of Lemma \ref{lem:Ndeltabd} at $\mathfrak{K} = \log N$).
   
   Recall that we assume $(N, N_1, N_2)$, $T$, and $\zeta$ are given as in Assumption \ref{ass:NT_assumption}. First, let $\mathfrak{N} = \floor{ e^{\mathfrak{K}}}$, and let $\mathfrak{N}_1 < N_1  < N_2 < \mathfrak{N}_2$ be integers such that $(\mathfrak{N}, \mathfrak{N}_1, \mathfrak{N}_2)$ satisfy the leftmost condition in the display \eqref{eqn:N1N2}. Let $\bm{\mathfrak{L}}(s) = [\mathfrak{L}_{ij}(s)]_{i,j \in \llbracket \mathfrak{N}_1, \mathfrak{N}_2 \rrbracket}$ denote the $\mathfrak{N} \times \mathfrak{N}$ time $s$ Lax matrix, whose entries satisfy $\mathfrak{L}_{i,i+1}(s) = \mathfrak{L}_{i+1,i}(s) = \mathfrak{a}_i(s)$ and $\mathfrak{L}_{ii}(s) = \mathfrak{b}_i(s)$, where for $t \geq 0$, $\{( \mathfrak{a}_i(t), \mathfrak{b}_i(t))\}_{i \in \llbracket \mathfrak{N}_1, \mathfrak{N}_2  \rrbracket}$ are defined from the dynamics \eqref{eqn:Flaschka_ev} with initial data defined as follows: 
   \begin{equation}\label{eqn:mfaibiN1N2}
       \mathfrak{a}_i(0) = a_i(0) \quad \text{for } i \in \llbracket N_1, N_2-1 \rrbracket, \qquad
       \mathfrak{b}_i(0) = b_i(0) \quad \text{for } i \in \llbracket N_1, N_2 \rrbracket.
   \end{equation}
   and otherwise, they are given as in \eqref{eqn:equil} (but with $(\mathfrak{a}_i, \mathfrak{b}_i)$ here playing the role of $(a_i, b_i)$ here)
   $$\{\mathfrak{b}_i(0)\}_{i \in \llbracket \mathfrak{N}_1, \mathfrak{N}_2  \rrbracket \setminus \llbracket N_1, N_2 \rrbracket} \cup \{\mathfrak{a}_i(0)\}_{i \in \llbracket \mathfrak{N}_1, \mathfrak{N}_2 -1 \rrbracket  \setminus \llbracket N_1, N_2 - 1 \rrbracket},$$ which are mutually independent of each other and of 
   $$ \{b_i(0)\}_{i \in \llbracket N_1, N_2 \rrbracket} \cup \{a_i(0) \}_{i \in \llbracket N_1, N_2 - 1 \rrbracket}  .$$
   Note that (upon setting $\mathfrak{a}_{\mathfrak{N}_2}(0) = 0$) $\{(\mathfrak{a}_i(0),\mathfrak{b}_i(0))\}_{i  \in \llbracket \mathfrak{N}_1, \mathfrak{N}_2  \rrbracket}$ has the joint distribution \eqref{eqn:equil} with $N = \mathfrak{N}$, so $\bm{\mathfrak{L}}(0)$ is sampled from thermal equilibrium. 
   
   Let $\{\lambda_{i;\mathfrak{N}}\}_{i=1}^{\mathfrak{N}}$ denote the eigenvalues of $\bm{\mathfrak{L}}(0)$ (and thus also of $\bm{\mathfrak{L}}(t)$ for any $t \geq 0$). In addition, let $\mathfrak{z} = e^{-100(\log \mathfrak{N})^{3/2}}$, and let $\varphi_{s}^{(\mathfrak{N})} : \llbracket 1, \mathfrak{N} \rrbracket \rightarrow \llbracket \mathfrak{N}_1, \mathfrak{N}_2 \rrbracket$ be a $\mathfrak{z}$-localization center bijection for $\bm{\mathfrak{L}}(s)$, and in addition, let 
   \begin{equation}
       \Lambda_{i; \mathfrak{N}}(s) \coloneqq \lambda_{(\varphi_{s}^{(\mathfrak{N})})^{-1}(i);\mathfrak{N}}, \qquad \text{for } i \in \llbracket \mathfrak{N}_1, \mathfrak{N}_2 \rrbracket
   \end{equation}
   as defined in Definition \ref{def:quasi_intro}.

   Now, our plan is to use a coupling argument. We already have defined a coupling between $\L(s)$ and $\bm{\mathfrak{L}}(s)$ via our coupling between $\L(0)$ and $\bm{\mathfrak{L}}(0)$. We claim that, on an event $\mathsf{E}_1$ satisfying 
   \begin{equation}\label{eqn:PE1bd}
       \mathbb{P}(\mathsf{E}_1^{\complement}) \leq \mathfrak{c}_1^{-1} e^{-c_1 \mathfrak{K}^2},
   \end{equation}
   
    \noindent for some constant $c_1>0$, we have the bound 
   \begin{equation}\label{eqn:Lscouple}
    \max_{i \in \llbracket N_1 + T^2, N_2 - T^2 \rrbracket}\bigl(|\mathfrak{a}_i(s) - a_i(s)| + |\mathfrak{b}_i(s) - b_i(s)|\bigr) \leq e^{-T^2/4}.
   \end{equation}
   Indeed, we define $\mathsf{E}_1$ as the event on which 
   \begin{equation}\label{eqn:aibibd}
       \max(|\mathfrak{a}_i(0)|, |\mathfrak{b}_i(0)|) \leq \mathfrak{K}, \qquad \max(|a_i(0)|, |b_i(0)|) \leq \mathfrak{K}
   \end{equation} 
   for all $i \in \mathbb{Z}$ for which either $\mathfrak{a}_i, \mathfrak{b}_i$ (in the first expression) or $a_i, b_i$ (in the second) above is well-defined. By Lemma \ref{lem:bd_lem}, this event has probability at least $1- \mathfrak{c}_1^{-1} e^{-\mathfrak{c}_1 \mathfrak{K}^2}$. Then, on $\mathsf{E}_1$, we may apply \cite[Lemma 4.5]{Agg25a} with the parameters $(\tilde{N}_1, \tilde{N}_2, \tilde{N}; N_1, N_2, N;N_1',N_2', N' ; T; A, K)$ there given by $(\mathfrak{N}_1, \mathfrak{N}_2, \mathfrak{N}; N_1, N_2, N; N_1, N_2, N; T \log N; \mathfrak{K}, T^2)$ here. Note that with this choice of $A, K$, we have $K \geq 200 A T \log N$, since $\mathfrak{K} \leq N^{1/200}= T^{1/2}$, so the assumptions of \cite[Lemma 4.5]{Agg25a} are satisfied on $\mathsf{E}_1$. Therefore, by this lemma, \eqref{eqn:Lscouple} holds on $\mathsf{E}_1$.

In what follows, we will utilize Lemma \ref{lem:Lax_eig_coupling}. Define $\L^{aug}(s) = [L_{ij}^{aug}]_{i, j \in \llbracket \mathfrak{N}_1, \mathfrak{N}_2 \rrbracket} $ as the $\mathfrak{N} \times \mathfrak{N}$ matrix which agrees with $\L(s)=[L_{ij}(s)]_{i, j}$ for $i, j \in \llbracket N_1, N_2 \rrbracket$, and which has all other entries set to $0$. Note that if $i \in \llbracket 1, N \rrbracket$, then $\varphi_s(i) \in \llbracket N_1, N_2 \rrbracket$ is a $\mathfrak{z} = e^{-100(\log \mathfrak{N})^{3/2}}$-localization center for the matrix $\L^{aug}(s)$. Thus, on the event $\mathsf{E}_2$, defined as the intersection of the $(\mathbf{L}, \tilde{\L}) = (\bm{\mathfrak{L}}(s), \L^{aug}(s))$ case of Lemma \ref{lem:Lax_eig_coupling} and of $\mathsf{SEP}_{\L(0)}(e^{-\mathfrak{K}^2})$, there exists a bijection $\psi : \llbracket N_1, N_2 \rrbracket \rightarrow  \llbracket N_1, N_2 \rrbracket $ such that, for $i \in \llbracket N_1 + 2 T^2, N_2 - 2T^2 \rrbracket$,
\begin{equation}\label{eqn:LLi_ssi_dif}
    |\Lambda_i-\Lambda_{\psi(i); \mathfrak{N}}| \leq e^{-c (\log \mathfrak{N})^3}, \qquad \text{ and } |\psi(i) - i| \leq (\log \mathfrak{N})^2.
\end{equation}
To obtain the above, we have used the fact that $(\log \mathfrak{N})^3 \leq T^2$, and the bound \eqref{eqn:Lscouple}. By Lemmas \ref{lem:sep_lem} and \ref{lem:Lax_eig_coupling}, the event $\mathsf{E}_2$ satisfies $ \mathbb{P}(\mathsf{E}_2^{\complement}) \leq c^{-1} e^{-c \mathfrak{K}^2}$.

We will next show that, on an event $\mathsf{E}$ of probability at least $1-\mathfrak{c}^{-1} e^{-\mathfrak{c} \mathfrak{K}^2}$, for all $k \in \llbracket N_1 + T^4, N_2 - T^4 \rrbracket$, the bound
   \begin{equation}\label{eqn:FGconc_Dif}
       \left| \sum_{i=\mathfrak{N}_1}^{\mathfrak{N}_2}F(\Lambda_{i;\mathfrak{N}}(s)) G(q_{k;\mathfrak{N}}(s) -q_{i;\mathfrak{N}}(s)) - \sum_{i=N_1}^{N_2}F(\Lambda_i(s)) G(q_k(s)  -q_i(s)) \right| \leq A B (\log \mathfrak{N})^{9}, 
    \end{equation}
    holds. We exhibit the event $\mathsf{E}$ on which \eqref{eqn:FGconc_Dif} holds by following the argument proving \cite[Equation (4.12)]{Agg25a}, with $\bm{\mathfrak{L}}$ here replacing $\tilde{\L}$ there. For the rest of the proof, we set
\[
K \coloneqq \lceil S (\log \mathfrak{N})^{9/2} \rceil.
\]

Let $\mathsf{E}$ denote the intersection of the following events:

\begin{enumerate}
\item[(i)] $\mathsf{E}_1$ from above, so \eqref{eqn:Lscouple} and \eqref{eqn:aibibd} hold. We shrink $\mathsf{E}_1$, so that in addition $\mathsf{E}_1$ contains the events
\begin{equation}\label{eqn:ailower}
\bigcap_{i \in \llbracket \mathfrak{N}_1, \mathfrak{N}_2 -1 \rrbracket} \{|\mathfrak{a}_i(s)| > e^{-(\log \mathfrak{N})^2}\}
\end{equation}
and
\begin{equation}\label{eqn:LamNupper}
\bigcap_{i \in \llbracket1, \mathfrak{N} \rrbracket} \{ |\lambda_{i;\mathfrak{N}}| \leq \mathfrak{K} \} \bigcap \bigcap_{i \in \llbracket1, N \rrbracket} \{ |\lambda_{i}| \leq \mathfrak{K} \}.
\end{equation}
We also shrink $\mathsf{E}_1$ so that $\bigcap_{t \geq 0}\mathsf{BND}_{\L(t)}(\mathfrak{K})$ holds. We may arrange that \eqref{eqn:PE1bd} still holds by shrinking $c_1$; indeed, the probability bounds for \eqref{eqn:ailower}, \eqref{eqn:LamNupper}, and $\bigcap_{t \geq 0}\mathsf{BND}_{\L(t)}(\mathfrak{K})$ follow from Lemma \ref{lem:bd_lem}, using that $\bm{\mathfrak{L}}(0)$ is sampled from thermal equilibrium, and then we recall the definition of $\mathfrak{K}$ (and that $\mathfrak{K} \geq \log N$).
\item[(ii)] $\mathsf{E}_2$ from the paragraph containing \eqref{eqn:LLi_ssi_dif}.
\item[(iii)] The event $\mathsf{E}_3$ on which
\begin{equation}
\bigl|(q_{k;\mathfrak{N}}(s)-q_{i;\mathfrak{N}}(s))-\alpha(k-i)\bigr|
\le |k-i|^{1/2} (\log \mathfrak{N})^2,  \quad \text{for } i,k \in \llbracket \mathfrak{N}_1+ T (\log \mathfrak{N})^4, \mathfrak{N}_2 -T (\log \mathfrak{N})^4 \rrbracket \label{eqn:E5b}
\end{equation}
and on which  
\begin{equation}\label{eqn:E5a2}
\sgn(k-i) (q_{k;\mathfrak{N}}(s)-q_{i;\mathfrak{N}}(s))
\ge \frac{\alpha}{2} |k-i| ,  \qquad \text{for } i,k \in \llbracket \mathfrak{N}_1 , \mathfrak{N}_2 \rrbracket \text{ with } |i-k| > T (\log \mathfrak{N})^6
\end{equation}
and 
\begin{equation}\label{eqn:E5a1}
 \bigl| (q_j(0)-q_i(0))-\alpha(j-i)\bigr|
\le  |j-i|^{1/2} (\log \mathfrak{N})^2 ,  \qquad \text{for } i,j \in \llbracket N_1 , N_2 \rrbracket.
\end{equation}

By Lemma \ref{lem:q_spacing} with $N$ there given by $\mathfrak{N}$ here, \eqref{eqn:E5b} and \eqref{eqn:E5a2} have probability at least $1-\mathfrak{c}_3^{-1} e^{-\mathfrak{c}_3 \mathfrak{K}^2}$. Moreover, the probability of \eqref{eqn:E5a1} has the same lower bound (after possibly decreasing $\mathfrak{c}_3$), by \cite[Lemma 3.10]{Agg25}, which is a Chernoff bound for random walks.
\end{enumerate}

By the probability bounds already recorded for $\mathsf{E}_1$ and $\mathsf{E}_2$, and the one above for $\mathsf{E}_3$, we have $\mathbb{P}(\mathsf{E}^\complement)\le c^{-1} e^{-c\mathfrak{K}^2}$. We restrict to $\mathsf{E}$ below. 

\medskip

\noindent\textbf{Step 1: truncation to $|i-k|\le K$.}
If $i \in \llbracket N_1 + T^2, N_2 - T^2 \rrbracket$, then on $\mathsf{E}$, by \eqref{eqn:E5b}, \eqref{eqn:E5a2}, \eqref{eqn:ailower}, and \eqref{eqn:Lscouple}, we have
\begin{equation}\label{eqn:E5a3}
\sgn(k-i) (q_{k}(s)-q_{i}(s))
\ge \frac{\alpha}{4} |k-i| \geq 2 S ,  \qquad \text{for }  |i-k| > K.
\end{equation}
If instead $i \notin \llbracket N_1 + T^2, N_2 - T^2 \rrbracket$, then since $k \in \llbracket N_1 + T^4, N_2 - T^4 \rrbracket$ we have $|i-k| > T^4/2 \geq T (\log \mathfrak{N})^6/2$. Moreover, by our restriction to $\bigcap_{t \geq 0}\mathsf{BND}_{\L(t)}(\mathfrak{K})$ and the bound $s \le T \log N$, we have $|q_i(s)-q_i(0)|, |q_k(s)-q_k(0)| \le T \log N \cdot \log \mathfrak{N}$, and hence \eqref{eqn:E5a1} gives
\[
|q_k(s)-q_i(s)| \ge |q_k(0)-q_i(0)| - 2 T (\log \mathfrak{N})^2 > T (\log \mathfrak{N})^3 > 2S.
\]
Therefore, if $|i-k|\ge K\ge S(\log \mathfrak{N})^{9/2}$, then in either case \eqref{eqn:E5a3} holds, and similarly, if $|i-k|\ge K\ge S(\log \mathfrak{N})^{9/2}$, then by \eqref{eqn:E5b} and \eqref{eqn:E5a2}, $|q_{k;\mathfrak{N}}(s)-q_{i;\mathfrak{N}}(s)|>2S$.
Therefore, since $\supp G \subseteq [-S,S]$, both $G(q_k(s)-q_i(s))$ and $G(q_{k;\mathfrak{N}}(s)-q_{i;\mathfrak{N}}(s))$ vanish whenever $|i-k|\ge K$. Hence
\begin{align}
&\Bigl|\sum_{i=\mathfrak{N}_1}^{\mathfrak{N}_2} F(\Lambda_{i;\mathfrak{N}}(s))\,G(q_{k;\mathfrak{N}}(s)-q_{i;\mathfrak{N}}(s))
-\sum_{i=N_1}^{N_2} F(\Lambda_i(s))\,G(q_k(s)-q_i(s))\Bigr| \nonumber\\
&\qquad =
\left|\sum_{|i-k|\le K}
\Bigl(F(\Lambda_{i;\mathfrak{N}}(s))\,G(q_{k;\mathfrak{N}}(s)-q_{i;\mathfrak{N}}(s))
- F(\Lambda_i(s))\,G(q_k(s)-q_i(s))\Bigr)\right|. \label{eqn:splitwindow}
\end{align}

\medskip

\noindent\textbf{Step 2: split into an $F$-difference and a $G$-difference term.}
Using $xy-x'y'=(x-x')y+x'(y-y')$, we may bound \eqref{eqn:splitwindow} by
\begin{equation}
\left|\sum_{|i-k|\le K}
\bigl(F(\Lambda_{i;\mathfrak{N}}(s))-F(\Lambda_i(s))\bigr)\,G(\Delta_i)\right|
+ A \sum_{|i-k|\le K}|G(\Delta_{i}^{(\mathfrak{N})})-G(\Delta_i)|. \label{eqn:FGsplit}
\end{equation}
where we abbreviate
\[
\Delta_{i}^{(\mathfrak{N})}\coloneqq q_{k;\mathfrak{N}}(s)-q_{i;\mathfrak{N}}(s),
\qquad
\Delta_i\coloneqq q_k(s)-q_i(s).
\]

\medskip

\noindent\textbf{Step 3: control the $F$-difference term.}
Set
\[
\mathcal{S}_F\coloneqq \sum_{|i-k|\le K}
\bigl(F(\Lambda_{i;\mathfrak{N}}(s))-F(\Lambda_i(s))\bigr)\,G(\Delta_i).
\]
On $\mathsf{E}$, by \eqref{eqn:LamNupper} all eigenvalues $\lambda_{i;\mathfrak{N}}$ of $\bm{\mathfrak{L}}(s)$ and all eigenvalues $\lambda_i$ of $\L(s)$ lie in $[-\mathfrak{K},\mathfrak{K}]$, and \eqref{eqn:LLi_ssi_dif} gives
$|\Lambda_i(s)-\Lambda_{\psi(i);\mathfrak{N}}(s)|\le e^{-c(\log \mathfrak{N})^3}
\le e^{-\mathfrak{K}^{5/2}}$ for $\mathfrak{K}$ large.
Thus, by the second assumption on $F$,
\[
|F(\Lambda_i(s))-F(\Lambda_{\psi(i);\mathfrak{N}}(s))|\le A e^{-\mathfrak{K}^2}.
\]
Since $\psi$ is a bijection and $|\psi(i)-i|\le (\log\mathfrak{N})^2$, the symmetric difference between
$\psi(\{i:|i-k|\le K\})$ and $\{i:|i-k|\le K\}$ has size at most $2(\log\mathfrak{N})^2$, and hence
\[
\left|\sum_{|i-k|\le K} F(\Lambda_{i;\mathfrak{N}}(s))\,G(\Delta_i)
- \sum_{|i-k|\le K} F(\Lambda_{\psi(i);\mathfrak{N}}(s))\,G(\Delta_{\psi(i)})\right|
\le 4AB(\log\mathfrak{N})^2 \le AB(\log\mathfrak{N})^3.
\]
Therefore,
\begin{equation}\label{eqn:FDIF}
\begin{aligned}
|\mathcal{S}_F|
&\le \left|\sum_{|i-k|\le K}
\Bigl(F(\Lambda_{\psi(i);\mathfrak{N}}(s))\,G(\Delta_{\psi(i)})
- F(\Lambda_i(s))\,G(\Delta_i)\Bigr)\right|
+ AB(\log\mathfrak{N})^3 \\
&\le 2B\sum_{|i-k|\le K}|F(\Lambda_{\psi(i);\mathfrak{N}}(s))-F(\Lambda_i(s))| \\
&\qquad + A\sum_{|i-k|\le K}|G(\Delta_{\psi(i)})-G(\Delta_i)|
+ AB(\log\mathfrak{N})^3 \\
&\le 2(2K+1)AB e^{-\mathfrak{K}^2}  + A\sum_{|i-k|\le K}|G(\Delta_{\psi(i)})-G(\Delta_i)|
+ AB(\log\mathfrak{N})^3 \\
&\le A\sum_{|i-k|\le K}|G(\Delta_{\psi(i)})-G(\Delta_i)|
+ 2AB(\log\mathfrak{N})^3.
\end{aligned}
\end{equation}
Above, we have used \eqref{eqn:LamNupper} to bound $|F(\Lambda_i(s))| \leq A$.

\medskip

\noindent\textbf{Step 4: control the $G$-difference terms.}
By our assumptions on $G$,
\begin{equation}\label{eqn:Gsummand}
\begin{aligned}
|G(\Delta_{i}^{(\mathfrak{N})})-G(\Delta_i)|
\le & B S^{-1}|\Delta_{i}^{(\mathfrak{N})}-\Delta_i| \\
&+ B \cdot \mathbbm{1}_{\{\Delta_{i}^{(\mathfrak{N})}\ge 0\ge \Delta_i\}}
+ B \cdot \mathbbm{1}_{\{\Delta_i\ge 0\ge \Delta_{i}^{(\mathfrak{N})}\}}.
\end{aligned}
\end{equation}
First, on $\mathsf{E}$ we have
\begin{equation}\label{eqn:Delta_i_dif}
|\Delta_{i}^{(\mathfrak{N})}-\Delta_i|
=|(q_{k;\mathfrak{N}}(s)-q_k(s))-(q_{i;\mathfrak{N}}(s)-q_i(s))|.
\end{equation}
Using $q_\ell-q_{\ell+1}=2\log a_\ell$, the bounds \eqref{eqn:Lscouple} and \eqref{eqn:ailower}, and the definition of $\mathfrak{K}$ in terms of $\mathfrak{N}$, we obtain
\[
|\Delta_{i}^{(\mathfrak{N})}-\Delta_i|
\le e^{-c(\log \mathfrak{N})^3}.
\]
Hence
\[
\sum_{|i-k|\le K} B S^{-1}|\Delta_{i}^{(\mathfrak{N})}-\Delta_i|
\le (2K+1)\,B S^{-1}\,e^{-c(\log\mathfrak{N})^3}\;\le\;B.
\]

Second, the indicator terms in \eqref{eqn:Gsummand} can only be nonzero if $\Delta_i$ and $\Delta_{i}^{(\mathfrak{N})}$
have opposite signs.
But by \eqref{eqn:E5b} and \eqref{eqn:E5a1}, for $|i-k|\ge (\log\mathfrak{N})^5$ both
$\Delta_i$ and $\Delta_{i}^{(\mathfrak{N})}$ are nonzero and have the same sign as $\alpha(k-i)$
(since both are within $ |k-i|^{1/2}(\log\mathfrak{N})^2 \leq  (\log \mathfrak{N})^{9/2}$ of $\alpha(k-i)$, which has size at least $\alpha (\log N)^5$). Therefore, these indicator functions vanish unless
$|i-k|\le (\log\mathfrak{N})^5$, which yields at most $2(\log\mathfrak{N})^5+1$ indices.
Consequently,
\[
\sum_{|i-k|\le K}
\bigl(\mathbbm{1}_{\{\Delta_{i}^{(\mathfrak{N})}\ge 0\ge \Delta_i\}}
+\mathbbm{1}_{\{\Delta_i\ge 0\ge \Delta_{i}^{(\mathfrak{N})}\}}\bigr)
\;\le\; 5(\log\mathfrak{N})^5.
\]
It follows that
\[
\sum_{|i-k|\le K}|G(\Delta_{i}^{(\mathfrak{N})})-G(\Delta_i)|
\le B + 5B(\log\mathfrak{N})^5
\le 6 B(\log\mathfrak{N})^{5}.
\]

For the first term in the last line of \eqref{eqn:FDIF}, we similarly have
\[
|G(\Delta_{\psi(i)})-G(\Delta_i)|
\le B S^{-1}|\Delta_{\psi(i)}-\Delta_i|
+ B \cdot \mathbbm{1}_{\{\Delta_{\psi(i)}\ge 0\ge \Delta_i\}}
+ B \cdot \mathbbm{1}_{\{\Delta_i\ge 0\ge \Delta_{\psi(i)}\}}.
\]
By \eqref{eqn:E5a1}, we obtain
\[
\begin{aligned}
|\Delta_{\psi(i)}-\Delta_i|
&= |q_i(s)-q_{\psi(i)}(s)| \\
&\le
\alpha |\psi(i)-i|
+ |\psi(i)-i|^{1/2}(\log \mathfrak{N})^2 \\
&\le C(\log\mathfrak{N})^3.
\end{aligned}
\]
Moreover, since $K \coloneqq \lceil S (\log \mathfrak{N})^{9/2} \rceil$,
\[
\sum_{|i-k|\le K} B S^{-1}|\Delta_{\psi(i)}-\Delta_i|
\le C B S^{-1}S(\log\mathfrak{N})^5 (\log\mathfrak{N})^3
\le C B (\log\mathfrak{N})^{8}.
\]
Also, the indicator terms can be nonzero only if $\Delta_i$ and $\Delta_{\psi(i)}$ have opposite signs. The bound on $|\Delta_{\psi(i)}-\Delta_i|$ together with \eqref{eqn:E5a1} shows that this can happen for at most $C(\log\mathfrak{N})^5$ indices $i$. Hence
\[
\sum_{|i-k|\le K}|G(\Delta_{\psi(i)})-G(\Delta_i)|
\le C B (\log\mathfrak{N})^{8}.
\]

\noindent Combining the result of steps 3--4 in \eqref{eqn:FGsplit} gives
\[
\left| \sum_{i=\mathfrak{N}_1}^{\mathfrak{N}_2}F(\Lambda_{i;\mathfrak{N}}(s)) G(q_{k;\mathfrak{N}}(s) -q_{i;\mathfrak{N}}(s)) - \sum_{i=N_1}^{N_2}F(\Lambda_i(s)) G(q_k(s)  -q_i(s)) \right| \leq A B (\log \mathfrak{N})^{9}
\]
on the overwhelmingly probable event $\mathsf{E}$. Lemma \ref{lem:Ndeltabd} then follows from combining this bound with \cite[Proposition 4.2]{Agg25} with $(N, \L, t)$ there given by $(\mathfrak{N}, \bm{\mathfrak{L}}, s)$ here.
\end{proof}

\subsection{Proof of Lemma \ref{lem:xi_bd_outline}}\label{app:xi_bd_outline}

We will show the two statements \eqref{eqn:small_Lambda_bd_outline} and \eqref{eqn:mathfrakX_bd_outline} of Lemma \ref{lem:xi_bd_outline} separately.

\begin{proof}[Proof of \eqref{eqn:small_Lambda_bd_outline}]
This follows immediately from Item \ref{item:conc3} of Lemma \ref{lem:no_q_conc}, except for the uniformity in $t $. The overwhelmingly probable extension to all $t \in [0, T \log N]$ follows from a union bound over a lattice mesh, together with a continuity estimate. We omit further details on this argument, as it is entirely analogous to the proofs of Lemmas \ref{lem:Xi_brownian_cont_allparam_pl} and \ref{lem:gen_brownian_cont_allparam_prelim}.  
\end{proof}

\begin{proof}[Proof of \eqref{eqn:mathfrakX_bd_outline}]
Set $A_2:=\log N$ and fix an integer $A_1>(\log N)^5$. For $\tau\in[0,\log N]$ and $\mathfrak{q}\in\mathbb{R}$ define the event
\[
\mathsf{E}_{\mathfrak{q},\tau}^{(A_1)}
:=\left\{\sup_{|\Lambda|\le A_1}\big|\mathfrak{X}\big(\Lambda, (\mathfrak{q}-\tau\ve(\Lambda)) / \alpha,\tau\big)\big|
\le A_1A_2^2\right\}.
\]
By Lemma~\ref{lem:Xbd1} (applied with this choice of $A_2$ and $A_1$), for every $\tau\in[0,\log N]$ and $\mathfrak{q}\in\mathbb{R}$,
\begin{equation}\label{eqn:PEqtau}
\mathbb{P}\bigl((\mathsf{E}_{\mathfrak{q},\tau}^{(A_1)})^{\complement}\bigr)
\le c_3^{-1}e^{-c_3A_1^2}.
\end{equation}

\noindent Set
\begin{equation}\label{eqnAdf}
A := C_0 A_1 \log N,
\end{equation}
where $C_0>0$ is a constant to be chosen sufficiently large. For
$|\mathfrak{q}|\le (\log N)^4$, $\tau\in[0,\log N]$, and $|\Lambda|\le A_1$ we have
\begin{align*}
\left|\frac{\mathfrak{q}-\tau \ve(\Lambda)}{\alpha}\right|
\le \alpha^{-1} (|\mathfrak{q}|+\tau\,|\ve(\Lambda)|)
& \le C\Bigl((\log N)^4+(\log N)\sup_{|\lambda|\le A_1}|\ve(\lambda)|\Bigr) \\
&\le C\big((\log N)^4 + A_1\log N\big)
\le C_0 A_1\log N = A,
\end{align*}
where we used the growth bound $|\ve(\lambda)|\le C(|\lambda|+1)$ for $|\lambda|\le A_1$
(and $A_1\ge (\log N)^5$) from Lemma \ref{lem:ve_tail}. Consequently, whenever
$|\mathfrak{q}|\le (\log N)^4$, $\tau\in[0,\log N]$, and $|\Lambda|\le A_1$,
the point $(\Lambda, (\mathfrak{q}-\tau \ve(\Lambda)) / \alpha,\tau)$ lies in the domain $K_0$ of Lemma \ref{lem:holder_norm} with this choice of $A$. Define the event $\mathsf{F}^{(A_1)}$ as the outcome of Lemma \ref{lem:holder_norm} with this choice of $A$ and with $u:=A^5$; this event  satisfies
\[
\mathbb{P}\bigl((\mathsf{F}^{(A_1)})^{\complement}\bigr)\le c^{-1}e^{-cA_1^2}.
\]

Let $\delta:=A^{-500}$ with $A$ as in \eqref{eqnAdf}, and let $\mathcal{M}$ be a $\delta$-mesh covering of $[-(\log N)^4, (\log N)^4] \times [-\log N, \log N]$. Thus $|\mathcal{M}|\le A_1^{C}\,(\log N)^C$. 
 Define
\[
\mathsf{E}^{(A_1)}:=\bigcap_{(\mathfrak{q}_0,\tau_0)\in\mathcal{M}}\mathsf{E}_{\mathfrak{q}_0,\tau_0}^{(A_1)}.
\]
By a union bound and \eqref{eqn:PEqtau}, 
\begin{equation}\label{eqn:PEbdXbdpf}
\mathbb{P}\bigl((\mathsf{E}^{(A_1)})^{\complement}\bigr)\le c^{-1}e^{-cA_1^2}.
\end{equation}

For any $(\mathfrak{q},\tau)\in\llbracket -(\log N)^4,(\log N)^4\rrbracket\times[0,\log N]$
denote as $(\mathfrak{q}_0,\tau_0)\in\mathcal{M}$ a point satisfying $\|(\mathfrak{q},\tau)-(\mathfrak{q}_0,\tau_0)\|_2\le C\delta$.
On $\mathsf{E}^{(A_1)}\cap\mathsf{F}^{(A_1)}$, for all $|\Lambda|\le A_1$,
\begin{align*}
\big|\mathfrak{X}\big(\Lambda, (\mathfrak{q}-\tau\ve(\Lambda)) / \alpha,\tau\big)\big|
&\le \big|\mathfrak{X}\big(\Lambda, (\mathfrak{q}-\tau\ve(\Lambda)) / \alpha,\tau\big)
-\mathfrak{X}\big(\Lambda, (\mathfrak{q}_0-\tau_0\ve(\Lambda)) / \alpha,\tau_0\big)\big|\\
& \qquad + 
\big|\mathfrak{X}\big(\Lambda, (\mathfrak{q}_0-\tau_0\ve(\Lambda)) / \alpha,\tau_0 \big)\big| \le A_1A_2^2 + C  u\,\delta^{2/5}A_1^{2/5}
\le 2A_1A_2^2,
\end{align*}
where in the last step we used $u=A^5$ and $\delta=A^{-500}$, so $u \delta^{2/5} A_1^{2/5}\le  A_1A_2^2$. Therefore, on $\mathsf{E}^{(A_1)}\cap\mathsf{F}^{(A_1)}$,
\begin{equation}\label{eqn:XbdfixedA1}
\sup_{|\Lambda|\le A_1}\ \sup_{|\mathfrak{q}|\le (\log N)^4}\ \sup_{\tau\in[0,\log N]}
\big|\mathfrak{X}\big(\Lambda, (\mathfrak{q}-\tau\ve(\Lambda)) / \alpha,\tau\big)\big|
\le 2 A_1A_2^2.
\end{equation}

Finally, take the intersection of $\mathsf{E}^{(A_1)}\cap\mathsf{F}^{(A_1)}$ over dyadic integers $A_1=2^\ell$ such that $A_1\ge(\log N)^5$. Since
\[
\sum_{\ell:\ 2^\ell\ge(\log N)^5} e^{-c\,2^{2\ell}}
\le e^{-c(\log N)^{10}/2}
\le e^{-c(\log N)^2},
\]
this intersection holds with probability at least $1-c^{-1}e^{-c(\log N)^2}$.

On this event, given any $\Lambda\in\mathbb{R}$ with $|\Lambda|$ larger than $2 (\log N)^5$, choose $A_1=2^\ell$ so that $A_1/2\le |\Lambda|+1\le A_1$. Then \eqref{eqn:XbdfixedA1} yields
\[
\sup_{|\mathfrak{q}|\le (\log N)^4}\sup_{\tau\in[0,\log N]}
\big|\mathfrak{X}\big(\Lambda, (\mathfrak{q}-\tau\ve(\Lambda)) / \alpha,\tau\big)\big|
\le 2A_1A_2^2
\le C(\log N)^{C}(|\Lambda|+1)^2,
\]
after using the facts that $A_1\le 2(|\Lambda|+1)$ and $A_2 = \log N$. This yields \eqref{eqn:mathfrakX_bd_outline} when $|\Lambda| \ge 2 (\log N)^5$. Since this bound \eqref{eqn:mathfrakX_bd_outline} also holds if $|\Lambda| \leq 2 (\log N)^5$ by Lemma \ref{lem:Xbd1}, we deduce \eqref{eqn:mathfrakX_bd_outline} in general.
\end{proof}

\section{Variance of linear statistic}\label{app:lin_var}
Below we will deal with a random collection of numbers $\mathbf{a} = (a_1,\dots, a_{N-1})$ and $\mathbf{b} = (b_1,\dots, b_N)$, where $N>0$ is an integer. In all cases, will denote by $\L =\L(\mathbf{a}, \mathbf{b})$ the symmetric tridiagonal $N \times N$ matrix with $L_{i,i+1} = L_{i+1,i} = a_i$ for $i \in \llbracket 1, N-1 \rrbracket$, and $L_{i i} = b_i$, for $i = \llbracket 1, N \rrbracket$ (and $L_{ij} = 0$ for all $(i,j)$ not of the above form). In all that follows, will also let $V(x) $ denote a polynomial potential function, which we always assume has an even degree $d \geq 2$, and positive leading coefficient.
Throughout this appendix, $\theta>0$ will be a real number.

\subsection{Density and continuity estimates}
\label{app:lin_var_apriori}

In this section we discuss probability densities that minimize a certain variational problem, and recall certain facts about them (from the prior works \cite{GM22,GM23}) or prove several for them. 

Let $\mathcal{F}^{(2)}$ denote the asymptotic free energy of the \emph{high temperature beta ensemble}
\begin{equation}\label{eqn:HTmeas}
d \mu_{N, HT}^{\theta, V}(\mathbf{a}, \mathbf{b}) = Z_N^{HT}(\theta, V)^{-1} \cdot e^{-\Tr V(L)} \prod_{j=1}^{N-1} \mathbbm{1}_{a_j>0} \cdot a_j^{2 \theta (1-j/N) - 1}   d\mathbf{a} d\mathbf{b},
\end{equation}

\noindent where $Z_N^{HT}(\theta, V)$ is a normalization constant so that $\mu_{N, HT}^{\theta, V}$ is a probability measure. The associated free energy is defined as
\begin{equation}
\label{F2}
\mathcal{F}^{(2)}(\theta, V) = -\lim_{N \rightarrow \infty} N^{-1} \log Z_N^{HT}(\theta, V).
\end{equation}

In the following discussion, we will match the notation of \cite{GM22} and recall several results from there. Let $\L = \L(\mathbf{a}, \mathbf{b})$ with $(\mathbf{a}, \mathbf{b})$ sampled from \eqref{eqn:HTmeas}, and let $\eig \mathbf{L} = (\lambda_1, \lambda_2, \ldots , \lambda_N)$. The \emph{equilibrium measure} for this ensemble is the function $\rho_{\theta, V}^{(2)} : \mathbb{R} \rightarrow \mathbb{R}_{\ge 0}$, defined to satisfy 
\begin{equation}\label{eqn:HTequil}
\lim_{N \rightarrow \infty} \frac{1}{N} \sum_{i=1}^N g(\lambda_i) = \int_{-\infty}^{\infty}g(x) \rho_{\theta, V}^{(2)}(x) dx,
\end{equation}
almost surely for all bounded continuous functions $g$. By \cite[Lemma 3.2]{GM22}, the equilibrium measure is the unique probability density function minimizing the free energy functional 
\begin{equation}\label{eqn:functional}
\rho \mapsto f_{\theta}^V (\rho) = \int_{-\infty}^{\infty}V(x) \rho(x) dx - \theta \int_{-\infty}^{\infty}\int_{-\infty}^{\infty}\log|x-y| \rho(x) \rho(y) d x d y + \int_{-\infty}^{\infty}\rho(x) \log \rho(x) dx.
\end{equation}

\noindent Also by \cite[Lemma 3.2]{GM22}, there exists a real number $\lambda_{\theta}^V \in \mathbb{R}$ such that $\rho_{\theta, V} = \rho_{\theta, V}^{(2)}$ satisfies
\begin{equation}\label{eqn:EL}
  V(x)  -2  \theta \int_{-\infty}^{\infty}\log|x-y| \rho_{\theta, V} (y)  d y + \log \rho_{\theta, V}(x) = \lambda_{\theta}^V.
\end{equation}
In fact, $\lambda_{\theta}^V$ is determined by the relation
\begin{equation}\label{eqn:lambda_theta_V}
 \lambda_{\theta}^V =  f_{\theta}^V (\rho_{\theta,V}) - \theta \int_{-\infty}^{\infty} \int_{-\infty}^{\infty}\log |x-y| \rho_{\theta, V}(x) \rho_{\theta, V}(y) dx dy,
\end{equation}

 \noindent by integrating \eqref{eqn:EL} against $\rho_{\theta, V}(x)$ and comparing the resulting expression with \eqref{eqn:functional}.

Throughout this section, we denote $\langle \cdot, \cdot \rangle_{2}$ as the usual $L^2(\mathbb{R})$ inner product. Recalling $\mathbf{T}$ from \eqref{operatort}, for two densities $\rho, \tilde \rho$, define the distance (which is well-defined by \cite[Equation (30)]{GM22})
\begin{flalign}\label{eqn:D}
\begin{aligned} 
D(\rho, \tilde \rho) & \coloneqq \left( - \int_{-\infty}^{\infty}\int_{-\infty}^{\infty}\log |x- y| \cdot (\rho(x)-\tilde{\rho}(x)) (\rho(y)-\tilde{\rho}(y) ) dx d y \right)^{1/2} \\
& = \bigg(-\frac{1}{2} \cdot \langle \rho-\tilde \rho, \T (\rho - \tilde \rho) \rangle_{2} \bigg)^{1/2} .
\end{aligned} 
\end{flalign}

 Let us record an intermediate lemma due to \cite{GM23}.\footnote{While the uniformity of the constant $\mathfrak{c}$ in Lemma \ref{lem:rho_tail_bd} over the family $(V_t)_{|t| \le \varepsilon}$ and $\theta \in (\delta,\theta_0)$ was not explicitly stated in \cite[Lemma 2.2]{GM23}, it quickly follows from the proof of \cite[Lemma 2.2]{GM23} and \cite[Lemma 3.2]{GM22}.}
  \begin{lem}[{\cite[Lemma 2.2]{GM23}}]\label{lem:rho_tail_bd}
Let $V_t = V_0 + t g$ be a one-parameter family of nonconstant polynomials of even degree. Let $I_{\varepsilon}$ denote either the interval $[-\varepsilon, \varepsilon]$ or the interval $[0, \varepsilon]$. Assume for some real numbers $C_0 \in \mathbb{R}$ and $a, \varepsilon > 0$, for all $x \in \mathbb{R}$ we have  
 \begin{equation}\label{eqn:V_growth}
 \displaystyle\inf_{t \in I_{\varepsilon}} V_t(x) >  a x^2 + C_0.
 \end{equation}
 Then for any real numbers $\theta_0 > \delta>0$, there exists a constant $\mathfrak{c} > 0$, depending only on $(a, C_0, \theta_0, \varepsilon,\delta)$, such that
\begin{equation}\label{eqn:rho_tail_bd}
     \displaystyle\sup_{t \in I_{\varepsilon}} \displaystyle\sup_{\theta \in (\delta, \theta_0)} \rho_{\theta,V_t}^{(2)} (x) \leq \mathfrak{c}^{-1} e^{-\mathfrak{c} x^2}.
 \end{equation}
 \end{lem}

This next lemma provides a continuity estimate for $\rho_{\theta,V}^{(2)}$ in the potential $V$.
\begin{lem}\label{lem:rho_apriori}
 Let $V$ be a polynomial of degree $2 k $ with positive leading coefficient, and let $g$ be any polynomial such that $V + \varepsilon g$ has positive leading coefficient for $\varepsilon > 0$ sufficiently small. Let $\rho = \rho_{\theta, V  }^{(2)}$, $\tilde V = V + \varepsilon g$, $\tilde{\rho} = \rho_{\theta, \tilde V }^{(2)}$, and $\Delta \rho = \rho - \tilde \rho$. There exists a constant $\mathfrak{C} > 0$ dependent on $(V, g, \theta)$, such that for all $\varepsilon>0$ sufficiently small and all $x \in \mathbb{R}$ we have
 \begin{equation}\label{eqn:rho_apriori_Tbd_stmt}
 |\T(\Delta \rho)(x)| \leq \mathfrak{C} \varepsilon |\log \varepsilon|^{\mathfrak{C}},
 \end{equation}
 \begin{equation}\label{eqn:rho_apriori_lambda_stmt}
 |\lambda_\theta^V - \lambda_\theta^{\tilde V}| \leq \mathfrak{C} \varepsilon |\log \varepsilon|^{\mathfrak{C}},
 \end{equation}
 \begin{equation}\label{eqn:rho_apriori_pw_stmt}
 |\rho(x) - \tilde{\rho}(x)| \leq 
 \rho(x) \cdot \big|1- \exp\big(\lambda_\theta^{\tilde V} - \lambda_\theta^{V} - \varepsilon g(x)  - \theta \T(\Delta \rho)(x) \big) \big|,
 \end{equation}
 and
 \begin{equation}\label{eqn:rho_ptwise_apriori}
 \sup_{x \in \mathbb{R}} |\rho(x) - \tilde \rho(x)| \leq \mathfrak{C} \varepsilon |\log \varepsilon|^{\mathfrak{C}}.
 \end{equation}
 \end{lem}

\begin{proof}

We first derive an initial pointwise bound on $D(\rho, \tilde \rho) $ given by \eqref{eqn:D}. For this, we compute analogously to the end of the proof of \cite[Lemma 3.2]{GM22}. Denoting $\Delta \rho \coloneqq (\rho- \tilde{\rho})$, we obtain 
\begin{align}
0 & \geq f_{\theta}^V(\rho)- f_{\theta}^V(\tilde \rho) \notag \\
&= \int_{-\infty}^{\infty}V(x) \Delta \rho(x) \, dx - \theta \int_{-\infty}^{\infty}\int_{-\infty}^{\infty}\log|x-y| \rho(x) \rho(y) \, dx \, d y \notag \\
&\qquad + \theta \int_{-\infty}^{\infty}\int_{-\infty}^{\infty}\log|x-y| \tilde{\rho}(x) \tilde{\rho}(y) \, dx \, d y 
+ \int_{-\infty}^{\infty}\bigl(\rho(x) \log \rho(x) - \tilde \rho(x) \log \tilde \rho(x)\bigr)  \, d x   \notag\\
&= \int_{-\infty}^{\infty}V(x) \Delta \rho(x) \, dx - \theta \langle \T \tilde \rho, \Delta \rho \rangle_2 - \frac{\theta}{2} \langle \T \Delta \rho,  \Delta \rho \rangle_2
+ \int_{-\infty}^{\infty}\bigl(\rho(x) \log \rho(x) - \tilde \rho(x) \log \tilde \rho(x)\bigr)  \, d x \notag\\
&= \int_{-\infty}^{\infty}\bigl(V(x)-\tilde V(x)\bigr) \Delta \rho(x) \, dx  - \frac{\theta}{2} \langle \T \Delta \rho,  \Delta \rho \rangle_2
+ \int_{-\infty}^{\infty}\rho(x) \log \frac{\rho(x)}{\tilde \rho(x)} \, d x. \label{eqn:big_display}
\end{align}
To obtain the final line, we replaced $V  = \tilde V + V - \tilde V$ and also integrated $ \Delta \rho$ against \eqref{eqn:EL} (with $V$ there replaced by $\tilde{V}$). Now, $\int_{-\infty}^{\infty}\rho(x) \log (\rho(x)/\tilde \rho(x)) d x \geq 0$, and $- \theta\langle \mathbf{T} \Delta \rho,  \Delta \rho \rangle_2 / 2 = \theta D(\rho, \tilde{\rho})^2 $. So we obtain 
\begin{equation}\label{eqn:Dbd}
 D(\rho, \tilde{\rho})^2 \leq \frac{1}{\theta} \int_{-\infty}^{\infty}\bigl(\tilde{V}(x)- V(x)\bigr) \Delta \rho(x) \, dx.
\end{equation}

Denote by $\varphi_{\varepsilon}$ a bump function supported on $[- |\log \varepsilon|, |\log \varepsilon|]$ and equal to $1$ on $[1- |\log \varepsilon|, |\log \varepsilon|-1]$. Plugging in $\tilde V - V = \varepsilon g$, and using the Plancherel theorem, we may write 
\begin{multline}\label{eqn:split}
\frac{1}{\theta} \int_{-\infty}^{\infty}\bigl(\tilde V(x)-V(x)\bigr) \Delta \rho(x) \, dx = \frac{1}{\theta} \int_{-\infty}^{\infty}\varphi_{\varepsilon}(x) \bigl(\tilde V(x)-V(x)\bigr) \Delta \rho(x) \, dx \\
+ \frac{1}{\theta} \int_{-\infty}^{\infty}\bigl(1- \varphi_{\varepsilon}(x)\bigr) \bigl(\tilde V(x)-V(x)\bigr) \Delta \rho(x) \, dx  \\
= \varepsilon \frac{1}{\theta} \int_{-\infty}^{\infty}\widehat{\varphi_{\varepsilon}  g}(t) \sqrt{|t|} \frac{1}{\sqrt{|t|}} \widehat{\Delta \rho}(t) \, d t  + \varepsilon \frac{1}{\theta} \int_{-\infty}^{\infty}\bigl(1- \varphi_{\varepsilon}(x)\bigr) g(x) \Delta \rho(x) \, dx.
\end{multline}
We bound the integral in the first term. Using Cauchy-Schwarz, and the fact that $|\widehat{\varphi_{\varepsilon}  g}(t) | \leq C |\log \varepsilon|^C(1+|t|)^{-4}$, we may bound 
\begin{equation}\label{eqn:first_int}
\left| \int_{-\infty}^{\infty}\widehat{\varphi_{\varepsilon}  g}(t) \sqrt{|t|} \frac{1}{\sqrt{|t|}} \widehat{\Delta \rho}(t) d t  \right| \leq 
C  |\log \varepsilon|^C \left( \int_{-\infty}^{\infty}|t|^{-1} |\widehat{\Delta \rho}(t)|^2 \, d t \right)^{1/2} \leq C |\log \varepsilon|^C D(\rho, \tilde{\rho}).
\end{equation}
To obtain the last inequality above, we have used the fact (which is a quick consequence of the Plancherel theorem; see, for example, \cite[Equation (30)]{GM22}) that $( \int_{-\infty}^{\infty}|t|^{-1} |\widehat{\Delta \rho}(t)|^2 \, d t )^{1/2} $ is the Fourier representation of the logarithmic energy from \eqref{eqn:D}, up to a positive constant.

Using Lemma \ref{lem:rho_tail_bd}, for some $c > 0$, both $|\rho(x)| \leq c^{-1} e^{-c |x|^2}$ and $|\tilde{\rho}(x)| \leq c^{-1}e^{-c |x|^2}$, so the second integral in the right hand side of \eqref{eqn:split} can be upper bounded in absolute value as
\begin{equation}\label{eqn:second_int}
\left| \int_{-\infty}^{\infty}\bigl(1- \varphi_{\varepsilon}(x)\bigr) g(x) \Delta \rho(x) \, dx \right| \leq c^{-1} e^{-c |\log \varepsilon|^2},
\end{equation}
or $c > 0$ small enough. We claim that for some $C = C(g)$ and $m = m(g)$,
\begin{equation}\label{eqn:Dbd2}
 D(\rho, \tilde{\rho}) \leq \varepsilon \frac{1}{\theta} C |\log \varepsilon|^C .
\end{equation}
We have two cases: Either $D(\rho, \tilde{\rho}) \leq c^{-1} e^{-c |\log \varepsilon|^2}$, in which case \eqref{eqn:Dbd2} holds, or alternatively, we must have, by \eqref{eqn:Dbd}, \eqref{eqn:split}, and \eqref{eqn:first_int} (bounding the second integral in the right hand side of \eqref{eqn:split} by $D(\rho, \tilde{\rho})$ using \eqref{eqn:second_int}),
\begin{equation}
 D(\rho, \tilde{\rho})^2 \leq \varepsilon \frac{1}{\theta} C \left( |\log \varepsilon|^C + 1\right) D(\rho, \tilde{\rho}).
\end{equation}
\noindent In this case, \eqref{eqn:Dbd2} follows. It then follows from \eqref{eqn:split}, \eqref{eqn:first_int}, \eqref{eqn:second_int}, \eqref{eqn:Dbd2}, and \eqref{eqn:big_display}, we obtain the bound 
\begin{equation}\label{eqn:rel_entbd}
\int_{-\infty}^{\infty}\rho(x) \log \frac{\rho(x)}{\tilde \rho(x)} \, d x \leq C \frac{1}{\theta} \varepsilon^2 |\log \varepsilon|^C.
\end{equation}

Next we will use the pointwise relation \eqref{eqn:EL}, which is equivalent to
\begin{equation}\label{eqn:rho_ptwise}
 \rho(x)  = e^{\lambda_\theta^V - V(x) + \theta \T(\rho)(x)} ,
\end{equation}
to deduce a bound on $|\rho(x) - \tilde{\rho}(x)|$. First, we must bound $|\lambda_\theta^V - \lambda_\theta^{\tilde V}| $. We note by \eqref{eqn:lambda_theta_V},
\begin{equation}\label{eqn:lambda_dif}
|\lambda_\theta^V - \lambda_\theta^{\tilde V}|  \leq | f_{\theta}^V(\rho)  -f_{\theta}^{\tilde V}(\tilde{\rho})| +  \left| \theta \int_{-\infty}^{\infty}\int_{-\infty}^{\infty}\log |x-y| \bigl(\rho(x) \rho(y) - \tilde \rho(x) \tilde \rho(y)\bigr) \, dx \, dy \right|.
\end{equation}
Then, by \eqref{eqn:big_display} and the bounds \eqref{eqn:rel_entbd}, \eqref{eqn:Dbd2}, \eqref{eqn:first_int}, and \eqref{eqn:split} above, and by Lemma \ref{lem:rho_tail_bd},
\begin{equation}\label{eqn:fbd}
|f_{\theta}^V(\rho)  -f_{\theta}^{\tilde V}(\tilde{\rho})| \leq |f_{\theta}^V(\rho)  -f_{\theta}^{V}(\tilde{\rho})|  + |f_{\theta}^V(\tilde{\rho})  -f_{\theta}^{\tilde V}(\tilde{\rho})| 
\leq C \varepsilon \frac{1}{\theta} |\log \varepsilon|^C.
\end{equation}
Moreover, the second difference in \eqref{eqn:lambda_dif} can be upper bounded as follows
  \begin{multline}
  \label{integralrhoestimate0}
 \left| \theta \int_{-\infty}^{\infty}\int_{-\infty}^{\infty}\log |x-y| \bigl((\rho(x) - \tilde \rho(x)) \rho(y) + \tilde \rho(x) (\rho(y)-  \tilde \rho(y))\bigr) \, dx \, dy \right| \\
  \leq  \frac{\theta}{2} \left( \left| \langle \T \Delta \rho, \rho \rangle_{L^2} \right| + \left| \langle \T \Delta \rho, \tilde \rho \rangle_{L^2} \right| \right),
\end{multline}
so it suffices to bound $\T (\Delta \rho)$. For this, we note that, for small $t_0 > 0$ and large $T>0$, we have
\begin{flalign}\label{eqn:ptwise}
\begin{aligned}
\frac{1}{2} \left| \T (\Delta \rho)(x) \right| &=\left | \int_{-\infty}^{\infty}\log|x-y| \Delta \rho(y) d y \right| \\
&= \left| \int_{-\infty}^{\infty}\widehat{\log|x-\cdot |}(t) \widehat{\Delta \rho}(t) d t \right| \\
&\leq \pi  \int_{-\infty}^{\infty}\frac{1}{|t|} | \widehat{\Delta \rho}(t)| d t  \\
&\leq  2\pi\left( \int_0^{t_0} \frac{1}{|t|} | \widehat{\Delta \rho}(t)| d t + \left( \int_{t_0}^{T} \frac{1}{t}  dt \right)^{1/2} \left( \int_{t_0}^T t^{-1} |\widehat{\Delta \rho}(t)|^2 \, d t  \right)^{1/2} + \int_{T}^{\infty} \frac{1}{|t|} | \widehat{\Delta \rho}(t)| d t \right) \\
&\leq C t_0 + (\log T + |\log t_0|)^{1/2} D(\rho, \tilde{\rho}) + \frac{C}{T^{1/2}} \|\Delta \rho\|_{L^2} .
\end{aligned}
\end{flalign}
Here we have used that $\widehat{\Delta \rho}(0)=0$ and, by Lemma \ref{lem:rho_tail_bd}, the first moment of $\Delta \rho$ is finite, so $|\widehat{\Delta \rho}(t)| \leq C |t|$; we also used the Fourier-side representation of $D(\rho,\tilde\rho)$ from \eqref{eqn:D} (as explained in the comment after \eqref{eqn:first_int}). To bound the last term we used Cauchy-Schwarz and the Plancherel theorem. If we choose $t_0 = \varepsilon$, $T = \varepsilon^{-2}$, and bound $\|\Delta \rho\|_{L^2} \leq C$ using \eqref{eqn:rho_tail_bd}, we obtain for each $x \in \mathbb{R}$,
 \begin{equation}
\bigl| \T (\Delta \rho)(x) \bigr| = \left| \int_{-\infty}^{\infty}\log|x-y| \Delta \rho(y) \, d y \right|  \leq C \varepsilon |\log \varepsilon|^C.
\end{equation}
This proves the bound \eqref{eqn:rho_apriori_Tbd_stmt}.

Thus, using \eqref{eqn:lambda_dif}, \eqref{eqn:fbd},  \eqref{integralrhoestimate0}, \eqref{eqn:rho_apriori_Tbd_stmt}, and Lemma \ref{lem:rho_tail_bd}, we obtain a bound of the form 
\begin{equation}
|\lambda_\theta^V - \lambda_\theta^{\tilde V}|  \leq C \varepsilon |\log \varepsilon|^C.
\end{equation}
This proves the bound \eqref{eqn:rho_apriori_lambda_stmt}.

Now, finally, we return to \eqref{eqn:rho_ptwise}; by substituting \eqref{eqn:rho_apriori_lambda_stmt} and the pointwise bound \eqref{eqn:rho_apriori_Tbd_stmt} we obtain
\begin{multline}\label{eqn:finalrhotilderhobd}
|\rho(x) - \tilde{\rho}(x)| \leq 
\rho(x) \cdot \big|1- \exp\big(\lambda_\theta^{\tilde V} - \lambda_\theta^{V} - \varepsilon g(x)  - \theta \T(\Delta \rho)(x) \big) \big| \\
\leq   \rho(x) \cdot \left|1- \exp\left(O (\varepsilon |\log \varepsilon|^C ) - \varepsilon g(x)  \right) \right|  
\leq C \varepsilon |\log \varepsilon|^C.
\end{multline}
\noindent The first inequality is \eqref{eqn:rho_apriori_pw_stmt}. It remains to justify the final inequality. If $|x| \leq |\log \varepsilon|$, we indeed obtain the bound above, since $\varepsilon g(x) = O(\varepsilon |\log \varepsilon|^C)$. On the other hand, if $|x| > |\log \varepsilon|$, then Lemma \ref{lem:rho_tail_bd} implies that
 $$\left|\rho(x) - \tilde{\rho}(x) \right| \leq \left|\rho(x) \right| + \left|\tilde{\rho}(x) \right|  \leq  c^{-1} e^{-c |\log \varepsilon|^2}  \leq C \varepsilon ,$$
 so we obtain \eqref{eqn:finalrhotilderhobd} in this case as well. This proves the bound \eqref{eqn:rho_ptwise_apriori}, and hence completes the proof of the lemma.
\end{proof}

\subsection{Covariances through free energies} 
\label{app:lin_var_res}

In this section we recall results from \cite{GM22,MM24} that express covariances between charges of the Toda lattice in terms of the free energies given by \eqref{F2} above and \eqref{1f} below. 

To that end, consider the \emph{generalized Gibbs ensemble} for the Toda Lattice associated with $V$, which is the probability measure on $(\mathbf{a}, \mathbf{b})$ defined by  
\begin{equation}\label{eqn:Toda_meas}
d \mu_{N, Toda}^{\theta, V}(\mathbf{a}, \mathbf{b}) = Z_N^{Toda}(\theta, V)^{-1} \cdot e^{-\Tr V(\L)} \prod_{j=1}^{N-1}  \mathbbm{1}_{a_j>0}  \cdot a_j^{2 \theta -1}d\mathbf{a} d\mathbf{b},
\end{equation}

\noindent where $Z_N^{Toda}(\theta, V)$ is a normalization constant so that $\mu_{N, Toda}^{\theta, V}(\mathbf{a}, \mathbf{b})$ is a probability measure. The associated free energy is defined as 
\begin{equation}
\label{1f}
\mathcal{F}^{(1)}(\theta, V) = -\lim_{N \rightarrow \infty} N^{-1} \log Z_N^{Toda}(\theta, V).
\end{equation}

 \noindent Let $\rho_{\theta,V}^{(1)} $ denote the associated equilibrium measure, defined to satisfy \eqref{eqn:HTequil} where, sampling $\L$ under \eqref{eqn:Toda_meas}, $(\lambda_1, \lambda_2, \ldots , \lambda_N)$ are the eigenvalues of $\L$.

The following is a corollary of \cite[Theorems 5 and 9]{MM24}, specialized to the setting of the Toda lattice. In fact, the final parts of the statement are not explicitly stated there, but are corollaries of the proofs (as we outline below). 

\begin{lem}[{\cite{GM22,MM24}}] \label{thm:mm}

 Let $V(x)$ be a polynomial of even degree with positive leading coefficient. Then, the following hold.
\begin{enumerate}
\item For any (real) polynomial $g$, we have
\begin{equation}\label{eqn:diffF_id}
\int_{-\infty}^{\infty} g(x) d\rho^{(1)}_{\theta, V }(x) =\i \partial_{t} \mathcal{F}^{(1)}(\theta, V + \i t g(x) )\big|_{t= 0 } ,
\end{equation}
and 
\begin{equation}\label{eqn:difftheta_id}
\int_{-\infty}^{\infty} g(x) \rho^{(1)}_{\theta, V }(x)  dx= \partial_{\theta}\left( \theta  \int_{-\infty}^{\infty} g(x) \rho^{(2)}_{\theta, V }(x)  dx \right).
\end{equation}
\label{item:D2}

\item Suppose $(\mathbf{a}, \mathbf{b})$ are sampled from $d \mu_{N, Toda}^{\theta, V}$, and let $(\lambda_1, \lambda_2, \ldots , \lambda_N)$ denote the eigenvalues of the corresponding tridiagonal matrix $\L = \L(\mathbf{a}, \mathbf{b})$. For any integers $m, n >0$, the limit 
\begin{equation}
\mathcal{C}_{m, n} \coloneqq \lim_{N \rightarrow \infty} \frac{1}{N} \Cov \left( \sum_{j=1}^N \lambda_j^m, \sum_{j=1}^N \lambda_j^n   \right) 
\end{equation}
 exists, and we have
\begin{equation}
\mathcal{C}_{m, n}  = \partial_{t_1} \partial_{t_2} \mathcal{F}^{(1)}(\theta, V + \i t_1 x^m+\i t_2 x^n )\big|_{t_1=t_2=0} .
\end{equation}
\label{item:D3}

\item Let g be a polynomial such that $\deg g < \deg V$. Then there exists a constant $\varepsilon > 0$ such that the function $t \mapsto \mathcal{F}^{(1)}(\theta, V + t g  ) $ is holomorphic for $|t| \leq \epsilon$.
\label{item:D4}

\item For any integer $n > 0$, the limit 
\begin{equation}\label{eqn:C0n}
\mathcal{C}_{0, n} = \mathcal{C}_{n, 0} \coloneqq \lim_{N \rightarrow \infty} \frac{1}{N} \Cov \left( \sum_{j=1}^N \lambda_j^n, -2 \sum_{j=1}^{N-1} \log a_j  \right)
\end{equation}
exists. We have 
\begin{equation}\label{eqn:C0n_formula}
\mathcal{C}_{n, 0} = -\i \partial_{\theta} \partial_{t} \mathcal{F}^{(1)}(\theta, V + \i t x^n )|_{t= 0 }   .
\end{equation} 
Moreover, 
\begin{equation}\label{eqn:C00_formula}
\mathcal{C}_{0, 0} = \var( - 2 \log a_j ) =   \lim_{N \rightarrow \infty} \frac{1}{N} \Cov \left( -2 \sum_{j=1}^{N-1} \log a_j, -2 \sum_{j=1}^{N-1} \log a_j  \right) = \partial_{\theta}^2 \mathcal{F}^{(1)}(\theta, V ).
\end{equation}

\label{item:D5}

\item Let g be a polynomial such that $V + \delta g $ has positive leading coefficient for $\delta > 0$. Then there exists a constant $\varepsilon > 0$ such that the quantity
$
\partial_{t_1} \partial_{t_2} \mathcal{F}^{(1)}(\theta, V + \delta g + \i t_1 x^m+\i t_2 x^n )|_{t_1=t_2=0}
$
is continuous in $\delta$ for all $\delta \in [0, \varepsilon]$.  \label{item:D6}
\end{enumerate}
\end{lem}

\begin{proof}[Proof of Lemma \ref{thm:mm} (Outline)]

\textit{Item \ref{item:D2}.} This holds by \cite[Lemma 3.6]{GM22}, and \cite[Theorem 1.5]{MM24}.

\textit{Item \ref{item:D3}.} This is \cite[Theorem 9]{MM24}.

\textit{Item \ref{item:D4}.} This follows from the fact that, for $|t|$ sufficiently small, the prelimit $-N^{-1} \log Z_N^{Toda}(\theta, V + t g)$ is analytic in $t$, and the convergence to the limit is uniform in $t$ (as follows from the proof of \cite[Theorem 6.7]{MM24}).

\textit{Item \ref{item:D5}.} This can be obtained in the same way as \cite[Theorem 1.9]{MM24}, by a slight modification of the proof. To explain this, in the definition of $J_1^{(N)}$ in \cite[Theorem 3.1]{MM24}, we substitute $U(\mathbf{a}, \mathbf{b}) =  \sum_{i=1}^{N-1}2 \log a_i + \tilde{U}(\mathbf{a}, \mathbf{b})$, where $\tilde{U}(\mathbf{a},\mathbf{b})$ is a function of the form already considered there (note $\mathbf{a}$ and $\mathbf{b}$ are interchanged there compared to here; in this discussion, we adhere to our notation). Let us briefly outline why the proof of this theorem goes through with such a choice of $U$. Following the terminology of \cite[Definition 1.1]{MM24}, the seed for $\sum_{i=1}^{N-1}2 \log a_i$ can be taken to be $2 \log a_i$, and the weed can be taken  to be $0$. Note that in this case $|u| e^{-h}$ is no longer bounded, as assumed in \cite[Lemma 6.6]{MM24}. However, the only place this is used is in the proof of \cite[Lemma B.1(2,3)]{MM24}), when bounding the derivatives $\partial_\alpha k_{\alpha, t}$ and $\partial_t^a k_{\alpha, t}$ of the kernel $ k_{\alpha, t}$ from \cite[Equation (6.1)]{MM24}. Since multiplying by $\log a_i$ leads to a function that remains integrable at $0$ (see \cite[Remark 1.7]{MM24} there), the desired $L^2$ bounds on these derivatives still hold (and uniformly in $t$). Therefore, the same proof goes through with the above choice of $U$.

 Consequently, with the above choice of $U$, we may take derivatives of $J_1^{(N)}$ in $t$ in order to extract the variance of the corresponding statistic, which in this case is taken to be 
\begin{equation}
   \sum_{i=1}^N g(\lambda_i) - \sum_{i=1}^{N-1} 2 \log a_i
\end{equation}
for some polynomial $g$. This gives \eqref{eqn:C00_formula}, which implies (by polarization) \eqref{eqn:C0n_formula}, except the derivatives in $\theta$ are taken in the imaginary direction. This is remedied by observing that both $\mathcal{F}^{(1)}(\theta, V + \i t x^n )|_{t= 0 }$ and $\partial_t \mathcal{F}^{(1)}(\theta, V + \i t x^n )|_{t= 0 }  $ are analytic in $\theta$ for $\theta$ in a complex neighborhood of $\theta_0$ (as they are uniform limits of analytic functions), which yields Item \ref{item:D5} of the lemma.

\textit{Item \ref{item:D6}.} For this, it suffices to know that the leading eigenvalue $\tilde{\lambda}(\alpha,t,\delta) $ of the transfer operator $\mathcal{L}_{\alpha,t,\delta}$ utilized in the proof of \cite[Theorem 6.7]{MM24} (which we now take to also depend on $\delta$), and its derivatives up to order $2$ in $t$, are continuous in $\delta$. Indeed, fixing any polynomial $p$, by the asymptotic expansion from \cite[Theorem 6.7]{MM24} and the proof of \cite[Lemma 5.2]{MM24}, for  $|t|$ small enough one has
\[
-\partial_t^2 \mathcal{F}^{(1)}(\theta, V + \delta g + \i t p ) =  \frac{1}{k} \partial_t^2 \log \tilde{\lambda}(\alpha,t,\delta),
\]
where $k$ is the integer associated with the underlying $k$-circular observable used in \cite{MM24}, so that $N = kM+\ell$ in \cite[Theorem 6.7]{MM24}.
Hence, continuity of $\lambda(\alpha,t,\delta)$ and of its first two $t$-derivatives in $\delta$ implies continuity in $\delta$ of $\partial_t^2\mathcal{F}^{(1)}(\theta, V + \delta g + \i t p )|_{t=0}$. Since the polynomial $p$ is arbitrary, the desired continuity of
\[
\partial_{t_1} \partial_{t_2} \mathcal{F}^{(1)}(\theta, V + \delta g + \i t_1 x^m+\i t_2 x^n )|_{t_1=t_2=0}
\]
then follows by polarization. By \cite[Proposition 2.3]{Gou15}, the desired continuity in $\delta$ of the second derivative of $\tilde{\lambda}$ follows from the fact that the corresponding transfer operator $\mathcal{L}_{\alpha,t,\delta}$ has continuous partial derivatives $\partial_t^a \partial_\delta^b \mathcal{L}_{\alpha,t,\delta}$ up to order two in both $\delta$ and $t$, i.e. for $a,b \in \{0,1,2 \}$, for all $(\delta, t) \in [0,\epsilon) \times (-\epsilon, \epsilon)$ (here we assume $\epsilon>0$ is small enough, and at $\delta = 0$ we refer to a one-sided derivative). This differentiability can in turn be derived in a manner similar to the arguments in the proof of \cite[Theorem 6.7]{MM24}. Indeed, if $h$ is a seed for $\Tr V(\L)$, $\tilde{h}$ is a (compatible) seed for $ \Tr g(\L)$, and $u$ is a (compatible) seed for $\Tr p(\L)$, then $|u|^a |\tilde{h}|^b e^{-h - \delta \tilde{h}}$ is bounded above by a constant for $a, b \in \{0,1,2\}$ and uniformly for $\delta \in [0, \epsilon]$. 
By the transfer kernel's formula in \cite[Definition 6.5]{MM24}, together with the Hilbert--Schmidt estimates used in \cite[Lemma 6.6]{MM24}, this ensures that the derivatives $\partial_t^a \partial_\delta^b \mathcal{L}_{\alpha,t,\delta}$ exist and depend continuously on $(\delta,t)$.
\end{proof}

\subsection{Proof of Lemma \ref{lem:lim_var_formula}}
\label{app:lin_var_pf}

Before proving Lemma \ref{lem:lim_var_formula}, we have the following lemma, which contains the essential computations of this section. In the process of proving it, we will show that $\rho_{\theta, V}^{(1)}$ is given by an explicit formula. The arguments below closely follow the discussion in \cite[Section 4]{Spo20}, using the bounds provided by Lemma \ref{lem:rho_apriori} to justify them. For any function $g : \mathbb{R} \rightarrow \mathbb{R}$, we set 
\begin{flalign*} 
\mu_g = \int_{-\infty}^{\infty} g(\lambda)\rho_{\theta, V}^{(1)}(\lambda) d\lambda,
\end{flalign*} 

\noindent whenever the right hand side converges. In addition, we denote by 
\begin{equation}\label{eqn:rho1inner} 
\langle f, g \rangle_{\rho_{\theta, V}^{(1)}} \coloneqq \int_{-\infty}^{\infty} f(\lambda) \overline{g(\lambda)}\rho_{\theta, V}^{(1)}(\lambda) d\lambda, 
\end{equation} 
the inner product on single variable functions associated to $\rho_{\theta, V}^{(1)}$.

 \begin{lem}\label{lem:spo}

 For any integer $n \ge 0$ and real number $\beta>0$, there exists a constant $\mathfrak{c}>0$ such that the following holds for any real number $\delta \in (0, \mathfrak{c})$. Let $V(x) = \beta x^2/2 + \delta x^{ 2n}$ and $g$ be a polynomial with $\deg g < n = \deg V$. We have
\begin{equation}\label{eqn:var_thm}
- \partial_{t}^2 \mathcal{F}^{(1)}(\theta, V +  t g )\big|_{t=0} =   \langle (1-\theta \mathbf{T}  \boldsymbol{\rho_{\theta, V}^{(2)}})^{-1} (g - \mu_{g}\varsigma_0) , (1-\theta \mathbf{T} \boldsymbol{\rho_{\theta, V}^{(2)}})^{-1} (g - \mu_{g}\varsigma_0) \rangle_{\rho_{\theta, V}^{(1)}};
\end{equation}
\begin{equation}\label{eqn:cov_0nthm}
-\partial_{\theta} \i \partial_{t} \mathcal{F}^{(1)}(\theta, V + \i t g ) \big|_{t= 0 }  = -\alpha(\theta, V) \langle (1-\theta \mathbf{T} \boldsymbol{\rho_{\theta, V}^{(2)}})^{-1} \varsigma_0 , (1-\theta \mathbf{T} \boldsymbol{\rho_{\theta, V}^{(2)}})^{-1} (g - \mu_{g}\varsigma_0) \rangle_{\rho_{\theta, V}^{(1)}};
\end{equation}
and 
\begin{equation}\label{eqn:cov_00thm}
\partial_{\theta}^2 \mathcal{F}^{(1)}(\theta, V ) \big|_{t= 0 }  = \alpha(\theta,V)^2 \langle (1-\theta \mathbf{T} \boldsymbol{\rho_{\theta, V}^{(2)}})^{-1} \varsigma_0 , (1-\theta \mathbf{T} \boldsymbol{\rho_{\theta, V}^{(2)}})^{-1} \varsigma_0\rangle_{\rho_{\theta, V}^{(1)}} .
\end{equation}
\end{lem}

 \begin{remark}
For \eqref{eqn:var_thm}, \eqref{eqn:cov_0nthm}, or \eqref{eqn:cov_00thm} to make sense, one must verify that the operator $(1 -\theta \T  \boldsymbol{\rho_{\theta, V}^{(2)}} )^{-1} $ exists and is bounded on a Hilbert space containing the one associated with the inner product \eqref{eqn:rho1inner}. When the potential $V(x)$ is equal to $V_0(x) = \beta x^2/2$, the existence and boundedness of $(1 -\theta \T  \boldsymbol{\rho_{\theta, V}^{(2)}} )^{-1}$ holds on $\mathcal{H}$ (from Definition \ref{def:inner_prod}). For our purposes, we only need the above results when $V(x) = V_{\delta}(x) = V_0(x) + \delta x^{2 n}$, for a fixed integer $n>0$ and a small real number $\delta>0$ (though it can likely be proven in greater generality with further effort). In this setting, it will suffice to work on the Hilbert space $\mathcal{H}$. Indeed, let us denote by $\mathcal{H}_{\delta}$ the Hilbert space associated to the inner product \eqref{eqn:rho1inner} when $V = V_{\delta}$. By Lemma \ref{lem:rho1_formula} below, $\mathcal{H} \subseteq \mathcal{H}_{\delta}$. That $(1 -\theta \T  \boldsymbol{\rho_{\theta, V}^{(2)}} )^{-1} $ exists and is bounded on $\mathcal{H}$ holds by Lemma \ref{lem:dress_existence} below. 
\end{remark}

\begin{lem}\label{lem:dress_existence}
Let $V(x)$ be a polynomial satisfying the hypotheses of Lemma \ref{lem:spo}. Then, the operator $ (1 - \theta \T \boldsymbol{\rho_{\theta, V}^{(2)}}) : \mathcal{H} \rightarrow \mathcal{H}$ is invertible with a bounded inverse.
\end{lem}

\begin{proof}

Let $\varrho_{\beta} = \rho_{\theta, V_0}^{(2)}$ (as in Definition \ref{def:alpha_Laxdos}), and let $\Delta \rho^{(2)} =   \rho_{\theta, V_\delta}^{(2)}-\varrho_{\beta} $. First, we observe that the second inequality in \eqref{eqn:finalrhotilderhobd} (with $g(x) =  x^{ 2n}$ and $\epsilon = \delta$) from the proof of Lemma \ref{lem:rho_apriori} implies 
\begin{equation}\label{eqn:deltarho2bd}
| \Delta \rho^{(2)}(x)| \leq \delta |\log \delta|^C (x^{2 n}+1) \varrho_{\beta}(x) \le C \delta |\log \delta|^C  (x^{2 n}+1) \varrho (x),
\end{equation}

\noindent where in the last statement we used the first statement of \eqref{rho2} with the third part of Lemma \ref{lem:ve_tail}. Denoting the operator norm on $\mathcal{H}$ by $\| \cdot \|_{\text{op},\mathcal{H}}$, this implies 
\begin{equation}\label{eqn:TDeltaRhoOp}
\| \T \boldsymbol{ \Delta \rho^{(2)} } \|_{\text{op},\mathcal{H}} \leq C \delta |\log \delta|^C,
\end{equation}

\noindent We then claim that, for $\delta>0$ small enough, we have (as operators on $\mathcal{H}$) that
\begin{align}
1 -\theta \T  \boldsymbol{\rho_{\theta, V_{\delta}}^{(2)}} &=  (1 -\theta \T \boldsymbol{\varrho_{\beta}}  )\big( 1 -\theta (1 -\theta \T \boldsymbol{\varrho_{\beta}} )^{-1} \T  \boldsymbol{\Delta \rho^{(2)}} \big) ; \label{eqn:rhodelta_dress}  \\
\big( 1 -\theta \T  \boldsymbol{\rho_{\theta, V_{\delta}}^{(2)}} \big)^{-1} &= \sum_{n=0}^{\infty} \big( \theta (1 -\theta \T \boldsymbol{\varrho_{\beta}} )^{-1} \T  \boldsymbol{\Delta \rho^{(2)}} \big)^n (1 -\theta \T \boldsymbol{\varrho_{\beta}}  )^{-1}. \label{eqn:inv_series}
\end{align}
Indeed, \eqref{eqn:rhodelta_dress} holds by the definition of $\Delta \rho^{(2)}$. Moreover, for $\delta > 0$ small enough, \eqref{eqn:TDeltaRhoOp} yields that the series on the right hand side of \eqref{eqn:inv_series} has operator norm upper bounded by $\sum_{n=0}^{\infty} C^n \delta^n |\log \delta|^{C n}  < \infty$; hence, this series is convergent, so \eqref{eqn:rhodelta_dress} implies that \eqref{eqn:inv_series} holds as an equality of operators on $\mathcal{H}$. This proves existence and boundedness $ (1 -\theta \T  \boldsymbol{\rho_{\theta, V}^{(2)}}  )^{-1} : \mathcal{H} \rightarrow \mathcal{H}$.
\end{proof}

We will use the following lemma (shown at the end of this section) to prove Lemma \ref{lem:spo}. 
\begin{lem}\label{lem:rho1_formula}
Let $V(x)$ be a polynomial satisfying the hypotheses of Lemma \ref{lem:spo}. Denoting $\varsigma_0^{\dr,\delta} \coloneqq (1 -\theta \T  \boldsymbol{\rho_{\theta, V}^{(2)}})^{-1} \varsigma_0 \in \mathcal{H}$, set 
$$\alpha(\theta, V) = \Bigg( \int_{-\infty}^{\infty}\theta \rho_{\theta, V}^{(2)}(x) \varsigma_0^{\dr, \delta}(x) dx \Bigg)^{-1}.$$ 

\noindent Then, for each $x \in \mathbb{R}$, we have
\begin{equation}\label{eqn:rho1rho2}
\rho^{(1)}_{\theta, V }(x) = \partial_{\theta} (\theta \rho^{(2)}_{\theta, V }(x)),
\end{equation}
and 
\begin{equation}\label{eqn:rho1eqns}
\rho_{\theta, V}^{(1)} = \theta \cdot \alpha(\theta, V)  \cdot  (1 -\theta \boldsymbol{\rho_{\theta, V}^{(2)}} \T )^{-1}  \rho_{\theta, V}^{(2)} = \theta \cdot \alpha(\theta, V) \cdot \rho_{\theta, V}^{(2)} (1 -\theta \T \boldsymbol{\rho_{\theta, V}^{(2)}} )^{-1}  \varsigma_0  .
\end{equation}

\noindent In addition, there exist constants $\mathfrak{c}>0$ and $\mathfrak{C}>1$ such that the following holds for all $x \in \mathbb{R}$. First,
\begin{equation}\label{eqn:s0drdf_lem}
\varsigma_0^{\dr,\delta}(x) > \mathfrak{c}, \qquad  |\varsigma_0^{\dr,\delta}(x) - \varsigma_0^{\dr}(x)| \leq \mathfrak{C}  \delta |\log \delta|^{\mathfrak{C}}  \log(|x|+2).
\end{equation}

\noindent Second, denoting by $\mathcal{H}_{\delta}$ the Hilbert space associated with the inner product \eqref{eqn:rho1inner}, we have  
\begin{equation}\label{eqn:H0Hd}
\rho_{\theta, V}^{(1)}(x)\leq \mathfrak{C} \varrho(x); \qquad \mathcal{H} \subseteq \mathcal{H}_{\delta}.
\end{equation}
\end{lem}

 \begin{remark}\label{rmk:dress_inv2}

 Let us explain why \eqref{eqn:rho1eqns} is well defined. First, it is possible to define a suitable space of functions (containing $\rho_{\theta, V}^{(2)} $) that is preserved by $1 -\theta \boldsymbol{\rho_{\theta, V}^{(2)}} \T$, on which $(1 -\theta \boldsymbol{\rho_{\theta, V}^{(2)}} \T )^{-1} $ exists. For example, we may take this space to be $\boldsymbol{\varrho} \cdot \mathcal{H} = \{ \varrho \cdot h : h \in \mathcal{H} \}$; it follows from the bound \eqref{eqn:rhovarrho_rat} below that $1 -\theta \boldsymbol{\rho_{\theta, V}^{(2)}} \T$ preserves this space. On this space, we have the identity of operators 
$$(1 -\theta \boldsymbol{\rho_{\theta, V}^{(2)}} \T )^{-1}  = 1+ \theta \boldsymbol{\rho_{\theta, V}^{(2)}}  (1 -\theta  \T \boldsymbol{\rho_{\theta, V}^{(2)}})^{-1} \T.$$
Moreover, as operators $\mathcal{H} \rightarrow \boldsymbol{\varrho} \cdot \mathcal{H}$, we have
\begin{equation}\label{eqn:DrrDid}
(1 -\theta \boldsymbol{\rho_{\theta, V}^{(2)}} \T )^{-1} \boldsymbol{\rho_{\theta, V}^{(2)}} = \boldsymbol{\rho_{\theta, V}^{(2)}} (1 -\theta \T \boldsymbol{\rho_{\theta, V}^{(2)}}  )^{-1} .
\end{equation}
The formulas above are the standard identities $(1-BA)^{-1}=1+B(1-AB)^{-1}A$ and $(1-BA)^{-1}B=B(1-AB)^{-1}$ with $A=\theta \T$ and $B=\boldsymbol{\rho_{\theta, V}^{(2)}}$. Note that the second equality in \eqref{eqn:rho1eqns} can be deduced from the formula \eqref{eqn:DrrDid} above. In addition, for any functions $f,g \in \mathcal{H}$,
\begin{equation}\label{eqn:Dadj}
\langle g, (1-\theta \boldsymbol{\rho_{\theta, V}^{(2)}} \T)^{-1} (\varrho f) \rangle_{2} = \langle  (1-\theta \T  \boldsymbol{\rho_{\theta, V}^{(2)}} )^{-1} g, \varrho f  \rangle_{2}.
\end{equation}
 \end{remark}

\begin{proof}[Proof of Lemma \ref{lem:spo}]
The proof begins with the relation 
\begin{equation}\label{eqn:spo_start}
-\partial_{t} \mathcal{F}^{(1)}(\theta, V + t g ) = \int_{-\infty}^{\infty} \rho^{(1)}_{\theta, V + t g }(x) g(x) dx
\end{equation}
which follows from Items \ref{item:D2} and \ref{item:D4} in Lemma \ref{thm:mm} above. We will proceed by taking $|t|$ small, expanding 
\begin{equation}\label{eqn:rho1_var}
\rho^{(1)}_{\theta, V + t g } = \rho^{(1)}_{\theta, V  } + t (\delta \rho^{(1)}_{\theta, V  } ) (g)  + o(t)
\end{equation}
and then inserting the result into \eqref{eqn:spo_start} to compute the second derivative of the free energy $\partial_t^2 \mathcal{F}^{(1)}$. In \eqref{eqn:rho1_var}, $(\delta \rho^{(1)}_{\theta, V  } )(g)$ denotes the first order variation of $\rho^{(1)}_{\theta, V  } $ as the potential is varied in the direction $g$. To compute $(\delta \rho^{(1)}_{\theta, V  }) (g)$, we will use the fact \eqref{eqn:rho1eqns} that $\rho_{\theta, V + t g }^{(1)}$ can be obtained from $\rho_{\theta, V + t g }^{(2)}$ by applying the dressing operator $(1 - \theta \mathbf{T} \bm{\rho_{\theta,V}^{(2)}})^{-1}$, together with the fact that  $\rho_{\theta, V + t g }^{(2)}$ solves the Euler--Lagrange equation \eqref{eqn:EL}.

Indeed, we first subtract the two Euler--Lagrange equations \eqref{eqn:EL} for $\rho = \rho^{(2)}_{\theta, V  }$ and $\tilde \rho = \rho^{(2)}_{\theta, V + t g }$ and multiply the difference by $\rho$. Denoting $\Delta \rho = \tilde \rho -  \rho $, we obtain 
\begin{equation}\label{eqn:ELsub}
t \rho(x) g(x) - \theta \rho(x) \T(\Delta \rho)(x) +\rho(x)( \log \tilde \rho(x) - \log \rho(x)) =  \rho(x)(\lambda_{\theta}^{\tilde V} -  \lambda_{\theta}^{ V}).
\end{equation}

\noindent We then claim that 
\begin{align}
\rho(x) (\log \tilde \rho(x) - \log \rho(x)) &=\Delta \rho(x) + \text{err}(x),  \label{eqn:errdef}\\
\text{where} \qquad |\text{err}(x)| &\leq C t^2 |\log t|^{C} \varrho(x), \label{eqn:errbd}
\end{align}
where $\varrho$ is as in Definition \ref{def:alpha_Laxdos}.

Before showing the bound \eqref{eqn:errbd}, we record two estimates. Recall that $\varrho =\alpha \theta \varrho_{\beta} \varsigma_0^{\dr} = \alpha \theta \rho_{\theta, \beta x^2/2}^{(2)} \varsigma_0^{\dr}$ (by \eqref{rho2}). Note in addition that $\varrho_{\beta}(x) \leq C \varrho(x)$ by the lower bound $\varsigma_0^{\dr}(x)>c$ in Lemma \ref{lem:ve_tail}. Recall $V(x) = \beta x^2/2 +\delta x^{2 n}$. We apply the estimates \eqref{eqn:rho_apriori_Tbd_stmt} and \eqref{eqn:rho_apriori_lambda_stmt} from Lemma \ref{lem:rho_apriori}, once with $\rho =  \varrho_{\beta}$ and $\tilde{\rho} = \rho_{\theta, V + t g}^{(2)}$, and once with $\rho = \varrho_{\beta}$ and $\tilde{\rho} = \rho_{\theta, V }^{(2)}$, in each case together with \eqref{eqn:rho_ptwise}, to obtain the following estimates, valid for all $t \leq \delta$ (for $\delta$ sufficiently small), and $x \in \mathbb{R}$, 
\begin{equation}\label{eqn:rhovarrhobeta_rat}
 \frac{\rho_{\theta, V + t g}^{(2)}(x)}{ \varrho_{\beta}(x)} \leq C e^{-  \delta x^{2 n}/2},\qquad  \frac{\rho_{\theta, V }^{(2)}(x)}{ \varrho_{\beta}(x)} \leq C e^{- \delta x^{2 n}/2}.
\end{equation}
Since $\varrho_{\beta}(x) \leq C \varrho(x)$, this implies
\begin{equation}\label{eqn:rhovarrho_rat}
 \frac{\rho_{\theta, V + t g}^{(2)}(x)}{ \varrho(x)} \leq C e^{-  \delta x^{2 n}/2},\qquad  \frac{\rho_{\theta, V }^{(2)}(x)}{ \varrho(x)} \leq C e^{- \delta x^{2 n}/2}.
\end{equation}

Now we will show \eqref{eqn:errbd}, and as stated above that display we set $\rho = \rho^{(2)}_{\theta, V  }$ and $\tilde \rho = \rho^{(2)}_{\theta, V + t g }$. With $B(x)=\lambda_\theta^{\tilde V}-\lambda_\theta^V-t g(x)+\theta \T(\Delta \rho)(x)$, by \eqref{eqn:ELsub} we have $\log \tilde \rho(x)-\log \rho(x)=B(x)$, so $\Delta \rho(x)=\rho(x)(e^{B(x)}-1)$. If $|x| \le |\log t|$, then \eqref{eqn:rho_apriori_lambda_stmt}, \eqref{eqn:rho_apriori_Tbd_stmt}, and the polynomial growth of $g$ imply $|B(x)| \le C |t| |\log t|^C$, so 
$$|\text{err}(x)| = \rho(x)|e^{B(x)}-1-B(x)|\leq C \rho(x) |t|^2 |\log t|^C,$$ 
which satisfies \eqref{eqn:errbd} by \eqref{eqn:rhovarrho_rat}. If instead $|x|>|\log t|$, then we proceed as follows. If $t g(x)\ge 0$, then $B(x)\le C |t| |\log t|^C \leq C$ (here we again use the bounds \eqref{eqn:rho_apriori_lambda_stmt} and \eqref{eqn:rho_apriori_Tbd_stmt}, which hold for all $x$, and that $t g(x)$ appears with a negative sign in the definition of $B(x)$) and so by the second bound in \eqref{eqn:rhovarrho_rat},
\begin{flalign*}
\left| \rho(x)(\log \rho(x)-\log \tilde \rho(x))-\Delta \rho(x) \right| =\rho(x)| e^{B(x)}-1-B(x) | & \leq C e^{-\delta x^{2n}/2} \varrho(x) |x|^C \\
& \leq C t^2 |\log t|^{C} \varrho(x).
\end{flalign*}

\noindent If instead $t g(x)<0$, then rewriting $\Delta \rho(x)=\tilde \rho(x)(1-e^{-B(x)})$ and using the first bound in \eqref{eqn:rhovarrho_rat}  to bound $\tilde{\rho}(x)$ gives
\[
\left| \rho(x)(\log \rho(x)-\log \tilde \rho(x))-\Delta \rho(x) \right|=\tilde \rho(x)\bigl|1-(1+B(x))e^{-B(x)}\bigr| \leq C t^2 |\log t|^{C} \varrho(x).
\]
In both cases above, we have also used that $g$ has polynomial growth to bound $|B(x)| \leq C |x|^C$. Thus, we have shown \eqref{eqn:errbd}.

So, denoting $\Delta \lambda =  \lambda_{\theta}^{\tilde V}-\lambda_{\theta}^V $ and recasting \eqref{eqn:ELsub} in operator notation, we obtain
\begin{equation}\label{eqn:elsub_expansion}
  (1 - \theta \boldsymbol{\rho} \T) \Delta \rho =  \rho \Delta \lambda - t \rho g - \text{err}
\end{equation}
which, upon multiplying both sides on the left by $1 + \theta \boldsymbol{\rho} (1-\theta \T \boldsymbol{\rho})^{-1} \T$ (see Remark \ref{rmk:dress_inv2} and in particular \eqref{eqn:DrrDid}),  implies 
\begin{equation}\label{eqn:Deltarho_eq}
 \Delta \rho =   \Delta \lambda \cdot  \bm{\rho} (1 - \theta  \T \boldsymbol{\rho} )^{-1} \varsigma_0 - t \bm{\rho} (1 - \theta  \T \boldsymbol{\rho} )^{-1} g + O(|t|^2 |\log t|^C).
\end{equation}
By \eqref{eqn:errbd}, the error term on the right hand side above, denoted $O(|t|^2 |\log t|^C)$, is a function of the form $\varrho \cdot h$, for some $h \in \mathcal{H}$ which satisfies
\begin{equation}\label{eqn:hHbd}
\left\| h \right\|_{\mathcal{H}} \leq C  |t|^2 |\log t|^C .
\end{equation}

Next, we must specify $\Delta \lambda$. We observe that it is uniquely determined by the property that $\int_{-\infty}^{\infty}\Delta \rho(x) dx = 0$: By \eqref{eqn:Deltarho_eq}, \eqref{eqn:hHbd}, and \eqref{eqn:s0drdf_lem}, 
\begin{equation}\label{eqn:deltalambda}
\Delta \lambda =  t \cdot \frac{\int_{-\infty}^{\infty}(1 - \theta  \T \boldsymbol{\rho} )^{-1} ( g)(x) \rho(x) dx}{\int_{-\infty}^{\infty}(1 - \theta \T \boldsymbol{\rho})^{-1} (\varsigma_0)(x) \rho(x) dx} + O(|t|^2 |\log t|^C)= t \delta \lambda(g) + O(|t|^2 |\log t|^C)  ,
\end{equation}
where we have denoted the coefficient of $t$ in the first term on the right hand side of the first equation as $\delta \lambda(g)$. We obtain (substituting back into \eqref{eqn:Deltarho_eq})
\begin{equation}\label{eqn:rho2var}
 \Delta \rho = t \left( \delta \lambda(g) \cdot \bm{\rho} (1 - \theta  \T \boldsymbol{\rho})^{-1} \varsigma_0 - \bm{\rho}  (1 - \theta  \T \boldsymbol{\rho})^{-1}  g \right) + \varrho \cdot h, \qquad \|h\|_{\mathcal{H}} \leq C t^2 (\log t)^C
\end{equation}
where we have used the stated property of the error term in \eqref{eqn:Deltarho_eq}, as well as the second bound in \eqref{eqn:rhovarrho_rat} and the fact that $(1 - \theta  \T \boldsymbol{\rho})^{-1} \varsigma_0 \in \mathcal{H}$. From this point, we follow the calculations of \cite[Section 4]{Spo20} to complete the proof.

Finally, we are in a position to use \eqref{eqn:rho1eqns} to compute the variation $\delta \rho^{(1)}(g)$ in Equation \eqref{eqn:rho1_var}. Denoting by $(\delta \rho^{(2)})(g)$ the coefficient of $t$ in \eqref{eqn:rho2var} (which is the variation of $\rho_{\theta, V}^{(2)}$ in direction $g$), we claim for any $f \in \mathcal{H}$ that
\begin{multline}\label{eqn:rho2_var}
\left\langle ( 1 - \theta \boldsymbol{\rho_{\theta, V+t g}^{(2)}} \T )^{-1} \theta \rho_{V +t g}^{(2)} - ( 1 - \theta \boldsymbol{\rho_{\theta, V}^{(2)}} \T )^{-1} \theta \rho_{V }^{(2)} , f\right\rangle_2\\
=t  \left\langle (1 -\theta \T \boldsymbol{\rho_{\theta, V}^{(2)}} )^{-1} \varsigma_0,  \theta \boldsymbol{(\delta \rho^{(2)})(g)} (1-  \theta \T \boldsymbol{\rho_{\theta, V}^{(2)}} )^{-1} f \right\rangle_2 + O(t^2 |\log t|^{C})
\end{multline}
where the constant in the $O(t^2 |\log t|^{C})$ error may depend on $f$ and $g$, and where $\boldsymbol{(\delta \rho^{(2)})(g)}$ denotes the multiplication operator corresponding to the function $(\delta \rho^{(2)})(g)$ (which is the coefficient of $t$ in \eqref{eqn:rho2var}). To justify \eqref{eqn:rho2_var}, we perform a computation and invoke an expansion of the form \eqref{eqn:inv_series} and the bound on the error in \eqref{eqn:rho2var} in order to bound lower order terms. Specifically, defining $\mathbf{D}_V \coloneqq (1-  \theta \T \boldsymbol{\rho_{\theta, V}^{(2)}} )^{-1} $, and $\mathbf{D}_{V+t g} \coloneqq (1-  \theta \T \boldsymbol{\rho_{\theta, V+ t g}^{(2)}} )^{-1} $, we have
\begin{align*}
 &\theta \boldsymbol{\rho_{\theta, V +t g}^{(2)}}( 1 - \theta  \T \boldsymbol{\rho_{\theta, V+t g}^{(2)}})^{-1} - \theta \boldsymbol{\rho_{\theta, V }^{(2)}}  ( 1 - \theta  \T \boldsymbol{\rho_{\theta, V}^{(2)}} )^{-1}  \\
&=  \theta \boldsymbol{\Delta  \rho}( 1 - \theta  \T \boldsymbol{\rho_{\theta, V}^{(2)}})^{-1} +  \theta \boldsymbol{\rho_{\theta, V+t g }^{(2)}}  \big( ( 1 - \theta  \T \boldsymbol{\rho_{\theta, V+t g}^{(2)}} )^{-1} -  ( 1 - \theta  \T \boldsymbol{\rho_{\theta, V}^{(2)}} )^{-1}  \big) \\
&=\theta \boldsymbol{\Delta  \rho}  \mathbf{D}_{V}+   \theta \boldsymbol{\rho_{\theta, V+t g }^{(2)}} \mathbf{D}_{V+ t g} \left(    \theta \T \boldsymbol{\Delta \rho} \right) \cdot  \mathbf{D}_V\\
&=\theta \boldsymbol{\Delta  \rho}  \mathbf{D}_{V}+  ( 1 - \theta \boldsymbol{\rho_{\theta, V+t g}^{(2)}} \T )^{-1} \big(  \theta \boldsymbol{\rho_{\theta, V+t g }^{(2)}}   \theta  \T  \boldsymbol{\Delta \rho} \big) \cdot  \mathbf{D}_V \\
&=  ( 1 - \theta \boldsymbol{\rho_{\theta, V+t g}^{(2)}} \T )^{-1} \left(  \theta  \boldsymbol{\Delta \rho} \right)   \mathbf{D}_V .
\end{align*}
By applying \eqref{eqn:Dadj} and \eqref{eqn:rhovarrho_rat} to move the operators onto $f$ in the left hand side of \eqref{eqn:rho2_var}, and then using the display above and \eqref{eqn:Dadj} and \eqref{eqn:rhovarrho_rat} again, the left hand side of \eqref{eqn:rho2_var} is given by $ \langle \theta  \boldsymbol{\Delta \rho} \mathbf{D}_V f, ( 1 - \theta \T \boldsymbol{\rho_{\theta, V + t g}^{(2)}}  )^{-1} \varsigma_0  \rangle_2 $. This expression in turn by \eqref{eqn:rho2var} satisfies
\begin{equation}\label{eqn:rho2_var2}
\left \langle \theta  \boldsymbol{\Delta \rho} \mathbf{D}_V f, ( 1 - \theta \T \boldsymbol{\rho_{\theta, V + t g}^{(2)}}  )^{-1} \varsigma_0 \right \rangle_2 
= \theta t \left \langle   \delta \rho^{(2)} (g)  \mathbf{D}_V f, ( 1 - \theta \T \boldsymbol{\rho_{\theta, V}^{(2)}}  )^{-1} \varsigma_0 \right \rangle_2 +  O(|t|^2 |\log t|^C).
\end{equation}
To obtain final equation above, we have performed two steps: First, we have used an expansion similar to \eqref{eqn:inv_series}, to see (as an equality of operators on $\mathcal{H}$),
\begin{equation}\label{eqn:series_inv_t}
\mathbf{D}_{V + t g} =( 1 -\theta  \T \boldsymbol{\rho_{\theta, V+ t g}^{(2)}} )^{-1} = \sum_{n=0}^{\infty} \big( ( 1 -\theta  \T \boldsymbol{\rho_{\theta, V}^{(2)}} )^{-1}   \T \boldsymbol{\Delta \rho}  \big)^n \cdot ( 1 -\theta  \T  \boldsymbol{\rho_{\theta, V}^{(2)}} )^{-1}.
\end{equation}

\noindent By \eqref{eqn:rho2var}, we may bound the operator norm of the sum of the terms corresponding to indices $n \geq 1$ by $C |t|$, and thus, 
\begin{equation}\label{eqn:series_inv_t_s0bd}
\bigg\| \sum_{n=1}^{\infty} \big( ( 1 -\theta  \T \boldsymbol{\rho_{\theta, V}^{(2)}}  )^{-1}   \T \boldsymbol{\Delta \rho}  \big)^n \cdot ( 1 -\theta  \T  \boldsymbol{\rho_{\theta, V}^{(2)}} )^{-1} \varsigma_0 \bigg\|_{\mathcal{H}} \leq C t.
\end{equation}
Moreover, since $f \in \mathcal{H}$, both $\mathbf{D}_V$ and $\mathbf{D}_{V + t g} $ are bounded on $\mathcal{H}$, and the error $\varrho \cdot h$ in \eqref{eqn:rho2var} satisfies the last bound stated there, we have 
\begin{align*}
\left \langle \theta  \left( \Delta \rho - t \delta\rho(g) \right) \mathbf{D}_V f, \mathbf{D}_{V+t g} \varsigma_0 \right \rangle_2 &\leq C t^2 |\log t|^C \\
t \left \langle \theta \delta\rho(g) \mathbf{D}_V f,\left(\mathbf{D}_{V+t g}- \mathbf{D}_{V} \right)\varsigma_0 \right \rangle_2 &\leq C t^2.
\end{align*}
In the second bound above we used \eqref{eqn:series_inv_t_s0bd}. The two bounds above show \eqref{eqn:rho2_var2} and thus \eqref{eqn:rho2_var}.

Note that $\delta \lambda(g) = \mu_g$ by \eqref{eqn:Dadj}, \eqref{eqn:deltalambda}, and \eqref{eqn:rho1eqns}: Indeed, by \eqref{eqn:Dadj} and \eqref{eqn:DrrDid} the quotient in \eqref{eqn:deltalambda} may be rewritten as
\[
\frac{\langle g^{\dr,\delta}, \theta \rho_{\theta,V}^{(2)} \varsigma_0 \rangle_2}{\langle \varsigma_0, \theta \rho_{\theta,V}^{(2)} \varsigma_0^{\dr,\delta}\rangle_2}
= \frac{\langle g, \theta \rho_{\theta,V}^{(2)} \varsigma_0^{\dr,\delta} \rangle_2}{\langle \varsigma_0, \theta \rho_{\theta,V}^{(2)} \varsigma_0^{\dr,\delta}\rangle_2} ,
\]
and \eqref{eqn:rho1eqns} identifies this ratio with $\int_{-\infty}^{\infty} g(x)\rho_{\theta,V}^{(1)}(x)\,dx=\mu_g$. Thus, as a consequence of \eqref{eqn:rho2_var} and the formula for $\delta \rho^{(2)}(g)$ from \eqref{eqn:rho2var}, using the notation $f^{\dr,\delta} \coloneqq \mathbf{D}_{V} f$ to distinguish it from the dressing $f^{\dr}$ used throughout the rest of the paper, we have
\begin{equation}
\frac{d}{d t} \langle f, (1-\theta \boldsymbol{\rho_{\theta, V+t g}^{(2)}} \T)^{-1} \theta \rho_{\theta, V+t g}^{(2)} \rangle_2 \big|_{t=0} = \langle  \varsigma_0^{\dr,\delta},  \theta \rho_{\theta, V}^{(2)} (\mu_g  \varsigma_0^{\dr,\delta} - g^{\dr,\delta} ) f^{\dr,\delta} \rangle_2 . \label{eqn:deriv_calc}
\end{equation}

\noindent  So, applying \eqref{eqn:deriv_calc} once with $f = g$ and once with $f = \varsigma_0$, and using \eqref{eqn:DrrDid} with \eqref{eqn:rho1eqns}, we obtain
\begin{align}
\frac{d}{d t} \langle g ,  \rho_{\theta, V+t g}^{(1)}  \rangle_2 \big|_{t=0} &= \frac{d}{d t} \bigg( \frac{ \langle g, (1-\theta \boldsymbol{\rho_{\theta, V+t g}^{(2)}} \T)^{-1} \theta \rho_{\theta, V+t g}^{(2)}  \rangle_2}{\langle \varsigma_0, (1-\theta \boldsymbol{\rho_{\theta, V+t g}^{(2)}} \T)^{-1} \theta \rho_{\theta, V + t g}^{(2)} \rangle_2} \bigg) \bigg|_{t=0} \notag \\
&=  \frac{\langle \varsigma_0^{\dr,\delta}\theta \rho_{\theta, V}^{(2)},   (\mu_g \varsigma_0^{\dr,\delta}- g^{\dr,\delta}) g^{\dr,\delta} \rangle_2}{\langle \varsigma_0^{\dr,\delta}\theta \rho_{\theta, V }^{(2)} , \varsigma_0  \rangle_2} \notag \\
& \qquad  - \frac{\langle \varsigma_0^{\dr,\delta} \theta \rho_{\theta, V}^{(2)}, g \rangle_2}{\langle \varsigma_0^{\dr,\delta} \theta \rho_{\theta, V }^{(2)},   \varsigma_0 \rangle_2^2 }  \langle \varsigma_0^{\dr,\delta},  \theta \rho_{\theta, V}^{(2)} (\mu_g \varsigma_0^{\dr,\delta}- g^{\dr,\delta}) \varsigma_0^{\dr,\delta} \rangle_2. \label{eqn:apply_deriv_calc}
\end{align}
Upon simplifying and using Lemma \ref{lem:rho1_formula}, the final expression in \eqref{eqn:apply_deriv_calc} above is quickly verified to match the negative of the right hand side of \eqref{eqn:var_thm}. This, with \eqref{eqn:spo_start}, proves \eqref{eqn:var_thm}.

The equations \eqref{eqn:cov_0nthm} and \eqref{eqn:cov_00thm} can be proved in a similar way (by following \cite[Equations (4.3) and (4.7)]{Spo20}); we omit further details. 
\end{proof}

Next, we outline a proof of Lemma \ref{lem:rho1_formula}.

\begin{proof}[Proof of Lemma \ref{lem:rho1_formula} (Outline)]
First, we outline how to show \eqref{eqn:rho1eqns}; it will follow from arguments similar to those in the proof of Lemma \ref{lem:spo} above (and following the computations in \cite{Spo20}). Then we show \eqref{eqn:s0drdf_lem} and \eqref{eqn:H0Hd}, and next \eqref{eqn:rho1rho2}.

More specifically, as in the proof of Lemma \ref{lem:spo}, one works on $\mathcal{H}$ from Definition \ref{def:inner_prod}. Let $\theta_1 = \theta$, and let $|\theta_1-\theta_2|$ be small, and denote $\Delta \lambda \coloneqq \lambda_{\theta_2}^V- \lambda_{\theta_1}^V$, and $\Delta \rho(x) \coloneqq  \rho_{\theta_2, V}^{(2)} (x) -  \rho_{\theta_1, V}^{(2)} (x)$. Following the proof of Lemma \ref{lem:spo} up until \eqref{eqn:elsub_expansion}, we obtain
\begin{equation}\label{eqn:ELdftheta}
\rho_{\theta_1, V}^{(2)} (\theta_2-\theta_1) \T \rho_{\theta_1, V}^{(2)} + \rho_{\theta_1, V}^{(2)} (\lambda_{\theta_2}^V- \lambda_{\theta_1}^V) = (1-\theta_1  \rho_{\theta_1, V}^{(2)} \T) \Delta \rho + \varrho \cdot h,
\end{equation}
where $h$ is a function bounded by $|h(x)| \leq C |\theta_1-\theta_2|^2 \cdot \left|\log |\theta_1- \theta_2| \right|^C (1+|x|)^{C}$ (where $C$ may depend on $\theta_1$). Rearranging \eqref{eqn:ELdftheta} (after first multiplying both sides by $(1-\theta_1  \bm{\rho_{\theta_1, V}^{(2)}} \T)^{-1}$), dividing by $\theta_2-\theta_1$, and using \eqref{eqn:DrrDid}, one obtains 
\begin{multline}\label{eqn:deltarho_theta_formula}
\frac{\Delta \rho}{\theta_2-\theta_1} \\
= \rho_{\theta_1,V}^{(2)} (1-\theta_1 \T \boldsymbol{\rho_{\theta_1,V}^{(2)}})^{-1}\T \rho_{\theta_1,V}^{(2)}
+ \frac{\lambda_{\theta_2}^V-\lambda_{\theta_1}^V}{\theta_2-\theta_1}\,\rho_{\theta_1,V}^{(2)} (1-\theta_1 \T \boldsymbol{\rho_{\theta_1,V}^{(2)}})^{-1}\varsigma_0
- \frac{(1-\theta_1  \rho_{\theta_1,V}^{(2)} \T)^{-1}(\varrho \cdot h)}{\theta_2-\theta_1}.
\end{multline}
Note that upon testing both sides of this identity against a continuous function of polynomial growth, the last term is upper bounded by $C |\theta_1-\theta_2| \cdot \left|\log |\theta_1- \theta_2| \right|^C$ by the stated bound on $h$. 

Now observe that the quantity $\lambda_{\theta_2}^V-\lambda_{\theta_1}^V$ is obtained by integrating both sides of the equation and using the fact that $\int_{-\infty}^{\infty}\Delta \rho(x)\,dx = 0$. Thus, we obtain 
\begin{flalign} 
\label{lambda_lambda} 
\frac{\lambda_{\theta_2}^V-\lambda_{\theta_1}^V}{\theta_2-\theta_1}
=
-\frac{\int_{-\infty}^{\infty}\rho_{\theta_1,V}^{(2)}(x)\bigl((1-\theta_1 \T \boldsymbol{\rho_{\theta_1,V}^{(2)}})^{-1}\T \rho_{\theta_1,V}^{(2)}\bigr)(x)\,dx}
{\int_{-\infty}^{\infty}\rho_{\theta_1,V}^{(2)}(x)\bigl((1-\theta_1 \T \boldsymbol{\rho_{\theta_1,V}^{(2)}})^{-1}\varsigma_0\bigr)(x)\,dx}
+ O( |\theta_1-\theta_2| |\log |\theta_1-\theta_2||^C).
\end{flalign} 

\noindent Thus, the coefficient of
\[
\rho_{\theta_1,V}^{(2)} (1-\theta_1 \T \boldsymbol{\rho_{\theta_1,V}^{(2)}})^{-1}\varsigma_0
\]
has a limit as $\theta_2$ tends to $\theta_1$, meaning that the left-hand side of \eqref{eqn:deltarho_theta_formula} does as well (after testing against a function). Substituting \eqref{lambda_lambda} for $(\lambda_{\theta_2}^V-\lambda_{\theta_1}^V)/(\theta_2-\theta_1)$ into 
\begin{flalign} 
\label{thetarho} 
\frac{\theta_2 \rho_{\theta_2,V}^{(2)}-\theta_1 \rho_{\theta_1,V}^{(2)}}{\theta_2-\theta_1}
= \rho_{\theta_2,V}^{(2)} + \theta_1 \frac{\Delta \rho}{\theta_2-\theta_1},
\end{flalign} 

\noindent which, after testing against a polynomial $g$, gives the right side of \eqref{eqn:difftheta_id} as $\theta_2$ tends to $\theta_1$. Using this and
\[
(1-\theta_1 \T \boldsymbol{\rho_{\theta_1,V}^{(2)}})^{-1}\varsigma_0
= \varsigma_0 + \theta_1 (1-\theta_1 \T \boldsymbol{\rho_{\theta_1,V}^{(2)}})^{-1}\T \rho_{\theta_1,V}^{(2)};
\]
 testing the resulting expression for $(\theta_2 \rho_{\theta_2,V}^{(2)}-\theta_1 \rho_{\theta_1,V}^{(2)}) / (\theta_2-\theta_1)$ against $g$; letting $\theta_2$ tend to $ \theta_1$; and using the fact that, if $\nu$ denotes the limit of the constant in \eqref{lambda_lambda}, then $1+\theta_1 \nu = \theta_1 \alpha(\theta_1,V)$, we obtain \eqref{eqn:rho1eqns} (first upon testing against any polynomial $g$ and thus pointwise).

Next, let us verify \eqref{eqn:s0drdf_lem}. Denoting $\varsigma_0^{\dr,\delta} \coloneqq (1 -\theta \T  \boldsymbol{\rho_{\theta, V}^{(2)}})^{-1} \varsigma_0$, and $\varsigma_0^{\dr}(x)$ the dressing of $\varsigma_0$ at $\delta = 0$. In what follows it will be useful to record the bound 
\begin{equation}\label{eqn:vsdrc}
\varsigma_0^{\dr}(x) > c \log(|x|+2),
\end{equation}
valid for all $x \in \mathbb{R}$. This follows from the lower bound $\varsigma_0^{\dr}(x) > c $ in Lemma \ref{lem:ve_tail}, together with the fact that $\varsigma_0^{\dr}(x) = 1 + 2 \theta \int_{-\infty}^{\infty} \log|x-y| \varsigma_0^{\dr}(y) \varrho_{\beta}(y) d y $.
 
Specifically, denoting $\mathbf{D} = (1 -\theta \T  \boldsymbol{\varrho_{\beta}})^{-1}$ the $\delta = 0$ dressing operator and $\Delta \rho = \rho_{\theta, V_{\delta}}^{(2)} - \varrho_{\beta}$, we have the following equalities in $\mathcal{H}$ (see \eqref{eqn:rhodelta_dress} and \eqref{eqn:inv_series} in the proof of Lemma \ref{lem:dress_existence}):
\begin{align}
\varsigma_0^{\dr,\delta}  &= (1 -\theta \mathbf{D} \T \boldsymbol{\Delta \rho} )^{-1} \mathbf{D} \varsigma_0 \notag \\
&= \left(\sum_{n=0}^{\infty}\left( \theta \mathbf{D} \T \boldsymbol{\Delta \rho}  \right)^n \right)\mathbf{D} \varsigma_0 = \varsigma_0^{\dr} + \theta \mathbf{D} \T \boldsymbol{\Delta \rho}  \left(\sum_{n=0}^{\infty}\left( \theta \mathbf{D} \T \boldsymbol{\Delta \rho}  \right)^n  \right)  \varsigma_0^{\dr}. \label{eqn:sdrd}
\end{align}
We now bound the second term on the right-hand side above. To do this, we first explain how to bound $ \theta \mathbf{D} \T \boldsymbol{\Delta \rho^{(2)}} \varsigma_0^{\dr}$. First, 
$$
\int_{-\infty}^{\infty} 2 \log|x-y| \Delta \rho^{(2)}(y) \varsigma_0^{\dr}(y) dy \leq C  \delta |\log \delta|^C  \log(|x|+2);
$$
we have used $|\varsigma_0^{\dr}(x)| \leq C \log(|x|+2)$ (by the first part of Lemma \ref{lem:ve_tail}), and the bound \eqref{eqn:deltarho2bd}. Using the bound in the display above together with Lemma \ref{lem:ve_tail}, $\left|\mathbf{D} \T \boldsymbol{\Delta \rho^{(2)}} \varsigma_0^{\dr}(x) \right| \leq C'  \delta |\log \delta|^C  \log(|x|+2)$. Bounding the terms involving higher powers of $ \theta \mathbf{D} \T \boldsymbol{\Delta \rho^{(2)}}$ is similar, and leads to higher powers of $C' \delta |\log \delta|^C$. Ultimately, using the obtained bounds in \eqref{eqn:sdrd}, we obtain 
\begin{equation}\label{eqn:s0drdf}
|\varsigma_0^{\dr,\delta}(x) - \varsigma_0^{\dr}(x)| \leq C  \delta |\log \delta|^C  \log(|x|+2),
\end{equation}
which is the second bound in \eqref{eqn:s0drdf_lem}. The bound \eqref{eqn:s0drdf}, together with \eqref{eqn:vsdrc}, implies that $\varsigma_0^{\dr,\delta}(x)  > c $, which is the first bound in \eqref{eqn:s0drdf_lem}. 

Next, we prove the bound in \eqref{eqn:H0Hd} (which implies the second inclusion). This follows from multiplying both sides of \eqref{eqn:rhovarrhobeta_rat} by $\varsigma_0^{\dr, \delta}/\varsigma_0^{\dr}$; indeed, \eqref{eqn:rho1eqns} implies that $\rho_{\theta, V}^{(1)}= \alpha(\theta, V) \theta \rho_{\theta, V}^{(2)} \cdot \varsigma_{0}^{\dr, \delta} $, and then the claim follows from \eqref{eqn:vsdrc} and the combination of \eqref{eqn:s0drdf} with the $f = \varsigma_0$ case of Lemma \ref{lem:ve_tail} (to bound $|\varsigma_0^{\dr, \delta}(x)| \leq C \log(|x|+2)$).

Finally, we show the pointwise relation \eqref{eqn:rho1rho2}. Having identified $\rho_{\theta_1, V}^{(1)}$ via \eqref{eqn:rho1eqns}, in order to show \eqref{eqn:rho1rho2} it remains by \eqref{thetarho} to verify that the last term in \eqref{eqn:deltarho_theta_formula} tends to $0$ pointwise, as $\theta_2$ tends to $\theta_1$. To that end, we claim for all  $f \in \mathcal{H}$ that 
\begin{equation}\label{eqn:fdrd_ptwisebd}
|f^{\dr, \delta}(x)| \leq C \big( f(x) + \log(|x|+2) \cdot \|f\|_{\mathcal{H}} \big);
\end{equation}
 
 \noindent observe at $\delta = 0$ that this follows from Lemma \ref{lem:ve_tail}. The bound \eqref{eqn:fdrd_ptwisebd} follows from the estimates
 \begin{align*}
 |f^{\dr, \delta}(x)| &\leq | f(x)| +2  \theta \int_{-\infty}^{\infty} \left| \log|x-y| f^{\dr,\delta}(y) \right|  \rho_{\theta, V}^{(2)}(y) dy \\
 &\leq |f(x)| +C \left( \int_{-\infty}^{\infty} \left| \log|x-y| \right|^2 \varrho(y) dy\right)^{1/2}\left( \int_{-\infty}^{\infty} \left| f(y) \right|^2 \varrho(y) dy\right)^{1/2};
 \end{align*}
 the first inequality is by the definition of $f^{\dr, \delta}$; the second is by Cauchy-Schwarz, \eqref{eqn:H0Hd}, \eqref{eqn:rho1eqns}, the bound $\varsigma_0^{\dr, \delta}(x) > c$ (from \eqref{eqn:s0drdf_lem}), and boundedness of $(1- \theta  \mathbf{T} \boldsymbol{\rho_{\theta, V}^{(2)}})^{-1} $ on $\mathcal{H}$ (by Lemma \ref{lem:dress_existence}). Thus, \eqref{eqn:fdrd_ptwisebd} follows by bounding $( \int_{-\infty}^{\infty} \left| \log|x-y| \right|^2 \varrho(y) dy )^{1/2} \leq C \log(|x|+2)$. 

 The estimate \eqref{eqn:fdrd_ptwisebd} applied with $f$ equal to the function $\mathbf{T}(\varrho h)$ (where $h$ is as in \eqref{eqn:ELdftheta} and \eqref{eqn:deltarho_theta_formula}) implies that for each $x$ the second order term of \eqref{eqn:ELdftheta} tends to $0$ upon multiplying both sides of the identity by $(1- \theta \boldsymbol{\rho_{\theta, V}^{(2)}} \mathbf{T})^{-1} =  1+ \theta \boldsymbol{\rho_{\theta, V}^{(2)}}  (1 -\theta  \T \boldsymbol{\rho_{\theta, V}^{(2)}})^{-1} \T$, dividing by $\theta_2-\theta_1$, and letting $\theta_2$ tend to $\theta_1$. This confirms \eqref{eqn:rho1rho2} and thus the lemma.
\end{proof}

Now let us show Lemma \ref{lem:lim_var_formula}.

\begin{proof}[Proof of Lemma \ref{lem:lim_var_formula}]

Set $n = \deg p$. We consider a family of potentials which perturb the potential of thermal equilibrium from Definition \ref{def:inf_thermal_eq}. Specifically, for $\delta \geq 0$, let $V_{\delta} (x) = \beta x^2/2 + \delta x^{ 2n}$ (this will be done in order to apply Lemma \ref{lem:spo} at $g = p$, which requires that $\deg V > \deg p = n$). Letting $\rho_{\theta,V_\delta}^{(1)}$ and $\rho_{\theta,V_{\delta}}^{(2)}$ denote the equilibrium measures defined below \eqref{1f} and \eqref{F2}, define
\begin{equation}\label{eqn:mupdf}
\mu_{p,\delta} = \int_{-\infty}^{\infty} p(\lambda)\rho_{\theta, V_{\delta}}^{(1)}(\lambda) d\lambda.
\end{equation}

\noindent In what follows, we further let $\varrho_{\beta}=\rho_{\theta,\beta x^2/2}^{(2)}$ denote the density associated with $V = \beta x^2/2$ from Definition \ref{def:alpha_Laxdos} and $\varrho=\rho_{\theta,\beta x^2/2}^{(1)}$ is the corresponding equilibrium measure of the Lax matrix under thermal equilibrium (also from Definition \ref{def:alpha_Laxdos}). 

We start by proving the formula \eqref{eqn:var_explicit1}
 for the limiting variance. We will proceed by first adding a small perturbation to $\beta x^2/2$ so that we may invoke Item \ref{item:D3} of Lemma \ref{thm:mm}, and then utilizing Lemma \ref{lem:spo}. By the holomorphicity of $\mathcal{F}^{(1)}(\theta, V_{\delta} + \mathrm{i} t p)$ in $t$ when $|t| \le c$ (recall the third part of Lemma \ref{thm:mm}), \eqref{eqn:var_thm} yields
 \begin{equation}\label{eqn:fe_var_ID}
     \partial_t^2 \mathcal{F}^{(1)}(\theta, V_{\delta} + \mathrm{i} t p)|_{t=0} =  \langle (1-\theta \mathbf{T}  \boldsymbol{\rho_{\theta,V_{\delta}}^{(2)}})^{-1} (p - \mu_{p,\delta}) , (1-\theta \mathbf{T}  \boldsymbol{\rho_{\theta,V_{\delta}}^{(2)}})^{-1} (p - \mu_{p,\delta}) \rangle_{\rho_{\theta,V_{\delta}}^{(1)}}
 \end{equation}
 thus obtained holds for any $\delta > 0$ small enough.
 
Next, we let $\delta$ tend to $0$, for which we must establish continuity in $\delta$ of both sides of the identity \eqref{eqn:fe_var_ID}. Item \ref{item:D6} of Lemma \ref{thm:mm} provides continuity of its left-hand side. To verify continuity of the right-hand side of \eqref{eqn:fe_var_ID}, we first claim that
\begin{align}
\big\| (1-\theta \mathbf{T}  \boldsymbol{\rho_{\theta,V_{\delta}}^{(2)}})^{-1}  - (1-\theta \mathbf{T}  \boldsymbol{\rho_{\theta,V_{0}}^{(2)}})^{-1}  \big\|_{\text{op}, \mathcal{H}} &\leq C \delta |\log \delta|^C  \label{eqn:DdD0}\\
\big| \rho_{\theta,V_{\delta}}^{(1)}(x) - \varrho(x) \big| &\leq C \delta |\log \delta|^C,  \quad \text{ for all } x \in \mathbb{R},  \label{eqn:rhodrho0pw}
\end{align}

\noindent where $\| \cdot \|_{\text{op}, \mathcal{H}}$ denotes the operator norm on $\mathcal{H}$. Indeed, \eqref{eqn:DdD0} follows from the series expansion \eqref{eqn:inv_series} and the bound \eqref{eqn:TDeltaRhoOp} shown in the proof of Lemma \ref{lem:dress_existence}. The bound \eqref{eqn:rhodrho0pw} follows quickly from \eqref{eqn:rho1eqns}, upon using \eqref{eqn:deltarho2bd}, \eqref{eqn:s0drdf_lem}, and the bounds $|\varsigma_0^{\dr}(x)|, |\varsigma_0^{\dr,\delta}(x)| \leq C \log(|x|+2)$ (which hold by Lemma \ref{lem:ve_tail} and \eqref{eqn:s0drdf_lem}). 

Now let use \eqref{eqn:DdD0} and \eqref{eqn:rhodrho0pw} to deduce imply continuity of the right-hand side of \eqref{eqn:fe_var_ID}. By \eqref{eqn:DdD0} and \eqref{eqn:H0Hd}, we have for $\gamma \in \{0, \delta\}$ that
\begin{multline}\label{eqn:DdeltaD0_dif}
 \langle (1-\theta \mathbf{T}  \boldsymbol{\rho_{\theta,V_{\gamma}}^{(2)}})^{-1} (p - \mu_{p,\delta}) , (1-\theta \mathbf{T}  \boldsymbol{\rho_{\theta,V_{\delta}}^{(2)}})^{-1} (p - \mu_{p,\delta}) - (1-\theta \mathbf{T}  \boldsymbol{\rho_{\theta,V_{0}}^{(2)}})^{-1} (p - \mu_{p,\delta}) \rangle_{\rho_{\theta,V_{\delta}}^{(1)}} \\
 \leq C \delta |\log \delta|^C \cdot  \| p - \mu_{p,\delta} \|_{\mathcal{H}}^2 \leq  C^2 \delta |\log \delta|^C.
\end{multline}
To obtain the first inequality in \eqref{eqn:DdeltaD0_dif}, we used \eqref{eqn:DdD0} and the fact that $\| (1-\theta \mathbf{T}  \boldsymbol{\rho_{\theta,V_{\delta}}^{(2)}})^{-1}\|_{\mathcal{H}}  \leq C$ for all $\delta \in [0, c]$ (by Lemma \ref{lem:dress_existence}). To obtain the last inequality, we used that $|\mu_{p,\delta}| \leq C$ (by \eqref{eqn:H0Hd}) and hence $\|p - \mu_{p,\delta} \|_{\mathcal{H}}^2 \leq C$ for all $\delta \in [0, c]$. 

Now, by \eqref{eqn:DdeltaD0_dif}, in order to show that 
\begin{align}
& \big| \langle (1-\theta \mathbf{T}  \boldsymbol{\rho_{\theta,V_{\delta}}^{(2)}})^{-1} (p - \mu_{p,\delta}) , (1-\theta \mathbf{T}  \boldsymbol{\rho_{\theta,V_{\delta}}^{(2)}})^{-1} (p - \mu_{p,\delta}) \rangle_{\rho_{\theta,V_{\delta}}^{(1)}} \notag \\
&\qquad \qquad - \langle (1-\theta \mathbf{T}  \boldsymbol{\rho_{\theta,V_{0}}^{(2)}})^{-1} (p - \mu_{p,0}) , (1-\theta \mathbf{T}  \boldsymbol{\rho_{\theta,V_{0}}^{(2)}})^{-1} (p - \mu_{p,0}) \rangle_{\varrho} \big| \le C \delta |\log \delta|^C,\label{eqn:rhsd0tobd}
\end{align}
 it suffices to bound the quantity, with $f = (1-\theta \mathbf{T}  \boldsymbol{\rho_{\theta,V_{0}}^{(2)}})^{-1} (p - \mu_{p,0}) $,
 \begin{equation}
 \big| \langle f, f \rangle_{\varrho} -   \langle f, f \rangle_{\rho_{\theta,V_{\delta}}^{(1)}}\big|  \leq C \delta |\log \delta|^C. \label{eqn:rho0rhod}
 \end{equation}
  Indeed, the bound \eqref{eqn:rhodrho0pw} together with the tail bounds $\varrho(x) \leq c^{-1} e^{-c x^2}$ and  $\rho_{\theta,V_{\delta}}^{(1)}(x) \leq c^{-1} e^{-c x^2} $ (here $c>0$ can be taken independent of $\delta$, and the former follows from Lemma \ref{lem:varrho_bd} and the latter from \eqref{eqn:H0Hd}) implies that $ |\mu_{p,0} - \mu_{p,\delta}|   \leq C \delta |\log \delta|^C $, and putting this together with \eqref{eqn:DdeltaD0_dif} gives 
  \begin{align*}
& \Bigl| \langle (1-\theta \mathbf{T}  \boldsymbol{\rho_{\theta,V_{\delta}}^{(2)}})^{-1} (p - \mu_{p,\delta}) , (1-\theta \mathbf{T}  \boldsymbol{\rho_{\theta,V_{\delta}}^{(2)}})^{-1} (p - \mu_{p,\delta}) \rangle_{\rho_{\theta,V_{\delta}}^{(1)}} \\
&\qquad \qquad - \langle (1-\theta \mathbf{T}  \boldsymbol{\rho_{\theta,V_{0}}^{(2)}})^{-1} (p - \mu_{p,0}) , (1-\theta \mathbf{T}  \boldsymbol{\rho_{\theta,V_{0}}^{(2)}})^{-1} (p - \mu_{p,0}) \rangle_{\rho_{\theta,V_{\delta}}^{(1)}} \Bigr|   \leq C \delta |\log \delta|^C .
\end{align*}
  Then, using bound \eqref{eqn:rho0rhod} with $f = (1-\theta \mathbf{T}  \boldsymbol{\rho_{\theta,V_{0}}^{(2)}})^{-1} (p - \mu_{p,0}) $ completes the proof of \eqref{eqn:rhsd0tobd}.

So, to show continuity of the right-hand side of \eqref{eqn:fe_var_ID}, it remains to show \eqref{eqn:rho0rhod}. Note that $f$ has polynomial growth in either case of \eqref{eqn:rho0rhod}. Indeed this holds when $f = \varsigma_0$ and, when $f = (1-\theta \mathbf{T}  \boldsymbol{\rho_{\theta,V_{0}}^{(2)}})^{-1} (p - \mu_{p,0}) $ it holds by Lemma \ref{lem:ve_tail} (applied with $f$ there equal to $p - \mu_{p,0} \varsigma_0$ here). Therefore,
\begin{align}
    \label{frhorho1} 
    \begin{aligned} 
&\int_{-\infty}^{\infty} |f(x)|^2 |\rho_{\theta,V_{\delta}}^{(1)} - \varrho(x)| dx \\
&\qquad \qquad \leq C \delta |\log \delta|^C \int_{-|\log \delta|}^{|\log \delta| } |f(x)|^2  dx + \int_{|x|>|\log \delta|} |f(x)|^2 \bigl( |\rho_{\theta,V_{\delta}}^{(1)}(x)| + \varrho(x)| \bigr) dx \\
&\qquad \qquad \leq 2C \delta |\log \delta|^{C} + c^{-1} \int_{|x|>|\log \delta|} |f(x)|^2 e^{-c x^2} dx \leq  3C \delta |\log \delta|^{C}.
\end{aligned} 
\end{align}
To obtain the first inequality above, we split the integral into two parts, and we used \eqref{eqn:rhodrho0pw} on the first term and the triangle inequality on the second. To obtain the second, we used the bound $|f(x)| \leq C x^C $ on the first term and the bound \eqref{eqn:H0Hd} together with Lemma \ref{lem:varrho_bd} on the second. To obtain the third, we again used the polynomial growth of $f$. Since \eqref{frhorho1} implies \eqref{eqn:rho0rhod}, this shows continuity of the right-hand side of \eqref{eqn:fe_var_ID} at $\delta = 0$ (assuming \eqref{eqn:DdD0} and \eqref{eqn:rhodrho0pw}).

 Now, letting $\delta$ tend to $0$ in \eqref{eqn:fe_var_ID} establishes this identity at $\delta = 0$. Together with the Item \ref{item:D3} of Lemma \ref{thm:mm}, this completes the proof of \eqref{eqn:var_explicit1}. The proofs of \eqref{eqn:s0p_var} and \eqref{eqn:s0_var} are entirely analogous (and therefore omitted), as they are deduced in a very similar way from Item \ref{item:D5} of Lemma \ref{thm:mm}, together with \eqref{eqn:cov_0nthm} and \eqref{eqn:cov_00thm}, respectively.
\end{proof}

	\section{Quasi-particles in infinite volume} 
	
	\label{ParticlesInfinite} 

    Recall that quasi-particles were defined in Section \ref{FluctuationsQ} (see Definition \ref{ass:NT_assumption}) for the Toda lattice on a long, but finite, interval $\llbracket N_1, N_2 \rrbracket$. In this section we first define these quasi-particles directly on the infinite line $\mathbb{Z}$ (the analog of this for the box-ball system was done in \cite{SDBS}). We then use Theorem \ref{thm:quasi_fluct_intro} to deduce their fluctuation scaling limit. Throughout, we fix real numbers $\beta,\theta > 0$. 
	
	\subsection{Quasi-particles on $\mathbb{Z}$} 
	
	\label{ParticlesInfinite0} 
	
	Let $(\mathbf{a}(t); \mathbf{b}(t))$ denote the Flaschka variables of the Toda lattice on $\mathbb{Z}$, initialized under thermal equilibrium $\mu_{\beta,\theta;\infty}$ (Definition \ref{def:inf_thermal_eq}), where $\mathbf{a}(t) = (a_i(t))$ and $\mathbf{b}(t) = (b_i(t))$. Also, let $(\mathbf{p} (t); \mathbf{q} (t))$ denote the associated state space dynamics (recall Section \ref{subsec:model_and_objs}), and let $\mathbf{L} (t) = [L_{ij} (t)]$ denote the associated Lax matrix as in \eqref{eqn:Flaschka_vars_inf} (with the $(N_1,N_2)$ there equal to $(-\infty,\infty)$ here), which is an operator on the Hilbert space $\ell^2 (\mathbb{Z}) = \{ (v_i)_{i \in \mathbb{Z}} \in \mathbb{R}^{\mathbb{Z}}: \sum_{i \in \mathbb{Z}} v_i^2 < \infty \}$. By \cite[Remark 1.10]{OINL}, almost surely each eigenvalue of $\mathbf{L}(t)$ has multiplicity one, for all $t \ge 0$. Let $\eig \mathbf{L}(t)$ denote the (countable) set of eigenvalues of $\mathbf{L}(t)$. For any $\lambda \in \eig \mathbf{L}(t)$, let $\mathbf{u}^{\lambda} (t)= (u^{\lambda}(i;t)) \in \ell^2 (\mathbb{Z})$ denote a unit eigenvector of $\mathbf{L}(t)$ with eigenvalue $\lambda$.

    The following lemma is shown in Section \ref{ParitclesN}. Its first part is an analog of Lemma \ref{lem:lax_eig} for the Toda lattice $(\mathbf{a}(t); \mathbf{b}(t))$ on $\mathbb{Z}$; its second will be used to define quasi-particles for this model.
    
	\begin{lem} 
		
		\label{lt010}
		
		The following two statements almost surely hold, for any real number $t \ge 0$. 
		
		\begin{enumerate} 
			\item The operator $\mathbf{L}(t)$ has pure point spectrum, and $\eig \mathbf{L}(t) = \eig \mathbf{L}(0)$.  
			\item For any $\lambda \in \eig \mathbf{L}(t)$, we have absolute convergence of the sum
			\begin{flalign}
			    \label{qlambda0t} 
				\mathfrak{Q}^{\lambda} (t) = \displaystyle\sum_{i \in \mathbb{Z}} u^{\lambda} (i;t)^2 \cdot q_i (t).
			\end{flalign} 
			
		 	\noindent Moreover, the set $\{ \mathfrak{Q}^{\nu} (0) : \nu \in \eig \mathbf{L}(0) \} \subset \mathbb{R}$ has no accumulation points. 
		 \end{enumerate} 
	\end{lem} 

    As indicated by the following definition, the $\mathfrak{Q}^{\lambda} (t)$ from \eqref{qlambda0t} will be the quasi-particles for the Toda lattice on $\mathbb{Z}$. The bijection $\Phi$ below orders the quasi-particles from their starting locations; so, it may be viewed as the analog (for the infinite model) of the function $\varphi_0^{-1}$ from Definition \ref{def:quasi_intro}. 
    
    Observe that, while the quasi-particles for the Toda lattice on a finite interval (see \eqref{eqn:quasi_part_intro}) are permutations of the original Toda particle locations $(q_i(t))$, the $\mathfrak{Q}^{\lambda} (t)$ are sums of these locations, weighted by the eigenvector entries. We will see that these entries are localized on only a few coordinates, so the $\mathfrak{Q}^{\lambda} (t)$ will closely approximate the $(q_i(t))$ (see \eqref{qtin3} below).

	\begin{definition} 
		
		\label{qlambdanlimit0} 
		
		Recalling $\mathfrak{Q}^{\lambda} (t)$ from \eqref{qlambda0t}, let $\Phi : \mathbb{Z} \rightarrow \eig \mathbf{L}(0)$ be a bijection so that $\mathfrak{Q}^{\Phi (i)} (0) \le \mathfrak{Q}^{\Phi (j)} (0)$ if $i \le j$ and $\mathfrak{Q}^{\Phi (0)} \ge 0 > \mathfrak{Q}^{\Phi (-1)}$; such a bijection exists by the second part of Lemma \ref{lt010}. For any $i \in \mathbb{Z}$, set $\Lambda_i = \Phi (i)$, and define the \emph{quasi-particle} $\mathfrak{Q}_i (t) = \mathfrak{Q}^{\Lambda_i} (t)$ for any $t \in \mathbb{R}_{\ge 0}$.
		
	\end{definition} 

    The following result is the counterpart of Theorem \ref{thm:quasi_fluct_intro} on the infinite line, stating that the scaling limit of the above quasi-particles is given by the dressed L\'{e}vy--Chentsov field from Definition \ref{def:Z_intro}. Its proof is given in Section \ref{ParitclesN}, as a consequence of Theorem \ref{thm:quasi_fluct_intro} with a truncation argument.  
	
	\begin{cor}
		
		\label{q200} 
		
		For any $\beta>0$, there exists a constant $\theta_0 (\beta)>0$ such that the following holds whenever $\theta < \theta_0 (\beta)$. Given a real number $T \ge 1$, denote the quasi-particle fluctuations by
		\begin{equation*}
			Z_{i}^{\mathfrak{Q}} (t) \coloneqq  T^{-1/2} \cdot ( \mathfrak{Q}_{i}(t) - \mathfrak{Q}_i(0)- t \ve(\Lambda_i)), \qquad \text{for each $i \in \mathbb{Z}$ and $t \ge 0$}. 
		\end{equation*}
		
		\noindent Then there exists a constant $\mathfrak{c}>0$ and a coupling between $\mathbf{L}(0)$ and $\mathcal{Z}$ such that the following holds with probability at least $1-\mathfrak{c}^{-1}e^{-\mathfrak{c} (\log T)^2}$. For any integer $k \in \llbracket -T \log T, T \log T \rrbracket$ and real number $t \in [0, T \log T]$, denoting $\kappa = kT^{-1}$ and $\tau = tT^{-1}$, we have
		\begin{equation}
			\label{zqktz001} 
			\big| Z_{k}^{\mathfrak{Q}}(t) - \mathcal{Z}(\Lambda_k, \alpha \kappa + \tau \ve(\Lambda_k), \tau) \big| \leq T^{-\mathfrak{c}}.
		\end{equation}
		
	\end{cor}
	
	The following corollary provides the analog of Corollary \ref{zqi2} on the infinite line. Given Corollary \ref{q200}, its proof is essentially a consequence of that of Corollary \ref{zqi2} given Theorem \ref{thm:quasi_fluct_intro}. 
	
	\begin{cor}
		
		\label{zqi200}
		
		For any $\beta>0$, there exists a constant $\theta_0 (\beta)>0$ such that the following holds whenever $\theta < \theta_0 (\beta)$. For any real number $\delta > 0$, there is a constant $\mathfrak{c} = \mathfrak{c}(\delta) > 0$ such that the below holds with probability at least $1 - \mathfrak{c}^{-1} e^{-\mathfrak{c} (\log T)^2}$. For any integers $k, k' \in \llbracket -T \log T, T \log T \rrbracket$ satisfying $|k-k'| \le T^{1-\delta}$ and $|\Lambda_k - \Lambda_{k'}| \le T^{-\delta}$, we have that
		\begin{flalign*}
			\displaystyle\sup_{t \in [0, T \log N]} |Z_k^{\mathfrak{Q}} (t) - Z_{k'}^{\mathfrak{Q}} (t)| \le T^{-\mathfrak{c}}.
		\end{flalign*}
		
	\end{cor}
	
	\begin{proof}
		This follows from Corollary \ref{q200} and \eqref{zlambdak00estimate}. 
	\end{proof}

	\subsection{Approximations by finite Lax matrices and quasi-particles}
	
	\label{ParitclesN} 

    In this section we show Lemma \ref{lt010} and Corollary \ref{q200}, by approximating the Lax operator $\mathbf{L}(t)$ by finite-dimensional ones. The analysis here is reminiscent of that in \cite{FVAM} studying the Anderson model through its truncations. 
	
	For any integer $N \ge 1$, let $(\mathbf{a}_N(t); \mathbf{b}_N(t))$ denote the Flaschka variables of the Toda lattice on $\llbracket -N, N \rrbracket$, initialized under thermal equilibrium $\mu_{\beta,\theta;-N,N}$, where $\mathbf{a}_N(t) = ((a_{i;N} (t))$ and $\mathbf{b}_N (t) = (b_{i;N} (t))$. Also let $(\mathbf{p}_N (t); \mathbf{q}_N (t))$ denote the associated state space dynamics, where $\mathbf{p}_N (t) = (p_{i;N} (t))$ and $\mathbf{q}_N (t) = (q_{i;N}(t))$. Let $\mathbf{L}_N (t) = [L_{ij;N} (t)]$ denote the associated Lax matrix, and denote $\eig \mathbf{L}_N (t) = (\lambda_{1;N}, \ldots , \lambda_{2N+1;N})$, which does not depend on $t$ by Lemma \ref{lem:lax_eig}. For any $j \in \llbracket 1, 2N+1 \rrbracket$ let $\mathbf{u}_{j;N} (t) = (u_{j;N} (-N; t), \ldots , u_{j;N} (N; t))$ denote a unit eigenvector of $\mathbf{L}_N (t)$ with eigenvalue $\lambda_{j;N}$; set $u_{j;N} (k;t) = 0$ for any integer $k \notin \llbracket -N, N \rrbracket$. We couple the $(\mathbf{a}_N(t); \mathbf{b}_N(t))$ over all pairs $(-N,N)$ to share the same initial data with the Toda lattice $(\mathbf{a}(t); \mathbf{b}(t))$ on $\mathbb{Z}$ from Section \ref{ParticlesInfinite}, namely, $a_{i;N} (0) = a_{i}(0)$ for all $i \in \llbracket -N,N-1 \rrbracket$, and $b_{i;N} (0) = b_{i}(0)$ for all $i \in \llbracket -N,N \rrbracket$. 	

    The following lemma, shown in Section \ref{Prooflambdan} (using Lemma \ref{lem:Lax_eig_coupling}), states that eigenvalues and eigenvectors of a Lax matrix $\mathbf{L}_N (t)$ approximate those of larger Lax matrices $\mathbf{L}_{N'} (t)$. 
    
	\begin{lem}
		\label{nnq0}
		There exist constants $\mathfrak{c}_0, \mathfrak{c}_1 >0$ such that, for any integer $N \ge 1$, the following holds with probability at least $1 - \mathfrak{c}_1^{-1} e^{-\mathfrak{c}_1 (\log N)^2}$. Let $N' \ge 2N$ be any integer and $t \in [0, N(\log N)^{-3}]$ be any real number. Set $\zeta = e^{-150(\log N)^{3/2}}$, and let $\lambda_{i;N} \in \eig \mathbf{L}_N (t)$ denote an eigenvalue of $\mathbf{L}_N (t)$ that has a $\zeta$-localization center $\varphi \in \llbracket -N/2, N/2 \rrbracket$ with respect to $\mathbf{L}_N (t)$. There exists an index $i' = \psi_{N;N'} (i) \in \llbracket 1,2N'+1 \rrbracket$ such that the below three statements hold.
		\begin{enumerate}
			\item We have  $|\lambda_{i;N} - \lambda_{i';N'}'| \le e^{-\mathfrak{c}_0 (\log N)^3}$.
			\item For any integer $k \in \llbracket -N', N' \rrbracket$, we have  $|u_{i;N} (k;t) - u_{i';N'} (k;t)| \le e^{-\mathfrak{c}_0 (\log N)^2}$. In particular, $\varphi$ is a $e^{-150(\log N')^{3/2}}$ localization center for $\lambda_{i';N'}$ with respect to $\mathbf{L}_{N'}(t)$.
			\item For any index $i'' \in \llbracket 1, 2N'+1 \rrbracket$ with $i'' \ne i'$, we have  $|\lambda_{i;N} - \lambda_{i'';N'}'| > e^{-\mathfrak{c}_0 (\log N)^2/2}$.
		\end{enumerate}
	\end{lem}

		By Lemma \ref{nnq0} and the Borel--Cantelli lemma, there almost surely exists a constant $\mathfrak{c}_0>0$ and an integer $N_0 \ge 1$ such that the event in Lemma \ref{nnq0} holds for all $N \ge N_0$. By Lemma \ref{lem:sep_lem}, we may also assume that $\mathsf{SEP}_{\mathsf{L}_N (0)} (e^{-\mathfrak{c}_0 (\log N)^2/10})$ holds for all $N \ge N_0$. By Lemmas \ref{lem:fin_inf_t_coupling}, \ref{lem:q_spacing}, and \ref{lem:inf_vol_bds}, given a real number $T \ge 1$, we may additionally assume for all $N \ge N_0$ (if $N_0 \ge T^4$) that 
		\begin{flalign}
			\label{qabnt00} 
			\begin{aligned} 
			& \displaystyle\sup_{t \in [0,T]} \displaystyle\inf_{|i| \le N - N^{1/2}} \big( |a_{i;N} (t) - a_i (t)| + |b_{i;N} (t) - b_i (t)| \big) \le e^{-(\log N)^2}; \\ 
			& \displaystyle\sup_{t \in [0,T]} \displaystyle\inf_{|i| \le N - N^{1/2}}  |q_{i;N} (t) - q_i (t)| \le e^{-(\log N)^2}; \qquad \displaystyle\max_{|i| \le N} |i|^{-1/2} |q_{i;N} (0) - \alpha i| \le \log N; \\ 
			& \displaystyle\max_{|i|,|j| \le 9N/10} |i-j|^{-1/2} \cdot |q_{i;N} (t) - q_{j;N} (t) - \alpha (i-j)| \le (\log N)^2; \\
			& \displaystyle\sup_{t \in [0,T]} \displaystyle\max_{|i| \le  N} (|a_i (t)| + |b_i (t)| ) \le 6 \log N.
			\end{aligned} 
		\end{flalign} 
		
		\noindent By Lemma \ref{lem:loc_cent_dif}, we may further assume for all $N \ge N_0$ that if, for some real number $t \in [0, N (\log N)^{-6}]$ and integers $\varphi \in \llbracket t (\log N)^4 - N, N - t (\log N)^4 \rrbracket$ and $\varphi' \in \mathbb{Z}$, we have 
		\begin{flalign}
			\label{ncenter0} 
			 \quad \min \{ |u_{i;N} (\varphi;t)|, |u_{i,N} (\varphi';t)| \} \ge e^{-150 (\log N)^{3/2}}, \quad \text{then} \quad |\varphi - \varphi'| \le 2 (\log N)^3.
		\end{flalign}
		
		\noindent By the probability bounds in Lemmas \ref{nnq0}, \ref{lem:sep_lem}, \ref{lem:fin_inf_t_coupling}, \ref{lem:q_spacing}, \ref{lem:inf_vol_bds}, and \ref{lem:loc_cent_dif}, we have 
		\begin{flalign}
			\label{kprobabilityn0} 
			\mathbb{P}[N_0 \le K] \ge 1 - \mathfrak{c}_0^{-1} e^{-\mathfrak{c}_1 (\log K)^2},
		\end{flalign}
		
		\noindent for some constant $\mathfrak{c}_1 > 0$ and any real number $K \ge T^{10}$. 

        The next corollary, shown in Section \ref{Prooflambdan}, indicates that eigenvalues of $\mathbf{L}(t)$ can be obtained from those of $\mathbf{L}_N (t)$ by a limiting procedure. Below, we recall the functions $\psi_{N;N'}$ from Lemma \ref{nnq0}. 
		
		\begin{cor}
			\label{psin}
			Let $N \ge N_0$ be an integer, let $t \ge 0$ be a real number, and let $i \in \llbracket 1, 2N+1 \rrbracket$ be an index such that $\lambda_{i;N}$ admits an $e^{-150(\log N)^{3/2}}$-localization center $\varphi \in \llbracket -N/2, N/2 \rrbracket$ with respect to $\mathbf{L}_N (t)$. Then the following statements hold.
			
			\begin{enumerate} 
				\item For any integers $N' \ge 2N$, and $N'' \ge 2N'$, we have that $\psi_{N';N''} (\psi_{N;N'} (i)) = \psi_{N;N''} (i)$. 
				\item For any integer $k \ge 0$, the limits
			\begin{flalign}
				\label{lambdann00} 
				\lambda_{i;N;\infty} = \displaystyle\lim_{N' \rightarrow \infty} \lambda_{\psi_{N;N'} (i); N'}; \qquad u_{i;N;\infty} (k;t) = \displaystyle\lim_{N' \rightarrow \infty} u_{\psi_{N;N'} (i); N'} (k;t),
				\end{flalign}
			
			\noindent exist and satisfy
			\begin{flalign}
				\label{lambdann002} 
				|\lambda_{i;N;\infty} - \lambda_{i;N}| \le 2\mathfrak{c}_0^{-1} e^{-\mathfrak{c}_0 (\log N)^2/2}; \quad \displaystyle\max_{k \in \mathbb{Z}}  |u_{i;N;\infty} (k;t) - u_{i;N} (k;t)| \le e^{-(\log N)^{3/2}}.
			\end{flalign}
			
			\item For any integers $N' \ge 2N$ and $k \in \mathbb{Z}$, we have that $u_{\psi_{N;N'}(i);N';\infty} (k;t) = u_{i;N;\infty} (k;t)$ and $\lambda_{\psi_{N;N'} (i); N';\infty} = \lambda_{i;N;\infty}$. Moreover, if $|k - \varphi| >  (\log N)^3$, then  
			\begin{flalign} 
				\label{ukin0} 
				|u_{i;N;\infty} (k,t)| \le 2e^{-(\max \{ \log |k - \varphi|, \log N \})^{3/2}}. 
			\end{flalign} 
			
			\item We have that $\lambda_{i;N;\infty} \in \eig \mathbf{L}(t)$, and $\mathbf{u}_{i;N;\infty} (t) = (u_{i;N;\infty} (k;t))_{k \in \mathbb{Z}}$ is a unit eigenvector of $\mathbf{L}(t)$ with eigenvalue $\lambda_{i;N;\infty}$.
			
			\end{enumerate} 
					
		\end{cor}

        The next corollary, shown in Section \ref{Prooflambdan}, bounds the smallest dimension of a Lax matrix $\mathbf{L}_N (t)$ associated with a given eigenvalue of $\mathbf{L}(t)$, in terms of its localization center.

		\begin{cor} 
			
			\label{lambdalambda2}

			 Fix $\lambda \in \eig \mathbf{L}(t)$, and let $N \ge N_0$ is the minimal integer (if it exists) for which $\lambda = \lambda_{i;N;\infty}$, for some index $i \in \llbracket 1, 2N+1 \rrbracket$. If $\varphi \in \llbracket -N, N \rrbracket$ is an $e^{-150(\log N)^{3/2}}$-localization center for $\lambda_{i;N}$ with respect to $\mathbf{L}_N (t)$, then $N \le 10 \max \{ N_0, \varphi \}$. 
		\end{cor} 
		
		\begin{proof}[Proof of Lemma \ref{lt010}]
			
			For any integer $N \ge 2N_0$, let $\mathfrak{S}_N = (\lambda_{i;N;\infty}; \mathbf{u}_{i;N;\infty} (t))_i$, where $i \in \llbracket 1, 2N+1 \rrbracket$ ranges over all indices such that $\lambda_{i;N}$ admits an $e^{-150(\log N)^{3/2}}$-localization center $\varphi \in \llbracket -N/2, N/2 \rrbracket$ with respect to $\mathbf{L}_N(t)$.	By the third part of Corollary \ref{psin}, the family $((\mathfrak{S}_N)_N, (\psi_{N,N'})_{N' \ge 2N})$ forms an inverse system and therefore admits an inverse limit, which we denote by $\mathfrak{S}_{\infty}$. 
			
			By the third part of Corollary \ref{psin}, for any $(\lambda, \mathbf{u}^{\lambda} (t))$, we have that $\lambda$ is an eigenvalue of $\mathbf{L}(t)$ with unit eigenvector $\mathbf{u}^{\lambda} (t)$. Moreover, the $(\mathbf{u}^{\lambda}(t))$ comprise an orthonormal eigenbasis of $\ell^2 (\mathbb{Z})$. Indeed, they are orthonormal since the $(\mathbf{u}_{i;N} (t))$ are. To verify that they form a basis, it suffices to show for any $k \in \mathbb{Z}$ that $\sum_{\lambda} u^{\lambda} (k;t)^2 = 1$; since $\sum_{i=1}^{2N+1} u_{i;N} (k;t)^2 = 1$ for any $N \ge 1$, this is a quick consequence of \eqref{lambdann002} and \eqref{ukin0}. Thus, $\mathbf{L}(t)$ has pure point spectrum for all $t \ge 0$ and, for any eigenvalue $\lambda \in \eig \mathbf{L}(t)$, there exists an integer $N(\lambda) \ge N_0$ such that the following holds. For all $N \ge N(\lambda)$, there is an index $i$ satisfying the assumptions of Corollary \ref{psin} such that \eqref{lambdann002} holds. This, together with the fact that $\eig \mathbf{L}_N (0) = \eig \mathbf{L}_N (t)$, implies that $\eig \mathbf{L}(0) = \mathbf{L}(t)$ for all $t \ge 0$. This confirms the first statement of the lemma. The second follows from the decay of $|u_{i;N;\infty} (k,t)|$ from \eqref{ukin0}, with the approximation for $(q_i(t))$ given by the second and third bounds in \eqref{qabnt00}. 
            
            To verify the third, fix $\lambda \in \eig \mathbf{L}(t)$ and let $N \ge N_0$ denote an integer such that there exists an index $i \in \llbracket 1, 2N+1 \rrbracket$ for which $\lambda_{i;N;\infty} = \lambda$. Then $\lambda_{i;N}$ admits an $e^{-150(\log N)^{3/2}}$-localization center $\varphi \in \llbracket -N/2, N/2 \rrbracket$ with respect to $\mathbf{L}_N (t)$ in the interval. Denoting 
			\begin{flalign*}
				\mathfrak{Q}_{i;N} (t) = \displaystyle\sum_{k=-N}^{N} u_{i;N} (k;t)^2 \cdot q_{k;N} (t),
			\end{flalign*}
			
			\noindent it quickly follows from \eqref{ukin0}, the second bound in \eqref{lambdann002}, and the second and third bounds in \eqref{qabnt00} that 
			\begin{flalign} 
			\label{qtin2} 
			|\mathfrak{Q}_{i;N} (t) - \mathfrak{Q}^{\lambda} (t)| \le e^{-c(\log N)^{3/2}}.
			\end{flalign} 
			
			 Additionally, we have 
			\begin{flalign*}
				|\mathfrak{Q}_{i;N} (t) - q_{\varphi;N} (t)| & \le \displaystyle\sum_{k=-\lfloor 2(\log N)^3 \rfloor}^{\lfloor 2(\log N)^3 \rfloor} u_{i;N} (\varphi+k;t)^2 \cdot | q_{\varphi+k;N} (t) - q_{\varphi;N} (t) | \\
				& \qquad + 2e^{-150(\log N)^{3/2}} \displaystyle\max_{j \in \llbracket -N, N \rrbracket} |q_{j;N} (t)| \le (\log N)^4
			\end{flalign*}
	
			\noindent where in the first inequality we used \eqref{ncenter0} and the fact that $\sum_{k=-N}^N u_{i;N} (k;t)^2 = 1$; in the second we used the last three bounds in \eqref{qabnt00} (and that $\varphi \in \llbracket -N/2, N/2 \rrbracket$). With \eqref{qtin2}, this yields
			\begin{flalign}
			\label{qtin3} 
			|\mathfrak{Q}^{\lambda} (t) - q_{\varphi;N} (t)| \le (\log N)^5.
			\end{flalign}

			\noindent Moreover, by the last three bounds in \eqref{qabnt00}, we have that 
			\begin{flalign}
				\label{qnt0t2}  
				|q_{\varphi;N} (t) - \alpha \varphi | \le |q_{\varphi;N} (0) - \alpha \varphi| + |q_{\varphi;N} (0) - q_{\varphi;N} (t)| \le (|\varphi|^{1/2} + t) (\log N)^2.
			\end{flalign}  
			
			Now let $\mathfrak{a} \ge 0$ be a real number. By \eqref{qtin3} and \eqref{qnt0t2}, if $|\mathfrak{Q}^{\lambda}(t)| \le \mathfrak{a}$, then $|\varphi| \le C(|\mathfrak{a}|+t+1) (\log N)^5$. By Corollary \ref{lambdalambda2}, if $N$ is the minimal integer for which $\lambda = \lambda_{i;N;\infty}$ for some $i \in \llbracket 1, 2N+1 \rrbracket$, then $N \le 10 \max \{ N_0, \varphi \} \le C \max \{ N_0, (|\mathfrak{a}|+t+1) (\log N)^5 \}$. Since $N$ satisfies this inequality only if $N \le (|\mathfrak{a}| + t + N_0) (\log N)^6$ or $N \le CN_0$, there are only finitely many $\mathfrak{Q}^{\lambda}(t) \in [-\mathfrak{a}, \mathfrak{a}]$. So, $(\mathfrak{Q}^{\lambda}(t))$ has no accumulation points, as in the third statement of the lemma.
		\end{proof}

		\begin{proof}[Proof of Lemma \ref{q200}]

			We will deduce Lemma \ref{q200} from Theorem \ref{thm:quasi_fluct_intro}. By \eqref{kprobabilityn0}, we may assume $N_0 \le T^{20}$. Let $\mathsf{T} \in \llbracket T (\log T)^{20}, 2T(\log T)^{20} \rrbracket$ be an odd integer, and let $\mathsf{N} = (\mathsf{T}^{100} - 1)/2 \ge 10N_0$. Set $\zeta = e^{-50 (\log \mathsf{N})^{3/2}}$, and let $\varphi_t : \llbracket 1, 2\mathsf{N}+1 \rrbracket \rightarrow \llbracket - \mathsf{N}, \mathsf{N} \rrbracket$ denote a $\zeta$-localization center bijection for $\mathbf{L}_{\mathsf{N}} (t)$, for each real number $t \ge 0$. Then Assumption \ref{ass:NT_assumption} is satisfied, with the $(T,N)$ there equal to $(\mathsf{T}, 2\mathsf{N}+1)$ here. Following Definition \ref{def:quasi_intro}, set $Q_{j;\mathsf{N}} (t) = q_{\varphi_t (j);\mathsf{N}} (t)$ for each $j \in \llbracket 1, 2\mathsf{N}+1 \rrbracket$ and $t \ge 0$. As in \eqref{eqn:quasi_part_fluct_intro}, for any $i \in \llbracket 1, 2\mathsf{N}+1 \rrbracket$, define 
			\begin{flalign*} 
				Z_{\varphi_0 (i);\mathsf{N}}^{\mathcal{Q}} (t) = T^{-1/2} \cdot (Q_{i;\mathsf{N}}(t) - Q_{i;\mathsf{N}} (0) - t \ve (\lambda_{i;N})).
			\end{flalign*} 
			
			\noindent  Throughout this proof, we let $t = T \tau \in [0, T \log T]$ be a real number and recall the notation in the proof of Lemma \ref{lt010}, with the $N$ there equal to $\mathsf{N}$ here. We further restrict to the event on which Theorem \ref{thm:quasi_fluct_intro}, with the $(T,N_1,N_2)$ there equal to $(\mathsf{T}, -\mathsf{N}, \mathsf{N})$ here, applies. 
			
			Ranging $\varphi$ in \eqref{qtin3} and \eqref{qnt0t2} over $\llbracket -T (\log T)^9, T(\log T)^9 \rrbracket$, and using the fourth bound in \eqref{qabnt00}, we find that there are at least $T \log T$ distinct values of $\mathfrak{Q}^{\lambda} (0)$, both in the interval $[-C T (\log T)^9, 0]$ and in the interval $[0, CT (\log T)^9]$. Hence, for $|k| \le T \log T$, we have that $|\mathfrak{Q}_k (0)| \le C T (\log T)^9$. Moreover, as shown at the end of the proof of Lemma \ref{lt010} for $\mathfrak{a} = CT (\log T)^9$, for any $\lambda \in \eig \mathbf{L}(0)$ for which $|\mathfrak{Q}^{\lambda} (0)| \le CT (\log T)^9$, we have since $\mathsf{N} \ge T^{50} > (|\mathfrak{a}| + t + N_0) (\log N)^6$ that $\lambda = \lambda_{i(\lambda);\mathsf{N};\infty}$ for some $i(\lambda) \in \llbracket 1, 2 \mathsf{N} + 1 \rrbracket$. Hence, if $\lambda = \Phi(k)$ (recall Definition \ref{qlambdanlimit0}) for some $k \in \llbracket -T \log T, T \log T \rrbracket$, then $\Phi (k) = \lambda_{i;\mathsf{N};\infty}$, where we have abbreviated $i = i_k = i (\Phi(k))$. Then, by \eqref{lambdann002} and \eqref{qtin3} (and the fact that $\mathfrak{Q}^{\lambda}(t) = \mathfrak{Q}_k (t)$), we have 
			\begin{flalign}
				\label{lambda2q2k00} 
				|\lambda - \lambda_{i;\mathsf{N}}| \le c^{-1} e^{-c(\log T)^2}; \qquad 	|\mathfrak{Q}_k (t) - Q_{i;\mathsf{N}} (t)| \le (\log T)^6.
			\end{flalign} 
			
			Moreover, applying \eqref{zkqtestimate00}, with the  $(k;T,N)$ there equal to $(\varphi_0 (i); \mathsf{T}, \mathsf{N})$ here (using the fact that $|\varphi_0 (i)| \le T (\log T)^{10} \le \mathsf{T}$ by \eqref{lambda2q2k00}, the fact that $|\mathfrak{Q}_k (0)| \le CT (\log T)^9$, and \eqref{qabnt00}), we obtain the existence of a coupling  between $\mathbf{L}(0)$ and $\mathcal{Z}$ such that 
			\begin{flalign}
			    \label{zqtz} 
				\big| Z_{\varphi_0 (i)}^{\mathcal{Q}}(t) - \mathcal{Z}(\lambda_{i;\mathsf{N}}, \alpha T^{-1} \varphi_0 (i) + \tau \ve(\lambda_{i;\mathsf{N}}), \tau) \big| \leq T^{-\mathfrak{c}}.
			\end{flalign}
			
			\noindent Here, we also used the fact that the law of $(\mathsf{T} T^{-1})^{1/2} \cdot \mathcal{Z}(\lambda_{i;\mathsf{N}}, \alpha \mathsf{T}^{-1} \varphi_0 (i) + t \mathsf{T}^{-1} \ve(\lambda_{i;\mathsf{N}}), t \mathsf{T}^{-1})$ (which is what would arise from applying \eqref{zkqtestimate00}, with the $(T,N)$ there equal to $(\mathsf{T}, \mathsf{N})$ here) coincides with that of $\mathcal{Z}(\lambda_{i;\mathsf{N}}, \alpha T^{-1} \varphi_0 (i) + \tau \ve(\lambda_{i;\mathsf{N}}), \tau)$, jointly as processes in the parameters $(i, \lambda_{i;\mathsf{N}}, \tau)$ (due to the general fact that, for any $r>0$, the law of the process $r^{-1/2} \cdot \mathcal{Z} (\lambda, r\kappa, r\tau)$ coincides with that of $\mathcal{Z} (\lambda, \kappa, \tau)$, jointly in the parameters $(\lambda,\kappa,\tau)$, as a quick consequence of \eqref{eqn:intro_Z}, \eqref{eqn:Xdef}, \eqref{eqn:logphi_intro}, \eqref{wf00}, and Definition \ref{def:white_noise_dress}). Combining \eqref{zqtz} with \eqref{lambda2q2k00}, Definition \ref{def:Z_intro}, and Lemmas \ref{lem:ve_tail}, \ref{lem:holder_cont}, and \ref{lem:bd_lem}, this yields with probability at least $1 - c^{-1} e^{-c(\log T)^2}$ that (recalling $\Lambda_k = \Phi(k) = \lambda$) 
			\begin{flalign}
				\label{zqktz00} 
				\big| Z_{k}^{\mathfrak{Q}}(t) - \mathcal{Z}(\Lambda_k, \alpha T^{-1} \varphi_0 (i) + \tau \ve(\Lambda_k), \tau) \big| \leq T^{-\mathfrak{c}}.
			\end{flalign}
			
			\noindent It remains to replace $T^{-1} \varphi_0 (i)$ with $\kappa$ in \eqref{zqktz00}. To do so, observe from \eqref{qtin3} and \eqref{qnt0t2} that $|\mathfrak{Q}_k (0) - \alpha \varphi_0 (i_k)| \le T^{1/2} (\log T)^8$. Since the $(\mathfrak{Q}_k (0))$ are increasing and $\mathfrak{Q}_0 (0) \ge 0 > \mathfrak{Q}_{-1} (0)$, this implies that $|k - \varphi_0 (i_k)| \le T^{1/2} (\log T)^C$ and thus that $|T^{-1} \varphi_0 (i) - \kappa| \le T^{-c}$. This, together with \eqref{zqktz00}, yields \eqref{zqktz001} (again using  Definition \ref{def:Z_intro}, and Lemmas \ref{lem:ve_tail}, \ref{lem:holder_cont}, and \ref{lem:bd_lem}).
		\end{proof}

	\subsection{Proofs of Lemma \ref{nnq0}, Corollary \ref{lambdann00}, and Corollary \ref{lambdalambda2}} 
	
	\label{Prooflambdan} 
	
	In this section we prove first Corollary \ref{lambdann00}, then Corollary \ref{lambdalambda2}, and next Lemma \ref{nnq0}. 
	
	\begin{proof}[Proof of Corollary \ref{lambdann00}] 
		
		Applying the first statement of Lemma \ref{nnq0} twice (with the $(i,N,N')$ there equal to $(i,N,N')$ and $(\psi_{N';N} (i),N',N'')$ here) yields $|\lambda_{\psi_{N'';N'} (\psi_{N';N}(i)); N''} - \lambda_{i;N}|  \le 2e^{-\mathfrak{c}_0 (\log N)^3}$. With the third part of Lemma \ref{nnq0}, this yields $\psi_{N'';N'} (\psi_{N';N}(i)) = \psi_{N'';N} (i)$, as in the first part of the lemma. Since $(\lambda_{\psi_{N;N'} (i); N'})_{N' \ge 2N}$ and $(u_{\psi_{N;N'} (i); N'} (k;t))_{N' \ge 2N}$ form Cauchy sequences, by the first and second parts of Lemma \ref{nnq0}, the limits \eqref{lambdann00} exist; the latter also yield the rates of convergence \eqref{lambdann002}. 
		
		The first statement of the third part of the corollary follows from the first, together with \eqref{lambdann00} (applied with the $N$ there equal to both $N$ and $N'$ here). If $|k| > 2N$, then taking $N' = |k|-1$ gives $u_{i;N;\infty} (k;t) = u_{i';|k|-1;\infty} (k;t)$ for $i' = \psi_{N;|k|-1}(i)$; together with \eqref{lambdann002} (with the $N$ there equal to $|k|-1$ here) and the fact that $u_{i';|k|-1} (k;t) = 0$, this yields \eqref{ukin0}. If instead $|k| \le 2N$, then we have $|u_{i;N} (k,t)| \le e^{-150 (\log N)^{3/2}}$ by \eqref{ncenter0} (and the fact that $\varphi$ is an $e^{-150(\log N)^{3/2}}$-localization center for $\lambda_{i;N}$ with respect to $\mathbf{L}(t)$), which with \eqref{lambdann002} again verifies \eqref{ukin0}. This confirms the third part of the corollary. Furthermore, for any $k \in \mathbb{Z}$, we have by \eqref{qabnt00}, \eqref{lambdann002}, and \eqref{ukin0} that
		\begin{flalign*}
			( \mathbf{L}(t) \cdot \mathbf{u}_{i;N;\infty} (t))_k = \displaystyle\sum_{j \in \mathbb{Z}} L_{kj} (t) \cdot u_{i;N;\infty} (j;t) & = \displaystyle\lim_{N' \rightarrow \infty}  \displaystyle\sum_{j \in \mathbb{Z}} L_{kj;N'} (t) \cdot u_{i;N;N'} (j;t) \\
			& = \displaystyle\lim_{N' \rightarrow \infty} \lambda_{i;N;N'} \cdot u_{i;N;N'} (k,t) = \lambda_{i;N;\infty} \cdot  u_{i;N;\infty} (k;t),
		\end{flalign*}
		
		\noindent which shows the fourth part of the corollary.
	\end{proof}
	
	\begin{proof}[Proof of Corollary \ref{lambdalambda2}] 
		
		Assume to the contrary that $N > 10 \max \{ N_0, \varphi \}$. Let $\hat{\varphi}$ denote an $N^{-1}$-localization center for $\lambda_{i;N}$ with respect to $\mathbf{L}_N (t)$. Then $|\varphi - \hat{\varphi}| \le 2 (\log N)^3$ by \eqref{ncenter0}, so setting $\tilde{N} = \lfloor N/2 \rfloor$ we have $\tilde{N} \ge 4 \max \{ N_0, |\hat{\varphi}| \}$. Let $\mathcal{S}$ denote the set of indices $j \in \llbracket 1, 2\tilde{N}+1 \rrbracket$ such that $\lambda_{j;\tilde{N}}$ admits an $e^{-150(\log \tilde{N})^{3/2}}$-localization center in the interval $\llbracket -\tilde{N}/2, \tilde{N}/2 \rrbracket$. By \eqref{ncenter0} (and since the $(\mathbf{u}_{j;N}(t))$ form an orthonormal basis), we have $\sum_{j \in \mathcal{S}} u_{j;\tilde{N}} (\hat{\varphi};t)^2 \ge 1 - c^{-1} e^{-c(\log N)^{3/2}}$, which with \eqref{lambdann002} yields $\sum_{j \in \mathcal{S}} u_{\psi_{\tilde{N};N} (j);N} (\hat{\varphi};t)^2 > 1 - N^{-2}$. So, since $\hat{\varphi}$ is an $N^{-2}$-localization center of $\lambda_{i;N}$ with respect to $\mathbf{L}_N (t)$, we have $i = \psi_{\tilde{N};N} (j)$ for some $j \in \mathcal{S}$. Hence, $\lambda = \lambda_{i;N;\infty} = \lambda_{j;\tilde{N};\infty}$, where the last equality holds by the third statement of Lemma \ref{lambdann00}. This contradicts the minimality of $N$, so $N > 10 \max \{ N_0, \varphi \}$, establishing the corollary.
	\end{proof}

	\begin{proof}[Proof of Lemma \ref{nnq0} (Outline)] 
		
		We only outline the proof of this lemma, since it largely follows from results in \cite[Section 5]{Agg25a}. It suffices to show it for any integer $N' \in \llbracket 2N, 4N \rrbracket$, as then the lemma would follow by induction on $\lfloor \log_2 (N'/N) \rfloor$. Its first part follows from Lemma \ref{lem:Lax_eig_coupling}, applied with the $(\mathbf{L}; \tilde{\mathbf{L}}; \mathcal{D})$ there equal to $(\mathbf{L}_{N'} (0); \mathbf{L}_{N} (0); \llbracket -N', N' \rrbracket \setminus \llbracket -N, N \rrbracket)$ here. In what follows, we restrict to the event $\mathsf{E} = \mathsf{SEP}_{\mathbf{L}_N (0)} (e^{-\mathfrak{c}_0 (\log N)^2/4}) \cap \mathsf{SEP}_{\mathbf{L}_{N'}(0)} (e^{-\mathfrak{c}_0 (\log N)^2/4})$ (as we may, by Lemma \ref{lem:sep_lem}). Then, the third part of the lemma follows from the first. 
		
		It remains to verify the second part of the lemma. We will only do so for a fixed choice of $t \ge 0$, as then the proof that it holds for all $t$ simultaneously follows from a union bound over a lattice mesh, with a continuity estimate. We omit further details on the latter, as it is entirely analogous to the proofs of Lemmas \ref{lem:Xi_brownian_cont_allparam_pl} and \ref{lem:gen_brownian_cont_allparam_prelim} (using the fact that $|\partial_t u_{i;N} (k;t)| \le CN^2$ with probability at least $1-c^{-1} e^{-c(\log N)^2}$, by \cite[Equation (5.4)]{Agg25a} and Lemma \ref{lem:bd_lem}). 
		
		Define the matrix $\mathbf{M} = [M_{ij}]$, with rows and columns indexed by $i,j \in \llbracket -N', N' \rrbracket$, by setting $M_{ij} = L_{ij} (t)$ for all $i$ and $j$. Then, $\mathbf{M}$ has the same law as $\mathbf{L}_{N'} (0)$, by first part of Lemma \ref{lem:fin_inf_t_coupling}. By \cite[Proposition 4.4]{Agg25a}, the pairs $(\mathbf{M}; \mathbf{L}_N (t))$ and $(\mathbf{M}; \mathbf{L}_{N'} (t))$ satisfy \cite[Assumption 5.3]{Agg25a} for the $(\bm{L}; \tilde{\bm{L}})$ there, with the $(\mathcal{D}; \delta)$ there equal to $(\llbracket -N', N' \rrbracket \setminus \llbracket -9N/10, 9N/10 \rrbracket; e^{-N/50})$ here. For any $z \in \mathbb{C}$, set $z = \lambda_{i;N} +  \mathrm{i} \eta$ where $\eta = e^{-(\log N)^{5/2}}$ and denote the resolvents $\mathbf{G}_N (z) = (\mathbf{L}_N (t) - z)^{-1}$, $\mathbf{G}_{N'} (z) = (\mathbf{L}_{N'} (t) - z)^{-1}$, and $\bm{\mathfrak{G}} (z) = (\mathbf{M} - z)^{-1}$. Then, by \cite[Lemma 5.4]{Agg25a} (with Lemma \ref{lem:bd_lem}),  we have with probability at least $1 - c^{-1} e^{-c(\log N)^2}$ that
		\begin{flalign}
			\label{gkk2} 
			\begin{aligned} 
				\displaystyle\max_{k \in \llbracket -4N/5, 4N/5 \rrbracket} | & G_{kk;N} (z) - G_{kk;N'} (z)| \\
				&  \le \displaystyle\max_{k \in \llbracket -4N/5, 4N/5 \rrbracket} \big( |G_{kk;N} (z) - \mathfrak{G}_{kk} (z)| + |G_{kk;N'} (z) - \mathfrak{G}_{kk} (z)|  \big) \le  e^{-cN}.
			\end{aligned} 
		\end{flalign}
		
		\noindent Next observe that
		\begin{flalign*}
			\big| G_{kk;N} (& z) - G_{kk;N'} (z) - \eta^{-1} (  u_{i;N} (k;t)^2 -  u_{i';N'} (k;t)^2) \big| \\
			&  \le \bigg| \displaystyle\frac{\eta^{-1}(\lambda_{i';N'} - \lambda_{i;N})^2  }{(\lambda_{i';N'} - \lambda_{i;N})^2 + \eta^2} \bigg| \cdot u_{i';N'} (k;t)^2 +  \eta \displaystyle\sum_{j \ne i}  \displaystyle\frac{u_{j;N} (k;t)^2}{(\lambda_{j;N} - \lambda_{i;N})^2} + \displaystyle\sum_{j \ne i'}  \displaystyle\frac{u_{j;N'} (k;t)^2}{(\lambda_{j;N'} - \lambda_{i;N})^2}  \\
			& \le e^{-(\log N)^2},
		\end{flalign*} 
		
		\noindent where the first bound holds by \eqref{gijuij} and the second by the first part of the lemma, the definition of $\eta = e^{-(\log N)^{5/2}}$, and our restriction to $\mathsf{E}$. This, with \eqref{gkk2}, yields $|u_{i;N} (k;t)^2 - u_{i';N'} (k;t)^2| \le e^{-(\log N)^2}$ if $|k| \le 4N/5$. Since $|\varphi| \le N/2$ and $\varphi$ is an $e^{-150(\log N)^{3/2}}$-localization center for $\lambda_{i;N}$ with respect to $\mathbf{L}_N (t)$, this in particular implies that $ u_{i';N'} (\varphi;t)^2 \ge u_{i;N} (\varphi;t)^2 - e^{-(\log N)^2} \ge e^{-150 (\log N')^{3/2}}$, so $\varphi$ is an $e^{-150(\log N')^{3/2}}$-localization center for $\lambda_{i';N'}$ with respect to $\mathbf{L}_N (t)$. 
			
			It remains to show that $|u_{i;N} (k;t)^2 - u_{i';N'} (k;t)^2| \le e^{-c (\log N)^2}$ for $|k| > 4N/5$. In this case, $|k-\varphi| > N/4 > T (\log N)^2$. Hence, by  \cite[Equation (5.1)]{Agg25a}, we have with probability at least $1 - c^{-1} e^{-c(\log N)^2}$ that $u_{i;N} (k;t)^2 \le e^{-(\log N)^2}$ and $u_{i';N'} (k;t)^2 \le e^{-(\log N)^2}$, establishing the lemma. 			
		\end{proof}

\bibliographystyle{abbrv}
\bibliography{bibliotek}

\end{document}